Étienne Ghys

# A singular mathematical promenade

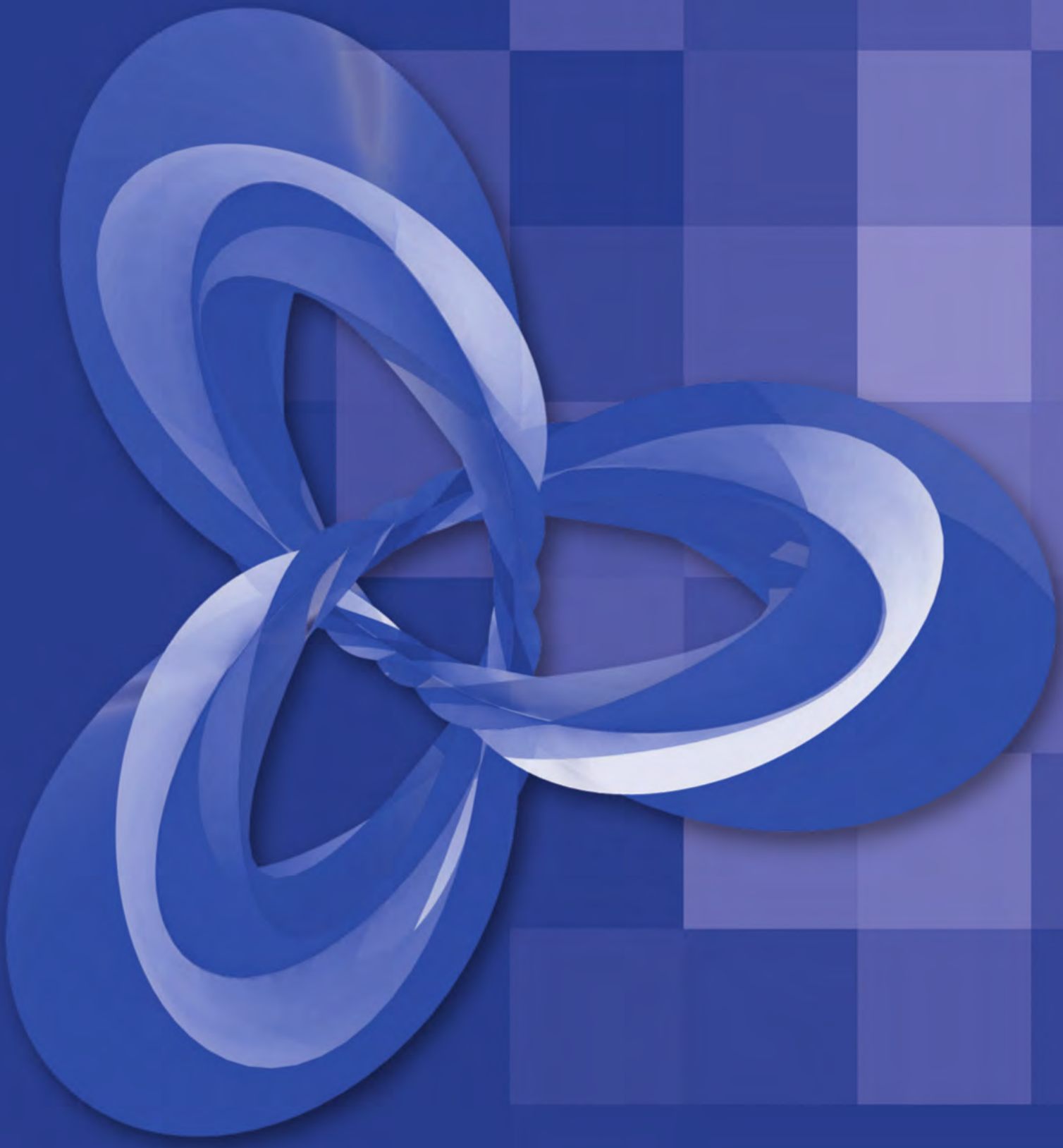

A singular mathematical promenade

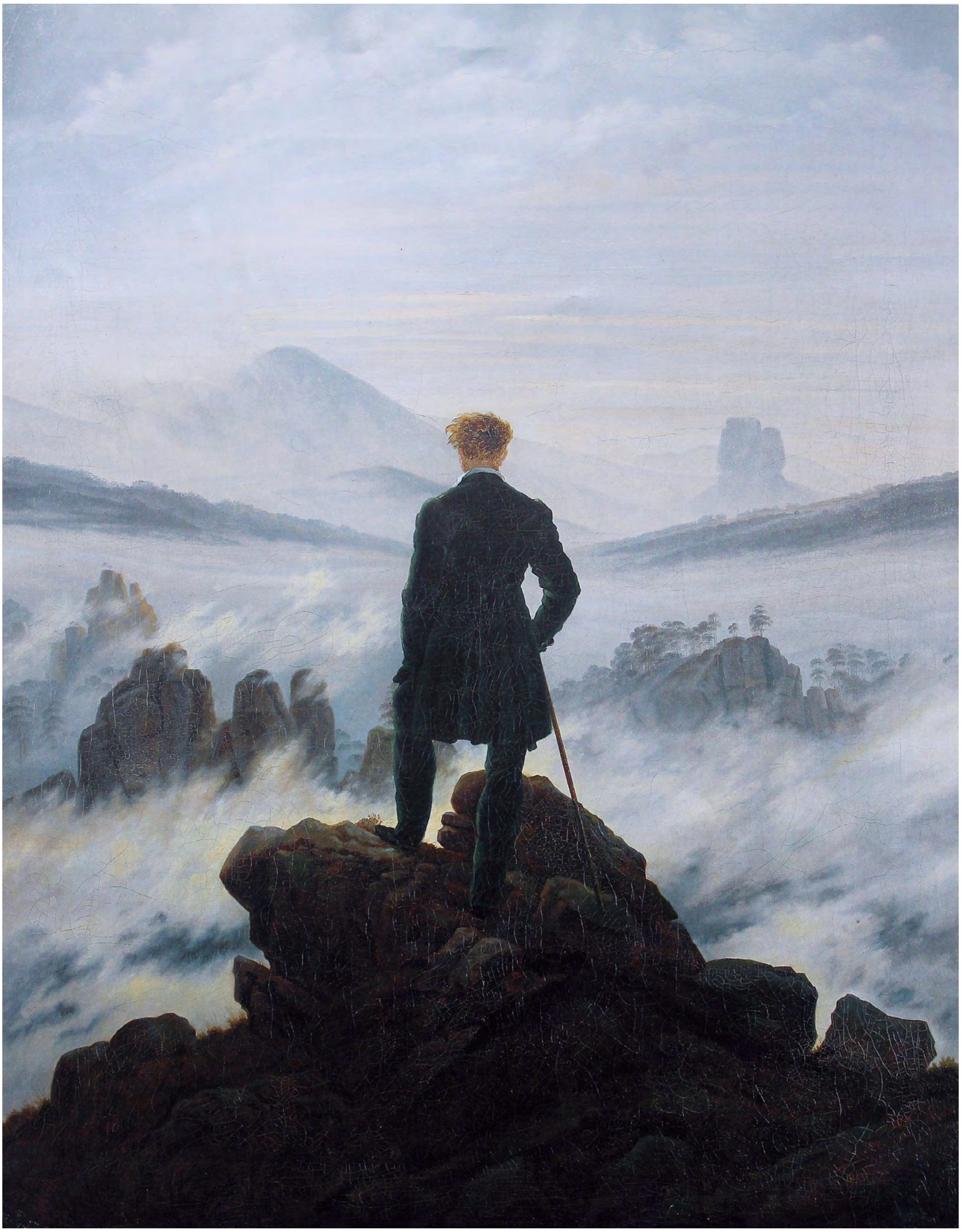

Caspar David Friedrich
"Wanderer above
the sea of fog" (1818) ©

ÉTIENNE GHYS

# A SINGULAR MATHEMATICAL PROMENADE

Contact address:

Étienne Ghys
CNRS-ENS Lyon, UMPA
46 Allée d'Italie
69364 Lyon
France
etienne.ghys@ens-lyon.fr
perso.ens-lyon.fr/ghys

ENS Éditions 2017
École Normale Supérieure de Lyon
15 parvis René Descartes
BP 7000
69342 Lyon cedex 07
France
www.ens-lyon.fr/editions/catalogue

*Pour Martine, qui m'a toujours encouragé à écrire "un petit livre".*

# Contents

# *Preface*

In March 2009, I attended an administrative meeting and the colleague sitting next to me was even more bored than I was. Obviously Maxim Kontsevich had something else in his mind. Suddenly, he passed me a Parisian métro ticket containing a scribble and a single word: “impossible”. That was the new theorem he wanted to share with me! It took me a few minutes and some whispering before I could guess the statement of the theorem and a few more minutes to find the proof. Here is the statement.

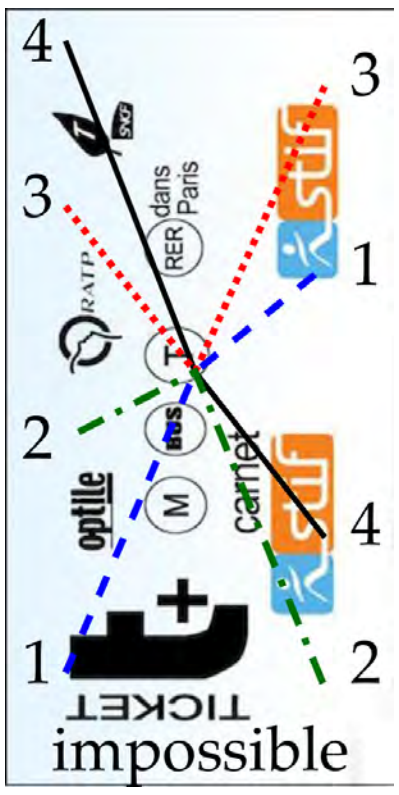

**Theorem.** *Four polynomials $P_1, P_2, P_3, P_4$ of a real variable $x$ cannot satisfy*

- $P_1(x) < P_2(x) < P_3(x) < P_4(x)$ *for small* $x < 0$,
- $P_2(x) < P_4(x) < P_1(x) < P_3(x)$ *for small* $x > 0$.

The *relative position* of the graphs of four real polynomials is subject to some constraints. I was fascinated: a new elementary result on four polynomials in 2009!

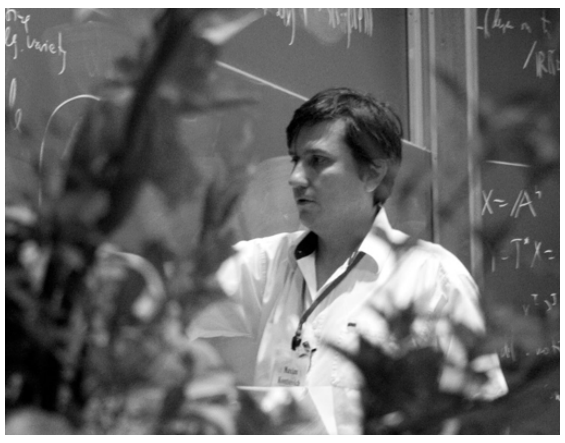

Later on, I tried to put this in a more general context, to study the situation when we have more than four polynomials etc. The result was a pleasant journey, with a lot of detours, in surprisingly different mathematical fields, in different periods of the history of mathematics. As usual, this led to open problems that I could only solve partially.

The purpose of this little book is to invite the reader on this mathematical promenade. I didn’t choose the most efficient way to reach a specific goal and actually there is no goal to this text. Almost all chapters are basically independent and you are welcome to skip as many of them as you want. If you find a section too arduous, or too flat, you can easily bypass it. We pay a visit to Hipparchus, Newton and Gauss, but also to many contemporary mathematicians. We play with a bit of algebra, topology, geometry, complex analysis, combinatorics, and computer science. A stroll in the mathematical world.

However, in order to reach some kind of goal and not to transform this promenade into a completely random walk, let me quote a result that will be proved in one of the last chapters. This is probably the only new result in this work.

Let us consider a point $p$ on a planar curve $\mathcal{C}$.

If $\mathcal{C}$ is *smooth*, the local picture is not so interesting.

If $\mathcal{C}$ is *singular* at $p$, the picture might be more complicated, like for instance a cuspidal point $x^2 = y^3$. Let us restrict our study to *algebraic* curves defined by some equation $F(x,y) = 0$ where $F$ is a polynomial in $x, y$ with real coefficients. It turns out that, in the neighborhood of one of its points, such a curve is the union of a finite number of irreducible pieces, usually called *branches*. The nature of these branches has been the subject of a lot of debate in the past, and we discuss this topic in detail. The main result is that *topologically* the branches are smooth! More precisely for every branch, there is a local *homeomorphism* of the plane sending it to a straight line. Every branch intersects a small circle centered at $p$ in exactly two points.

The relative position of the many branches of a curve is much more subtle. In the neighborhood of a singular point, the topology is described by an even number of points on a circle paired two by two: the pairing is given by the branches. We get $2n$ points on a circle grouped in $n$ pairs, each pair having a color, or a letter.

I can now state a theorem that will be more or less our final destination, some kind of lighthouse showing some direction.

**Theorem.** *There is no singularity of a real algebraic curve in the plane consisting of five branches $A, B, C, D, E$ intersecting a small circle as in the picture in the margin.*

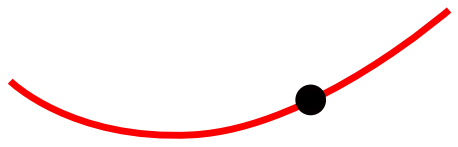

A smooth curve.

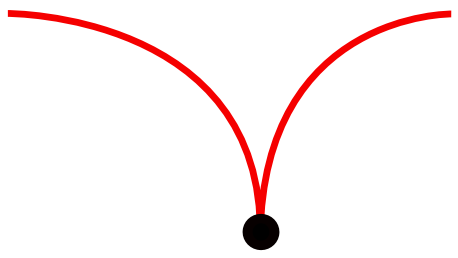

A cuspidal point.

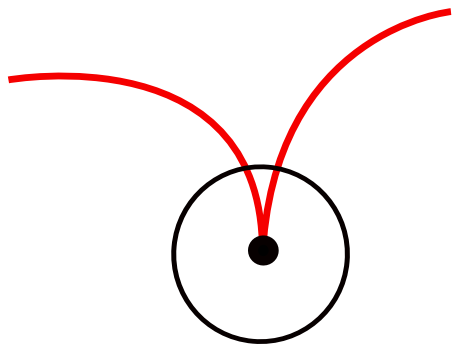

A small circle intersects the curve in two points.

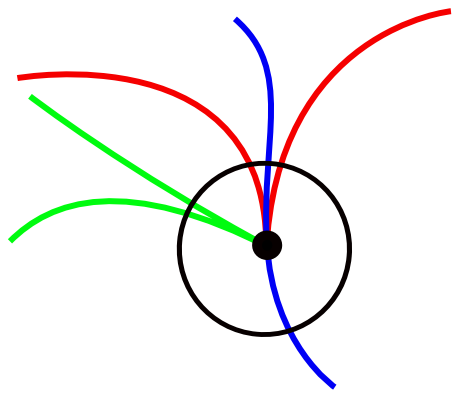

A curve with three branches.

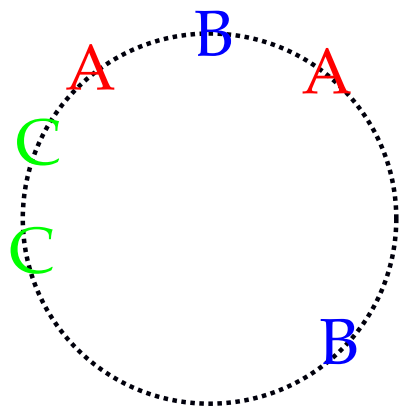

The associated cyclic word $ABACCB$.

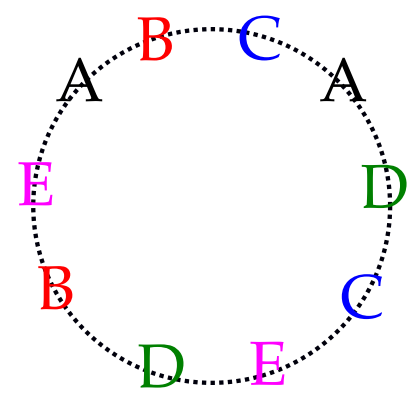

Impossible five branches.

Actually, I will prove a much more precise theorem giving a complete description of all possible topological configurations of the branches of an analytic curve.

I wrote this "petit livre" with one specific reader in mind: myself, when I was an undergraduate... To be very specific, I limited the prerequisites to my own background when I passed the "agrégation" examination, exactly forty years ago ☺! I vividly remember that I had (and I still have) great difficulties reading long mathematical treatises, full of technical details, and that I preferred looking at pictures. I have now learned that precision and details are frequently necessary in mathematics, but I am still very fond of promenades. I did try to imagine what could have been my own reactions faced with this book, as a beginner. This "conversation" between the two "versions of myself" has been interesting and reminded me of the short story "El Otro" by Borges. Was it a dream? A reconstruction of the past?

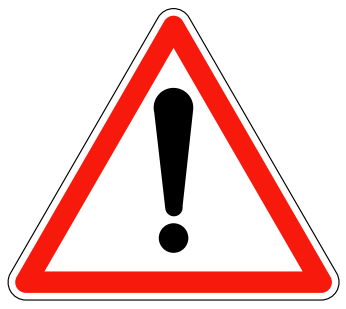

©

A word of caution is in order: this is not a fully fledged textbook with a definition-theorem-proof structure. You have to be prepared to get lost from time to time, like in many promenades. I know that you will grumble about me because of the lack of precise definitions, and indeed you will have to accept half baked definitions... Of course, textbooks are necessary and I provide many references in the margins. However, I am convinced that mathematical ideas and examples precede formal proofs and definitions. As d'Alembert said once: "*Just go on... and faith will catch up with you!*". You may see every now and then a beautiful panorama emerging from the mist, like the one on the frontispiece of this essay, by Caspar David Friedrich (Der Wanderer über dem Nebelmeer (1818)): a suggestion of the mathematical world?

*"Allez en avant, et la foi vous viendra."*

I hope some motivated undergraduates of today will enjoy some of these panoramas.

We can now embark for our voyage.

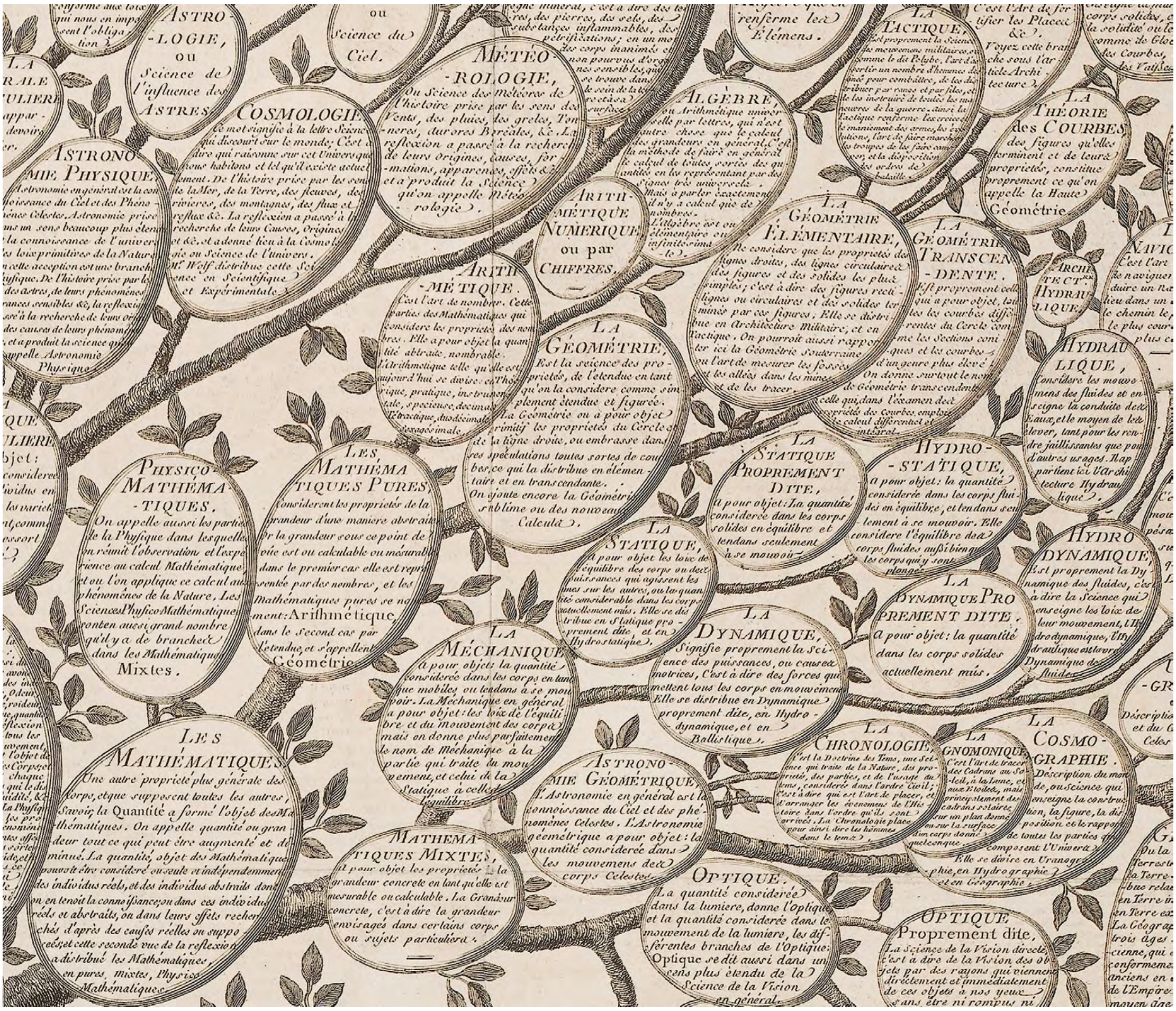

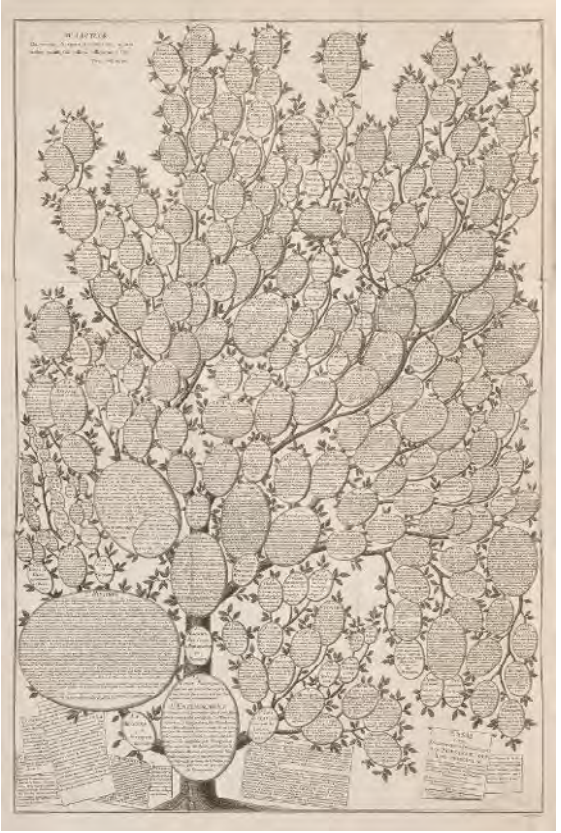

A detail from "Essai d'une distribution généalogique des sciences et des arts principaux" (Chrétien Fréderic Guillaume Roth, 1769). It was included as a frontispiece of the famous Encyclopédie by Diderot and D'Alembert. "Mathematics" are located in the lower left corner and the "theory of curves" is in the right upper corner. Is human knowledge organized as a tree?

# *Road map*

SINCE WE WILL DEFINITELY NOT FOLLOW A DIRECT ROUTE, and since you should be prepared for some optional detours, a rough outline of our itinerary might be useful, like in the promotional presentation of a touristic package by a tour operator.

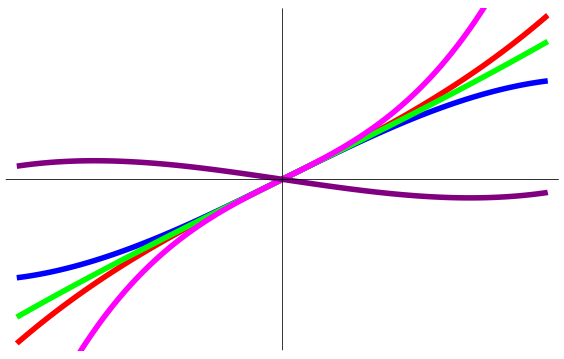

A permutation defined by 5 polynomials.

The first four chapters discuss the *relative positions of the graphs of a family of real polynomials* $P_1, ..., P_n$, in the spirit of the theorem of Kontsevich that I mentioned in the preface. Comparing the values of $P_i(x)$ for small negative and small positive values of $x$ yields a permutation of $\{1, ..., n\}$ which describes the local picture in the neighborhood of 0. I will give a fairly precise characterization of these permutations. It turns out that they have already been considered in a different disguise by combinatorists, under the name "*separable permutations*". We then examine *push and pop stacks*, as presented by Donald Knuth in *The art of computer programming*. We also *count* the number of separable permutations, and this will be an opportunity to discover that these numbers have already been considered by Hipparchus, more than two millenniums ago.

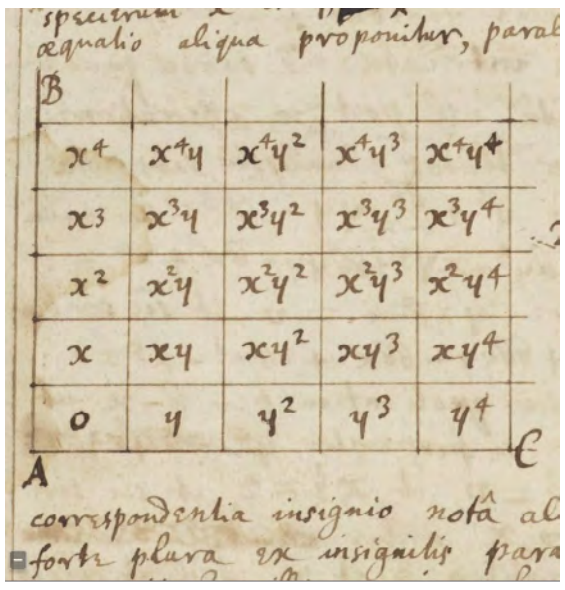

| | | | | |
|---|---|---|---|---|
| $x^4$ | $x^4y$ | $x^4y^2$ | $x^4y^3$ | $x^4y^4$ |
| $x^3$ | $x^3y$ | $x^3y^2$ | $x^3y^3$ | $x^3y^4$ |
| $x^2$ | $x^2y$ | $x^2y^2$ | $x^2y^3$ | $x^2y^4$ |
| $x$ | $xy$ | $xy^2$ | $xy^3$ | $xy^4$ |
| 0 | $y$ | $y^2$ | $y^3$ | $y^4$ |

From Newton's *de methodis*. ©

We then try to generalize the problem from graphs of polynomials to *planar curves*, implicitly defined by some real polynomial equation $F(x, y) = 0$. This requires the understanding of the topology of an algebraic (or analytic) curve in the neighborhood of a singular point. The first important results are due to Newton in 1669, in an extraordinary paper entitled *Tractatus de methodis serierum et fluxionum*, that we study over two chapters. This paper contains a detailed presentation of the famous

*Newton's method* for finding approximations of the roots of polynomials. It also introduces the related idea of *Newton's polygons*. Strictly speaking, Newton did not provide proofs, but he did understand that locally an analytic curve consists of a finite number of *branches*, which are "graphs" of *formal power series with rational exponents*. An additional chapter — that I called formal algebra — explains Newton's results in modern terminology and offers proofs.

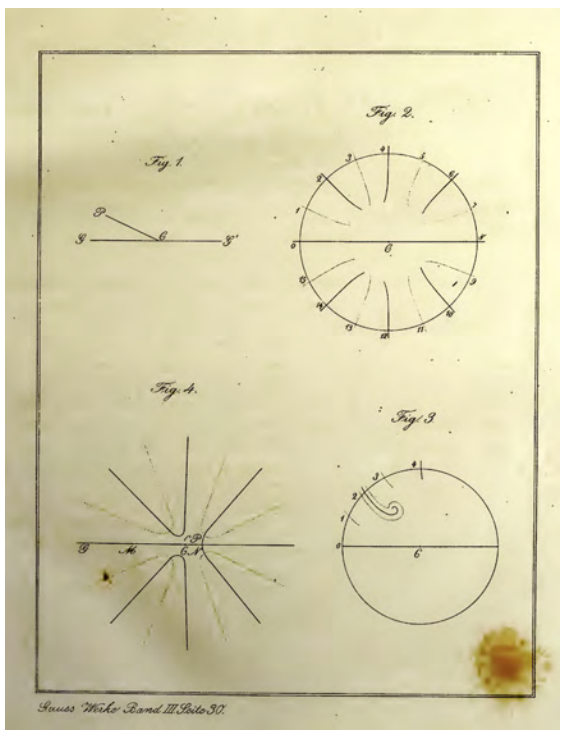

From Gauss's doctoral dissertation. ©

Up to this point, the discussion will be purely algebraic. We then review Gauss's first proof of the fundamental theorem of algebra — his doctoral dissertation in 1799 — using arguments of topological nature, which were revolutionary at that time. This is based on the *unproved claim that an algebraic curve entering a disc has to get out*. The proof of this claim is more subtle than one could imagine and two mathematicians sharing the same name could not prove it in the 19th century.

Euler, Cauchy and Poincaré were great masters in the manipulation of series. Two chapters deal with their discoveries. At the end of the second one, using the Calcul des limites de Cauchy, we finally get the proof of the *convergence of Newton's series*. This enables us to show that a small circle around a singularity of a plane real analytic curve intersects the curve in an even number of points and defines a *chord diagram*, i.e. $2n$ points cyclically ordered on a circle and grouped in pairs.

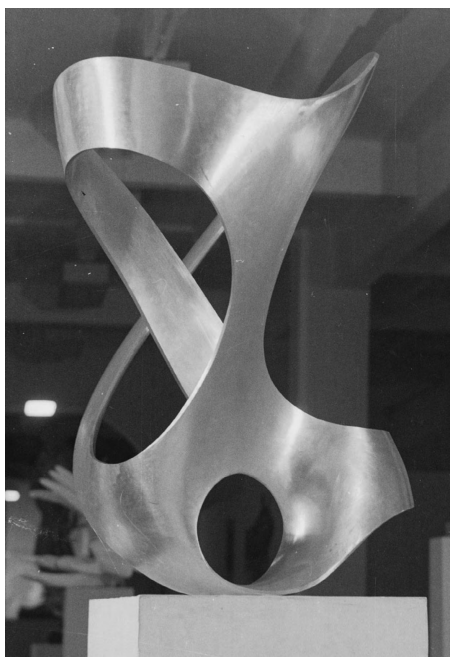

A disc blown-up three times. ©

The three following chapters are concerned with the topology of singularities of analytic planar curves. We explain the *blowing up* method, which is a kind of microscope enabling us to look deeply into the singularity. Topologically, this introduces a Moebius band, or Moebius necklaces if the microscope is used several times. The blowing up operation will turn out to be a powerful tool in the resolution of singularities.

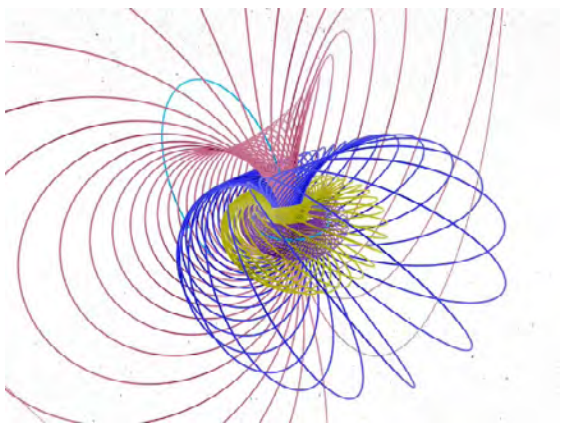

Hopf fibers. ©

The local pictures for *complex* planar curves are beautiful and worth a visit. Since $\mathbb{C}^2$ has real dimension 4, we intersect the curve with small 3-dimensional spheres around the singularity. From this viewpoint, even straight lines produce remarkable objects, like the Hopf fibration.

More complicated singularities, like for example the cusp $x^3 - y^2 = 0$ are described by *knots* and *links*. In order to understand the general case, we pay a visit to Victor Puiseux, who proposed in 1850 a completely new approach to Newton's series. In 1968, Jack Milnor used these ideas to give a complete topological picture, but still over the *complex numbers*.

Interestingly, we discover a link between separable permutations and the associahedron. This is a family of convex polytopes introduced by Tamari and Stasheff in order to understand the meaning of "associativity up to homotopy". Using his polytopes, Jim Stasheff was able to give a characterization of spaces having the same homotopy type as topological groups. It turned out that this was the starting point of operad theory, which plays a fundamental role in modern homotopy theory and algebraic topology. Operads are very general algebraic structures and they are perfectly adapted to our situation. Typical examples are given by trees, braids, configuration spaces etc. We will see that the collection of all singularities, up to homeomorphisms, can be seen as a singular operad and this helps understanding the global picture.

Just for fun, we examine a short note of Gauss, concerning closed loops in the plane, with ordinary double points. Going around the loop, each double point is visited twice, so that this defines some chord diagram. Can we characterize this kind of diagrams?

We finally reach our loose goal: the complete characterization, in two chapters, of the chord diagrams which are associated to singularities of real analytic planar curves.

Two additional chapters conclude the book. One on Gauss's approach to linking numbers and a final one, with no proof, on Kontsevich's universal invariant for knots. The main purpose of this final chapter is to encourage the reader to continue the exploration.

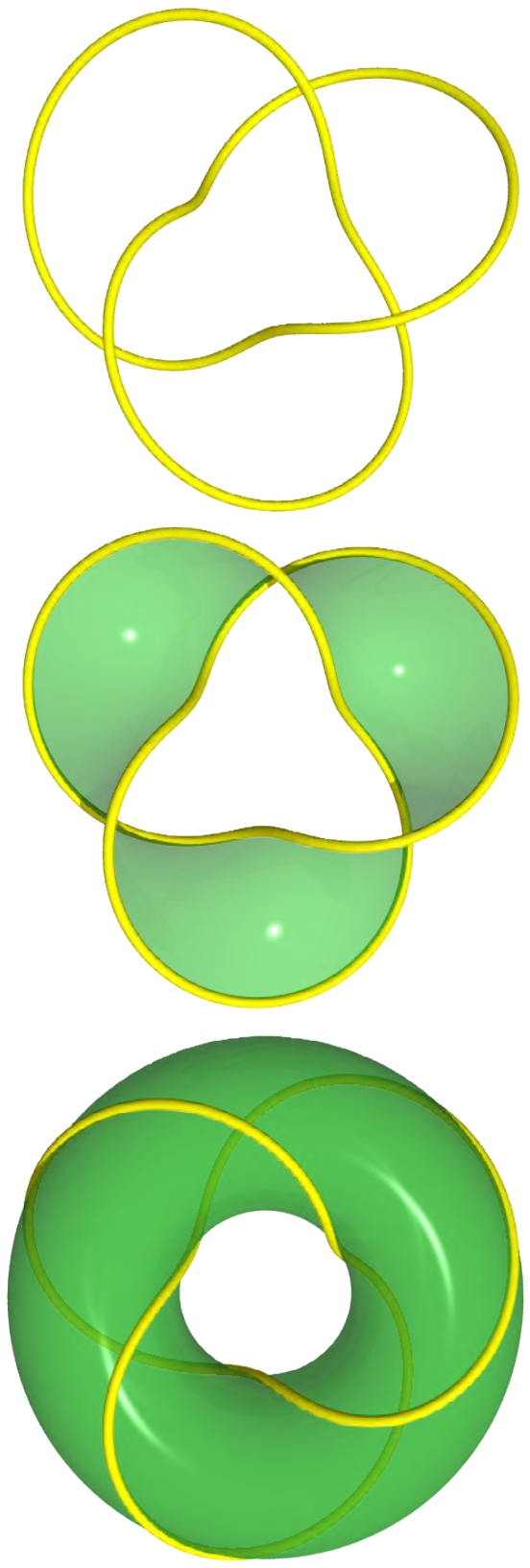

The trefoil knot. ©

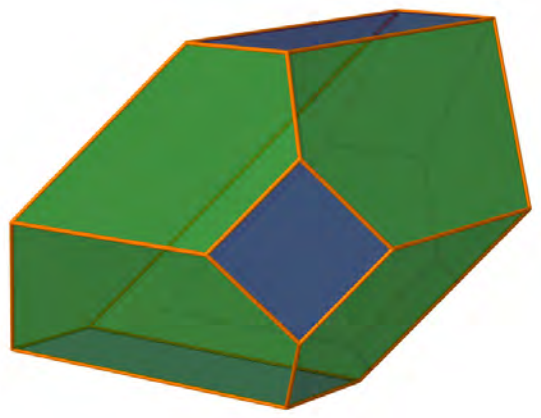

The associahedron. ©

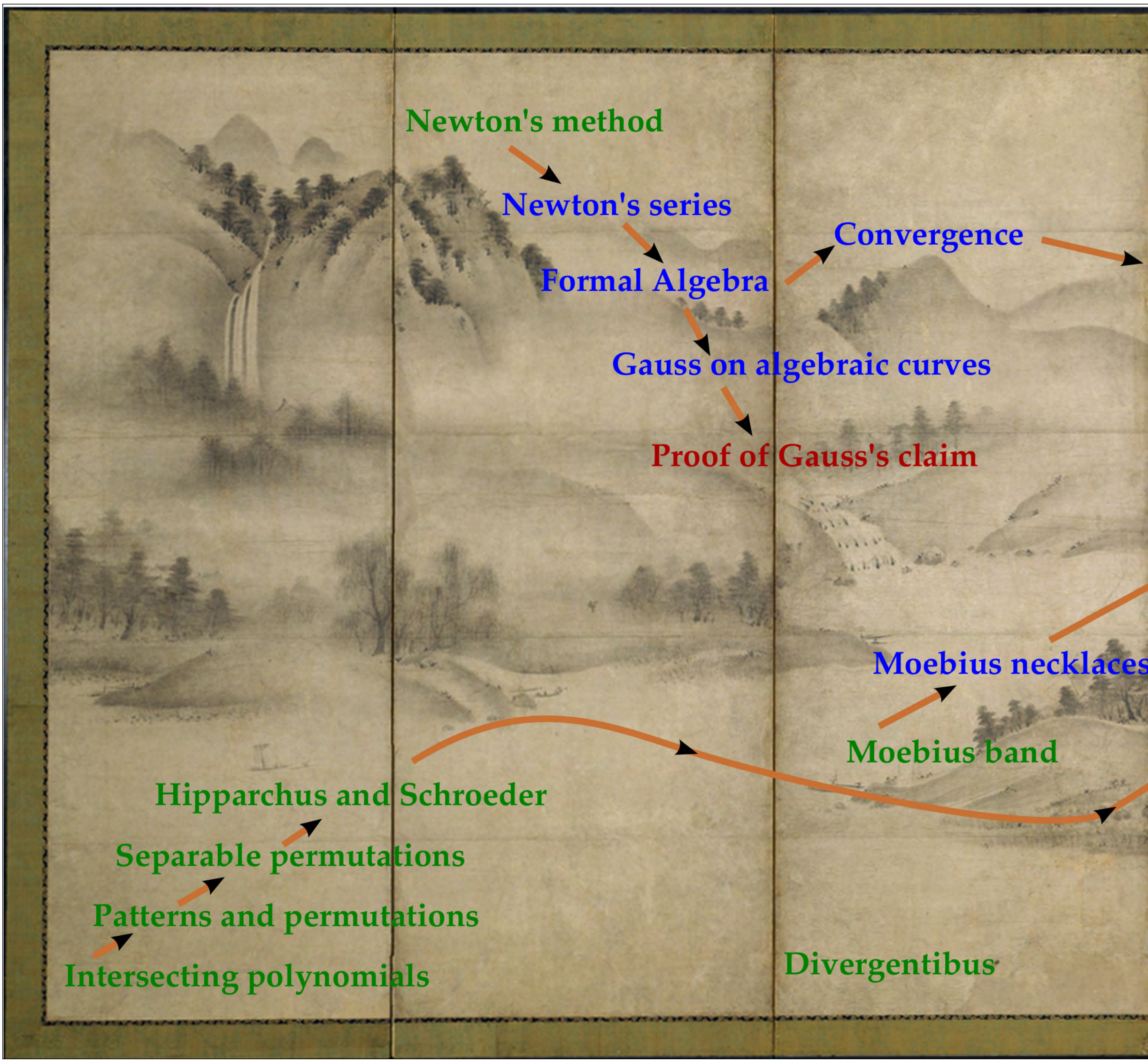

Colors (green, blue, red, and black) give a very subjective idea of the difficulty.

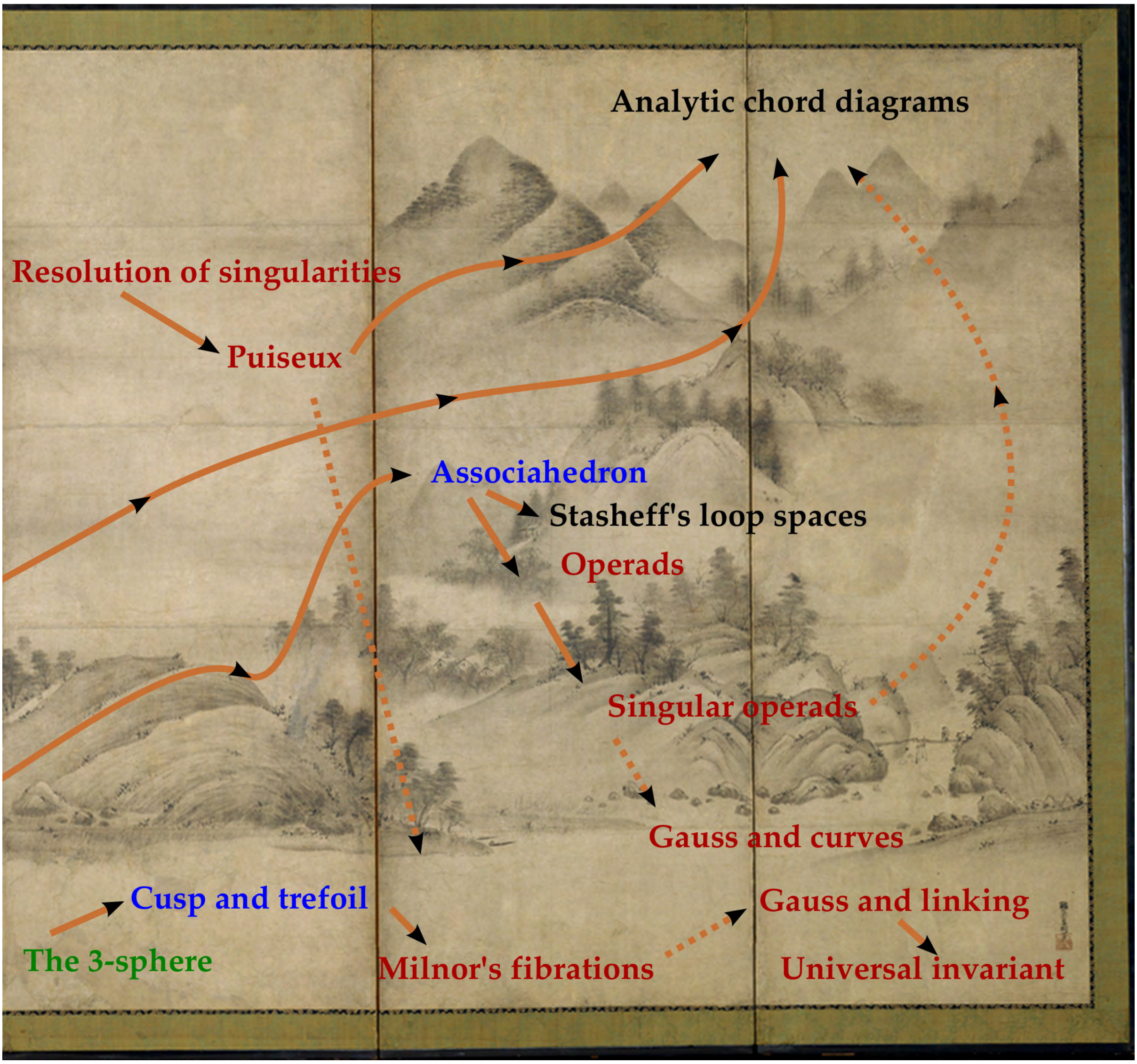

Landscape of the Four Seasons (Eight Views of the Xiao and Xiang Rivers), by Soami, early 16th century. ©

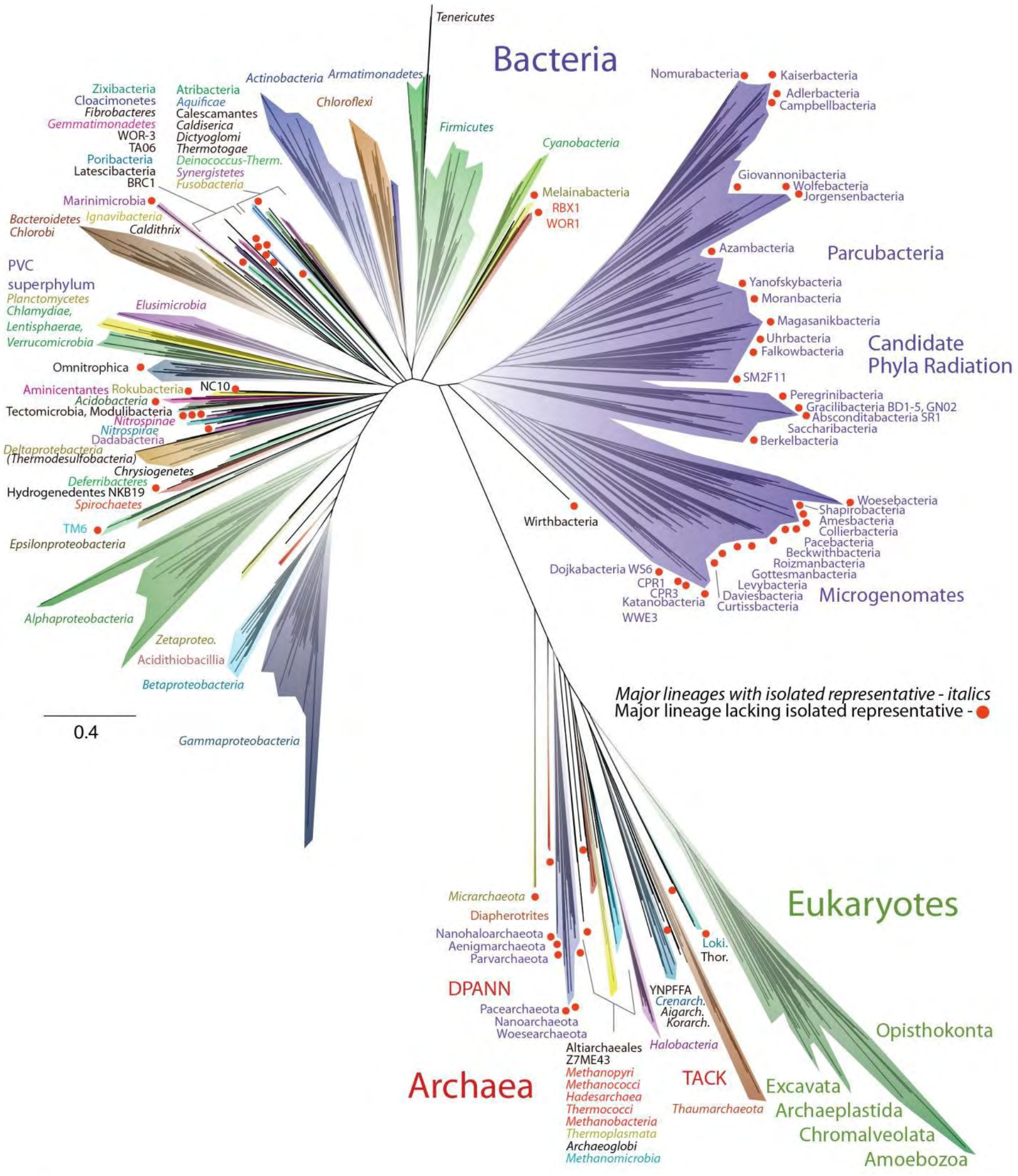

From “A new view of the tree of life” Nature Microbiology 1, (2016).
Can these branches be made graphs of polynomials?

# *Intersecting polynomials*
# *Maxim Kontsevich*

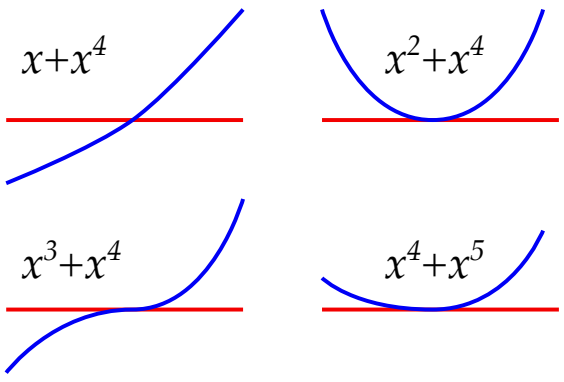

## *Polynomial interchanges*

BEFORE I PROVE KONTSEVICH'S THEOREM, let me begin with a much more elementary observation. Consider the position of the graph of a single real nonzero polynomial $P(x)$ with respect to the $x$-axis, in the neighborhood of 0.

There are two possibilities. Either the graph of $P$ crosses the $x$-axis at 0, or it stays on the same side. To distinguish between these two cases, I introduce the following definition.

**Definition.** Let $P(x) = a_0 + a_1x + a_2x^2 + \cdots$ be a polynomial (or a formal series). The *valuation* $v(P)$ of $P$ (at 0) is the smallest integer $k$ such that $a_k \neq 0$. By convention, the valuation of the zero polynomial is $\infty$.

Clearly, the graph of $P$ crosses the $x$-axis at 0 if and only if the valuation $v(P)$ is an odd integer.

If we are given *two* distinct polynomials $P_1, P_2$, the sign of $P_1(x) - P_2(x)$ changes at 0 if and only if $v(P_1 - P_2)$ is odd.

Suppose now that we have *three* polynomials $P_1, P_2, P_3$ and let us look at the possible pictures in the neighborhood of the origin. The six figures in the margin show that all six permutations of $\{1, 2, 3\}$ can be realized if we choose conveniently the polynomials.

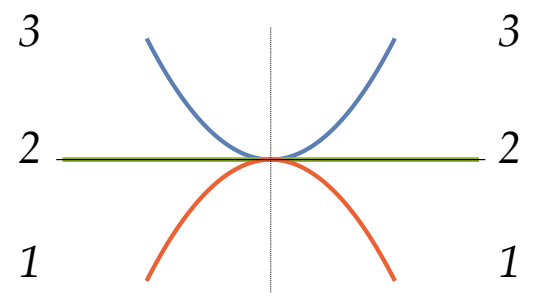

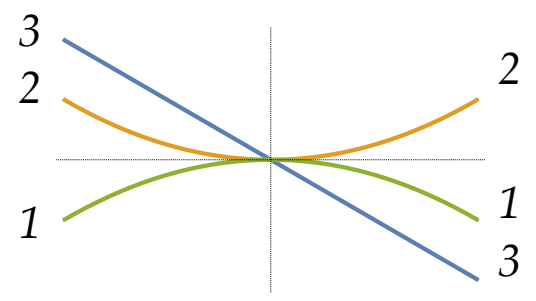

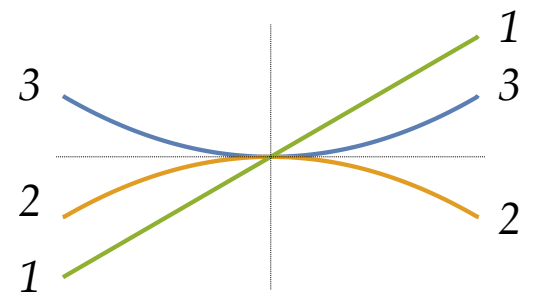

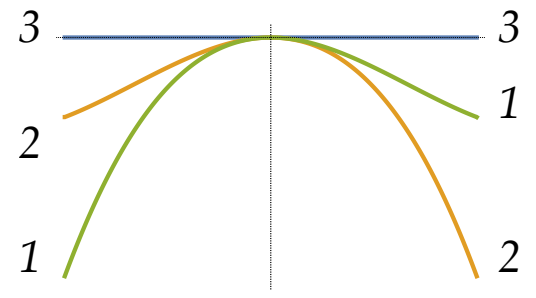

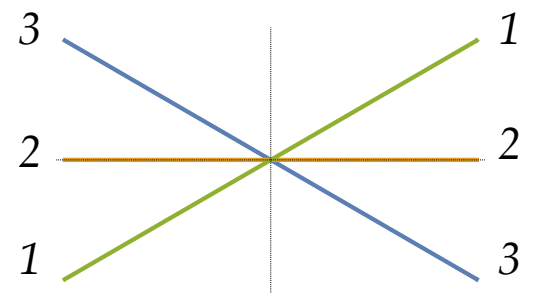

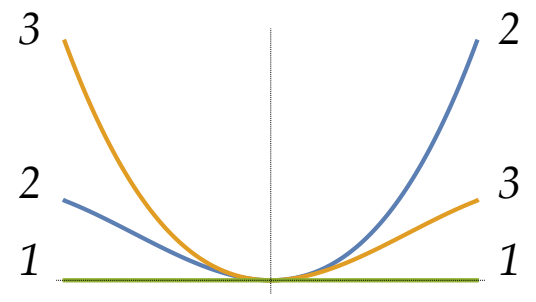

For instance:

$$\begin{array}{lclcl}
P_1(x) = -x^2 & ; & P_2(x) = 0 & ; & P_3(x) = x^2 \\
P_1(x) = -x^2 & ; & P_2(x) = x^2 & ; & P_3(x) = -x \\
P_1(x) = x & ; & P_2(x) = -x^2 & ; & P_3(x) = x^2 \\
P_1(x) = -x^2 + x^3 & ; & P_2(x) = -x^2 - x^3 & ; & P_3(x) = 0 \\
P_1(x) = x & ; & P_2(x) = 0 & ; & P_3(x) = -x \\
P_1(x) = 0 & ; & P_2(x) = x^2 + x^3 & ; & P_3(x) = x^2 - x^3.
\end{array}$$

Hence Kontsevich's phenomenon begins with *four* polynomials.

Note that all the previous pictures may have suggested that I assumed $P_i(0) = 0$ but this is not necessary. This is only due to the fact that this book mainly discusses local properties, in the neighborhood of a single point $(0,0)$.

The polynomials $x^2P_1(x), \ldots, x^2P_n(x)$ are ordered as $P_1(x), \ldots, P_n(x)$ and they all vanish at $x = 0$.

We can now prove the métro ticket theorem mentioned in the preface.

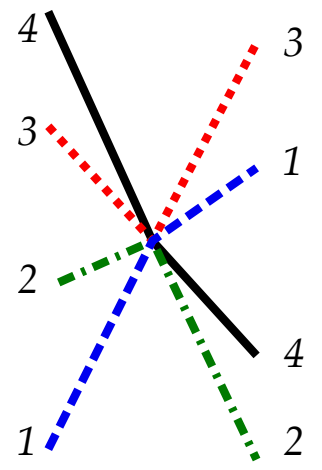

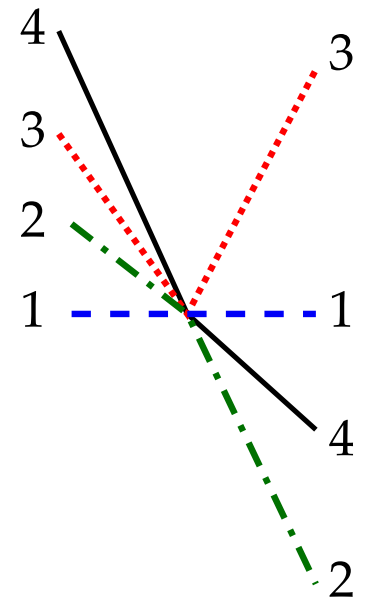

By contradiction, we assume that there exist four polynomials $P_1, P_2, P_3, P_4$ such that:

1. $P_1(x) < P_2(x) < P_3(x) < P_4(x)$ for small $x < 0$,
2. $P_2(x) < P_4(x) < P_1(x) < P_3(x)$ for small $x > 0$.

Replacing $P_i$ by $P_i - P_1$, we can assume that $P_1 = 0$.

Since $P_2$ and $P_4$ change sign at the origin, their valuations $v(P_2), v(P_4)$ are odd.

Since $P_3$ does not change sign, its valuation $v(P_3)$ is even.

From $0 < P_2(x) < P_3(x) < P_4(x)$ for small negative $x$, we deduce that $v(P_2) \geq v(P_3) \geq v(P_4)$.

Similarly, from $|P_4(x)| < |P_2(x)|$ for small positive $x$, we deduce that $v(P_4) \geq v(P_2)$.

That would force the three valuations to be equal, but two of them are odd and the third is even!

Contradiction. ⊡

I use the symbol ⊡ at the end of a proof. Will my astute reader guess why I put a · in a □? Hint: think in French.

Note that the same proof applies to real analytic functions but does not apply to smooth functions. Indeed my reader will easily find four $C^\infty$ functions $P_i$ crossing at the origin according to the "forbidden" permutation.

Why?

Changing orientations along the $x$-axis, we see that the inverse permutation is also forbidden. As an exercise, I recommend showing that the remaining 22 permutations of $\{1,2,3,4\}$ occur for suitable choices of the $P_i$'s ($i = 1,2,3,4$).

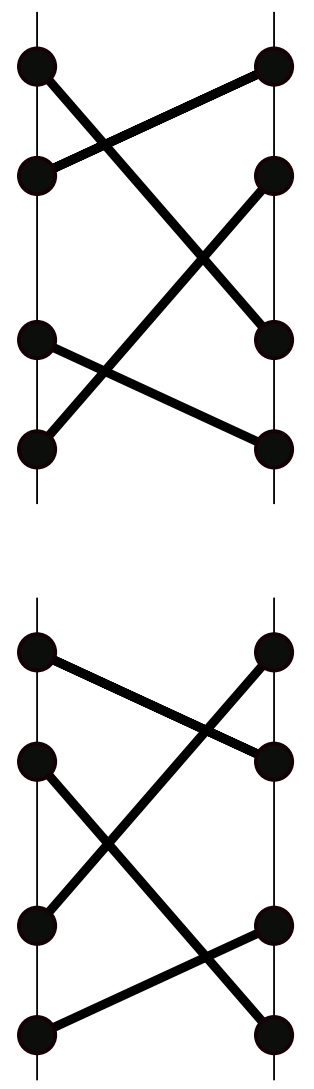

The forbidden permutations.

Let us now try to analyze the situation with any number of polynomials.

**Definition.** Let $n \geq 2$ be some integer and $\pi$ some permutation of $\{1,2,\ldots,n\}$. We say that $\pi$ is a *polynomial interchange* if there exist $n$ polynomials $P_1,\ldots,P_n$ such that:

- $P_1(x) < P_2(x) < \ldots < P_n(x)$ for small negative $x$.
- $P_{\pi(1)}(x) < P_{\pi(2)}(x) < \ldots < P_{\pi(n)}(x)$ for small positive $x$.

Our goal is to give a fairly precise description of polynomial interchanges.

## *Separable permutations*

**Definition.** Let $n \geq 2$ be some integer and $\pi$ some permutation of $\{1,2,\ldots,n\}$. We say that $\pi$ is *separable* if it does not "contain" one of the two forbidden permutations, i.e. if there do not exist four indices $1 \leq i_1 < i_2 < i_3 < i_4 \leq n$ such that $\pi(i_2) < \pi(i_4) < \pi(i_1) < \pi(i_3)$ or $\pi(i_3) < \pi(i_1) < \pi(i_4) < \pi(i_2)$.

The reason for the terminology "separable" will become clear in the next chapter.

In other words, a permutation is separable if it does not contain one of the two Kontsevich's permutations on four letters. It should be clear that a polynomial interchange is necessarily separable. In this section, we prove the converse.

Let us begin with a lemma which seems to be "folklore" in the combinatorics literature[1].

[1] P. Bose, J. F. Buss, and A. Lubiw. Pattern matching for permutations. *Inform. Process. Lett.*, 65(5):277–283, 1998.

**Lemma.** *Let $\pi$ be some separable permutation of $\{1,2,\ldots,n\}$ (for $n \geq 3$). Then there is a* proper *interval $I = \{k, k+1, \ldots, k+l\}$ (with $1 \leq k \leq k+l \leq n$) of length $l+1 \geq 2$ whose image by $\pi$ is an interval.*

We can assume that $\pi(1) < \pi(2)$ since otherwise we could replace $\pi$ by the "reverse" permutation $\overline{\pi}(k) = n+1-\pi(k)$.

Observe that if $\pi$ is a polynomial interchange, so is $\overline{\pi}$ (multiply all polynomials by $-x$). Similarly, by the very definition of separable permutations, $\pi$ and $\overline{\pi}$ are simultaneously separable.

If $\pi(2) = \pi(1)+1$ we are done since the image of $\{1,2\}$ is the interval $\{\pi(1), \pi(2)\}$. Hence we assume that $\pi(2) > \pi(1)+1$.

Consider the smallest integer $k$ such that $\pi(\{2,\ldots,k\})$ contains the interval $J$ between $\pi(1)+1$ and $\pi(2)$. Observe that $\pi(k)$ is in $J$ so that $\pi(1) < \pi(k) < \pi(2)$.

If the image $\pi(\{2,\ldots,k\})$ is exactly equal to the interval $J$, we found a nontrivial interval whose image by $\pi$ is an interval.

Otherwise, choose an element $l$ between 2 and $k$ whose image by $\pi$ is "outside" $J$. We have $\pi(l) < \pi(1)$ or $\pi(l) > \pi(2)$.

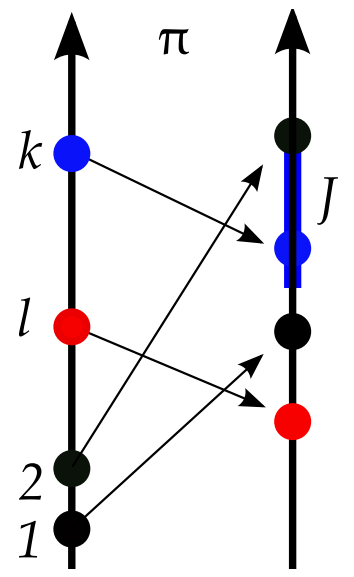

If $\pi(l) < \pi(1)$, the four elements $1, 2, l, k$ satisfy $1 < 2 < l < k$ and $\pi(l) < \pi(1) < \pi(k) < \pi(2)$. Therefore they are ordered as a "forbidden permutation" which is impossible, by definition of a separable permutation.

We can therefore assume that all elements of $\pi(\{2,\ldots,k\})$ are greater than or equal to $\pi(1)$.

We can also assume that $\pi(\{2,\ldots,k\})$ is not an interval since otherwise we are done. Therefore there is at least one "gap" in $\pi(\{2,\ldots,k\})$, which must be greater than $\pi(2)$. So there exists $m$ such that $k < m$ and $l$ such that $2 < l < k$ and $\pi(m) < \pi(l)$. The four elements $2, l, k, m$ are such that $2 < l < k < m$ and $\pi(k) < \pi(2) < \pi(m) < \pi(l)$ so that they are ordered as the other "forbidden permutation" which is impossible.

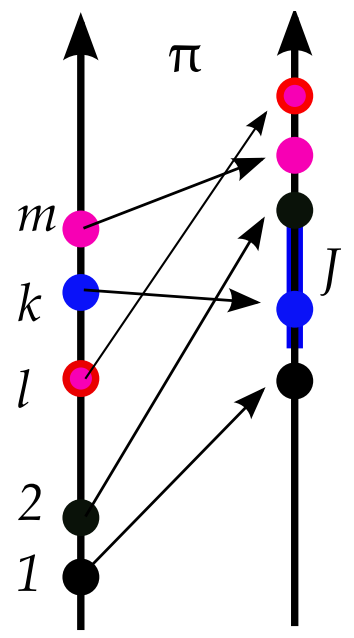

The lemma is proved. ⊡

It is easy to improve the lemma:

**Lemma.** *Let $\pi$ be some separable permutation of $\{1, 2, \ldots, n\}$. Then there are two consecutive integers whose images are consecutive.*

The proof is obvious by induction since any proper interval whose image is an interval defines another separable permutation with a smaller value of $n$, which therefore contains two consecutive elements with consecutive images. ⊡

We can now prove the main result of this chapter. The Kontsevich counter-example is somehow the only one.

**Theorem.** *A permutation is a polynomial interchange if and only if it is separable.*

We have already noticed that polynomial interchanges are separable: this is Kontsevich's observation.

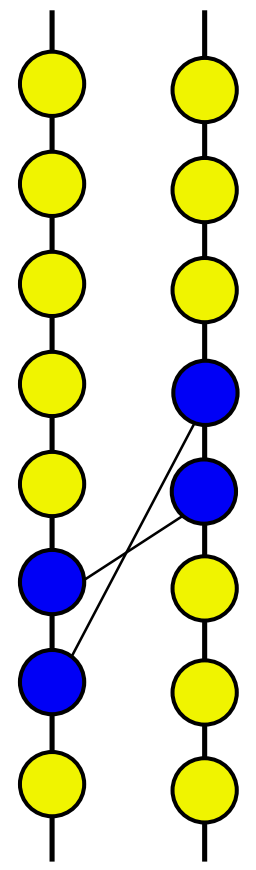

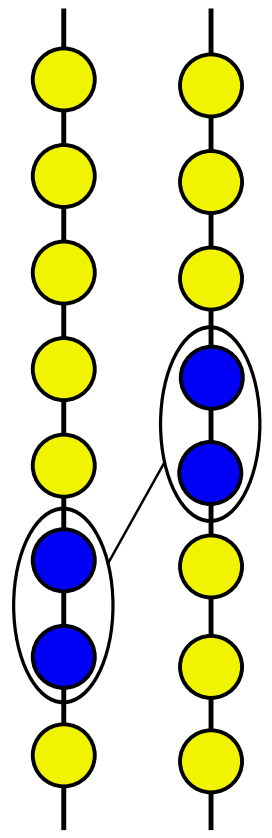

Again by induction on $n$, we show that every separable permutation is a polynomial interchange. Let $\pi$ be a separable permutation of $\{1,2,\ldots,n\}$. We know that there are two consecutive integers $i, i+1$ with consecutive images $\pi(i), \pi(i+1)$.

If $\{i, i+1\}$ and $\{\pi(i), \pi(i+1)\}$ are "collapsed" into single points, we produce a permutation $\pi'$ on $n-1$ objects which is obviously separable, and therefore a polynomial interchange by induction. It follows that there are $n-1$ polynomials

$$P_1, \ldots, P_{n-1}$$

which intersect at the origin according to $\pi'$. The only thing that remains to be done is to split the $i$-th polynomial $P_i$ in order to produce $n$ polynomials

$$P_1, \ldots, P_{i-1}, P_i', P_i'', P_{i+1}, \ldots, P_{n-1}$$

which intersect according to $\pi$. It suffices to set

$$P_i'(x) = P_i(x) \quad ; \quad P_i''(x) = P_i(x) + (-x)^N$$

for a sufficiently large value of $N$, even or odd, according to whether $\pi(i+1) > \pi(i)$ or $\pi(i+1) < \pi(i)$. ⊡

Now that we have identified the polynomial interchanges, our next duty is to understand the structure of those separable permutations.

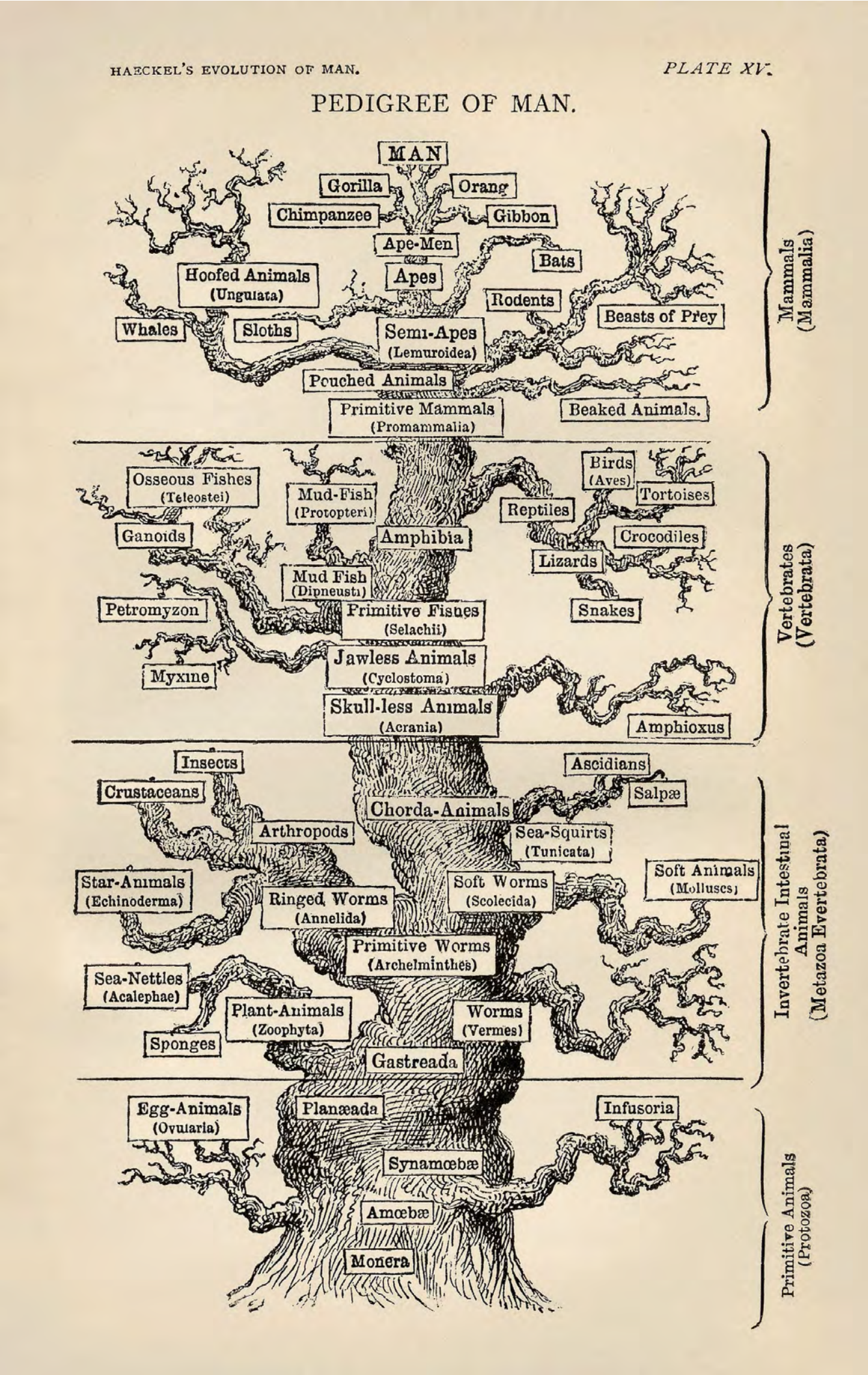

Ernst Haeckel's "tree of life" (1879).
Man on top of the tree of life? ©

# *Patterns and permutations*
# *Donald Knuth*

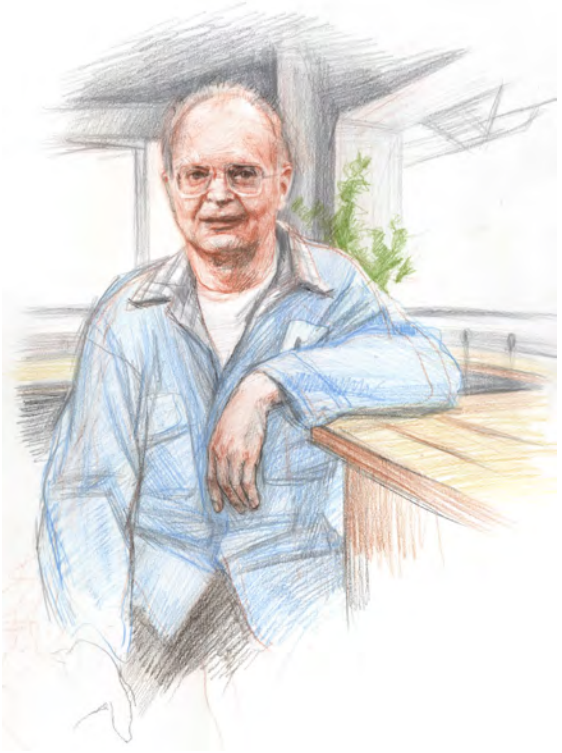

Donald Knuth. ©

## *Permutations*

IF THEY HAVE BEEN MATHEMATICALLY TRAINED AS I WAS, many of my readers may have felt some discomfort in the previous chapter. After all, permutations are usually defined as bijections from a set to *itself* and their *raison d'être* is that they constitute a group. Instead, we manipulated permutations in a strange way when we used the expression "The permutation $\pi$ *contains* one of the two Kontsevich's forbidden permutations" to mean that there are four indices $1 \le i_1 < i_2 < i_3 < i_4 \le n$ such that

$$\pi(i_2) < \pi(i_4) < \pi(i_1) < \pi(i_3), \text{ or}$$
$$\pi(i_3) < \pi(i_1) < \pi(i_4) < \pi(i_2).$$

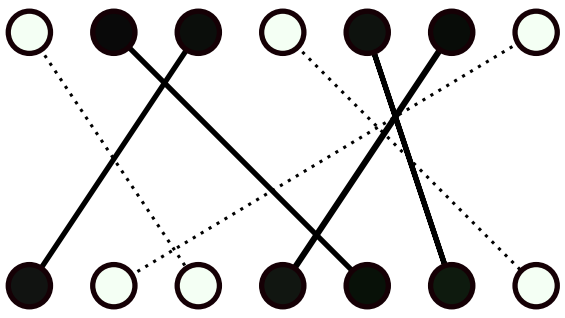

It certainly does not mean that the set $\{i_1, i_2, i_3, i_4\}$ is invariant under $\pi$. We are *not* taking the restriction to an invariant subset.

We are going to use the word "permutation" from a slightly different perspective, closer to computer science. This approach is in good part due to Donald Knuth in his great book *The art of computer programming*[2]. The more recent book[3] is a good source of information and shows that this area is currently blossoming.

Consider a finite set $E$ equipped with two total orderings $\ll$ and $\lll$. Order its elements using the first ordering

$$x_1 \ll x_2 \ll \ldots \ll x_n.$$

[2] D. E. Knuth. *The art of computer programming. Vol. 1: Fundamental algorithms.* Second printing. Addison-Wesley Publishing Co., Reading, Mass.-London-Don Mills, Ont, 1969.

[3] S. Kitaev. *Patterns in permutations and words.* Monographs in Theoretical Computer Science. An EATCS Series. Springer, Heidelberg, 2011. With a foreword by Jeffrey B. Remmel.

Look now at the way they are ordered under $\lll$. This defines a permutation $\pi$ of $\{1, 2, \dots, n\}$ such that

$$x_{\pi(1)} \lll x_{\pi(2)} \lll \dots \lll x_{\pi(n)}.$$

We will adopt this point of view: a permutation is a comparison between two total orderings.

For instance, the set $\{1, 2, 3, 4\}$ can be equipped with the orderings $1 \ll 2 \ll 3 \ll 4$ and $2 \lll 4 \lll 1 \lll 3$; we will denote by $(2, 4, 1, 3)$ the associated permutation.

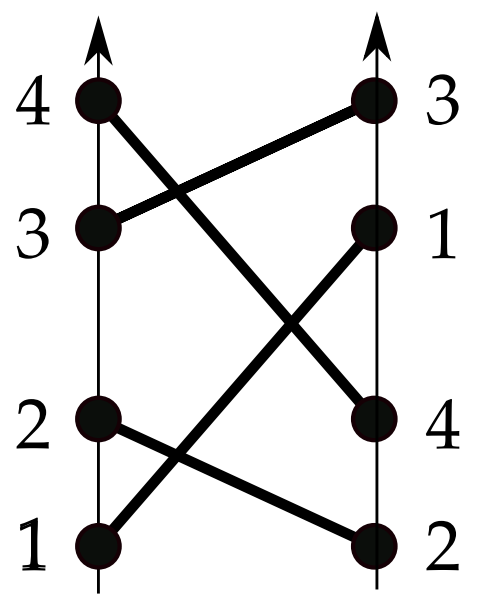

Be careful. In this figure, $1 \ll 2 \ll 3 \ll 4$ and $2 \lll 4 \lll 1 \lll 3$ and $(\pi(1), \pi(2), \pi(3), \pi(4)) = (2, 4, 1, 3)$, so that $\pi$ is actually the *inverse* of the permutation that you see following the edges, from the dots on the left to the dots on the right.

Any finite set of real polynomials $\{P_i(x)\}$ can be ordered in at least two ways: by comparing the values of $P_i(x)$ for small negative or for small positive values of $x$. This leads to polynomial interchanges.

We can certainly *restrict* orderings to subsets and this defines the concept of containment for permutations.

**Definition.** Let $\pi$ be the permutation of $\{1, \dots, n\}$ associated to two total orderings $\ll$ and $\lll$ on a set $E$ with $n$ elements. Let $F \subset E$ be a subset with $p$ elements. The restrictions of $\ll$ and $\lll$ to $F$ define a permutation $\sigma : \{1, \dots, p\} \to \{1, \dots, p\}$. In such a situation, we will say that $\sigma$ is *contained* in $\pi$ and we will write $\sigma \le \pi$.

For instance, every permutation different from the identity contains $(2, 1)$.

Denote by $\Sigma_n$ the set (not seen as a group) of permutations of $\{1, \dots, n\}$ and $\Sigma_\infty$ the disjoint union of the $\Sigma_n$'s, for $n \ge 1$. This defines a partial ordering $\le$ on $\Sigma_\infty$. Understanding this ordering is called *pattern recognition*, as one also says that $\sigma$ *is a pattern in* $\pi$ when $\sigma \le \pi$.

The website Database of Permutation Pattern Avoidance contains a huge number of examples.

A subset $\mathcal{C} \subset \Sigma_\infty$ is called a *permutation class* if $\pi \in \mathcal{C}$ and $\sigma \le \pi$ implies $\sigma \in \mathcal{C}$. For such a permutation class, we can consider its *basis* $\mathcal{B}$ consisting of those permutations $\pi$ which are *not* in $\mathcal{C}$ but such that any $\sigma \le \pi$, different from $\pi$, is in $\mathcal{C}$. So, a permutation $\pi$ is in $\mathcal{C}$ if and only if it does not contain an element of $\mathcal{B}$. I will write $\mathcal{C} = Av(\mathcal{B})$ and say that $\mathcal{C}$ consists of permutations *avoiding* $\mathcal{B}$.

Try and prove the Erdös-Szekeres theorem: every permutation $\pi \in \Sigma_n$ with $n > (p-1)(q-1)$ contains $(1, 2, 3, \dots, p)$ or $(q, q-1, \dots, 2, 1)$.

For instance the set $Inter \subset \Sigma_\infty$ of polynomial interchanges is obviously a permutation class. We have seen that its basis consists of two elements $(2, 4, 1, 3), (3, 1, 4, 2)$.

$Av((1, 2, 3))$ consists of those permutations that can be written as the union of two decreasing sequences.

The following are typical questions in the theory. Given a permutation class $\mathcal{C}$:

- Can we determine its basis? When is it finite?
- Can we count the number of elements in $\mathcal{C} \cap \Sigma_n$? Or at least, can we estimate this number?
- Can we decide algorithmically if a given permutation $\pi$ is in $\mathcal{C}$? What is the complexity of such an algorithm?

Find some example of a permutation class with an infinite basis.

We will answer all these questions in due time for the class of polynomial interchanges/separable permutations.

## *Stack-sortable permutations*

The theory of permutation patterns received a strong impetus from one exercise in volume 1 of *The art of computer programming*. Donald Knuth had the idea of attributing a degree of difficulty to the exercises in his book.

A "0" means that the reader should be able to solve it instantly.

A "10" requires one minute.

A "20" may require several hours, etc.

The scale is indeed logarithmic and even seems to have some pole in the neighborhood of 50...

The exercise that I want to discuss is labeled [M28]. The M means that it is aimed at mathematically inclined readers and the 28 is an indication of the time required to solve it (in the previously explained logarithmic scale). Today, this is not so hard but it turns out that this exercise had a lasting influence on combinatorics.

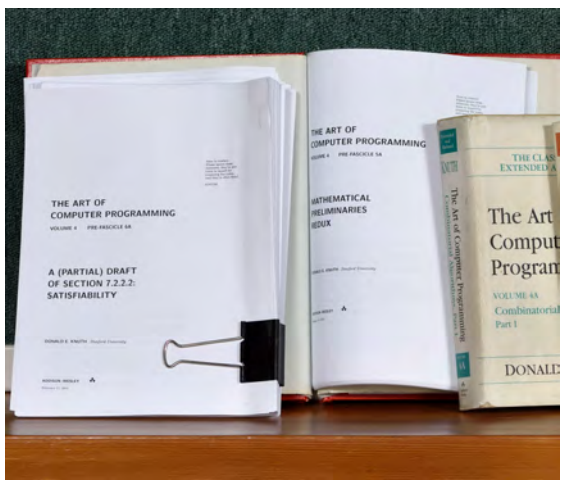

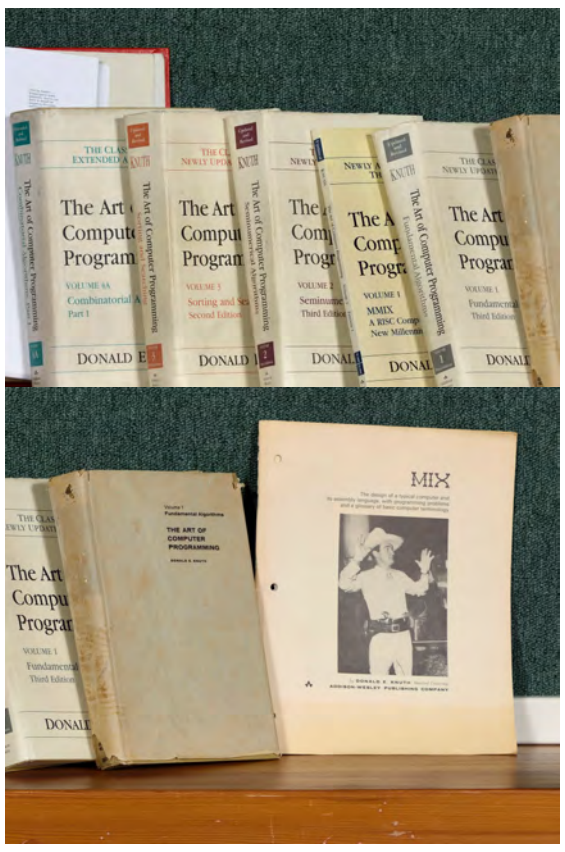

The art of computer programming.

I will describe a permutation class which is defined by a *stack*.

Imagine $n$ objects labeled $1, 2, \ldots, n$ lined on some horizontal line, in this order from left to right: $1 \ll 2 \ll \ldots \ll n$. On the right of $n$, there is a *stack*. This is some kind of well in which the objects can be piled on top of each other.

Initially, the stack is empty. Then select the object $n$ and push it in the stack. Then, there are two options. Either we *push* the last element on the line onto the top of the stack. Or we *pop* the top element of the stack to the right.

Look at the figure in the margin, and the evolution of the objects under a sequence Push, Push, Pop, Push, Pop, Pop, Push, Push, Push, Pop, Pop, Pop.

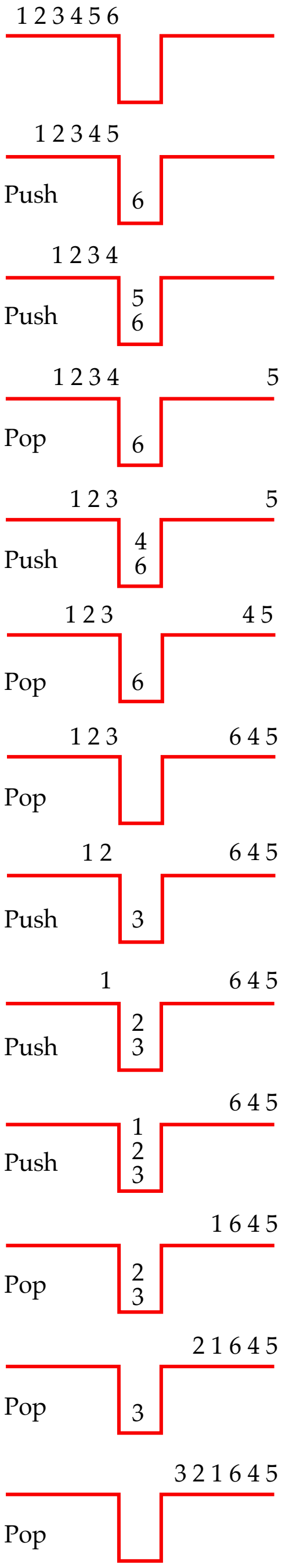

At the end of the operation the sequence $(1, 2, 3, 4, 5, 6)$ has been transformed in $(3, 2, 1, 6, 4, 5)$. This could be seen either as a permutation or as two orderings: $1 \ll 2 \ll 3 \ll 4 \ll 5 \ll 6$ (on the left) and $3 \lll 2 \lll 1 \lll 6 \lll 4 \lll 5$ (on the right).

**Definition.** A permutation $\pi$ is *stack-sortable* if it is the result of a sequence of Push and Pop's applied to $\{1, 2, \ldots, n\}$.

Exercise 5 from Knuth chapter 2, evaluated as [M28], is the following:

**Theorem.** *A permutation is stack-sortable if and only if it does not contain* $(2, 3, 1)$.

Let us solve this exercise.

Start with a permutation, for example $(3, 2, 1, 6, 4, 5)$. The last element is 5. If we want to sort this permutation with a stack, there is no choice: we have to push all elements in the list $(1, 2, 3, 4, 5, 6)$ until 5 is available on top of the stack, so that we can pop it and put it in its proper place, at the end of the output list. Then we look for the second to last, that is to say 4, and we continue pushing until 4 is on top of the stack etc. So, if we want to sort a permutation, there is only one way to do it.

We only have to understand when the sorting might go wrong. This will happen precisely when it would be time to pop some object $a$ which is unfortunately already in the stack but not on top, below some object $b \ll a$. If $b$ has been already pushed in the stack, this is because we had previously to pop some other object $c \ll b \ll a$.

We have $a \lll c$ since $c$ has already been popped and we are trying to pop $a$. Similarly, we have $b \lll a$ since we don't want to pop $b$ but $a$. The subset $\{c, b, a\}$ in $\{1, 2, \ldots, n\}$ therefore gives rise to the containment $(2, 3, 1) \le \pi$ as we had to show. ⊡

I strongly encourage the reader to do *all* exercises in Knuth's book.

The train track in the margin produces the stack-sortable permutations. A train consisting of cars $(1, 2, \ldots, n)$ arrives from the left. Cars can then be disassembled and each one has to follow the tracks in the direction given by the arrows. The train is assembled again on the exit side, on the right of the picture.

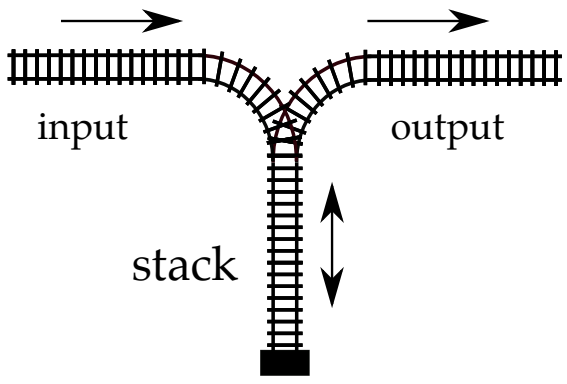

This simple train track produces stack-sortable permutations.

Knuth also defines *deques* (a combination of deck and queue). They are produced by the more complicated train track pictured below.

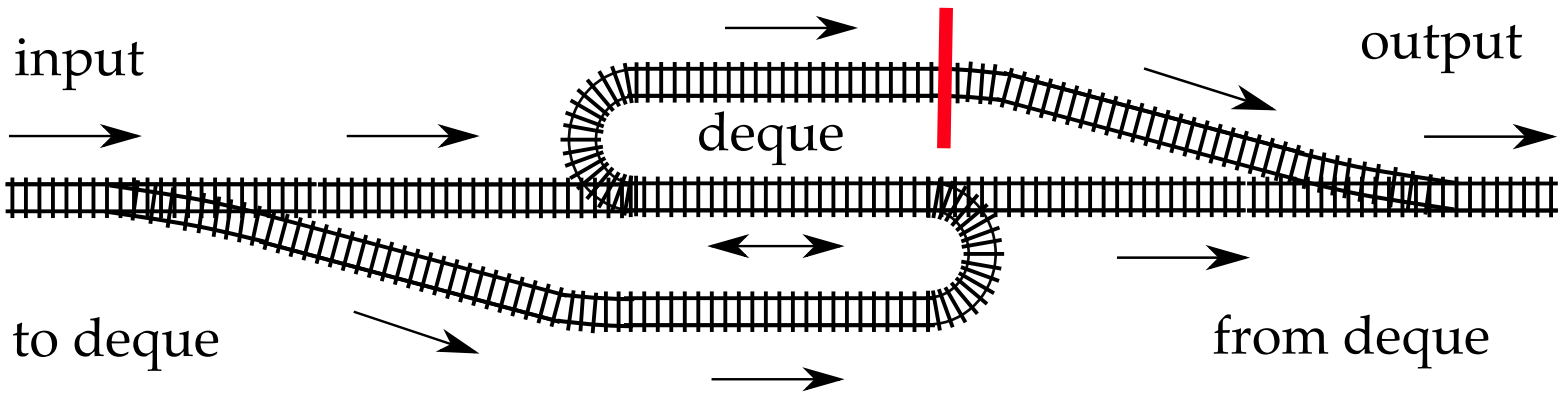

This is a deck-queue = a deque.

Beware ! the difficulty of this question might be around 60!

What are the deque-sortable permutations?

If the the red door is closed, we get an *output-restricted deque*. The associated permutations are exactly those which do not contain $(4,2,3,1), (4,1,3,2)$. We are getting closer to the characterization of polynomial interchanges, which avoid $(2,4,1,3)$ and $(3,1,4,2)$.

For much more about this fascinating field, I recommend the above mentioned book by Kitaev.

## *Ubiquitous Catalan*

Exercise 4 in the same chapter of Knuth's book is rated [M34]. However, it is easier to solve after having solved exercise 5.

The problem is to count the number of stack-sortable permutations of length $n$.

This is the famous $n$-th *Catalan number* $C_n$, that appears almost everywhere in mathematics.

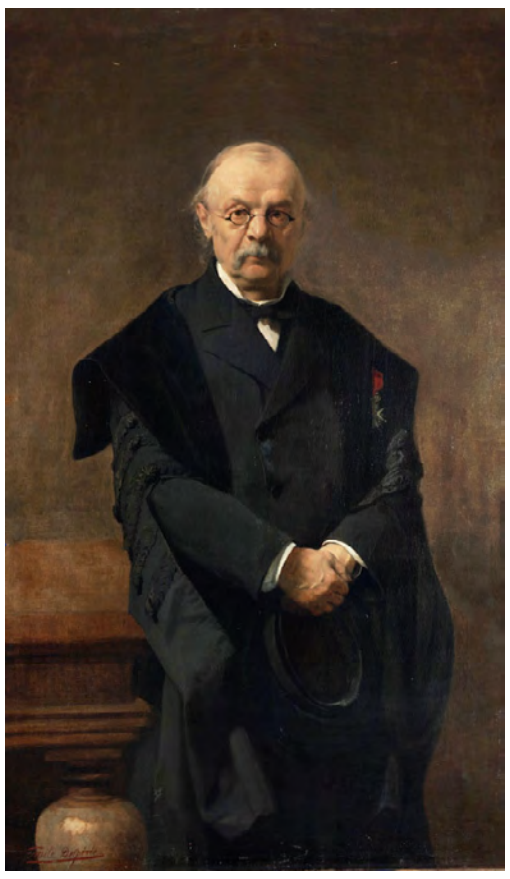

Eugène Catalan was born in 1814, in Bruges, then belonging to the Napoleonic French empire. ©

The first values are:

1,2, 5, 14, 42, 132, 429, 1430, 4862, 16796, 58786, 208012, 742900, 2674440, 9694845, 35357670, etc. (sequence A000108 in OEIS).

Stack-sortable permutations are described *uniquely* by a sequence of $2n$ Push's and Pop's. Conversely, a sequence of Push's and Pop's defines a permutation under the only condition that we are not forced to pop an empty stack. Said differently, every initial segment of the sequence should contain at least as many Push's as Pop's. This could also be described by looking at the evolution of the number of elements in the stack. The stack is empty at time $0$ and $2n$, changes by $+1$ or $-1$ steps for each Push and Pop. This is called a *Dyck word* of length $2n$. Here is an example of a Dyck word of length 24.

Walther von Dyck (1856–1934) was "the first to define a mathematical group, in the modern sense". This is at least what can be found in *Wikipedia*. The question is much more subtle and many more names should be mentioned. However, he was indeed one of the first to manipulate presentations of groups, with generators and relations.

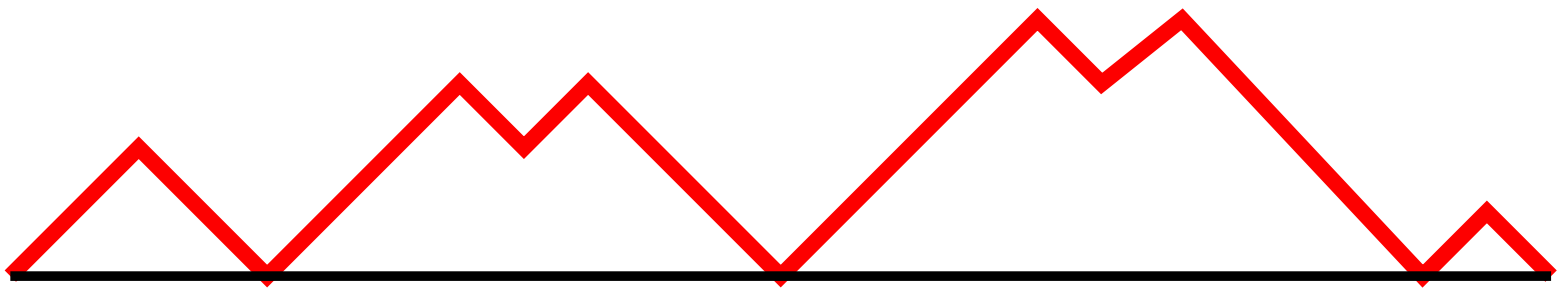

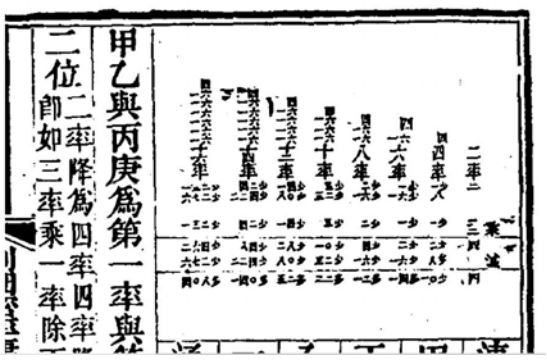

The first appearance of Catalan's numbers in Mingantu's book "The Quick Method for Obtaining the Precise Ratio of Division of a Circle" around 1730. ©

The number of these Dyck words is one of the many definitions of the $n$-th Catalan number. Alternatively, we could think of Push as an open parenthesis "(" and Pop as a closing one ")". The condition that we never pop an empty stack is now equivalent to the fact that the sequence of parentheses is correctly balanced. This means recursively that every open "(" is coupled with a closed ")" which enclose together a correctly balanced sequence of parenthesis. For instance, for $n = 3$, there are $C_3 = 5$ sequences

$$(())() \quad ; \quad ()(()) \quad ; \quad ((())) \quad ; \quad ()()() \quad ; \quad (()()).$$

There is also an interpretation in terms of *rooted planar trees*. A picture is worth a thousand words. For some strange reason, mathematicians and computer scientists tend to draw trees upside down: the root is on top and the leaves are on the bottom.

The picture in the margin is an example of such a tree. It has one root, 3 internal nodes and 4 leaves. The tree is *planar* because the children of its nodes are ordered from left to right. Equivalently, we will say that a tree is planar if its leaves have been totally ordered in such a way that the descendants of any node define some interval.

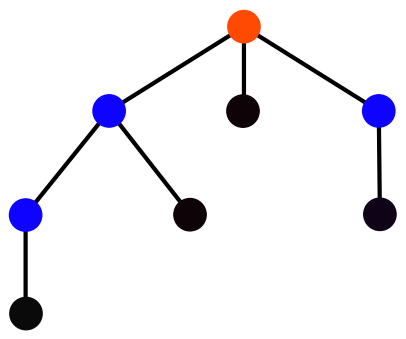

Such a tree defines a Dyck word. Just start from the root and follow the tree externally, going counter-clockwise. At each step you get further or closer to the root: this gives the sequence of +1 and −1, or "Push" and "Pop".

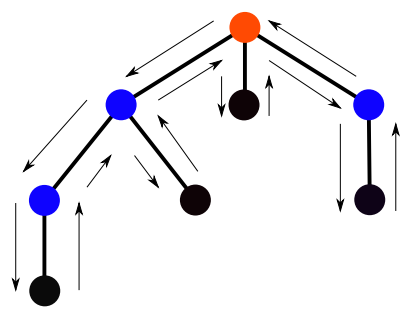

In this example, we get the sequence $+ + + - - + - - + - + + - -$. Conversely, a Dyck word can be transformed into a rooted planar tree.

All in all, there are bijections between:

- Stack-sortable permutations of $n$ objects.
- Dyck words of length $2n$.
- Balanced bracketings with $2n$ parentheses ($n$ open and $n$ closed).
- Rooted planar trees with $n$ edges.

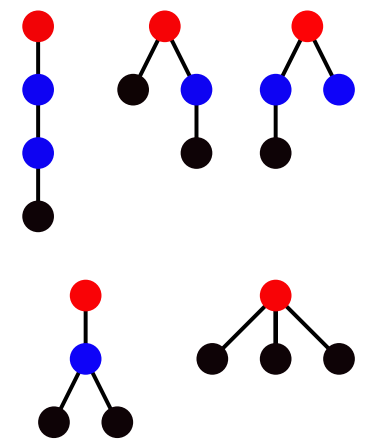

5 rooted planar trees with 3 edges.

The cardinality of any of these sets is the $n$-th Catalan number $C_n$.

Given a stack-sortable permutation, look at the last element $n$ of the list of objects (on the left, in black in the picture) and at its position $k$ after the sorting process. The permutation maps the (red) interval $\{k, \ldots, n-1\}$ to $\{k+1, \ldots, n\}$ and the yellow one $\{1, 2, \ldots, k-1\}$ to itself. Therefore, it induces a stack-sortable permutation on these two intervals. Hence we get the recurrence relation:

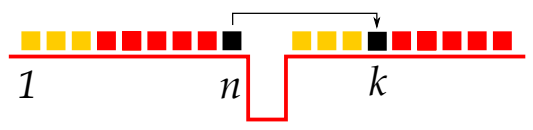

$$C_n = \sum_{k=1}^{n} C_{k-1} C_{n-k}.$$

Look at this recurrence relation in terms of Dyck words, bracketings and rooted planar trees.

This is Catalan's characteristic signature: it is found in many different contexts.

For instance, consider rooted planar *binary trees*. Their definition depends on the authors but let me define them as planar rooted trees such that every node has no children or two children, one being "to the left" and the other being "to the right". If such a tree has $n$ internal nodes, it has $n+2$ leaves, and $2n+2$ edges. If you remove its root, you get two rooted planar binary trees. Conversely you can add a common root to two rooted planar binary trees to produce a bigger rooted planar binary tree. This shows, after a moment of reflection, that the number

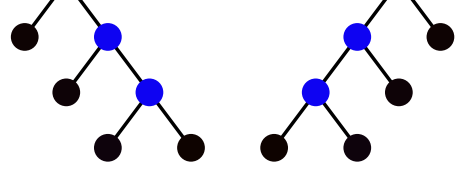

5 rooted planar binary trees with 4 leaves.

of rooted planar binary trees with $n+1$ leaves satisfies Catalan's recurrence relation. One can check that there are $1,1,2,5$ rooted planar binary trees with $1,2,3,4$ leaves, and we therefore get by induction that *$C_n$ is also the number of rooted planar binary trees with $n+1$ leaves.*

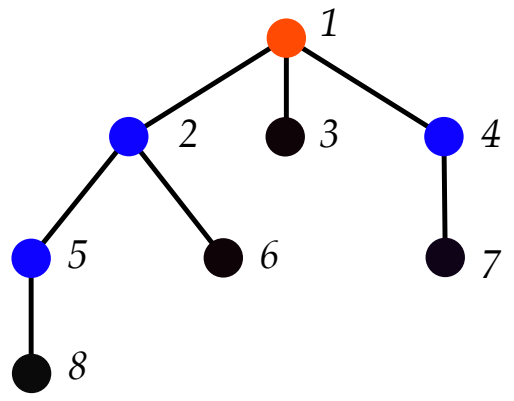

This suggests that there should be some correspondence between rooted planar trees and rooted planar binary trees. This is indeed the case. Let me present a slight variation on the so-called *Knuth transform* or *first child-next sibling representation*. Starting with some rooted planar tree $T$ with $n$ edges (first picture), we construct a rooted planar binary tree $T_{bin}$ with $n+1$ leaves (last picture). I first construct an auxiliary tree $T'$ (second and third pictures). The set of nodes of $T'$ is the same as the set of nodes of $T$. The root is the same. Every node $v$ of $T'$ has at most two children. The first is the eldest child of $v$ in $T$, if it exists. The second is the next sibling of $v$ in $T$, that is to say the eldest among siblings younger than $v$, if it exists. Then I transform $T'$ in a rooted planar binary tree $T_{bin}$ in the following way. First I delete the root and the edge going out of it. The new root of $T_{bin}$ is the eldest child of the root of $T$. For every node of $T'$, I add a left child if this node has no children in $T$ and a right child if this node has no younger sibling. Thus, if the node is a leaf of $T'$ (i.e. it has no children and no younger sibling in $T$), I add two children in $T_{bin}$ (see the green dots in the fourth picture). Check that this gives a bijection between rooted planar trees with $n$ edges and rooted planar binary trees with $n+1$ leaves ([M15]).

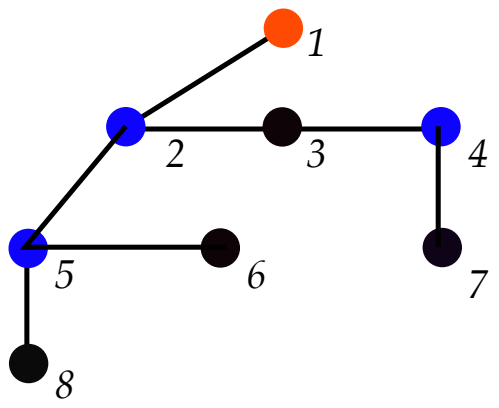

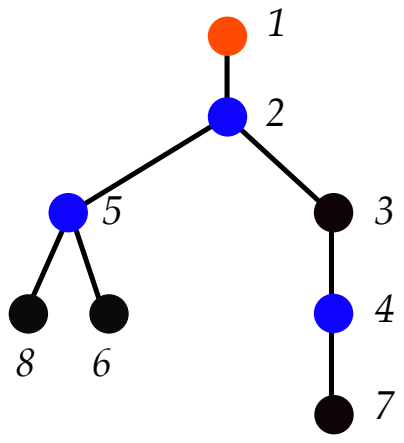

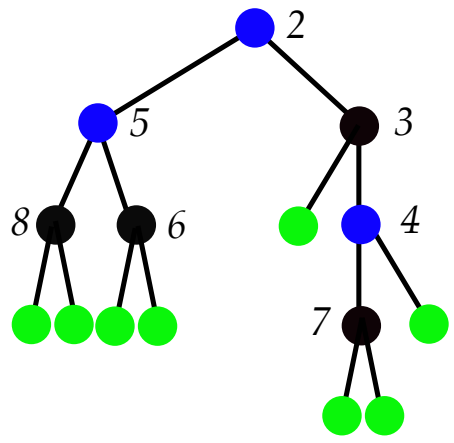

As usual in combinatorics, this sequence $C_n$ is encoded by its formal *generating series*

$$C(t) = \sum_{n\geq 0} C_n t^n$$

and the recurrence becomes:

$$C(t) = tC^2(t) + 1.$$

Remembering secondary school quadratic equations, we get

Why did I choose the − sign in front of the square root?

$$C(t) = \frac{1-\sqrt{1-4t}}{2t}.$$

It follows from this formula that the radius of convergence of $C(t)$ is 1/4. So, by the Cauchy-Hadamard theorem, the growth of $C_n$ can be estimated:

$$\limsup_{n\to\infty} \frac{1}{n} \log C_n = \log 4.$$

This formula can also be used to get a neat expression for $C_n$: just use Newton's binomial series.

$$\sqrt{1-4t} = \sum_{n=0}^{\infty} \binom{1/2}{n} (-4t)^n.$$

Comparing the coefficients of $t^n$:

$$C_n = -\frac{1}{2}\binom{1/2}{n+1}(-4)^{n+1}.$$

We now clean a bit:

$$\begin{aligned} C_n &= -\tfrac{1}{2}\tfrac{1}{(n+1)!}\left(\tfrac{1}{2}\right)\left(\tfrac{1}{2}-1\right)\cdots\left(\tfrac{1}{2}-n\right)(-1)^{n+1}2^{2(n+1)} = \ldots \\ &= \tfrac{1}{n+1}\tbinom{2n}{n}. \end{aligned}$$

The second "…" mean that you are encouraged to do the computation yourself!

The bibles of the subject of Catalan numbers are the books by Stanley[4], [5]. The Catalan website "Catalan Numbers", maintained by Igor Pak, is a remarkable source of information. The book by Flajolet and Sedgewick[6] provides a wider perspective (see in particular chapter 6 on trees).

[4] R. P. Stanley. *Enumerative combinatorics. Vol. I*. The Wadsworth & Brooks/Cole Mathematics Series. Wadsworth & Brooks/Cole Advanced Books & Software, Monterey, CA, 1986. With a foreword by Gian-Carlo Rota.

[5] R. P. Stanley. *Catalan numbers*. Cambridge University Press, New York, 2015.

[6] R. S. P. Flajolet. *An introduction to the analysis of algorithms*. Addison-Wesley, 2nd ed edition, 2013.

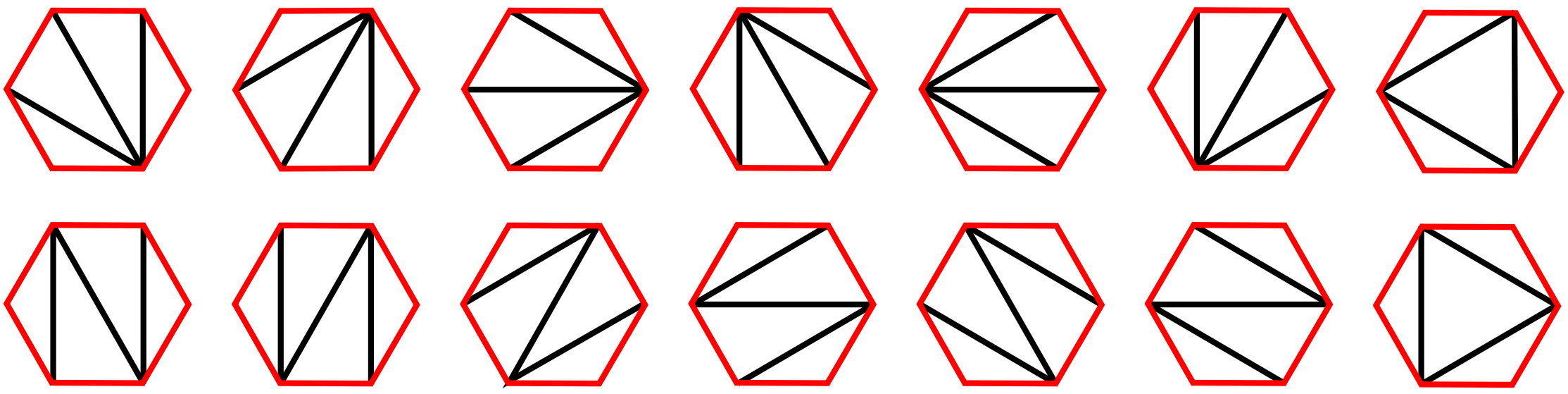

There are $C_4 = 14$ ways to subdivide a hexagon in triangles.

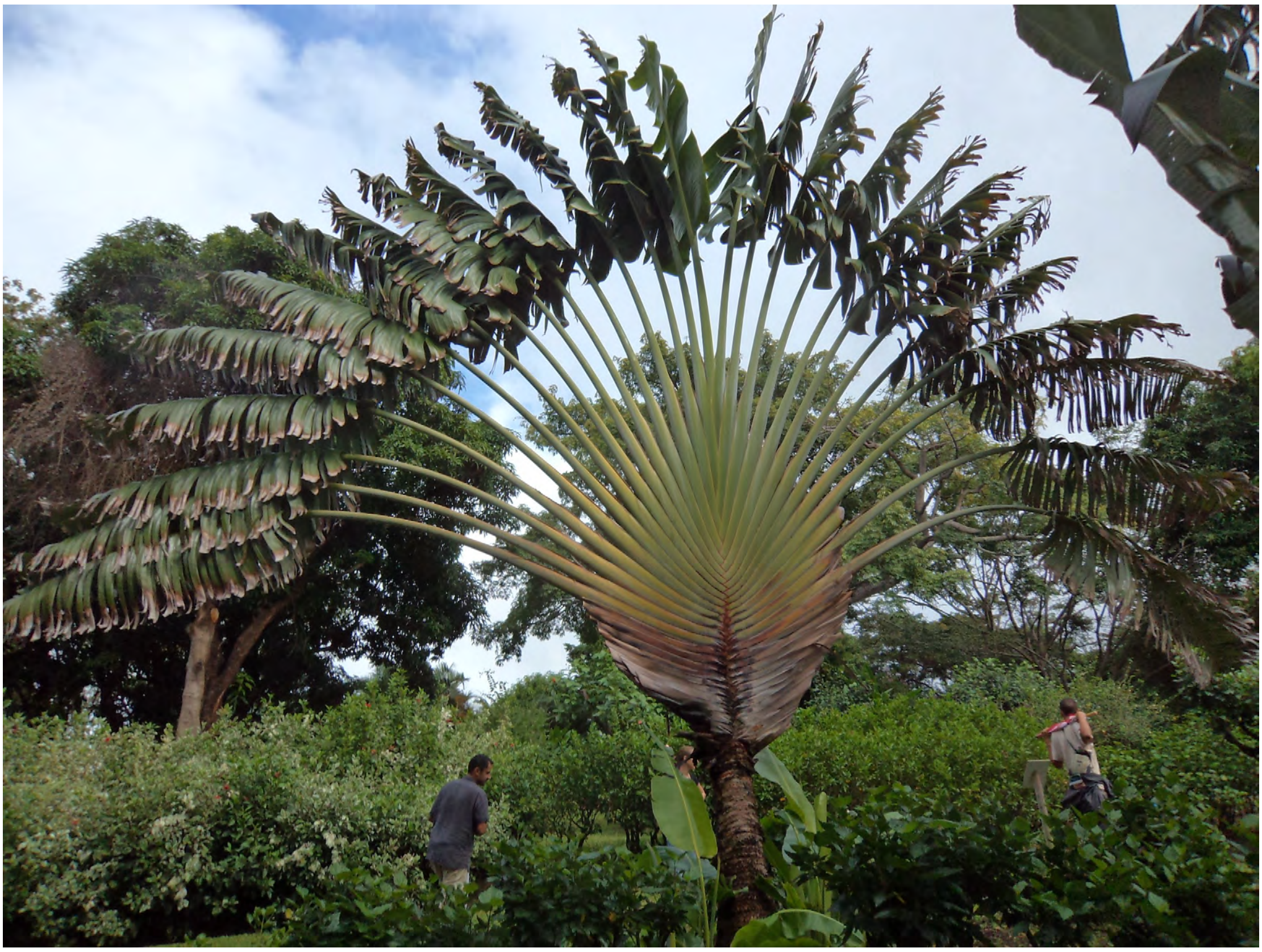

A planar tree in Guadeloupe, commonly known as traveller's palm.

# *Separable permutations*

## *From polynomials to trees*

The ring of polynomials $K[x]$ with coefficients in some field $K$ of characteristic 0 is equipped with a valuation $v$, given by the degree of the first nonzero coefficient, and a natural *ultrametric distance*, defined in terms of $v$:

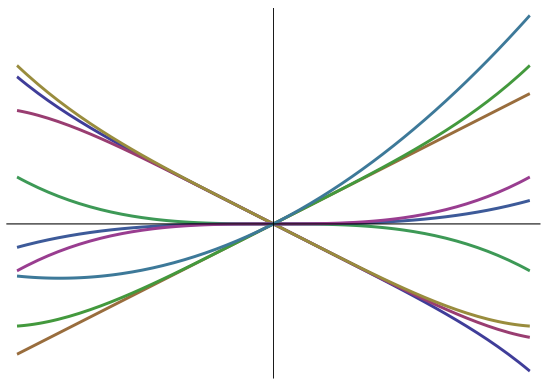

9 polynomials whose associated tree is below. Can you propose 9 equations?

$$dist(P, Q) = \exp(-v(P - Q)).$$

In plain words, two polynomials are close if their first $k$ derivatives at 0 coincide, for a large value of $k$. The ultrametric property for a distance means precisely that for every $\epsilon > 0$ the relation

$$dist(P, Q) < \epsilon$$

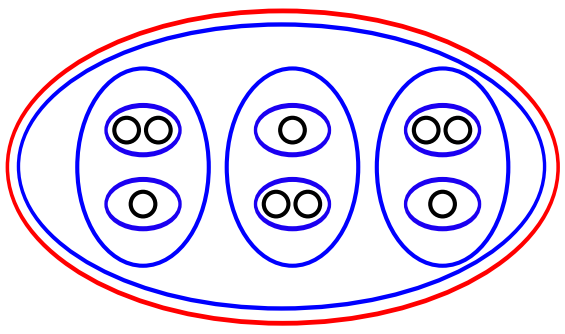

Nested equivalence relations.

is an equivalence relation. As $\epsilon$ decreases, these equivalence relations get finer and their intersection is trivial.

*Consider a finite set E of polynomials.* Define a tree in the following way. The root is labeled by the set $E$. The children of the root are labeled by the equivalence classes of the relation $v(P - Q) \geq 1$. The grandchildren of the root are labeled by the equivalence classes of the relation $v(P - Q) \geq 2$. And in general, the $k$-th generation corresponds to equivalence classes of the relation $v(P - Q) \geq k$. This tree is infinite but the equivalence classes stabilize to singletons when $k$ is large. We can therefore do some pruning in order to get a finite tree whose leaves are labeled by the elements of $E$. Conversely, the valuation structure can be recovered from the tree. Given two elements $P, Q$ of $E$, seen as

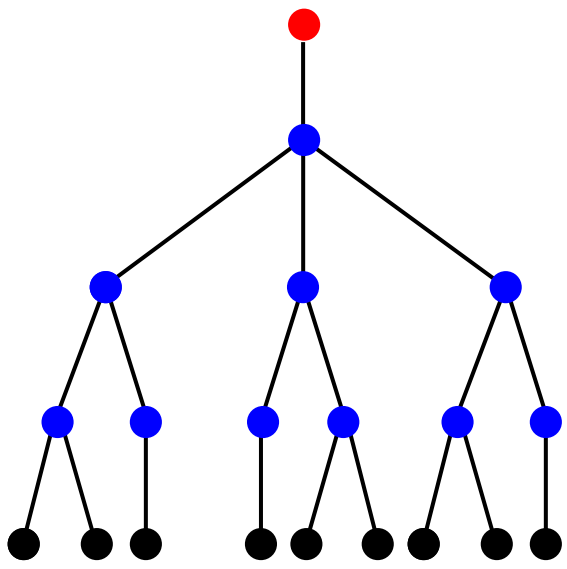

Associated tree.

leaves, we look in the tree for their closest common ancestor. The valuation of $P - Q$ is the *level* of this ancestor, defined as the length of the path connecting it to the root.

*Now suppose that $K$ is the field of real numbers.* As already observed, any finite set of real polynomials is equipped with two total orderings, comparing values for small negative and small positive $x$. Both provide an ordering on the set of leaves of the tree given by the valuation. Note that the descendants of a node, i.e. an equivalence class at some level, define an interval in each of these orderings. *Our tree is therefore a* planar *tree in two ways*.

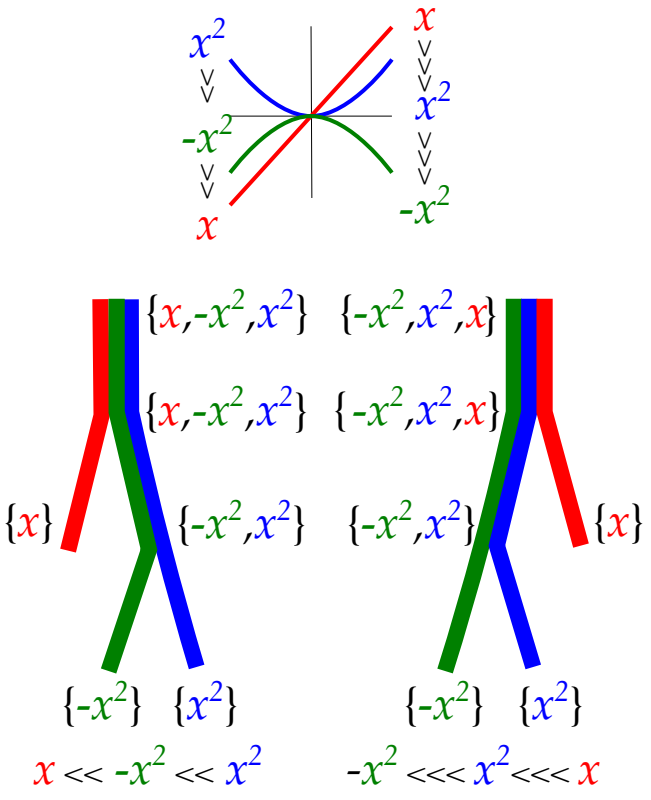

The two planar trees associated to the polynomials $\{x, -x^2, x^2\}$, looking at small negative (left $\ll$) and small positive (right $\lll$) values of $x$.

By convention, *let us choose the first order (i.e. for small $x < 0$) and let us associate to our set $E$ of polynomials the corresponding planar tree*. The comparison between the two orderings defines a permutation $\pi$ that we called a *polynomial interchange*.

If $P(x) < R(x) < Q(x)$ for small positive $x$ (or for small negative $x$), then $val(P - Q) \geq val(R - Q)$. It follows that all equivalence classes $val(P - Q) \geq k$ are intervals in both orderings: our trees are indeed planar.

The set of leaves of *any* rooted planar tree is equipped with two canonical orderings. The first, denoted $\ll$, is simply the order given by the definition of planarity. The second, denoted $\lll$, is defined in the following way. Given two leaves $a$ and $b$, denote by $a \vee b$ their closest common ancestor. Then $a \ll b$ and $a \lll b$ hold true simultaneously if and only if the level of $a \vee b$ is even. One has to check that this defines indeed an ordering $\lll$, in other words that $a \lll b \lll c \lll a$ is not possible. We can assume that $a \ll b \ll c$ or $c \ll b \ll a$ and the second case reduces to the first by symmetry. If we had $a \lll b \lll c \lll a$, the levels of $a \vee b$ and $b \vee c$ should be even, and the level of $c \vee a$ should be odd. The pictures in the margin show that this is not possible.

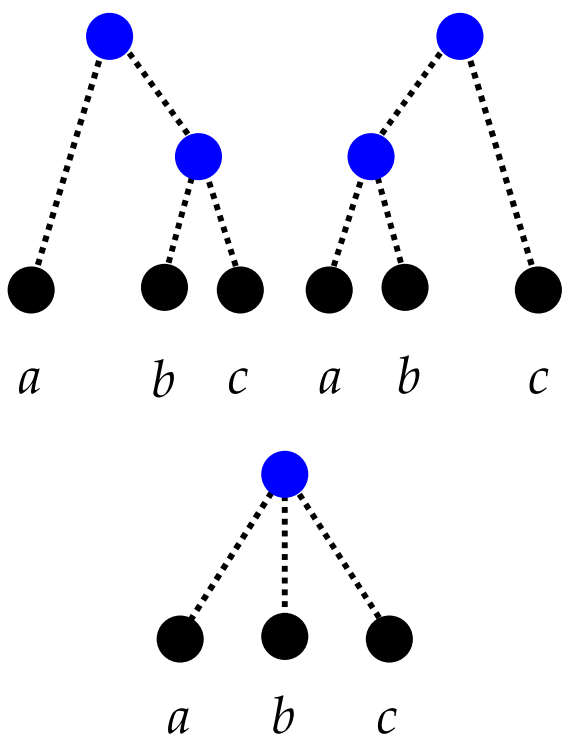

Let us sum up.

• A finite set of real polynomials produces a rooted planar tree.

• A rooted planar tree defines two orderings in its set of leaves, and therefore a permutation of the leaves.

• The permutation associated to the planar tree which is associated to a finite set of polynomials is simply the corresponding polynomial interchange.

Our trees contain too much information and we will prune their edges.

First, if an internal node has only one leaf among its descendants, we can delete all its descendants without changing the permutation (green path in the margin).

Suppose now that two internal nodes $P, Q$ are connected in the tree by some path such that all vertices between $P$ and $Q$ are non-ramified, i.e. have only one child (red and blue paths). If the number of edges in this path is even, just delete it and identify the two endpoints as a single node. If the number of edges in this path is odd, just delete it and connect the two endpoints by a single edge. This produces a new tree. In this process, the levels of some nodes have changed, but only by an even number. Therefore, if we compute the valuation in the new tree, the parity did not change and this parity is the only information that matters in order to construct the polynomial interchange. Note that the pruned tree has the property that all its internal nodes have at least two children.

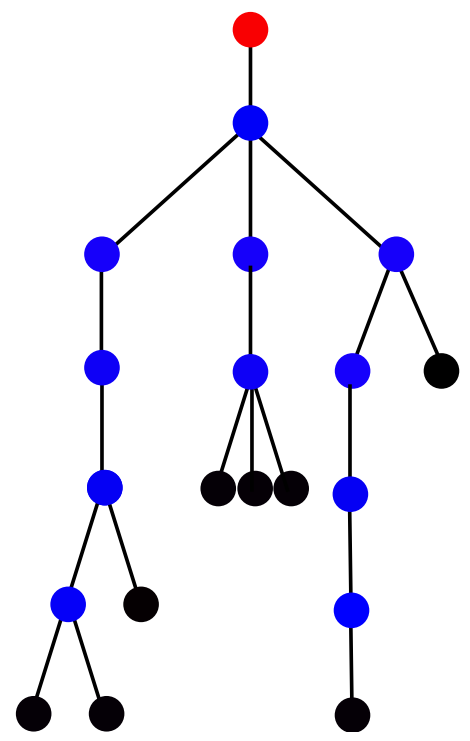

Tree before pruning.

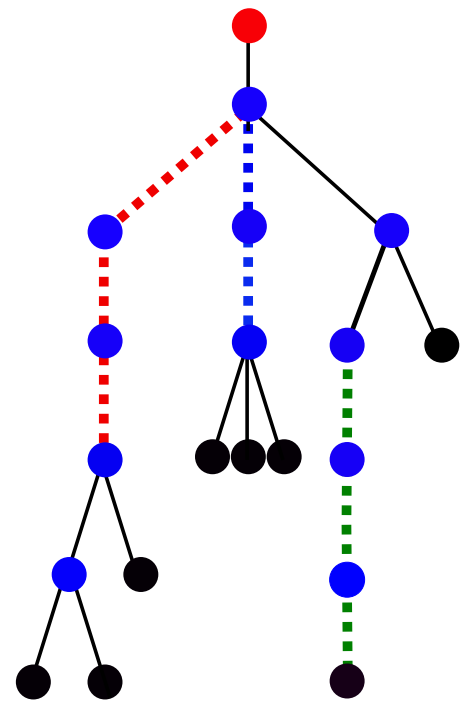

Deleting edges.

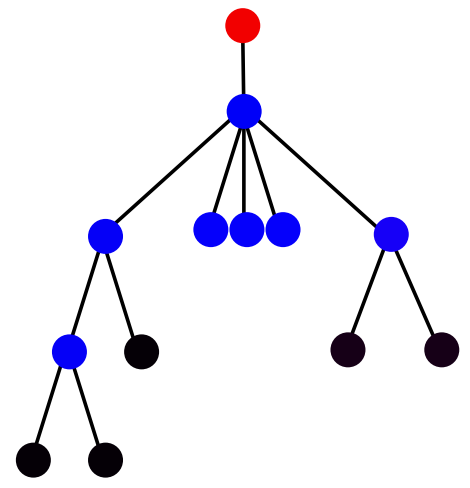

Reconnecting.

In summary, given $n$ polynomials, we constructed a *rooted planar tree* such that:

- The root can have any number of children.
- Every internal node has at least two children.
- There are exactly $n$ leaves, labeled by the $n$ polynomials.

Say that a planar tree is *pruned* if it satisfies these properties. It should be clear that for any pruned tree, one can find $n$ polynomials such that the associated pruned tree is the given one. We have seen that for any finite set of polynomials, the associated polynomial interchange can be read from the tree. In particular, the number of polynomial interchanges is less than or equal to the number of pruned trees.

## *From a permutation to a tree*

We will now show that the number of polynomial interchanges is *equal* to the number of pruned trees. The issue is to prove that *two different pruned trees produce different permutations.*

Let $T$ be a pruned tree with leaves $1 \ll 2 \ll \ldots \ll n$, from left to right, and let $\pi$ be the associated permutation ($n \geq 2$).

**Lemma.** *The images by $\pi$ of two consecutive integers $i, i+1$ are consecutive if and only if the corresponding leaves $i, i+1$ are siblings, i.e. have a common parent.*

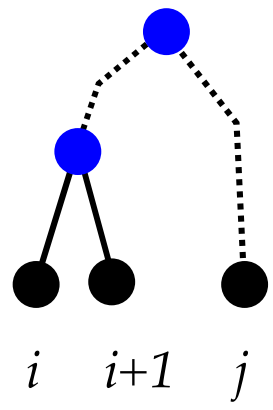

Indeed, if $i$ and $i+1$ have a common parent, one has $i \vee j = (i+1) \vee j$ for every $j \neq i, i+1$,. It follows that $i \lll j \lll i+1$ and $i+1 \lll j \lll i$ are impossible. Said differently, $i, i+1$ are also consecutive for $\lll$.

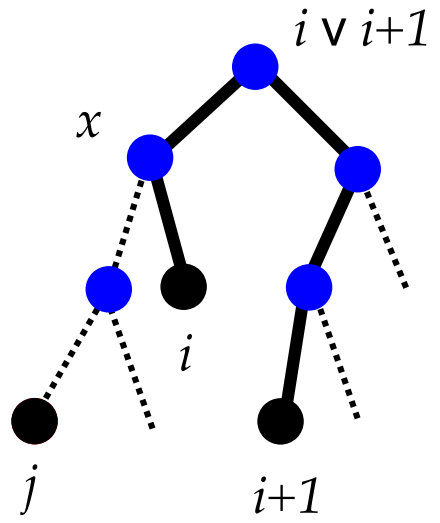

Conversely, suppose that $i \vee (i+1)$ is not a parent. Then the path connecting $i$ and $i+1$ in the tree has length at least 3 and contains therefore some internal node $x$ such that the levels of $x$ and $i \vee (i+1)$ are different. Choose a leaf $j$ which is a descendant of $x$ different from $i, i+1$ (which exists since $T$ is pruned). It follows that $j$ is between $i$ and $i+1$ for the order $\lll$. ⊡

This gives another proof that any polynomial interchange contains at least two consecutive integers with consecutive images. It suffices to consider an internal node with the highest level in $T$: it has at least two children which must be leaves, and therefore siblings.

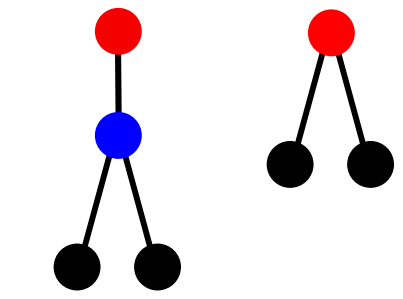

The 2 pruned trees with 2 leaves, defining the transposition and the identity.

*We can now prove, by induction on $n \geq 2$, that there is at most one pruned tree producing a given permutation $\pi$.* This is of course trivial for $n = 2$. If $T_1$ and $T_2$ have $n$ leaves and define the same permutation $\pi$, the previous lemma shows that there is some pair or consecutive leaves $i, i+1$ which are siblings, for both $T_1$ and $T_2$. Delete the leaf $i+1$ from $T_1$ and $T_2$, producing trees $T_1'$ and $T_2'$ with $n-1$ leaves. Clearly $T_1'$ and $T_2'$ define the same permutation $\pi'$ on $n-1$ leaves. We should be careful however that $T_1'$ or $T_2'$ might not be pruned. This happens precisely when $i$ and $i+1$ are the only children of some internal node in $T_1$ or $T_2$. The induction hypothesis shows that the pruned trees $T_1''$ and $T_2''$ are equal to some $T''$. Our trees $T_1$ and $T_2$ are obtained from $T''$ under one of the following two operations: adding two children to the $i$-th leaf of $T''$, or adding a sibling to the $i$-th leaf of $T''$. By assumption $T_1$ and $T_2$ produce the same permutation so that the parents of $i$ and $i+1$ in $T_1$ and $T_2$ have levels of the same parity. It follows that $T_1$ and $T_2$ are obtained from $T''$ by performing the same operation. Hence $T_1 = T_2$, as desired. ⊡

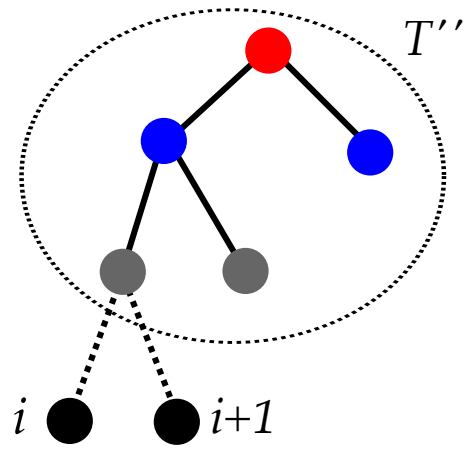

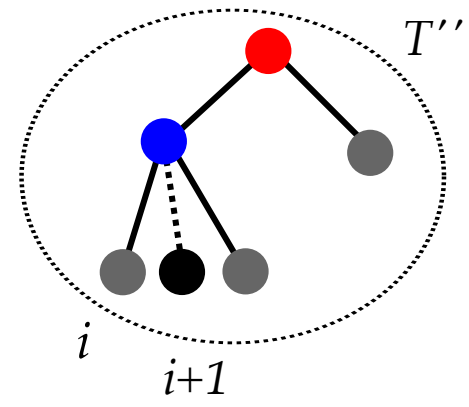

*Therefore the number of polynomial interchanges of size $n$ is equal to the number of pruned trees with $n$ leaves.*

## *From a pruned tree to a polynomial interchange and a separable permutation*

This leads us to the original definition of *separable permutations*. Given two permutations $\pi_1$ and $\pi_2$ of $n_1$ and $n_2$ ordered objects, we can think of two ways to produce a permutation of $n_1 + n_2$ objects. Number the first $n_1$ objects $\{1, 2, \ldots, n_1\}$ and the next $n_2$ as $\{n_1 + 1, n_1 + 2, \ldots, n_1 + n_2\}$. Denote by $\pi_1 \oplus \pi_2$ the permutation defined by

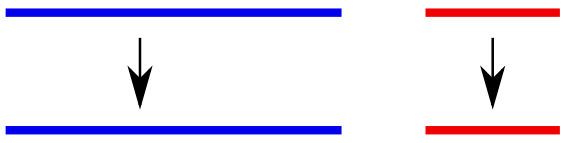

$\pi_1 \oplus \pi_2$.

$$\begin{aligned} \pi_1 \oplus \pi_2(k) &= \pi_1(k) && \text{if } 1 \le k \le n_1 \\ &= n_1 + \pi_2(k - n_1) && \text{if } n_1 + 1 \le k \le n_1 + n_2. \end{aligned}$$

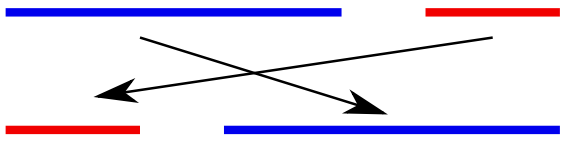

$\pi_1 \ominus \pi_2$.

Then define

$$\begin{aligned} \pi_1 \ominus \pi_2(k) &= \pi_2(k) + n_1 && \text{if } 1 \le k \le n_2 \\ &= \pi_1(k - n_2) && \text{if } n_2 + 1 \le k \le n_1 + n_2. \end{aligned}$$

In the definition of 1998[7], a permutation is separable if it is obtained from several copies of the trivial permutation on one object by successive $\oplus$ and $\ominus$ operations. We indeed have a fairly good understanding of these permutations[8].

[7] P. Bose, J. F. Buss, and A. Lubiw. Pattern matching for permutations. *Inform. Process. Lett.*, 65(5):277–283, 1998.

Note that if a permutation $\pi$ is separable, so is its reverse $\overline{\pi}(k) = n + 1 - \pi(k)$.

[8] É. Ghys. Intersecting curves (variation on an observation of Maxim Kontsevich). *Amer. Math. Monthly*, 120(3):232–242, 2013.

**Theorem.** *Let $\pi$ be a permutation of $\{1, \ldots, n\}$. The following conditions are equivalent.*

1. *$\pi$ is the polynomial interchange associated to n distinct polynomials $P_1, \ldots, P_n$.*
2. *$\pi$ does not contain $(2,4,1,3)$ or $(3,1,4,2)$.*
3. *$\pi$ is the permutation defined by some pruned tree.*
4. *$\pi$ is obtained from several copies of the trivial permutation on one object by successive $\oplus$ and $\ominus$ operations.*

*These permutations have already been defined as* separable.

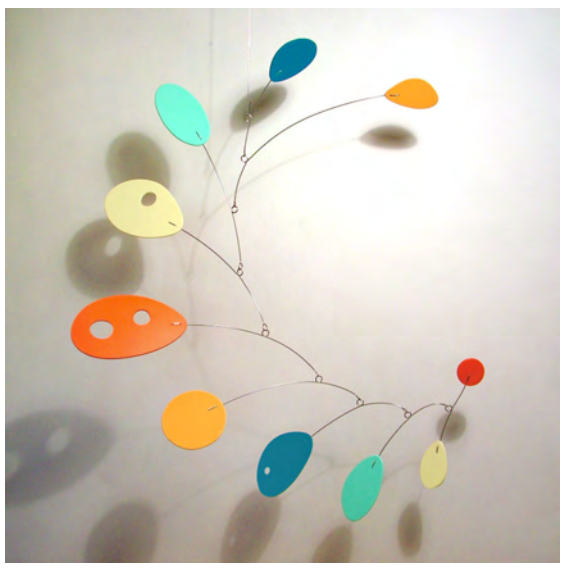

A mobile, à la Alexander Calder. Of course, these mobiles are not meant to be planar but any two of their planar realizations differ by a separable permutation of their leaves. ©

We already proved everything except the equivalence between 3 and 4. Let us prove this equivalence by induction.

Let $\pi$ be described by a pruned tree $T$. The root of this tree might have a unique child. If this is the case, the descendants of this unique child define another pruned tree $\overline{T}$ in which the new root has several children. The permutation associated to $\overline{T}$ is the reverse $\overline{\pi}$. Observe that $\pi_1 \ominus \pi_2 = \overline{\overline{\pi_1} \oplus \overline{\pi_2}}$. We can therefore assume that the root has several children. These children define pruned trees and polynomial interchanges $\pi_1, \pi_2, \ldots, \pi_k$ whose $\oplus$ sum is $\pi$. By induction, all permutations $\pi_1, \pi_2, \ldots, \pi_k$ are obtained from the trivial permutation by $\oplus$ and $\ominus$ operations so that the same is true for $\pi$.

The converse is just as easy. We have to show that the $\oplus$ and $\ominus$ sums of two permutations associated to pruned trees are also associated to a pruned tree. Just join the two roots of two pruned trees by a common parent, or add a grandparent, depending on the sign. ⊡

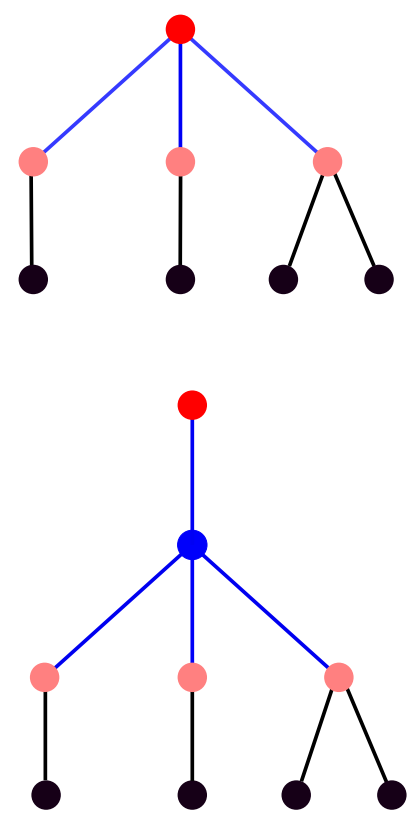

Adding an extra parent or a grandparent, to join several trees.

## *Train tracks, stacks, floorplans and permutons*

Imagine a (mathematical) train consisting of $n$ cars. Insert a mathematical *turntable*, one of these devices that one sees sometimes in railways, enabling 180 degrees rotations. The turntable is mathematical since it can fit any number of consecutive cars. We also assume that once a segment of cars has been reversed, these cars are hooked in a permanent way: they cannot be disconnected in future. Allow the train to move several times in the turntable. The final permutation of the cars is separable, almost by definition.

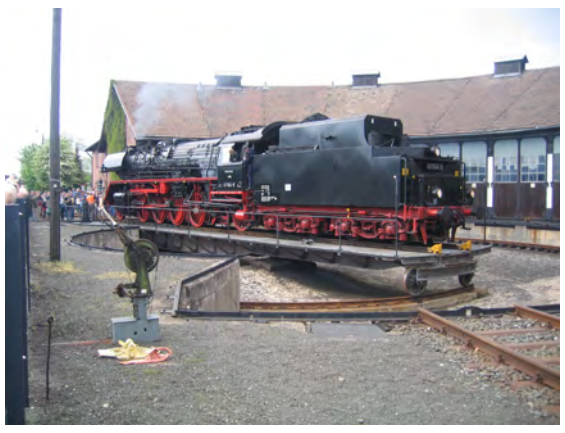

A railway turntable. ©

The same idea can also been expressed using *pop stacks in series*, as in the previous chapter. Imagine an indefinite sequence of stacks aligned on the right of the sequence $1, 2, \ldots, n$. The rules of the sorting game are different from the single stack case. At each step, we are allowed to push an element of the list on top of the first stack, or to pop the *full content* of some stack to the next one.

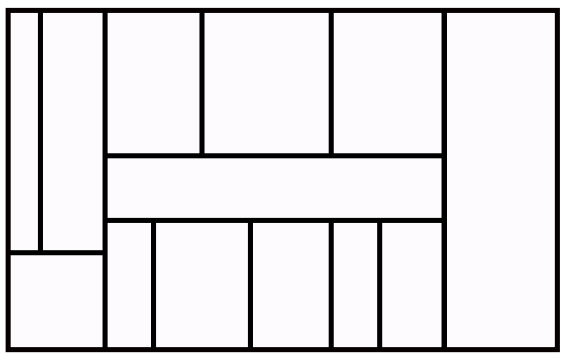

Here is another occurrence of separable permutations. Start with a rectangle and decompose it in several *rectangular rooms* by successive vertical or horizontal slicing. This is a *slicing floorplan*, as in the margin. Find a good definition for equivalent slicing floorplans and find a bijection with separable permutations[9].

[9] E. Ackerman, G. Barequet, and R. Y. Pinter. A bijection between permutations and floorplans, and its applications. *Discrete Appl. Math.*, 154(12):1674–1684, 2006.

To conclude this chapter let me mention a recent preprint[10] describing the probabilistic behavior of separable permutations when $n$ tends to infinity. Given a permutation $\pi$ of $\{1,\ldots,n\}$, consider its graph: this is the subset $\{(i,\pi(i))\} \subset \{1,\ldots,n\}^2$. Rescale this picture to draw it in the unit square $[0,1]^2$. To each permutation $\pi$, we associate the probability measure $\mu_\pi$ in the square which is the sum of the $n$ Dirac masses of weight $1/n$ located at $(i/n, \pi(i)/n)$. The space $Prob([0,1]^2)$ of probability measures on the square is compact (for the weak topology) so that we can study the accumulation points of the $\mu_\pi$'s in $Prob([0,1]^2)$. It is easy to see that any accumulation point $\mu$ is a *permuton*: a probability measure on the square whose two marginals (its two projections on the axes), are the Lebesgue measure on $[0,1]$. The preprint by Bassino et al. describes the limits of separable permutations. For each $n$, choose a separable permutation at random (uniformly distributed among all separable permutations). This produces a random probability distribution on $Prob([0,1]^2)$ for each $n$. The authors prove that this sequence of probabilities converges to a well defined probability in the space $Prob([0,1]^2)$. This limit is a random probability on the space of permutons: the *separable random permuton*. The two pictures in the margin, extracted from this preprint show typical graphs of separable permutations for large values of $n$.

[10] F. Bassino, M. Bouvel, V. Féray, L. Gerin, and A. Pierrot. The Brownian limit of separable permutations. 2016.

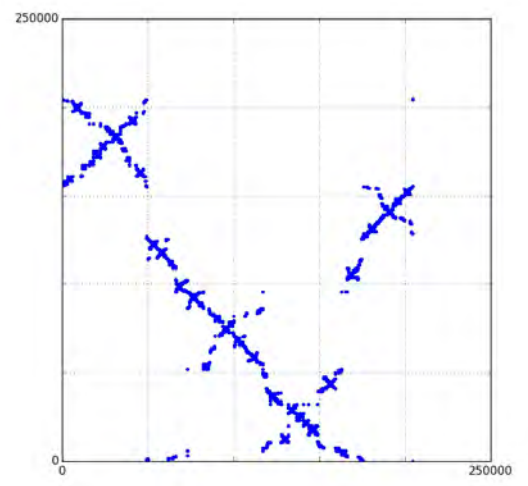

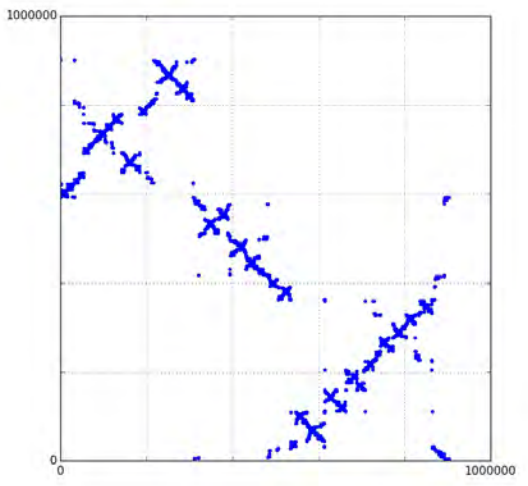

For much more on the combinatorics of separable permutations, see[11] or[12]. Let me however propose one exercise.

*Show that there is an algorithm deciding whether or not a given permutation of* $\{1,2,\ldots,n\}$ *is separable in linear time (in n).*

Note that there is an obvious algorithm in *polynomial time*: for each 4-tuple $1 \le i_1 \le i_2 \le i_3 \le i_4 \le n$ check whether their images are ordered like one of the two forbidden permutation. Going from polynomial time to linear time might be important since $\log n!$ grows faster than a linear function in $n$, but slower than a quadratic function, by Stirling's formula. Therefore, if you can find a linear time algorithm, you prove in particular that the number of separable permutations is small when compared with $n!$. As a hint for this exercise, read again the proof of the bijection between separable permutations and pruned trees.

[11] S. Kitaev. *Patterns in permutations and words*. Monographs in Theoretical Computer Science. An EATCS Series. Springer, Heidelberg, 2011. With a foreword by Jeffrey B. Remmel.

[12] M. Bóna. *Combinatorics of permutations*. Discrete Mathematics and its Applications (Boca Raton). CRC Press, Boca Raton, FL, second edition, 2012. With a foreword by Richard Stanley.

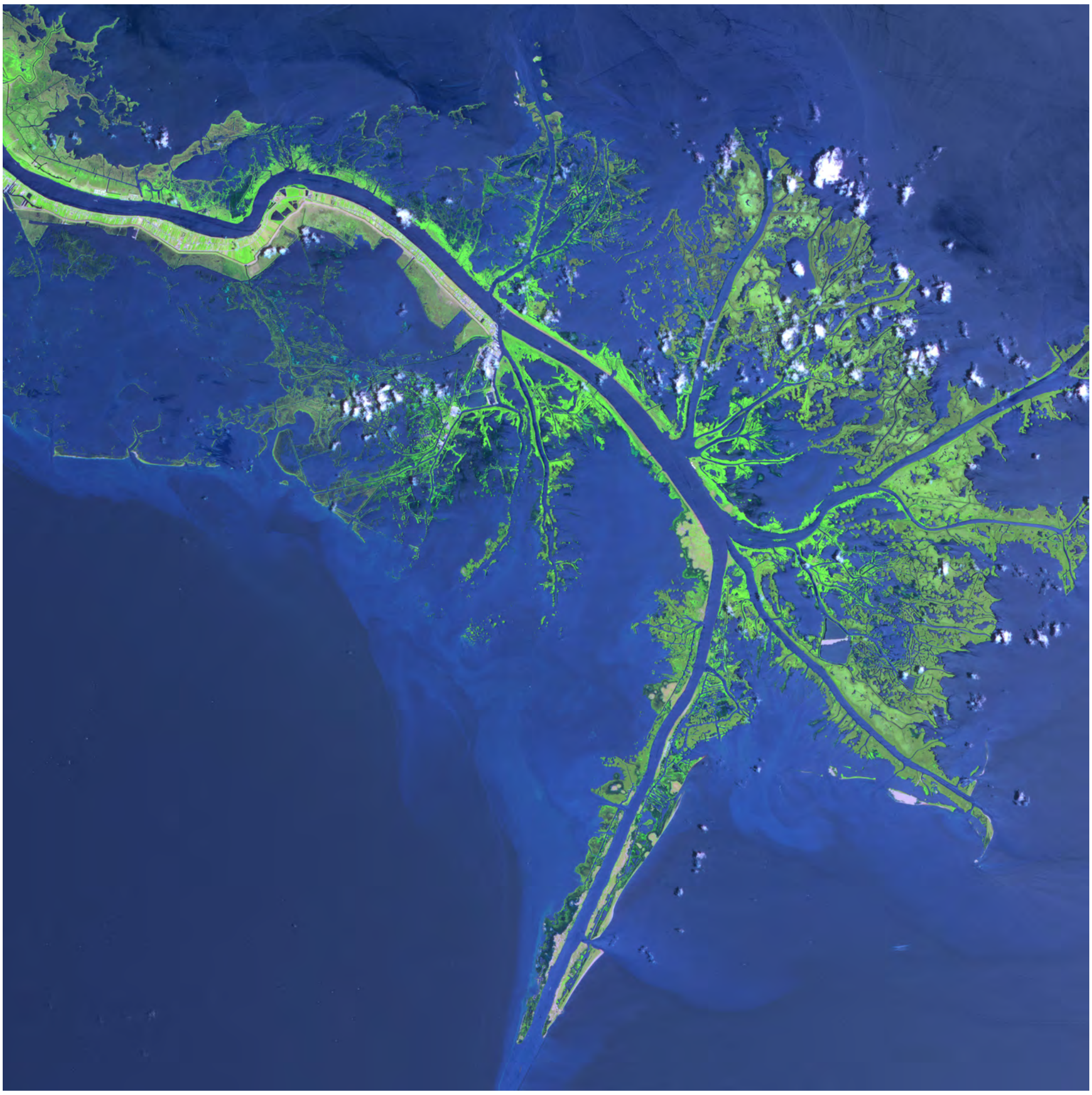

Mississippi River Delta. ©

# *Hipparchus and Schroeder*

## *Let us count trees*

WE ARE GOING TO COUNT the number $a(n)$ of polynomial interchanges (or separable permutations) on $n$ objects.

Let $b(n)$ be the number of pruned trees with $n$ leaves which are such that the root has *at least two children* if $n \geq 2$ (and has no children, if $n = 1$). For such a tree, we can create a new root which has the original root as its unique child. It follows that $a(n) = 2b(n)$ for $n \geq 2$. The first values of $b$ are:

$b(1) = 1$: a tiny tree whose root is also its unique leaf.

$b(2) = 1$: a tiny tree with two branches and two leaves.

$b(3) = 3$.

It is very tempting to establish a recurrence relation for $b(n)$.

Start with a pruned tree with $n$ leaves such that the root has at least two children. If we delete the root and the adjacent branches, we get a certain number of trees, having a total of $n$ leaves. Conversely, starting with an ordered set of at least two pruned trees having $n$ leaves in total, we can add a new root and connect it to the previous roots, in order to construct a pruned tree with $n$ leaves.

Therefore, the following relation holds:

$$b(n) = \sum_{i_1, i_2, \ldots, i_k; i_1 + \cdots + i_k = n} b(i_1) b(i_2) \cdots b(i_k).$$

We now use the classical method of *generating series*. Define the

The first values of $a(n)$:

| | |
|---|---|
| 1 | 1 |
| 2 | 2 |
| 3 | 6 |
| 4 | 22 |
| 5 | 90 |
| 6 | 394 |
| 7 | 1806 |
| 8 | 8558 |
| 9 | 41586 |
| 10 | 206098 |
| 11 | 1037718 |
| 12 | 5293446 |
| 13 | 27297738 |
| 14 | 142078746 |
| 15 | 745387038 |
| 16 | 3937603038 |
| 17 | 20927156706 |
| 18 | 111818026018 |
| 19 | 600318853926 |
| 20 | 3236724317174 |
| 21 | 17518619320890 |
| 22 | 95149655201962 |
| 23 | 518431875418926 |
| 24 | 2832923350929742 |
| 25 | 15521467648875090 |
| 26 | 85249942588971314 |
| 27 | 469286147871837366 |
| 28 | 2588758890960637798 |
| 29 | 14308406109097843626 |
| 30 | 79228031819993134650 |
| 31 | 439442782615614361662 |

formal power series $H$ by:

$$H(t) = \sum_{n=1}^{\infty} b(n)t^n = t + t^2 + 3t^3 + \cdots.$$

Let us square $H$:

$$H(t)^2 = t^2 + 2t^3 + 7t^4 + \cdots.$$

The coefficient of $t^n$ in this new series is $\sum_{i_1+i_2=n} b(i_1)b(i_2)$, which is equal to the number of pruned trees with $n$ leaves such that the root has exactly two children. Using $H(t)^3$, we would count the number of trees whose root has three children, etc.

The infinite series

$$H(t)^2 + H(t)^3 + \cdots$$

counts therefore all trees, except the only one which has a single leaf. Hence, this infinite sum is $H(t) - t$.

Therefore

$$H(t) - t = H(t)^2 + H(t)^3 + \cdots.$$

Summing the geometric series, we get:

$$H(t) - t = \frac{H(t)^2}{1 - H(t)}$$

or

$$2H(t)^2 - (1+t)H(t) + t = 0,$$

which yields:

$$H(t) = \sum_{n=1}^{\infty} b(n)t^n = (1 + t - \sqrt{1 - 6t + t^2})/4.$$

As a function of a complex variable, $(1 + t - \sqrt{1 - 6t + t^2})/4$ is well defined and holomorphic in the disc of center 0 whose radius is the smallest of the two roots of $1 - 6t + t^2 = 0$, which is $t = 3 - 2\sqrt{2}$. The radius of convergence of $H(t)$ is therefore $3 - 2\sqrt{2}$. In other words

$$\limsup_{n\to\infty} \frac{1}{n} \log b(n) = \log(3 + 2\sqrt{2}).$$

The reader will easily show that the lim sup can be replaced by a lim. The $a(n)$'s are the *large Schroeder numbers*, and the

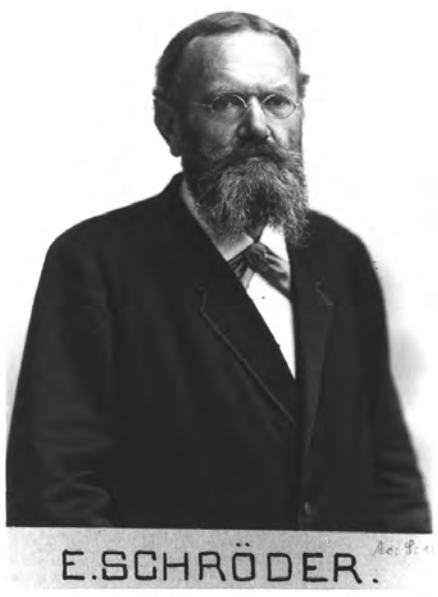

Friedrich Wilhelm Karl Ernst Schroeder (1841-1902) had many sporting hobbies: cycling, hiking, swimming, ice-skating, horseback riding, and gardening. Because he was always seen riding his bicycle around Karlsruhe, he was known locally as the 'Bicycle-professor' (see MacTutor History of Mathematics archive). ©

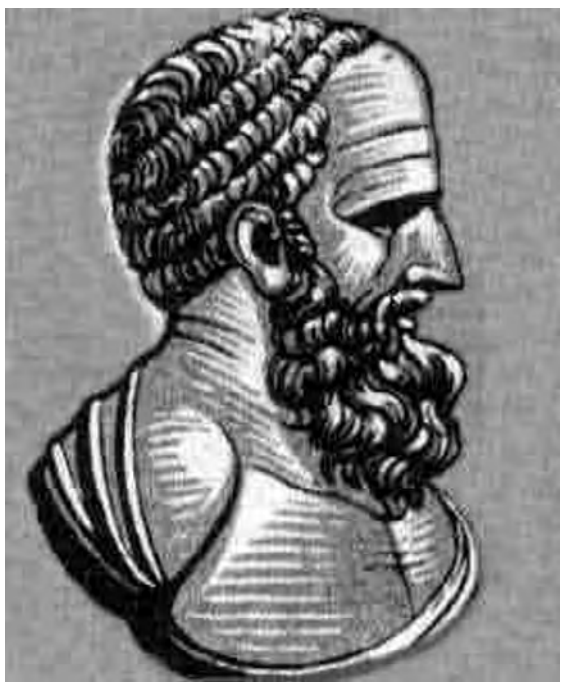

Hipparchus (c. 190 – c. 120 bc.) ©

$b(n)$'s are the *small Schroeder numbers*. Do not forget that $a(n) = 2b(n)$ for $n \geq 2$.

The amazing *On-Line Encyclopedia of Integer Sequences* has several pages dedicated to these two sequences (among many other pages) and contains probably too much information! For instance, one finds the asymptotic estimate:

$$a(n) \sim \frac{(3+2\sqrt{2})^n}{2n\sqrt{2\pi n}\sqrt{3\sqrt{2}-4}(1-\frac{9\sqrt{2}+24}{32n}+\cdots)}.$$

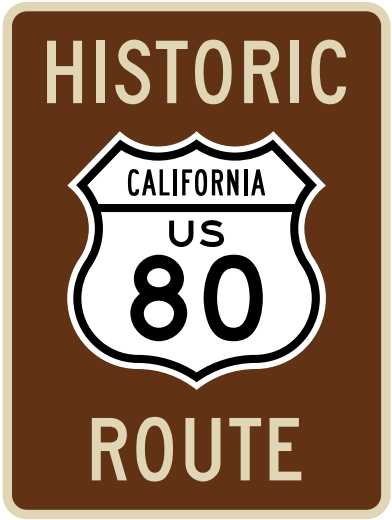

©

## *Hipparchus and Schroeder*

*Arnold's Principle* asserts that

> If a notion bears a personal name, then this is not the name of the discoverer.

and its complement, *Berry's Principle*:

> The Arnold Principle is applicable to itself[13].

[13] V. I. Arnold. On teaching mathematics.

In a different context, this is also called Stigler's law of eponymy. A. N. Whitehead reportedly said "Everything of importance has been said before by somebody who did not discover it".

This applies in particular to the discovery of Schroeder numbers. Ernst Schroeder was an important German logician who explained that his aim was[14]:

> To design logic as a calculating discipline, especially to give access to the exact handling of relative concepts, and, from then on, by emancipation from the routine claims of natural language, to withdraw any fertile soil from 'cliché' in the field of philosophy as well. This should prepare the ground for a scientific universal language that, widely differing from linguistic efforts like Volapük [a universal language like Esperanto, very popular in Germany at the time], looks more like a sign language than like a sound language.

[14] V. Peckhaus. 19th century logic between philosophy and mathematics. *Bull. Symbolic Logic*, 5(4):433–450, 1999.

Given his viewpoint on logic, it was a very natural question for him to count the number of correct bracketings on a word of length $n$. This is the purpose of his 1870 paper[15].

[15] E. Schröder. Vier combinatorische probleme. *Bull. Symbolic Logic*, 15:361–376, 1870.

For a word of length 2, there are two possibilities:

$$ab \quad \text{and} \quad (ab).$$

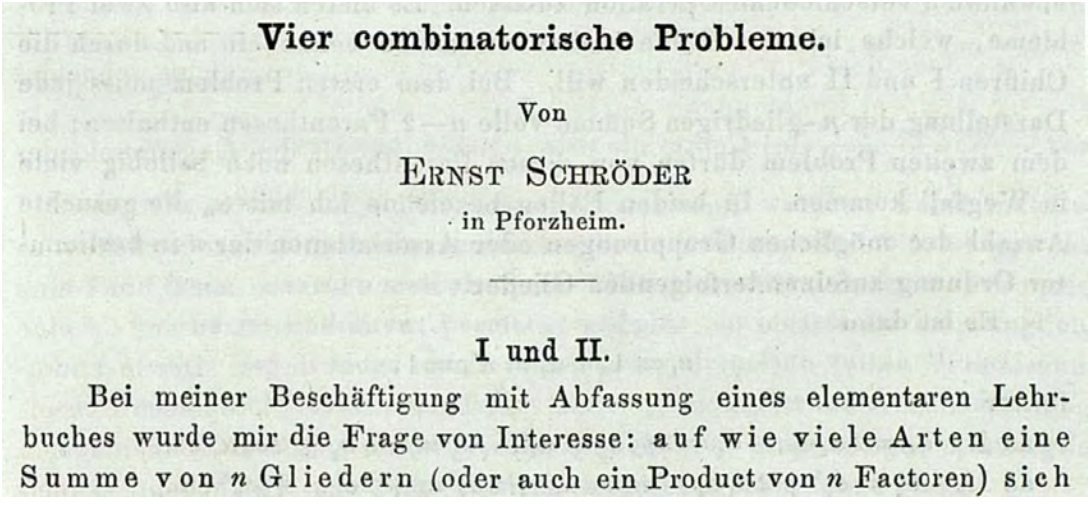

**Vier combinatorische Probleme.**

Von

Ernst Schröder

in Pforzheim.

**I und II.**

Bei meiner Beschäftigung mit Abfassung eines elementaren Lehrbuches wurde mir die Frage von Interesse: auf wie viele Arten eine Summe von $n$ Gliedern (oder auch ein Product von $n$ Factoren) sich

Six possibilities for a word of length 3:

$$abc \quad (ab)c \quad a(bc)$$
$$(abc) \quad ((ab)c) \quad (a(bc)).$$

| | |
|---|---|
| $abcd$ | $(abcd)$ |
| $(ab)cd$ | $((ab)cd)$ |
| $a(bc)d$ | $(a(bc)d)$ |
| $(ab(cd)$ | $((ab(cd))$ |
| $(ab)(cd)$ | $((ab)(cd))$ |
| $(abc)d$ | $((ab)cd)$ |
| $a(bcd)$ | $(a(bcd))$ |
| $((ab)c)d$ | $(((ab)c)d)$ |
| $(a(bc))d$ | $((a(bc))d)$ |
| $a((bc)d$ | $(a((bc)d)$ |
| $(a(b(cd))$ | $((a(b(cd)))$ |

The rules of the game are that a single letter cannot be enclosed in parentheses like $(a)$ and one should not duplicate parentheses like $((ab))$. The full word can be inside a single pair of parentheses or not (and this is why large Schroeder numbers are even). Note that a pair of parentheses can enclose more than two letters. There are 22 possibilities for a word of length 4.

It should be clear to the reader that these 22 expressions are nothing more than the list of the 22 pruned trees with 4 leaves. Indeed, parenthesized words can be described by pruned trees, as shown in the margin. Schroeder was simply counting pruned trees, alias polynomial interchanges, alias separable permutations. His paper contains the recurrence relation and the generating function, described above.

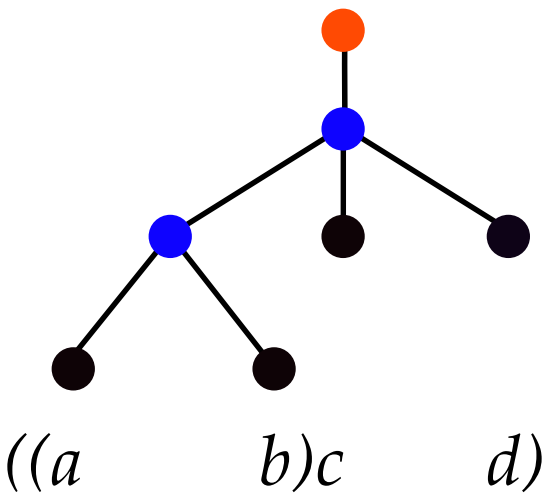

In 1994, David Hough, a graduate student at George Washington University (USA), was reading exercise 1.45 in Stanley's book[16]:

> The following quotation is from Plutarch's *Table Talk* VIII.9.732 'Chrysippus says that the number of compound propositions that can be made from only ten simple propositions exceeds a million'. (Hipparchus, to be sure, refuted this by showing that on the affirmative side there are 103,049 compound statements, and on the negative side 310,952.
>
> According to Y. Heath, A History of Greek Mathematics. vol 2, p. 245; "it seems impossible to make anything of these figures." [Heath also notes that a variant reading of 103,049 is 101,049.]
>
> Can in fact any sense be made of Plutarch's statement?

[16] R. P. Stanley. *Enumerative combinatorics. Vol. I.* The Wadsworth & Brooks/Cole Mathematics Series. Wadsworth & Brooks/Cole Advanced Books & Software, Monterey, CA, 1986. With a foreword by Gian-Carlo Rota.

Hough noticed that 103,049 is the tenth small Schroeder number $b(10)$ and it could not be a coincidence.

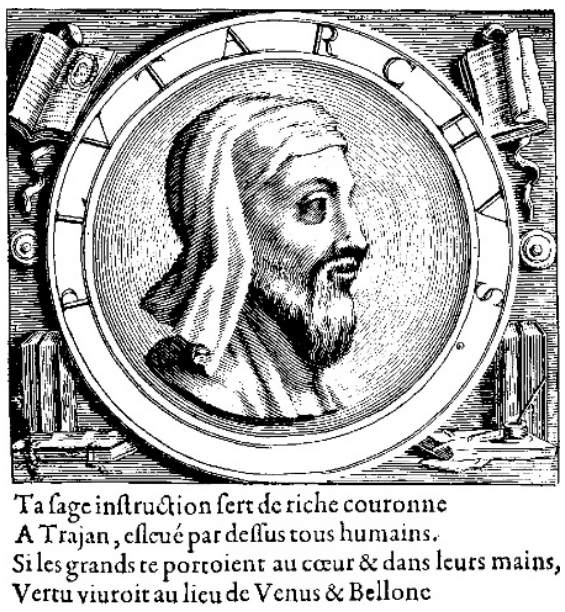

Plutarch
(c. AD 46 – AD 120). ©

Plutarch was a Greek historian and biographer, whose role (in our story) is limited to the retranscription of a quotation by Hipparchus two hundred years earlier. It is hard to imagine that this number 103,049 could have been remembered during such a long time without having been preserved in a book in the possession of Plutarch.

Hipparchus is probably the most important astronomer of ancient times. He is well known for his discovery of the precession of equinoxes but above all for the construction of a coherent scientific description of the motion of planets. His successor, Ptolemy, three hundred years later, is famous for the *Ptolemaic geocentric system* which became the astronomical dogma until Copernicus introduced the heliocentric system, many centuries later. Ptolemy owes a lot to Hipparchus and does not always acknowledge his debt. But that is not related to our story.

So, according to Hough, Hipparchus, under the transmission of Plutarch, was counting parenthesized words of length 10. Several historical papers have been written about this discovery of Schroeder numbers by Hipparchus[17] [18].

[17] R. P. Stanley. Hipparchus, Plutarch, Schröder, and Hough. *Amer. Math. Monthly*, 104(4):344–350, 1997.

[18] F. Acerbi. On the shoulders of Hipparchus: a reappraisal of ancient Greek combinatorics. *Arch. Hist. Exact Sci.*, 57(6):465–502, 2003.

An article in MathPages[19] provides a slightly more elaborate explanation in terms of Stoic logic (some pre-Aristotelean logic taught in particular by Chrysippus, and criticized by Hipparchus).

[19] K. Brown. Hipparchus on compound statements. MathPages.

Given a certain number of logical assertions $a_1, a_2, \dots, a_k$ there are at least two ways to combine them by conjunction or disjunction:

- $a_1$ OR $a_2$ OR $\cdots$ OR $a_n$, which is an $n$-ary function $OR(a_1, \dots, a_n)$,
- $a_1$ AND $a_2$ AND $\cdots$ AND $a_n$, which is an $n$-ary function $AND(a_1, \dots, a_n)$.

In modern Boolean notation, we use + for OR and a dot, or just concatenation, for AND. Now, consider a word of length $n$ (for example $abcd$ for $n = 4$). For each of the $n - 1$ spaces between letters, choose a "+" or a ".". There are $2^{n-1}$ possibilities (8 in our example).

| | |
|---|---|
| $abcd$ | $abc + d$ |
| $ab + cd$ | $ab + c + d$ |
| $a + bcd$ | $a + bc + d$ |
| $a + b + cd$ | $a + b + c + d$ |

We are used to give priority to multiplication above addition but if we want to specify an order to evaluate this logical function, parentheses are necessary.

This can be described by a pruned tree. Associate the symbol OR to a node of odd level and AND to a node of even level. Each node acts accordingly on the set of its children, these children being embraced in a pair of parentheses. If the root has a single child, it is not necessary to label it with AND since it acts on a singleton.

The reader is encouraged to show, as an exercise, that two different expressions, i.e. two different pruned trees, define two different Boolean functions $\{0,1\}^n \to \{0,1\}$ when evaluated at $a_i = 0$ or $1$ (false or true). So, Hipparchus was right: there are $a(10) = 2 \times 103{,}049$ ways of combining 10 assertions, using OR or AND, in the sense just described. One could ask why he mentioned $b(10)$ and not $a(10)$. He may have noticed the natural involution among *compound propositions* given by negation, which basically permutes AND and OR?

There is a related open question, called the Dedekind problem. There are $2^{2^n}$ Boolean functions, that is to say functions from $\{0,1\}^n$ to $\{0,1\}$. It is easy to see that any such Boolean function can be written by some formula using OR, AND and NOT. Those functions that can be described by formulas which are not involving NOT are called *monotone Boolean functions* (but note that we do not impose that each variable appears once in the formula, like in the case of Hipparchus). *The question of computing the number of monotone Boolean functions is open for* $n > 8$. This number has also a nice topological interpretation: it is the number of simplicial complexes whose vertices are $\{1,2,\ldots,n\}$.

Most mathematicians, including myself, have a naive idea about Greek mathematics. We believe that it only consists of Geometry, in the spirit of Euclid. The example of the computation by Hipparchus of the tenth Schroeder number may be a hint that the Ancient Greeks had developed a fairly elaborate understanding of combinatorics: this is the theme of the article by Acerbi quoted above.

| | |
|---|---|
| $abcd$ | $abc+d$ |
| $ab(c+d)$ | $a(bc+d)$ |
| $ab+cd$ | $a(b+c)d$ |
| $a(b+cd)$ | $(ab+c)d$ |
| $ab+c+d$ | $a(b+c)+d$ |
| $a(b+c+d)$ | $a+bcd$ |
| $(a+b)cd$ | $(a+bc)d$ |
| $a+bc+d$ | $(a+b)(c+d)$ |
| $(a+b)c+d$ | $a+b(c+d)$ |
| $a+(b+c)d$ | $a+b+cd$ |
| $(a+b+c)d$ | $a+b+c+d$ |

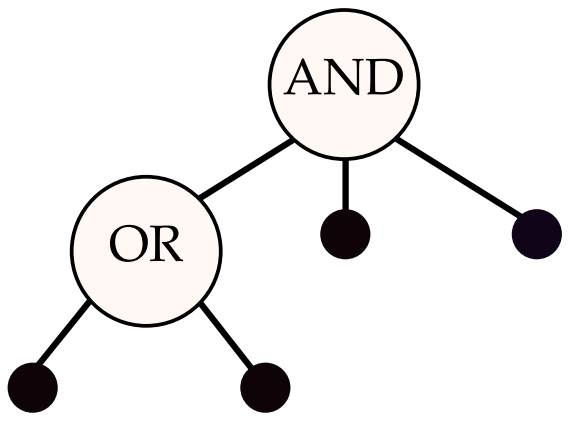

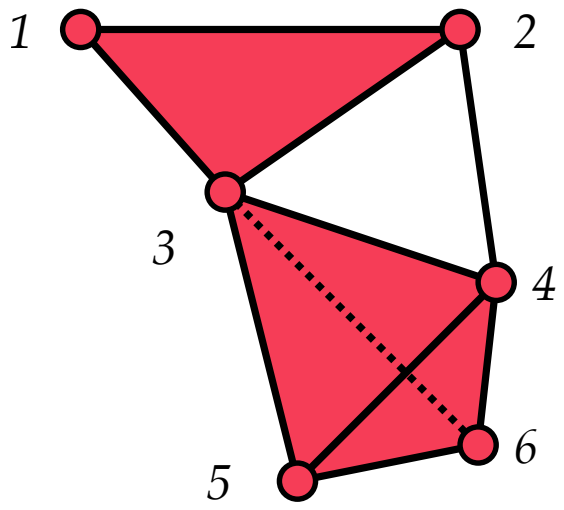

A simplicial complex with vertex set $V = \{1,\ldots,6\}$. The function associating 0 to a subset $X \subset V$ if $X$ is a simplex, and 1 otherwise, is a monotone Boolean function.

The book by Netz[20] offers new perspectives on this history. The first chapter discusses Greek combinatorics and in particular the number 103,049. It also contains a description of another combinatorial puzzle, found in the famous Archimedes palimpsest (the reader is urged to read[21] like a detective story). This is made out of 14 polygonal pieces and is similar to a *Tangram* game.

Netz *"asked his colleague at Stanford, Persi, a noted combinatorist, to help him solve what he assumed to be a simple question: how many ways are there to put together the square? [. . . ] It took Diaconis a couple of months and collaborative work with three colleagues to come up with the number of solutions: 17,142."*

Did Archimedes know the answer?

[20] R. Netz. *Ludic proof.* Cambridge University Press, Cambridge, 2009. Greek mathematics and the Alexandrian aesthetic.

[21] R. Netz and W. Noel. *The Archimedes codex*. Phoenix, London, 2008. Revealing the secrets of the world's greatest palimpsest.

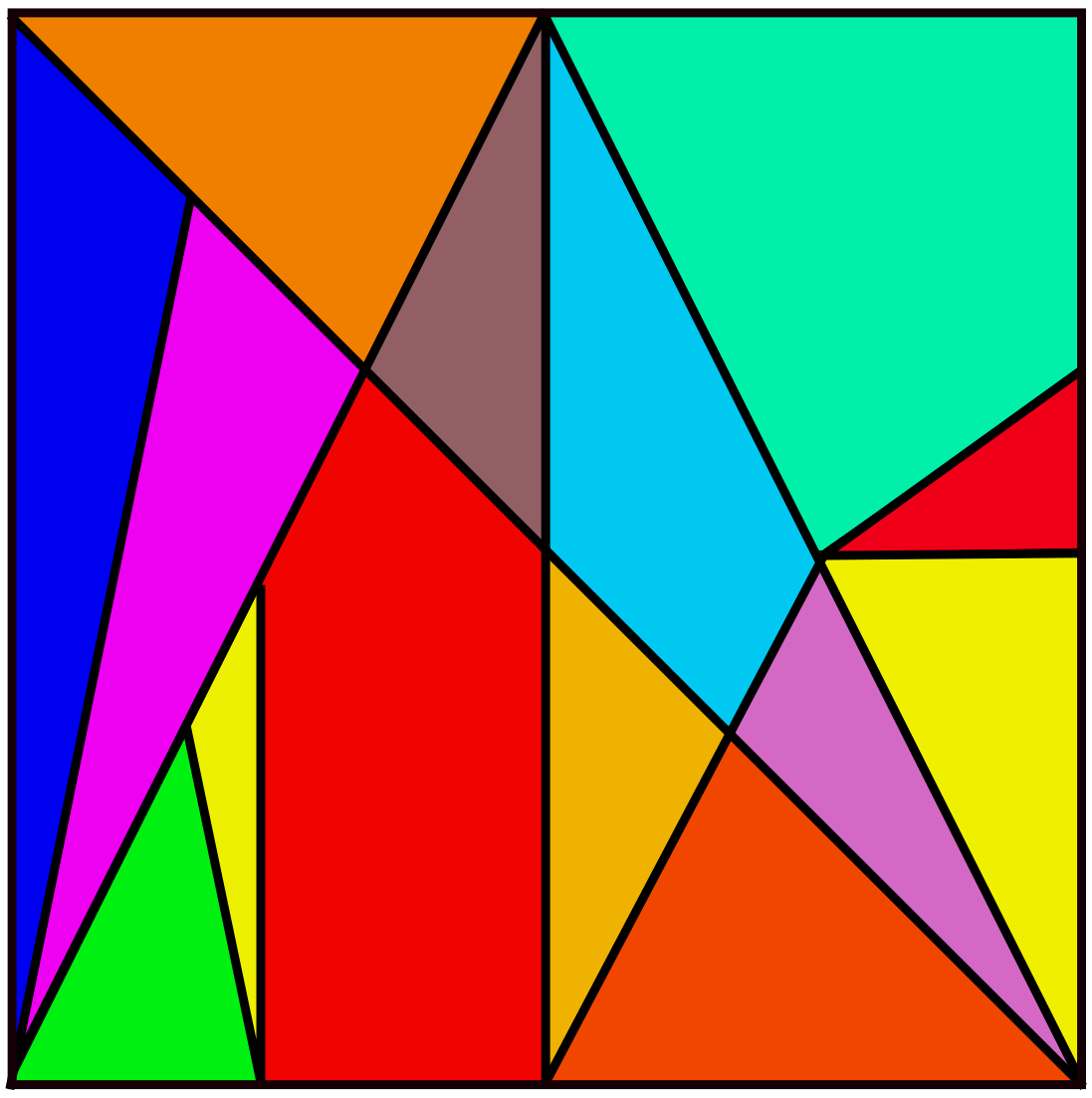

Archimedes' Stomachion.

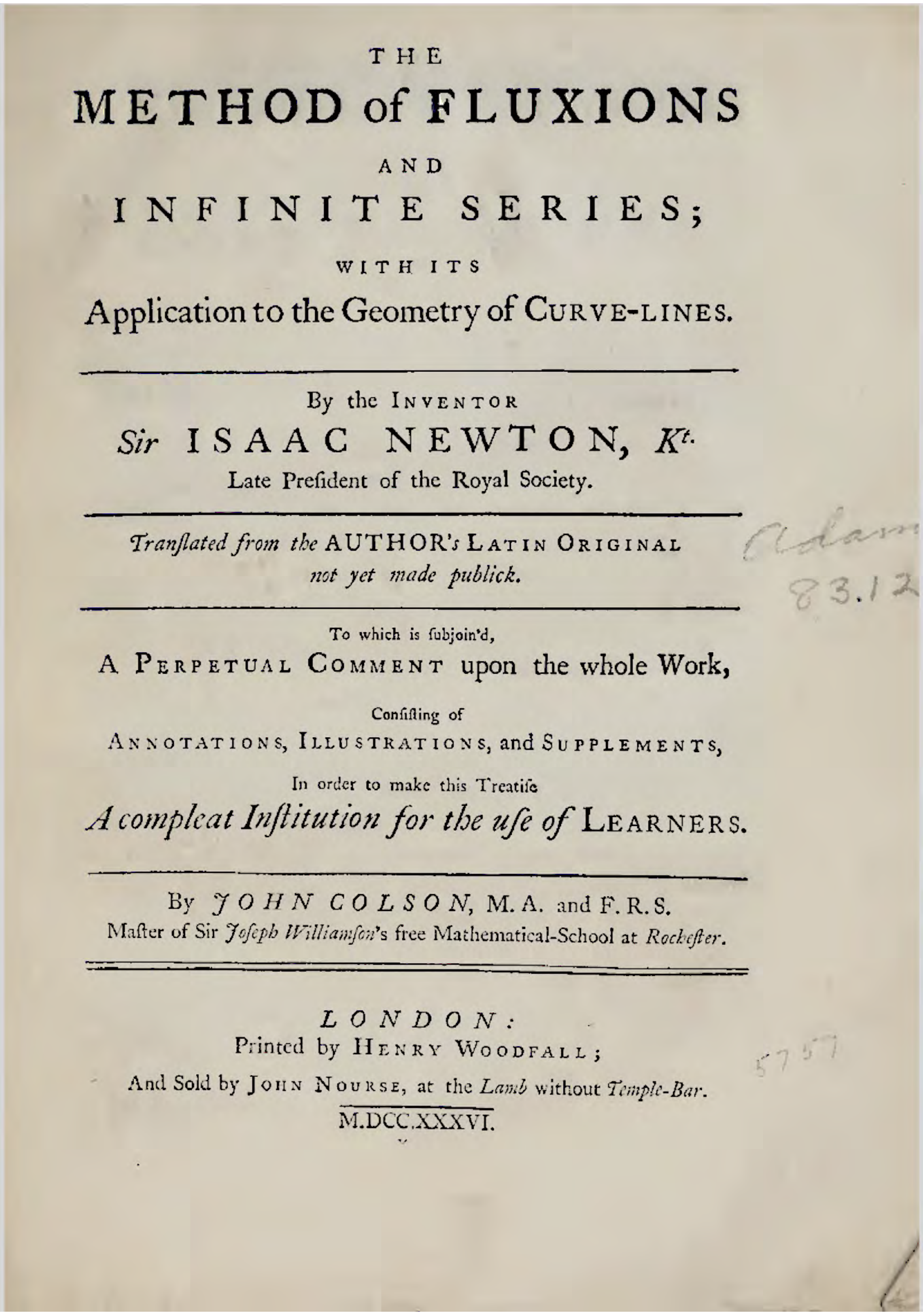

THE

METHOD of FLUXIONS

AND

INFINITE SERIES;

WITH ITS

Application to the Geometry of CURVE-LINES.

By the INVENTOR

*Sir* ISAAC NEWTON, *K^t.*

Late Prefident of the Royal Society.

*Tranflated from the* AUTHOR's LATIN ORIGINAL
*not yet made publick.*

To which is fubjoin'd,

A PERPETUAL COMMENT upon the whole Work,

Confifting of

ANNOTATIONS, ILLUSTRATIONS, and SUPPLEMENTS,

In order to make this Treatife

*A compleat Inftitution for the ufe of* LEARNERS.

By *JOHN COLSON*, M. A. and F. R. S.
Mafter of Sir *Jofeph Williamfon*'s free Mathematical-School at *Rochefter*.

*LONDON:*
Printed by HENRY WOODFALL;
And Sold by JOHN NOURSE, at the *Lamb* without *Temple-Bar*.
M.DCC.XXXVI.

# *De methodis serierum et fluxionum*
# *Newton's method*

## *Algebraic curves*

SINCE THE INTRODUCTION OF COORDINATES BY RENÉ DESCARTES, the study of planar curves, especially *planar algebraic curves* (defined by some polynomial equation $P(x,y)=0$), has become a central theme in mathematics and continues to be so. Of course, equations of degree 1 and 2 (lines and conics) were very familiar. When XVII-th century mathematicians looked at higher degree curves, they found a jungle, consisting of many different shapes that they tried to tame. For instance, Isaac Newton wrote a long memoir on curves of degree 3, decomposing them in a great number of "species". See the discussion in[22] or in[23].

[22] J. Stillwell. *Mathematics and its history*. Undergraduate Texts in Mathematics. Springer, New York, third edition, 2010.

[23] W. W. R. Ball. On Newton's Classification of Cubic Curves. *Proc. London Math. Soc.*, S1-22(1):104–143, 1890.

Very quickly, it appeared clearly that singular points play a central role in the understanding of the geometry of these curves. A point $(x_0,y_0)$ is *singular* if it lies on the curve, i.e. if $P(x_0,y_0)=0$, and the partial derivatives $\partial P/\partial x$ and $\partial P/\partial y$ vanish at $(x_0,y_0)$. In a neighborhood of a *regular* (i.e. non-singular) point, a modern mathematician has no difficulty applying the implicit function theorem: in suitable smooth coordinates around such a point, the curve is a straight line.

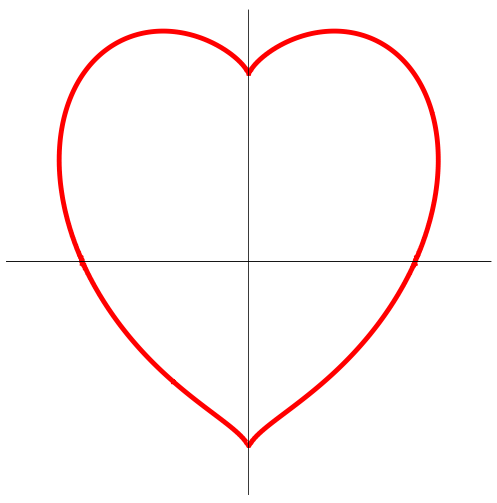

An algebraic curve with two singular points.

Singular points might however be much more complicated and deciphering their nature took a long time.

In this chapter, I am going to describe one of the major steps toward this understanding, following Newton's book *De methodis serierum et fluxionum*.

I do not want to go into any historical detail about the rivalry between Newton and Leibniz concerning the invention of calculus. Let me recommend specifically for our purpose the excellent Newton's biography by Westfall[24].

[24] R. S. Westfall. *Never at rest*. Cambridge University Press, Cambridge, 1980. A biography of Isaac Newton.

Here is the description by Newton himself of his *Annus mirabilis* (see[25]):

[25] R. S. Westfall. Newton's marvelous years of discovery and their aftermath: myth versus manuscript. *Isis*, 71(256):109–121, 1980.

> In the beginning of the year 1665 I found the Method of approximating series and the Rule for reducing any dignity of any Binomial into such a series. The same year in May I found the method of Tangents of Gregory and Slusius, and in November had the direct method of fluxions and the next year in January had the theory of Colours and in May following I had entrance into ye inverse method of fluxions. And the same year I began to think of gravity extending to ye orb of the Moon. All this was in the two plague years of 1665–1666. For in those days I was in the prime of my age for invention and minded Mathematicks and Philosophy more then at any time since.

In July 1669, based on his 1665 ideas, Newton had written *De Analysi per aequationes numero terminorum infinitas*.

In 1671, he wrote *De methodis serierum et fluxionum* but did not publish it.

In 1676, he wrote two famous letters to Leibniz (through Oldenburg, as an intermediary): *espistola prior* and *epistola posterior*.

The English translation (by Colson) of *De methodis* appeared in 1736 (therefore 9 years after Newton's death). A French translation of the English translation, by Buffon, appeared in 1740.

It seems that Colson did not accept to show the latin manuscript to Buffon.

All these papers contain a rather precise description of singular points of algebraic curves, in terms of what is called today *Puiseux series*, following once again Arnold's principle.

We want to study a curve $P(x,y) = 0$, where $P$ is a polynomial with *complex coefficients*. One should understand first that Newton is not thinking of this as a *curve*, but as a *function*: given $x$, he wants to solve the equation $P(x,y) = 0$ and to find $y$ as a

Later, we will study the case of real coefficients as well as functions $P$ which are only assumed to be analytic.

function $y(x)$. His main result is that it is indeed possible, as soon as one is willing to consider $y(x)$ as an infinite series *in rational powers* of $x$. Let me state a theorem that will be made precise later on, and that Newton "almost" proved.

**Theorem.** *Any polynomial equation $P(x,y) = 0$ (where $P$ is not divisible by $x$) such that $P(0,0) = 0$ is equivalent, in the neighborhood of $(0,0)$, to a finite number of equalities $y = f_i(x)$ (with $i = 1,\dots,n$) where $f_i$ is a* Puiseux series *of the form:*

$$f_i(x) = \sum_{k=1}^{\infty} a_{i,k} x^{\frac{k}{m_i}}$$

*for some complex coefficients $a_{i,k}$ and some positive integers $m_i$.*

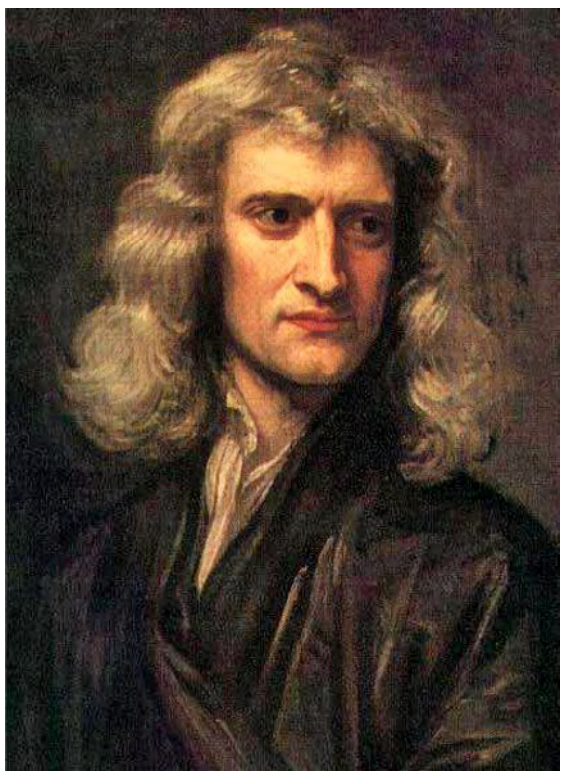

Newton in 1689, by Godfrey Kneller. ©

In other words, $\{P = 0\}$ is the union of a finite number of *graphs* of series $f_i$. We are in a position similar to Kontsevich's original question and it will be natural to ask ourselves what is the topological nature of these graphs. However, before we study this question, there are many details to be fixed, since in particular these $f_i$'s are not really *functions*. Think for instance of the "graph" of the square root.

We will look closely at the first part of *De methodis serierum et fluxionum*. The frontispiece of this important book is on the first page of this chapter. In order to simplify my readers' task, I shall follow the English translation.

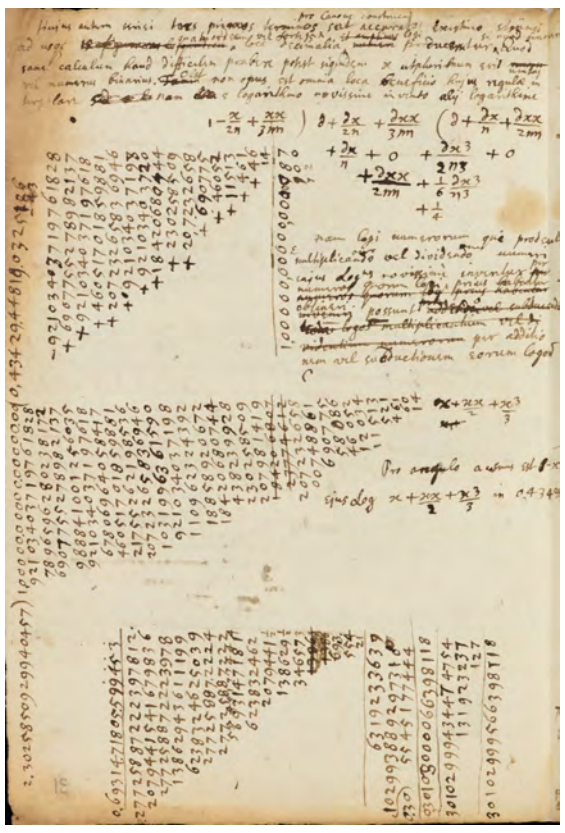

A page from *De methodis.* ©

## *Newton's method*

Let us start reading Newton.

> Since there is a great conformity between the Operations in Species, and the same Operations in common Numbers; nor do they seem to differ, except in the Characters by which they are represented, the first being general and indefinite, and the other definite and particular: I cannot but wonder that no body has thought of accommodating the lately discover'd Doctrine of Decimal Fractions in like manner to Species [. . . ] especially since it might have open'd a way to more abstruse Discoveries.

*Explanation*: By *common number*, Newton means…a common number, that is to say what we call today a *complex number*. Note that very few mathematicians at that time would consider these numbers as “common”. By *species*, he means a polynomial in $x$, or an entire series, or what is called today a *Laurent series*, or maybe a *Puiseux series*, i.e. a series in rational powers of $x$. In any case, in Newton’s words a *species* is some kind of function.

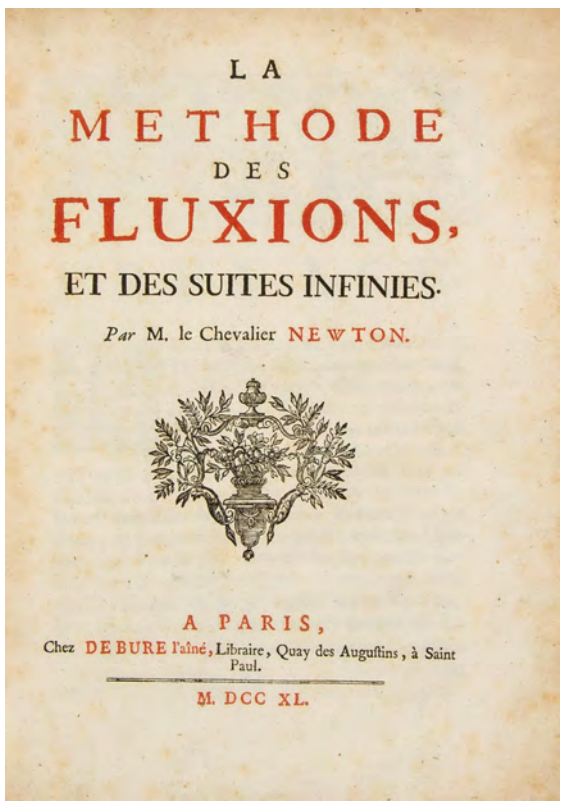

Cover page of the French translation by Buffon. ©

> But since this Doctrine of Species, has the same relation to Algebra, as the Doctrine of Decimal Numbers has to common Arithmetick; the Operations of Addition, Subtraction, Multiplication, Division, and Extraction of Roots, may easily be learned from thence, if the Learner be but skill’d in Decimal Arithmetick, and the Vulgar Algebra, and observes the correspondence that obtains between Decimal Fractions and Algebraick Terms infinitely continued.

*Explanation*: Newton observes that series can be manipulated just in the same way as numbers, for which we have four operations $(+,-,\times,/)$. In modern terminology, he observes that *common numbers* and *Laurent series* are both fields.

> For as in Numbers, the Places towards the right-hand continually decrease in a Decimal or Subdecuple Proportion; so it is in Species respectively, when the Terms are disposed, (as is often enjoin’d in what follows) in an uniform Progression on infinitely continued, according to the Order of the Dimensions of any Numerator or Denominator.

*Explanation*: Again in modern anachronistic terminology, Newton tells us about the topology of these two fields. Two real numbers are close if their decimal expansions agree until a large rank and, analogously, two polynomials in $x$, or two series, are close in the neighborhood of 0 if the valuation of their difference is large.

At this stage, we can guess Newton’s strategy. He will teach us a way of solving polynomial equations $P(y) = 0$ where $P$ is a polynomial with coefficients in some field, which could consist either of *common numbers* or of *species*. This will therefore apply to equations of the form $P(x,y) = 0$ where $P$ is a polynomial in two variables, seen as a polynomial $P(x)(y)$ in one variable $y$ with coefficients in the field of rational functions $\mathbb{C}(x)$, or the field $\mathbb{C}((x))$ of Laurent series.

In a very pedagogical presentation, Newton gives several examples showing the analogy between species and common numbers. First he shows how to expand $a^2/(b-x)$ as a series in $x$. Just as we would explain in primary school that dividing 1 by $0.9 = 1 - 0.1$ yields $1.11111111\ldots$. This is easy and must have been easy also for his readers.

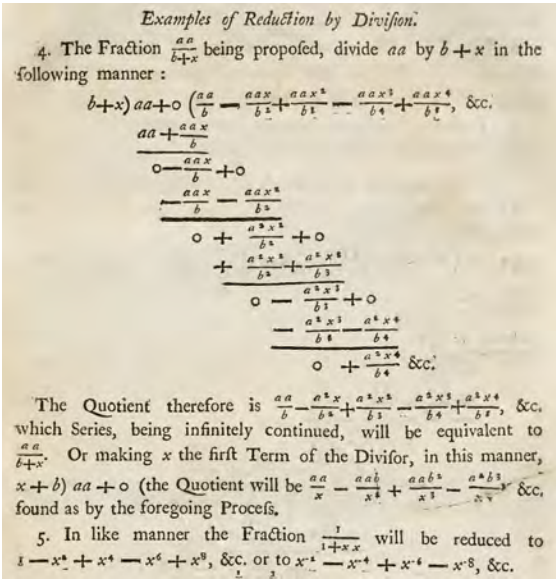
*Examples of Reduction by Division.*

4. The Fraction $\frac{aa}{b+x}$ being propofed, divide $aa$ by $b+x$ in the following manner:

$$b+x)\,aa+0\left(\frac{aa}{b}-\frac{aax}{b^2}+\frac{aax^2}{b^3}-\frac{aax^3}{b^4}+\frac{aax^4}{b^5}, \text{ \&c.}\right.$$

The Quotient therefore is $\frac{aa}{b}-\frac{a^2x}{b^2}+\frac{a^2x^2}{b^3}-\frac{a^2x^3}{b^4}+\frac{a^2x^4}{b^5}$, &c. which Series, being infinitely continued, will be equivalent to $\frac{aa}{b+x}$. Or making $x$ the firft Term of the Divifor, in this manner, $x+b)\,aa+0$ (the Quotient will be $\frac{aa}{x}-\frac{aab}{x^2}+\frac{aab^2}{x^3}-\frac{a^2b^3}{x^4}$, &c. found as by the foregoing Procefs.

5. In like manner the Fraction $\frac{1}{1+xx}$ will be reduced to $1-x^2+x^4-x^6+x^8$, &c. or to $x^{-2}-x^{-4}+x^{-6}-x^{-8}$, &c.

Then he explains the meaning of rational exponents, which must also have been familiar to most of his readers. He can then present his famous binomial formula through the example of the square root of $a^2+x^2$, as an infinite series in $x$:

$$(a^2+x^2)^{\frac{1}{2}} = a + \frac{x^2}{2a} - \frac{x^4}{8a^3} + \frac{x^6}{16a^5} \;\&c.$$

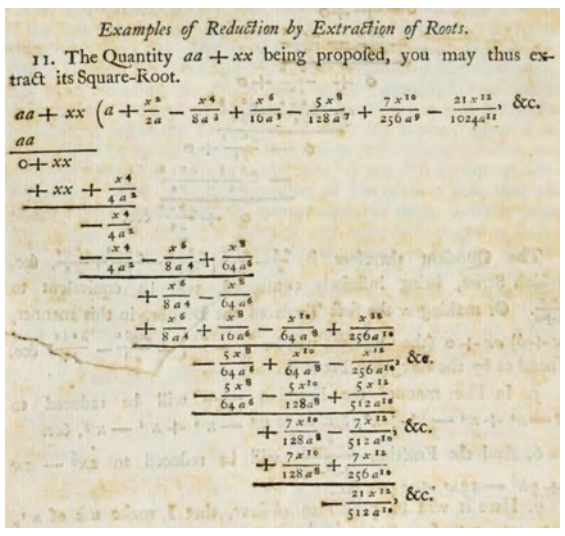
*Examples of Reduction by Extraction of Roots.*

11. The Quantity $aa+xx$ being propofed, you may thus extract its Square-Root.

$$aa+xx\left(a+\frac{x^2}{2a}-\frac{x^4}{8a^3}+\frac{x^6}{16a^5}-\frac{5x^8}{128a^7}+\frac{7x^{10}}{256a^9}-\frac{21x^{12}}{1024a^{11}}, \text{ \&c.}\right.$$

We now come to the part which is the most interesting for us.

He would like to solve what he calls "affected equations" which are polynomial equations whose coefficients are Species, that is to say equations $P(x,y) = P(x)(y) = 0$. Again in a very pedagogical way he declares that he will begin by solving *ordinary* equations in common numbers of the form $P(y) = 0$, where $P$ is a polynomial in $\mathbb{C}[y]$.

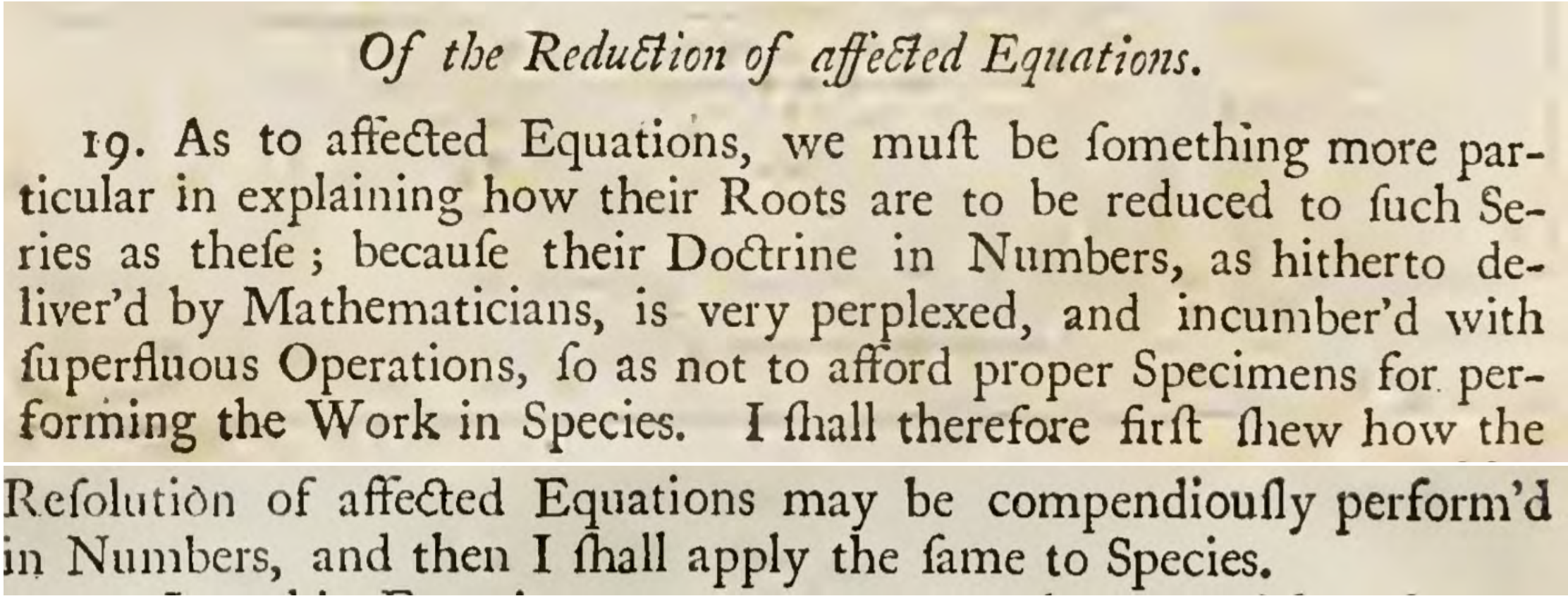
*Of the Reduction of affected Equations.*

19. As to affected Equations, we muft be fomething more particular in explaining how their Roots are to be reduced to fuch Series as thefe; becaufe their Doctrine in Numbers, as hitherto deliver'd by Mathematicians, is very perplexed, and incumber'd with fuperfluous Operations, fo as not to afford proper Specimens for performing the Work in Species. I fhall therefore firft fhew how the Refolution of affected Equations may be compendioufly perform'd in Numbers, and then I fhall apply the fame to Species.

This is the famous *Newton's method* which is one of the most fundamental tools in analysis.

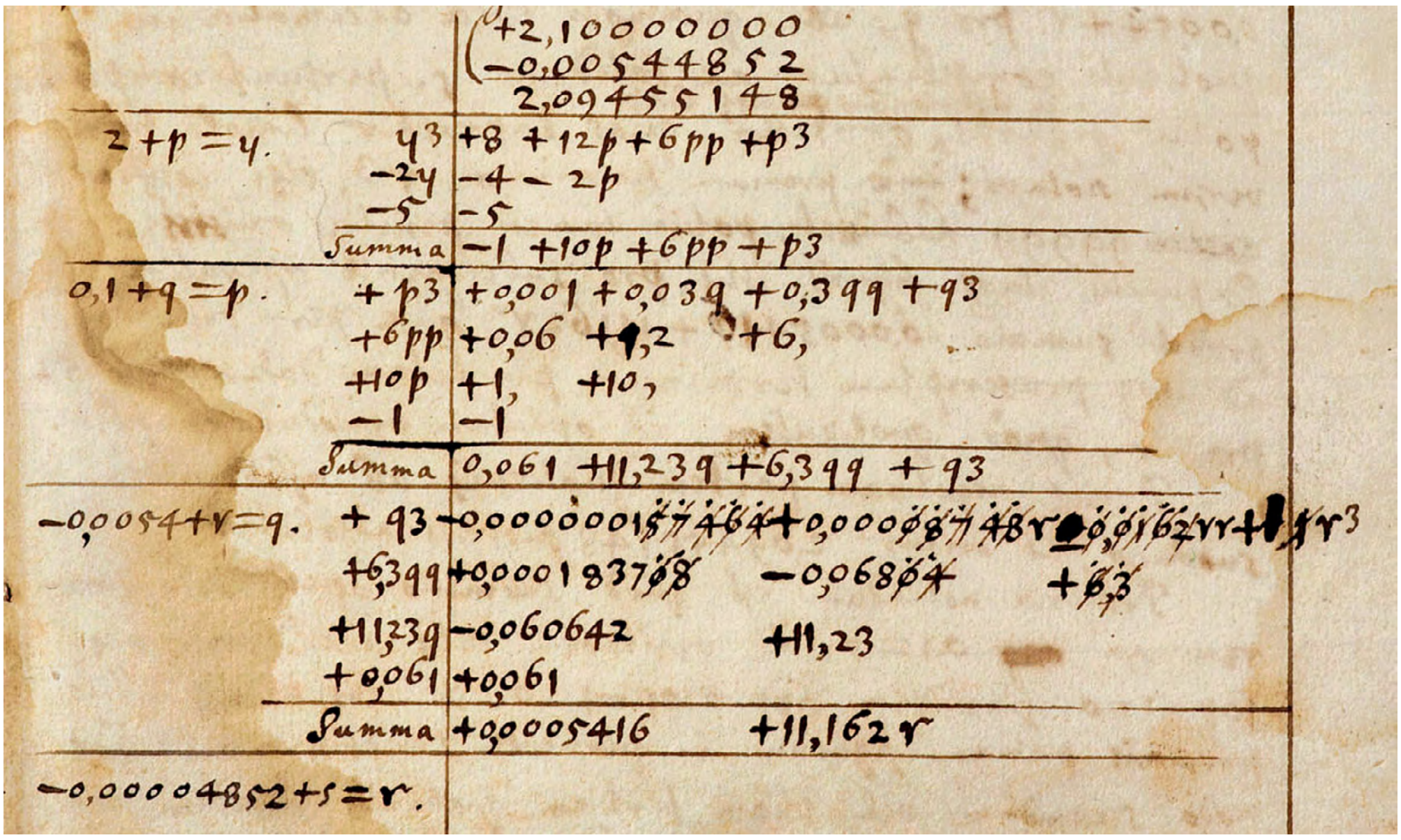

Look at the way he presents the computations. His example is the cubic equation

$$y^3 - 2y - 5 = 0.$$

He observes, by trial and error, that there is a root which is not very different from 2. Therefore he looks for $y$ as $y = 2 + p$ with a small $p$. Substituting in the original equation, he finds

$$p^3 + 6p^2 + 10p - 1 = 0.$$

He can now "reject $p^3 + 6p^2$ because of its smallness" to obtain

$$10p - 1 \simeq 0$$

so that $p$ is close to $1/10$. He can then set $p = 0.1 + q$ and substitute in the equation to get

$$q^3 + 6.3q^2 + 11.23q + 0.061 = 0$$

and "since $11.23q + 0.061 = 0$ is near the truth" he knows that $q$ is close to $-0.061/11.23 \simeq -0.0054$.

Writing $q = -0.0054 + r$, he can substitute as before and continue the operation "as far as I please". He finally gets the solution close to 2.09455148.

The next paragraph shows, if necessary, that Newton was incredibly gifted for computations. "The work can be most abbreviated" indeed. He explains what all undergraduate students (should) know: that at every step the number of correct decimals is essentially doubled and that it is therefore not necessary to compute exactly the $p, q, r, s, etc$. This is why in his table, some digits are barred: this is not a blunder, this is a clever simplification.

In 1690, Raphson (1648–1715) (fellow of the Royal Society, and therefore knowing very well Newton) published a method for solving equations in *Analysis aequationum universalis*. Start with an approximate solution $y_0$ of $P(y) = 0$ and consider the sequence $(y_k)_{k\geq 0}$ defined by

$$y_{k+1} = y_k - \frac{P(y_k)}{P'(y_k)}$$

which, if everything works fine, converges to a solution. Raphson does not mention Newton. Some historians claim that the two methods are *very* different. In the case of Raphson, one keeps the same equation and computes the sequence $y_k$. In the case of Newton, at each step one computes a new equation. The two methods give exactly the same result and are formally identical, but clearly if one computes by hand, Newton's presentation is much more efficient. One could say that Raphson is using an *iteration* and Newton a *recursion*[26]. Some mathematicians claim that Raphson understood the role of the derivative of $P$ and that Newton was only linearizing the equation. Well, who could say that Newton, the inventor of derivative, could not have noticed that the linear part is the derivative? As far as I am concerned, I will continue speaking of *Newton's method* and not of Newton-Raphson's method.

As a final comment, needless to say that Newton does not discuss at all any question about the convergence. Note also that his example only involves real roots of real polynomials.

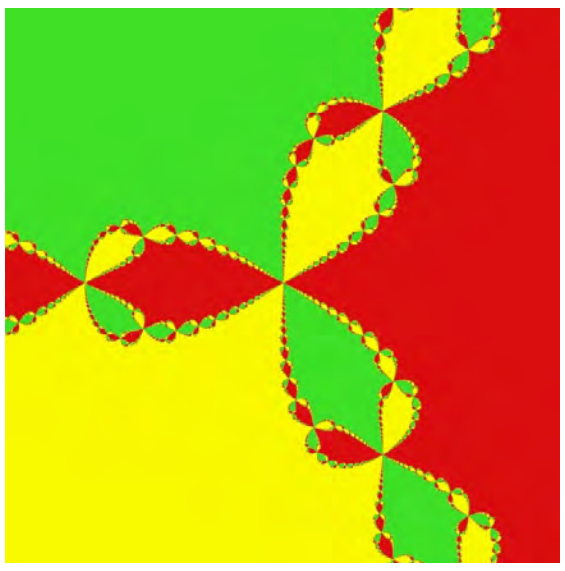

©

Newton's method can be used for finding roots of polynomials $P(z)$ with complex coefficients. Starting with some $z_{init}$, we hope that the iteration of Newton's algorithm will converge to a root. The plane (or at least the set of $z_{init}$ for which the method works) can therefore be decomposed in several domains, according to the limiting root.

In 1880, Cayley asked for a description of this decomposition. He wrote that the question is easy in degree 2 (exercise for the reader) and that for degree 3 it is "anything but obvious". Indeed, it is known today that this decomposition has a fractal nature. This is known as *Newton's rabbit*.

[26] C. Christensen. Newton's method for resolving affected equations. *College Math. J.*, 27(5):330–340, 1996.

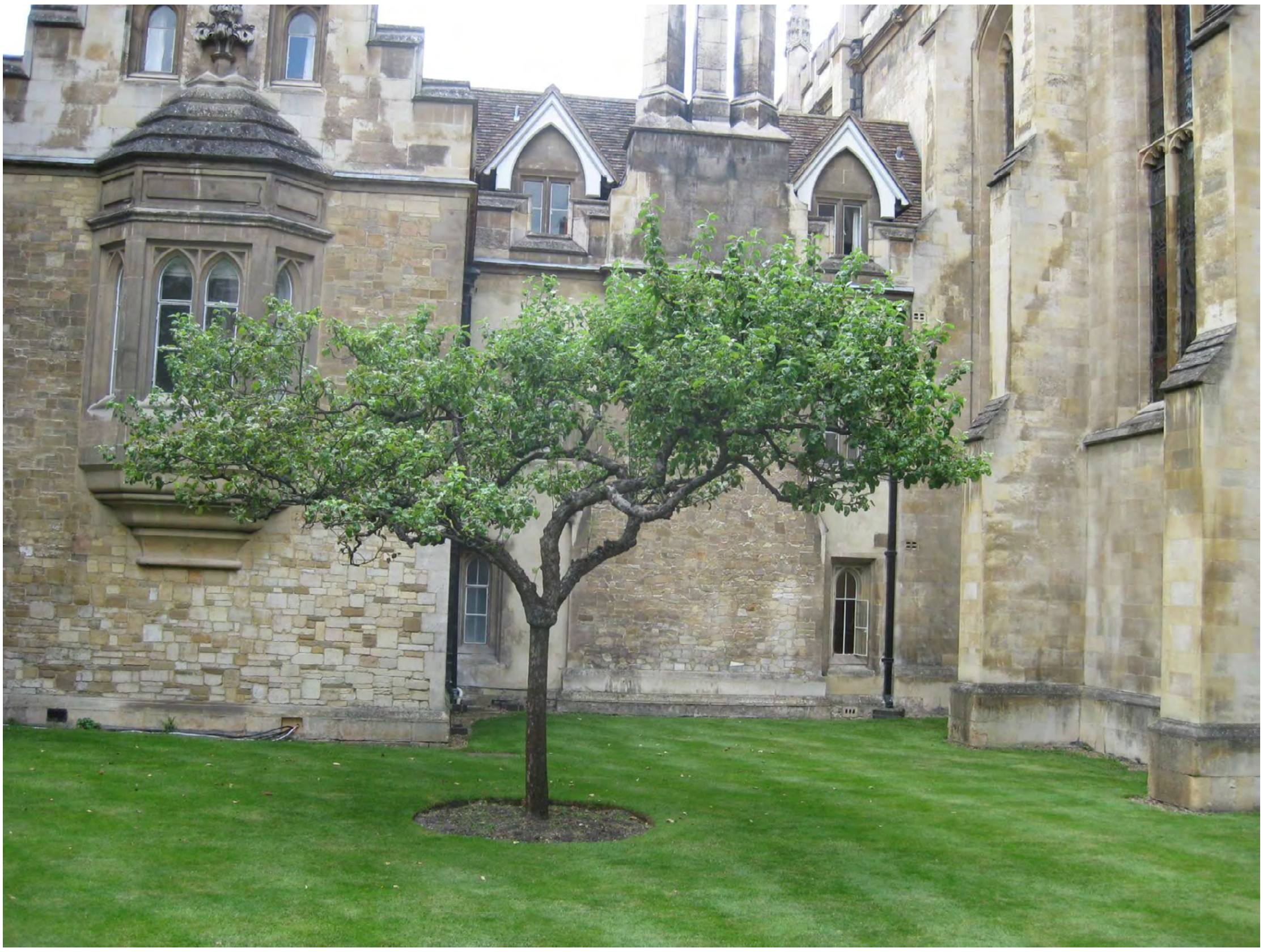

Newton's apple tree in Trinity college. A myth is sometimes circulated that this was the tree from which the apple dropped onto Isaac Newton. In fact, he was not in Cambridge during his *Annus mirabilis*.

# *De methodis serierum et fluxionum* Newton’s series

## *Affected equations*

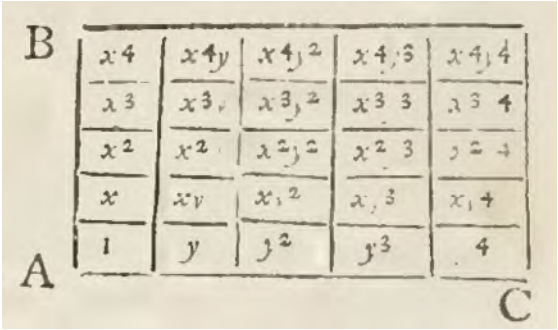

©

Newton can now solve “affected equations”, whose coefficients are *functions* of $x$. His example is

$$y^6 - 5xy^5 + (x^3/a)y^4 - 7a^2x^2y^2 + 6a^3x^3 + b^2x^4 = 0.$$

In this equation, $a, b$ are some parameters. Note that Newton takes great care to write *homogeneous* equations. For simplicity, I will be less careful and choose $a = b = 1$. Newton looks first for an *approximate solution* of the form $y = ux^\alpha$ where $u$ is some unknown nonzero constant and $\alpha$ is some unknown rational number. Substituting, he finds

Do not forget that Newton was *also* a physicist.

$$u^6x^{6\alpha} - 5u^5x^{1+5\alpha} + u^4x^{3+4\alpha} - 7u^2x^{2+2\alpha} + 6x^3 + x^4 = 0.$$

This is an expression involving “monomials" in rational powers of $x$. The exponents are $6\alpha, 1+5\alpha, 3+4\alpha, 2+2\alpha, 3, 4$. If we study the situation in the neighborhood of 0, the largest term corresponds to the smallest of these exponents. For a generic choice of $\alpha$, the six exponents are different. In this is the case and if we wish to express the fact that the dominant term vanishes, that forces $u = 0$: this is certainly not what we want to do.

Therefore, we have to choose $\alpha$ such that at least two of the six exponents are equal and moreover such that they are the smallest. Newton expresses this condition using his famous *polygon*. He draws a kind of checker board subdivided into squares (that he calls *parallelograms*). For each nonzero monomial $a_{ij}x^jy^i$ $(i, j \geq 0)$ in the original equation, he marks a star in the box $(i, j)$. In his example there are six stars.

Choosing $\alpha$ and comparing the exponents $j + i\alpha$ can be interpreted in a geometric way, which is clearly explained by Newton. Place a ruler on the checker board and move it until it touches the marked stars.

> Then, when any Equation is proposed, mark such of the Parallelograms as correspond to all its Terms, and let a Ruler be apply'd to two, or perhaps more, of the Parallelograms so mark'd, of which let one be the lowest in the left-hand Column at AB, the other touching the Ruler towards the right-hand; and let all the rest, not touching the Ruler, lie above it. Then select: those Terms of the Equation which are represented by the Parallelograms that touch the Ruler.

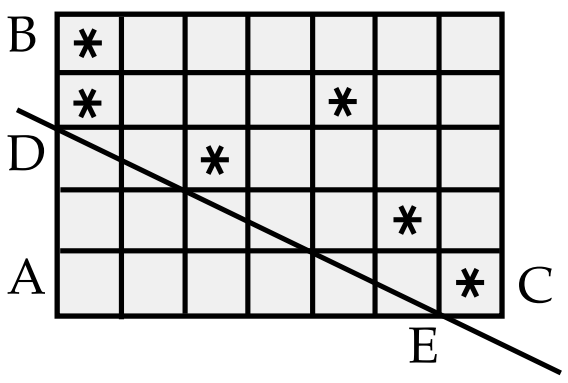

For some reason, Newton marks the monomials $x^j$ on the vertical axis and $y^i$ on the horizontal. ©

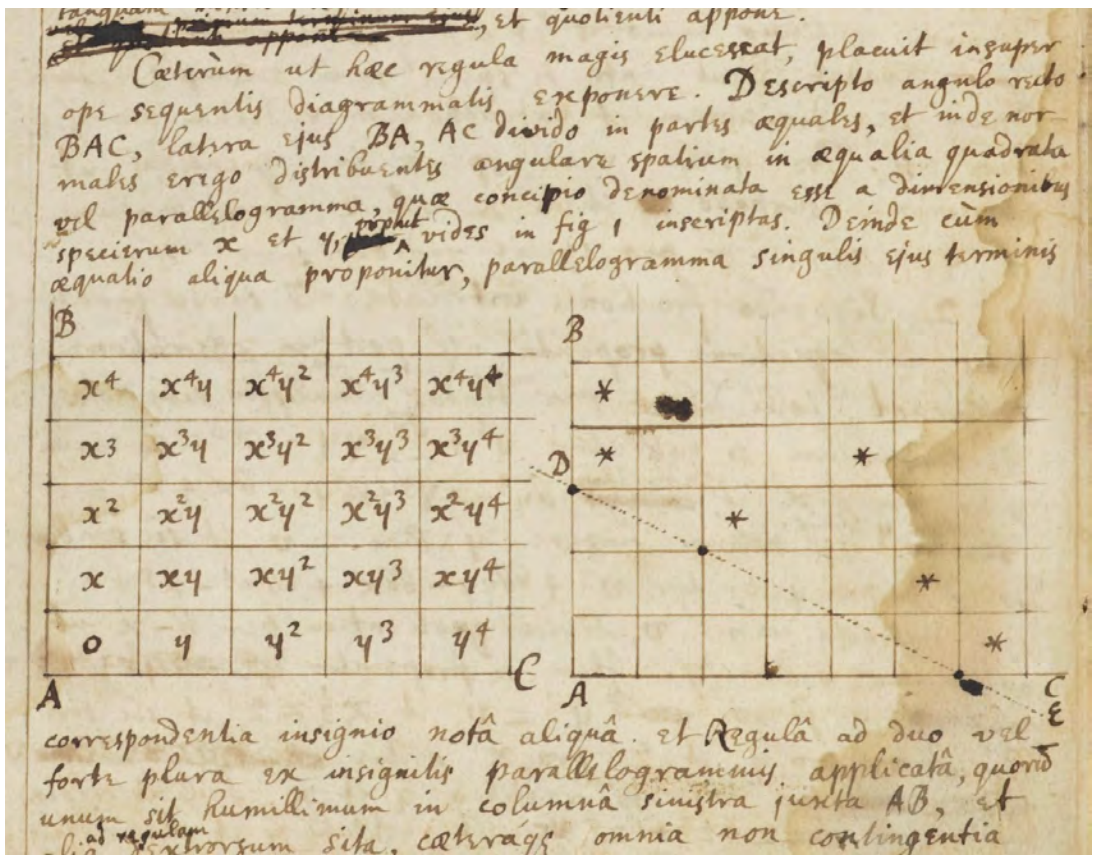

So, in his case, the coefficient $\alpha$ is chosen to be equal to $1/2$ (the slope of the line $DE$) and the three *dominant monomials* $x^3, x^2y^2$ and $y^6$ are chosen. Indeed, for $\alpha = 1/2$, the equation becomes, ordering in increasing powers of $x$:

$$(u^6 - 7u^2 + 6)x^3 - 5u^5x^{7/2} + x^4 + u^4x^5 = 0.$$

We are therefore led to choose $u$ as a solution of the equation

$$u^6 - 7u^2 + 6 = 0$$

which contains three monomials since the ruler touches three stars. There are six solutions

$$u = \pm 1 \quad ; \quad \pm\sqrt{2} \qquad ; \qquad \pm\sqrt{-3}.$$

Newton seems to ignore the last two imaginary solutions. He may be only interested by the real solutions but even if this is the case, this is a mistake, as will be seen later in this section.

He then chooses the first solution. He can write $y = \sqrt{x} + p$, as in his method with *common numbers*. Then, it suffices to "continue the process at pleasure".

However, with no explanation, he abandons suddenly his first example and switches to other numerical examples for which he "exhibits the praxis of his resolution".

Let me show how to continue Newton's first example. For simplicity, I slightly change his presentation.

Instead of improving the first approximate solution $y \simeq \sqrt{x}$ by adding an unknown $p$, let us set

$$x = x_1^2 \quad ; \quad y = x_1(1 + y_1).$$

We substitute these values in the original equation and simplify by $x_1^6$.

$$\begin{gathered} -5x_1 + x_1^2 + x_1^4 - 8y_1 - 25x_1y_1 + 4x_1^4y_1 \\ +8y_1^2 - 50x_1y_1^2 + 6x_1^4y_1^2 + 20y_1^3 \\ -50x_1y_1^3 + 4x_1^4y_1^3 + 15y_1^4 - 25x_1y_1^4 \\ +x_1^4y_1^4 + 6y_1^5 - 5x_1y_1^5 + y_1^6 = 0. \end{gathered}$$

In this new equation, the coefficients of $x_1$ and $y_1$ are not zero. So Newton's ruler passes now through $(0,1)$ and $(1,0)$. This is another way of saying that the new equation is not singular at the origin. So the dominant terms are linear

$$-5x_1 - 8y_1$$

which yields

$$y_1 \simeq -\frac{5}{8}x_1.$$

Continuing the process we set:

$$x_1 = x_2 \quad ; \quad y_1 = -\frac{5}{8}x_2(1 + y_2)$$

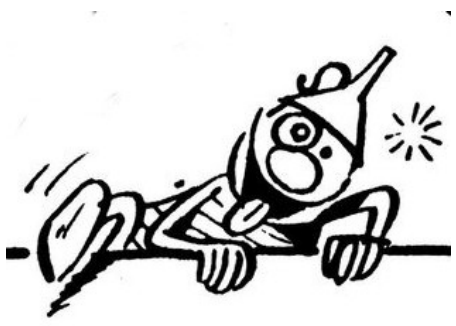

On October 24, 1676, Newton sent a letter to Leibniz, "describing" his contribution to *calculus*. At the end of the letter, he writes "[...] inverse problems of tangents are within our power, and others more difficult than those, and to solve them I have used a twofold method of which one part is neater, the other more general. At present, I thought fit to register them both in transposed letters...", and then he conceals his method in an anagram:

*5accd10effh11i4l3m9n6oqq*

*r8s11t9y3x : 11ab3cdd10eg1*

*0ill4m7n6o3p3q6r5s11t8vx, 3a*

*c4egh5i4l4m5n8oq4r3s6t4v,*

*aaddcecceiijmmnnooprrrss*

*ssttuu.*

Poor Leibniz! He must have struggled to find the meaning of the anagram.

Poor reader! Even if I give the solution, (s)he will have to translate from latin to his own language, and (s)he will then understand that the content is not so clear!

*"Una methodus consistit in extractione fluentis quantitatis ex aequatione simul involvente fluxionem ejus: altera tantum in assumptione seriei pro quantitate qualibet incognita ex qua caetera commode derivari possunt, et in collatione terminorum homologorum aequationis resultantis, as eruendos terminos assumptae seriei".*

and so on, we would find an expansion of $y$ as the product of $\sqrt{x}$ and a series in integral powers of $x$.

With the help of *Mathematica*, my computer finds

$$\begin{aligned} y(x) = \ & x^{1/2} - 5\cdot 2^{-3}x + 79\cdot 2^{-5}x^{3/2} - 14185\cdot 2^{-10}x^2 \\ & + 3118083\cdot 2^{-15}x^{5/2} - 189696965\cdot 2^{-18}x^3 \\ & + 24625187405\cdot 2^{-22}x^{7/2} - 1670815928565\cdot 2^{-25}x^4 + \cdots \end{aligned}$$

For the other solution $u = \sqrt{2}$, the same computer claims that

$$\begin{aligned} y(x) = \ & \sqrt{2}x^{1/2} + 2x - 13\sqrt{2}\cdot 2^{-2}5^{-4}x^{3/2} + 3825^{-2}x^2 \\ & - 267229\sqrt{2}\cdot 2^{-5}5^{-3}x^{5/2} + 903813\cdot 2^{-1}5^{-4}x^3 \\ & - 1661176381\sqrt{2}\cdot 2^{-7}5^{-5}x^{7/2} + 777992628\cdot 5^{-6}x^4 + \cdots \end{aligned}$$

A final comment on the motivation of Newton. Since he "proved" that any "function" $y(x)$ defined by some implicit relation $P(x,y) = 0$ can be expanded as power series of $x$ (at the cost of using rational exponents) and since he, of course, knows very well the derivative and primitive of any power $x^\alpha$, he can use his technique to compute derivatives and primitives of any series. In other words, he is able to compute the derivative and the primitive of "any" function. The rest of *De methodis serierum et fluxionum* is devoted to many applications of this method.

## *A mistake of Newton?*

It is amazing to realize that Newton missed a root of the equation

$$u^6 - 7u^2 + 6 = 0.$$

One might believe that he thought that the imaginary roots $\pm\sqrt{-3}$ would lead to imaginary solutions for $y(x)$. But this is not so and I believe that this is indeed a mistake.

Discovering in 2016 a mistake in an important paper written by Newton around 1669 is an interesting experience. Looking at the original manuscript, we see that Newton had to fix a blunder and to glue a piece of paper above the original page. I suspect that the library of Trinity College would not agree to peel off the precious manuscript to see what is beneath. One should use X-rays.

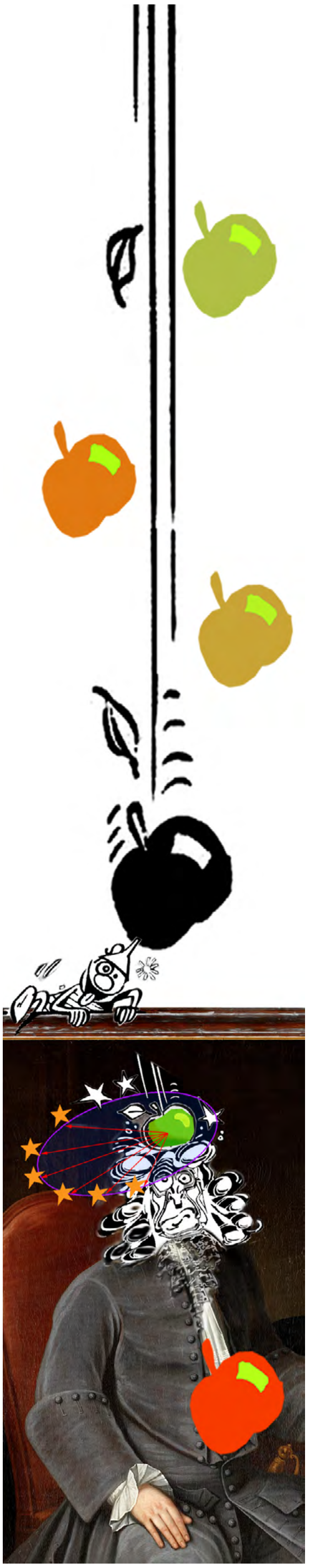

Newton at the precise moment of the mistake? (Homage to Gotlib). ©

A mistake? ©

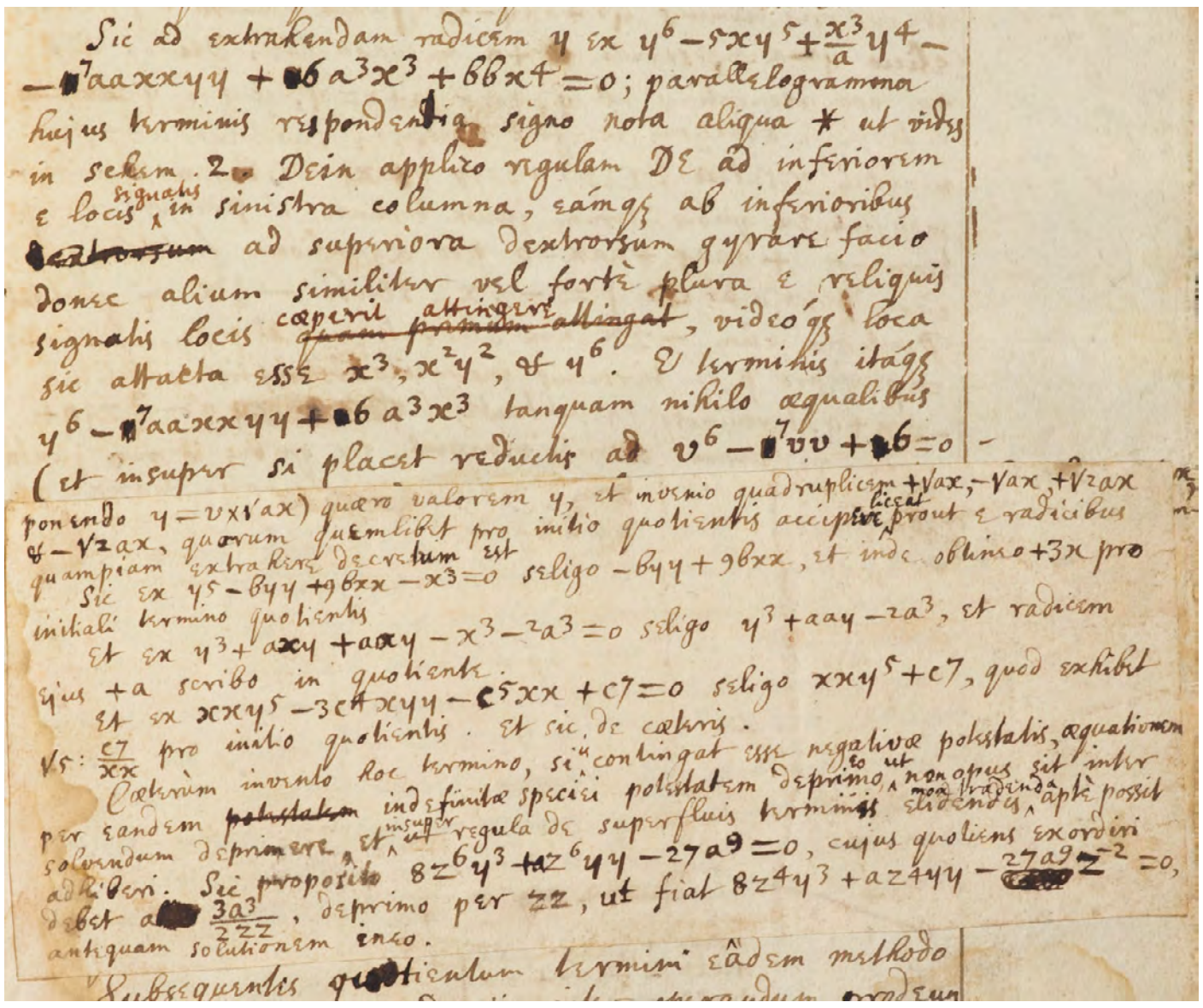

Actually, there could be another interpretation. In his commentary of the *Epistola posterior*, Turnbull[27] (note 68, page 159) mentions another "error": according to him, the "roof" of the square root sign is not long enough and Newton wrote mistakenly $\sqrt{2}x$ instead of $\sqrt{2x}$. Then he comments that "Newton rejects the imaginary cases given by $v^2 + 3 = 0$". It is indeed possible that Newton made a mistake with $\sqrt{2}x$ which led him to think of $\sqrt{-3}x$ as imaginary, and then to reject it. If he had written $\sqrt{-3x}$ he would have seen that this solution is not imaginary at all when $x < 0$. We'll never know.

[27] I. Newton. *The correspondence of Isaac Newton, Vol. II: 1676–1687*. Published for the Royal Society. Cambridge University Press, New York, 1960.

Indeed, $\pm\sqrt{-3}$ is imaginary but the approximate function $y \simeq \pm\sqrt{-3}\sqrt{x}$ is real if $x$ is a *negative real number* so that it should not have been discarded.

For the real root $u = 1$, we have set

$$x = x_1^2 \quad ; \quad y = x_1(1 + y_1).$$

For the imaginary root $\sqrt{-3}$, we set

$$x = -x_1^2 \quad ; \quad y = \sqrt{3}x_1(1 + y_1)$$

and we proceed as before. We finally get a third *real* solution.

$$
\begin{aligned}
y(x) = \; & -3^{1/2}(-x)^{1/2} - 9\cdot 2^{-3}(-x) - 721\cdot 2^{-6}5^{-1}3^{-1/2}(-x)^{3/2} - 36543\cdot 2^{-10}5^{-2}(-x)^2 \\
& -27986569\cdot 2^{-15}3^{-3/2}5^{-3}(-x)^{5/2} - 96025589\cdot 2^{-18}5^{-4}(-x)^3 \\
& +169264391911\cdot 2^{-22}3^{-5/2}5^{-5}(-x)^{7/2} + 1398151100829\cdot 2^{-25}5^{-6}(-x)^4 + \cdots
\end{aligned}
$$

One may ask why we found three solutions and not six since the equation $u^6 - 7u^2 + 6 = 0$ has indeed six solutions. This is simply because opposite roots give rise to the same solution. Do not forget that Newton considers $\sqrt{x}$ as a 2-valued function, so that for him $\sqrt{x}$ and $-\sqrt{x}$ are "the same". I agree that writing $\sqrt{x} = -\sqrt{x}$ might lead to contradictions, but not under the pen of Newton. We are wise to teach our students that $\sqrt{x}$ is the positive root for $x$ real and positive, and to choose some principal determination for $x \in \mathbb{C} \setminus \mathbb{R}_-$. In modern terminology, the two parameterized curves $(t^2, t)$ and $(t^2, -t)$ are the same curves, with different parameterizations.

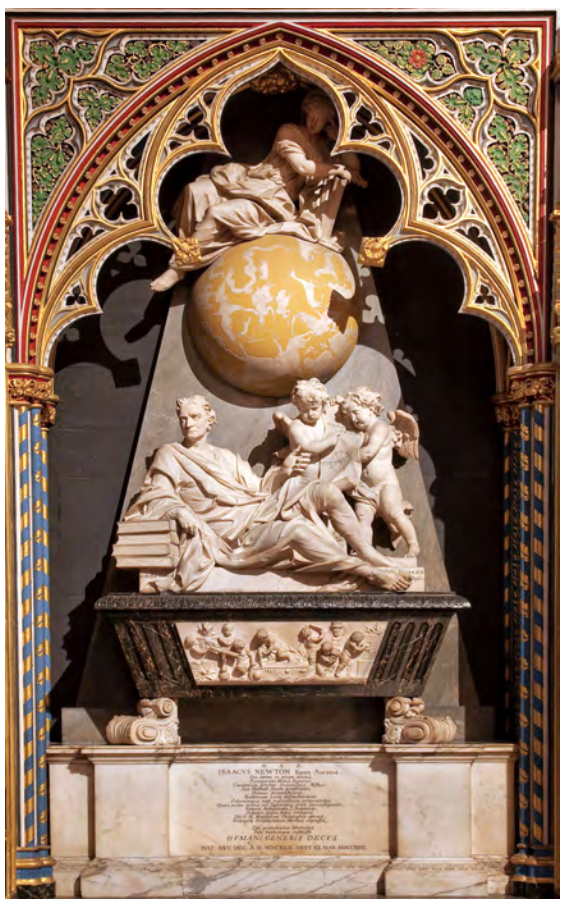

Newton's monument in Westminster Abbey.

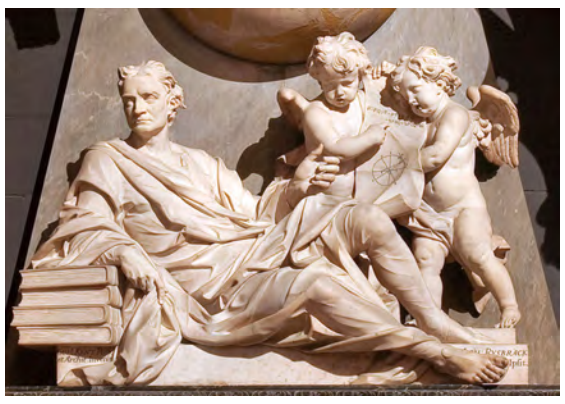

Right hand on his four main books.

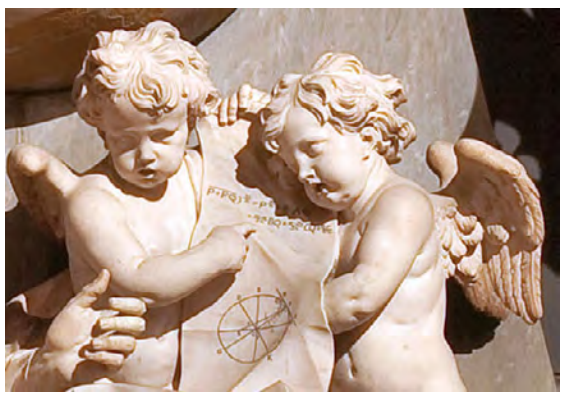

Left hand pointing to two angels showing... the binomial series.

## *What Newton did not prove*

The definition of "convergence" was not at Newton's disposal. However, his numerical computations suggest that his series are indeed convergent and he even uses the terminology *convergent*. To be honest, one could say that he only shows that his series give *asymptotic expansions*. In practice, a series

$$a_1x^{\alpha_1} + a_2x^{\alpha_2} + \cdots$$

(where $\alpha_1 < \alpha_2 < \ldots$ are rational exponents) is *asymptotic* to a function $f(x)$ if for every $n \geq 1$:

$$f(x) - \sum_{i=1}^{k} a_i x^{\alpha_i} = o(x^{\alpha_k}).$$

This does not imply that $f_k(x) = \sum_{i=1}^{k} a_i x^{\alpha_i}$ converges to $f$, but is frequently as useful, and sometimes even more useful, than a usual convergence.

Another aspect that he does not discuss concerns the nature of the rational exponents that appear in his series. At each step, a new rational number appears and it is not clear that this sequence of exponents converges to infinity. Even less clear is the fact that all the denominators are *bounded*. However, Newton does observe that his method is not restricted to polynomial equations $P(x, y) = 0$ but works perfectly for "aequationes numero terminorum infinitas" of the form $\sum_{i,j \geq 0} a_{ij} x^i y^j = 0$ (with $a_{00} = 0$), involving what are called today *analytic functions*.

To conclude this chapter, let us look at Newton's original curve.

$$P(x, y) = y^6 - 5xy^5 + x^3y^4 - 7x^2y^2 + 6x^3 + x^4 = 0.$$

If I ask my computer to plot this curve in a $[-50, +20] \times [-50, +50]$ box, I get the first plot in the margin. This may look surprising since we only see two branches in the neighborhood of the origin. Let us zoom and look in a smaller box $[-1, 1] \times [-2, 2]$ (second plot), we can guess another branch. Zooming more in $[-.1, .1] \times [-.4, .4]$ (third plot), this is easier to see. The local situation is completely clear in $[-.01, .01] \times [-.2, .2]$ (fourth plot). Finally, in $[-.001, .001] \times [-.05, .05]$, we do see three branches asymptotic to $\pm\sqrt{x}, \pm\sqrt{2x}, \pm\sqrt{-3x}$ as predicted by Newton (except that he forgot the third, for $x < 0$).

The polynomial $P(x, y)$ is prime: it does not split non-trivially as a product of two polynomials. However, as a convergent power series in $x, y$, in the neighborhood of the origin, it does split as a product of three factors, corresponding to the three branches.

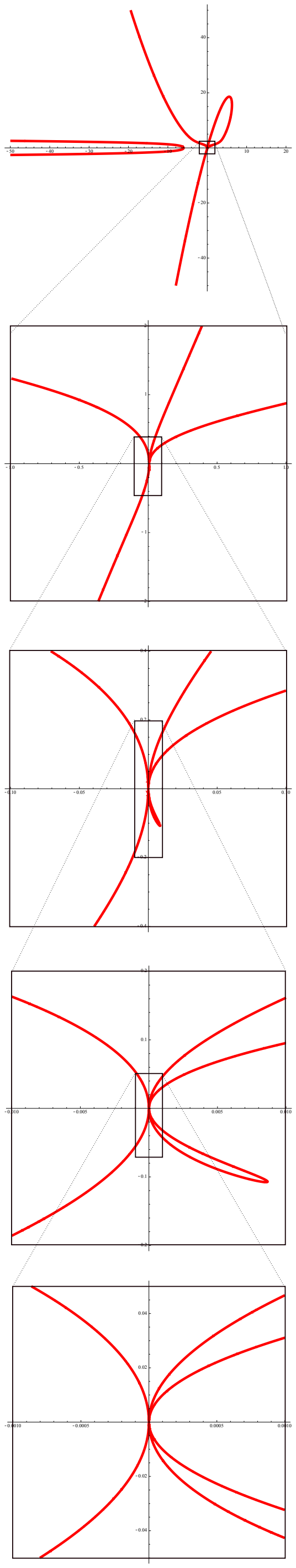

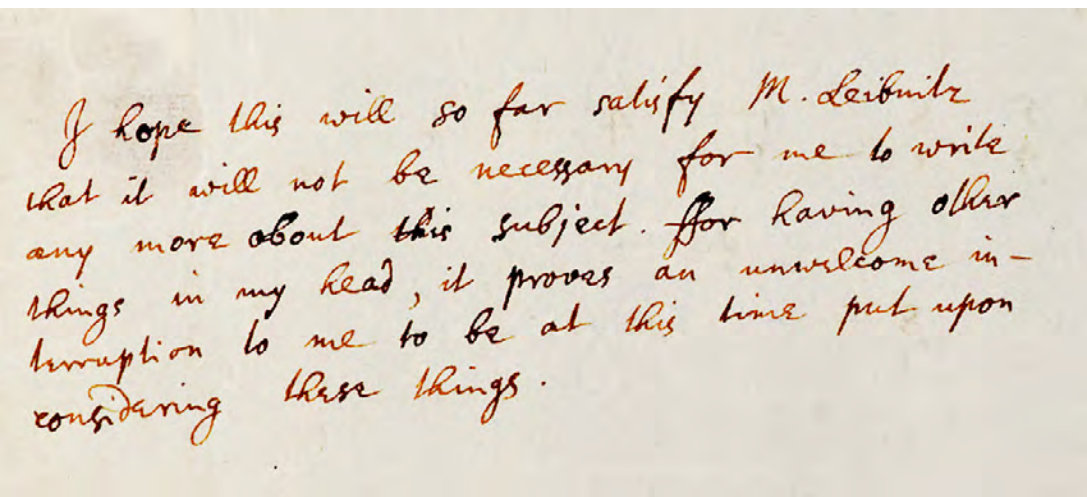
I hope this will so far satisfy M. Leibnitz
that it will not be necessary for me to write
any more about this subject. For having other
things in my head, it proves an unwelcome in-
terruption to me to be at this time put upon
considering these things.

When Newton asked Oldenbourg to forward his *Epistola Posterior* to Leibniz, he added this P.S. Yes, he had indeed other things in his head.

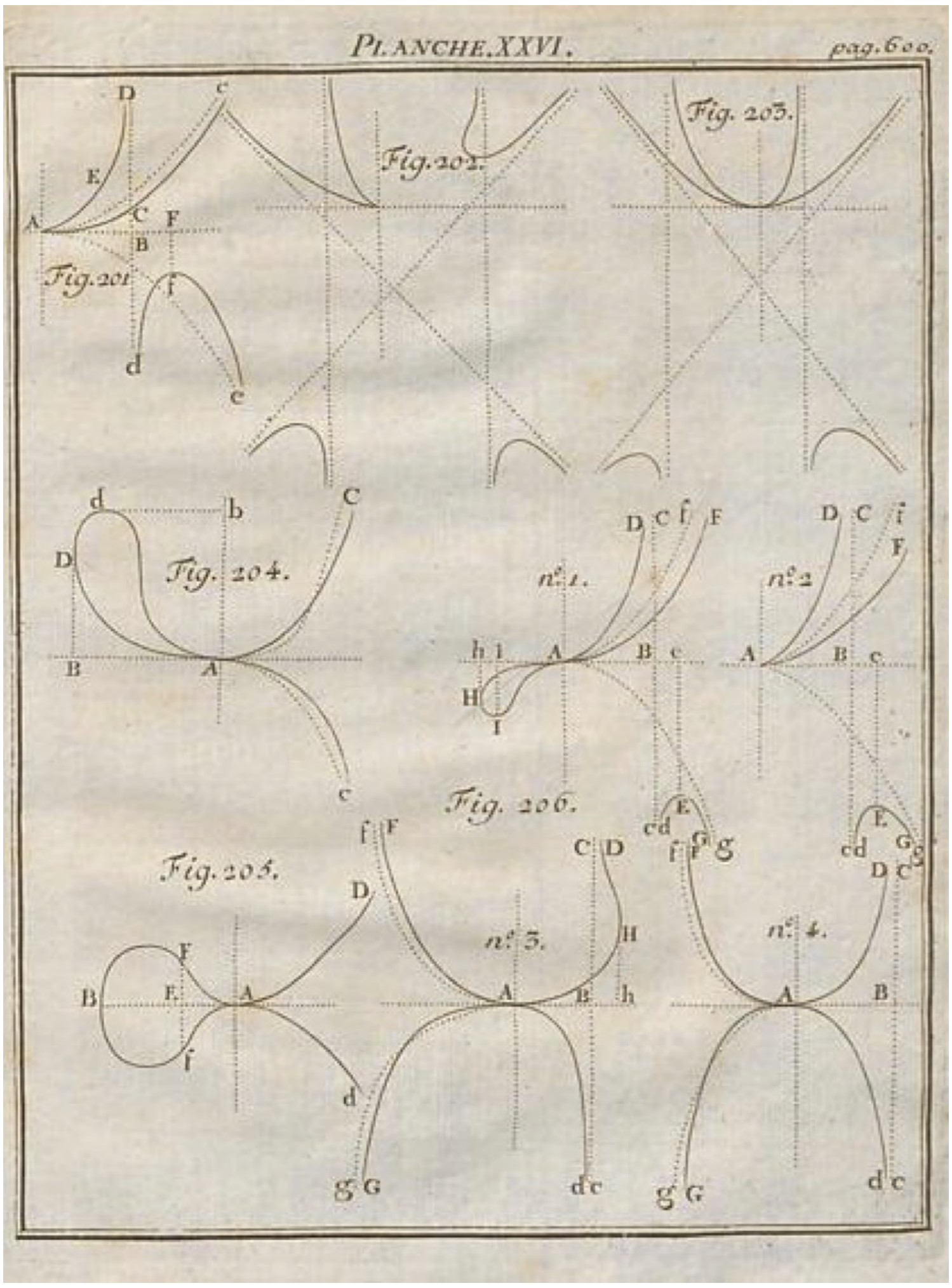

A plate from Cramer's book on curves. ©

# *Some formal algebra*

The algebra in this chapter will be "formal" since we will consider formal series.

## *Finding one solution*

Let me repeat Newton's arguments, expressing them in a more modern algebraic terminology. Attributing all of this to Newton requires infinite imagination and extrapolation. I will emphasize however an important contribution of Cramer. Usually a good part of what follows is attributed to Puiseux, but this would require at least as much of extrapolation. I'll describe Puiseux's contribution in due course.

There are excellent books on this topic. I recommend in particular Walker[28], Brieskorn and Knörrer[29], Wall[30], and Casas-Alvero[31].

Let $K$ denote some *algebraically closed field of characteristic zero*. The main example that I have in mind is of course the field $\mathbb{C}$ of complex numbers.

Some notation:

– $K[x]$ is the ring of *polynomials* in $x$ with coefficients in $K$.

– $K(x)$ is the field of *rational functions* in $x$ with coefficients in $K$: this is the quotient field of $K[x]$.

– $K[[x]]$ is the ring of *formal series* in $x$: expressions of the form $\sum_{i=0}^{\infty} a_i x^i$ where the $a_i$'s are in $K$, without any reference to convergence matters.

[28] R. J. Walker. *Algebraic curves*. Springer-Verlag, New York-Heidelberg, 1978. Reprint of the 1950 edition.

[29] E. Brieskorn and H. Knörrer. *Plane algebraic curves*. Modern Birkhäuser Classics. Birkhäuser/Springer Basel AG, Basel, 1986.

[30] C. T. C. Wall. *Singular points of plane curves*, volume 63 of *London Mathematical Society Student Texts*. Cambridge University Press, Cambridge, 2004.

[31] E. Casas-Alvero. *Singularities of plane curves*. London Mathematical Society Lecture Note Series. Cambridge University Press, 1 edition, 2000.

– In a similar way, in two variables, we define $K[x,y]$ (polynomials), $K(x,y)$ (rational functions) and $K[[x,y]]$ (formal series).

We can now state Newton's theorem in a precise form.

**Theorem** (Newton-Cramer). *Let $F(x,y)$ be a formal series in $K[[x,y]]$ vanishing at the origin and not divisible by $x$. Then there is an integer $m \geq 1$ and a formal power series $f(t) \in K[[t]]$ vanishing at 0 such that $F(t^m, f(t))$ vanishes identically. In other words, the equation $F(x,y) = 0$ has at least one "solution" of the form $y = f(x^{1/m})$.*

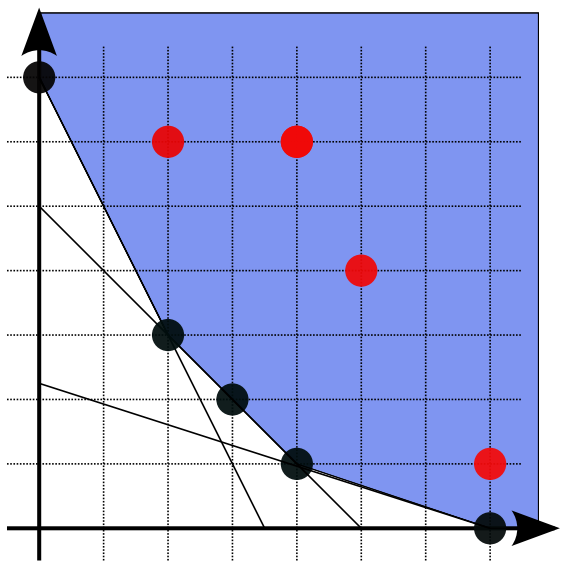

A Newton polygon. Terms in the series with nonzero $a_{ij}$ are represented by dots. There are three supporting lines. Note that I do not follow Newton's strange idea of writing the $i$'s on the vertical axis. For the benefit of the reader, I follow the tradition and use $x$ and $y$ for the horizontal and vertical axes.

We have already seen the general structure of the proof.

Let me write $F_0(x_0, y_0)$ instead of $F(x,y)$ since I will describe some iterative construction involving some $x_k, y_k$'s. So, let $F_0(x_0, y_0) = \sum_{i,j} a_{ij} x_0^i y_0^j$ (with $a_{00} = 0$) be a formal series. For each $(i,j)$ such that $a_{ij} \neq 0$, consider the quarter plane $x \geq i; y \geq j$. The *Newton polygon* is the convex hull of the union of these quarter planes. The picture in the margin shows the polygon for

$$F_0 = y_0^7 - x_0^2 y_0^3 + x_0^2 y_0^6 + x_0^3 y_0^2 \\ + x_0^4 y_0 + x_0^4 y_0^6 + x_0^5 y_0^4 + x_0^7 + x_0^7 y_0.$$

We can always assume that $F_0$ is not divisible by $y_0$ or by $x_0$ since we could divide by some monomial $y_0^j$ or $x_0^i$ without changing the problem. In other words, each of the axis intersects the Newton polygon. The boundary of this polygon, away from the axis, consists of a finite number of segments, included in *supporting lines* whose equations have the form $\alpha i + \beta j = \gamma$ where $\alpha, \beta$ are positive integers that we can assume relatively prime. Choose one of these lines $\alpha_0 i + \beta_0 j = \gamma_0$ and select the finite number of coefficients $a_{ij}$ such that $(i,j)$ lies on this line. This defines a *"dominant polynomial"*

$$F_{dom}(x_0, y_0) = \sum_{\alpha_0 i + \beta_0 j = \gamma_0} a_{ij} x_0^i y_0^j.$$

Example: Choose the supporting line

$$2i + j = 7,$$

so that $(\alpha, \beta, \gamma) = (2, 1, 7)$. $F_{dom}(x_0, y_0) = y_0^7 - x_0^2 y_0^3$ and the polynomial $p$ is

$$p(u) = u^7 - u^3.$$

Let us look for an approximate solution of our equation $F_0(x_0, y_0) = 0$ of the form

$$x_0 = t^{\alpha_0} \quad ; \quad y_0 = u t^{\beta_0}.$$

Choose $u_0 = 1$ so that the approximate solution is

$$x_0 = t^2; y_0 = t.$$

By "approximate", I mean that it solves the dominant part of our equation

$$F_{dom}(t^{\alpha_0}, u t^{\beta_0}) = 0.$$

We get a polynomial equation in $u$:

$$p(u) = \sum_{\alpha_0 i + \beta_0 j = \gamma_0} a_{ij} u^j = 0.$$

Since $K$ is algebraically closed, there is at least one nonzero solution $u_0$. We then come back to the original equation $F_0(x_0, y_0) = 0$ and substitute $x_1^{\alpha_0}$ for $x_0$ and $u_0 x_1^{\beta_0}(1 + y_1)$ for $y_0$. This produces a new formal series in $(x_1, y_1)$ which, by construction, is divisible by $x_1^{\gamma_0}$. Dividing by $x_1^{\gamma_0}$, we get another equivalent equation $F_1(x_1, y_1) = 0$

...and the process can be continued "at pleasure" (in Newton's words), producing a sequence of equations $F_k(x_k, y_k)$ $(k \geq 1)$ and of integers $\alpha_k, \beta_k, \gamma_k$.

One important property is missing and was *not* discussed by Newton. We have to show that after a finite number of steps, the coefficients $\alpha_k$ are always equal to 1, which means that the slopes of all supporting lines are inverses of integers and not only rational numbers. This is important since each step implies the introduction of a root $x_{k+1} = x_k^{1/\alpha_k}$ and we would get into trouble if we had to do that an infinite number of times.

This is analyzed in full detail in Chapter VII of Cramer's excellent book *Introduction a l'analyse des lignes courbes algébriques*[32], published in 1750. The author gives a full credit to Newton but explains that:

> La vraye Méthode des Séries est fondée sur le Parallélogramme de Mr. Newton, invention excellente, mais dont l'Auteur n'a pas donné la Démonstration, dont il semble même n'avoir pas senti tout le prix.

**Definition.** If $F$ is a formal power series in $K[[x, y]]$ not divisible by $x$, the *multiplicity*, denoted $mult(F)$, is the valuation of $F(0, y)$ as a series in $y$. This is also the smallest height of a point of the Newton polygon of $F$ on the vertical axis.

Note that by convexity any supporting line intersects the $j$-axis *below* $mult(F)$. In particular the degree of the polynomial $p(u)$ is at most $mult(F)$.

**Lemma.** $mult(F_1) \leq mult(F)$.

First step. In $F_0$ substitute

$$x_0 \to x_1^2; y_0 \to x_1(1 + y_1).$$

and divide by $x_1^7$. We get $F_1(x_1, y_1) = x_1 + 4y_1 + 2x_1^2 + 2x_1 y_1 + 18y_1^2 + 6x_1^2 y_1 + x_1 y_1^2 + 34y_1^3 + 10x_1^2 y_1^2 + 35y_1^4 + {}+10x_1^2 y_1^3 + 21y_1^5 + {}+5x_1^2 y_1^4 + 7y_1^6 + 3x_1^7 + x_1^2 y_1^5 + y_1^7 + x_1^8 + 10x_1^7 y_1 + x_1^8 y_1 + 21x_1^7 y_1^2 + 24x_1^7 y_1^3 + 16x_1^7 y_1^4 + 6x_1^7 y_1^5 + {}+x_1^7 y_1^6$, which has a nontrivial linear term in $y_1$ so that $y_1$ can be expanded as a power series in $x_1$.

[32] G. Cramer. *Introduction à l'Analyse des lignes courbes algébriques*. Frères Cramer et Cl. Philibert, 1750.

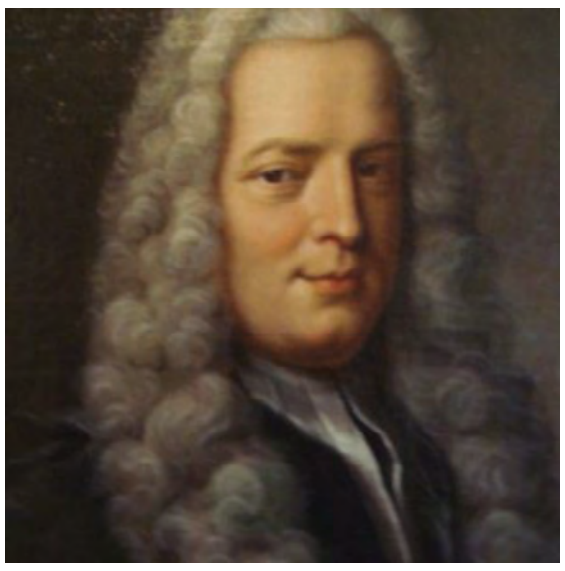

Gabriel Cramer (1704-1752). His book on curves contains, besides a serious analysis of Newton's series, a theory of linear equations in $n$ unknowns (the famous Cramer's rule) and the elements of elimination theory. I like the title of appendix 1: "De l'évanouissement des inconnues" which looks more enticing than "elimination".

By definition

$$F_1(x_1, y_1) = x_1^{-\gamma} \sum_{i,j} a_{ij} x_1^{\alpha_0 i + \beta_0 j} u_0^j (1 + y_1)^j.$$

In order to get $mult(F_1)$, we let $x_1 = 0$ and look at the valuation of $p(u_0(1+y_1))$ as a polynomial in $y_1$.

$$mult(F_1) \leq degree(p) \leq mult(F).$$

$\boxdot$

So under Newton's algorithm, the sequence of multiplicities $mult(F_k)$ is non-increasing. This inequality is strict unless $F_0$ has a very special form that we now analyze.

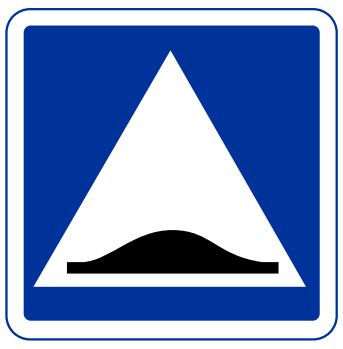

One has equality if and only if $p(u_0(1+y_1))$ contains only one monomial of degree $mult(F)$. This implies in particular that the degree of $p$ is equal to $mult(F)$.

Said differently, the root $y_1 = 0$ of $p(u_0(1+y_1)) = 0$ should be multiple of order $mult(F)$. This means that $p$ has the form

$$p(u) = C(u - u_0)^{mult(F)}.$$

This polynomial has nonzero coefficients in each degree from 0 to $mult(F)$. Hence, the segment of the boundary of the Newton polygon that we have chosen contains dots for each value of $j$ from $j = 0$ to $j = mult(F)$. This implies that the Newton polygon has only one side (other from the axes segment) and that $\alpha_0 = 1$.

Let's sum up.

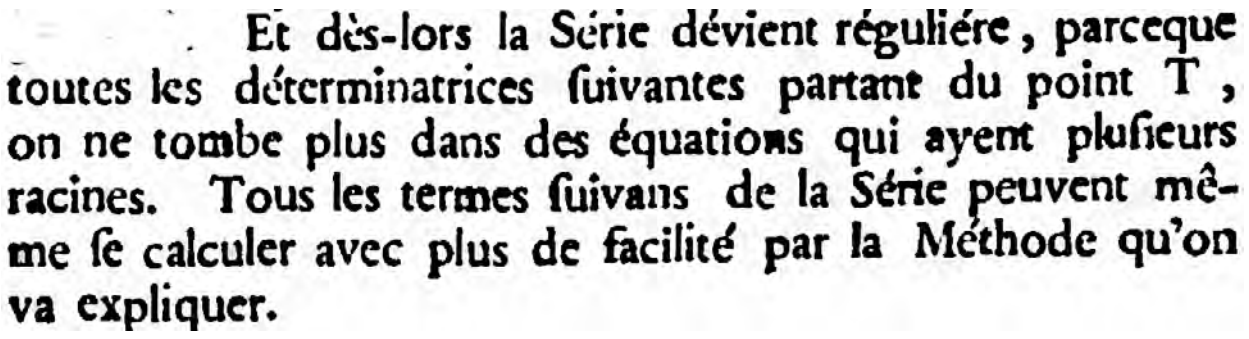

Et dès-lors la Série dévient réguliére, parceque toutes les déterminatrices suivantes partant du point T, on ne tombe plus dans des équations qui ayent plusieurs racines. Tous les termes suivans de la Série peuvent même se calculer avec plus de facilité par la Méthode qu'on va expliquer.

Extract of the proof by Cramer (page 200).

Along the algorithm, the multiplicities $mult(F_k)$ are non-increasing and therefore they have to be constant after some time. At this stage, all Newton polygons have $\alpha_k = 1$ (and moreover have the very special structure that was just described). We have

$$x_0 = x_1^{\alpha_0} = x_2^{\alpha_0 \alpha_1} = \ldots = x_k^{\alpha_0 \alpha_1 \cdots \alpha_{k-1}} = \ldots$$

and

$$y_0 = u_0 x_1^{\beta_1}(1 + y_1) = u_0 x_1^{\beta_1}(1 + u_1 x_2^{\beta_2}(1 + y_2)) = \dots$$

Since the $\alpha_k$'s are equal to 1 for large values of $k$, we can set $m$ to be the product of all the $\alpha_k$'s and call $t$ the value of $x_k$ for large $k$. We have $x_0 = t^m$ and more generally, each $x_k$ is a power of $t$ with some integral positive exponent.

The inductive formula

$$y_k = u_k x_{k+1}^{\beta_{k+1}}(1 + y_{k+1})$$

defines a sequence of polynomials $y_k(t)$ in the variable $t$. This sequence "converges" to a limiting series $f(t) \in K[[t]]$. This means that the valuation of $f(t) - y_k(t)$ goes to infinity when $k$ converges to infinity.

To be complete, we should check that this is indeed a solution to our problem, i.e. that $F(t^m, f(t))$ does vanish identically. I encourage the reader to check it. After all, the algorithm had only one goal: to find a solution.

This is the proof of Newton's theorem: every equation of the form $F(x, y) = 0$ has some solution, in a properly defined sense.

⊡

## *Algebraic closure*

It is time to give a precise definition of series with rational exponents.

Denote by $K[[x]][x^{-1}]$ the field of *formal Laurent series*, that is to say formal expressions of the form $\sum_{i \geq i_0}^{\infty} a_i x^i$ (where $i_0$ might be a negative integer). This is the quotient field of the ring $K[[x]]$.

More generally, if $m$ is a non-zero positive integer, we denote by $K[[x^{1/m}]][x^{-1/m}]$ the field of formal Laurent series in the variable $x^{1/m}$: formal expressions of the form $\sum_{i \geq i_0}^{\infty} a_i x^{i/m}$ (where $i_0$ might be a negative integer). The subfield consisting of series for which $a_i = 0$ whenever $i$ is not divisible by $m$ is canonically isomorphic to $K[[x]][x^{-1}]$ so that we can see $K[[x^{1/m}]][x^{-1/m}]$ as a field extension of $K[[x]][x^{-1}]$. The Galois group of this extension is easy to describe: it consists of the $m$-th roots of unity. The

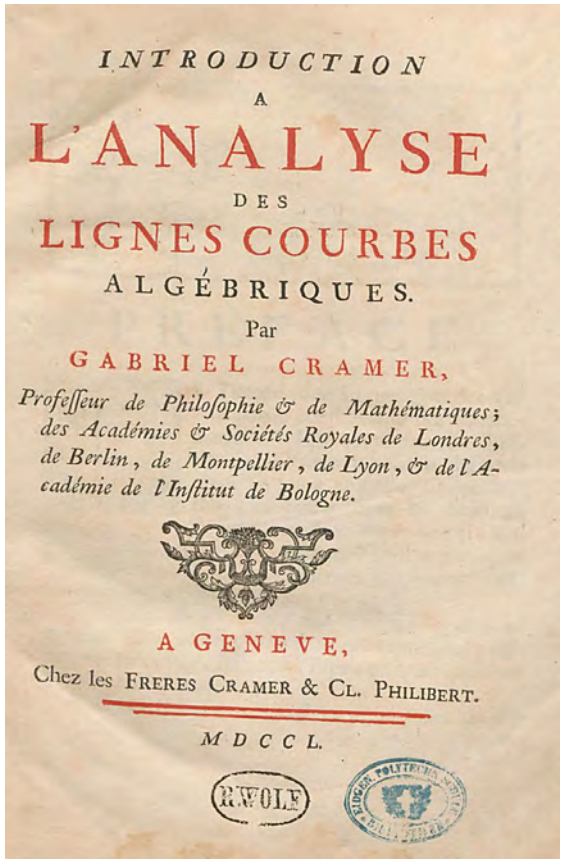
INTRODUCTION
A
L'ANALYSE
DES
LIGNES COURBES
ALGÉBRIQUES.
Par
GABRIEL CRAMER,
Professeur de Philosophie & de Mathématiques; des Académies & Sociétés Royales de Londres, de Berlin, de Montpellier, de Lyon, & de l'Académie de l'Institut de Bologne.
A GENEVE,
Chez les FRERES CRAMER & CL. PHILIBERT.
MDCCL.

©

Consider the algebraic closure $\overline{\mathbb{F}}_2$ of the field $\mathbb{F}_2$ with 2 elements. Try the algorithm on the polynomial $F = y^2 + x^2 y + x^2$ with coefficients in $\mathbb{F}_2$ and show that we get the solution

$$y = x\left(1 + x^{\frac{1}{2}} + \cdots + x^{1-2^{-k}} + \cdots\right).$$

This is *not* a Puiseux series since the exponents $1 - 2^{-k}$ do not have a common denominator and don't even tend to infinity. What happened ?

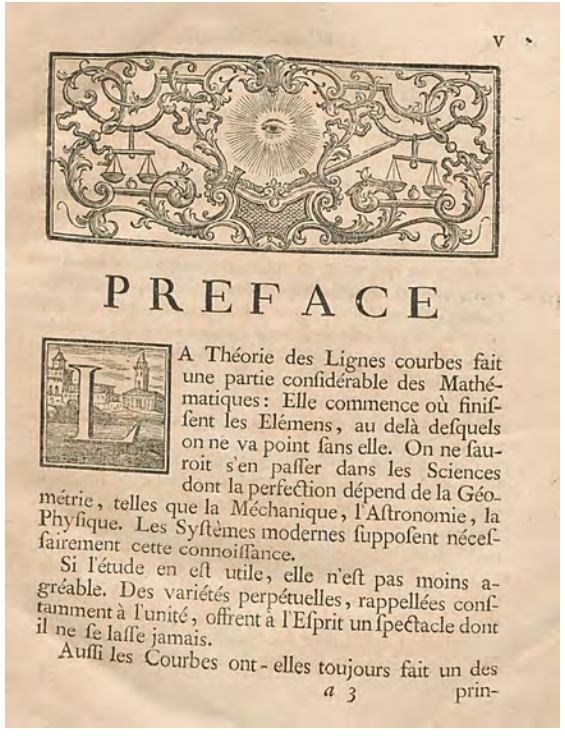
v

PREFACE

LA Théorie des Lignes courbes fait une partie considérable des Mathématiques: Elle commence où finissent les Elémens, au delà desquels on ne va point sans elle. On ne sauroit s'en passer dans les Sciences dont la perfection dépend de la Géométrie, telles que la Méchanique, l'Astronomie, la Physique. Les Systèmes modernes supposent nécessairement cette connoissance.

Si l'étude en est utile, elle n'est pas moins agréable. Des variétés perpétuelles, rappellées constamment à l'unité, offrent à l'Esprit un spectacle dont il ne se lasse jamais.

Aussi les Courbes ont-elles toujours fait un des

a 3 prin-

"Un spectacle dont on ne se lasse jamais".

action of such a root $\omega$ on

$$\sum_{i \geq i_0}^{\infty} a_i x^{i/m}$$

produces

$$\sum_{i \geq i_0}^{\infty} \omega^i a_i x^{i/m}.$$

This is a Galois extension: the elements of $K[\![x^{1/m}]\!][x^{-1/m}]$ which are invariant under the Galois group action are in $K[\![x]\!][x^{-1}]$.

In the same way if $m_1$ divides $m_2$, the field $K[\![x^{1/m_1}]\!][x^{-1/m_1}]$ is a subfield of $K[\![x^{1/m_2}]\!][x^{-1/m_2}]$. The direct limit of all these extensions of $K[\![x]\!][x^{-1}]$ is denoted by $K[\![x^{\star}]\!][x^{\star -1}]$.

*This is the field of Puiseux series*: series with rational exponents, having a common denominator. In down to earth terms, a Puiseux series is a formal expression of the form $\sum_{i \geq i_0}^{\infty} a_i x^{i/m}$ for some non-zero positive integer $m$. Puiseux series with $i_0 \geq 0$ constitute a ring, that we denote by $K[\![x^{\star}]\!]$.

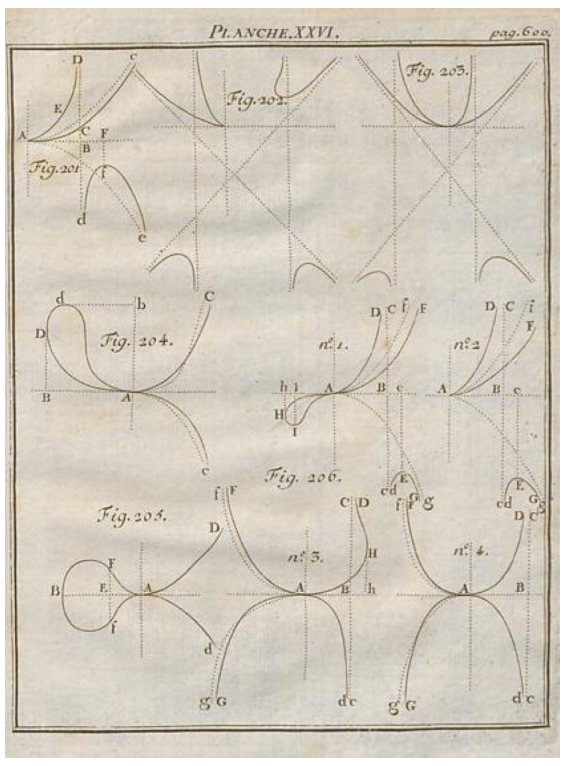

©

**Theorem** (Newton-Cramer). *The field of Puiseux series $K[\![x^{\star}]\!][x^{\star -1}]$ is algebraically closed. This is the algebraic closure of the field of Laurent series $K[\![x]\!][x^{-1}]$.*

This theorem is nothing but a restatement of the main theorem of this chapter.

The fact that $K[\![x^{\star}]\!][x^{\star -1}]$ is an algebraic extension of $K[\![x]\!][x^{-1}]$ is clear. Indeed any Puiseux series lies in some $K[\![x^{1/m}]\!][x^{-1/m}]$ and is therefore algebraic over $K[\![x]\!][x^{-1}]$.

Consider a (non-constant) polynomial equation with coefficients in $K[\![x^{\star}]\!][x^{\star -1}]$ and variable $y$. Defining $\overline{x} = x^{1/m}$ for some highly divisible $m$ and multiplying all coefficients by a high power of $\overline{x}$ we can assume that the coefficients of our polynomial are in $K[\![\overline{x}]\!]$. Our equation is therefore of the form $F(\overline{x}, y) = 0$ where $F$ is a formal power series. We know that such an equation has a solution as a series in $\overline{x}^{1/\overline{m}}$ for some $\overline{m}$, which is in particular in $K[\![x^{\star}]\!][x^{\star -1}]$. So, for any non-constant polynomial with coefficients in $K[\![x^{\star}]\!][x^{\star -1}]$ we found a root in $K[\![x^{\star}]\!][x^{\star -1}]$. $\boxdot$

## *Finding all solutions*

If we think of $F(x,y) = 0$ as an equation where the unknown is a series $y(x) \in K[\![x^\star]\!][x^{\star -1}]$, we can try, as we would do with a usual polynomial equation, to factor $F$ as a product of linear factors in the algebraic closure

$$F = A(x,y)(y - f_1(x))(y - f_2(x))\cdots(y - f_n(x))$$

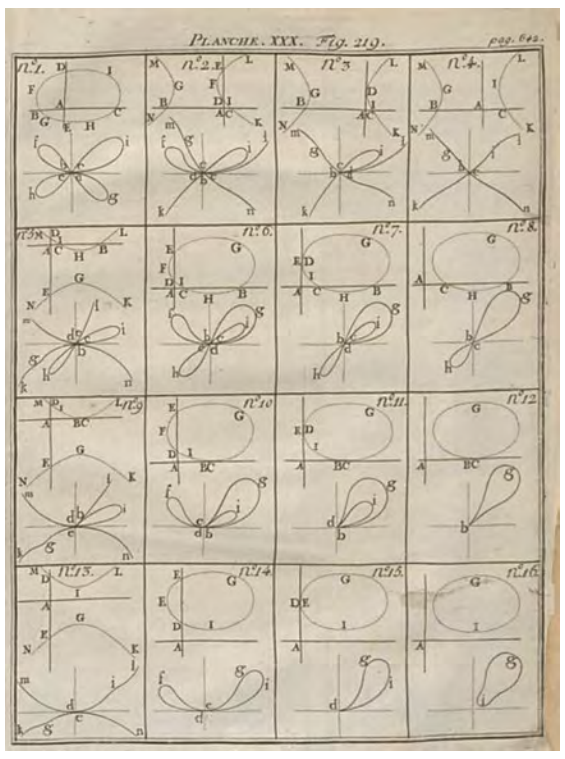

©

where $A(0,0) \neq 0$ and the $n$ solutions $f_i(x)$ are in $K[\![x^\star]\!][x^{\star -1}]$. That would be obvious if $F$ was a polynomial in the $y$ variable, but it is only a formal series. It is not even clear that our equation has a finite number of solutions.

As a matter of fact, Newton was right and our equations are indeed very close to being "standard" polynomial equations, as I explain now.

Let me begin with some elementary observations.

**Lemma.** *Suppose a formal series $y = f(x) \in K[\![x]\!]$ is a solution to the equation $F(x,y) = 0$ where $F \in K[\![x,y]\!]$. Then $F$ is divisible by $y - f(x)$ in $K[\![x,y]\!]$.*

This is obvious if $f(x) = 0$. Now the formal transformation $(x,y) \mapsto (x, y - f(x))$ induces an automorphism of $K[\![x,y]\!]$ sending $y$ to $y - f(x)$. ⊡

**Lemma.** *For $f \in K[\![x^\star]\!]$, define*

$$\overline{f}(x,y) = (y - f_1(x))(y - f_2(x))\cdots(y - f_m(x))$$

*where $f_1, \ldots, f_m$ are the Galois conjugates of $f$. Then $\overline{f}(x,y)$ is in $K[\![x,y]\!]$.*

Clear, since $\overline{f}(x,y)$ is a polynomial in $y$ whose coefficients are invariant under the Galois group. Note that this polynomial is the minimal polynomial of the element $f$ of $K[\![x^\star]\!][x^{\star -1}]$ as an algebraic extension of $K[\![x]\!][x^{-1}]$. ⊡

**Lemma.** *Suppose a formal Puiseux series $y = f(x) \in K[\![x^\star]\!]$ is a solution to the equation $F(x,y) = 0$ where $F \in K[\![x,y]\!]$. Then the associated series $\overline{f}(x,y) \in K[\![x,y]\!]$ divides $F(x,y)$ in the ring $K[\![x,y]\!]$.*

Since $f$ is a solution and the equation is invariant under the Galois group, all the conjugates are also solutions. We then show, using the first lemma $m$ times that $F$ is divisible by $\overline{f}(x,y)$ in $K[[x^{1/m},y]]$. Now the quotient $F/\overline{f}$ is Galois invariant so that it is actually in $K[[x,y]]$. ⊡

We can now prove the so-called *Weierstrass preparation theorem*, for formal series.

**Theorem.** *Let $F(x,y) \in K[[x,y]]$. Assume that $F$ is not divisible by $x$. Then $F$ can be written as a product $A(x,y)P(x,y)$ where $A, P$ are in $K[[x,y]]$ and*

- *$A(0,0) \neq 0$ so that $A$ is an invertible element.*
- *$P(x,y)$ is a polynomial in $y$.*

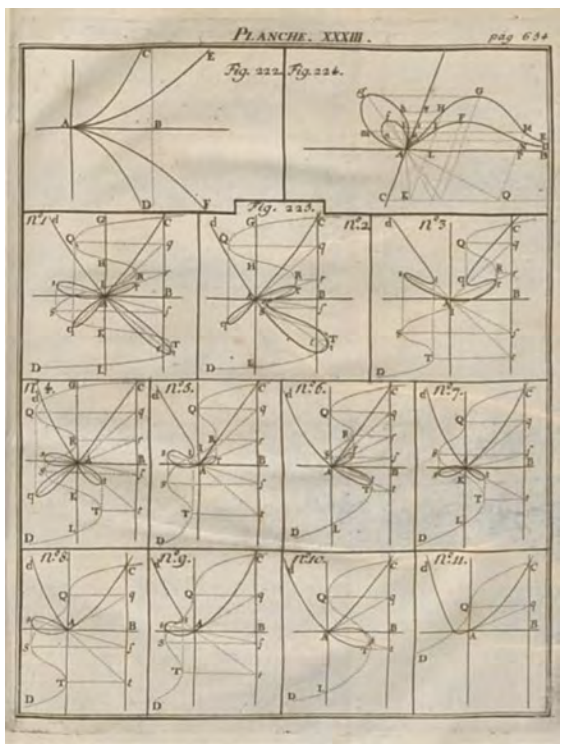

The proof is by induction on the multiplicity $mult(F)$. Note that the valuation of a product is the sum of the valuations and that $mult(F) = 0$ means precisely that $F(0,0) \neq 0$. If $mult(F) \geq 1$, we know that $F$ has at least one solution in $y = f(x) \in K[[x^\star]]$ and that $F$ is divisible by $\overline{f}(x,y) \in K[[x,y]]$. The quotient has a lower multiplicity. ⊡

Now, we can harvest and state two corollaries that follow easily from the previous theorem. The proofs are the same as in the classical case of polynomial rings over fields.

**Theorem.** *Any nonzero element $F$ of $K[[x,y]]$ can be split as*

$$F = A(x,y)x^r(y - f_1(x))(y - f_2(x))\cdots(y - f_n(x))$$

*where $A \in K[[x,y]]$ is such that $A(0,0) \neq 0$, the $n$ solutions $f_i(x)$ are in $K[[x^\star]]$, and $r \geq 0$.*

**Theorem.** *The ring $K[[x,y]]$ is a unique factorization domain.*

The irreducible factors are the $\overline{f_i} \in K[[x,y]]$.

We conclude this chapter with two exercises.

*Exercise.* Newton's algorithm produces solutions $y(x)$. At each step, we have to choose one of the supporting lines on the boundary of the polygon, and a root of the corresponding polynomial equation. Show that this algorithm produces *all solutions* $f_i(x)$ of $F(x,y) = 0$.

*Exercise.* Suppose we follow Newton's algorithm using some choices of segments, leading eventually to a solution $y = f(x)$. The process produces a sequence of formal series $F_k(x_k, y_k)$. We proved that the multiplicities $mult(F_k)$ of $F_k$ are eventually constant. Show that this "eventual multiplicity" is just the multiplicity of the root, that is to say the number of factors equal to $(y - f(x))$ in the above decomposition of $F$ as a product $A(x,y)x^r(y - f_1(x))(y - f_2(x))\cdots(y - f_n(x))$.

Enough algebra for the time being!

> Il est fâcheux que Mr. Newton se soit contenté d'étaler ses découvertes sans y joindre de démonstrations et qu'il ait préféré le plaisir de se faire admirer à celui d'instruire.

From Cramer's preface.

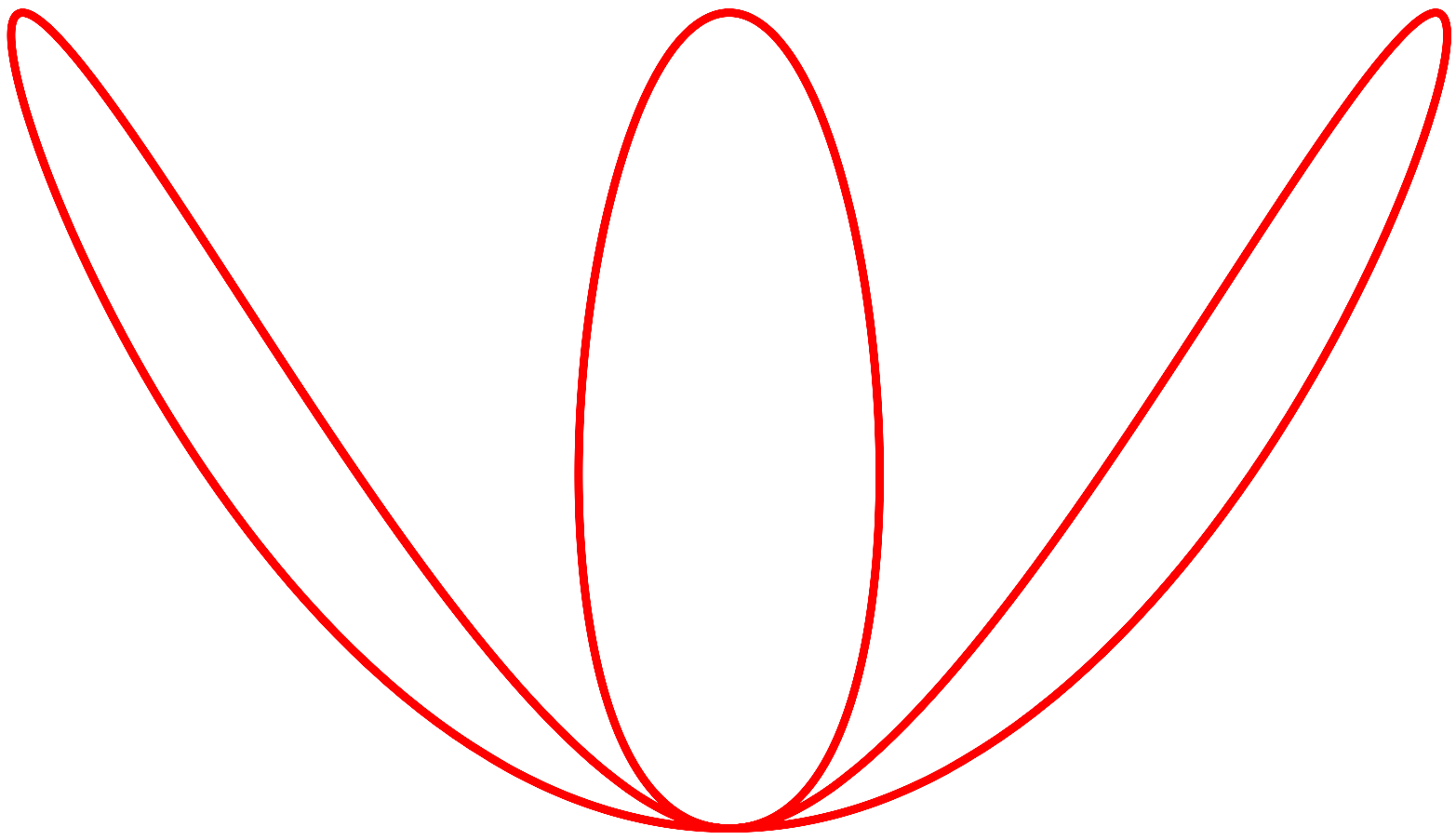

This is one of the sextic curves studied by Cramer, whose equation is

$$y^6 - (y - x^2)(y - 4x^2)^2 = 0.$$

The singular point has 3 branches. Can you draw these branches?

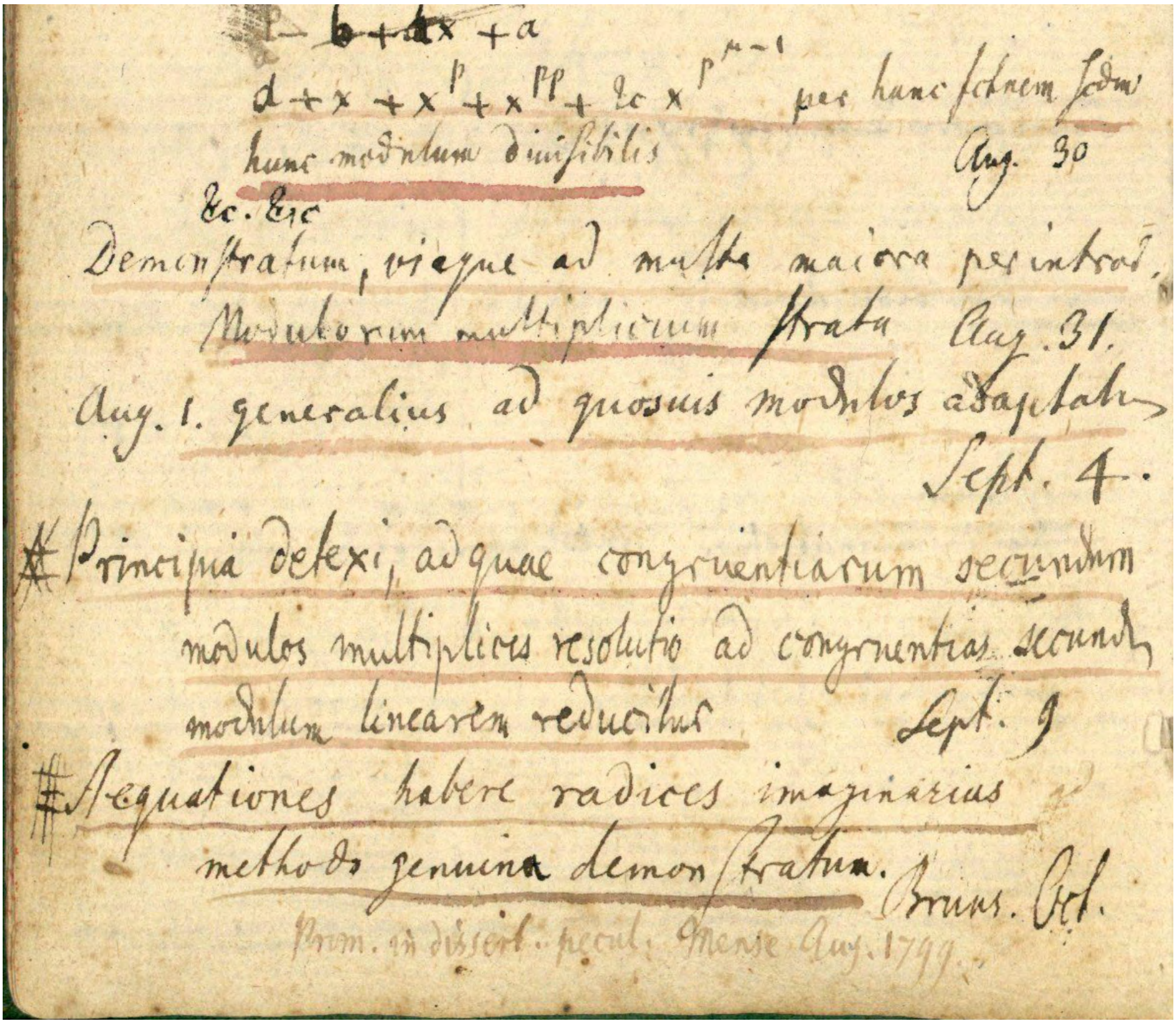
hunc modulum divisibilis Aug. 30

&c. &c.

Demonstratum, viaque ad multo maiora per introd. Modulorum multiplicium strata Aug. 31.

Aug. 1. generalius ad quosuis modulos adaptatum Sept. 4.

# Principia detexi, ad quae congruentiarum secundum modulos multiplices resolutio ad congruentias secundum modulum linearem reducitur Sept. 9

# Aequationes habere radices imaginarias methodo genuina demonstratum. Bruns. Oct.

Prom. in dissert. pecul. mense Aug. 1799.

Starting from 1796 (when he was 19 years old) Gauss recorded his mathematical discoveries in his famous Tagebuch. An impressive list of results. See Klein's commentaries[33] and [34] for an English translation. This page concerns August to October 1797. The last item *"Aequationes habere radices imaginarias methodo genuina demonstratum"* announces his proof of the fundamental theorem of algebra. Below this line, with a different ink, a later addition mentions that this was the theme of his dissertation: *"Prom[ulgatum] in dissert[atione] pecul[iari] mense Aug. 1799"* ©

[33] F. Klein. Gauß' wissenschaftliches Tagebuch 1796–1814. *Math. Ann.*, 57(1):1–34, 1903.

[34] J. J. Gray. A commentary on Gauss's mathematical diary, 1796–1814, with an English translation. *Exposition. Math.*, 2(2):97–130, 1984.

# *Curuam algebraicam neque alicubi subito abrumpi posse: Gauss on algebraic curves*

## *The fundamental theorem of Algebra*

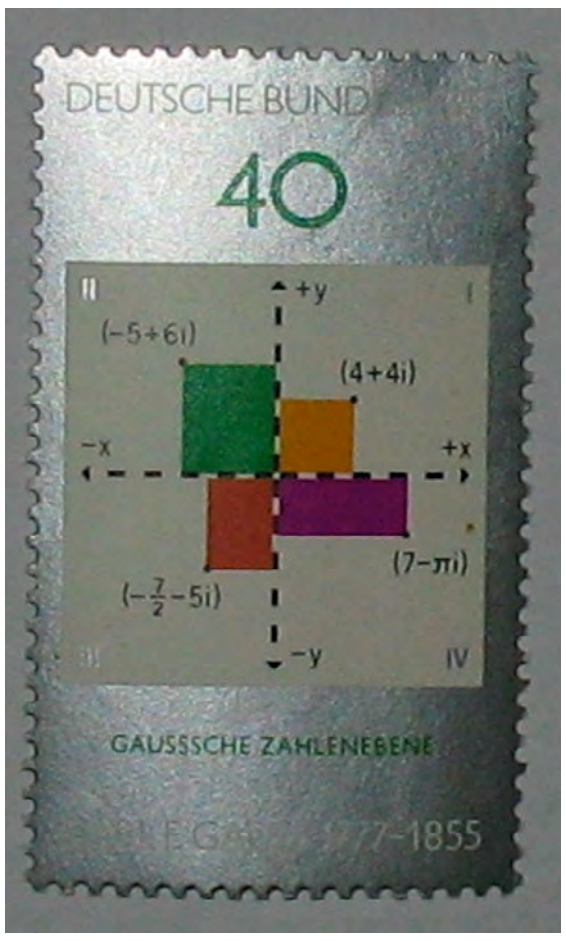

A stamp commemorating Gauss's complex plane. ©

Carl Friedrich Gauss was 22 years old when he defended his thesis in 1799. This is a remarkable piece of work[35] containing what may possibly be considered as the first "proof" of the fundamental theorem of algebra.

[35] C. F. Gauß. *Werke. Band III*. Georg Olms Verlag, Hildesheim, 1973. Reprint of the 1866 original.

*Any non-constant polynomial with complex coefficients has at least one root.*

In slightly different terminology, and not using the words "complex" or "imaginary", which were suspicious at that time, he proved that any real polynomial is a product of factors of degrees 1 or 2. In a different language, the title of his PhD is:

> DEMONSTRATIO NOVA THEOREMATIS OMNEM FVNCTIONEM ALGEBRAICAM RATIONALEM INTEGRAM VNIVS VARIABILIS IN FACTORES REALES PRIMI VEL SECUNDI GRADVS RESOLVI POSSE

Gauss received his degree from the university of Helmstedt. His formal advisor was Johann Friedrich Pfaff who read carefully the dissertation. However, this doctorate was *in absentia*: there was no oral presentation. The manuscript mentions that the main result was obtained in October 1797. An English translation of the thesis by Ernest Fandreyer is available online.

This is not a proof by today's standards, but I will present a slight variation on the same theme which is perfectly acceptable by 21st century mathematicians. It was not the first attempt of a proof. Among Gauss's predecessors, one might mention d'Alembert, Euler and Lagrange. None of these previous

"proofs" were solid, even at that time, but I will try to reconstruct d'Alembert's proof since he used Newton's polygon.

The first half of Gauss's thesis deals with a criticism of his predecessors. He carefully explains why the proofs of d'Alembert, Euler and Lagrange are flawed. It is hard to imagine a similar situation today of a very young man defending his PhD and beginning by a systematic destruction of great and respected Masters who had passed away only fifteen years earlier, or even were still alive (as in the case of Lagrange). Then, in a second part, Gauss gives his proof. Beautiful proof, indeed, but not totally exempt of "unproved facts". At a crucial moment, to be described later, he needs a fairly precise description of the local structure of a real algebraic curve. He then asserts, with no proof, that

Gauss's "unproved facts" have nothing to do with Trump's "alternative facts": after all they are true.

> But according to higher mathematics, any algebraic curve (or the individual parts of such an algebraic curve if it perhaps consists of several parts) either turns back into itself or extends to infinity. Consequently, a branch of any algebraic curve which enters a limited space, must necessarily exit from this space somewhere.

*"Iam ex geometria sublimori constat, quamuis curuam algebraicam, (siue singulas cuiusuis curuae algebraicae partes, si forte e pluribus composita sit) aut in se redientem aut vtrimque in infinitum excurrentem esse, adeoque si ramus aliquis curuae algebraicae in spatium definitum intret, eundem necessario ex hoc spatio rursus alicubi exire debere."*

In other words, Gauss is claiming that an algebraic curve cannot simply stop at some point. The "proof" is given in a footnote: it is a typical example of a proof by intimidation:

> It seems to have been proved with sufficient certainty that an algebraic curve can neither be broken off suddenly anywhere (as happens e.g. with the transcendental curve whose equation is $y = 1/\log x$) nor lose itself, so to say, in some point after infinitely many coils (like the logarithmic spiral). As far as I know, nobody has raised any doubts about this. However, should someone demand it then I will undertake to give a proof that is not subject to any doubt, on some other occasion.

*"Satis bene certe demonstratum esse videtur, curuam algebraicam neque alicubi subito abrumpi posse (vti e. g. euenit in curua transscendente, cuius aequatio $y = 1/\log x$), neque post spiras infinitas in aliquo puncto se quasi perdere (vt spiralis logarithmica), quantumque scio nemo dubium contra rem mouit. Attamen si quis postulat, demonstrationem nullis dubiis obnoxiam alia occasione tradere suscipiam."*

"Nobody has raised doubts and he will prove it on some other occasion" ☺ ! Actually, he never proved this fact (even though he published later three other proofs of the fundamental theorem of algebra, as if he was himself not convinced). What an arrogant (and brilliant) young man!

Gauss gives two examples of curves. The first is the graph of $1/\log(x)$ and the is the logarithmic spiral ($r = \exp(\theta)$ in polar coordinates). Both can be defined by some equation $F(x,y) = 0$ and both have some kind of *stopping point*. If we draw a small disk around this point, the curve enters this disk but does not exit. The (correct) claim of Gauss is that this is due to the transcendental nature of these curves and that this does not happen for algebraic curves for which $P(x,y)$ is a polynomial.

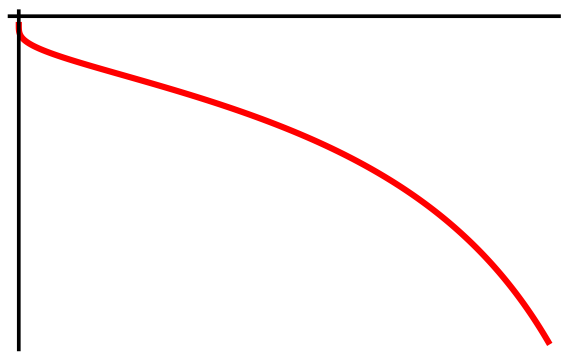

The curve $x\exp(-1/y) = 1$ (i.e. $y = 1/\log x$) has a "dead end" at $(0,0)$.

## *A reconstruction of the proof by Gauss*

My intention is certainly not to analyze this proof from a historical point of view. There would be much to be discussed: the concept of continuity, of curve, the topological arguments, and above all the geometrical use of complex numbers as points in a plane. I recommend the books by Dhombres and Alvarez[36] and Van der Waerden[37]. Let me only mention a lucid point of view expressed by Gauss, more than twenty years before Abel and Galois. In the following, what he calls a *pure equation* is an equation of the form $x^n = a$.

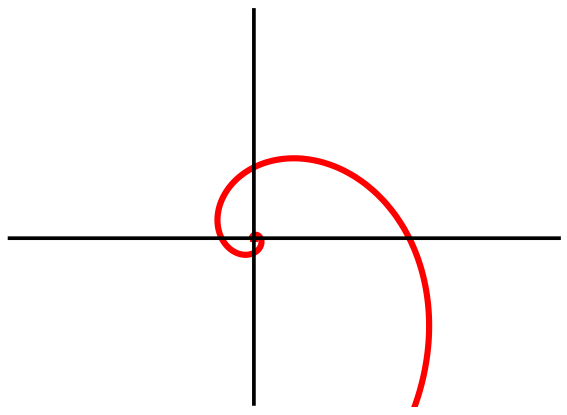

The logarithmic spiral with equation $y - x\log\tan(x^2 + y^2) = 0$ ($\rho = \exp(-\theta)$ in polar coordinates) has infinitely many coils as it converges to the origin.

> [...] after so much labor of such great mathematicians there is very little hope left ever to arrive at a general solution of algebraic equations. It seems more and more probable that such a solution is entirely impossible and contradictory. This must not at all be considered paradoxical, as that which is commonly called the solution of an equation is indeed nothing other than its reduction to pure equations. For the solution of pure equations is here not taught but presupposed; and if you express the roots of an equation $x^m = H$ by $\sqrt[m]{H}$, you have in no way solved it, and you have not done more than if you had devised some symbol to denote the root of an equation $x^h + Ax^{h-1} + etc. = 0$ and set the root equal to this.

[36] J. Dhombres and C. Alvarez. *Une histoire de l'invention mathématique : les démonstrations du théorème fondamental de l'algèbre dans le cadre de l'analyse réelle et de l'analyse complexe de Gauss à Liouville*. Hermann, 2013.

[37] P. D. B. L. van der Waerden. *A History of Algebra: From al-Khwarizmi to Emmy Noether*. Springer-Verlag Berlin Heidelberg, 1 edition, 1985.

My modest purpose is to propose a modern reconstruction of the proof, showing why Gauss needed some understanding of the local nature of algebraic curves. Let $P(z)$ be a monic polynomial of degree $n \geq 1$ with complex coefficients. The main idea is to think of $z = x + iy$ as a point in the plane and of $P(x+iy)$ as $p(x,y) + iq(x,y)$, defining two real polynomials in $(x,y)$. Proving the existence of a complex root of $P$ is equivalent to showing that the two algebraic curves $p(x,y) = 0$ and $q(x,y) = 0$

Note that the "simple idea" of thinking of a polynomial as a map from a plane to another plane was a new idea in 1797.

have a non empty intersection. We are going to analyze the qualitative behavior of these two curves in the neighborhood of infinity.

When the modulus of $z$ is large, $P(z)$ and $z^n$ are equivalent, so that we consider as a first approximation the curves

$$\Re(x+iy)^n = 0 \quad ; \quad \Im(x+iy)^n = 0.$$

These equations are easy to solve: they define radial lines

$$\arg z = \frac{(2k+1)\pi}{2n} \quad (0 \le k \le 2n-1); \quad \arg z = \frac{2k\pi}{2n} \quad (0 \le k \le 2n-1),$$

which intersect at the origin. These $2n$ lines intersect each circle $|z| = R$ at $4n$ points. The first claim of Gauss is the following:

**Lemma.** *When $R$ is large enough, each of the two algebraic curves $p(x,y) = 0$ and $q(x,y) = 0$ intersects the circle $|z| = R$ at $2n$ points which are close to the previous ones.*

The real and imaginary parts of $\frac{1}{R^n}P(R\exp(i\theta))$ are trigonometric polynomials of degree $n$ in the variable $\theta$ which are close to $\cos(n\theta)$ and $\sin(n\theta)$. Therefore each one vanishes at most $2n$ times and they do vanish $2n$ times by the intermediate value theorem. Elementary details are left to the reader. The proof of this point by Gauss is perfect. ⊡

Now comes the topological part of the proof.

Suppose first that the algebraic curves $p(x,y) = 0$ and $q(x,y) = 0$, that will be called the blue and the red curves, are smooth. Inside the disc $|z| \le R$ they consist of a finite number of arcs, each diffeomorphic to $[0,1]$ and a certain number of loops, diffeomorphic to a circle. This follows from the classification of compact one dimensional manifolds (see for instance[38] or [39]). There are $4n$ points on a circle, blue and red, with alternating colors. We will say that two points (of the same color) are *paired* if they are boundary points of one of these blue or red arcs, *inside the disc*. So our set of $4n$ points consists of $2n$ pairs.

Consider four distinct points on the circle, two of them colored in red and the other two in blue. From the topological point of view, there are two possibilities. They could be *linked* or *unlinked*. Going around the circle, we read alternate colors,

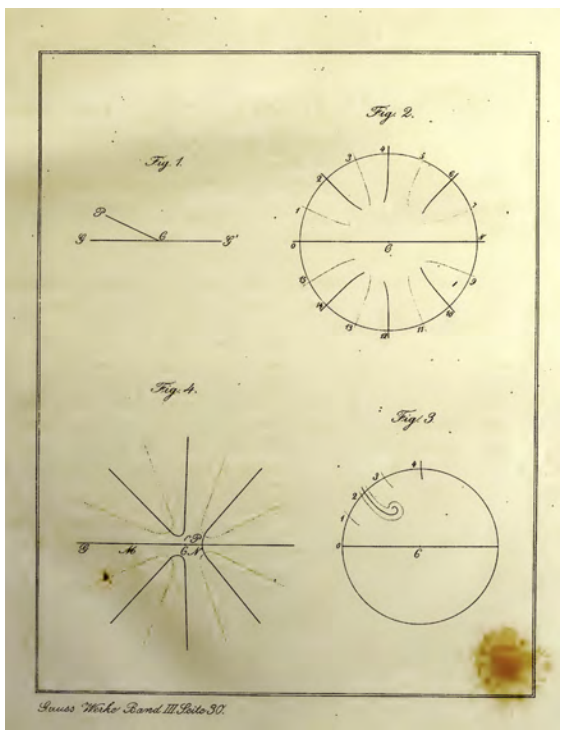

This is the only figure from Gauss's dissertation. ©

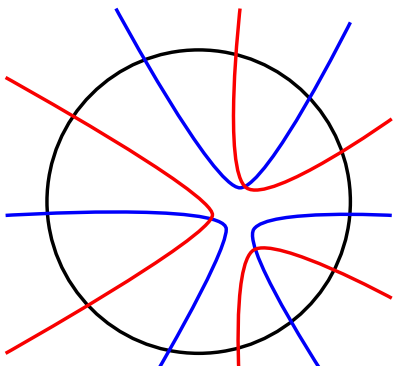

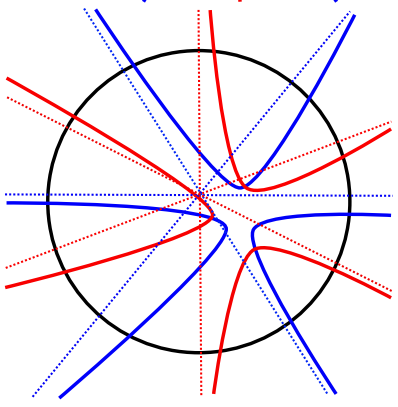

38 J. W. Milnor. *Topology from the differentiable viewpoint.* Princeton Landmarks in Mathematics. Princeton University Press, Princeton, NJ, 1997.

39 V. Guillemin and A. Pollack. *Differential topology.* AMS Chelsea Publishing, Providence, RI, 2010. Reprint of the 1974 original.

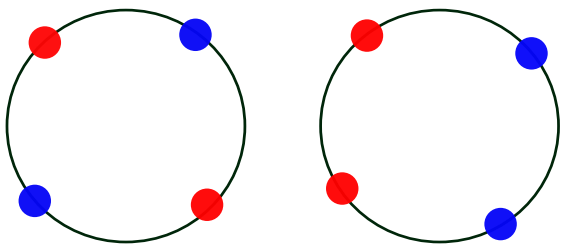

Linked and unlinked.

like "blue, red, blue, red" in the linked case, and "blue, blue, red, red" in the unlinked case. The crucial topological lemma, which is intuitive is the following.

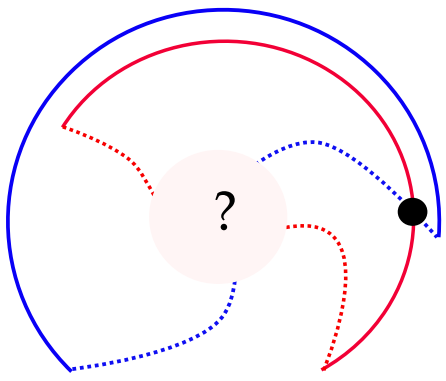

**Lemma.** *Let $b_0, b_1, r_0, r_1$ be four points on the circle such that $\{b_0, b_1\}$ and $\{r_0, r_1\}$ are linked. Let $b$ (resp. $r$) be a smooth arc in the disc connecting $b_0$ and $b_1$ (resp. $r_0$ and $r_1$). Then $b$ and $r$ intersect non-trivially.*

This follows from one of the very first theorems in algebraic topology (therefore not formally at the disposal of Gauss). *Two closed curves in the plane which intersect transversally have an* even number *of intersection points (see for instance Milnor's book).* If there could exist *disjoint* arcs connecting the blue and the red points *inside the disc*, we could construct two closed loops in the plane intersecting in exactly one point (see the figure in the margin). ⊡

The fact that *two transversal closed curves in the plane intersect in an even number of points* is more or less equivalent to Jordan's theorem: *"the complement of a closed embedded curve in the plane has exactly two connected components"*. Indeed, if $c_1, c_2$ are closed and transversal, one can first modify $c_1$ slightly, without changing its intersection with $c_2$, in such a way that $c_1$ becomes an immersion with transversal self-intersections. Then one modifies $c_1$ as in the picture below, again without changing the intersection with $c_2$, in order to replace it by a disjoint union of closed *embedded* curves. Now by Jordan's theorem, each time $c_2$ enters a connected component of the complement of (the modified) $c_1$, it has to exit, so that there is indeed an even number of intersections. Try to prove Jordan's theorem from the parity of intersection.

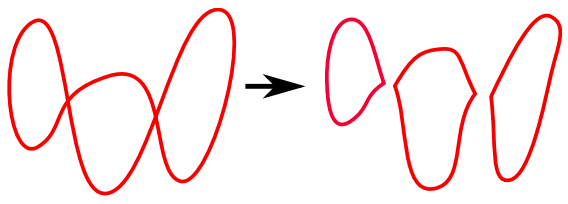

**Lemma.** *Suppose $2k$ persons sit around a table and they shake their hands two by two, without crossing arms! Then, at least two neighbors shake their hands.*

For $k = 2$, this is the previous lemma. Consider two persons shaking hands. They decompose the boundary of the table in two intervals. If one is empty, we are done. Otherwise, proceed by induction. ⊡

Still assuming that the two curves $p(x, y) = 0$ and $q(x, y) = 0$, blue and red, are smooth, we prove the fundamental theorem of algebra, following Gauss. By contradiction, assume that the blue and red arcs do not cross. By the previous lemma two neighbors on the circle are paired. This is impossible since consecutive points do not have the same color. ⊡

Now, we understand the difficulty for which "Nobody has raised doubts". If there were an algebraic curve with a dead end, an arc could penetrate *inside the disc* and stop there, without exiting and that would be fatal for the proof.

Let us make Gauss's claim precise.

**Theorem.** *Let $(x_0, y_0)$ be a point on some real algebraic curve $C$ defined by $F(x, y) = 0$ where $F$ is a real polynomial in $\mathbb{R}[x, y]$. Then there is a homeomorphism of some small disc centered in $(x_0, y_0)$ sending $C$ to the union of an* even number *of distinct radii.*

Prove that the number of "non-crossing pairings" of an even number of points on a circle is a Catalan number.

This claim is indeed true and "I will undertake to give a proof that is not subject to any doubt, on some other occasion."

You'll have to wait until the next chapter!

*Assuming this is true*, it is easy to finish the proof. If the blue and red curves $p(x,y) = 0; q(x,y) = 0$ are singular and disjoint, it suffices to modify them slightly, as in the margin, locally in disjoint small neighborhoods of all singular points, connecting the radii in pairs, so that they become *disjoint non-singular arcs*. We have seen that this is not possible. ⊡

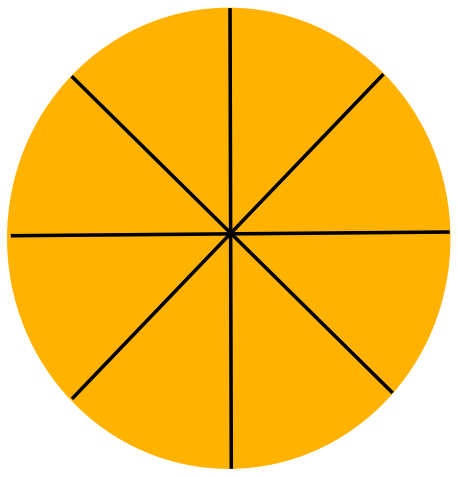

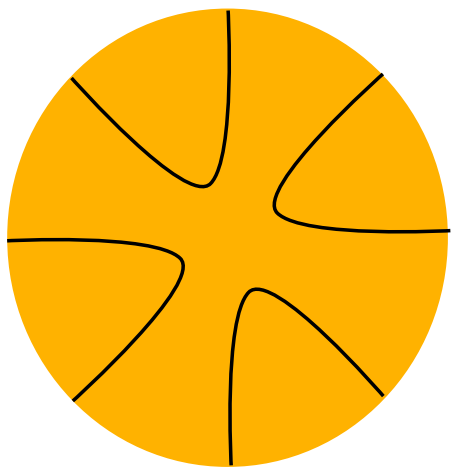

## *Comments on this proof*

Steve Smale presented this proof in a paper dealing with effective versions of the fundamental theorem of algebra[40]. He emphasized Gauss's unproved claim:

> But for the moment, I wish to point out what an immense gap Gauss's proof contained. It is a subtle point even today that a real algebraic plane curve cannot enter a disk without leaving.

He also comments on the endless debate about who gave the "first" accepted proof.

> One can understand the historical situation better perhaps from the point of view of Imre Lakatos[41]. Lakatos in the tradition of Hegel, on one hand, and Popper, on the other, sees mathematics as a development which proceeds as a series of 'proofs and refutations'.

There are many ways to "fix" the proof and to fill the "immense gap". First I should mention the long detailed paper by Ostroswki, dated 1920, fully dedicated to the proof of Gauss's claim[42]. The curves $p(x,y) = 0$ and $q(x,y) = 0$ used by Gauss are indeed algebraic curves, but they are very special algebraic curves. In modern terminology, these polynomials are real and imaginary parts of a holomorphic function $P(z)$ and are therefore *harmonic* polynomials. The detailed proof by Ostrowski actually deals with harmonic polynomials, which is sufficient for our present problem. With elementary notions on complex analysis, it is indeed easy to fill the details, as I show now.

[40] S. Smale. The fundamental theorem of algebra and complexity theory. *Bull. Amer. Math. Soc. (N.S.)*, 4(1):1–36, 1981.

[41] I. Lakatos. *Proofs and refutations*. Cambridge Philosophy Classics. Cambridge University Press, Cambridge, paperback edition, 2015. Originally published in 1976.

[42] A. Ostrowski. Über den ersten und vierten Gauss'schen Beweis des Fundamental satzes der Algebra. *Nachrichten der Gesellschaft der Wissenschaften Göttingen*, 1920.

Think of $P(z) = P(x+iy) = p(x,y) + iq(x,y)$ as a map from $\mathbb{C} \simeq \mathbb{R}^2$ to another copy of itself. The differential of this map $P$ at a point $z_0 = x_0 + iy_0$ can be seen either as a $2 \times 2$ real matrix or as the complex number $P'(z_0)$. Hence, critical points are simply the finitely many zeroes of the derivative $P'$. The blue and red curves are the inverse images of the two axes. Let us analyze the inverse image by $P$ of some line.

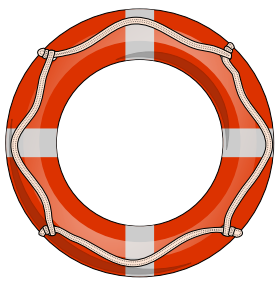

In case of emergency, the book[43] can be helpful to understand these pictures.

[43] T. Needham. *Visual complex analysis*. The Clarendon Press, Oxford University Press,

In the neighborhood of some point $z_0$, we have

$$P(z) - P(z_0) = c_k(z-z_0)^k + c_{k+1}(z-z_0)^{k+1} + \cdots + c_n(z-z_0)^n$$

for some $k \geq 1$ (the valuation of $P(z) - P(z_0)$ at $z_0$). Hence

$$P(z) - P(z_0) = \left( (z-z_0)\sqrt[k]{c_k}\sqrt[k]{1 + \frac{c_{k+1}}{c_k}(z-z_0) + \cdots} \right)^k = \phi(z)^k.$$

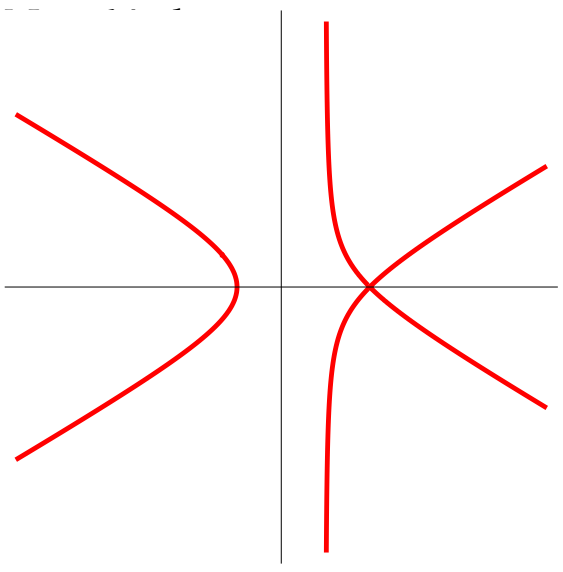

Inverse image of the vertical axis by $2z^3 - 3z^2 + 1 + i$. The critical points are $z = 0, 1$ and the critical values are $1 + i$ and $i$: one of them is on the vertical axis.

Here, $\sqrt[k]{c_k}$ is any choice of the $k$-th root and the second $k$-th root is a convergent power series by Newton's binomial theorem. The differential at $z_0$ of $\phi$ is not zero, so that $\phi$ is a local diffeomorphism. In short, $P(z)$ is the local composition of a diffeomorphism and of the map $z \mapsto P(z_0) + (z-z_0)^k$. It is therefore obvious that the inverse image by $P$ of a smooth curve going through $P(z_0)$ is the union of $k$ smooth curves through $z_0$. In particular, locally there are $2k$ half lines, and this proves Gauss's claim in the special case that he needed. This special case is indeed very special since these $k$ smooth curves make equal angles. ⊡

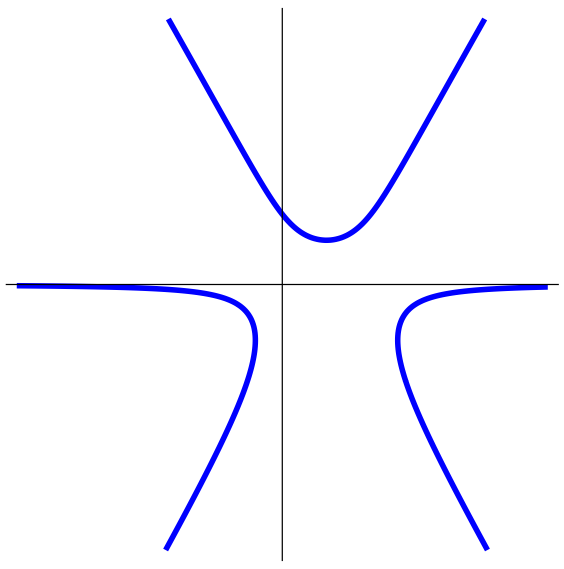

Inverse image of the horizontal axis by $2z^3 - 3z^2 + 1 + i$. There are no critical value on this axis.

*But do not forget that we still did not prove Gauss's claim in its full generality.*

There is another way to fill Gauss's "immense gap". Rotating the axis by an angle $\theta$, we can replace $P(z)$ by $\exp(i\theta)P(z)$. The curve $p(x,y) = 0$ (resp. $q(x,y) = 0$) is singular if and only if one of the critical values of $P$ is on the real (resp. imaginary) axis. Hence it suffices to rotate by a suitable $\theta$ to avoid this, so that Gauss could as well have started with the assumption that the blue and red curves are smooth. This easy argument was not easy in 1797.

Today, there are many proofs of the fundamental theorem of algebra. I recommend Eisermann's paper[44] for a lucid overview.

[44] M. Eisermann. The fundamental theorem of algebra made effective: an elementary real-algebraic proof via Sturm chains. *Amer. Math. Monthly*, 119(9):715–752, 2012.

This proof by Gauss is certainly neither the most direct nor the easiest. Cleaning it requires some subtle topological arguments but on the way we get ample rewards, as we understand much better complex polynomials as maps.

Let me present my favorite proof, in the spirit of Gauss's topological proof, that can be found in Smale's above mentioned paper. For me, this is the simplest[45] one. Choose a point $z_0$ in such a way that the segment $\gamma$ connecting $0$ to $P(z_0)$ does not contain one of the finitely many critical values of $P$. This is possible if $0$ is not a critical value but if this is the case $0$ is a value so that $P$ has a root. Look at the inverse image of $\gamma$ by $P$. This is a smooth compact manifold of dimension 1, with boundary. The point $z_0$ is a boundary point of one component. The other boundary point of this component is clearly a root of $P$. Voilà! ⊡

What do I mean by "simplest"? Probably not the shortest since this proof contains a lot of implicit facts that should be proven. Simplicity is a subtle and very personal concept in mathematics. In this special case, I would say that this is simple because I think I could not forget it.

[45] É. Ghys. Inner simplicity vs. outer simplicity. In R. Kossak and P. Ording, editors, *Simplicity: Ideals of Practice in Mathematics and the Arts Conferences (CUNY New York, 2013)*. Springer Verlag, 2017.

This simple proof actually gives much more. Away from the critical values, we can pull back the radial vector field $-x\partial/\partial x - y\partial/\partial y$ by the differential of $P$. We get a vector field in the plane, away from the critical points. The trajectories of this vector field are precisely the inverse images of the radial lines. Hence, starting from a point and solving this differential equation, we should arrive at the roots of $P$. One way to approximate the solutions of an ODE is to use the standard Euler method. It turns out that the Euler iterative scheme coincides with Newton's method. Newton, Gauss and Euler together!

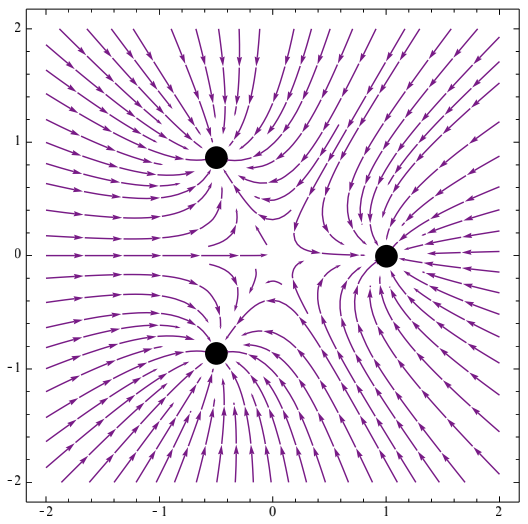

The phase portrait of the vector field $-P(z)/P'(z)$ for $P(z) = z^3 - 1$. Trajectories are mapped by $P$ on radial lines.

## *A proof by d'Alembert*

I also describe a proof by d'Alembert[46] for two reasons. The first is that in France the fundamental theorem of algebra is often called d'Alembert's theorem ☺. The second is that this is closely related to Newton's polygons that we analyzed earlier. See[47] and [48] for much more on this proof. D'Alembert does not mention Newton. How could a Frenchman acknowledge the contribution of an Englishman?

Let me present simplified version of his "proof". Suppose we

[46] J. D'Alembert. Recherches sur le calcul intégral. *Histoire de l'Acad. Royale Berlin*, pages 182–224, 1748.

[47] C. Baltus. D'Alembert's proof of the fundamental theorem of algebra. *Historia Math.*, 31(4):414–428, 2004.

[48] C. Gilain. Sur l'histoire du théorème fondamental de l'algèbre: théorie des équations et calcul intégral. *Arch. Hist. Exact Sci.*, 42(2):91–136, 1991.

want to solve

$$z^n + a_{n-1}z^{n-1} + \cdots + a_0 = 0.$$

Set $z = y/\varepsilon$. We get a strictly equivalent equation:

$$y^n + a_{n-1}\varepsilon y^{n-1} + \cdots + \varepsilon^n a_0 = 0.$$

Of course, $y = 0$ is a solution for $\varepsilon = 0$ and we want a solution for $\varepsilon \neq 0$. Consider the above equation as an equation of the form $F(\varepsilon, y) = 0$. By Newton and Cramer, there are non-trivial solutions $y(\varepsilon)$, at least for small values of $\varepsilon$, expressed as Puiseux series in $\varepsilon$. "Therefore", we found a root of our equation. ⊡

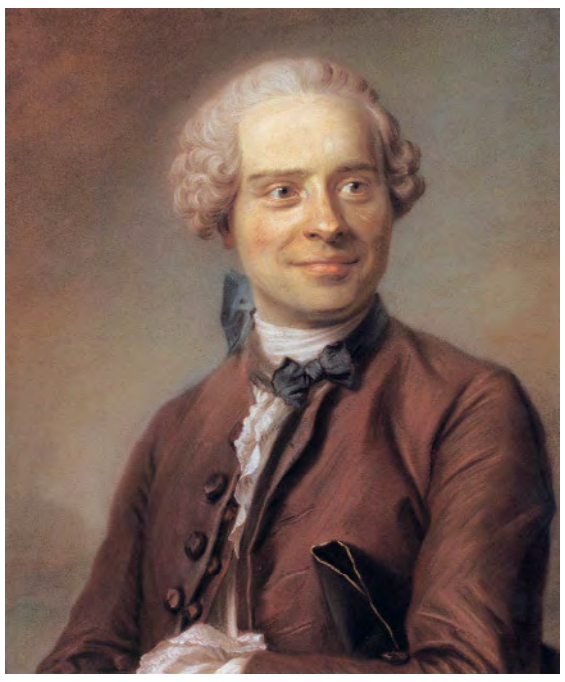

Jean Le Rond d'Alembert (1717–1783). ©

The previous "therefore" is subject to a lot of discussion. One of the main difficulties is that neither Newton, nor Cramer, nor d'Alembert proved the convergence of the series. Even worse, Newton's algorithm constructing the Puiseux series uses the fundamental theorem of algebra. A vicious circle.

Amazingly d'Alembert also published some version of this proof in a memoir dealing with *the cause of winds*[49].

[49] J. D'Alembert. *Réflexions sur la cause genérale des vents*. David l'ainé, Paris, 1747.

Let me conclude this chapter by an exercise, suggested by my former student Victor Kleptsyn. Look at the inverse images of the real axis (say in red) and the imaginary axis (in blue) by some complex polynomial $P(z)$. This produces some graph in a big disc. Each edge is colored in blue or red. Singular points of the blue (resp. red) graph are critical points of $P$ which are mapped to the real (resp. imaginary) axis: they present an even number of blue (resp. red) edges going out of a vertex. Generically, there is no such singular point. The local picture around the intersection of the two graphs has been described above: $4k$ edges going out of the vertex, cyclically alternating blue and red. These intersections do exist by the fundamental theorem of algebra. On large circles, we have alternation between red and blue.

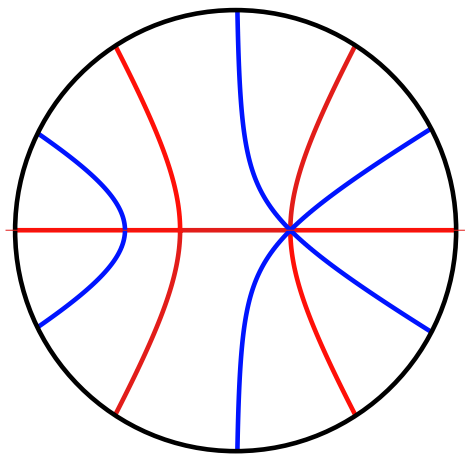

Inverse image of the real and imaginary axes by $z^3 - 3z + 2$. There are two critical values: 0 and 4.

The question concerns the converse. Suppose we have a colored graph in a disc presenting all the previous qualitative features. Under which conditions does there exist some polynomial $P(z)$ such that its associated colored graph is homeomorphic to the given graph, under some homeomorphism of the disc?

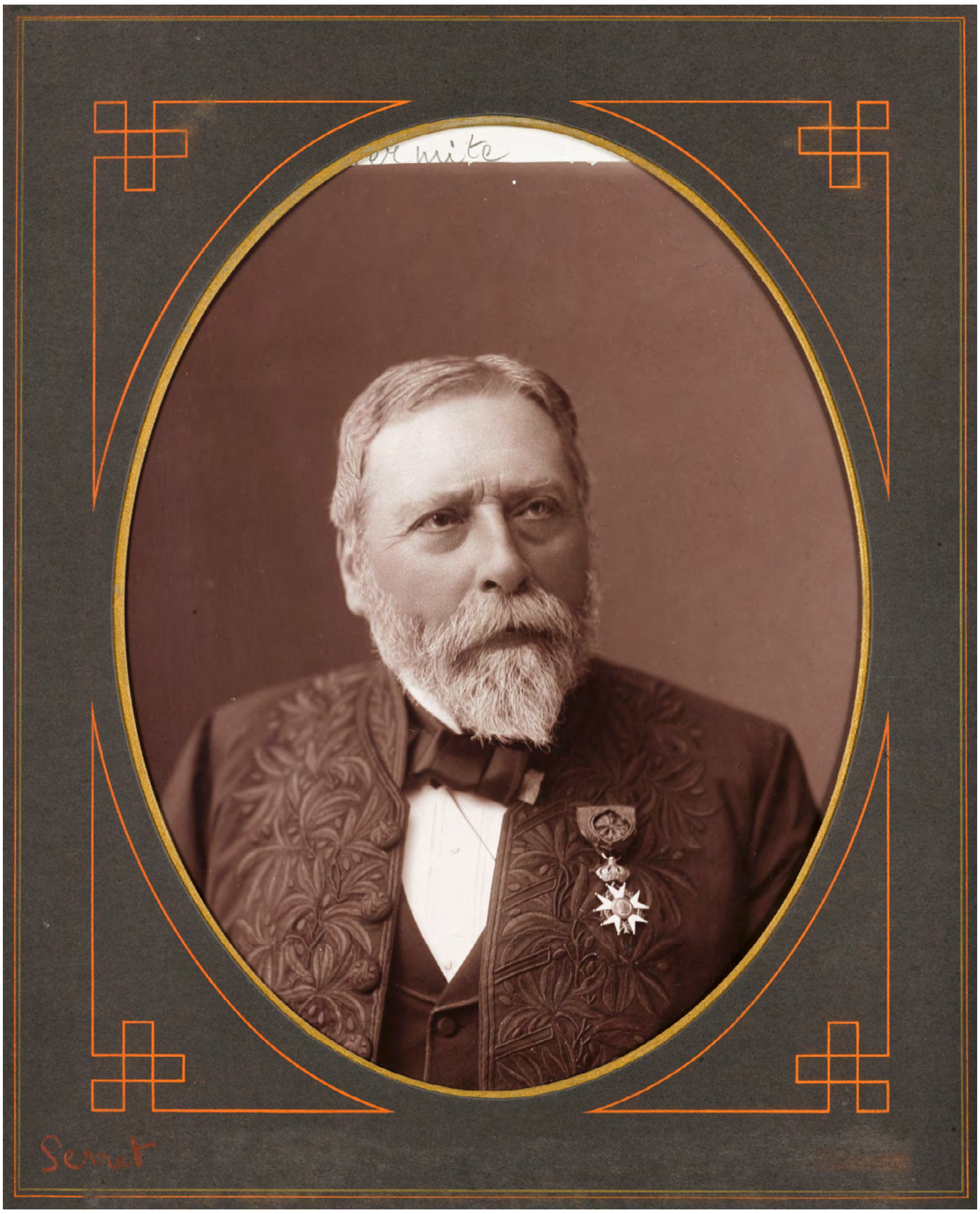

Joseph Alfred Serret (1819–1885). ©

# *Proof of Gauss's claim on singularities of algebraic curves: two papers by two Serret's*

IT IS TIME TO PROVE GAUSS'S ASSERTION: "the neighborhood of a point of a planar real algebraic curve is homeomorphic to an *even* number of radii in a disc".

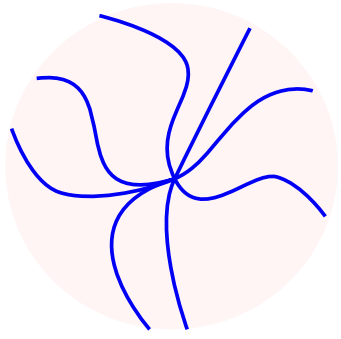

Local picture of an algebraic curve.

## *Insufficient proofs*

Joseph Alfred Serret (1819–1885) should not be confused with Paul Joseph Serret (1827–1898).

Joseph Alfred had a brilliant career. He signed his books as "Membre de l'Institut et du bureau des longitudes, Professeur au Collège de France et à la Faculté des sciences de Paris". In 1849, he published a very influential *Cours d'algèbre supérieure* in two volumes containing one of the first systematic expositions of Galois theory. He is also at the origin of the Frenet-Serret frame for curves in 3-space.

The younger Paul Joseph had a much more modest career. He signed his books "Docteur ès sciences, membre de la société philomatique". He taught in collège Sainte-Barbe in Paris. I could not find his portrait.

In 1847, Joseph Alfred wrote a paper[50] in which he "proves" an assertion from Newton:

> If a straight line is asymptotic to a branch of an algebraic curve, then it is also asymptotic to another branch.

[50] J. A. Serret. Théorème sur les courbes algébriques asymptotiques. *Nouvelles annales de mathématiques, journal des candidats aux écoles polytechnique et normale*, 6:217–218, 1847.

According to Joseph Alfred, *"Ce théorème est dû à Newton, et est énoncé, si je ne me trompe, dans son* Enumeratio Linearum tertio ordains*"*. *"This theorem is due to Newton and is stated, if I am not mistaken, in his* Enumeratio Linearum tertio ordains.*"*

Note that what Newton calls here a *branch* is one half of what we call a branch. As a simple example, look at Descartes's folium $x^3 + y^3 = 3xy$. Its asymptote is approached by the curve as $x$ tends to $+\infty$ *and* $-\infty$. This corresponds to *two* branches in Newton's terminology and to one branch at infinity in ours.

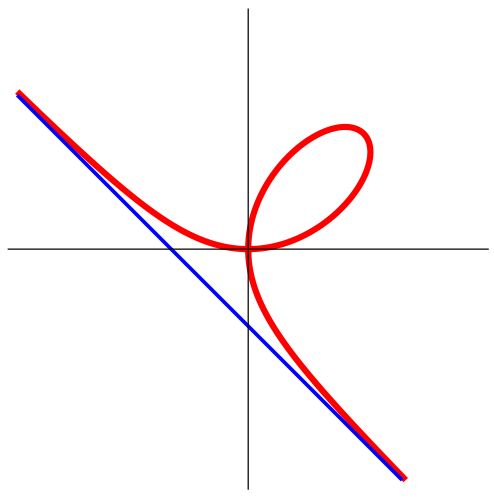

The asymptote to Descartes's folium $x^3 + y^3 - 3xy = 0$.

Joseph Alfred Serret's proof consists of the following. Let $F(x, y) = 0$ be the equation of the curve in a coordinate system so that $y = 0$ is the asymptote. Let us change $x$ in $1/x$. This produces a second algebraic curve $F_1(x, y) = 0$. Now if the original curve had a single branch asymptotic to $y = 0$, then the algebraic curve $F_1(x, y) = 0$ would have a "point d'arrêt" at $(0, 0)$, i.e. a *dead end*, which is impossible. Amazingly, Joseph Alfred takes for granted that such a stopping point is impossible. Clearly, this is not a proof in any sense of the term.

*"ce qui ne peut arriver pour une courbe algébrique".*

Ironically, he criticizes Euler for his lack of rigor. At the end of his short paper, he indeed quotes Euler's *Introductio in analysin infinitorum* (volume 2, chapter 7, section 174):

> **Quam ob** rem Linea curva duos habebit ramos in infinitum excurrentes inter se oppositos...

*"For this reason, the curve has two branches at infinity which will be opposite to each other..."*

The last sentence of Joseph Alfred's papers is: "This *quam ob rem* needed a proof". Did he really believe that Euler, or Newton, could not have thought of the change of variables $x \mapsto 1/x$?

Eighteen years later, Paul Joseph wrote another short paper[51], in the same journal, criticizing the earlier paper of his homonymous and prestigious colleague. He begins by asserting that Joseph Alfred's contribution is a reduction of the problem of asymptotes to the problem of stopping points of algebraic curves, but that this was "a priori obvious". Now — Paul Joseph insists — the main question remains open: one still has to prove that an algebraic curve cannot have a stopping point. He finally proposes the following proof.

[51] P. J. Serret. Note sur les courbes algébriques. *Nouvelles annales de mathématiques, journal des candidats aux écoles polytechnique et normale*, 4:311–313, 1865.

Let $(0, 0)$ be a point on an algebraic curve $F(x, y) = 0$. Let us intersect the curve with a small circle centered at the origin $x^2 + y^2 = r^2$. We use the well-known parameterization.

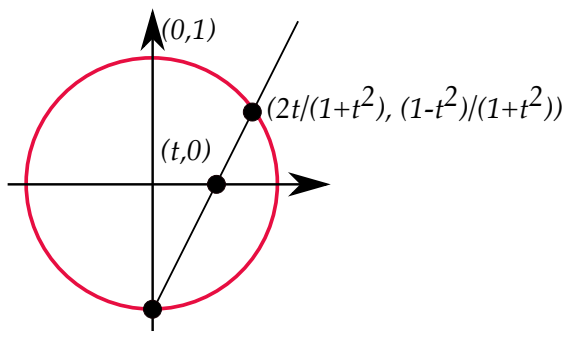

$$x = \frac{2rt}{1+t^2} \quad ; \quad y = \frac{(1-t^2)r}{1+t^2}.$$

By substitution in $F(x,y) = 0$ and multiplication by $(1+t^2)^d$ where $d$ is the degree of $F$ , we get an equation $\phi_{2d}(t) = 0$, where $\phi_{2d}$ is a polynomial of degree $2d$. Now, if the point $(0,0)$ happened to be a stopping point, the curve would intersect a small circle in a single point, so that an equation in $t$ of even degree $2d$ would have a single root, *"ce qui serait absurde"*.

By degree of $F$, I mean the maximum $i+j$ when $x^i y^j$ varies among the monomials with a non-trivial coefficient. We will check later that Paul Joseph was right and that the degree is indeed exactly $2d$ and not just $\leq 2d$.

Amazing. How could Paul Joseph not know that $t^2$ is of degree two and has a single root? This root is double but this is exactly our problem. We could imagine an algebraic curve going to some point and going back following the *same* path.

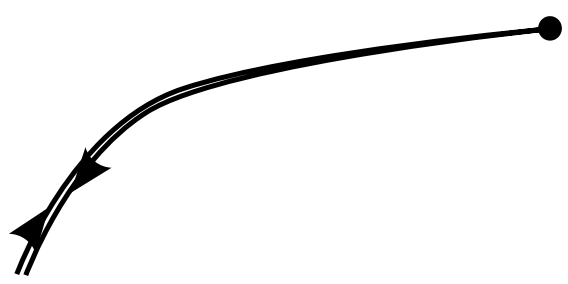

There is something to be proved.

## *Two important facts in commutative algebra*

I collect here two basic theorems on polynomials which will enable us to fix Paul Joseph's proof. See for example[52], [53] or [54] for much more about algebra. All rings will be assumed to be commutative. Some useful definitions are in the margins.

[52] S. Lang. *Algebra*, volume 211 of *Graduate Texts in Mathematics*. Springer-Verlag, New York, third edition, 2002.

[53] M. Artin. *Algebra*. Prentice Hall, Inc., Englewood Cliffs, NJ, 1991.

[54] N. Bourbaki. *Commutative algebra. Chapters 1–7*. Elements of Mathematics (Berlin). Springer-Verlag, Berlin, 1998.

**Theorem.** *Let $\mathcal{R}$ be a unique factorization domain. Then the polynomial ring $\mathcal{R}[x]$ is also a unique factorization domain.*

Say that a polynomial in $\mathcal{R}[x]$ is *primitive* if its coefficients are relatively prime. The key point is the so-called... Gauss's lemma.

Proved by Gauss in Article 42 of his *Disquisitiones Arithmeticae* in 1801, three years after his PhD.

**Lemma.** *The product of two primitive polynomials in $\mathcal{R}[x]$ is primitive.*

The (modern) proof is easy (but somehow indirect). If $p$ is prime in $\mathcal{R}$, the ring $\mathcal{R}/p$ is an integral domain. If $P_1(x)$ and $P_2(x)$ are two polynomials whose product is not primitive, then all coefficients of $P_1P_2$ are divisible by some prime $p$. We can reduce all coefficients modulo $p$ and get the following equality in $(\mathcal{R}/p)[x]$:

$$\overline{P}_1(x)\overline{P}_2(x) = 0.$$

Since the polynomial ring over an integral domain is an integral domain, we conclude that $\overline{P}_1(x)$ or $\overline{P}_2(x)$ is zero in $(\mathcal{R}/p)[x]$. ⊡

An *integral domain* is a ring in which the product of two nonzero elements is nonzero. A *unit* in a ring is an element which admits an inverse. Two elements $a,b$ in a ring are called *associated*, denoted $a \equiv b$, if there is a unit $u$ such that $b = ua$.

Define the *content* $cont(P)$ of a polynomial $P(x) \in \mathcal{R}[x]$ as the greatest common divisor of its coefficients. Clearly, every polynomial $P(x)$ can be written as the product $cont(P)\tilde{P}(x)$ where $\tilde{P}(x)$ is primitive. Gauss's lemma simply means that the content of a product is the product of the contents.

We can now prove the theorem. We are going to show that prime elements in $\mathcal{R}[x]$ are:

1. Prime elements of $\mathcal{R}$, seen as constant polynomials,
2. Primitive polynomials in $\mathcal{R}[x]$ which are prime when seen as polynomials over the quotient field $Quot(\mathcal{R})$ of $\mathcal{R}$.

The ring of polynomials over a field is Euclidean. This applies to $Quot(\mathcal{R})[x]$ so that any element $P(x)$ of $\mathcal{R}[x]$ can be written as a product of prime polynomials in $Quot(\mathcal{R})[x]$. Chasing denominators, we can write $P$ as a product of elements of the types 1/ and 2/:

$$P(x) = u \cdot r_1 \cdots r_k \cdot P_1(x) \cdots P_l(x).$$

Here, $u$ is a unit in $\mathcal{R}$, the $r_i$ are primes in $\mathcal{R}$ and the $P_i$'s are primitive and irreducible in $Quot(\mathcal{R})[x]$. By Gauss, the product $r_1 \cdots r_k$ is the content of $P$ and is therefore uniquely defined by $P$.

Since $\mathcal{R}$ is a unique factorization domain, the $r_i$ are uniquely defined by $P$ (up to units and up to permutation).

Since $Quot(\mathcal{R})[x]$ is also a unique factorization domain, the factors $P_i[x]$ are also uniquely defined up to permutation and units, *in* $Quot(\mathcal{R})[x]$. Now, an equality $Q(x) = aP(x)$ where $P(x), Q(x)$ are primitive in $\mathcal{R}[x]$ and $a$ is in $Quot(\mathcal{R})$ implies that $a$ is a unit in $\mathcal{R}$. ⊡

The immediate corollary is that *for any field K, the polynomial rings* $K[x_1, \ldots, x_n]$ *are unique factorization domains*. In this special case, the theorem means that any non-constant polynomial in $K[x_1, \ldots, x_n]$ can be written as a product of irreducible factors, in a unique way, up to permutation and multiplication by constant factors (in $K$).

The second algebraic result concerns the *resultant*. Let $P_1(x), P_2(x)$ denote two polynomials in the polynomial ring $\mathcal{R}[x]$ over some integral domain $\mathcal{R}$, of degrees $d_1, d_2 \geq 1$. Denote by $\mathcal{R}_d[x]$ the

A *prime element* $p$ in an integral domain $\mathcal{R}$ is an element such that the quotient ring $\mathcal{R}/p$ is an integral domain.
An element $a$ in an integral domain is *irreducible* if it is not the product of two non-units. Prime elements are irreducible. The converse does not hold in general.
A *unique factorization domain* (sometimes called factorial ring) is a ring in which every element is a product of prime elements, unique up to the ordering and units. Euclidean and principal rings — for instance the ring of polynomials over a field — are unique factorization domains. In this case, the concepts of primes and irreducible coincide and greatest common divisors are well defined.

$\mathcal{R}$-module of polynomials of degrees at most $d$, isomorphic to $\mathcal{R}^{d+1}$. Consider the map

$$\Phi : (A_1, A_2) \in \mathcal{R}_{d_2-1}[x] \times \mathcal{R}[x]_{d_1-1} \mapsto A_1P_1 - A_2P_2 \in \mathcal{R}_{d_1+d_2-1}[x].$$

This can be seen as a linear map from $\mathcal{R}^{d_1+d_2}$ into itself. Its determinant is called the *resultant* of $P_1$ and $P_2$, and denoted $Res(P_1, P_2)$. This element of $\mathcal{R}$ is a universal polynomial expression, with coefficients in $\mathbb{Z}$, in the coefficients of $P_1$ and $P_2$.

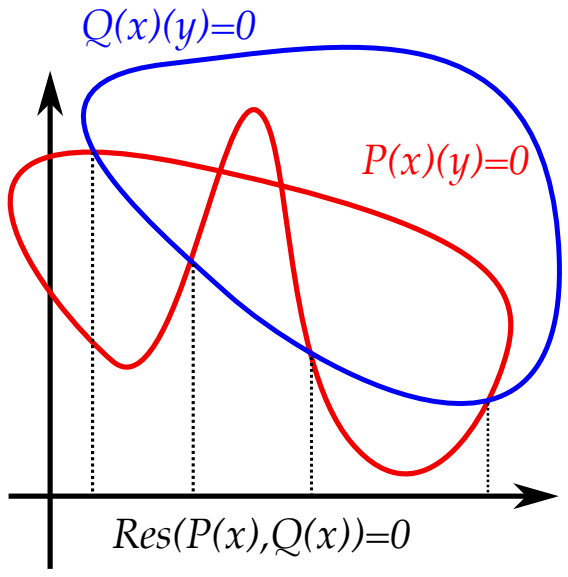

The projection of the intersection of two curves $P(x,y) = 0$ and $Q(x,y) = 0$ on the $x$ axis is given by the zeros of the resultant.

**Theorem.** *Suppose $\mathcal{R}$ is a unique factorization domain. The resultant $Res(P_1, P_2)$ is equal to zero if and only if $P_1$ and $P_2$ have a common non-trivial divisor in $\mathcal{R}[x]$.*

Indeed, if $P_1 = QQ_1$ and $P_2 = QQ_2$, the element $(Q_2, Q_1)$ is in the kernel of $\Phi$ so that the resultant vanishes.

Conversely, if the resultant vanishes, the kernel of $\Phi$ is not trivial so that there are non trivial elements $A_1, A_2$ in $\mathcal{R}_{d_1-1}[x]$ and $\mathcal{R}_{d_2-1}[x]$ such that $A_1P_1 = A_2P_2$. The conclusion follows from the fact that $\mathcal{R}[x]$ is a unique factorization domain: if $P_1$ and $P_2$ were relatively prime, $P_1$ would divide $A_2$ which is impossible since the degree of $A_2$ is less than the degree of $P_1$. ⊡

## *Proof of Gauss's claim*

We can now prove that the neighborhood of a point on a real algebraic curve consists in an *even* number of *arcs only intersecting at the origin.*

Let $F(x,y) = 0$ be the equation of our real algebraic curve passing through the origin $(0,0)$. Write $F$ as a product of irreducible factors:

$$F(x,y) = F_1(x,y) \cdots F_n(x,y).$$

The polynomials $F_i, F_j$ are non-associated when $F_i \not\equiv F_j$ for $i \neq j$: there is no constant $c$ such that $F_j = cF_i$.

Without changing the zero locus of $F$ in the neighborhood of $(0,0)$, we can delete some of the factors and assume that all the $F_i$'s vanish at $(0,0)$ and are non-associated irreducible factors.

The zero locus of $F$ in the neighborhood of $(0,0)$ is the union of the zero loci of the $F_i$'s.

In order to prove Gauss's claim, we prove two lemmas.

**Lemma** (1). *Let $P(x,y) \in \mathbb{R}[x,y]$ be an irreducible polynomial and $Q(x,y) \in \mathbb{R}[x,y]$ some polynomial. Suppose that the curves $P(x,y) = 0$ and $Q(x,y) = 0$ have an infinite number of intersection points in some small neighborhood of the origin. Then $P$ divides $Q$ in $\mathbb{R}[x,y]$. In particular, if $P$ and $Q$ are both irreducible and not associated, then the two corresponding curves can only intersect in isolated points.*

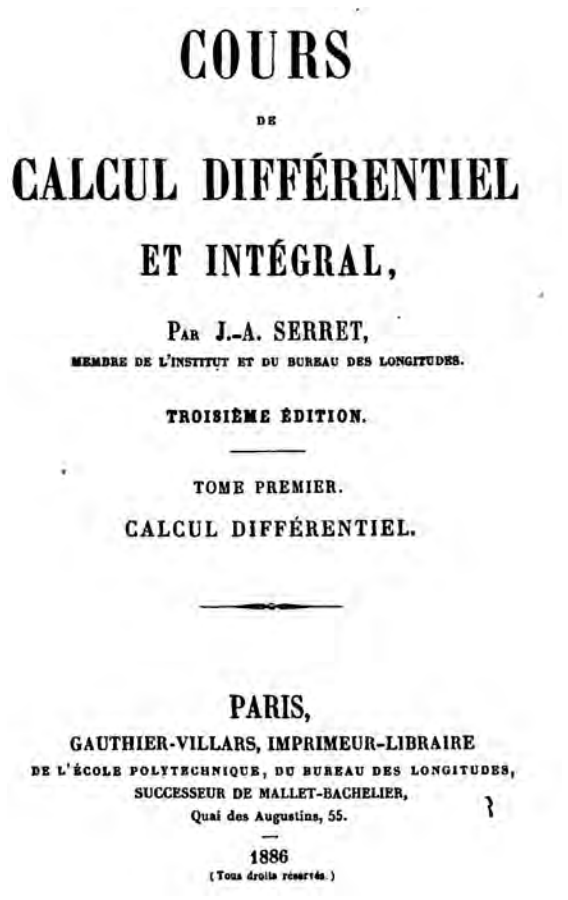
COURS
DE
CALCUL DIFFÉRENTIEL
ET INTÉGRAL,
PAR J.-A. SERRET,
MEMBRE DE L'INSTITUT ET DU BUREAU DES LONGITUDES.
TROISIÈME ÉDITION.
TOME PREMIER.
CALCUL DIFFÉRENTIEL.
PARIS,
GAUTHIER-VILLARS, IMPRIMEUR-LIBRAIRE
DE L'ÉCOLE POLYTECHNIQUE, DU BUREAU DES LONGITUDES,
SUCCESSEUR DE MALLET-BACHELIER,
Quai des Augustins, 55.
1886
(Tous droits réservés.)

**Lemma** (2). *If $P(x,y)$ is irreducible, its zero locus in the neighborhood of the origin consists of an even number of arcs converging to $(0,0)$.*

Let us begin by the first lemma. If $P(x,y) = 0$ contains an infinite number of points on the same vertical axis $x = x_0$, the polynomial $P(x,y)$ must be divisible by $(x - x_0)$ and since we assume that it is irreducible and vanishes at the origin, this implies that $P(x,y)$ is a constant multiple of $x$, for which the lemma is obvious. Without loss of generality, we can therefore assume that $P(x,y) = 0$ intersects every vertical line in a finite number of points.

If $P(x_0,y_0) = Q(x_0,y_0) = 0$, the two polynomials $P(x_0,y), Q(x_0,y)$, seen as elements of $\mathbb{R}[y]$, have a common root $y_0$ and therefore their resultant vanishes, as an element of $\mathbb{R}$.

Assume that the curves $P(x,y) = 0$ and $Q(x,y) = 0$ have an infinite number of intersection points in some small neighborhood of the origin. Let us look at the resultant of $P, Q \in \mathbb{R}[x][y]$ as an element of $\mathbb{R}[x]$. This resultant vanishes for an infinite number of values $x_0$ and therefore vanishes identically. We have seen that this implies that $P, Q$ have a common factor in $\mathbb{R}[x,y]$. Since $P$ is irreducible, this shows that $P$ divides $Q$. ⊡

We now prove the second lemma following Paul Joseph's idea. Let us set

$$F(x,y) = \sum_{i,j} a_{ij} x^i y^j.$$

Observe that when $F(x,y) = x$, the lemma is obvious.

Denote by $d$ the degree of $F$ (which is by definition the maximum value of $i + j$ for which $a_{ij} \neq 0$). Note that since $F$ is irreducible it is not divisible by $x$ (unless it is a constant multiple of $x$) so that one of the coefficients $a_{0j}$ is not 0. Fixing $r$, we get a parameterization by $t$ of the circle of radius $r$ (minus $(0,-r)$):

$$x = r\frac{2t}{1+t^2} \quad ; \quad y = r\frac{1-t^2}{1+t^2}.$$

Substitute in $F(x,y)$ and multiply the result by $(1+t^2)^d$:

$$\phi_{2d,r}(t) = \sum_{i,j} a_{ij} r^{i+j} (2t)^i (1-t^2)^j (1+t^2)^{d-i-j}.$$

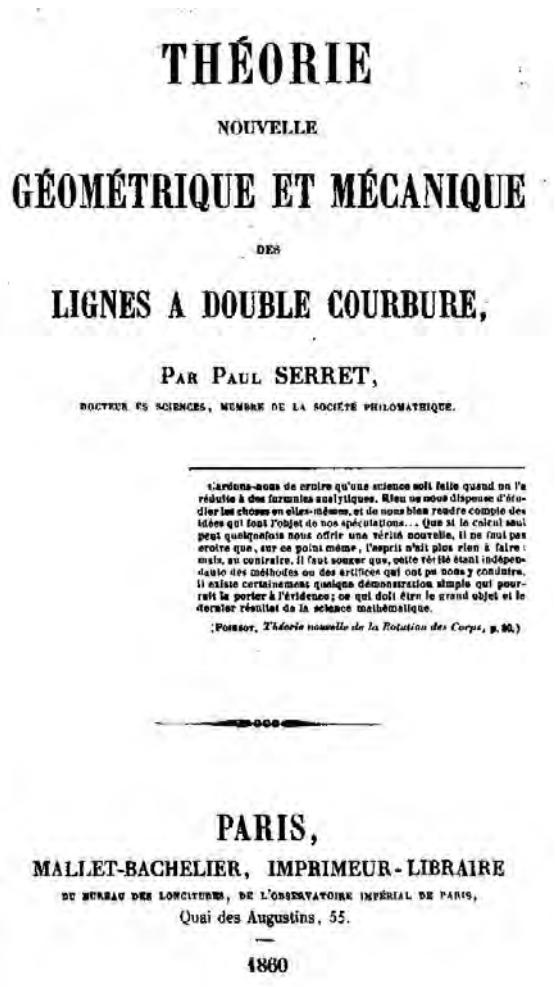

THÉORIE
NOUVELLE
GÉOMÉTRIQUE ET MÉCANIQUE
DES
LIGNES A DOUBLE COURBURE,
PAR PAUL SERRET,
DOCTEUR ÈS SCIENCES, MEMBRE DE LA SOCIÉTÉ PHILOMATHIQUE.

Gardons-nous de croire qu'une science soit faite quand on l'a réduite à des formules analytiques. Rien ne nous dispense d'étudier les choses en elles-mêmes, et de nous bien rendre compte des idées qui font l'objet de nos spéculations... Que si le calcul seul peut quelquefois nous offrir une vérité nouvelle, il ne faut pas croire que, sur ce point même, l'esprit n'ait plus rien à faire : mais, au contraire, il faut songer que, cette vérité étant indépendante des méthodes ou des artifices qui ont pu nous y conduire, il existe certainement quelque démonstration simple qui pourrait la porter à l'évidence ; ce qui doit être le grand objet et le dernier résultat de la science mathématique.

(POINSOT, *Théorie nouvelle de la Rotation des Corps*, p. 80.)

PARIS,
MALLET-BACHELIER, IMPRIMEUR-LIBRAIRE
DU BUREAU DES LONGITUDES, DE L'OBSERVATOIRE IMPÉRIAL DE PARIS,
Quai des Augustins, 55.
1860

This is a polynomial in $t$ whose highest monomial is

$$\left(\sum_j (-1)^j a_{0j} r^j\right) t^{2d}$$

which is certainly not zero for small $r \neq 0$. Hence, Paul Joseph is right and the degree of $\phi_{2d,r}(t)$ is equal to $2d$.

In order to complete the proof, we still have to show that the roots of $\phi_{2d,r}(t) = 0$ are simple for small $r \neq 0$, so that there is an even number of roots. At a double root $t_0$, the polynomial $\phi_{2d,r}(t)$ and its derivative vanish simultaneously. Geometrically, this means that the tangent to the circle at this point is also tangent to the curve $F(x,y) = 0$. Said differently, the double points that we want to exclude correspond to the intersection of $F(x,y) = 0$ and the curve $y\partial F/\partial x - x\partial F/\partial y = 0$. Since we assume that $F$ is irreducible, the first lemma implies that these two curves intersect in a finite number of points, unless $F$ divides $y\partial F/\partial x - x\partial F/\partial y$ . For degree reason, this can only happen if $y\partial F/\partial x - x\partial F/\partial y$ is a constant multiple of $F$, which means in turn that $F$ is constant on circles This implies that $F$ is a polynomial in $x^2+y^2$ and since it vanishes at the origin, it should be divisible by $x^2+y^2$. It is therefore a constant multiple of $x^2+y^2$, whose zero locus reduces to the origin. ⊡

The proof of Gauss's claim is essentially finished. The restriction of $F$ to each circle $x^2+y^2 = r^2$ for small nonzero $r$ has an even number of zeroes which are simple. Using the implicit function theorem in this very elementary situation, we conclude that these zeroes define an even number of disjoint curves converging to the origin. This does not say anything about the limiting directions of these curves: they might a priori converge to the origin without having a limiting tangent. ⊡

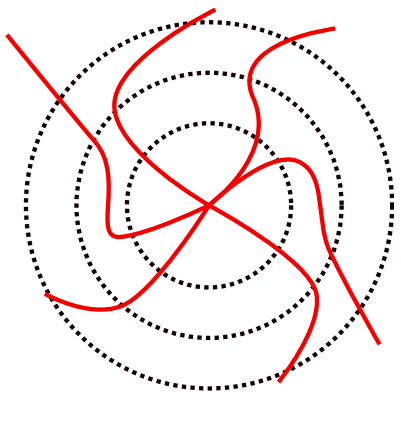

We will see that any branch of an algebraic curve does have a tangent.

Carl Friedrich Gauss and Paul Joseph Serret were right.

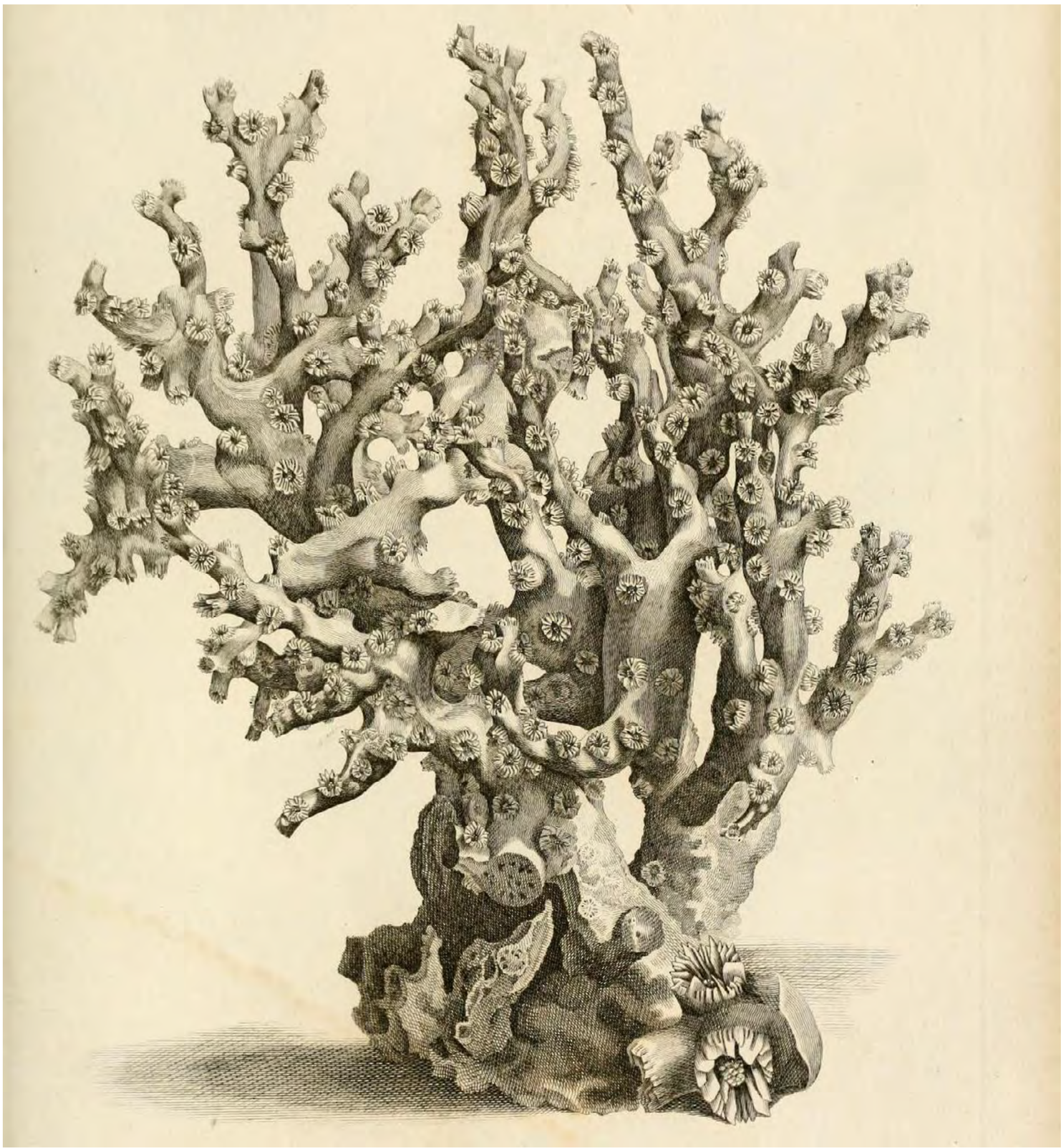

A plate from a book by J. Lamouroux, dated 1821, representing *Oculina Hirtella*. This book was in the library of *HMS Beagle*, which was also Darwin's cabin for five years. Ramis sparsis divergentibus! ©

# *De seriebus divergentibus: Euler, Cauchy and Poincaré*

## *Euler's seriebus divergentibus*

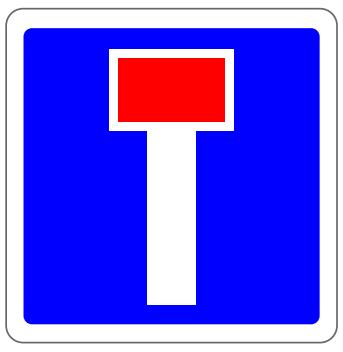

Dead end. This chapter is completely independent from the rest of the book. ©

NEWTON DID NOT LIMIT THE USE OF INFINITE SERIES TO EQUATIONS OF THE FORM $F(x, y) = 0$. He also used them in a systematic way to solve *differential equations*. His approach was essentially practical. He looks for a solution as a formal series and computes inductively a large number of terms of the series in order to get some "accuracy". There was no systematic understanding of the concept of convergence, but in all the cases that he treated the series were indeed convergent.

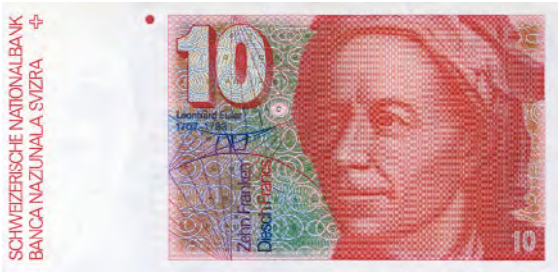

Leonhard Euler (1707-1783) ©

Later, Euler became the great Master of series. It is a common opinion among contemporary mathematicians that Euler was careless with series and that he manipulated series which "make no sense". For instance, his formula[55]

$$1 - 2^3 + 3^3 - 4^3 + 5^3 - 6^3 + etc \ldots = -\frac{1}{8}$$

is shocking for undergraduate students, who have been taught the definition of a convergent series very early and refuse to consider these horrors. Not so! Euler knows what he does. He discusses various procedures for attributing a sum to a series, even if it is divergent, and tries to compare these procedures. His series are not the most general: they are implicitly defined by some kind of algorithm, to use an anachronism. He is convinced that divergent series do represent "something" inherently linked with the nature of the series. His paper *De seriebus divergentibus*[56] is a pure gem and I recommend it to any mathematician.

[55] L. Euler. Remarques sur un beau rapport entre les séries des puissances tant directes que réciproques. *Mémoires de l'académie des sciences de Berlin*, 17:83–106, 1768. See the Euler Archive for English translations and comments.

[56] L. Euler. De seriebus divergentibus. *Novi Commentarii academiae scientiarum Petropolitanae*, 5:205–237, 1760. See the Euler Archive for English translations and comments.

One of his examples is famous:

$$S = 1 - 1! + 2! - 3! + 4! - 5! + \cdots.$$

Using five different methods of summation, Euler gets values which seem to indicate that $S$ should be close to 0.5963473621237. One of the most convincing methods uses the fact that the formal series

$$\hat{f}(x) = x - 1!x^2 + 2!x^3 - \cdots$$

is a solution of the *linear* differential equation

$$x^2 y' + y = x.$$

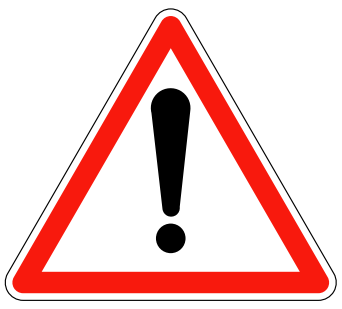

This is a very elementary equation for which one finds an explicit solution which is equal to 0 for $x = 0$:

$$f(x) = \exp\left(\frac{1}{x}\right) \int_0^x \frac{1}{t} \exp\left(-\frac{1}{t}\right) dt.$$

Be careful! I would certainly not like to give to my young reader the wrong feeling that *any* divergent series has a well defined sum. Some divergent series are indeed summable using different methods and produce *different* sums. As an easy example, think of $1 - 1 + 1 - 1 + \cdots$.

Somehow, we could say that $f(x)$ represents the value of the formal series $\hat{f}(x)$. The numerical value found by Euler 0.596347362123 is the value $f(1)$.

Just type the following in your computer

$$N[Exp[1] * Integrate[Exp[-1/t]/t, t, 0, 1], 100]$$

and get immediately

0.5963473623231940743410784993692793760741778601525487815734849104823272191148744174

not in complete agreement though with Euler's numerical result.

Then came the period of disgrâce for divergent series.

The new master was Augustin Cauchy who defined clearly the concept of *convergence* and who is usually associated with mathematical rigor. This is not completely wrong but this is without any doubt an exaggerated simplification. On the one hand, rigor did exist before Cauchy and on the other hand Cauchy did not reject divergent series[57]. Unfortunately, even today, many students are still convinced that divergent series come from the devil...

[57] A. Cauchy. Sur un emploi légitime des séries divergentes. *Comptes Rendus Académie des Sciences*, XVII:18–25, 1843.

In 1821, in the preface to his Cours d'analyse, Cauchy wrote that he was *forced* to abandon divergent series!

> I have been forced to admit some propositions which will seem, perhaps, hard to accept. For instance, that a divergent series has no sum.

*"J'ai été forcé d'admettre diverses propositions qui paraîtront peut-être un peu dures. Par exemple qu'une série divergente n'a pas de somme..."*

In a famous letter to Holmboe, Abel wrote in 1826:

> Divergent series are in general something fatal, and it is a disgrace to base any proof on them.

*"Les séries divergentes sont en général quelque chose de bien fatal et c'est une honte qu'on ose y fonder aucune démonstration."* volume 2 of Abel's collected papers.

## *Poincaré*

The next great master was Poincaré who clearly understood that divergent series are not only useful but also necessary to solve natural questions from celestial mechanics. I refrain from presenting these dynamical aspects, even though they are fascinating and connected to current research activity.

Let me quote from the second volume of the *Méthodes nouvelles de mécanique céleste*.

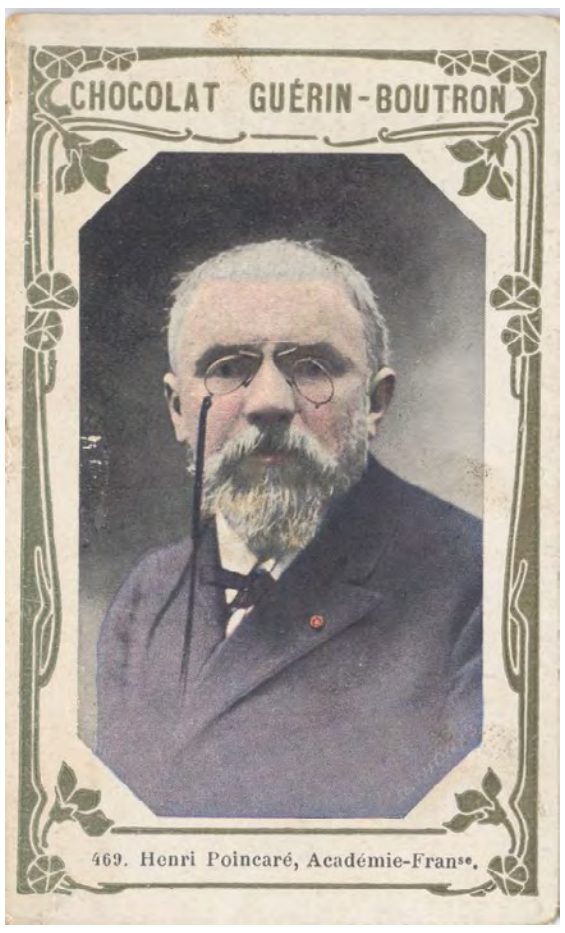

Henri Poincaré (1854-1912). How many mathematicians were so famous during their life time that their photograph was printed on chocolate bars?

> There is a kind of misunderstanding between the geometers and the astronomers, concerning the meaning of the word convergence. The geometers, concerned with absolute rigor and not bothered by the length of the inextricable computations that they conceive to be possible without trying to undertake them explicitly, would say that a series is convergent when the sum of the terms tends to a definite limit, even if the first terms decrease very slowly. On the contrary, the astronomers have the habit of saying that a series converges when, for instance, the first 20 terms decrease very rapidly, even if the remaining terms would grow forever. Thus, let us take a simple example and consider the two series which have as general term
>
> $$\frac{1000^n}{n!} \quad \text{and} \quad \frac{n!}{1000^n}.$$
>
> The geometers will say that the first series converges, and even that it converges fast [...]; and they will say that the second series diverges [...] On the contrary, the astronomers will consider the first series as divergent, [...] , and the second series as convergent. The two rules are legitimate: the first one in the theoretical researches; the second one in the numerical applications.

> Il y a entre les géomètres et les astronomes une sorte de malentendu au sujet de la signification du mot convergence. Les géomètres, préoccupés de la parfaite rigueur et souvent trop indifférents à la longueur de calculs inextricables dont ils conçoivent la possibilité, sans songer à les entreprendre effectivement, disent qu'une série est convergente quand la somme des termes tend vers une limite déterminée, quand même les premiers termes diminueraient très lentement. Les astronomes, au contraire, ont coutume de dire qu'une série converge quand les vingt premiers termes, par exemple, diminuent très rapidement, quand même les termes suivants devraient croître indéfiniment. Ainsi, pour prendre un exemple simple, considérons les deux séries qui ont pour terme général $\frac{1000^n}{n!}$ et $\frac{n!}{1000^n}$. Les géomètres diront que la première série converge, et même qu'elle converge rapidement, [...] mais ils regarderont la seconde comme divergente [...]. Les astronomes, au contraire, regardéront la première série comme divergente, [...] et la seconde comme convergente [...] Les deux règles sont légitimes : la première, dans les recherches théoriques ; la seconde, dans les applications numériques.

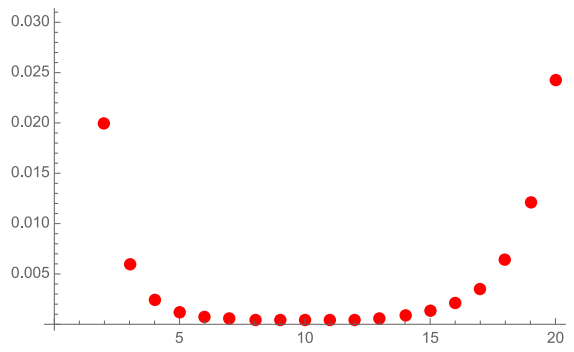

$\frac{n!}{10^n}$ for $n = 1, \ldots, 20$.

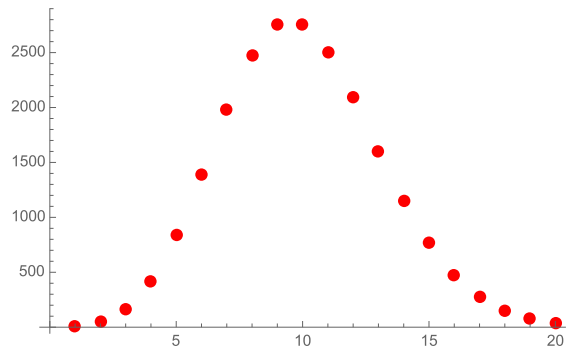

$\frac{10^n}{n!}$ for $n = 1, \ldots, 20$.

The example of Poincaré is perfect: look at the values of $\frac{n!}{1000^n}x^n$ and $\frac{1000^n}{n!}x^n$ for $x = 100$ and .01 respectively, and for $n = 1, \ldots, 20$.

## *The saddle-node and Euler's equation*

Let us look at a very simple example showing that there are no choices: in order to understand ordinary differential equations, even with polynomial coefficients, we *have* to deal with divergent series.

Consider the following simple system:

$$\frac{dx}{dt} = x^2 \quad ; \quad \frac{dy}{dt} = -y + x.$$

It is called a *saddle-node* because it looks like a saddle where $x > 0$ and a node when $x < 0$. This may look like a very degenerate situation but it appears in *codimension 1*: at the origin the linear part of the vector field has one vanishing eigenvalue. Therefore similar saddle-nodes should appear in generic one parameter families of vector fields in the plane.

The picture in the margin shows the phase portrait of this vector field. Clearly, we see a smooth invariant curve passing through the origin (and different from the $y$ axis). This is called the *central manifold*.

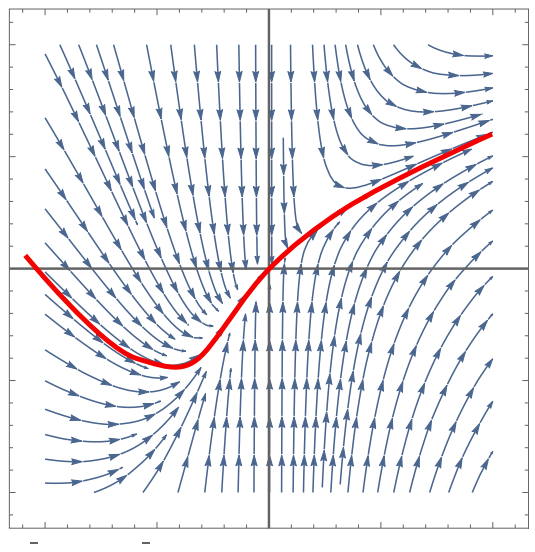

Phase portrait of the saddle-node.

Looking for this curve as a graph $y(x)$, we get immediately the Euler equation $x^2 \frac{dy}{dx} + y = x$. So the equation of the central manifold is the $C^\infty$ function defined by

$$y = f(x) = \exp\left(\frac{1}{x}\right) \int_0^x \frac{1}{t} \exp\left(-\frac{1}{t}\right) dt,$$

and we have to understand how this function is related to the formal divergent series $\hat{f}$.

## *Euler function, Stokes phenomenon etc.*

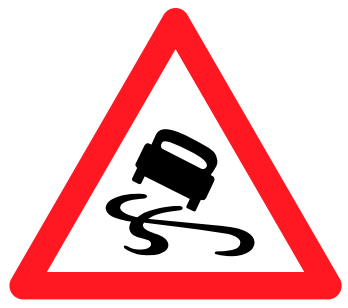

Caution: dangerous changes of variables! ©

I follow the presentation by Hardy[58].

Change variable and set $t = x/(1 + xw)$ so that

$$f(x) = \int_0^\infty \exp(-w) \frac{x}{1 + xw} \, dw.$$

[58] G. H. Hardy. *Divergent series*. Éditions Jacques Gabay, Sceaux, 1992. With a preface by J. E. Littlewood and a note by L. S. Bosanquet, Reprint of the revised (1963) edition.

This yields

$$\begin{aligned} f(x) &= \int_0^\infty \exp(-w) \left(x - x^2 w + x^3 w^2 - \cdots + (-1)^{n-1} x^n w^{n-1}\right) dw \\ &\quad + (-1)^n x^{n+1} \int_0^\infty \frac{\exp(-w) w^n}{1 + xw} \, dw \\ &= x - 1! x^2 + 2! x^3 - \cdots + (-1)^{n-1}(n-1)! x^n + R_n(x). \end{aligned}$$

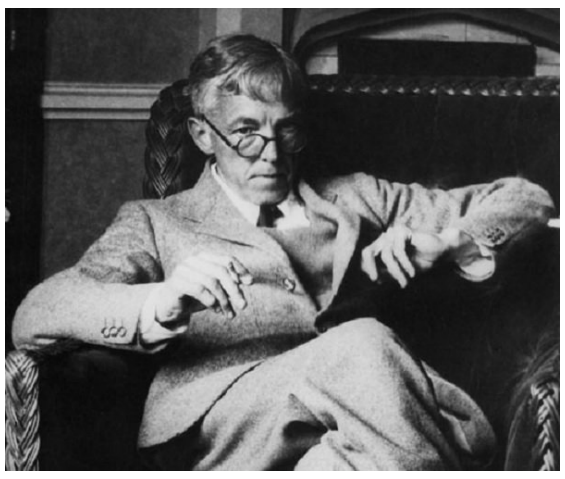

"Young men should prove theorems, old men should write books". (Hardy 1877-1947). ©

The term $R_n$ is easy to majorize. If $x, w > 0$, we have $1 + xw > 1$ so that

$$|f(x) - (x - 1!x^2 + 2!x^3 - \cdots + (-1)^{n-1}(n-1)!x^n)| \leq n! x^{n+1}.$$

In other words, *the formal series $\hat{f}$ is asymptotic to the $C^\infty$ function $f$.*

Recall that a series $\sum_k a_k x^k$ is *asymptotic* to a function $f(x)$ if for every $n$, we have $f(x) - \sum_{k=1}^n a_k x^k = o(x^n)$.

Actually, much more can be said. Suppose now that $x$ is a complex number which is not a negative real number. Then the formula defining $f$ makes perfect sense, so that *$f$ is a holomorphic function in $\mathbb{C} \setminus \mathbb{R}_-$.* Suppose now that $x$ belongs to a sector where its argument is in $[-\pi + \delta, \pi - \delta]$, for some $\delta > 0$. In this sector,

$|1+xw|$ is minorized ($w$ is still a positive real number) so that we get some inequality

$$|f(x)-(x-1!x^2+2!x^3-\cdots+(-1)^{n-1}(n-1)!x^n)| \leq C(\delta)n!|x|^{n+1}$$

for $x$ in this sector. Said differently, *the formal series $\hat{f}$ is asymptotic to the holomorphic function $f$ in any sector not containing the negative real line.*

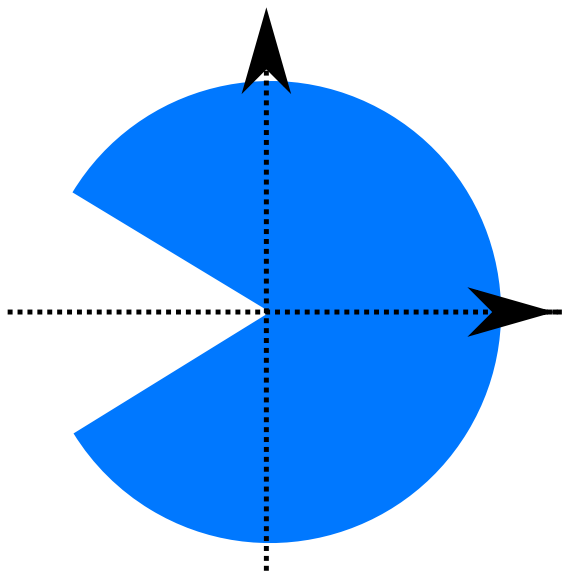

Pacman.

Still more can be said. Let us continue with the presentation by Hardy of Euler's manipulations.

$$\begin{aligned} f(x) &= \exp(\tfrac{1}{x})\int_0^x \frac{\exp(-\frac{1}{t})}{t}\,dt = \exp(\tfrac{1}{x})\int_{\frac{1}{x}}^{\infty}\frac{\exp(-u)}{u}\,du \\ &= -\exp(\tfrac{1}{x})li(\exp(-\tfrac{1}{x})) \end{aligned}$$

where $li$ is the *integral logarithm* defined for $0<v<1$ by

$$li(v) = \int_0^v \frac{dv}{\log v} = -\int_{\log\frac{1}{v}}^{\infty}\frac{\exp(-u)}{u}\,du.$$

$$\begin{aligned} -li(\exp(-y)) &= \int_y^{\infty}\frac{\exp(-u)}{u}\,du \\ &= \int_1^{\infty}\frac{\exp(-u)}{u}\,du - \int_0^1 \frac{1-\exp(-u)}{u}\,du \\ &\quad -\int_1^y \frac{du}{u} + \int_0^y \frac{1-\exp(-u)}{u}\,du \\ &= -\gamma - \log y + y - \tfrac{1}{2\cdot 2!}y^2 + \tfrac{1}{3\cdot 3!}y^3 - \cdots \end{aligned}$$

where $\gamma$ is... the Euler constant. It follows that

$$f(x) = \exp(\frac{1}{x})\log x + S(\frac{1}{x})$$

where

$$S(y) = -\exp(y)\left(\gamma - y + \frac{1}{2\cdot 2!}y^2 - \frac{1}{3\cdot 3!}y^3 + \cdots\right).$$

Note that $S(y)$ is an entire function, i.e. holomorphic and uniform in the full complex plane.

This provides some holomorphic extension of $f$ on the universal cover of $\mathbb{C}\setminus\{0\}$. As we go one turn around the origin, the function changes by $2i\pi\exp(\frac{1}{x})$.

Note that $li(x) = \int_0^{\infty}\frac{dx}{\log x}$ is also defined for $x>1$ as an improper integral:

$$\lim_{\varepsilon\to 0}\left(\int_0^{1-\varepsilon}+\int_{1+\varepsilon}^{\infty}\right)\frac{dx}{\log x}.$$

It is famous in number theory as it gives a very accurate estimate for the number $\pi(x)$ of prime numbers $\leq x$ when $x$ tends to infinity. In particular

$$|\pi(x)-li(x)| = o\left(\frac{x}{(\log x)^N}\right)$$

for every $N\geq 1$.

Since this chapter is about divergent series, it might be a good idea to mention the divergent asymptotic expansion $li(x) = \frac{x}{\log x}\left(\sum_{k=0}^{\infty}\frac{k!}{(\log x)^k}\right)$, that you should be able to prove yourself, integrating by parts.

Let us sum up the properties of $f$.

- *It is a multivalued holomorphic function which is defined in the whole plane, or more precisely a holomorphic function on the Riemann surface of the logarithm.*
- *In any sector of angle* $< 2\pi$*, the function* $f$ *is asymptotic to the formal series* $\hat{f}(x) = x - 1!x^2 + 2!x^3 - \cdots$.
- *The monodromy, that is to say the change in the value of* $f(x)$ *as* $x$ *goes around the origin, is* $2i\pi\exp(\frac{1}{x})$*, which is flat in any sector of angle* $< 2\pi$*. This means that any two single-valued determinations of* $f$ *in a sector have the same asymptotic expansion* $\hat{f}$.

A "holomorphic function on the Riemann surface of the logarithm" is an old fashioned way of speaking of $\phi(\log z)$ where $\phi$ is a holomorphic function defined on the complex plane. This is multivalued since $\log z$ is defined "up to $2i\pi$".

The divergence of the formal series $\hat{f}$ corresponds to the "multivalued up to a flat function" property of the function $f$. This is not the first time that a phenomenon in the real domain is explained by another one in the complex domain. This is called the *Stokes phenomenon* (discovered March 19th, 1857 at 3 a.m.).

Of course Euler's example is just an example. The remarkable fact is that this example is significant and that a beautiful theory has been developed. Analytic or even algebraic differential equations may have solutions which are divergent series, but one can give a perfectly well-defined meaning to their sum, as holomorphic *multivalued* functions.

I refrain from continuing in that direction since our promenade would not go where I plan to go. Even in promenades it is good to sail towards some kind of heading.

For a fascinating description of the historical development of the theory, I strongly recommend[59]. For a more systematic description, at an accessible level, these lecture notes[60] will be useful.

This paper[61] is a modern presentation of divergent series published for the 300th anniversary of Leonhard Euler's birth.

[59] J.-P. Ramis. Poincaré et les développements asymptotiques (première partie). *Gaz. Math.*, 133:33–72, 2012.; and J.-P. Ramis. Les développements asymptotiques après Poincaré: continuité et divergences. *Gaz. Math.*, 134:17–36, 2012.

[60] M. Loday-Richaud. Divergent series and differential equations. Part of a book. Submitted. 2014, hal-01011050.

[61] V. S. Varadarajan. Euler and his work on infinite series. *Bull. Amer. Math. Soc. (N.S.)*, 44(4):515–539, 2007.

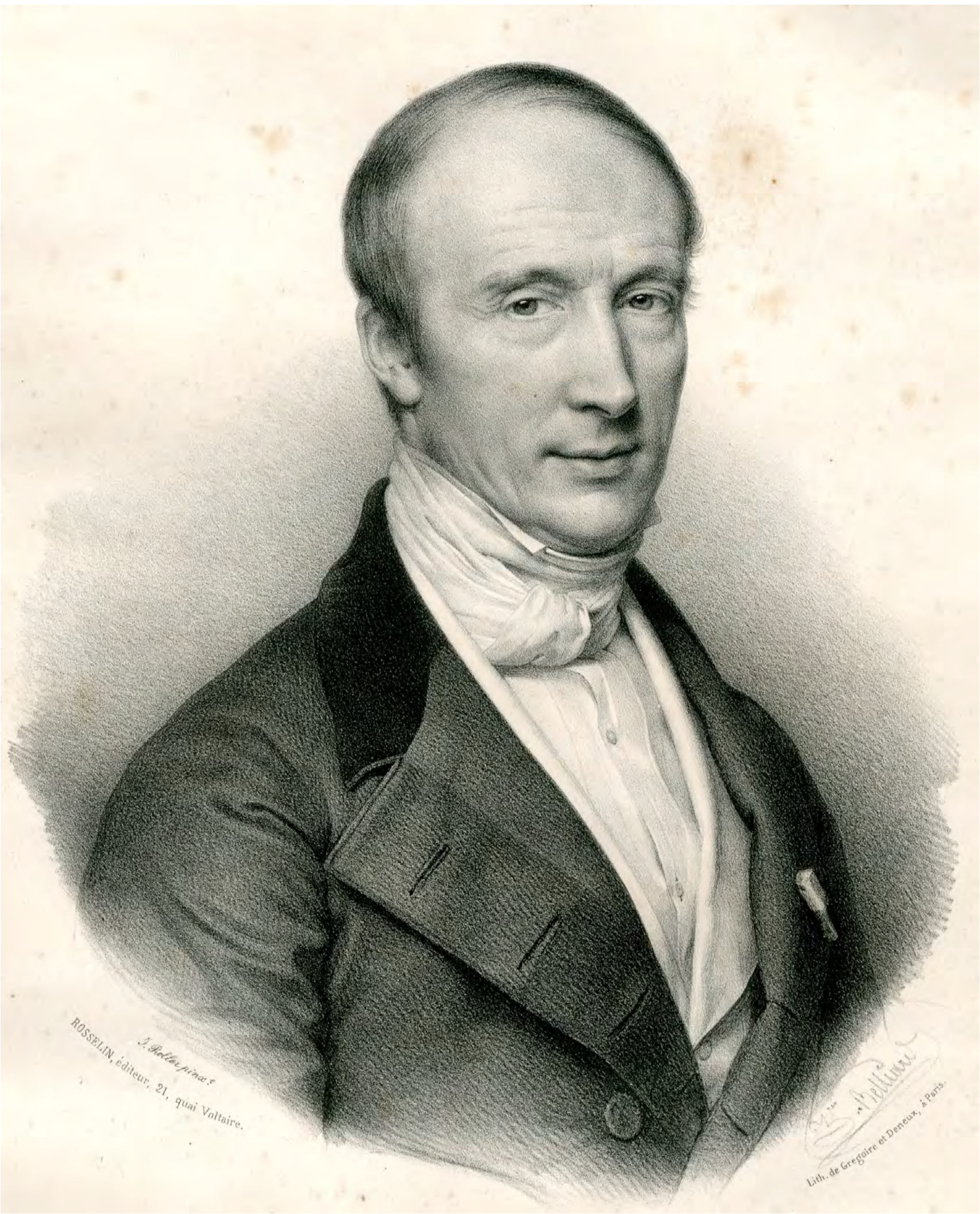

Augustin Cauchy.
(1789 –1857). ©

# *Convergence*
# *Le calcul des limites de Cauchy*

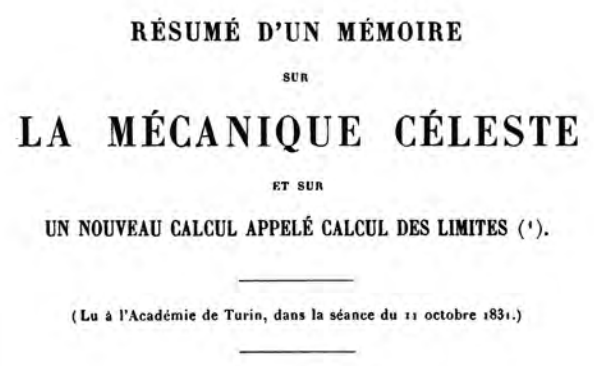
RÉSUMÉ D'UN MÉMOIRE

SUR

LA MÉCANIQUE CÉLESTE

ET SUR

UN NOUVEAU CALCUL APPELÉ CALCUL DES LIMITES (¹).

(Lu à l'Académie de Turin, dans la séance du 11 octobre 1831.)

I PROVE NOW THE SO-CALLED PUISEUX THEOREM giving a local parameterization of a complex algebraic curve in a neighborhood of a singular point, in terms of *convergent* power series. We will not follow Puiseux's original approach. Instead, I prefer to use a method introduced by Cauchy under the name *Calcul des limites*[62].

[62] A. Cauchy. *Œuvres complètes*, volume 12. Gauthier Villars, 1882.

The word *limite* should not be understood as limit but as *bound*. The present terminology is "majorant method".

## *The implicit function theorem*

I begin with a proof *à la Cauchy* of the classical implicit function theorem. This used to be the standard proof in old textbooks but is frequently ignored today and replaced by more powerful methods, based on fixed point theorems. It has nevertheless some advantages: it is very elementary and almost entirely combinatorial. I recommend this book[63] for an interesting historical approach.

[63] S. G. Krantz and H. R. Parks. *The implicit function theorem*. Birkhäuser Boston, Inc., Boston, MA, 2002. History, theory, and applications.

Let us denote by $K$ a field of characteristic 0, equipped with a norm, which is a map $x \in K \mapsto |x| \in \mathbb{R}_+$ such that $|xy| = |x||y|$ and $|x+y| \leq |x| + |y|$. Assume that $|x| = 0$ if and only if $x = 0$ and that $K$ equipped with $|\,|$ is complete: Cauchy sequences converge.

I assume moreover that the norm is not the trivial one for which $|x| = 1$ for all nonzero elements $x$.

I basically have in mind the case of $\mathbb{C}$ and $\mathbb{R}$ but there are many other examples ($p$-adic fields in particular).

Let us denote by $K\{x\}$ the ring of series

$$f(x) = \sum_{k \geq 0} u_k x^k$$

for which the $u_k$ are in $K$ and satisfy some inequality of the form

$$|u_k| \leq Cr^k$$

for some $C, r > 0$. Since $K$ is complete, this corresponds to series which are absolutely convergent in some neighborhood of 0 (germs of *analytic functions* for $\mathbb{C}, \mathbb{R}$). For simplicity, we will say that the elements of $K\{x\}$ are *convergent series*.

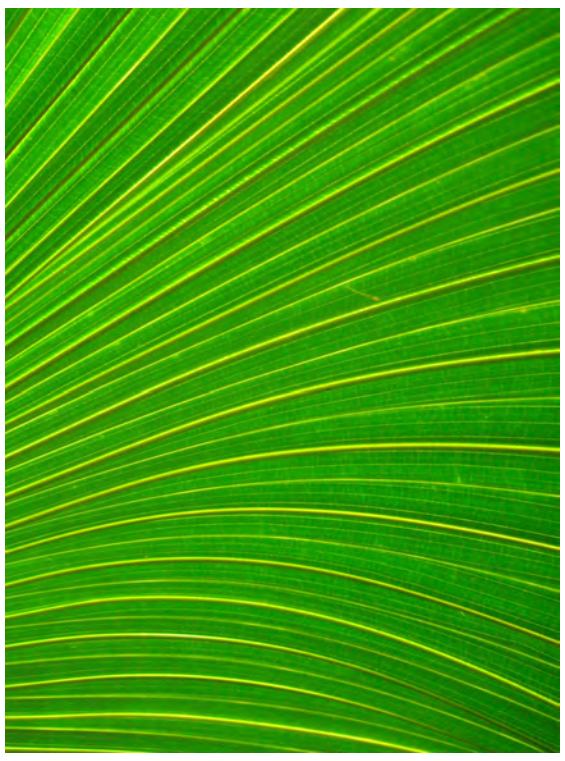

Level curves of a regular function $F(x, y)$.

Similarly, denote by $K\{x, y\}$ the ring of series

$$F(x, y) = \sum_{i,j \geq 0} a_{ij} x^i y^j$$

which are *convergent*, i.e. for which the $a_{ij}$'s are in $K$ and such there exist $C, r > 0$ such that for all $i, j$:

$$|a_{ij}| \leq Cr^{i+j}.$$

**Theorem** (Implicit function theorem). *Let $F \in K\{x, y\}$ be such that $F(0,0) = 0$ and $\partial F/\partial y(0,0) \neq 0$. Then there is a convergent series $f(x) \in K\{x\}$ such that $f(0) = 0$ and $F(x, f(x)) = 0$. The solutions $(x, y)$ to the equation $F(x, y) = 0$ in the neighborhood of $(0,0)$ in $K^2$ are precisely the pairs $(x, f(x))$.*

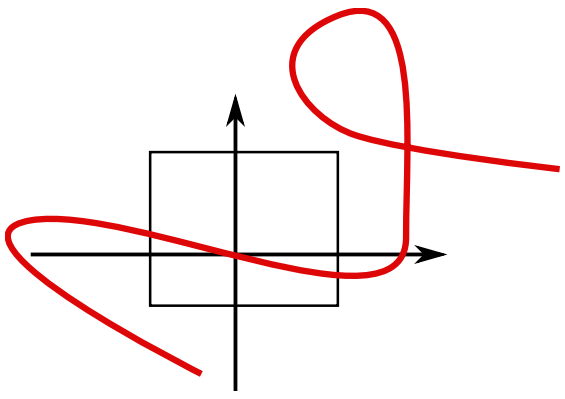

The proof is the following. If we substitute a formal series $y = \sum_{k \geq 1} u_k x^k$ in the formal series $\sum_{i,j \geq 0} a_{ij} x^i y^j$ (with $a_{00} = 0$), the result is a formal series $\sum_{l \geq 1} v_l x^l$ whose coefficients depend on the $u_k$'s and the $a_{ij}$'s. Let us compute the first terms in

$$\sum_{i,j} a_{ij} x^i \Big(\sum_{k \geq 1} u_k x^k\Big)^j = \sum_{l \geq 1} v_l x^l.$$

We find

$$v_1 = a_{10} + a_{01} u_1$$
$$v_2 = a_{20} + a_{11} u_1 + a_{02} u_1^2 + a_{01} u_2$$
$$v_3 = a_{30} + a_{21} u_1 + a_{12} u_1^2 + a_{03} u_1^3 + a_{11} u_2 + 2a_{02} u_1 u_2 + a_{01} u_3$$
$$v_4 = a_{40} + a_{31} u_1 + a_{13} u_1^3 + a_{04} u_1^4 + a_{21} u_2 + 2a_{12} u_1 u_2 + 3a_{03} u_1^2 u_2 + a_{02} u_2^2 + a_{11} u_3 + 2a_{02} u_1 u_3 + a_{01} u_4$$

*etc.*

Even if this is complicated, we can prove immediately by induction that $v_l$ is written as

$$v_l = G_l\left((a_{ij})_{i+j\le l}, (u_k)_{k\le l-1}\right) + a_{01}u_l$$

where $G_l$ is a polynomial expression with positive integral coefficients involving the $a_{ij}$'s for $i+j \le l$ and the $u_k$'s for $k \le l-1$.

*Exercise*: Read carefully this proof and produce a sharp version of the implicit function theorem. Suppose that $F$ is convergent in some bi-disc $\{|x| < \alpha\,;\,|y| < \beta\}$ and that $|F(x,y)| < M$ in this bi-disc. What is the maximum value of $\rho(\alpha,\beta,M,|a_{01}|) > 0$ such that one can guarantee that $f$ is convergent on the open disc of radius $\rho$?

Our problem is to show that given a *convergent* series $F$, there is a unique *convergent* $f(x)$ such that $F(x,f(x)) = 0$. In other words, we are given the $a_{ij}$'s such that $|a_{ij}| \le Cr^{i+j}$ and we want to show that the equations $v_l = 0$ with unknowns $u_l$ have a unique convergent solution.

By our hypothesis, $a_{01} \ne 0$, so that, multiplying $F$ by $-1/a_{01}$, we can suppose that $a_{01} = -1$. In the same way, changing $x, y$ by constant multiples, we can assume that $C = 1$ and $r = 1$. In other words, we assume that $|a_{ij}| \le 1$ for all $i,j \ge 0$.

Since $G_l$ only depends on the $u_k$'s for $k \le l-1$ (and the $a_{ij}$'s), the previous formulae define by induction a unique series $u_l$ (depending on the $a_{ij}$'s):

$$\begin{aligned}
&u_1 = a_{10}\\
&u_2 = a_{20} + a_{11}u_1 + a_{02}u_1^2 = a_{20} + a_{11}a_{10} + a_{02}a_{10}^2\\
&\ldots\\
&u_l = G_l\left((a_{ij})_{i+j\le l}, (u_k)_{k\le l-1}\right).
\end{aligned}$$

Our task is to show that this series $\sum_l u_l x^l$ is convergent.

Now comes Cauchy's simple and beautiful idea (see[64] for details). We are going to check the theorem in *one* specific example and then show that the general case follows at once.

[64] U. Bottazzini and J. Gray. *Hidden harmony—geometric fantasies*. Sources and Studies in the History of Mathematics and Physical Sciences. Springer, New York, 2013. The rise of complex function theory.

For this example, let us choose $\overline{F}$ such that $\overline{a}_{01} = -1$ and all other $\overline{a}_{ij} = 1$:

$$\overline{F}(x,y) = -y + x + x^2 + xy + y^2 + x^3 + x^2y + xy^2 + y^3 + \cdots.$$

Let $\overline{u}_l$ be the corresponding sequence associated to this choice of $\overline{F}$ defined by:

Note that $G_l$ does not involve $\overline{a}_{01}$ which is the only $\overline{a}_{ij}$ which is not equal to 1.

$$\overline{u}_l = G_l\left((1), (\overline{u}_k)_{k\le l-1}\right) \quad (k = 1, 2, \ldots).$$

It is easy to solve $\overline{F}(x,y) = 0$ since the equation

$$\frac{1}{(1-x)(1-y)} - 1 - 2y = 0$$

is equivalent to $y = \frac{1}{4}\left(1 \pm \sqrt{\frac{1-9x}{1-x}}\right)$. In the neighborhood of 0, we have to choose the − sign and there is a unique analytic solution

$$\begin{aligned} y &= f(x) = \frac{1}{4}\left(1 - \sqrt{\frac{1-9x}{1-x}}\right) \\ &= \overline{u}_1 x + \overline{u}_2 x^2 + \cdots \\ &= x + 3x^2 + 13x^3 + 71x^4 + 441x^5 + 2955x^6 + \cdots. \end{aligned}$$

The coefficients $\overline{u}_k$ obviously satisfy some inequality $\overline{u}_k \le c\rho^k$ (with $\rho < 9$) since $f$ is analytic in some neighborhood of 0.

Now let us study the case of *a general F*, for which we assumed $|a_{ij}| \le 1$. Since the polynomials $G_l$ have *positive integral coefficients*, it follows by induction that $|u_l| \le \overline{u}_l$. Indeed:

$$\begin{aligned} |u_{l+1}| &= |G_l\left((a_{ij})_{i+j\le l+1}, (u_k)_{k\le l}\right)| \\ &\le |G_l\left((1), (|u_k|)_{k\le l}\right)| \\ &\le |G_l\left((1), (\overline{u}_k)_{k\le l}\right)| \\ &\le \overline{u}_{l+1}. \end{aligned}$$

In particular, $|u_k| \le c\rho^k$ and the series $f(x) = \sum_k u_k x^k$ is convergent. The proof of the theorem is almost finished. We found a convergent solution $y = f(x)$ and we still have to show that *all* solutions of $F(x,y) = 0$ in the neighborhood of the origin are of the form $(x, f(x))$.

When $K = \mathbb{C}$, i.e. when $F$ is a holomorphic function defined in some neighborhood of $(0,0) \in \mathbb{C}^2$, Cauchy gives an integral formula for the implicit function:

$$f(x) = \frac{1}{2i\pi}\int_C \frac{y\frac{\partial F}{\partial y}}{F(x,y)}\,dy$$

where $C$ is a small circle around the origin. Prove it!

In the ring $K\{x,y\}$ it is clear that an element $F(x,y)$ is divisible by $y$ if and only if it vanishes when 0 is substituted for $y$. The transformation $(x,y) \mapsto (x, y - f(x))$ induces an automorphism of $K\{x,y\}$ mapping $y$ to $(y - f(x))$. We know that $y = f(x)$ is a solution to $F(x,y) = 0$ so that the previous remark implies that $F$ is divisible by $y - f(x)$ in $K\{x,y\}$. The quotient is nonzero at $(0,0)$ since

$$F(x,y) = -y + a_{10}x + \cdots \quad \text{and} \quad f(x) = a_{10}x + \cdots.$$

Therefore $F(x,y) = U(x,y)(y - f(x))$, where $U \in K\{x,y\}$ is such that $U(0,0) \neq 0$. In particular in the neighborhood of $(0,0)$, the equation $F(x,y) = 0$ is indeed equivalent to $y = f(x)$. The implicit function theorem is proved. ⊡

For a good description *du calcul des limites*, I recommend Hille's book[65].

[65] E. Hille. *Ordinary differential equations in the complex domain*. Dover Publications, Inc., Mineola, NY, 1997. Reprint of the 1976 original.

## *Puiseux theorem*

Recall that we have already solved implicit equations of the form $F(x,y) = 0$ when $F$ is a non-trivial *formal* series in $K[\![x,y]\!]$ where $K$ is an algebraically closed field of characteristic 0.

We showed (with the help of Newton and Cramer) that any nonzero element $F(x,y)$ in $K[\![x,y]\!]$ can be split as a product:

$$F(x,y) = A(x,y)x^r(y - f_1(x))(y - f_2(x))\cdots(y - f_n(x))$$

where $A(0,0) \neq 0$ and the $n$ solutions $f_i(x)$ are formal *Puiseux series* in $K[\![x^\star]\!]$.

Our goal now is to show that if $F(x,y)$ is a convergent series so are the $f_i(x)$'s.

*We assume that $K$ is an algebraically closed field equipped with a complete norm.*

Since every element $f(x)$ in $K[\![x^\star]\!]$ lies in a ring $K[\![x^{1/m}]\!]$ for some $m \geq 1$, i.e. is a series in the variable $x^{1/m}$, there is no difficulty in defining convergent Puiseux series.

Denote by $K\{x^\star\}$ and $K\{x^\star, y^\star\}$ the rings of convergent Puiseux series in one and two variables.

Even though series in $K\{x^\star\}$ are convergent by definition, we should be cautious: they are not actual functions defined in the neighborhood of 0. They are multivalued functions of $x$.

**Theorem** ("Puiseux theorem"). *Any nonzero element $F$ in $K\{x,y\}$ can be split as*

$$F = A(x,y)x^r(y - f_1(x))(y - f_2(x))\cdots(y - f_n(x))$$

*where $A(0,0) \neq 0$ and the $n$ solutions $f_i(x)$ are* convergent Puiseux series *in $K\{x^\star\}$.*

The proof might look a bit cumbersome but the reader should keep in mind that this theorem is just a slight generalization of the implicit function theorem.

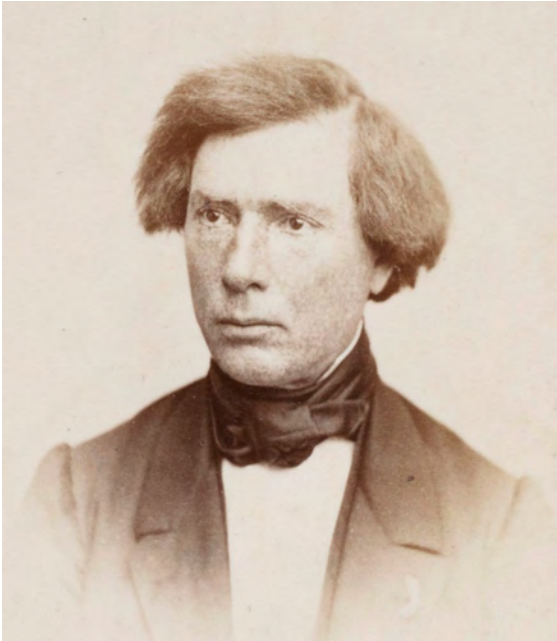

Victor Puiseux (1820 – 1883).

The completion of the algebraic closure of the field of $p$-adic numbers is a good example.

I hope that my reader has guessed the definition of the ring of convergent Puiseux series in two variables: just consider convergent power series in $(x,y)$ and replace formally $x$ and $y$ by $x^{1/m}$ and $y^{1/m}$ for some $m \geq 1$.

Recall the following facts.

1. If some formal Puiseux series $y = f(x)$ is a root of $F(x,y) = 0$, it is obtained by an application of Newton's algorithm making a choice of a side of the Newton polygon at each step.

2. At each step of the algorithm, we define $x_k = x_{k+1}^{\alpha_k}$ and $y_k = u_k x_k^{\beta_k}(1 + y_{k+1})$ for some positive integers $\alpha_k, \beta_k$, and we replace $F_k(x_k, y_k)$ by $F_{k+1}(x_{k+1}, y_{k+1}) = x_{k+1}^{-\gamma_k} F_k(x_k, y_k)$ (for some positive integer $\gamma_k$). Therefore the series $f(x)$ can also be described by series $y_k(x_k)$ ($k \geq 1$). Clearly, it is equivalent to prove the convergence of $f(x)$ or of anyone of the $y_k(x_k)$'s.

3. After a certain number of steps, the multiplicities of $F_k(x_k, y_k)$ (i.e. the valuations of $F_k(0, y_k)$) remain equal to some "ultimate multiplicity" *mult* $\geq 1$ (Cramer's theorem).

4. This "ultimate constant" *mult associated to a root* $y = f(x)$ *of* $F(x,y) = 0$ is also the multiplicity of the root, in other words the number of equal factors $(y - f(x))$ appearing in the splitting of $F$.

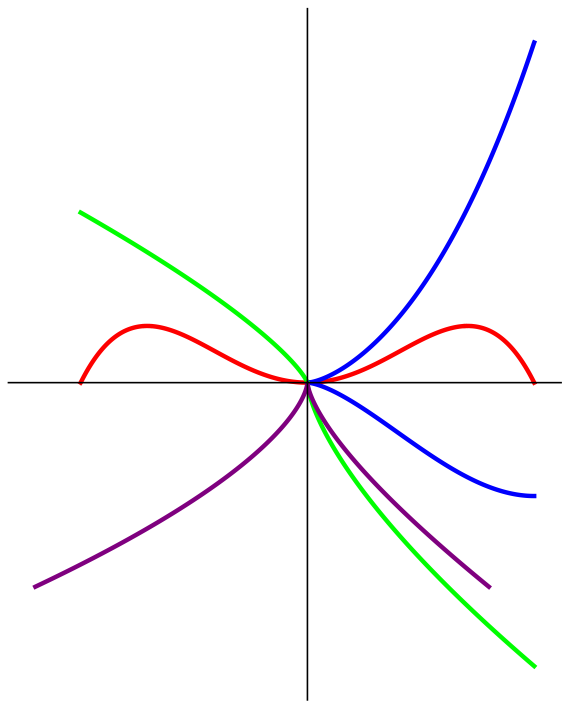

A singular point with four Puiseux branches.

We can finish the proof of Puiseux theorem.

Let $F$ be a convergent series in $K\{x,y\}$ and let

$$F(x,y) = A(x,y)x^r(y - f_1(x))(y - f_2(x))\cdots(y - f_n(x))$$

be its decomposition as a product of (a priori formal) Puiseux series. Choose $N$ such that all the $f_i(x)$'s belong to $K[[x^{1/N}]]$ and set $\overline{x} = x^{1/N}$ so that $F$ can also be seen as an element of $K\{\overline{x}, y\}$ and all the $f_i$'s as elements of $K[[\overline{x}]]$.

We are reduced to the case of $F$ in $K\{x,y\}$ such that all the $f_i(x)$'s are formal series in $K[[x]]$ and we have to prove that these $f_i$'s are actually convergent series, i.e. belong to $K\{x\}$.

If the "ultimate multiplicity" *mult* of a root $y = f(x)$ is equal to 1, the path followed by Newton's algorithm and leading to the solution $f(x)$ will eventually lead to some $F_k(x_k, y_k)$ with multiplicity 1. The implicit function theorem applied to the convergent $F_k$ shows that the solution $f(x)$ is also convergent.

If a polynomial has a multiple root, this root is also a root of its derivative. In our context, this means that if $f(x)$ is a formal series in $K[[x]]$ which is a solution of $F(x,y) = 0$ with multiplicity $mult \geq 2$, then the same series is a solution of $\partial F/\partial y(x,y) = 0$ with a smaller multiplicity. Of course, if $F(x,y)$ is convergent, so are its partial derivatives with respect to $y$. A simple induction finishes the proof. ⊡

## *Corollaries*

We have done most of the job. It is time for dessert.

First, we get the same corollaries that we had for formal series, with the same proofs.

*We continue to assume that $K$ is algebraically closed, of characteristic $0$, and equipped with a complete norm.*

**Theorem** (Weierstrass preparation theorem). *Let $F(x,y)$ be a convergent series in the ring $K\{x,y\}$ which is not divisible by $x$. Then $F$ can be written as a product $A(x,y)P(x,y)$ where $A, P$ are in $K\{x,y\}$ and*

- *$A(0,0) \neq 0$ so that $A$ is an invertible element.*
- *$P(x,y)$ is a polynomial in $y$.* ⊡

**Theorem.** *The ring $K\{x,y\}$ is a unique factorization domain.* ⊡

A very useful formulation of Puiseux theorem is given in terms of parameterization.

**Theorem** (Puiseux parameterization). *Let $F(x,y)$ be a nonzero convergent series in the ring $K\{x,y\}$, vanishing at the origin and not divisible by $x$. Then, there exist*

1. *integers $m_i \geq 1$,*
2. *open sets $U_i \subset K$ containing $0$ (for the topology defined by the norm),*
3. *series $g_i \in K\{x\}$ converging on $U_i$,*

*such that the intersection of the curve $F(x,y) = 0$ with a small neighborhood of $(0,0) \in K^2$ is the union of the images of the maps*

$$\phi_i : t \in U_i \mapsto (t^{m_i}, g_i(t)) \in K^2.$$

*Moreover these maps $\phi_i$ are injective and their images only intersect at the origin.*

Observe first that the set of roots $f_1(x), \ldots, f_n(x) \in K\{x^\star\}$ of our equation $F(x,y) = 0$ is invariant under the Galois group. For each root $f_i(x)$, let $m_i$ be the smallest integer such that $f_i(x) \in K\{x^{1/m_i}\}$. In particular $f_i(x)$ has $m_i$ distinct Galois conjugates. By reordering, we can assume that $f_1(x), \ldots, f_{m_1}(x)$ are the $m_1$ conjugates of $f_1$. None of these $f_i(x)$ $(i = 1, \ldots, m_1)$ is a "function" of $x$ in the usual sense. However, a choice of some $m_1$-th root of $x$ defines a specific value for each of the $f_i(x)$'s $(i = 1, \ldots, m_1)$. Changing the root of $x$ simply permutes the values for $f_i(x)$ $(i = 1, \ldots, m_1)$. Said differently, there is a convergent $g_1(t) \in K\{t\}$ such that the $m_1$ values of $f_1(x)$ are the $m_1$ values of $g_1(\sqrt[n]{x})$ for the $m_1$ possible choices of $\sqrt[n]{x}$. All these points are parameterized in some neighborhood of 0 by:

$$\phi_1 : t \in U \mapsto (t^{m_1}, g_1(t)) \in K^2.$$

This defines a finite number of $\phi_i$'s as in the theorem, one for each orbit of the Galois group on the set of roots. The images of the $\phi_i$'s cover the zero locus of $F$ (always in a neighborhood of the origin).

It remains to show that the $\phi_i$'s are injective and that their images only intersect at the origin.

Zeroes of an analytic function are isolated. The following lemma simply states that the same is true in $K\{x\}$ for a general $K$. I leave the proof as an exercise for the reader.

**Lemma.** *Let $h$ be a convergent series $K\{x\}$. If there is a sequence $(x_k)_{k\geq 0} \in K \setminus \{0\}$ converging to 0 such that $h(x_k) = 0$, then $h = 0$.* ⊡

Suppose now that $\phi_1$, for instance, is not injective in the neighborhood of 0. That would imply that there is some $m_1$-th root of unity $\omega_1$ such that the solutions to $g_1(\omega_1 t) = g_1(t)$ accumulate to 0. According to the lemma, $g_1(\omega_1 t)$ and $g_1(t)$ would be equal, identically, and that would contradict the fact that $m_1$ is the smallest integer such that $f_1(x) \in K\{x^{1/m_1}\}$.

The same argument shows that the intersections of the images

of

$$\phi_1 : t \in U \mapsto (t^{m_1}, g_1(t)) \in K^2 \quad ; \quad \phi_2 : t \in U_2 \mapsto (t^{m_2}, g_2(t)) \in K^2$$

is non trivial (i.e. accumulates to the origin) if and only if $m_1 = m_2 = m$ and there is some $m$-th root of the unity $\omega$ so that $g_2(t) = g_1(\omega t)$, identically. In this case, the two images actually coincide in the neighborhood of the origin. ⊡

The images of the $\phi_i$'s are usually called the *branches* of the curve $F(x,y) = 0$. The Puiseux-type parameterization of a branch is unique up to the Galois action.

In particular, a neighborhood of the origin in $\{F(x,y) = 0\}$ is homeomorphic to the union of a finite number of balls in $K$ intersecting in a single point. Note that "a ball" is an interval in $\mathbb{R}$, a disc in $\mathbb{C}$, and a Cantor set for the $p$-adic numbers.

## *Real numbers*

So far, we assumed that the field $K$ is algebraically closed. Let us study the case of real numbers which, after all, is at the origin of our promenade.

Let $F(x,y) \in \mathbb{R}\{x,y\}$ be a nonzero convergent series vanishing at the origin. Look at its zero set $\{(x,y) \mid F(x,y) = 0\}$ either as a complex curve in $\mathbb{C}^2$ or as a real curve in $\mathbb{R}^2$, in the neighborhood of $(0,0)$. Here, we are primarily interested in the description of the real curve.

For some $F(x,y) \in \mathbb{R}\{x,y\}$, it might happen that the real part of its zero set is reduced to the origin. The most obvious example is $x^2 + y^2 = 0$. Over the complex numbers this curve consists of two imaginary branches $y = ix$ and $y = -ix$, which only intersect at $(0,0)$. Of course, since we are only interested in the real part of the zero set of $F$, we simply discard all irreducible factors of $F$ whose zero sets reduce to the origin (over the real numbers).

Over the complex numbers, this zero set is the union of some branches parameterized by:

$$\phi_i : t \in U_i \mapsto (t^{m_i}, g_i(t)) \in \mathbb{C}^2.$$

Since $F(x,y)$ has real coefficients, its zero locus in $\mathbb{C}^2$ is globally invariant under complex conjugation. Since branches are disjoint away from the origin, a real point different from the origin has to belong to a branch which coincides with its conjugate. The complex conjugate of the image of $\phi_i$ is the image of

$$\overline{\phi_i} : t \in U_i \mapsto (t^{m_i}, \overline{g_i(\bar{t})}) \in \mathbb{C}^2.$$

Therefore, branches containing real points different from the origin are such that

$$\overline{g_i(\bar{t})} = g_i(\omega t)$$

for some $m_i$-th root of unity $\omega$. Writing

$$g_i(t) = \sum_{k \geq 1} a_k t^k,$$

this condition means

$$\overline{a_k} = a_k \omega^k.$$

Let $\mu$ be one of the two square roots of $\omega$ and set $t = \mu s$. Then

$$t^{m_i} = \mu^{m_i} s^{m_i} = \pm s^{m_i}$$

and

$$g_i(t) = \sum_{k \geq 1} a_k t^k = \sum_{k \geq 1} a_k \mu^k s^k = \sum_{k \geq 1} b_k s^k.$$

Now, the coefficients $b_k$ are real since

$$\overline{b_k} = \overline{a_k}\,\overline{\mu}^k = a_k \omega^k \mu^{-k} = a_k \mu^k = b_k.$$

Let us sum up this discussion.

**Theorem.** *Let $F(x,y) \in \mathbb{R}\{x,y\}$ be a nonzero converging series with real coefficients, vanishing at the origin and not divisible by $x$. Assume that the zero locus of $F$ in the neighborhood of $(0,0) \in \mathbb{R}^2$ is not reduced to the origin. Then this zero locus is the union of a finite number of curves of the form*

$$\phi_i : t \in ]-\epsilon_i, +\epsilon_i[ \mapsto (\pm t^{m_i}, g_i(t)) \in \mathbb{R}^2.$$

*where $g_i$ is a convergent series with real coefficients. The $\phi_i$ are injective and their images only intersect at the origin.* ⊡

It is easy to see that these curves $\phi_i$ are transverse to small circles centered at the origin. Indeed, tangent points correspond to the vanishing of

$$\frac{d}{dt}(t^{2m_i} + g_i^2(t)) = 2m_i t^{m_i - 1} + 2 g_i(t) g_i'(t)$$

whose zeroes are isolated. Note that this expression cannot be identically 0 otherwise the curve would be a circle.

We proved more than Gauss's claim: locally a real analytic curve is made out of a finite number of branches such that the following properties hold.

- Each branch is homeomorphic to $]-\epsilon,+\epsilon[$ and is transverse to small circles centered at the origin and intersects these circles in exactly two points (one for $t>0$ and one for $t<0$).

- Two different branches only intersect at the origin.

- Along a branch $y/x$ converges when $t$ tends to $0$ to some limit in $\mathbb{R}\cup\{\infty\}$. This means that every branch has a well-defined tangent at the origin.

In particular an algebraic curve cannot reach the origin as an infinite spiral.

Note that if $(0,0)$ is an isolated point in the zero locus of $F(x,y)\in\mathbb{R}\{x,y\}$, every complex branch only intersects its complex conjugate at the origin. In this case, there must be an even number of branches, conjugate two by two, and the multiplicity of $F$ is even. Here is a simple corollary, analogous to the fact that every odd degree real polynomial has a real root.

*Let $F(x,y)$ be a nonzero converging series with real coefficients, vanishing at the origin, not divisible by $x$, and with odd multiplicity. Then the real curve $F(x,y)=0$ is not reduced to the origin. For small real values of $x$, there is at least one real solution to $F(x,y)=0$.*

This simple fact has been transformed in a powerful tool by Poincaré who used it in numerous situations, like for instance for proving the existence of periodic orbits in the 3-body problem (see[66] page 70). This is his *continuity method*.

[66] H. Poincaré. *Les méthodes nouvelles de la mécanique céleste. Tome I.* Les Grands Classiques Gauthier-Villars. Librairie Scientifique et Technique Albert Blanchard, Paris, 1987. Reprint of the 1892 original.

## *Chord diagrams*

The local topology of an analytic curve in the neighborhood of a singular point suggests the following definition, which will be important in the rest of this book.

**Definition.** 1. A *chord diagram* is a set of $2n$ points on a circle equipped with some involution with no fixed points. In other words, a collection of $2n$ points paired two by two.

2. Two chord diagrams are *equivalent* if there is an orientation preserving homeomorphism of the circle mapping the first to the second, and commuting with the involution. In other words, we consider a cyclic word on $2n$ letters where each letter appears exactly twice. We can also draw chords connecting pairs. This is sometimes called a *Gauss word*, or a *matching*, or a *pairing*, depending on the context. I had to make a choice and I chose *chord diagram*.

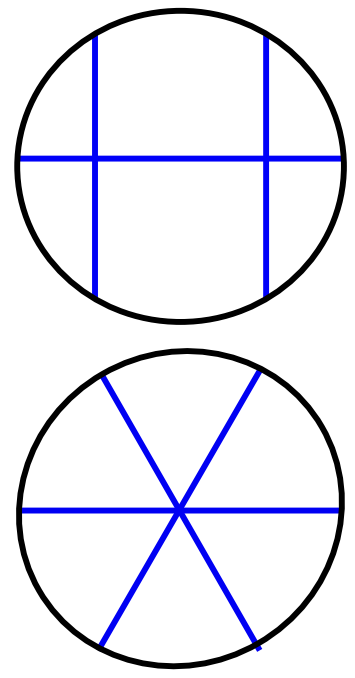

Two diagrams with three chords.

3. The chord diagram associated to an analytic curve at some (singular) point is the chord diagram obtained by intersecting the curve with a small circle around the point, where pairs of points correspond to branches. Such a chord diagram is called *analytic*.

We would like to understand analytic chord diagrams and the topology of real analytic curves.

Be patient! We'll get there.

## *A controversy concerning the shape of bird beaks?*

In 1751, Euler wrote a very interesting paper (in French) about the shape of analytic curves. In the introduction, he mentions that

> Even Geometry is not exempt from controversies and apparent contradictions, although we quite often maintain the contrary.

*"Même la géométrie n'est pas exemte [sic] de controverses, & des contradictions apparentes, quoi qu'on soutienne souvent le contraire."*

Should I give examples of modern controversies?

The controversy that Euler wanted to elucidate concerns the shape of cuspidal points[67]. There was a disagreement between Mr. le Marquis de l'Hôpital and Mr. Guà de Malves. Euler acted as a judge and dissipated the apparent contradictions in a brilliant way.

[67] L. Euler. Sur le point de rebroussement de la seconde espèce de M. le Marquis de l'Hôpital. *Mémoires de l'académie des sciences de Berlin*, 5:203–221, 1751. See the Euler archive for English translations and comments.

So far, our discussion only concerned the *topology* of branches in the neighborhood of a singular point. We did not say much about their *geometry*. We only mentioned that a branch has a tangent at the singular point.

L'Hôpital's book is entitled *Analyse des infiniment petits pour l'intelligence des lignes courbes* and was published in 1696. It is the first textbook on differential calculus. It contains a classification of singular branches of analytic curves in four categories.

Let me express this in modern terminology. Choose coordinates so that the tangent is $y = 0$. Locally, our branch is the union of two *half branches* which are graphs of two functions $h_1(x), h_2(x)$, defined in small intervals of the form $]-\epsilon, 0]$ or $[0, \epsilon[$. These functions are smooth, away from the origin.

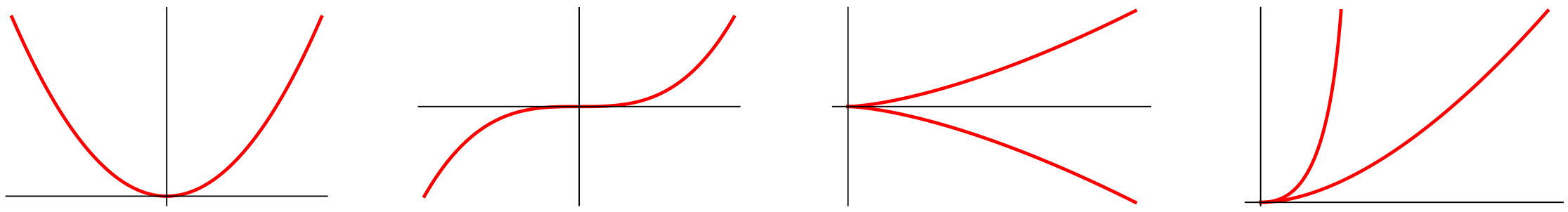

The four cases are:

1. $h_1$ is defined on $]-\epsilon, 0]$ and $h_2$ on $[0, \epsilon[$ and their second derivatives have the same sign. In this case, the curve is convex (or concave) and is on one side of its tangent.

2. An *inflexion point*. The same as before except that the second derivatives have different signs.

3. A standard *cuspidal point*. Here, $h_1$ and $h_2$ are defined on the same side of the origin and their second derivatives have opposite signs. So both half-branches have opposite convexities.

4. A *bird beak* ("point de rebroussement à bec d'oiseau") in which the second derivatives have the same sign on the two half branches.

The frontier between geometry and topology is unclear. Let me say that topology deals with properties invariant under homeomorphisms and geometry invariants under. . . smaller groups, like for instance euclidean isometries, projective automorphisms, or simply diffeomorphisms. For instance, I would consider the existence of a tangent to a curve as a geometric property.

It is very easy to find examples of the first three categories. As for the fourth category, l'Hôpital gave the following example. Wrap a thread on some curve with an inflection point, and attach it at some other point on the curve. When you unwrap it, keeping it tight, the end point will describe a curve (called the *involute*) which will present such a bird beak. I simply chose $y = x^3$ as an inflection curve and I asked my computer to draw l'Hôpital's curve. The result is in the margin. Indeed,

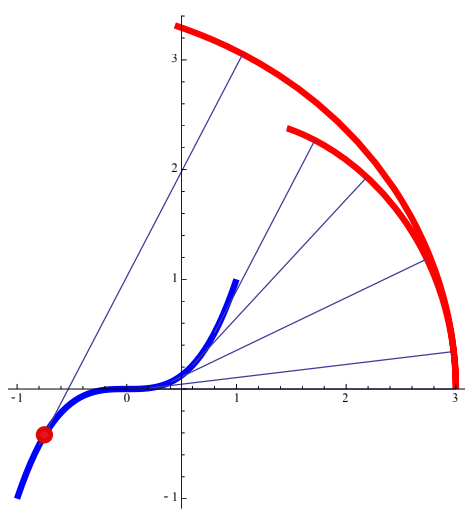

The length of the thread plus the curvilinear length along the curve is constant.

the end point of the thread describes a red curve which presents a bird beak when the thread is tangent at the inflection point, as claimed. The half-branches have the same concavity. This is "mechanically" obvious for l'Hôpital.

In 1740, Mr. Guà de Malves published an amazing book[68] whose purpose was to avoid Newton's techniques and to use only Descartes. I should recall the controversy between the English and the French during the eighteenth century around Descartes and Newton.

[68] J.-P. Gua de Malves. *Usages de l'analyse de Descartes pour découvrir, sans le secours du Calcul Différentiel, les Propriétés, ou affectations principales des lignes géométriques de tous les ordres*. Briasson, 1740.

As an illustration of this Anglo-French war, I recommend the *lettres sur Descartes et Newton*, by Voltaire.

> A Frenchman who arrives in London, will find philosophy, like everything else, very much changed there. He had left the world a plenum, and he now finds it a vacuum. At Paris the universe is seen composed of vortices of subtile matter; but nothing like it is seen in London. In France, it is the pressure of the moon that causes the tides; but in England it is the sea that gravitates towards the moon; so that when you think that the moon should make it high tide, those gentlemen fancy it should be low tide, which very unluckily cannot be proved. For to be able to do this, it is necessary the moon and the tides should have been inquired into at the very instant of the creation.

*"Un Français qui arrive à Londres trouve les choses bien changées en philosophie comme dans tout le reste. Il a laissé le monde plein ; il le trouve vide. À Paris, on voit l'univers composé de tourbillons de matière subtile ; à Londres, on ne voit rien de cela. Chez nous, c'est la pression de la lune qui cause le flux de la mer chez les Anglais, c'est la mer qui gravite vers la lune, de façon que, quand vous croyez que la lune devrait nous donner marée haute, ces Messieurs croient qu'on doit avoir marée basse ; ce qui malheureusement ne peut se vérifier, car il aurait fallu, pour s'en éclaircir, examiner la lune et les marées au premier instant de la création."*

Anyway, Guà's book is about a debate, still active today: should algebraic geometry use transcendental tools from differential geometry? Among the "theorems" in this book, one finds the claim that l'Hôpital is wrong and that bird beaks don't exist.

Guà is aware of l'Hôpital's example but he criticizes it in the following way. Suppose you look at two parabolas $y = x^2$ and $y = 2x^2$ but only for $x \geq 0$. You get two convex half parabolas whose union looks like a bird beak. Therefore, according to Guà, l'Hôpital's example of a beak is artificial: the complete analytic curve contains two smooth branches, as in the parabolas example, and the mechanical construction using the thread is missing one half of the curve. Convincing? Guà continues and "proves" that bird beaks are impossible for an analytic curve.

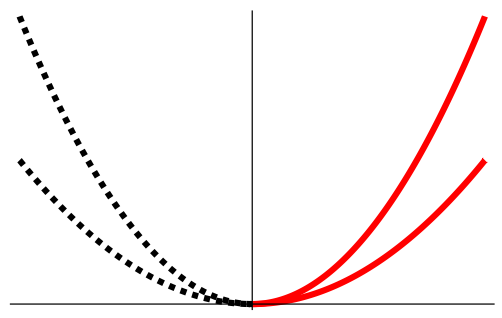

The "proof" goes more or less along the following lines. A branch has the form $y = ax^{p/q} + o(x^{p/q})$ for some pair of relatively prime integers $p, q$ with $p > q$ if $y = 0$ is the tangent at 0. If there is a beak, $q$ has to be even since otherwise $y$ would be defined for all $x$, positive or negative. The concavity is given by the sign of the second derivative, which is of the order of $a\left(\frac{p}{q}\right)\left(\frac{p}{q}-1\right)x^{\frac{p-2q}{q}}$. Since $p-2q$ is odd and $q$ is even, the two determinations of this second derivative have different signs and the two half branches have opposite convexity: this is not a bird beak.

Now comes the great Euler. His paper is very clear and unquestionable. Initially, he was convinced by Guà's argument but he found a mistake in 1744. In a column entitled "Did Euler prove Cramer's rule", Rob Bradley mentions a letter between Euler and Cramer related to this topic.

Euler, Guà, and l'Hôpital freely use Puiseux series and do not raise any doubt concerning their convergence. What is remarkable in Euler's paper is the description of the role of *complex* numbers to understand *real* algebraic curves (in 1751). Here is one of Euler's examples:

$$y = x^{1/2} \pm x^{3/4}.$$

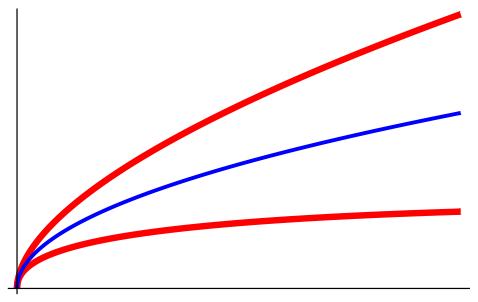

The graph in the margin does look like an eagle beak. How does one know that these two graphs, with ± signs, do belong to the same branch and cannot be completed as in our example with two parabolas? Euler gives a convincing argument using complex numbers. I strongly encourage my reader to find the mistake in Guà's "proof".

Eliminate the radicals. Euler finds

$$y^4 - 2xy^2 + x^2 - x^3 - 4yx^2 = 0$$

and you should draw the Newton polygon and check that there is indeed a single branch at the origin.

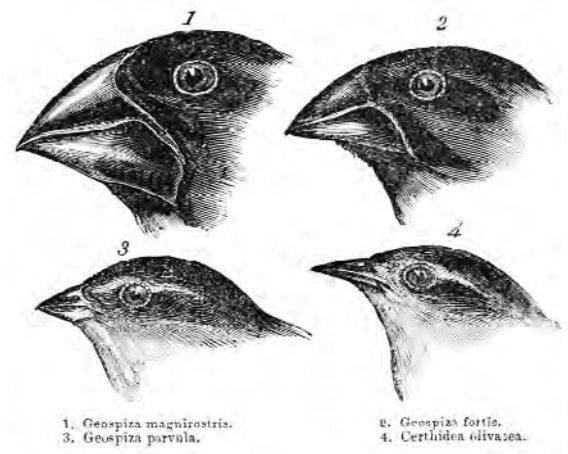

The shapes of finches beaks from the Galapagos islands were important in Darwin's discovery of evolution. ©

Today bird beaks are not mentioned anymore in math books. These points are now called *second order cusps*, in a more neutral way. Sometimes, one still sees the name *ramphoid curve*, from the greek "rampho" associated to the crooked beaks belonging to birds of prey.

To conclude this chapter, I have only one piece of advice: stop reading this book and read (some of) Euler's papers. Now!

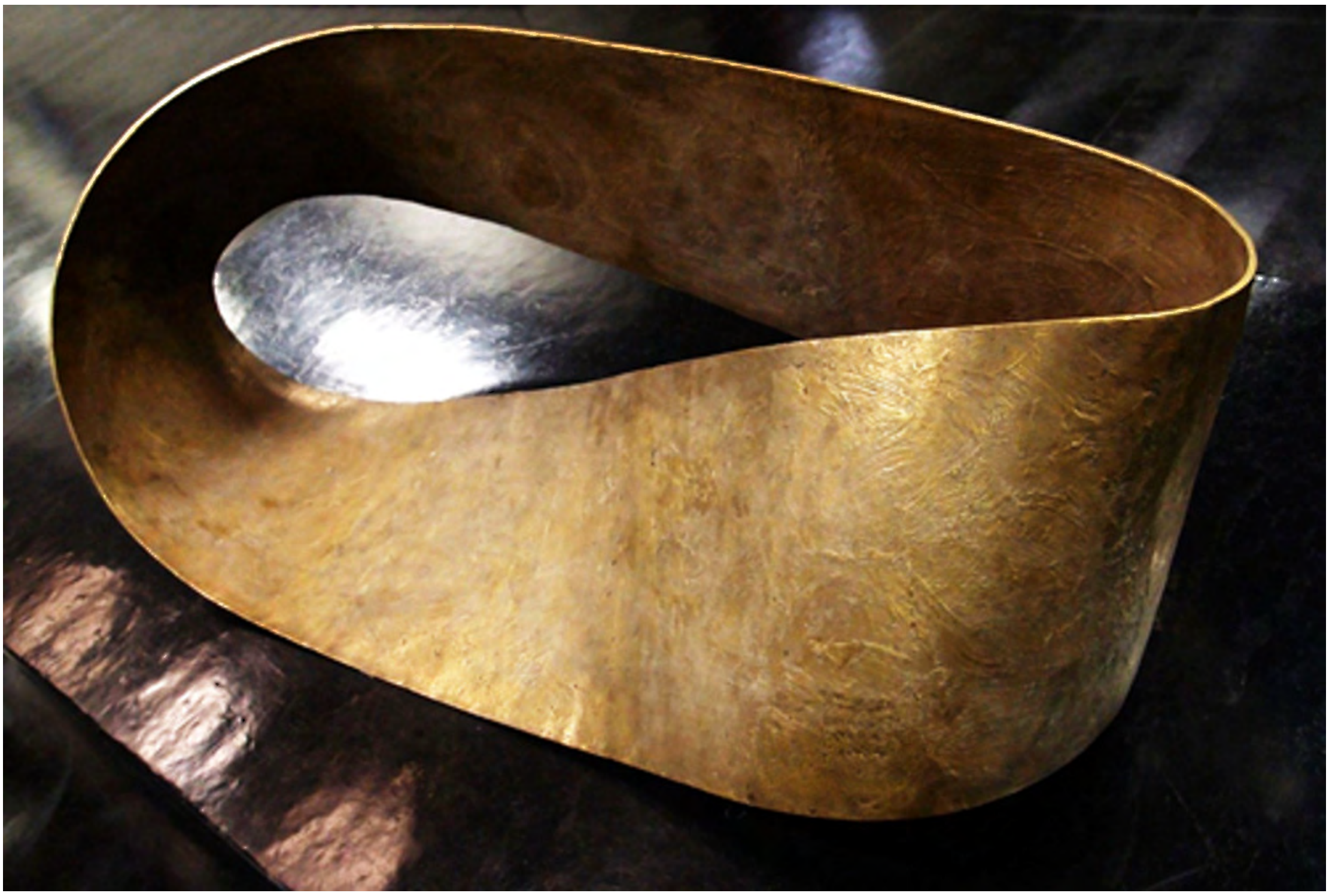

A Moebius band in the main hall of IMPA, where the first draft of this book was written.

# *Moebius and his band*

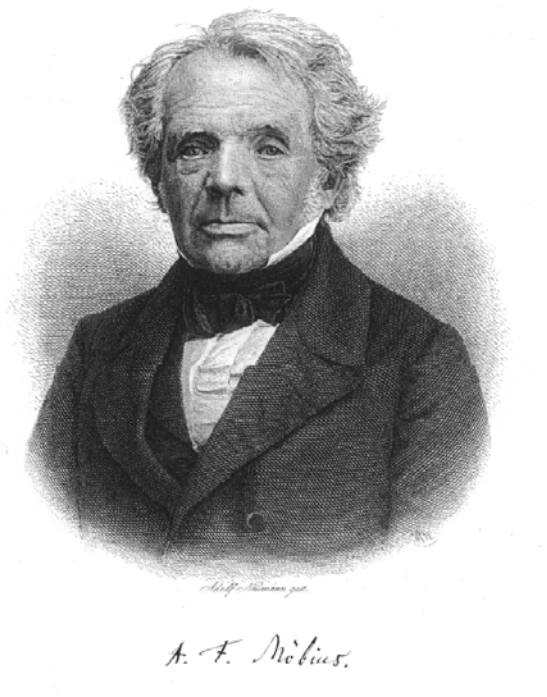

August Ferdinand Möbius (1825–1884).

THIS IS THE TITLE OF A BOOK[69] dedicated to German mathematics in the nineteenth century. In this chapter, we discuss the topology associated to the process of *desingularization* of analytic curves, leading to some beautiful necklaces made out of Moebius bands.

[69] R. F. John Fauvel, Robin Wilson. *Moebius and his Band: Mathematics and Astronomy in Nineteenth-Century Germany*. OUP, 1993.

## *Polar coordinates*

Look at the following familiar map from a cylinder to a plane:

$$\Phi : (\rho, \theta) \in \mathbb{R} \times \mathbb{R}/2\pi\mathbb{Z} \mapsto (\rho \cos\theta, \rho \sin\theta) \in \mathbb{R}^2.$$

It has the following properties:

1. $\Phi$ restricted to $\mathbb{R}_+^\star \times \mathbb{R}/2\pi\mathbb{Z}$ is a diffeomorphism onto the punctured plane $\mathbb{R}^2 \setminus \{(0,0)\}$.
2. $\Phi$ "collapses" the circle $\{0\} \times \mathbb{R}/2\pi\mathbb{Z}$ to the origin.
3. The inverse image by $\Phi$ of a point which is different from the origin contains precisely two points, of the form $(\rho, \theta)$ and $(-\rho, \theta + \pi)$.

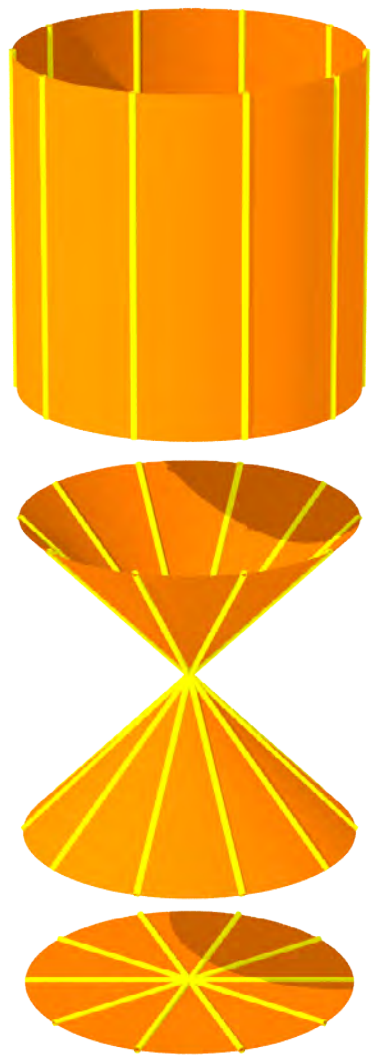

Property 3. is not very convenient for a coordinate system and this is the main reason why $\Phi$ will be slightly modified in a moment. Sometimes, one restricts $\Phi$ to $\mathbb{R}_+ \times \mathbb{R}/2\pi\mathbb{Z}$ but this introduces some artificial boundary.

Property 2. is interesting in the context of *desingularization*. In a small neighborhood of the origin, $\Phi^{-1}$ is behaving like a microscope: tiny circles $x^2 + y^2 = \epsilon^2$, of perimeter $2\pi\epsilon$, are mapped by $\Phi^{-1}$ to two big circles $\{\pm\epsilon\} \times \mathbb{R}/2\pi\mathbb{Z}$, of perimeter $2\pi$.

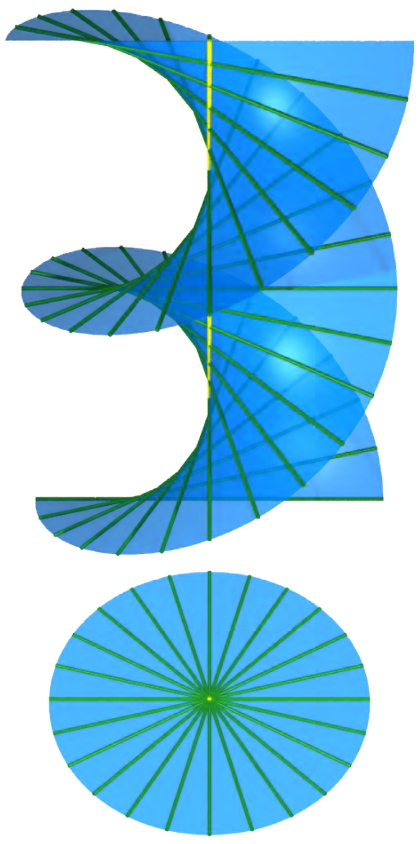

©

As a first naive example, consider a straight line $D$ passing through the origin. Its inverse image $\Phi^{-1}(D)$ consists of two "lines" $\theta = \alpha$ and $\theta = \alpha + \pi$ plus the circle $\{0\} \times \mathbb{R}/2\pi\mathbb{Z}$. Therefore, if two distinct lines $D_1$ and $D_2$ intersect at the origin, their inverse images become somehow disjoint. The operation $\Phi^{-1}$ has removed the intersection point. The "somehow" is due to the fact that $\Phi^{-1}(D)$ contains $\Phi^{-1}(0,0)$ so that the inverse images of two intersecting lines cannot be disjoint.

A better procedure is the following. Given a subset $X$ of the plane, let me denote by $\widehat{\Phi^{-1}}(X)$ the closure in $\mathbb{R} \times \mathbb{R}/2\pi\mathbb{Z}$ of $\Phi^{-1}(X \setminus \{(0,0)\})$. With this definition $\widehat{\Phi^{-1}}(D_1)$ and $\widehat{\Phi^{-1}}(D_2)$ are indeed disjoint. I will say that $\widehat{\Phi^{-1}}$ is the *strict transform*.

In order to visualize $\Phi$, look at the surface $S$ embedded in $\mathbb{R}^2 \times \mathbb{R}/2\pi\mathbb{Z}$ and defined by $x \sin\theta = y \cos\theta$. This is analogous to a double spiral staircase. The picture in the margin represents a simple staircase in $\mathbb{R}^2 \times [0, 2\pi[$. Note that $S$ is a *smooth* surface. Our map $\Phi$ corresponds to the projection onto the horizontal plane $\mathbb{R}^2$ and $\Phi^{-1}(0,0)$ is the vertical $\mathbb{R}/2\pi\mathbb{Z} \times \{(0,0)\}$.

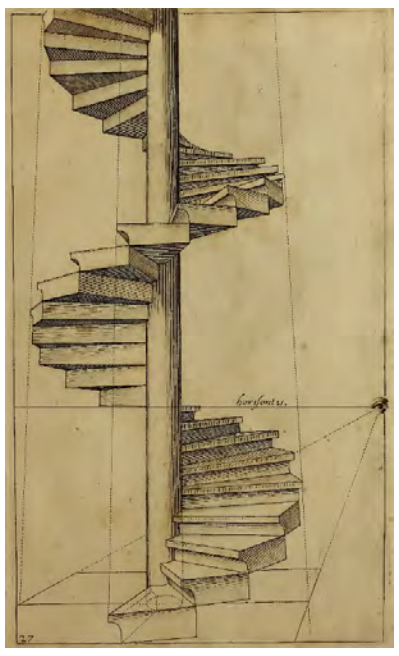

Plate 27 of Instruction en la science de perspective, H. Hondius (1625). ©

As a second simple example, consider the planar curve $x^3 = y^2$, having a cuspidal singular point at the origin. Its strict transform has equation $\rho = \sin^2\theta / \cos^3\theta$ (with two components, as it should be) and is not singular anymore. It is now smooth and tangent to the circle $\mathbb{R}/2\pi\mathbb{Z} \times \{(0,0)\}$.

The general idea is that the strict transform of a curve is "less singular" than the original curve at the origin. Repeating the operation a certain number of times, the singular curve will hopefully be transformed in a smooth curve.

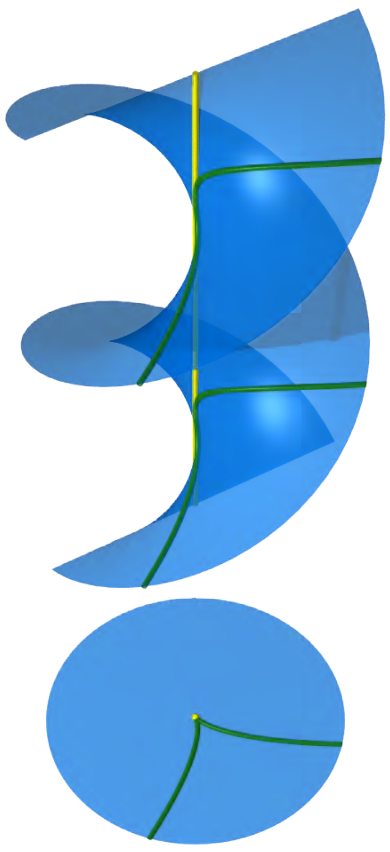

©

Before going on, we have to fix the problem that preimages by $\Phi$ contain two points. Iterating the process $n$ times, we get $2^n$ points and that is very difficult to handle.

## *The Moebius band*

Note that the involution sending $(\rho, \theta)$ to $(-\rho, \theta + \pi)$ has no fixed points. An easy way to get rid of the double preimages of $\Phi$ is to identify the points $(\rho,\theta)$ and $(-\rho,\theta+\pi)$ in $\mathbb{R}\times\mathbb{R}/2\pi\mathbb{Z}$. The corresponding quotient is a smooth surface. The involution reverses orientation since its Jacobian determinant is $-1$. The quotient surface is not orientable: this is the *famous Moebius band*[70].

The same thing could be expressed in the following way. The set of straight lines passing through the origin is a circle, can be parameterized either

– by its *slope* $t$ which is an element of the real projective line $P^1_{\mathbb{R}} \simeq \mathbb{R}\cup\{\infty\}$.

– or by its angle $\theta$, now *defined modulo* $\pi$.

Let $\mathcal{M}$ be the space of pairs $(p,D)$ where $p$ is a point in the plane and $D$ a line passing through the origin, *and through* $p$. It can be seen either as

$$\mathcal{M} = \{((x,y),t) \in \mathbb{R}^2 \times (\mathbb{R}\cup\{\infty\}) \mid y = tx\}$$

or

$$\mathcal{M} = \{((x,y),\theta)) \in \mathbb{R}^2 \times \mathbb{R}/\pi\mathbb{Z} \mid x\sin\theta = y\cos\theta\}.$$

The first presentation has the advantage of having a very simple equation and the disadvantage of not showing immediately that $\mathcal{M}$ is a smooth surface in the neighborhood of $t=\infty$. However, on a second thought, replacing $t$ by $t' = 1/t$, the equation becomes $x = t'y$ and the disadvantage disappears. The second presentation shows that $\mathcal{M}$ is indeed the quotient of $\mathbb{R}\times\mathbb{R}/2\pi\mathbb{Z}$ by the involution mentioned above.

Note that $x=y=0$ defines a circle $E$ embedded in $\mathcal{M}$: this is called the *exceptional divisor*.

Define the map

$$\Psi : ((x,y),t) \in \mathcal{M} \mapsto (x,y) \in \mathbb{R}^2.$$

It has exactly the desired properties:

1. $\Psi$ restricted to the complement of the exceptional divisor is a diffeomorphism onto the complement of the origin in the plane.

2. $\Psi$ "collapses" the exceptional divisor to the origin.

[70] P. Popescu-Pampu. La bande que tout le monde connaît. *Images des Mathématiques*, 2010.

In his 1895 memoir *Analysis Situs*, Henri Poincaré does not mention the Moebius band but "La surface unilatère que tout le monde connaît" (the unilateral surface that everybody knows).

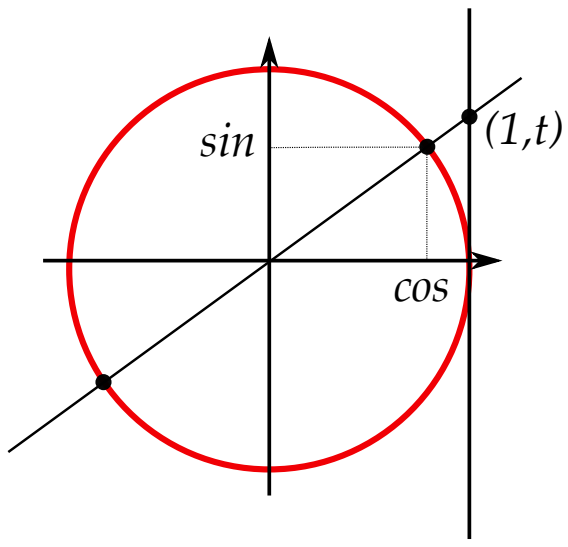

The circle of angles modulo $\pi$ is a real projective line.

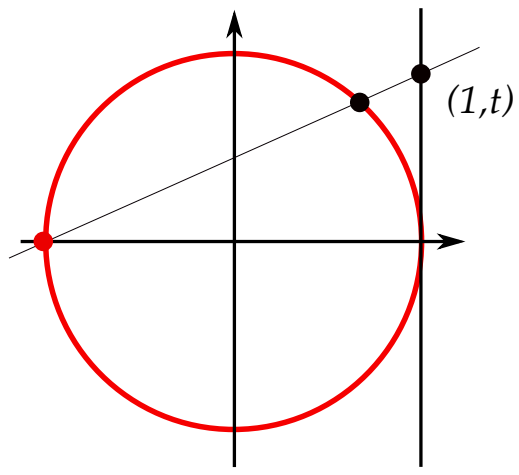

The circle of angles modulo $2\pi$ is also a real projective line.

The terminology *divisor* comes from algebraic geometry and is confusing since the exceptional divisor is a circle embedded in $\mathcal{M}$ which does *not* divide the surface in two components, unlike the core of an annulus.

For this reason, we will say that $\mathcal{M}$ has been obtained from the plane by *blowing up the origin*. Conversely, $\Psi$ is a *blowing down* map.

Since we want to work locally, it is usually useful to restrict to the compact surface with boundary

$$\overline{\mathcal{M}} = \{((x,y),\theta))|x^2+y^2 \leq 1\} \subset \mathcal{M}.$$

Clearly $\overline{\mathcal{M}}$ is obtained from $[-1,+1] \times [0,\pi]$, after gluing $(t,0)$ and $(-t,\pi)$. This is the familiar Moebius band: a rectangle where two opposite sides are identified after a twist.

It should be clear that that the boundary of $\overline{\mathcal{M}}$ is connected since it is mapped homeomorphically onto the boundary of a disc by $\Psi$.

It should be equally clear that the complement of the exceptional divisor in $\overline{\mathcal{M}}$ is connected. This is a complicated way to say that if we cut open the band along the middle circle, the result is a *bona fide* annulus. Indeed, the complement of $E$ is homeomorphic to a punctured disc.

Finally, the inverse image of a circle, say $x^2+y^2 = 1/2$ is a circle embedded in $\overline{\mathcal{M}}$ which disconnects it into two parts. The first is an untwisted annulus, mapped to $x^2+y^2 \geq 1/2$ by $\Psi$ and the other is a smaller Moebius band, mapped to $x^2+y^2 \leq 1/2$ by $\Psi$.

The terminology *exceptional* also comes from algebraic geometry and is more difficult to explain.Two transversal *complex* curves in a *complex* surface have a positive intersection number. Blowing up a point in a *complex* surface produces an exceptional divisor which is now a *complex* projective line. To compute its self-intersection number, we have to use non-holomorphic deformations. This self intersection turns out to be equal to −1, which was considered surprising and exceptional by the algebraic geometers of the past.

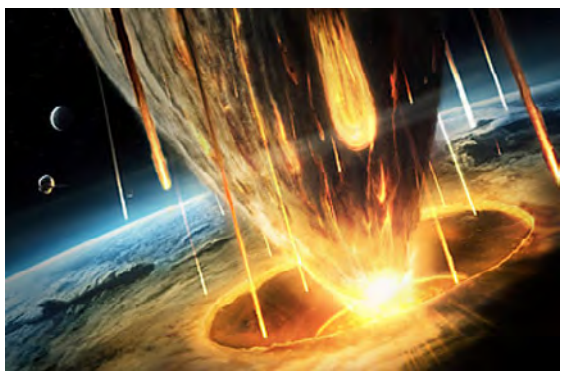

The website Impact Earth enables you to blow up our planet at any point. ©

## *Some pictures*

The Moebius band is undoubtedly one of the very few mathematical objects that are famous outside of the mathematical world. In Science Fiction, in Art, Philosophy etc.

Just for the fun of it, let me quote a few sentences from the famous psychoanalyst Jacques Lacan[71] in his seminar "l'Étourdit", in 1972:

[71] J. Lacan. *L'Etourdit*. Seuil, 1973.

> Le non-enseignable, je l'ai fait mathème de l'assurer de la fixion de l'opinion vraie, fixion écrite avec un x mais non sans ressource d'équivoque. Ainsi un objet aussi facile à fabriquer que la bande de Moebius en tant qu'elle s'imagine, met à portée de toutes mains ce qui est inimaginable dès que son dire à s'oublier, fait le dit s'endurer. D'où a procédé ma fixion de ce point doxa que

I am unable to translate into English (or even in understandable French).

> je n'ai pas dit, je ne le sais pas et ne peux donc - pas plus que FREUD - en rendre compte de ce que j'enseigne, sinon à suivre ses effets dans le discours analytique, effet de sa mathématisation qui ne vient pas d'une machine, mais qui s'avère tenir du machin une fois qu'il l'a produite.

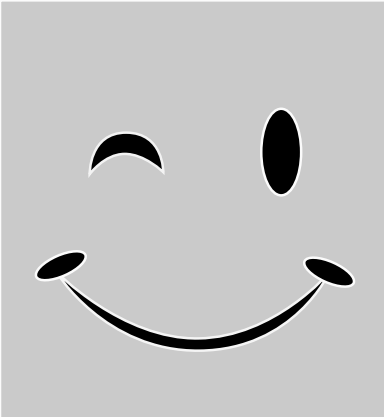

For great (and serious) comments on the Moebius band, I urge the reader to *look at* J. Scott Carter's book[72].

"look at" is more appropriate than "read" in this case.

The band is named after Moebius, who published it in 1865, but — as usual — he was not the first. Listing had described the same object in 1862.

I could easily produce a book full of Moebius bands, of different shapes, colors etc. Let me present here only a small sample.

Make a simple knot with a band of paper and tighten it. You get something like in the figure below. When you close your regular pentagon you produced a Moebius band.

[72] J. S. Carter. *How surfaces intersect in space: an introduction to topology*. K & E series on knots and everything 2. World Scientific, 2nd ed edition, 1993.

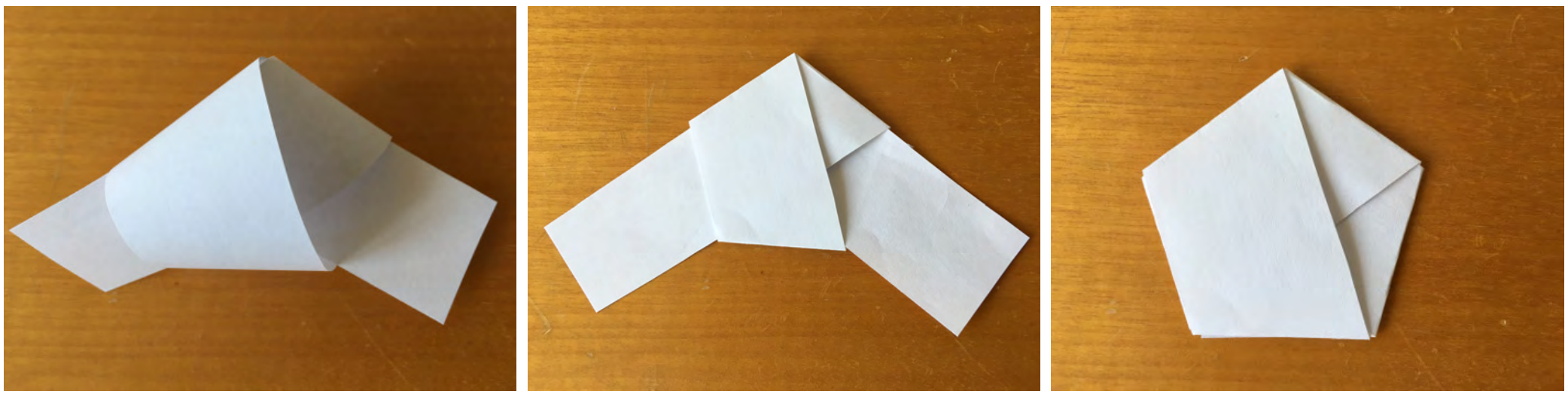

Chapter 4 of the beautiful *Topological picturebook*[73] is dedicated to the impossible *tribar*.

[73] G. K. Francis. *A topological picturebook*. Springer-Verlag, 2006.

Consider a disc, or the interior of an ellipse in the plane. Its exterior has the topology of an annulus. Now, think of this ellipse in the *real projective plane* where we have to add the line at infinity, which is topologically a circle: one point for each direction. Every line in the projective plane intersects infinity in exactly one point.

Do not confuse the real projective plane, obtained from the plane $\mathbb{R}^2$ by adding a real projective line (i.e. a circle) at infinity, with the complex projective line, obtained from $\mathbb{C}$ by adding a point at infinity.

Show that this implies that the complement of a disc in the projective plane is a Moebius band.

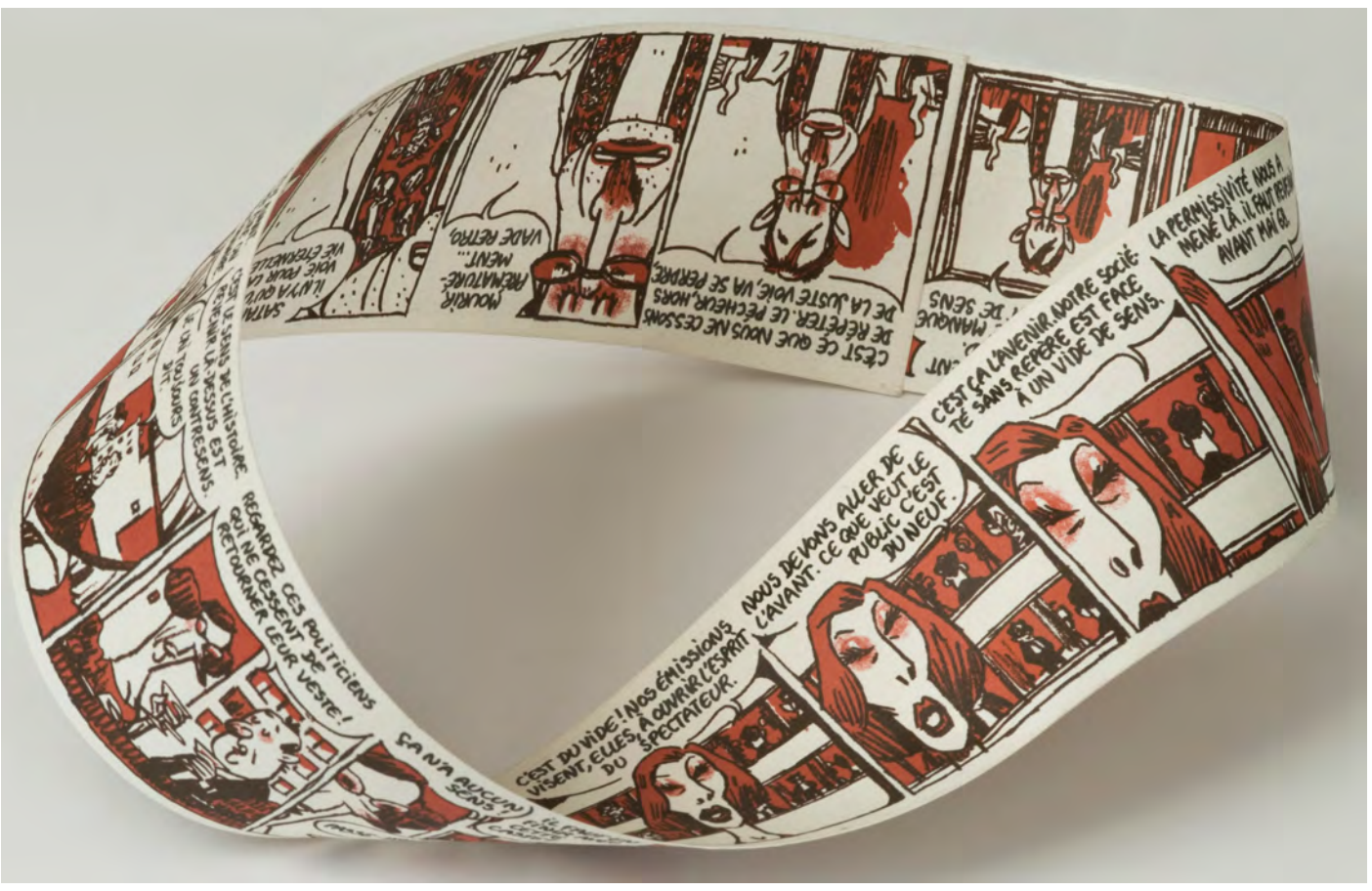

This is a Moebius cartoon, by Étienne Lecroart. You can read the story and come back to the beginning upside down and the story starts again! ©

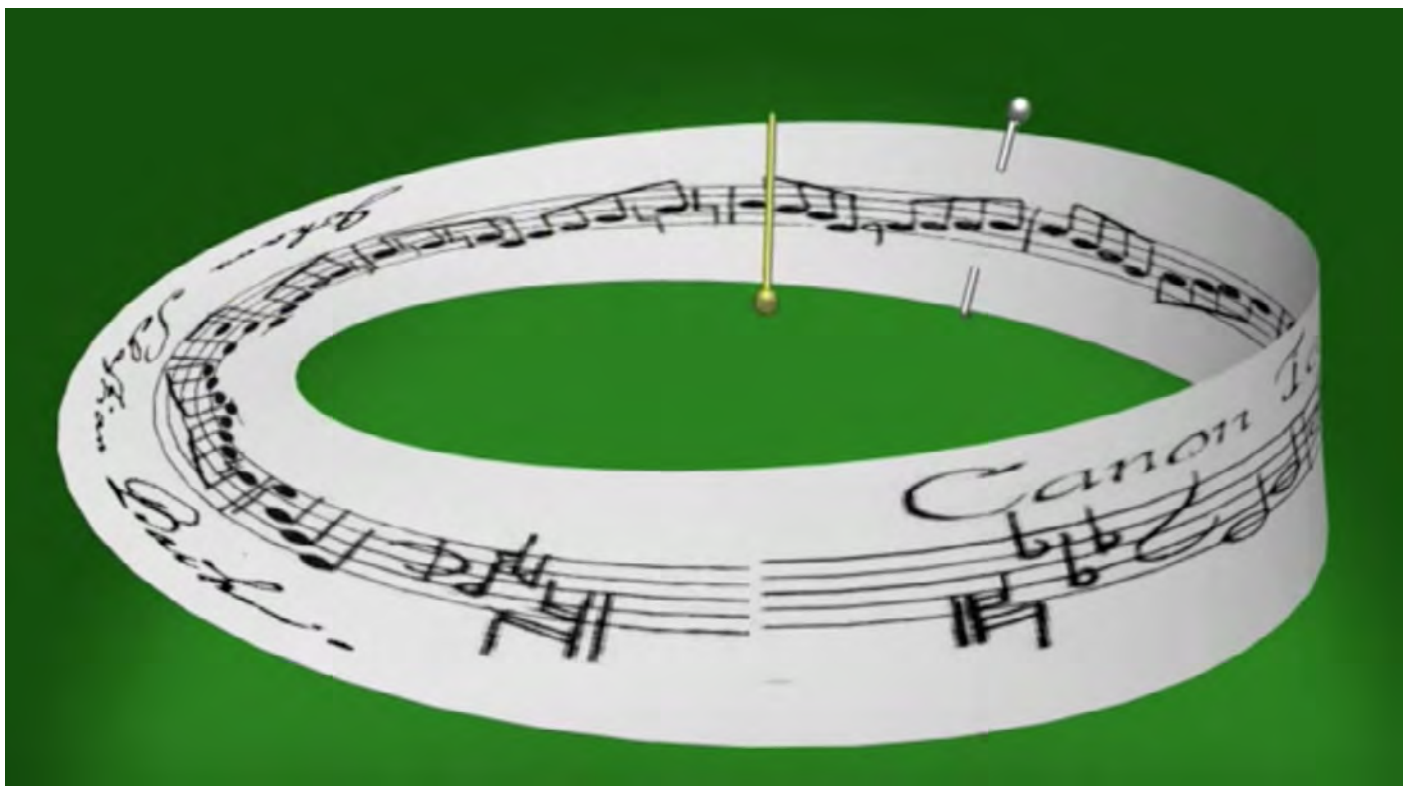

This is a tribute by J. Leys to J.S. Bach Bach crab cannon on a Moebius strip (1747), by Jos Leys. ©

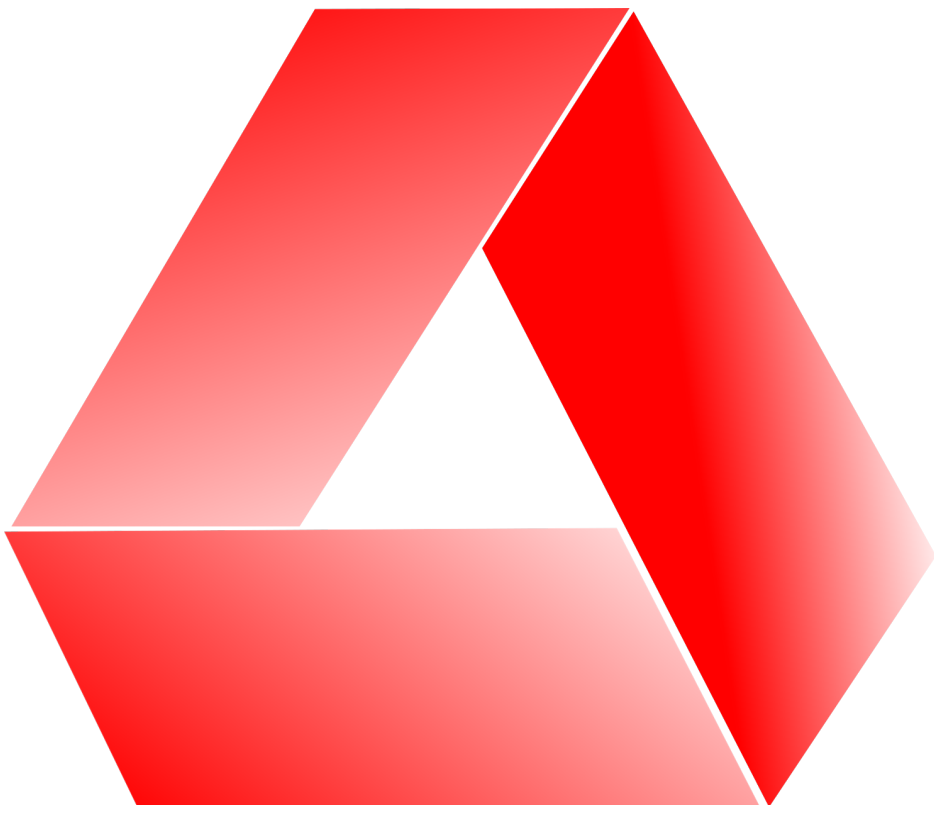

This may look misleading. Is it possible to realize it in such a way that it is made out of three *planar trapeziums* as it looks on the picture? If the three pieces are twisted, this object describes a Moebius band in space. Its boundary is a circle, as it should be, but this circle is *knotted* in space: it is a trefoil knot. This is different from the usual picture in which the boundary is unknotted.

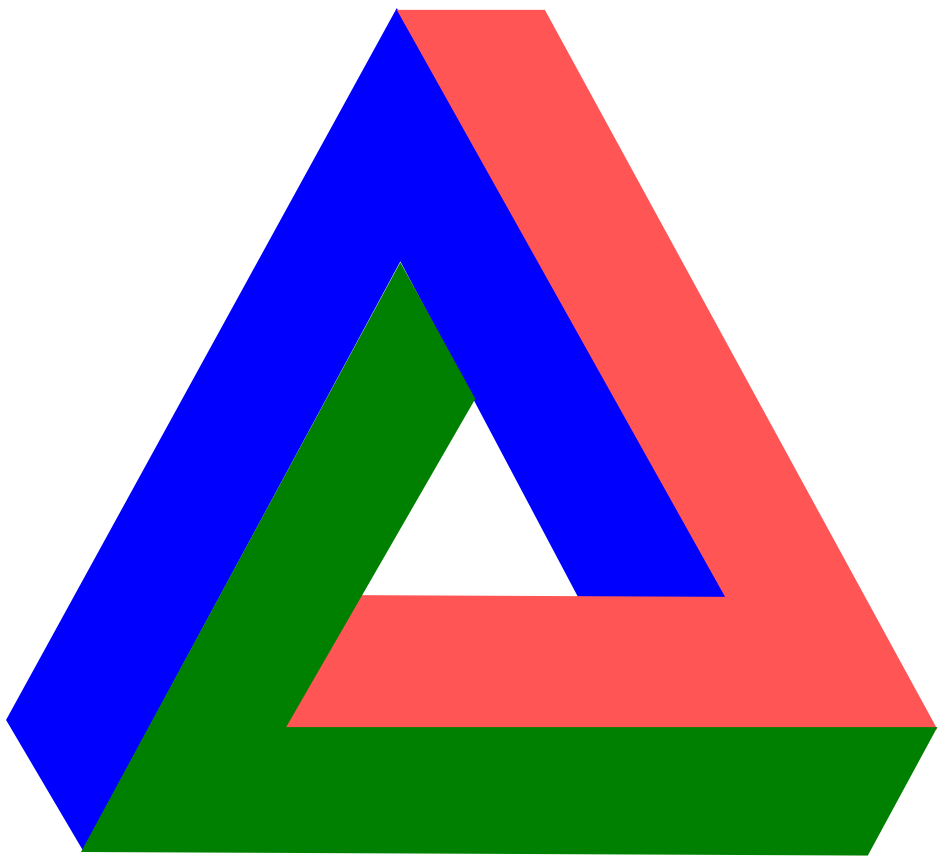

An *impossible object*. ©

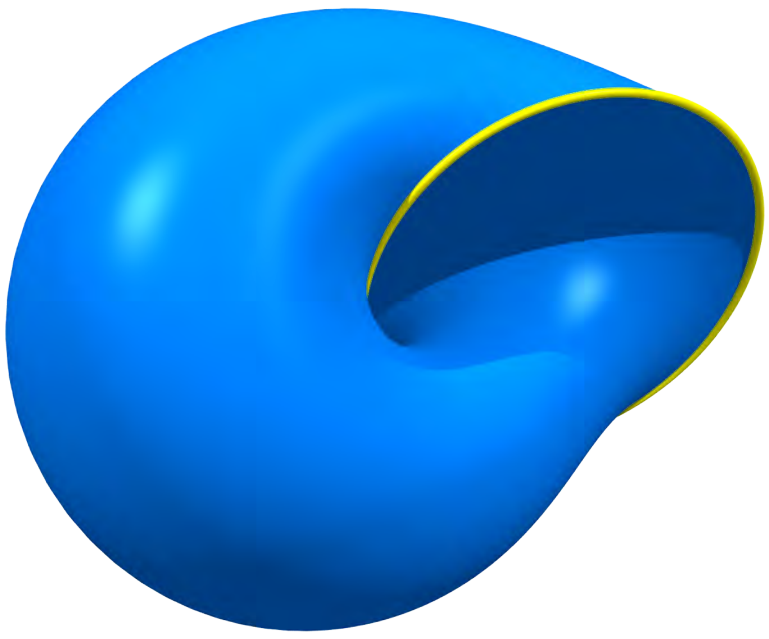

Starting from the usual Moebius band, and deforming until the boundary becomes a round circle, we get the *Moebius snail*. ©

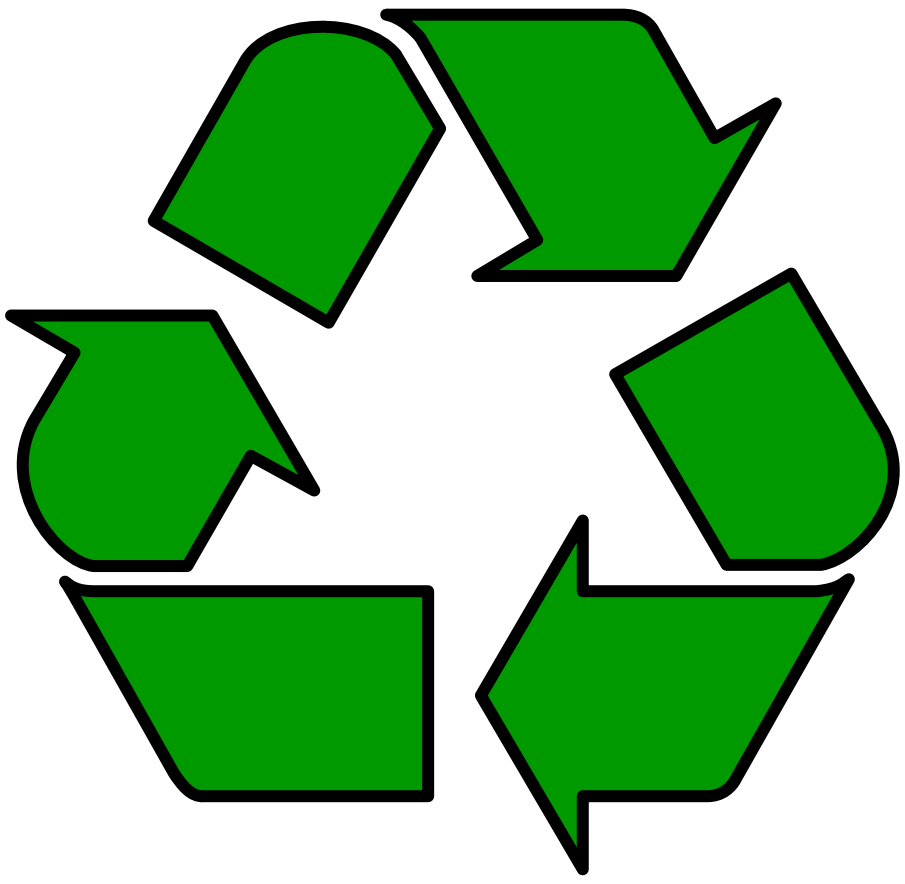

Have you noticed that the recycling symbol is a Moebius band? ©

Look carefully at the cauldron for the Rio 2016 Olympics. It consists of a large circle formed of many rotating hinges.

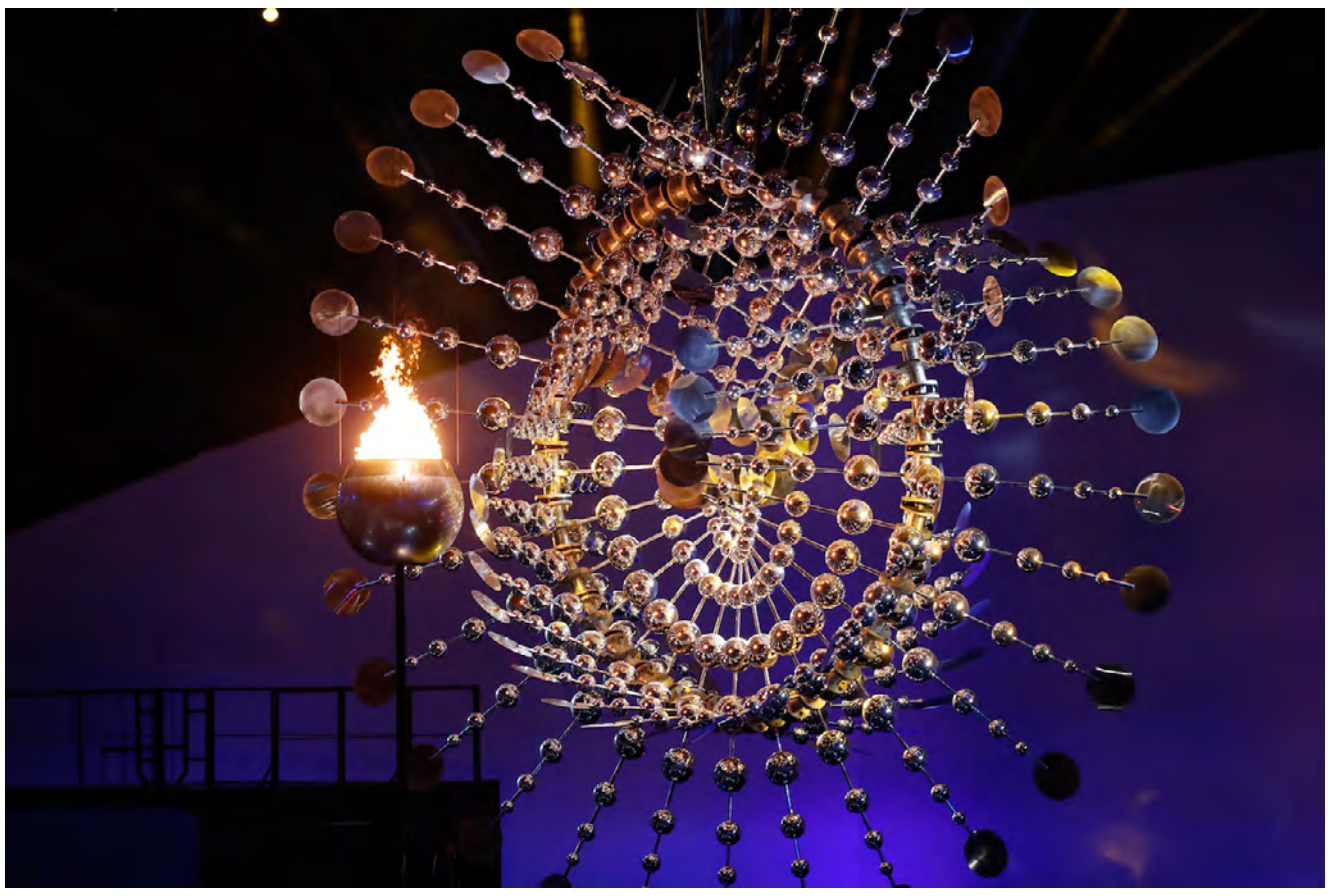

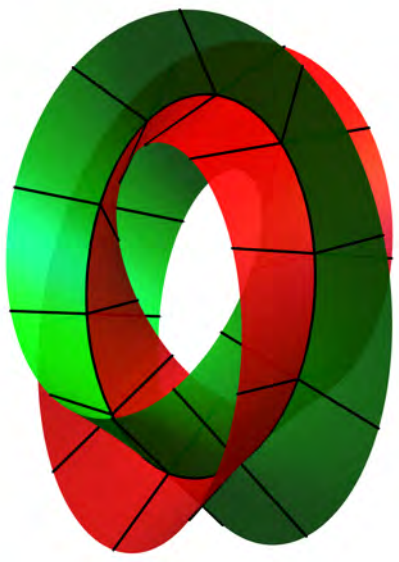

A kinetic sculpture by Anthony Howe. ©

Each hinge carries four arms. Hence, this is the union of many segments. Check that this represents two Moebius bands intersecting along their common core, as in the margin.

The following drawings show the collapse of the exceptional divisor in a Moebius band.

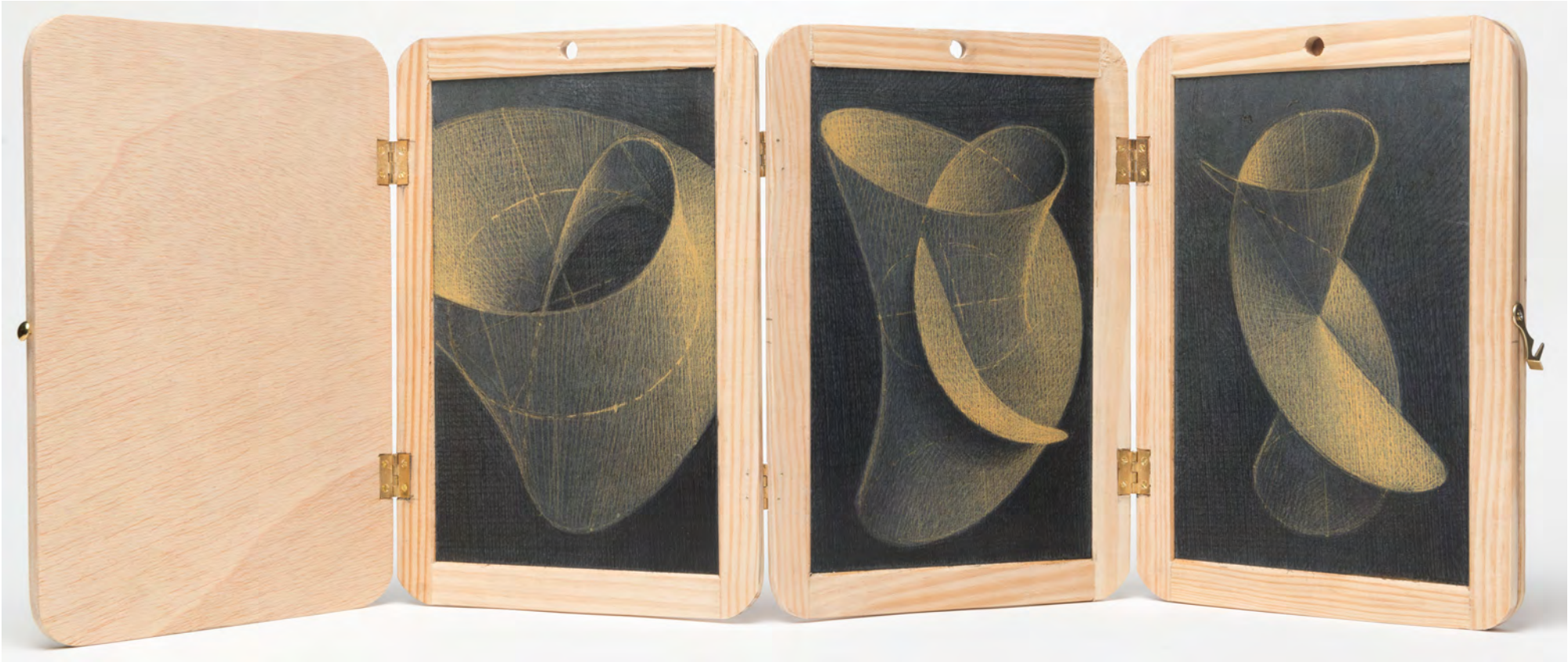

Reducing the "Möbius", crayons de couleur sur ardoises, bois, métal, by Sylvie Pic. ©

The final picture is a cone over some immersed closed loop in the sphere with two double points. A cone over a circle is a disc,

as was expected from a blowing down map. Slicing the cone with a plane, we find a Descartes folium so that the equation of the cone could be $x^3 + y^3 - 3xyz = 0$, as in the following wire models.

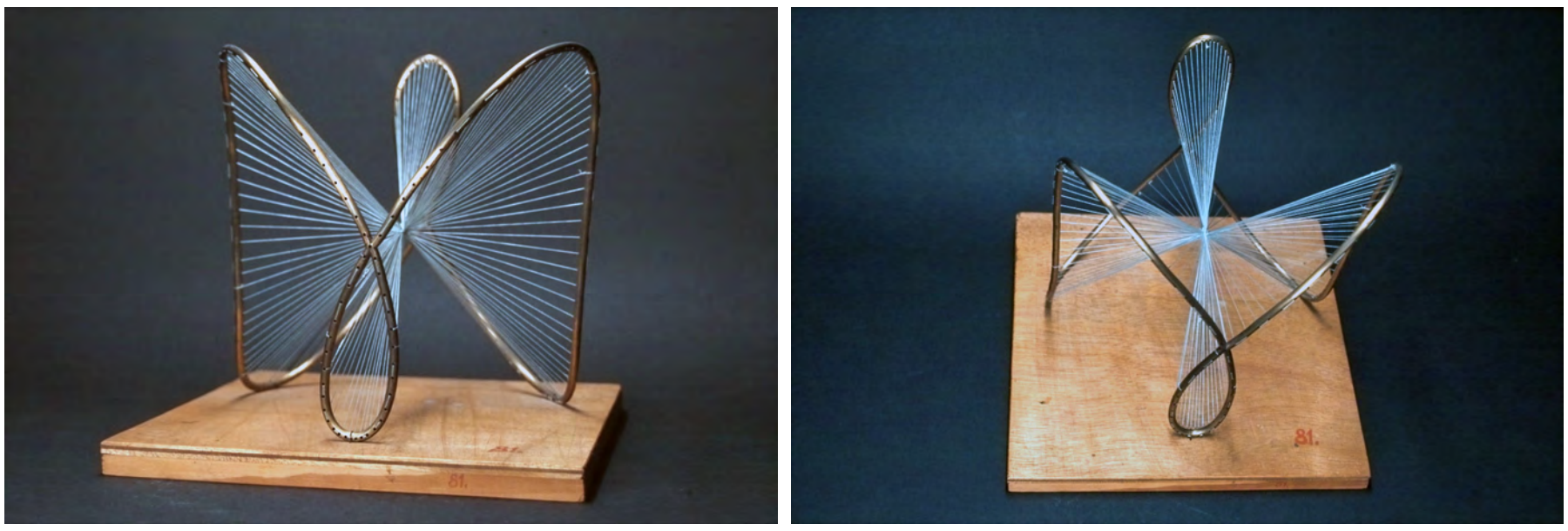

A third order cone. ©

## *Testing our microscope*

Let us test the efficiency of our microscope $\Psi^{-1}$.

If $X \subset \mathbb{R}^2$ is any set, its *strict transform* as the closure of $\Psi^{-1}(X \setminus \{(0,0)\})$ in $\mathcal{M}$.

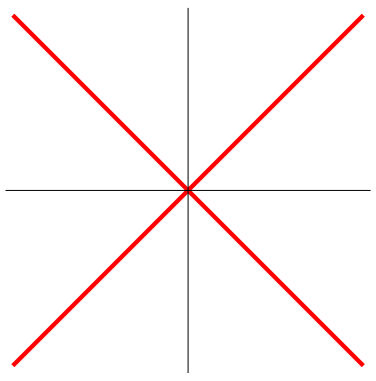

Let us try first with two intersecting lines $x^2 - y^2 = 0$. So, let us set $y/x = t$, which together with $y = \pm x$ gives $tx = \pm x$. Since we compute the strict transform, we work outside the origin and get $t = \pm 1$. Note that $t = \infty$ is not in the strict transform, as one easily sees in the chart $t' = 1/t$. Now, the closure of $t = \pm 1$ in $\mathcal{M}$ consists of two disjoint curves. Hence the strict transform of two smooth curves intersecting transversally at the origin produces *disjoint smooth curves*.

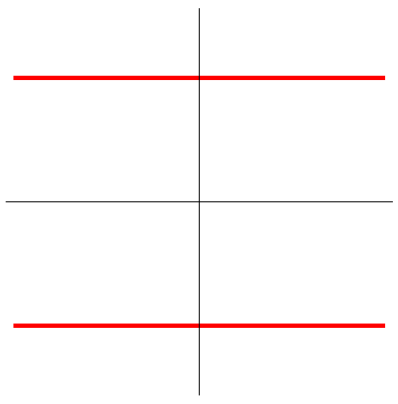

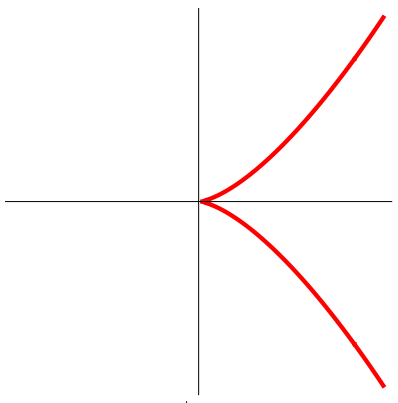

What about a cuspidal point $y^2 - x^3 = 0$? This gives $t^2x^2 - x^3 = 0$, which may be simplified by $x^2$ so that $x = t^2$. Therefore, in the coordinates $(x, t)$ of $\mathcal{M}$, the strict transform is a smooth parabola, tangent to the exceptional divisor $(x = 0)$.

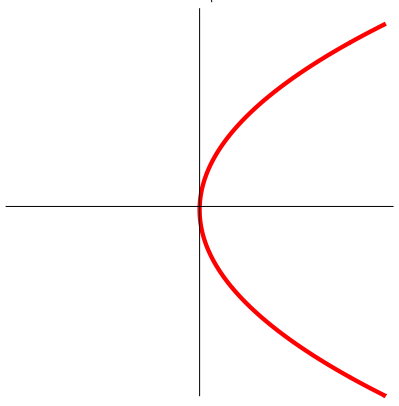

Now, let us consider $y^2 - x^5 = 0$. The strict transform is $t^2 = x^3$ and is therefore a cuspidal point.

It is clear that a single blowup will be insufficient and that we have to blow up again. Just in the same way that the Newton's algorithm does not always terminate at the first step.

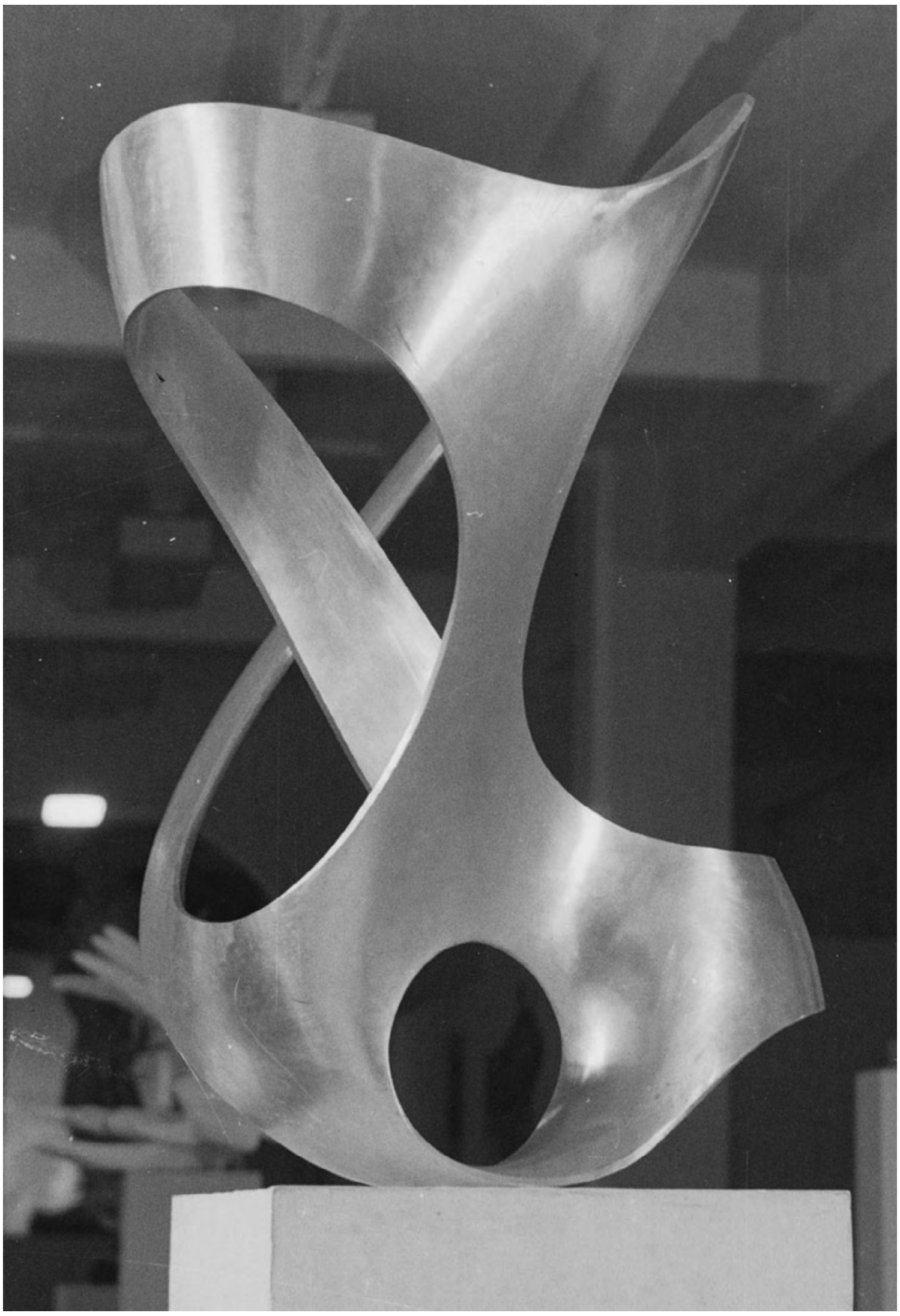

Max Bill, “Unité tripartite”, 1948-49, sculpture, MAC/USP, São Paulo, Brasil. Does the “tripartite” relate to the fact that this is the connected sum of three projective planes minus a disc?

# Moebius necklaces

## Blowing up several times

In order to blow up a point on a surface, delete a disc around that point, blow up this disc, and glue the boundary of the Moebius band to the boundary of the complement of the disc.

We have seen how to blow up a point in the plane. The construction can be generalized. Given a point $p$ on a *smooth surface* $S$, one blows up $S$ at $p$, producing another smooth surface $S_p$ and a blowing down map $\Psi_p : S_p \to S$ as before. The inverse image of $p$ is the exceptional divisor $E_p$: its elements are the tangent lines at $p$ and it is identified with the projective line $P^1(T_p(S))$ constructed from the tangent plane $T_p(S)$ of $S$ at $p$. Outside the exceptional divisor, $\Psi_p$ is a diffeomorphism onto $S \setminus \{p\}$.

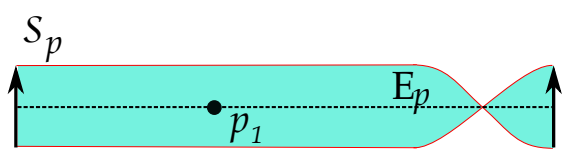

In these pictures one should glue the corresponding arrows. This is a Moebius band.

Let us iterate the process. Choose some point $p_1$ in the exceptional divisor $E_p = \Psi_p^{-1}(p)$ and blow up $S_p$ at $p_1$. The result is a smooth surface $S_{p,p_1}$ and a blowing down map $\Psi_{p_1}$ from $S_{p,p_1}$ to $S_p$ with an exceptional divisor $E_{p_1} \subset S_{p,p_1}$. The inverse image $(\Psi_p \circ \Psi_{p_1})^{-1}(p)$ consists of the union of $E_{p_1}$ and of the strict transform of $E_p$ under $\Psi_{p_1}$. This union is called the exceptional divisor of the composition $\Psi_p \circ \Psi_{p_1} : S_{p,p_1} \to S$. Outside this divisor, the map $\Psi_p \circ \Psi_{p_1}$ is a diffeomorphism onto $S \setminus \{p\}$.

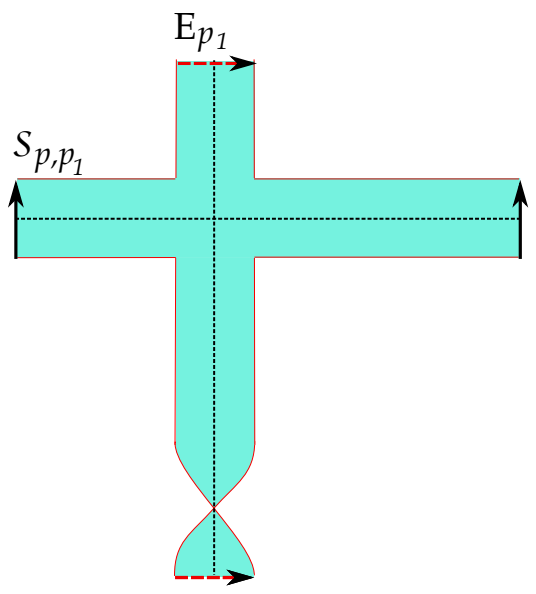

Blowing-up twice.

Choose a point $p_2$ in $(\Psi_p \circ \Psi_{p_1})^{-1}(p)$ and continue the process any finite number of times, “at pleasure”.

The final result is:

- a smooth surface $\overline{S}$,
- a smooth map $\overline{\Psi} : \overline{S} \to S$, which sends diffeomorphically $\overline{\Psi}^{-1}(S \setminus \{p\})$ to $S \setminus \{p\}$.

Note in particular that the boundary of $\overline{S}$ is connected. The inverse image $\overline{\Psi}^{-1}(p)$ is the *exceptional divisor*. It is a finite union of smooth circles embedded in $\overline{S}$. Any two of these circles are either disjoint or intersect transversally in a single point. Three different circles don't intersect. The picture is reminiscent of the Olympic games logo. The difference is that the olympic rings are disjoint, unlike in our situation where some circles do intersect.

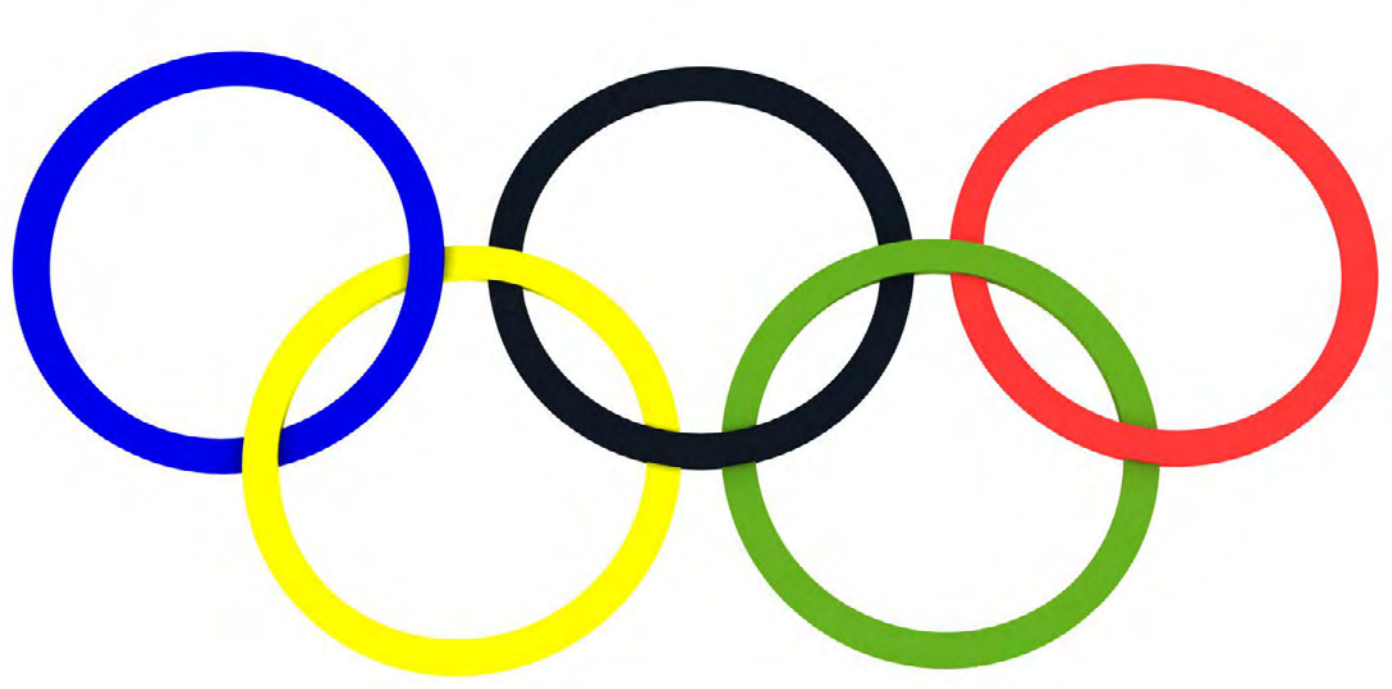

This composition of blowing ups is the multi-lens microscope that we will use and that will enable us to analyze all singular points.

## *The microscope*

Before we use our microscope, let us examine it. If we start with a disc $S$, the single step blowing up $\overline{S}$ is a Moebius band. We will illustrate the topology of $\overline{S}$ in the general case of a finite number of blowing ups.

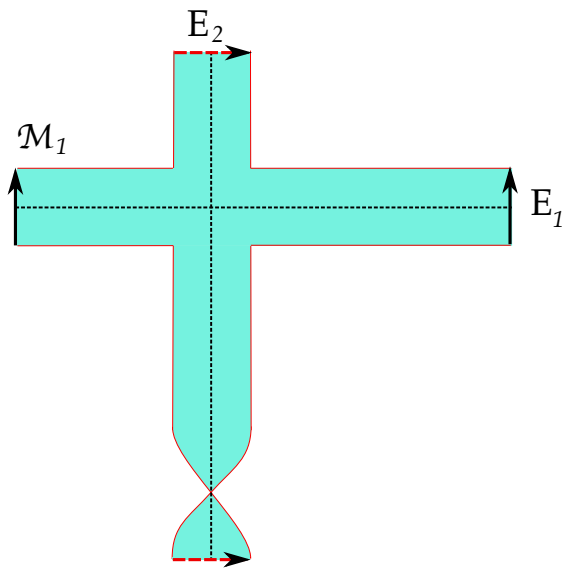

Let us begin with the two step blowing up. We have to picture the result of blowing up a point in a Moebius band. Start with a Moebius band $\mathcal{M}$ containing the exceptional divisor $E$ as its core. As before, pick a point $p_1$ in $E$ and blow up $\mathcal{M}$ at $p_1$. The result will be a surface $\mathcal{M}_1$ containing two circles $E_1$ and $E_2$ intersecting in one point. Here $E_1$ is the strict transform of $E$ and $E_2$ is the exceptional divisor of the second blow-up.

Let $\gamma$ be a loop on some surface $S$. Start with some orientation of the tangent space of $S$ at $\gamma(0)$ and follow it along $\gamma$. When the loop comes back to its origin, the orientation is either the original one or has been reversed. Accordingly, I will say that $\gamma$ is *orienting* or *disorienting*. Formally, this defines a homomorphism from the fundamental group of $S$ (or its first homology) to $\mathbb{Z}/2\mathbb{Z}$.

Have a look at the logo of the Olympic games in Rio de Janeiro and its orienting/disorienting loops.

Coming back to our blown up Moebius band $\mathcal{M}_1$ we will see that $E_1$ is orienting and $E_2$ is disorienting.

The fact that $E_2$ is disorienting should be clear. When we blew up $\mathcal{M}$ at $p_1$, we introduced a Moebius band whose core is $E_2$. As for $E_1$, it is the strict transform of the core $E$ of $\mathcal{M}$. Clearly $E$ is disorienting in $\mathcal{M}$ but this does not imply that its strict transform is disorienting as well. Quite the contrary, as we will see.

To construct $\mathcal{M}_1$, we dig a hole in the original $\mathcal{M}$ and glue another Moebius band to its boundary. Since we are dealing with topology, we can dig a "square hole".

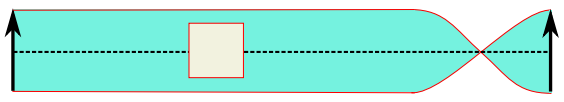

Deleting a small annular neighborhood of the boundary of $\mathcal{M}$, we can even imagine that the square hole crosses $\mathcal{M}$ "from side to side" (not forgetting that a Moebius band has only one boundary circle). In this case, the complement of the square in $\mathcal{M}$ is another square. So the construction of $\mathcal{M}_1$ can be done in another way. Start with a Moebius band, choose two disjoint intervals on its boundary, and glue the two opposite sides of a square to these two intervals. There remains a question. The gluing of the two sides can be done in two ways: with or without a twist.

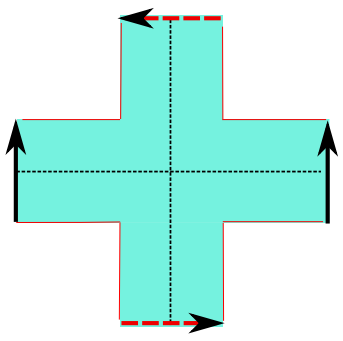

The previous construction can be visualized in the following manner. Consider a cross. Glue the top and bottom sides with a twist, so that the vertical part of the cross becomes a Moebius band. The vertical axis is the disorienting curve $E_2$.

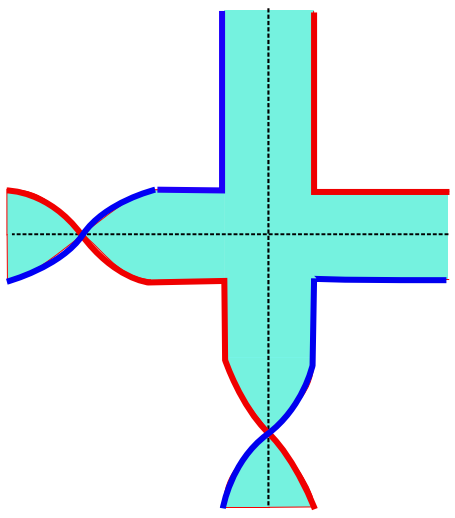

Now, we have to glue the left and right sides of the cross and we have to decide whether they should be twisted or not.

Let us try first with a twist. The boundary of the resulting surface is not connected: it cannot be our surface $\mathcal{M}_1$. Therefore,

the two sides have to be glued without twisting: the curve $E_1$ is indeed orienting.

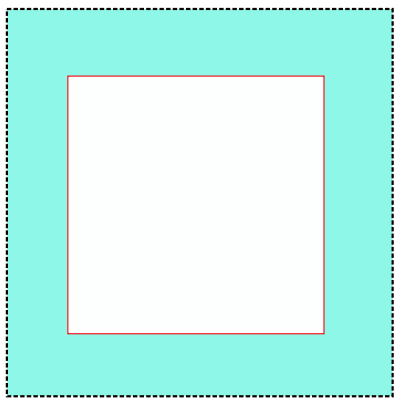

We have obtained a good picture for $\mathcal{M}_1$. A friend recommended that I show the best picture of $\mathcal{M}_1$ ☺! I drew in red and blue the exceptional divisor.

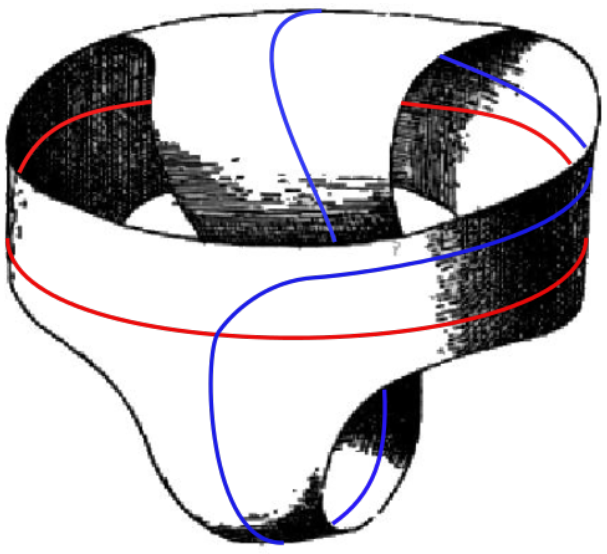

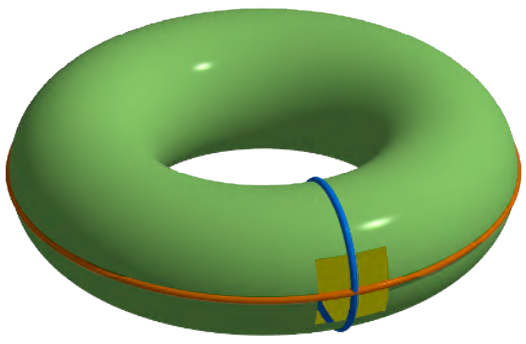

©

## *Interlaced hearts*

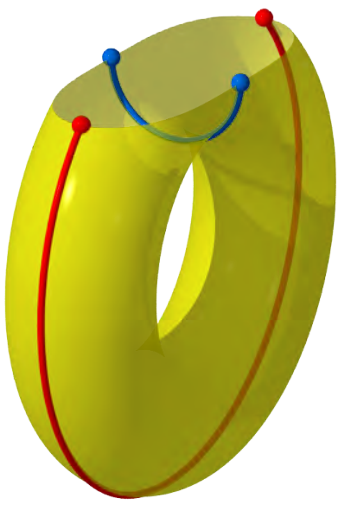

©

I encourage the reader to practice the following topological tricks.

Start with a cross, glue the opposite sides with no twist. Cut open the resulting surface along the central cross, in other words along the two circles. The result is a square frame. Would you have guessed it? Imagine a torus in 3-space and dig a square in it. Then cut it open along a meridian and a parallel. Clearly, what is left is a square with a square hole: a square frame.

Amazingly, this example of a cross with identified opposite sides has already been considered by Gauss under the name *Doppelring*. In his remarkable paper *Gauss als Geometer*, Stäckel[74] relates a conversation between Gauss and Moebius. Gauss observes that the Doppelring has a connected boundary. More interestingly, he observes that there are two *disjoint* arcs connecting two *linked* pairs of points on the boundary. I recall that the impossibility of such a configuration in a disc was the crucial point in his proof of the fundamental theorem of algebra.

[74] Stäckel. Literaturberichte: Materialien für eine wissenschaftliche Biographie von Gauß. *Monatsh. Math. Phys.*, 32(1):A5, 1922. Gesammelt von F. Klein, M. Brendel und L. Schlesinger. Heft V: C. F. Gauß als Geometer. In Kommission bei B. G. Teubner in Leipzig, 1918.

Take again the same cross and glue the opposite sides, now with a twist. Cut open the resulting surface along the two circles. The result is...

Well, it depends how you twisted it. As an abstract surface it is well defined: it consists of two connected components, each homeomorphic to an annulus. However, the way it is embedded in space depends on the twisting. Experiment!

The most impressive result occurs when both opposite sides of the cross are twisted, but in a different manner, left and right so to speak. This produces two interlaced hearts.

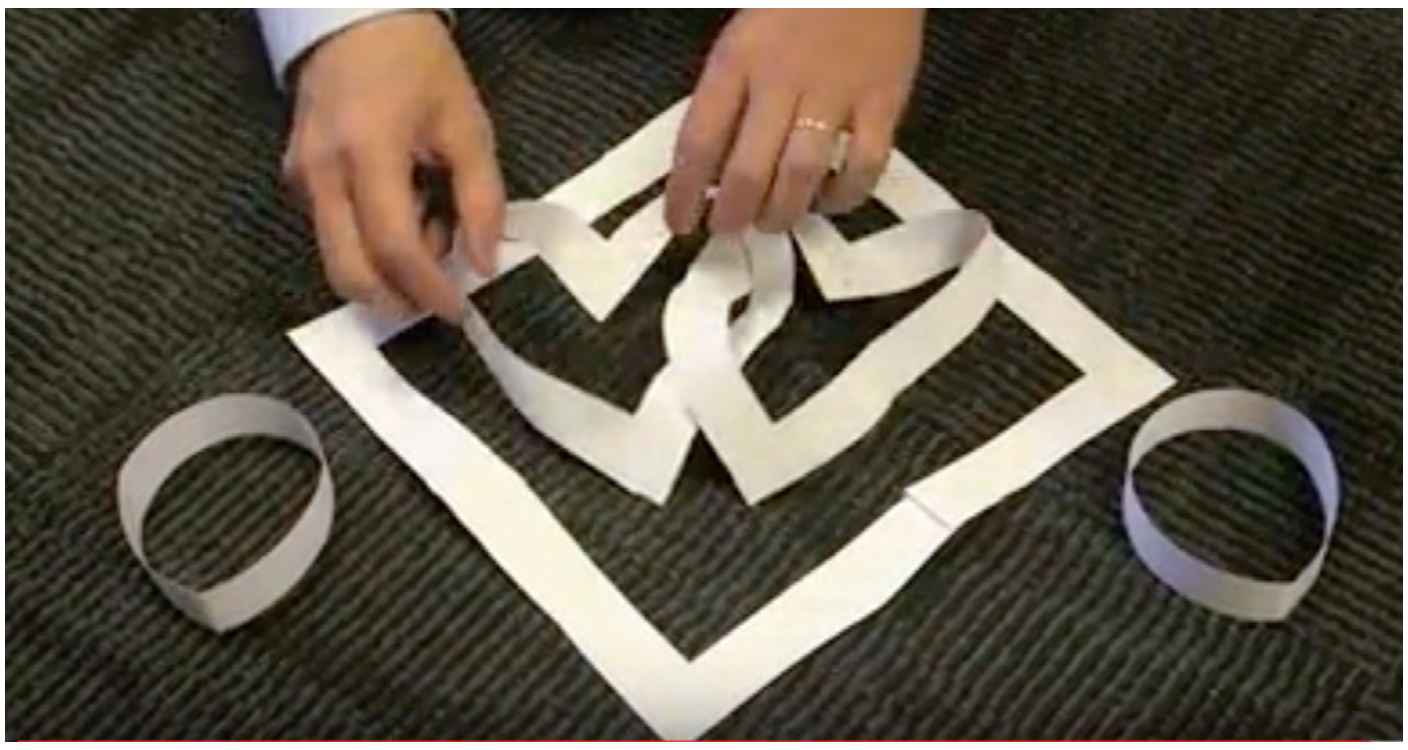

On this topic, the reader must look at Tadashi Tokieda's presentation Unexpected shapes on Youtube (in two parts). ©

Finally, gluing two of the opposite sides with a twist and not the other two, we get our blown up Moebius band, with a connected boundary. What happens when we cut it open along the two exceptional circles? This is easy since the blowing up is a homeomorphism outside the exceptional divisor. We get something homeomorphic to a punctured disc. Indeed, we also get a square frame. Also, the way this frame is embedded in space depends on the twisting. Practice these topomagical tricks!

## *Blowing up more points*

I now describe the situation when more points are blown up.

It is easy to describe the topology of the resulting surface. Blowing up a point amounts to digging a hole in the surface and to gluing a Moebius band on the boundary. Said differently, the operation of blowing up is equivalent to the *connected sum* with a projective plane.

Given two connected surfaces $M_1$ and $M_2$, their connected sum $M_1 \sharp M_2$ is obtained by deleting a disc from each of them and gluing them along the newly created boundary. Since we already observed that when we delete a disc in a projective plane we get a Moebius band, the topological effect of a blowing up is the connected sum with a projective plane.

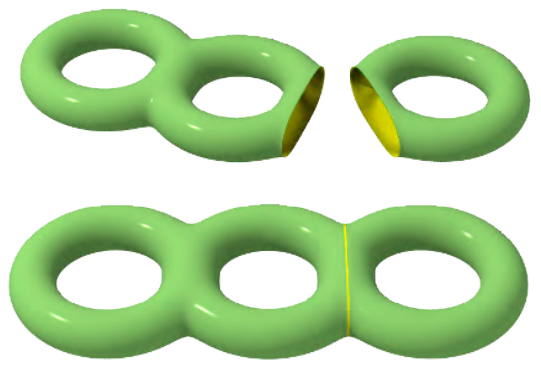

There is a subtle orientation question here. The two punctured surfaces could be glued in two different ways since the boundary circle has two orientations. However, orientable surfaces do have orientation reversing homeomorphisms. Check that this implies that the connected sum is indeed well defined among connected non-oriented surfaces, orientable or not.

Therefore, if we blow up a disc $k$-times in a row, the resulting surface is the connected sum of $k$ projective planes, minus a disc. Recall that any compact non orientable surface with a single boundary component is homeomorphic to such a surface, and that the number $1-k$ is known as the Euler-Poincaré characteristic of the surface. See for instance the books[75] and [76]. However, this is only a partial description of the result since we still have to describe the position and nature of the exceptional divisor. The latter does not only depend on $k$ but also on the choices of the $k$ successive points that have been blown up.

Look at Max Bill's beautiful sculpture illustrating this chapter. A paragraph of Ton Marar's paper[77] is dedicated to showing that this sculpture represents a connected sum of *three projective* planes (minus a disc). This is explained in the following pictures, extracted from this paper.

[75] S. Barr. *Experiments in topology*. Dover Publications, Inc., Mineola, NY, 1989. Reprint of the 1964 original.

[76] V. G. Boltyanskiĭ and V. A. Efremovich. *Intuitive combinatorial topology*. Universitext. Springer-Verlag, New York, 2001.

[77] T. Marar. Aspectos topológicos na arte concreta, 2004. II Bienal da Sociedade Brasileira de Matemática, Salvador, Universidade Federal da Bahia.

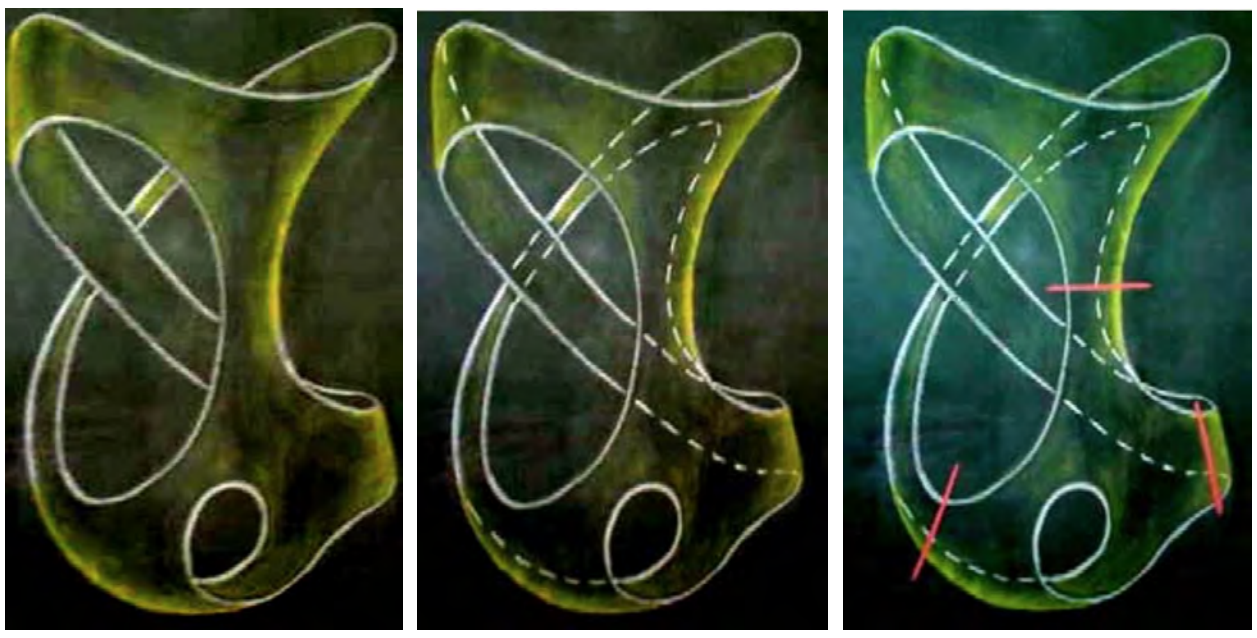

The same paper contains another version of the same surface, inspired from the already mentioned book by Francis (page 101).

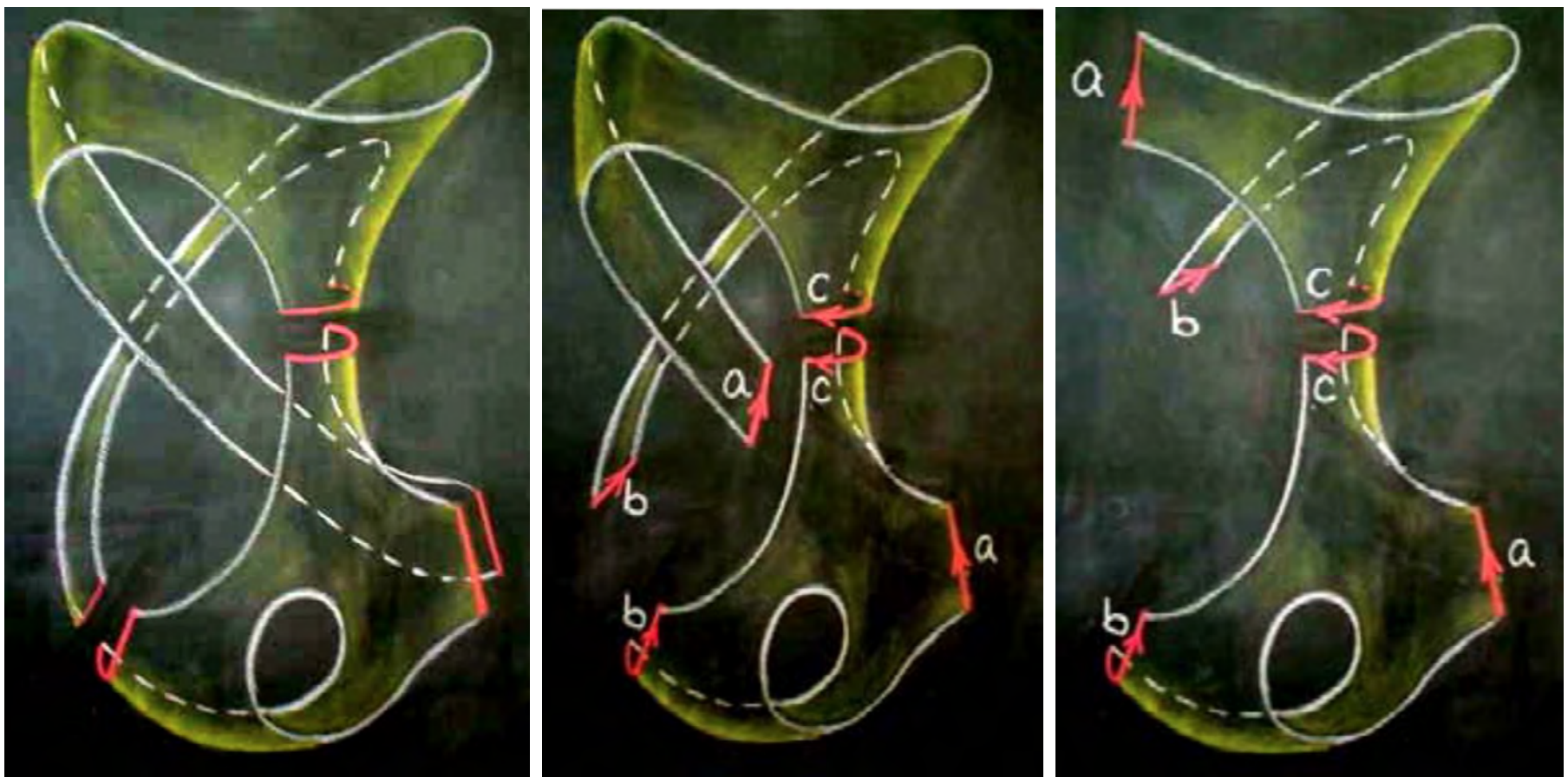

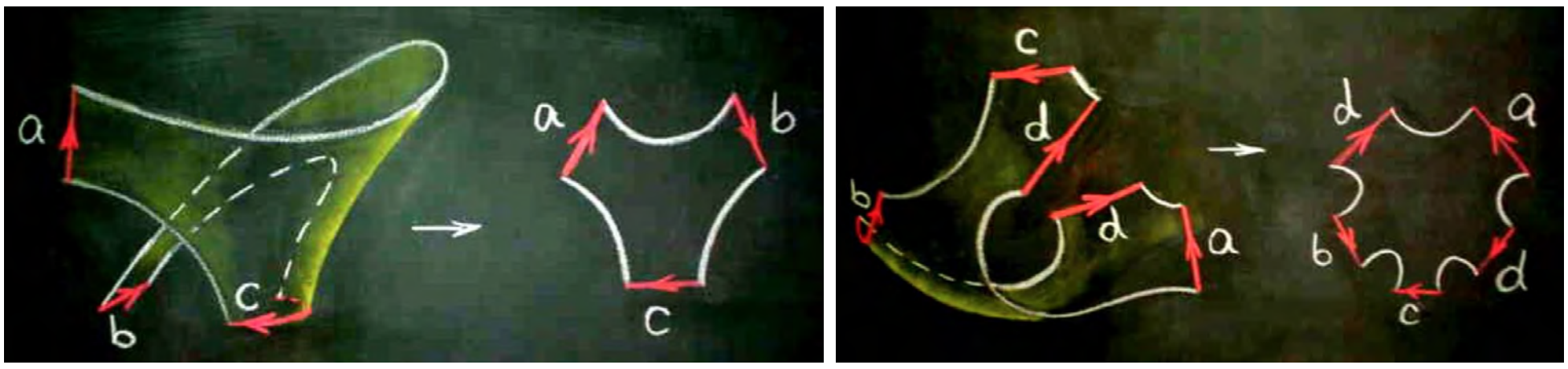

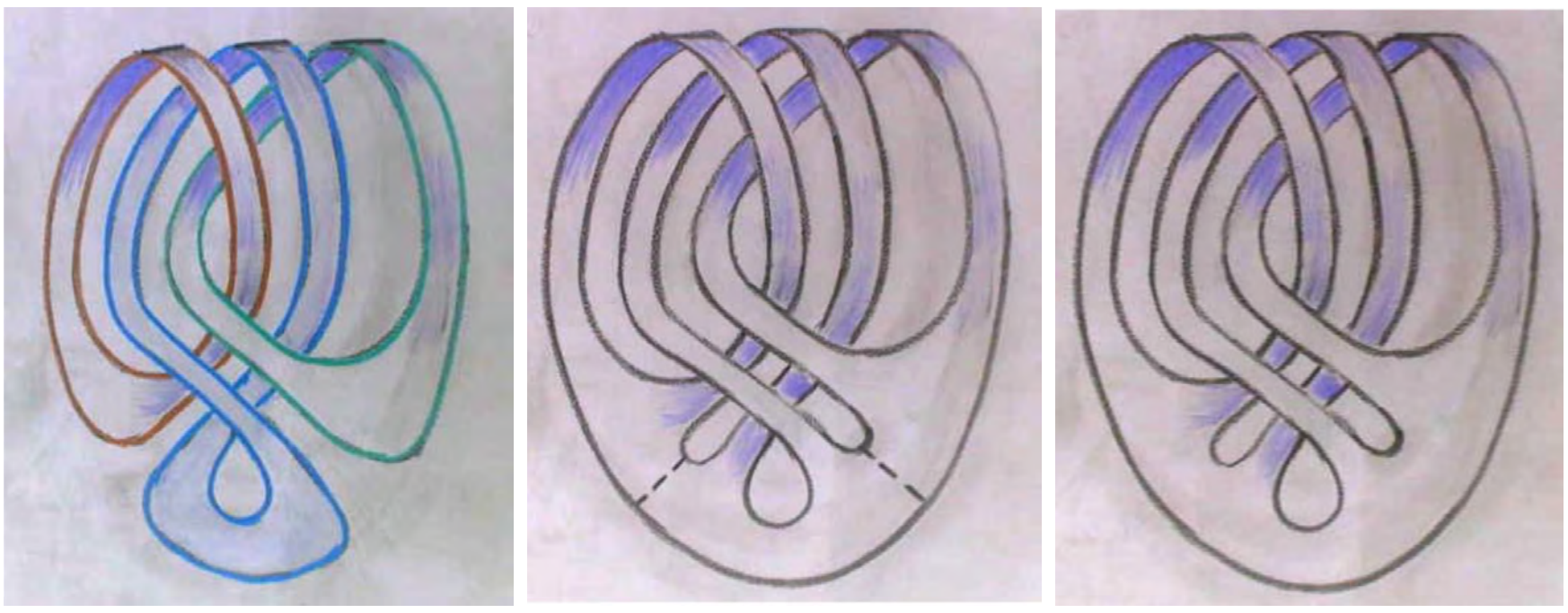

## *Necklaces of divisors*

We still have to describe the topology of the exceptional divisor *inside* the connected sum of projective planes.

At the first step, there is no surprise: the divisor is the core of the Moebius band.

At the second step, we blow up a point of the Moebius band. The case of interest is when we blow up a point on $E_1$ as discussed earlier. Algebraic geometers think of the projective line as a line... and draw it as a line, even though it is homeomorphic to a circle...

To forgive them, I should recall that the projective line over the complex numbers is homeomorphic to a 2-sphere, and to a Cantor set for the $p$-adic numbers.

When we come to the third blowing up, we may choose the point either on $E_1$, or on $E_2$, or at the intersection of $E_1$ and $E_2$. In all cases, the blown up surface is a connected sum of *three projective* planes (minus a disc), that is to say Max Bill's surface. However, the location of the exceptional divisor on this surface is not the same. As an exercise, the reader should try to (mentally) draw the three possible exceptional divisors, directly on the sculpture.

The general situation is now easy to describe. Combinatorially, the many components of the exceptional divisor are organized as a tree.

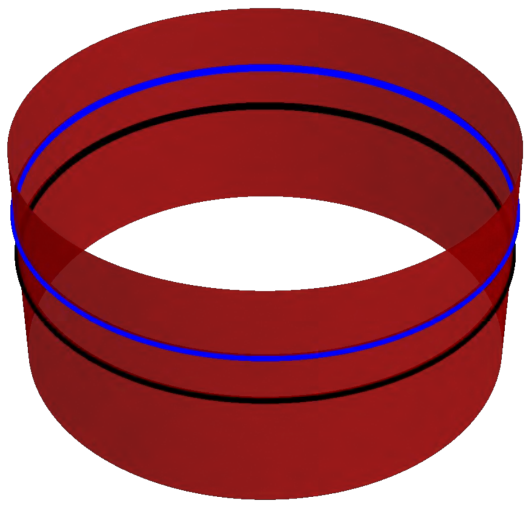

At each blowing up, one new Moebius band is attached on the previous necklace.

However, this changes the orientability of the band on which the new band is attached.

In order to prove that, consider a closed curve $\gamma$ on a surface $S$. Deform slightly $\gamma$ to some $\gamma'$, transversal to $\gamma$, and count the number of intersection points in $\gamma \cap \gamma'$ modulo 2. This is called the *self intersection* of $\gamma$. It is equal to 0 or to 1 according to whether $\gamma$ is orienting or disorienting.

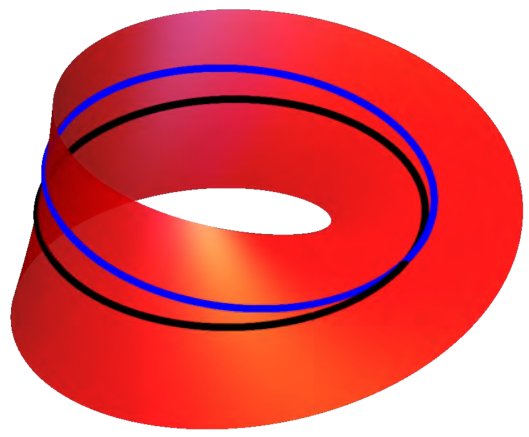

Now, let us blow up some point $p$ of $\gamma \subset S$. Choose $\gamma'$ which also passes through $p$ and let us blow up the picture at $p$. The two strict transforms $\bar{\gamma}$ and $\bar{\gamma}'$ have precisely one less intersection points than $\gamma$ and $\gamma'$ since the tangents at $p$ are different. It follows that the self intersection of $\bar{\gamma}$ is equal to the self intersection of $\gamma$ minus (or plus since we count modulo 2!) one.

Let us work out an example. The picture in the margin illustrates a succession of six blow ups. The thick points represent the centers of the blowing ups. The lines represent the projective lines (do not forget that they are actually circles). The dashed lines represent the new divisors which appear at each step. So the blow down maps, represented by downwards arrows, are collapsing these dotted lines to points of the same color. Double lines represent the components which are orienting. At the end of the process, the exceptional divisor consists of six circles.

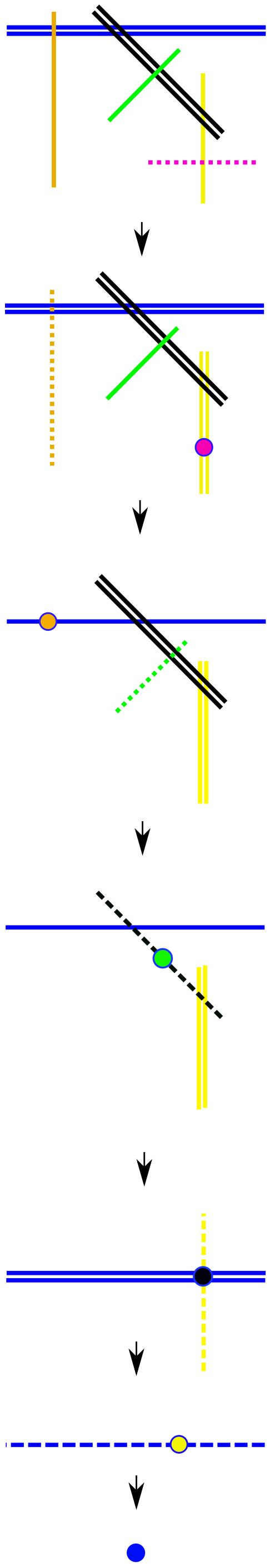

We can now draw the corresponding necklace, made of four Moebius bands and two annuli. The opposite sides of the six strips should be glued as suggested by the picture.

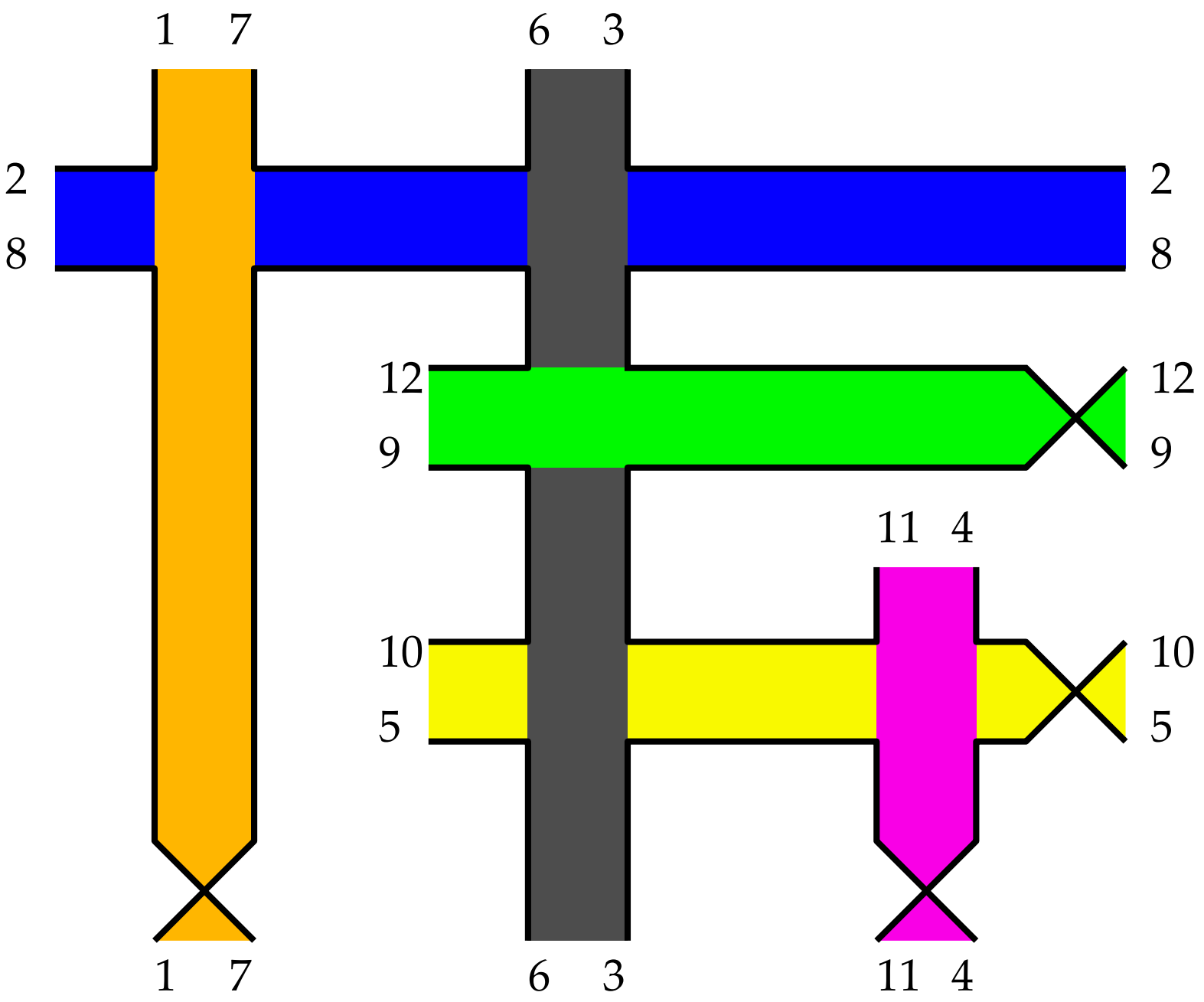

You should check that the boundary is indeed connected, as it should be. Go around the boundary, following the numbers from 1 to 12 and then back to 1.

## *Plumbing*

Here is another view of the previous example.

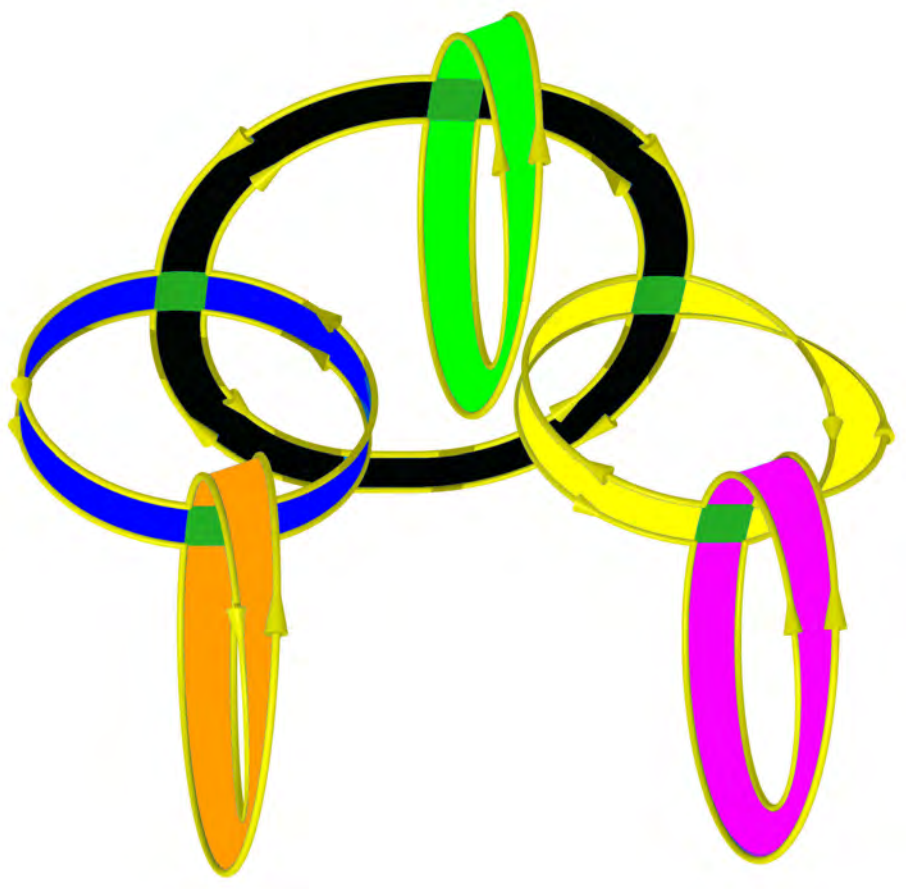

[78] B. Ozbagci and P. Popescu-Pampu. Generalized plumbings and Murasugi sums. *Arnold Math. J.*, 2(1):69–119, 2016.

A topologist would say that this surface is obtained by *plumbing* several Moebius bands and annuli. This operation is very simple. Suppose you have two surfaces $S_1, S_2$ with non-empty boundary. Choose two embeddings $i_1, i_2$ of the square $[-1,1]^2$ in $S_1$ and $S_2$ in such a way that the images $i_1(\{\pm 1\} \times [-1,1])$ and $i_2([-1,1] \times \{\pm 1\})$ lie in the boundary of $S_1$ and $S_2$ respectively. Now, for each $(x,y) \in [-1,1]^2$, identify $i_1(x,y)$ and $i_2(x,y)$. The result is the plumbing of $S_1$ and $S_2$ along $i_1, i_2$. This is a surface with boundary (and corners that can be smoothed easily). See[78] for a presentation of some variations around this construction.

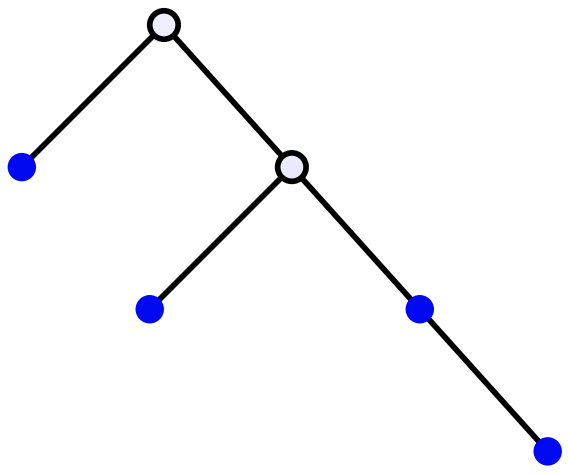

Blue dots correspond to Moebius bands and white ones to annuli.

Let us start now with a rooted planar tree. For each node, take an annulus or a Moebius band. Now, plumb all these bands together according to the blueprint given by the tree. Each band is plumbed to all the bands associated to its children in the tree, as in the picture in the margin. Note that the annulus and the Moebius band admit four homeomorphisms permuting opposite sides of the square so that the operation of plumbing such a band is well-defined. The final result of this plumbing is a surface $S$ with boundary.

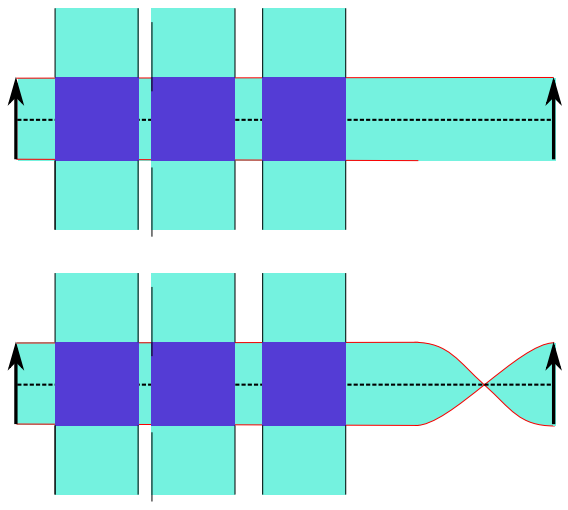

Each band (annulus or Moebius band) contains a circle as its core. The union of these circles defines a graph $E \subset S$ that we can call the divisor, even though our $S$ has not necessarily been constructed by a sequence of blowing ups. There is a projection $\pi$ of $S$ on $E$ such that the inverse image $\pi^{-1}(x)$ consists of one arc if $x$ is a regular point of $E$ and two intersecting arcs otherwise. Let us denote by $S/E$ the topological space obtained by collapsing $E$ to a single point. If $S$ is the result of a sequence of blowing ups, we know that $S/E$ is a closed disc and the projection of $S$ to $S/E$ is the blowing down map.

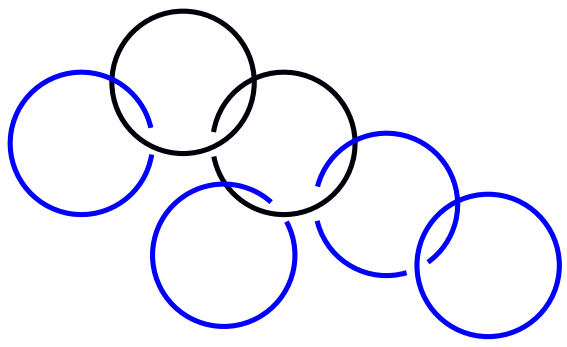

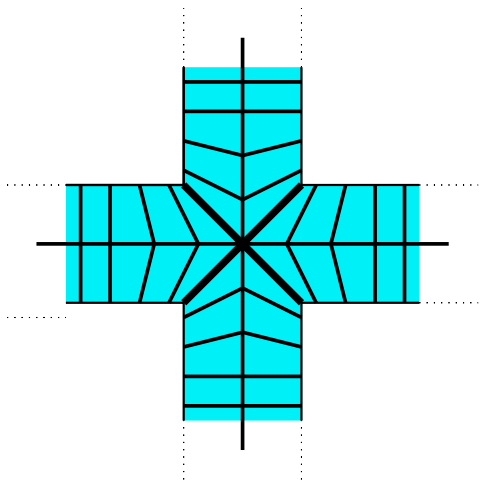

The fibers of the projection $\pi : S \to E$.
Do not confuse $\pi$ from $S$ to the divisor $E$ with the blowing down map $\Psi$ from $S$ to $S/E$.

*Exercise*: Show that the space $S/E$ is homeomorphic to a cone whose basis is the disjoint union of $k$ circles, where $k$ is the number of connected components of the boundary of $S$.

In particular, *the quotient space $S/E$ is homeomorphic to a disc if and only if the boundary of $S$ is connected.*

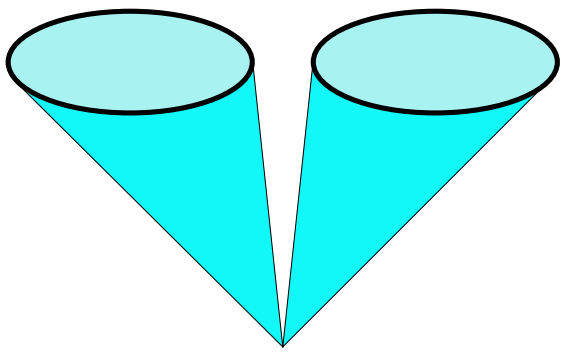

The cone on two circles.

The following exercise gives a simple criterion enabling us to check directly from the blueprint if the boundary of $S$ is connected. This is easier than drawing the picture on a sheet of paper and following carefully the boundary. The solution of this exercise requires some understanding of the homology of surfaces. Suppose the tree has $n$ nodes. Consider the symmetric $n \times n$ matrix $A$, with coefficients in $\mathbb{Z}/2\mathbb{Z}$, defined in the following way. Set $a_{ii} = 0$ if the node $i$ is an annulus and $a_{ii} = 1$ if it is a Moebius band. If $i \neq j$, set $a_{ij} = 1$ if the nodes $i, j$ are adjacent in the tree and 0 otherwise.

*Exercise*: Show that the boundary of $S$ is connected if and only if the matrix $A$ is invertible (over $\mathbb{Z}/2\mathbb{Z}$).

*Hint:* Check the following:

- The injection $E \subset S$ and the projection $\pi : S \to E$ induce inverse isomorphisms between $H_1(S, \mathbb{Z}/2\mathbb{Z})$ and $H_1(E, \mathbb{Z}/2\mathbb{Z})$.
- A basis of $H_1(E, \mathbb{Z}/2\mathbb{Z})$ is given by the cores of the $n$ bands.
- The symmetric intersection form on $H_1(E, \mathbb{Z}/2\mathbb{Z})$ is given by the matrix $A$.
- The kernel of the intersection form is the image of $H_1(\partial S, \mathbb{Z}/2\mathbb{Z})$ in $H_1(S, \mathbb{Z}/2\mathbb{Z})$.

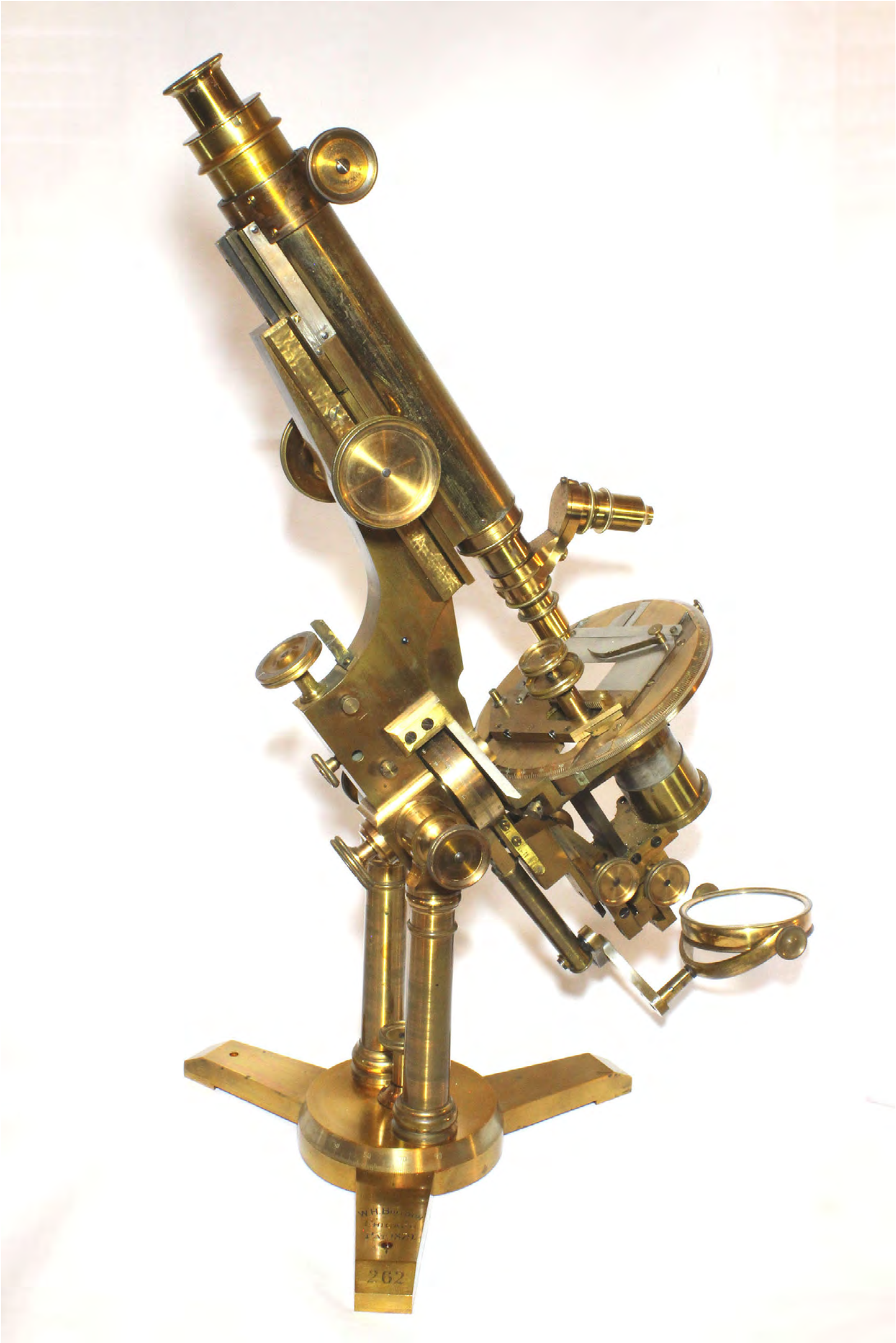

An 1882 microscope. ©

# *Resolution of singularities*

We now use our microscope to analyze the nature of singularities and prove a theorem which is essentially due to Max Noether[79].

[79] M. Noether. Rationale Ausführungen der Operationen in der Theorie der algebraischen Funktionen. *Math. Ann.*, 23:311–358, 1883.

## *Blowing up a branch*

Consider some singular point of a real analytic plane curve defined by some equation $F(x, y) = 0$.

Suppose that we have found a real branch of this curve, i.e. a solution of the form

$$x = \pm t^m \quad ; \quad y = \sum_{k \geq 1} a_k t^k.$$

Let us look at the set $I \subset \mathbb{N} \subset \mathbb{Z}$ of integers $k$ such that $a_k \neq 0$. We can always assume that the greatest common divisor of the elements of $I$ together with $m$ is 1. In other words, the subgroup of $\mathbb{Z}$ generated by $m$ and $I$ is $\mathbb{Z}$.

Let $\mu \geq 1$ be the smallest integer such that $a_\mu \neq 0$.

If $\mu < m$, the series $y/x$ "tends to infinity" as $t$ tends to 0, which means geometrically that the vertical axis $x = 0$ is tangent to the branch at the origin.

If $\mu > m$, the series $y/x$ "tends to 0" as $t$ tends to 0, which means geometrically that the horizontal axis $y = 0$ is tangent to the branch at the origin.

If $\mu = m$, the tangent at the origin is the line $y = a_m x$.

Until now, following Newton, we looked at $y$ as a "function"

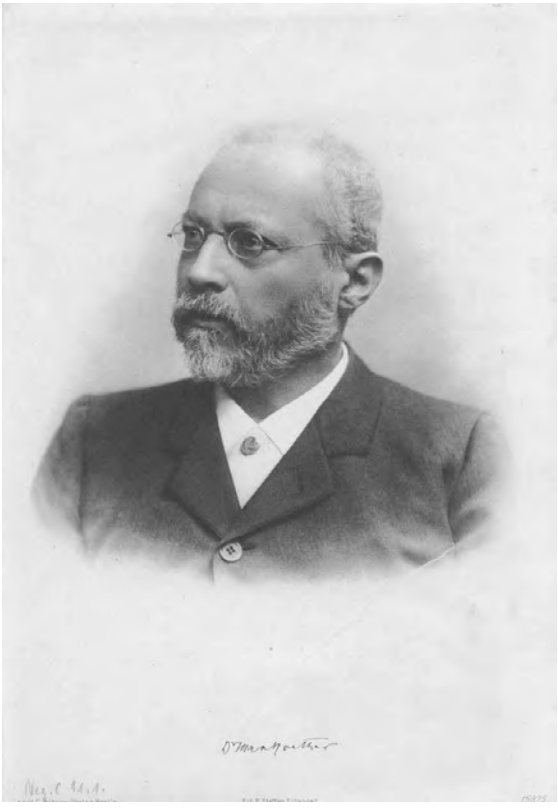

Max Noether (1844–1921). ©

of $x$. We are now more interested by the *curve* $F(x,y) = 0$ so that we can permute the roles of $x$ and $y$.

Hence we can always assume that $\mu \geq m$. Indeed, we can define $\tau$ as some $\mu$-th root of $\pm y = \pm \sum_{k \geq \mu} a_k t^k$ so that $\tau$ is a power series in $t$. Reverse the role of $x$ and $y$ so that we have now $y = \tau^\mu$ and $x$ is a series in integral powers of $\tau$.

If our branch is singular, i.e. if $m > 1$, proceed as follows.

1. Let $\beta_1$ be the smallest integer in $I$ which is not a multiple of $m$.

2. Let $\beta_2$ be the smallest integer in $I$ which is not in the group generated by $m$ and $\beta_1$.

Continue in this way until you obtain a family of integers in $I$ generating $\mathbb{Z}$. This defines a finite sequence of integers $m < \beta_1 < \beta_2 < \ldots < \beta_g$. This list is the *Puiseux characteristic of the branch*. Once again, Puiseux is not responsible for this definition, which was introduced later by Halphen and Smith[80]. Poor Puiseux!

Let us look at the effect of a blowing up on our branch. Recall that in practice this amounts to looking at the coordinates $(x, y_1)$ where $y_1 = y/x$ is the slope of the line passing through the origin and $(x,y)$. In these coordinates $(x,y_1)$, we have:

$$x = \pm t^m \quad ; \quad y_1 = \sum_{k \geq \mu} a_k t^{k-m}.$$

The Euclidean division of $\beta_1$ by $m$ gives

$$\beta_1 = mq + m_1 \quad \text{with} \quad (0 < m_1 < m)$$

and we have

$$x = \pm t^m \quad ; \quad y_1 = a_m + a_{2m} t^m + \cdots + a_{qm} t^{(q-1)m} + \sum_{k \geq \beta_1} a_k t^{k-m}.$$

Now, we translate in order to bring back the singularity to the origin. Said differently, we set $y_2 = y_1 - a_m$

$$x = \pm t^m \quad ; \quad y_2 = a_{2m} t^m + \cdots + a_{qm} t^{(q-1)m} + \sum_{k \geq \beta_1} a_k t^{k-m}.$$

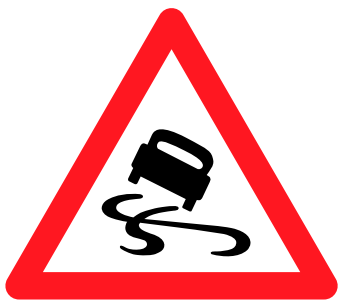

©

This procedure might look complicated. It is very similar to *Euclid's algorithm*. Given two positive integers $0 < a \leq b$, one subtracts $a$ to $b$ so that one has now $a, b-a$. If $0 < a \leq b-a$, one continues: $a, b-2a$. Continue while the first integer is smaller than the second. This is nothing more than the Euclidean division of $b$ by $a$. Then, permute the two integers and continue the process. The algorithm finishes after a finite number of steps, when the second integer is equal to 0. At this final step, the first number is the g.c.d of $a$ and $b$. For example $(6,9) \to (6,3) \to (3,6) \to (3,3) \to (3,0)$. In our more complicated situation, we proceed in the same way, blowing up as many times as necessary until we can permute the roles of $x$ and $y$ and continue...

[80] E. R. Garcí a Barroso, P. D. González Pérez, and P. Popescu-Pampu. Variations on inversion theorems for Newton–Puiseux series. *Math. Ann.*, 368(3-4):1359–1397, 2017.

and we can blow up and translate again if $q \geq 2$. After $q$ steps, this gives

$$x = \pm t^m \quad ; \quad y_{2q} = \sum_{k \geq \beta_1} a_k t^{k-qm}.$$

Since $\beta_1 - qm = m_1 < m$, the vertical axis is tangent to this last curve at the origin. As before, we permute the role of the two coordinates so that

$$y_{2q} = \pm \tau^{m_1} \quad ; \quad x = \sum_{k \geq 1} b_k \tau^k.$$

In other words, after a certain number of blowing ups, the new curve has multiplicity $m_1 < m$. Continuing in this way, after a finite number of steps, the curve is smooth.

We have therefore proved the following.

**Theorem.** *Let C be a branch of some analytic curve* $F(x,y) = 0$ *in the neighborhood of the origin. Then the strict transform of C by a suitable succession of blowing ups is a smooth curve.*

Strictly speaking, we have only defined blowing up maps for real surfaces. Similar maps could also be defined over the complex numbers. In this case, the result of a succession of blow ups is a holomorphic surface (i.e. of complex dimension 2), containing an exceptional divisor which is now a union of complex projective lines intersecting transversally. These blowing up operations could even be defined in the context of algebraic surfaces over any field. Most of this chapter would adapt verbatim to this general case.

In dimension 3 or more, the desingularization is much more subtle and blowing up points is not sufficient. Hironaka proved in 1964 that any algebraic variety over a field of zero characteristic can be desingularized. I have always been told that this proof is a *tour de force* which is very hard to digest. However, in his 2007 lectures on the resolution of singularities (Princeton University Press), J. Kollár writes that "The lingering perception that the proof of resolution is very hard gradually diverged from reality. ... it is feasible to prove resolution in the last two weeks of a beginning algebraic geometry course." So, why don't you try to read Kollár?

## *Blowing up all branches*

In the neighborhood of the origin, a curve $F(x,y) = 0$ consists of several branches. We learned how to desingularize each of these branches, but the many smooth curves that we obtain may be in a rather complicated relative situation. More blowing ups are necessary to untangle the strings. Using the previous theorem, we can desingularize all the branches, one by one. We have

– a blowing down map $\Psi$ from some surface $S$ to a neighborhood of the origin,

– an exceptional divisor $E \subset S$ mapped to the origin by $\Psi$,

so that the strict transform of our curve is the union of a certain number of smooth curves. Each of these curves intersects the exceptional divisor in a single point.

If all these points are distinct, our job is finished: our singular curve $F = 0$ has been "desingularized" as a union of disjoint smooth curves.

The only task that we still have to do is to deal with a certain number of *smooth curves* passing through the same point $p$ in the exceptional divisor. Note that some of these curves may be tangent to the divisor.

The point $p$ might belong to one or two components of $E$. Add this or these components to the list of our smooth curves passing through $p$. Choose local coordinates $(x, y)$ in the neighborhood of $p$ such that the equations of the smooth curves are $y = f_i(x)$ where the $f_i$'s are distinct convergent power series.

Then, blow up again, introducing a new projective line. The strict transforms of the curves will remain smooth and will intersect the new component of the divisor in a point corresponding to the derivatives of the $f_i$'s at the origin. Two curves might be tangent at the origin, i.e. two $f_i$'s might have the same derivative at 0, but we may blow up again. This process will separate the $f_i$'s by some of their Taylor polynomials. Eventually, the result is a collection of smooth curves which are disjoint and transverse to the exceptional divisor. We can even assume that the final collection of curves only intersect the divisor in regular points, i.e. not at the intersection of two circles.

This is Noether's theorem.

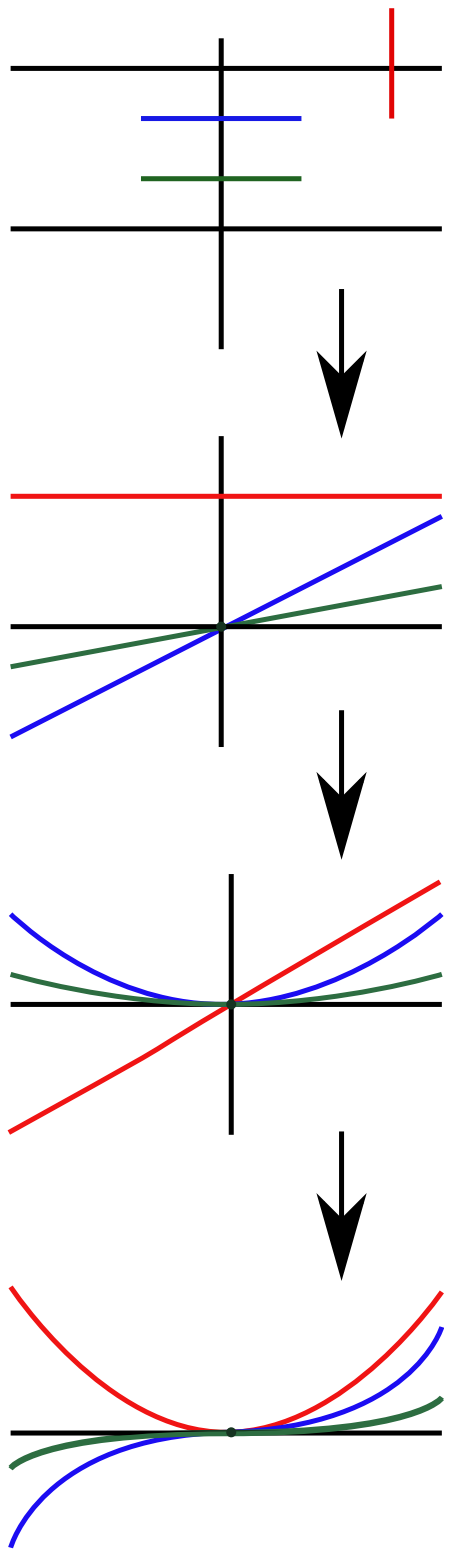

The red, blue and green curves are tangent to the black divisor. Three blowing ups make them transversal to the (new) divisor.

**Theorem.** *Let $C$ be some analytic curve in the neighborhood of the origin. Then the strict transform of $C$ under a suitable succession of blowing ups is a disjoint union of smooth curves transverse to the exceptional divisor.*

## *Quadratic transforms*

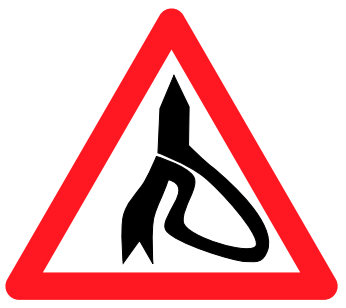

This is a detour. ©

Max Noether was working in the global context of *algebraic curves* and not of local singularities of analytic curves. His microscope was slightly different and is called the *quadratic transform*.

Let me introduce the *Cremona* group of the projective plane $P^2(K)$ over some field $K$[81]. It consists of $K$-*automorphisms* of the field $K(x, y)$ of rational functions in two variables. Such an automorphism is completely defined by two rational functions $f(x, y), g(x, y)$ which are the images of $x$ and $y$. Two functions $f, g$ define an element of the Cremona group if the transformation $(x, y) \dashrightarrow (f(x, y), g(x, y))$ is *birational*.

Projective transformations defined by elements of $PGL(3, K)$

[81] S. Cantat. The Cremona group in two variables. In *European Congress of Mathematics*, pages 211–225. Eur. Math. Soc., Zürich, 2013.

are birational isomorphisms of the plane, but the Cremona group is much bigger. A good example is the quadratic involution

$$\sigma : (x,y) \dashrightarrow (1/x, 1/y)$$

that can also be seen in homogeneous coordinates $[x : y : z]$ as

$$\sigma : [x : y : z] \dashrightarrow [yz : zx : xy].$$

It collapses the line containing two of the three points $[1 : 0 : 0]$, $[0 : 1 : 0]$ and $[0 : 0 : 1]$ to the third point. The involution is not defined at these points. Away from the three lines, $\sigma$ is a bijection, and even an involution. Note also that $[1 : 1 : 1]$ is a fixed point of $\sigma$.

If $A, B, C, M$ is a projective basis in $P^2(K)$ (no three are collinear), there is a projective transformation $\phi$ sending them to $[1 : 0 : 0]$, $[0 : 1 : 0]$, $[0 : 0 : 1]$ and $[1 : 1 : 1]$. The conjugate $\phi^{-1} \circ \sigma \circ \phi$ is the *quadratic transform* associated to the triangle $A, B, C$ (and fixed point $M$).

Max Noether used these maps instead of the blowing ups.

The advantage is that all the process is done in the projective plane without having to introduce a new surface. The drawback is that $\sigma$ is blowing up and down at the same time. It collapses lines and blows up points, so that while resolving some singularities, it creates new ones.

Noether claimed that $\sigma$ and $PGL(3, K)$ generate the full Cremona group. This is indeed true but his proof was not correct.

Start with an algebraic curve defined by some polynomial equation $P(x,y) = 0$. Choose a singular point $A$ and select two points $B, C$ which are not in the curve and such that lines $AB, BC, CA$ intersect the algebraic curve transversally (except in $A$ of course). Then, consider the image of the curve by a quadratic transform associated to $A, B, C$ (and some fixed point $M$ which is not relevant). The point $A$ is blown up. The other intersections of the curve with $AB, BC, CA$ produce smooth branches intersecting transversally. So, we can blow up as many times as necessary, at the cost of introducing multiple points where smooth curves intersect transversally.

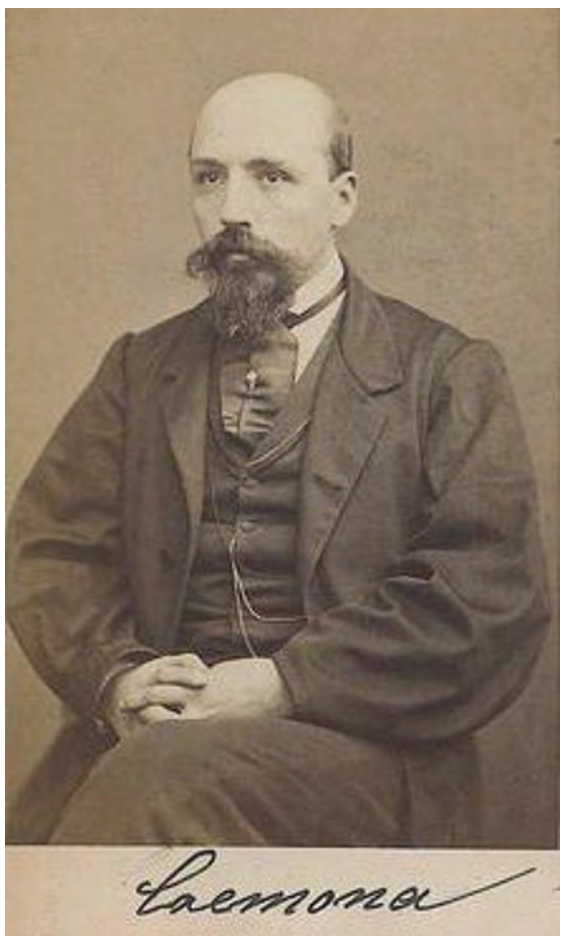

Luigi Cremona (1830–1903). ©

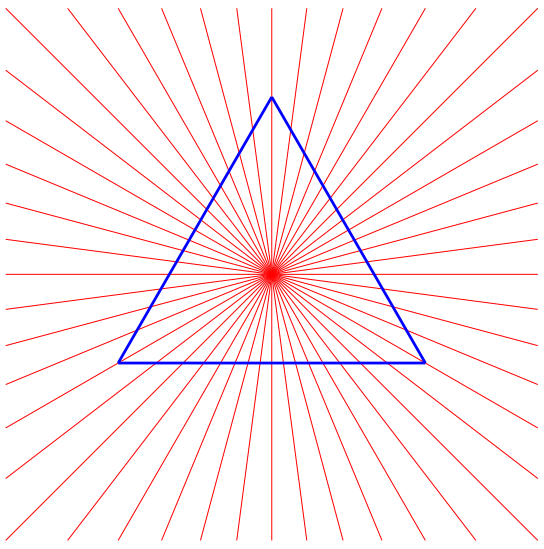

The quadratic involution maps the pencil of lines through $M$ to the pencil of conics passing through $A, B, C, M$.

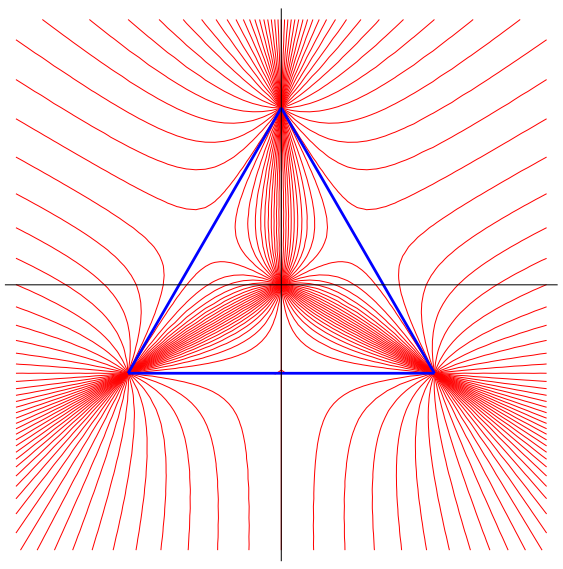

This is the way Noether expressed his theorem:

**Theorem.** *Any algebraic curve can be transformed under a suitable Cremona automorphism into another curve whose only singularities are* ordinary, *that is to say, consist of some smooth branches intersecting transversally.*

There is another famous involution in the (real) plane: the *inversion*. Sixty years ago, all secondary school students were familiar with it. Textbooks were full of exercises of the following style: take your favorite theorem in plane geometry, transform it by inversion and produce a new theorem. The definition is very simple. Choose a point $P$ in the Euclidean plane, called the *pole* of the inversion. Every point $Q$ is sent by inversion to the point $Q'$ such that $P, Q, Q'$ are on the same line and such that the product $\overline{PQ} \cdot \overline{PQ'} = 1$. This involution is not defined at the pole $P$, maps circles not containing $P$ to circles, and circles containing $P$ to straight lines not containing $P$. If $P$ is the origin of the complex plane, this is just the transformation $z \in \mathbb{C}^\star \mapsto 1/\overline{z} \in \mathbb{C}^\star$.

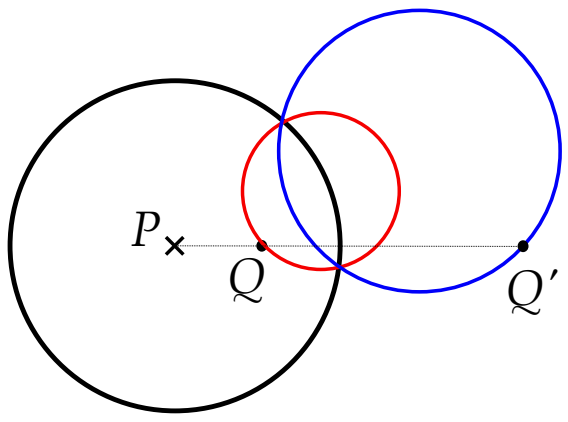

$\overline{PQ} \cdot \overline{PQ'} = 1.$

For instance, the theorem that all French kids call *Chasles' relation* states that for 3 points $Q', R', S'$ on an oriented line, the following holds true: $\overline{Q'R'} + \overline{R'S'} = \overline{Q'S'}$. Transform this by inversion and you discover Ptolemy's theorem: "Let a convex quadrilateral $PQRS$ be inscribed in a circle. Then the sum of the products of the two pairs of opposite sides equals the product of its two diagonals."

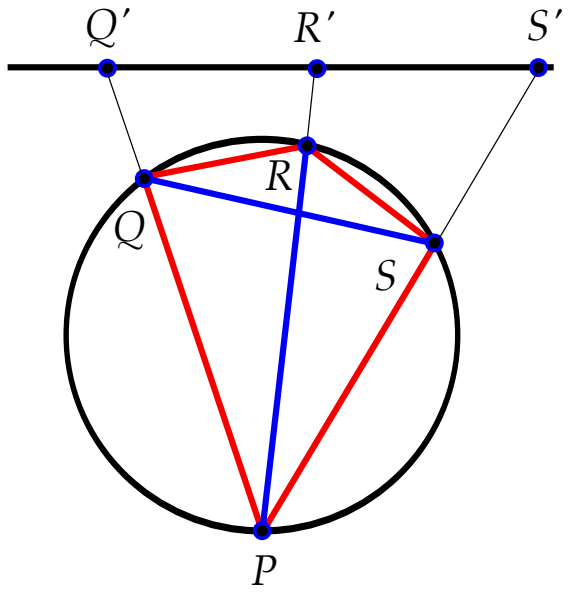

$QR.SP + PQ.RS = PR.QS.$

It turns out that the inversion is a special case of the quadratic transform. For the first vertex $A$ of our triangle, let us choose the point $[0 : 0 : 1]$ in $P^2(\mathbb{R})$, alias the origin of the plane $\mathbb{R}^2$, alias the point $0 \in \mathbb{C}$. For the second and third vertex, $B, C$, let us choose the so called *cyclic points*: those points which used to be famous among students, which are at the same time *at infinity* and *imaginary*. More precisely, they are the points $[1 : i : 0]$ and $[1 : -i : 0]$ ($i$ is $\sqrt{-1}$). They are called *cyclic* since all circles in the Euclidean plane pass through these points. For the fixed point $M$, choose for instance the point $[1 : 0 : 1]$, i.e. the point $1 \in \mathbb{C}$. I encourage my reader to show that the quadratic transform in this case is just the inversion. This can be checked by blind computation or by using classical projective geometry. Note that the image of a

I should speak of the *complexification* of circles, but this complexification was always implicit in the past. The three points $[0 : 0 : 1]$, $[1 : i : 0]$ and $[1 : -i : 0]$ are in $P^2(\mathbb{C})$ and not in $P^2(\mathbb{R})$, but the quadratic involution that they define in $P^2(\mathbb{C})$ preserves $P^2(\mathbb{R})$ and induces the inversion in the real euclidean plane, seen as the complement of the line at infinity in $P^2(\mathbb{R})$. Check!

(generic) straight line by a quadratic transformation is a conic passing through the three vertices of the triangle. Note also that any conic passing through the cyclic points is a circle. Enjoy the proof!

## *Let us work out an example*

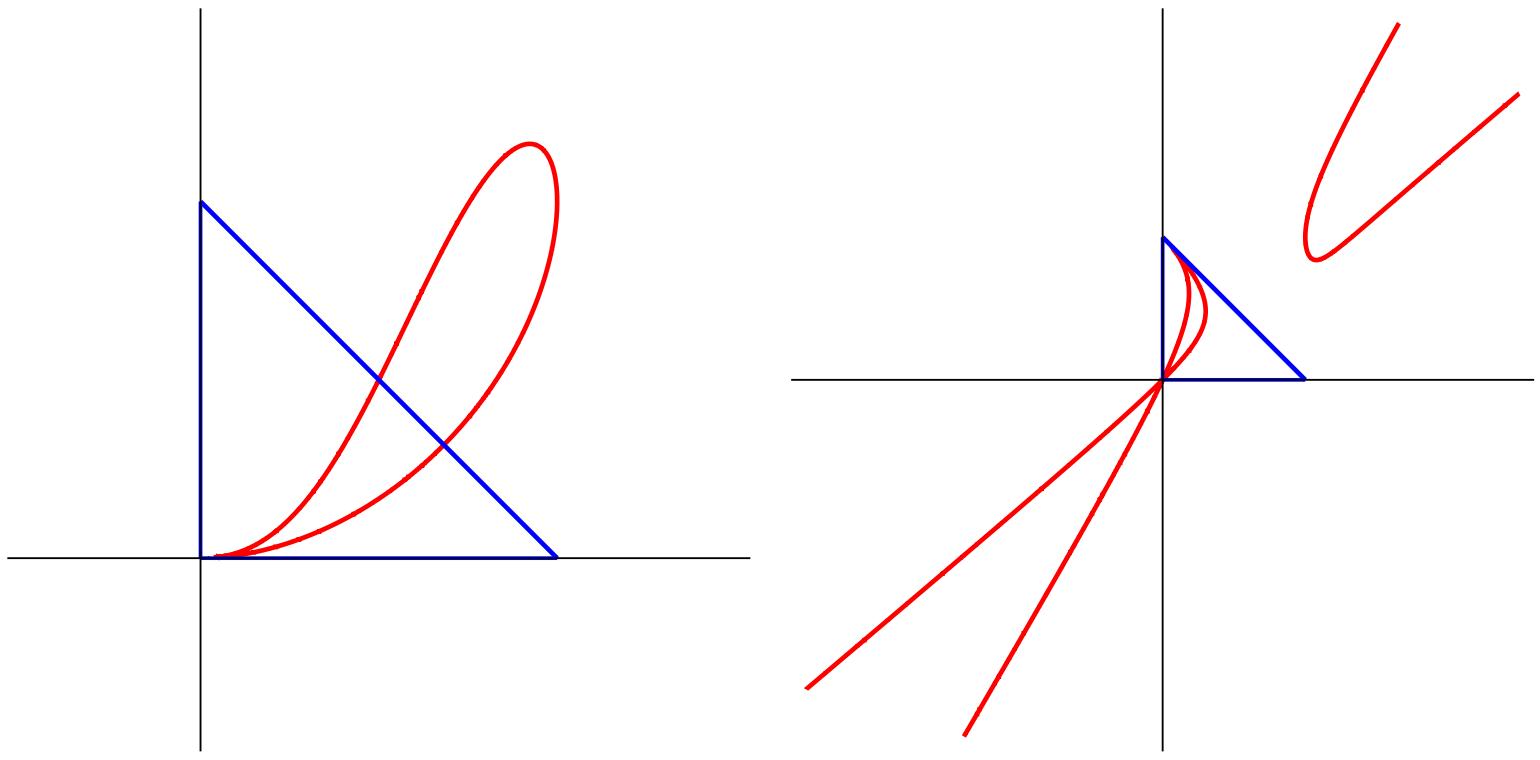

Look at the curve with equation $x^4 + x^2y^2 - 2x^2y - xy^2 + y^2 = 0$. This is an instance of Euler's *ramphoid curves*, discussed earlier, with a second order cusp.

Let us choose a (blue) triangle with one vertex at the singular point and transversal to the curve everywhere else. Let us perform a quadratic transform. The result is shown on the right picture above (in which I zoomed out). The singular point of the ramphoid curve being located at a vertex, the transformation behaves like a blowing up in the neighborhood of this point. This vertex is blown up to the opposite side of the triangle. However, the singularity is too deep to be resolved at the first step. The new curve still has a singular point at a (some other) vertex. Each edge on the triangle is collapsed to the opposite vertex and this creates a double point at the origin.

Choose some other (bigger green) triangle with a vertex on the singular point, as shown next page. Apply once more the corresponding quadratic transform. The result is shown on the right. The new curve is still singular at the lower left corner whereas the other vertices are ordinary double points.

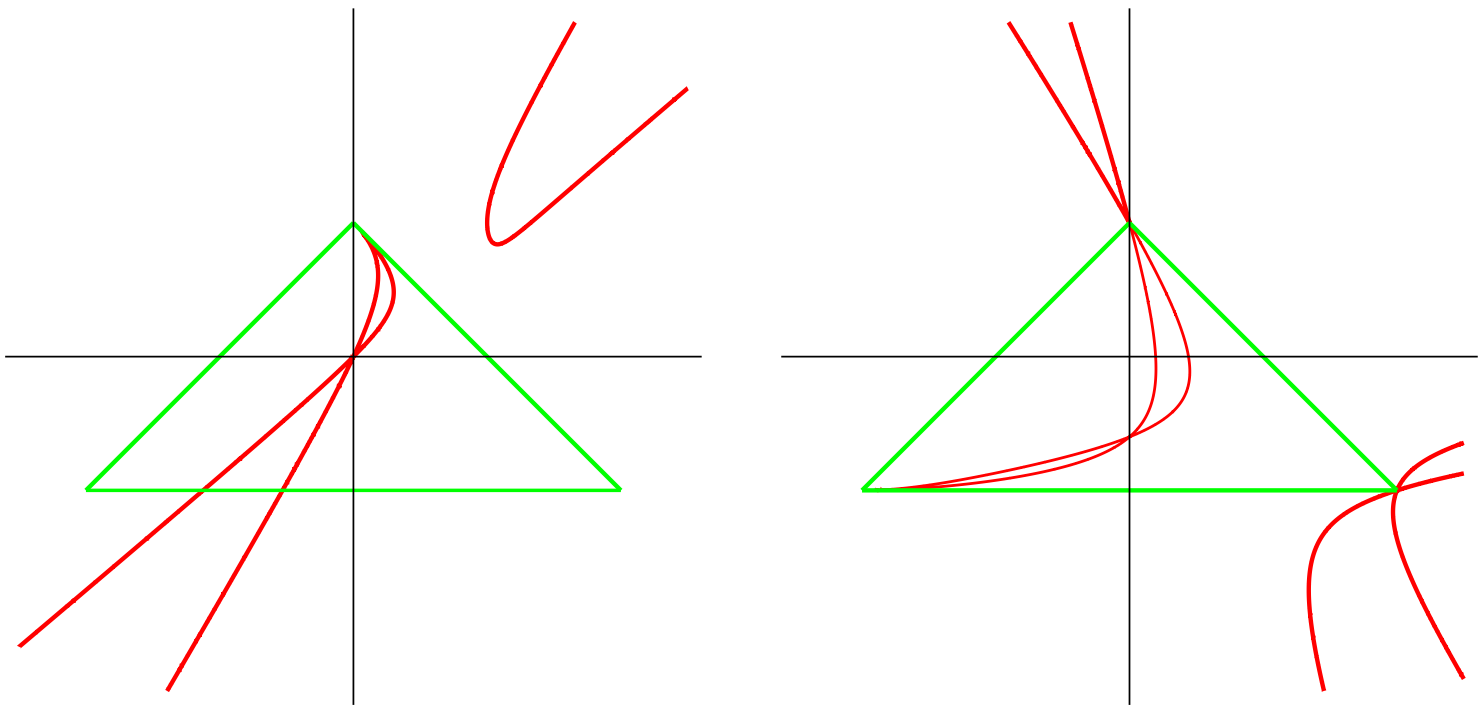

Choose some other (purple) triangle. One more quadratic transform leads finally to a curve whose only singularities are transversal intersections of smooth curves.

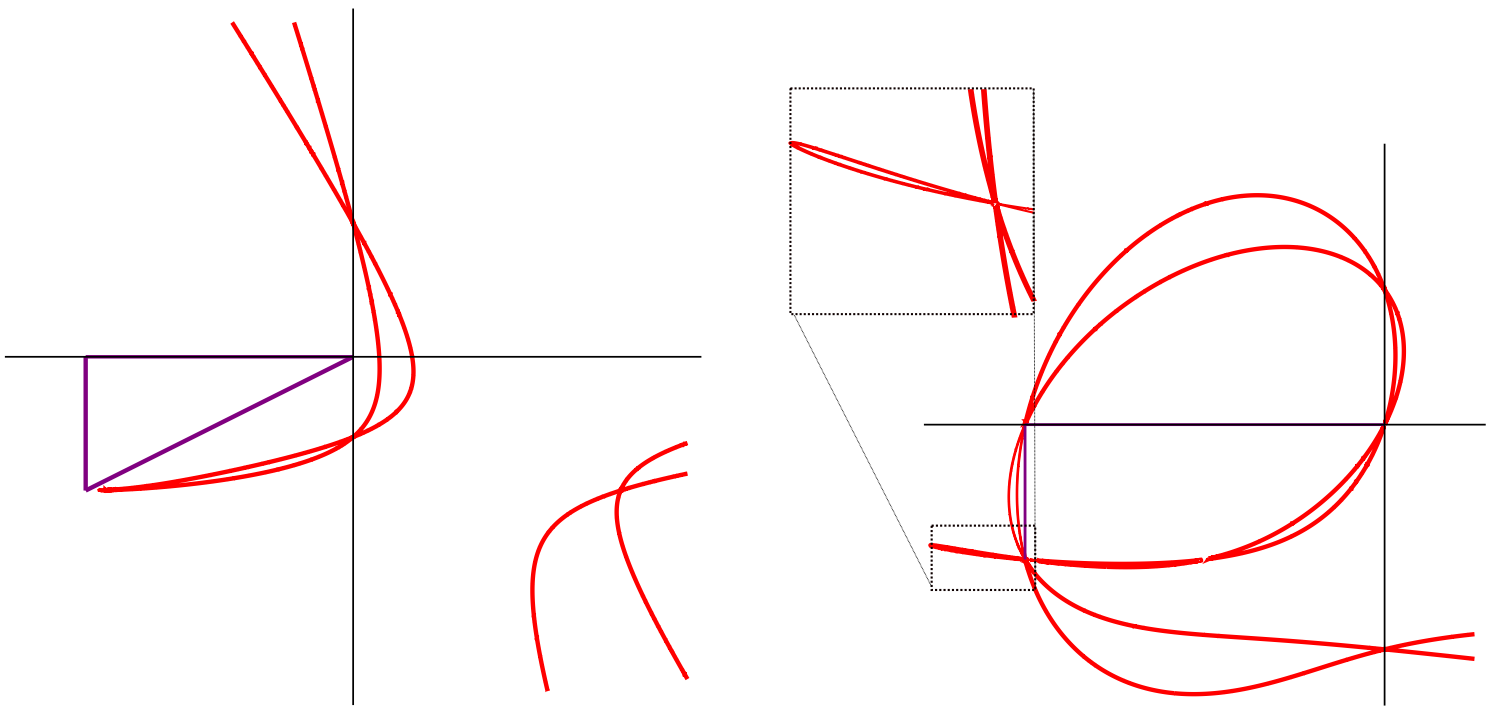

Noether's theorem is indisputably beautiful, but these *ordinary* singularities are not so simple after all. The following exercise shows that $n$ smooth curves intersecting transversally still contain too much information.

*Exercise (not so easy):* Suppose you have a finite number $n$ of *smooth* analytic curves intersecting transversally at one point. Show that for $n = 1, 2, 3, 4$ there is a local analytic diffeomorphism of the plane sending them to $n$ straight lines in the plane. Show that this is not necessarily true when $n \geq 5$. Can you describe the *moduli space* of $n$ transversal smooth curves, that is to say the quotient space under local diffeomorphisms?

There is another approach. Consider the tangent space to the projective space $P^d(K)$ of dimension $d$ over a field $K$. Projectivize this tangent space in order to produce an algebraic variety of dimension $2d-1$, which can therefore be embedded in some higher dimensional projective space $P^{2(2d-1)+1}(K)$. Given an algebraic curve $C$ in $P^d(K)$, we can look at the (Zariski) closure of the set of its tangent lines at regular points. This produces another algebraic curve $C_1$ in some other projective space of dimension $d_1$. Repeating the process, we finally get a smooth embedded curve $C_n$ in a projective space of some high dimension $d_n$. Now choose a *generic* projection to a curve in $P^2(K)$. The output is a curve $\overline{C}$ which is smooth with a finite number of *ordinary double points*.

A smooth projective algebraic variety of dimension $k$ can be, by definition, embedded in some projective space of some dimension. If we project it generically on some $2k+1$ dimensional projective subspace, this defines an embedding.

**Theorem.** *Any algebraic curve is birationally equivalent to another curve whose only singularities are ordinary double points where two smooth branches intersect transversally.*

One could be optimistic and expect that any planar algebraic curve is birationally equivalent to a *smooth* planar curve, but this is far from being true. The genus of a planar smooth curve of degree $d$ is $(d-1)(d-2)/2$ so that if an algebraic curve has a genus which is not an integer of this form, double points are compulsory.

One could be less optimistic and hope that any algebraic curve can be transformed to some curve whose singularities are ordinary double points using some Cremona transformation. Alas! This is not true either. The birational equivalence provided by the previous theorem might not be induced by some Cremona transformation (see[82] page 42).

For a modern presentation of all these concepts, I recommend Wall and Dolgachev's books[83,84] and, for a traditional version, the book by Semple and Roth[85].

[82] J. Kollár. *Lectures on resolution of singularities*, volume 166 of *Annals of Mathematics Studies*. Princeton University Press, Princeton, NJ, 2007.

[83] I. V. Dolgachev. *Classical algebraic geometry. A modern view*. Cambridge University Press, Cambridge, 2012.

[84] C. T. C. Wall. *Singular points of plane curves*, volume 63 of *London Mathematical Society Student Texts*. Cambridge University Press, Cambridge, 2004.

[85] J. G. Semple and L. Roth. *Introduction to algebraic geometry*. Oxford Science Publications. The Clarendon Press, Oxford University Press, New York, 1985. Reprint of the 1949 original.

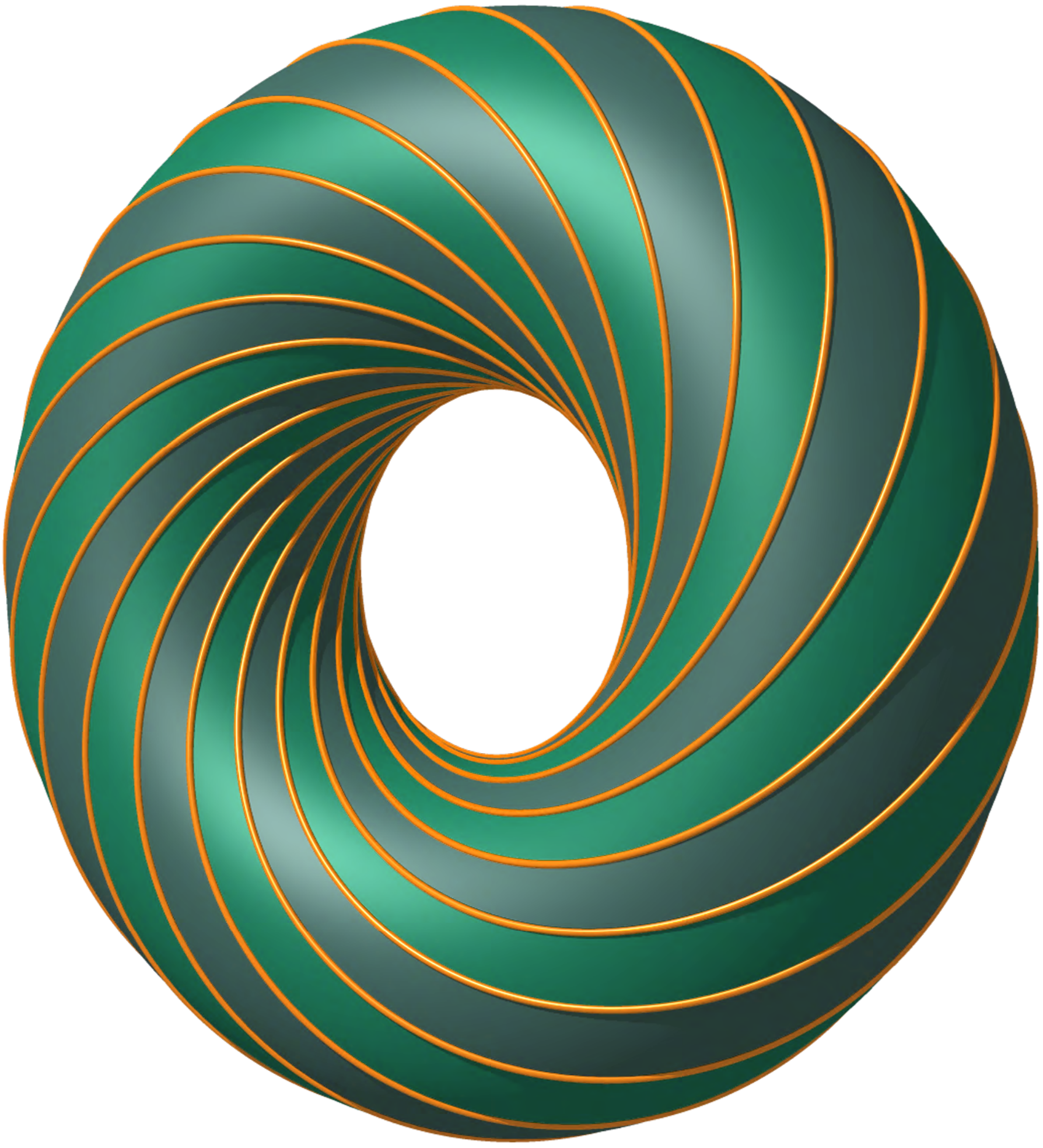

A Clifford torus filled with the so-called Villarceau circles. Each of these circles is the intersection of some line in $\mathbb{C}^2$ (passing through the origin) with the unit sphere $\mathbb{S}^3$ (and projected stereographically in 3-space).

# *The 3-sphere and the Hopf fibration*

## *A complex world?*

IT TOOK A LONG TIME BEFORE COMPLEX NUMBERS COULD BE ACCEPTED BY MATHEMATICIANS AS GENUINE NUMBERS.

> With the rise of algebra, the complex roots of real equations clamoured more and more insistently for recognition.

These are Coolidge's words in his wonderful book[86] describing the slow emergence of complex geometry in mathematics. As we have seen, Gauss was one of the most important pioneers, thinking of a complex number as a point in the plane. Visualizing $\mathbb{C}^2$ was much harder since it is 4-dimensional over the real numbers and only visionaries could imagine the fourth dimension during the nineteenth century. Many unsuccessful attempts are explained in Coolidge's book.

[86] J. L. Coolidge. *Geometry of the complex domain*. Clarendon Press, Oxford, 1924.

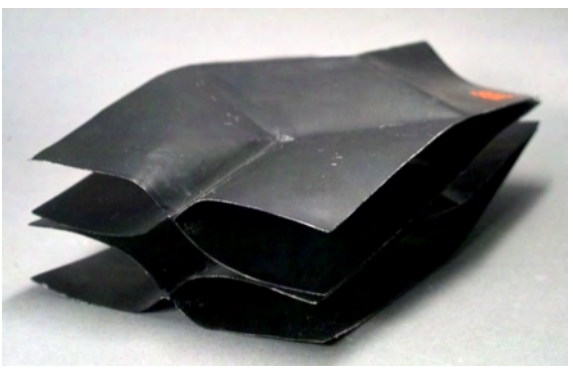

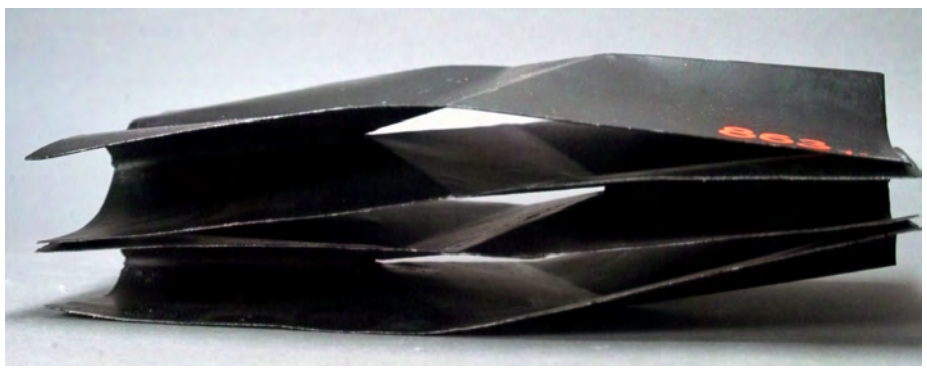

A model à la Riemann from the Göttingen Collection of Mathematical Models and Instruments. ©

Even Riemann, with his revolutionary concept now called *"Riemann surface"*, had to "see" them as some surfaces in the real 3 dimensional space, *spread over* $\mathbb{C}$ exhibiting some strange *cut lines* where the surface intersected itself, in some kind of virtual way. The least one can say is that the geometry over the complex numbers carried some air of mystery.

Nevertheless, it became progressively clear that complex geometry is not complex at all, and that it is of great help for understanding the real domain. The following quote by Paul Painlevé, in 1900[87], is a good example.

[87] P. Painlevé. *Œuvres de Paul Painlevé. Tome I.* Éditions du Centre National de la Recherche Scientifique, Paris, 1973. Analyse des travaux scientifiques, pages 72–73.

> It came to appear that, between two truths of the real domain, the easiest and shortest path quite often passes through the complex domain.

*"Il apparut que, entre deux vérités du domaine réel, le chemin le plus facile et le plus court passe bien souvent par le domaine complexe."*

Nowadays, complex geometry is better understood. Roughly speaking, there are two kinds of approaches.

The first consists in using complex numbers formally, as elements of some algebraically closed field, without any attempt to visualize them. This has been very efficient in modern algebraic geometry and indeed, the algebraic properties of $\mathbb{C}$ are amazingly powerful. The drawback is that the original questions, coming from real numbers, are usually forgotten. A famous algebraic geometer was once lecturing on (complex) Abelian varieties. At the end of his lecture, a question came about *real* Abelian varieties. The speaker was surprised and took some time before he answered, earnestly:

> Sorry, I never thought about reality!

The second approach consists in *drawing pictures*, projections, sections etc. More importantly, one tries to develop some intuition of high dimensional spaces, based on *analogy*. Modern topologists and geometers are no longer afraid by objects in $\mathbb{C}^2$ and they even consider them as very concrete. In this chapter, we try to develop some of this intuition.

According to a hoax circulating on the internet, Sophus Lie would have said:

> Life is complex because it has a real part and an imaginary part.

It is hard to believe that such a stern mathematician could have said such a thing.

## *The round* 3*-sphere*

Most of the time geometers draw a line on the blackboard when they mean $P^1(\mathbb{C})$ in $P^2(\mathbb{C})$, even if they know that $P^1(\mathbb{C})$ is a 2-dimensional (Riemann) sphere and that $P^2(\mathbb{C})$ is a non-contractible 4-dimensional manifold which does not have much in common with a blackboard.

They frequently draw a circle in the plane when they mean a 3-sphere in $\mathbb{C}^2$. They draw a real branch of a curve $P(x,y) = 0$ even though they do know that the actual topology over the complex numbers is much richer.

We use these "wrong pictures" since they are often the only possible approximation of the "reality" in the complex world.

Our goal is to give a description, as visual as possible, of the neighborhood of a point in an analytic curve $F(x,y) = 0$ in $\mathbb{C}^2$. Here, $x$ and $y$ are a complex numbers $x_1 + ix_2$ and $y_1 + iy_2$ and $x_1, x_2, y_1, y_2$ are real numbers. The curve is actually given by two equations

$$\mathfrak{R}(F(x_1 + ix_2, y_1 + iy_2)) = \mathfrak{I}(F(x_1 + ix_2, y_1 + iy_2)) = 0$$

in $\mathbb{R}^4$, so that from the point of view of real numbers, our curve is a surface. The very natural idea is to intersect our curve/surface with a small 3-dimensional sphere of radius $\epsilon$ and we hope to see something 1-dimensional (over the reals).

We therefore start with a description of the 3-sphere. The intersection of our curve with the sphere will be pictured later.

There are several ways to visualize the unit 3-sphere

$$\begin{aligned} \mathbb{S}^3 &= \{(x,y) \in \mathbb{C}^2 \mid |x|^2 + |y|^2 = 1\} \\ &= \{(x_1, x_2, y_1, y_2) \in \mathbb{R}^4 \mid x_1^2 + x_2^2 + y_1^2 + y_2^2 = 1\}. \end{aligned}$$

We could first use the *stereographic projection*.

Choose for instance the point $N = (0,0,1,0) \in \mathbb{R}^4$ as the *north pole* of $\mathbb{S}^3$ and project from $N$ to the tangent plane at the south pole $(0,0,-1,0) \in \mathbb{R}^4$. The point $(x_1, x_2, y_1, y_2) \in \mathbb{S}^3$ is mapped to $(u, v, -1, w)$ such that the points $N$, $(x_1, x_2, y_1, y_2)$ and $(u, v, -1, w)$ are aligned. In formula,

$$\Pi : (x_1, x_2, y_1, y_2) \in \mathbb{S}^3 \setminus \{N\} \mapsto \left(\frac{2x_1}{1 - y_1}, \frac{2x_2}{1 - y_1}, \frac{2y_2}{1 - y_1}\right) \in \mathbb{R}^3.$$

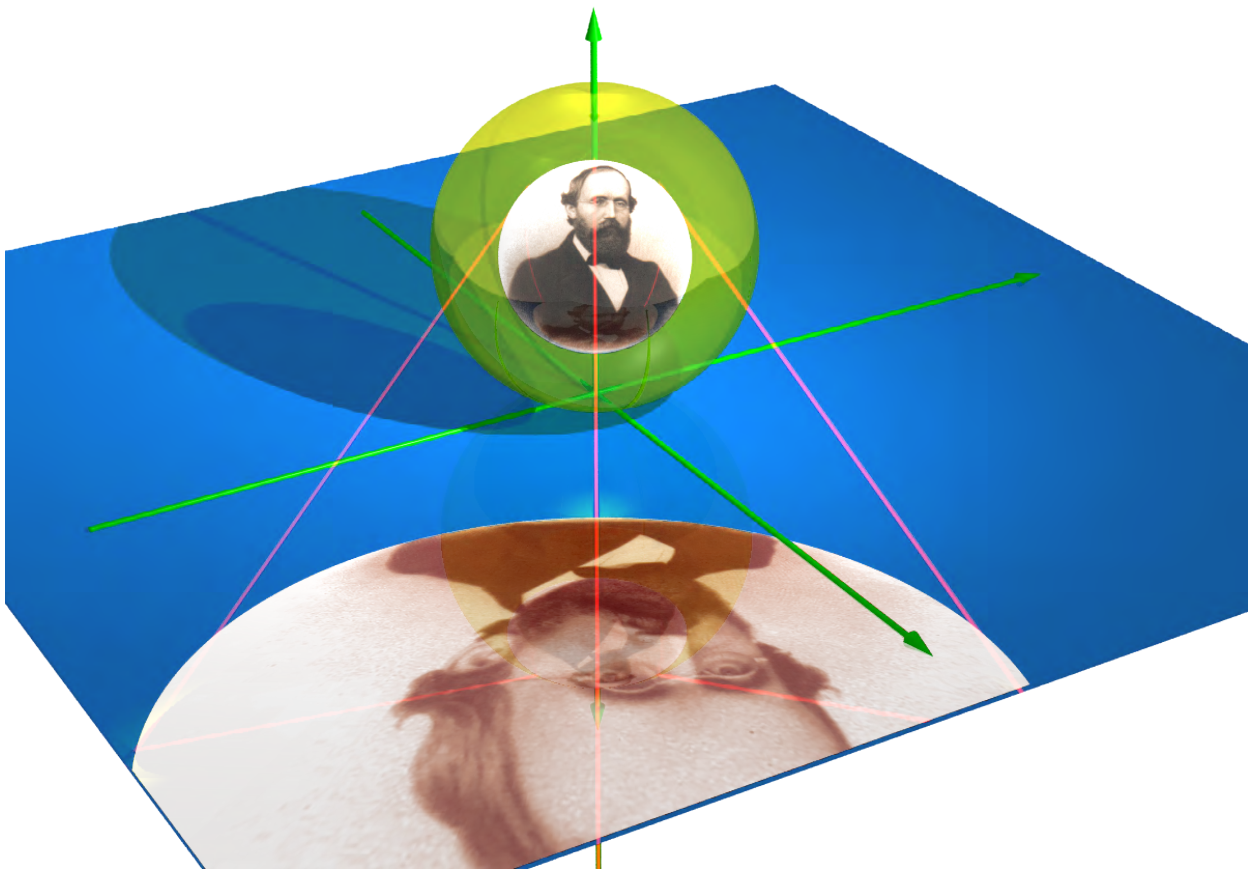

Projecting Bernhard Riemann stereographically. ©

The sphere *minus one point* can therefore be represented as the ordinary 3-space. Some symmetries are lost however since the north pole is completely arbitrary. This is a good opportunity to recommend the famous book[88] on imagination in geometry. I also recommend the movie Dimensions.

[88] D. Hilbert and S. Cohn-Vossen. *Geometry and the imagination*. Chelsea Publishing Company, New York, N. Y., 1952. Translated by P. Neményi.

The following properties of the stereographic projection are well known.

- the projection is *conformal*: its differential at any point is a similarity.
- the image of a circle on the 3-sphere is a circle in 3-space (or a straight line if the original circle passes through the north pole).

According to some historians, these properties were established (in the two dimensional case) by Hipparchus, whom we already met in this book.

The group $SO(4)$ of positive rotations of the sphere $\mathbb{S}^3$ can therefore be seen as a group of conformal diffeomorphisms of $\mathbb{R}^3 \cup \{\infty\}$.

The group of conformal diffeomorphisms of the $n$-sphere is actually much bigger than $SO(n+1)$ as it is non-compact. For instance, $\mathbb{R}^{n+2}$ can be equipped with the quadratic form of signature $(n+1,1)$ given by $q = x_1^2 + \cdots + x_{n+1}^2 - x_{n+2}^2$, so that the $n$-sphere can be interpreted as the intersection of the isotropic cone $q = 0$ with the hyperplane $x_{n+2} = 1$. Equivalently, the $n$-sphere can be thought as the space of isotropic lines. The non-compact group $SO(n+1,1)$ induces a conformal action on the $n$-sphere.

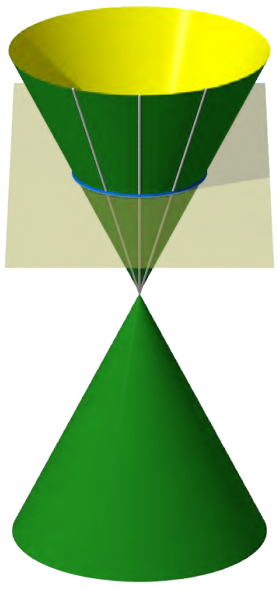

The conformal geometry of spheres is very rich.

Let me mention only two properties. Any conformal diffeomorphism between two connected open sets in a sphere of dimension at least 3 turns out to be the restriction of a global conformal diffeomorphism (Liouville's theorem). This is in strong contrast with the dimension 2 case where conformal diffeomorphisms coincide with holomorphic or anti-holomorphic diffeomorphisms, and the mathematical landscape would be much less beautiful if holomorphic maps would reduce to Moebius automorphisms $(az+b)/(cz+d)$ of the Riemann sphere.

If the conformal group of a Riemannian manifold is non-compact, this manifold is conformal to the sphere or to Euclidean space. This is the Obata and Lelong-Ferrand theorem.

I refrain from continuing in this direction since we could easily get lost and never come back from our mathematical promenade. I recommend the textbook by Berger[89] as well as his vast panorama[90]. For the reader interested in old fashioned presentations, the book by Coolidge[91] is beautiful.

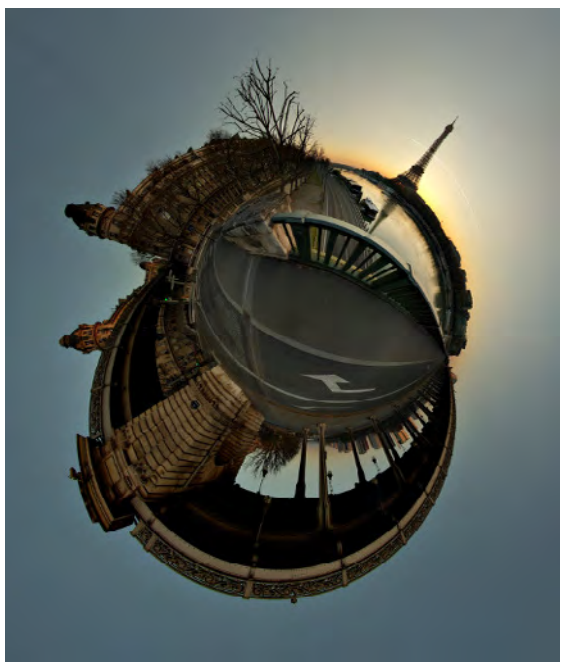

A conformal view of Paris. ©

[89] M. Berger. *Geometry I*. Universitext. Springer-Verlag, Berlin, 2009.

[90] M. Berger. *Geometry revealed. A Jacob's ladder to modern higher geometry*. Springer, Heidelberg, 2010.

[91] J. L. Coolidge. *A treatise on the circle and the sphere*. Clarendon Press, Oxford, 1916.

## *The "square" 3-sphere*

Squaring the circle?

Since $|x|^2+|y|^2=1$ on the 3-sphere, we can split it into two parts $T_1, T_2$ defined by

$$T_1=\{(x,y)\in\mathbb{S}^3 |\ |x|^2\le 1/2\}\quad,\quad T_2=\{(x,y)\in\mathbb{S}^3 |\ |y|^2\le 1/2\}.$$

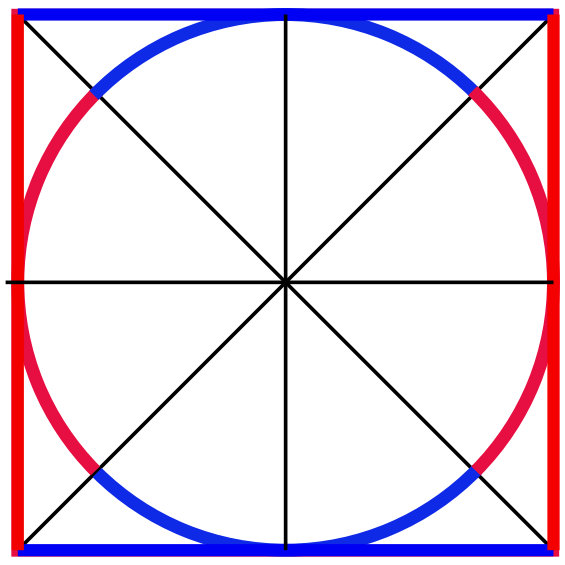

The intersection of $T_1$ and $T_2$ is a *Clifford torus* parameterized by

$$(\theta,\phi)\in(\mathbb{R}/2\pi\mathbb{Z})^2\mapsto\left(\frac{\sqrt{2}}{2}\exp(i\theta),\frac{\sqrt{2}}{2}\exp(i\phi)\right)\in\mathbb{S}^3.$$

As for $T_1$ and $T_2$, they are *solid tori*, parameterized by a product of a unit disc $D^2$ in $\mathbb{C}$ and a circle.

$$(z,\phi)\in D^2\times(\mathbb{R}/2\pi\mathbb{Z})\mapsto\left(\frac{\sqrt{2}}{2}z,\sqrt{1-\frac{|z|^2}{2}}\exp(i\phi)\right)\in T_1\subset\mathbb{S}^3$$

$$(\theta,z)\in(\mathbb{R}/2\pi\mathbb{Z})\times D^2\mapsto\left(\sqrt{1-\frac{|z|^2}{2}}\exp(i\theta),\frac{\sqrt{2}}{2}z\right)\in T_2\subset\mathbb{S}^3.$$

The 3-sphere is therefore the union of two solid tori, glued along their boundaries. The meridians of $\partial T_1$, that is to say the circles which bound a disc in $T_1$, are glued to parallels of $\partial T_2$, which do not bound a disc in $T_2$, and conversely.

We could also use the "square sphere"

$$D^2 \times D^2 = \{(x,y) \in \mathbb{C}^2 \,|\, |x| \leq 1 \text{ and } |y| \leq 1\}.$$

Its boundary consists of two solid tori $T_1' = \{|x| \leq 1 \text{ and } |y| = 1\}$ and $T_2' = \{|x| = 1 \text{ and } |y| \leq 1\}$. Using radial projection, the two solid tori $T_1, T_2$ are identified with $T_1', T_2'$. It is frequently more convenient to use the square sphere, since we can draw pictures in the solid torus without having to use stereographic projection. This simple but very useful idea is due to Kähler[92].

[92] E. Kähler. über die verzweigung einer algebraischen funktion zweier veränderlichen in der umgebung einer singulären stelle. *Math. Z.*, 30(1):188–204, 1929.

## *The* 3*-sphere is* very *round*

William Thurston[93], one of the masters of the visual aspect of mathematics, used to say that the 3-sphere is "rounder" than the other spheres. He had in mind the important fact that the group $SO(n+1)$ is *not* a simple group if and only if $n = 3$ (and of course $n = 0, 1$). This is related to what was called in the old literature *Clifford's parallelism*.

[93] W. P. Thurston. How to see 3-manifolds. *Classical Quantum Gravity*, 15(9):2545–2571, 1998. Topology of the Universe Conference (Cleveland, OH, 1997).

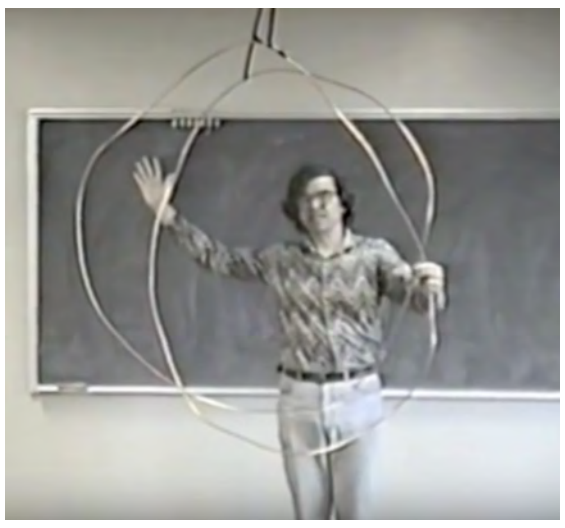

A screenshot from *Knots to Narnia*. William Thurston (1946—2012) shows that when he goes around a knot, he arrives "somewhere else". This vintage YouTube video is a must. ©

Recall that quaternions are formal expressions of the form $q = x_1 + ix_2 + jy_1 + ky_2$ where $x_1, x_2, y_1, y_2$ are real numbers and the formal symbols $i, j, k$ satisfy $ij = -ji = k; jk = -kj = i; ki = -ik = j$ and $i^2 = j^2 = k^2 = -1$. This defines a non-commutative division algebra $\mathbb{H}$. The conjugate $\bar{q}$ of $q$ is defined to be $x_1 - ix_2 - jy_1 - ky_2$ and the norm $N(q)$ is the product $q\bar{q} = x_1^2 + x_2^2 + y_1^2 + y_2^2$. This norm is multiplicative, i.e. $N(q_1q_2) = N(q_1)N(q_2)$ and the inverse of a nonzero quaternion is $q^{-1} = \bar{q}/N(q)$.

It follows that *the* 3*-sphere is identified with the group of unit quaternions* $\{q \in \mathbb{H} \,|\, N(q) = 1\}$. It is one of the great successes of the mathematical twentieth century to prove that the only spheres which can be equipped with the structure of a topological group are $\mathbb{S}^0 \simeq \mathbb{Z}/2\mathbb{Z}$, $\mathbb{S}^1 \simeq SO(2)$ and $\mathbb{S}^3$. A good starting point for this topic is *Numbers*[94].

[94] H.-D. Ebbinghaus, H. Hermes, F. Hirzebruch, M. Koecher, K. Mainzer, J. Neukirch, A. Prestel, and R. Remmert. *Numbers*, volume 123 of *Graduate Texts in Mathematics*. Springer-Verlag, New York, 1990.

But this is not the only reason why the 3-sphere is rounder.

Just as any group, it can be seen as homogeneous in two commuting ways, using the right and left actions. Given two unit quaternions $q_1, q_2$, the map $q \in \mathbb{H} \mapsto q_1 q q_2^{-1} \in \mathbb{H}$ is an isometry, and defines an element of $SO(4)$. It turns out that this homomorphism from $\mathbb{S}^3 \times \mathbb{S}^3$ to $SO(4)$ is onto and its kernel only contains $\pm(1,1)$. In other words, every rotation of the 3-sphere is the composition of a *left rotation* and a *right rotation* which commute. This situation is unique to dimension 3 as all other rotation groups are simple (with the obvious exception of $SO(2)$).

A (non-trivial) group is *simple* if it does not contain a proper normal subgroup. It is usual to say that a Lie group is simple if every normal Lie subgroup is either discrete or open. This is equivalent to say that its Lie algebra is simple, i.e. does not contain a proper ideal. The only non-proper normal subgroup of $SO(n)$ (for $n \neq 1$) is $\{\pm Id\}$ for $n$ even, so that $SO(n)$ is not simple as a group for $n$ even, but is simple as a Lie group.

## *The Hopf fibration*

We can now begin our description of the topology of algebraic curves in $\mathbb{C}^2$ and start with the simplest possible curve: *a straight line*.

Let us look at the intersection of the lines $x = 0$ and $y = 0$ with the unit sphere.

Under the stereographic projection, since the line $x = 0$ passes through the north pole, its image is simply a vertical straight line. The other line $y = 0$ is projected onto a circle which "goes around the vertical line $x = 0$".

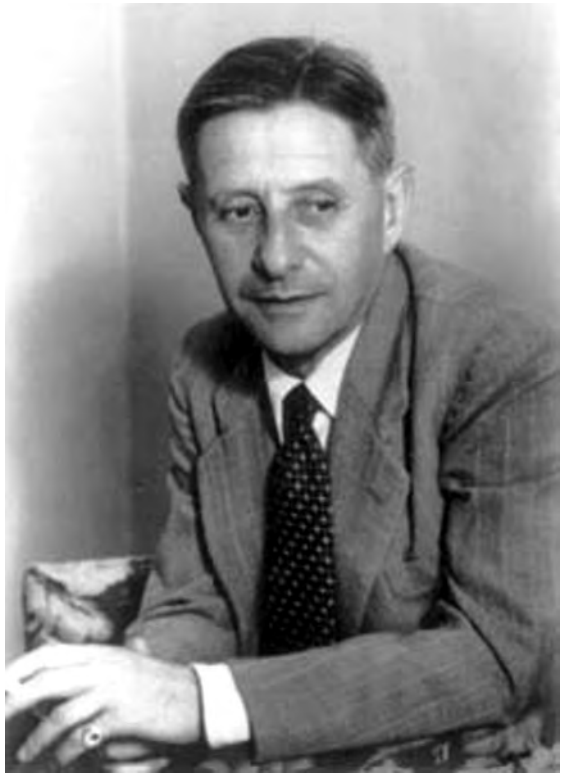

Heinz Hopf (1894–1971) should not be confused with Eberhard Hopf (1902–1983) (the only mathematician who moved from the US to Germany in 1936?). ©

In our decomposition in two solid tori, $x = 0$ becomes the circle $\{0\} \times (\mathbb{R}/2\pi\mathbb{Z})$ which is the core of $T_1$. Conversely, $y = 0$ becomes the circle $(\mathbb{R}/2\pi\mathbb{Z}) \times \{0\}$ which is the core of $T_2$. Note that these two circles are linked.

The line $y = x$ intersects the sphere on a circle which is in $T_1$ and $T_2$: it is neither a meridian nor a parallel but its homotopy class is $(1,1)$ in both $T_1$ and $T_2$.

All this structure is globally described by the so-called *Hopf fibration*. Every point $(x, y)$ of the punctured plane $\mathbb{C}^2 \setminus \{(0,0)\}$ belongs to a unique complex line passing through the origin, that is to say defines an element of $P^1(\mathbb{C})$. In other words, a line through the origin has an equation $y = \lambda x$ where $\lambda$ belongs to $\mathbb{C} \cup \{\infty\}$, identified with the Riemann sphere, or with a 2-dimensional sphere $\mathbb{S}^2$. This defines a map

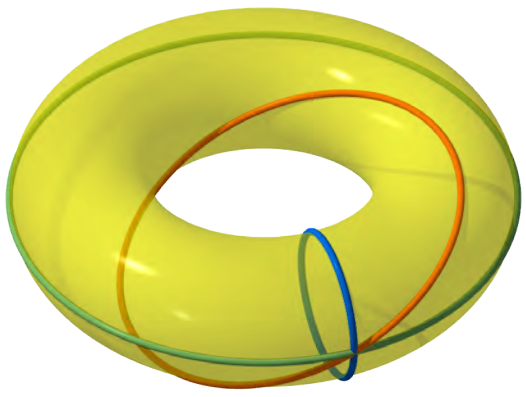

A meridian, a parallel, and a (1,1)-circle on a torus. ©

$$\pi : \mathbb{S}^3 \to \mathbb{S}^2$$

whose fibers are circles, intersections of complex lines with the sphere. Any two fibers are linked.

Of course, this map was not invented by Hopf! His contribution was to show that it is not homotopic to a constant map, but that is another story.

Here are some pictures of the Hopf fibration, under the stereographic projection, extracted from Dimensions.

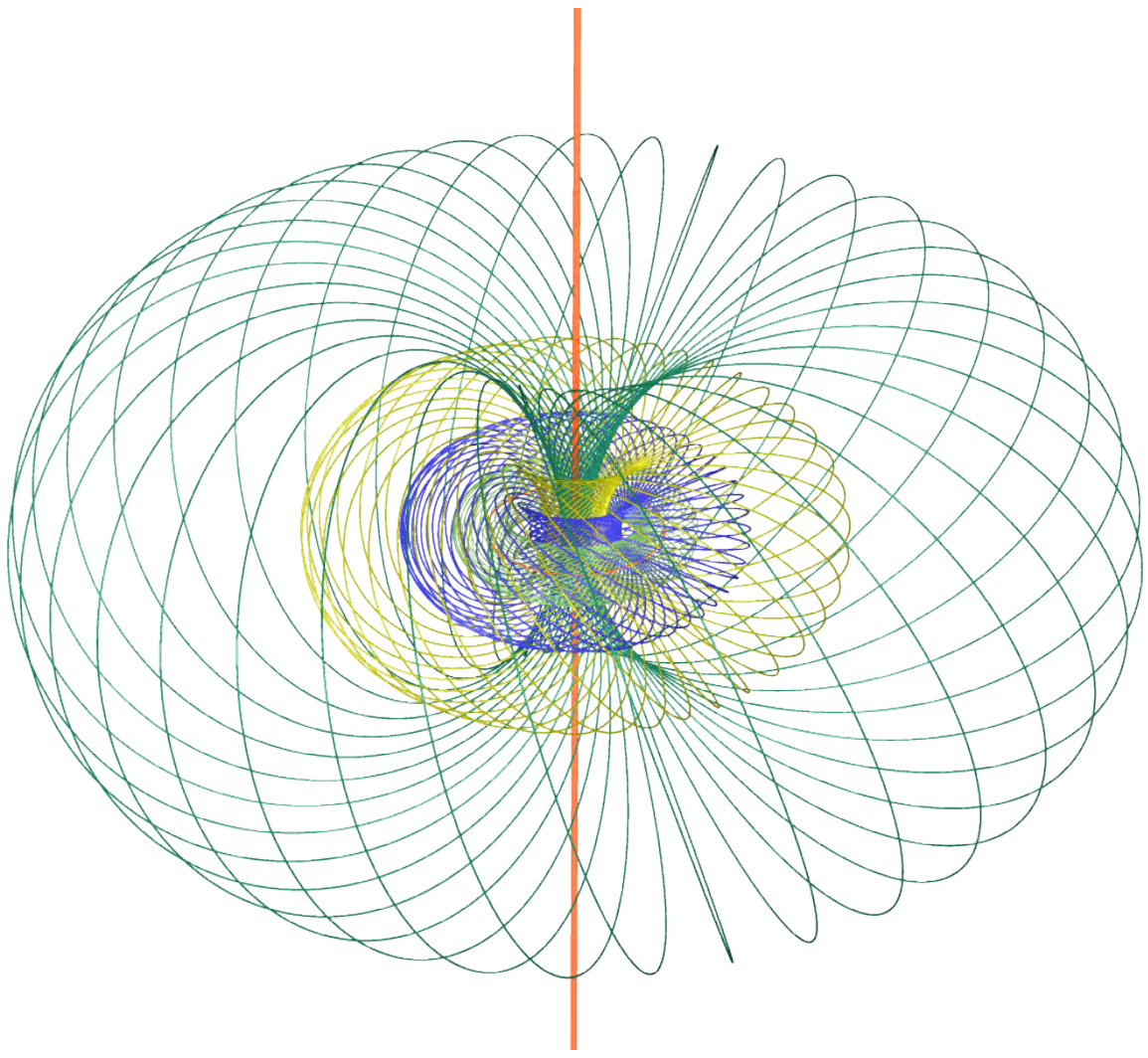

The Hopf fibration: each circle is a fiber of $\pi$. The inverse image of a circle by $\pi$ is a Clifford torus, which is a union of fibers (of the same color on the picture). ©

What is the real version of the Hopf fibration? It does exist but it is a little bit disappointing. Every point $(x, y)$ of $\mathbb{R}^2 \setminus \{(0,0)\}$ belongs to a unique line passing through the origin, and defines an element of $P^1(\mathbb{R})$. Such a line has the form $y = \lambda x$ where $\lambda$ belongs to $\mathbb{R} \cup \{\infty\}$ which is a 1-dimensional sphere, i.e. a circle $\mathbb{S}^1$. This defines a map

$$\pi : \mathbb{S}^1 \to \mathbb{S}^1$$

whose fibers are $\mathbb{S}^0$, intersections of real lines with the unit circle. This is just the multiplication by 2 in $\mathbb{R}/\mathbb{Z}$. Do not forget that a zero dimensional sphere is a pair of points.

## *Hopf links*

A Hopf fiber is just a round circle in the sphere, so there is not much to say about it. (This is not quite true: the geometry of the space of circles in 3-space is wonderful. Look at the modern book by Cecil[95] or Blaschke's classical Vorlesungen[96]).

[95] T. E. Cecil. *Lie sphere geometry. With applications to submanifolds*. Universitext. Springer-Verlag, New York, 1992.

[96] W. Blaschke. *Vorlesungen über Differentialgeometrie und geometrische Grundlagen von Einsteins Relativitätstheorie. Band I. Elementare Differentialgeometrie*. Dover Publications, New York, N. Y., 1945. 3d ed.

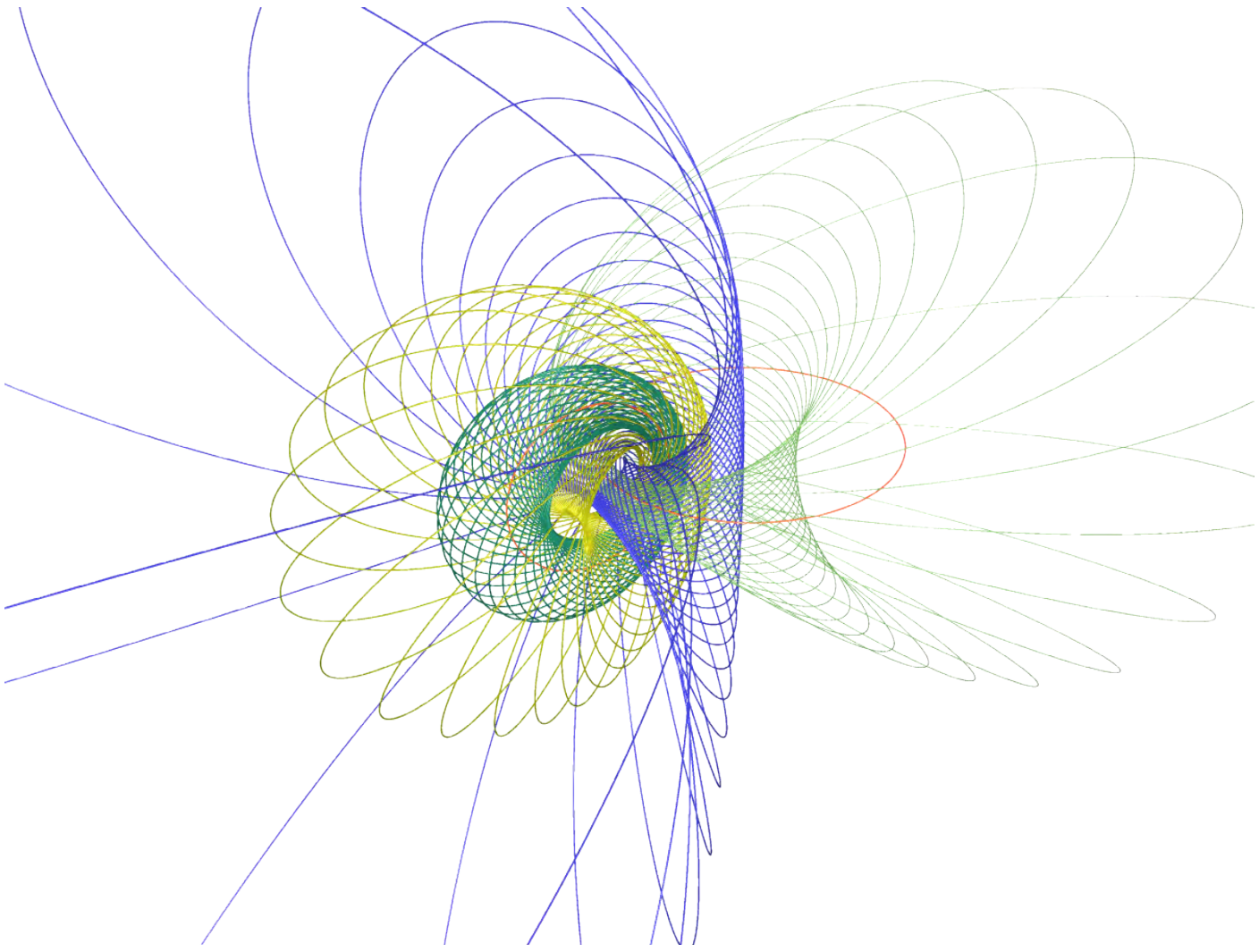

The same picture as above, after a rotation of the 3-sphere, which corresponds to a conformal map on $\mathbb{R}^3 \cup \{\infty\}$. ©

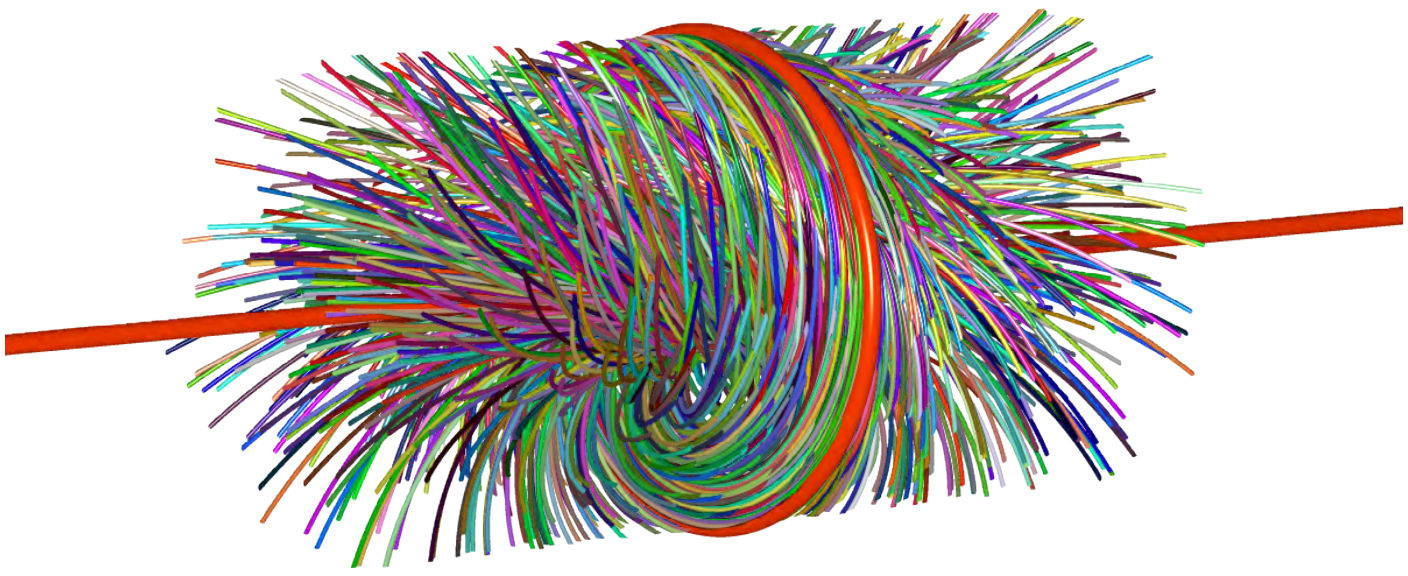

Hopf circles in the neighborhood of one of them, projected as a line in space (in red) . One of my readers told me that he sees this figure as "frightening". Do you agree with him ? ©

*Two Hopf circles* are more interesting since they define the simplest non-trivial link. Notice that even though they are linked, the two circles bound an annulus. Indeed, look at the pre-image by the Hopf fibration of some arc connecting two points: it is an annulus.

Two Hopf circles, bounding an annulus. ©

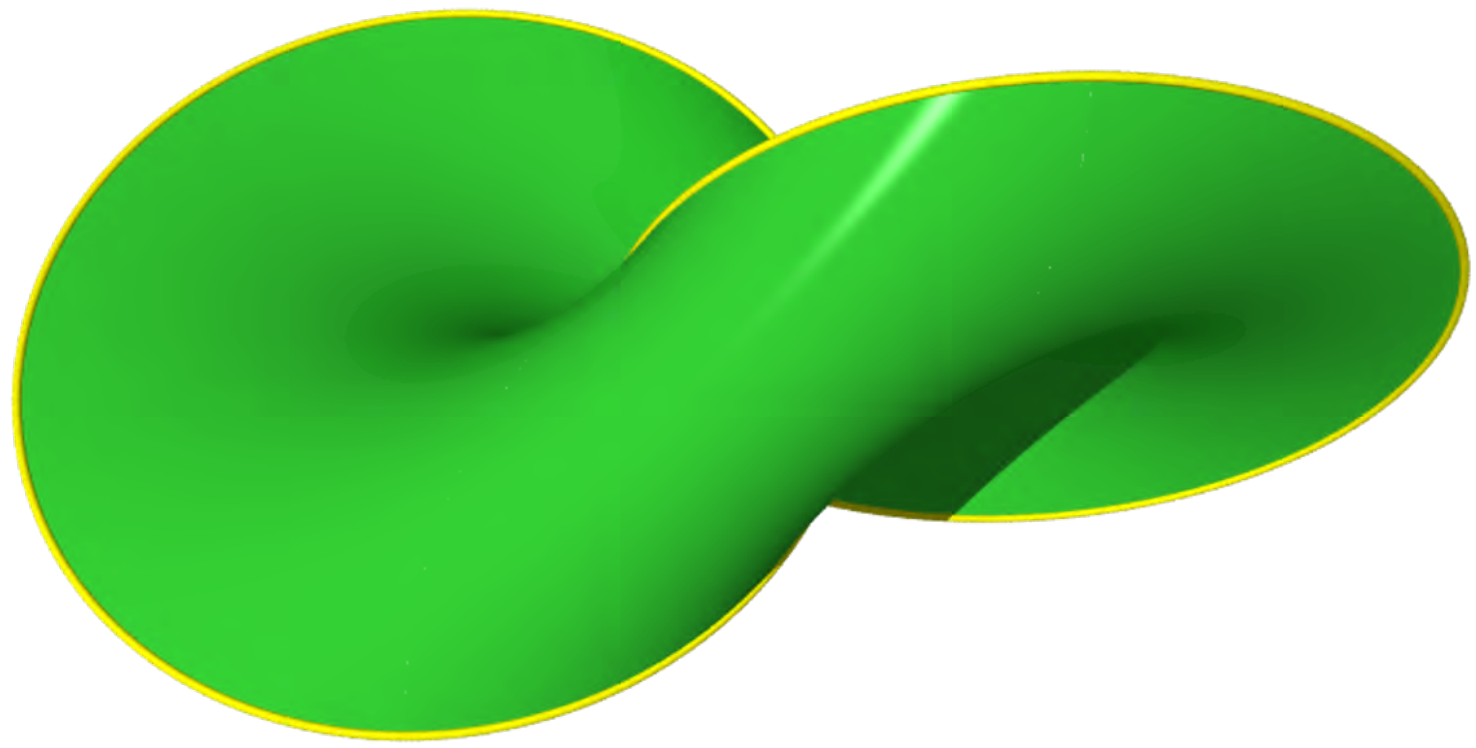

*Three Hopf circles* (or more) give a *Hopf link*. Each component is a circle and any two components are linked once. It is easy to find an orientable surface having such a link as boundary. Indeed, let us consider $n$ complex numbers $\lambda_1, \ldots, \lambda_n$ and the polynomial

$$F(x,y) = (y - \lambda_1 x)(y - \lambda_2 x) \cdots (y - \lambda_n x).$$

The intersection of the 3-sphere with the set of $(x,y)$ such that $F(x,y)$ is a positive real number is a surface whose boundary consists of $n$ Hopf circles. All this will be greatly generalized in the following chapters.

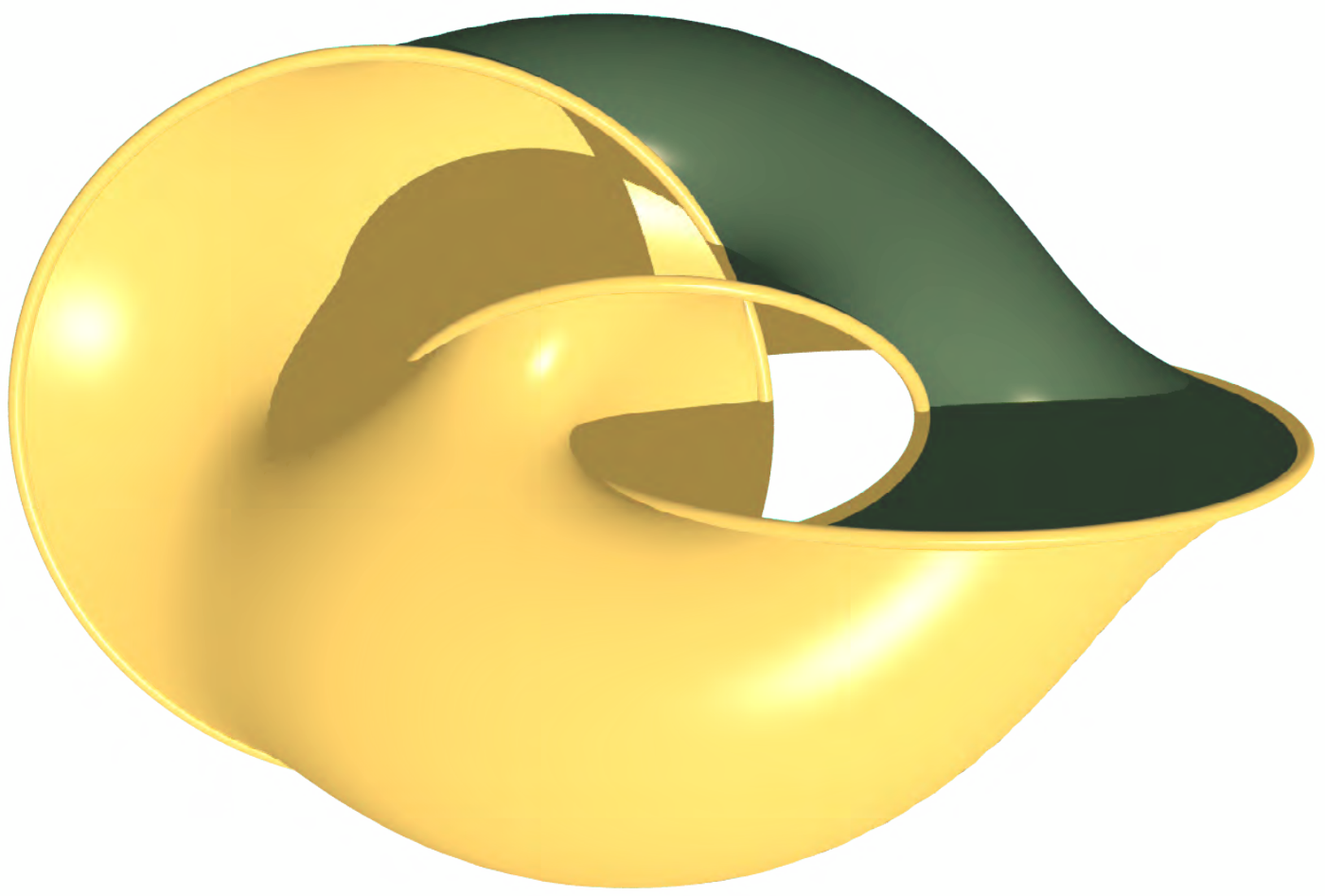

Three Hopf circles, bounding a surface. ©

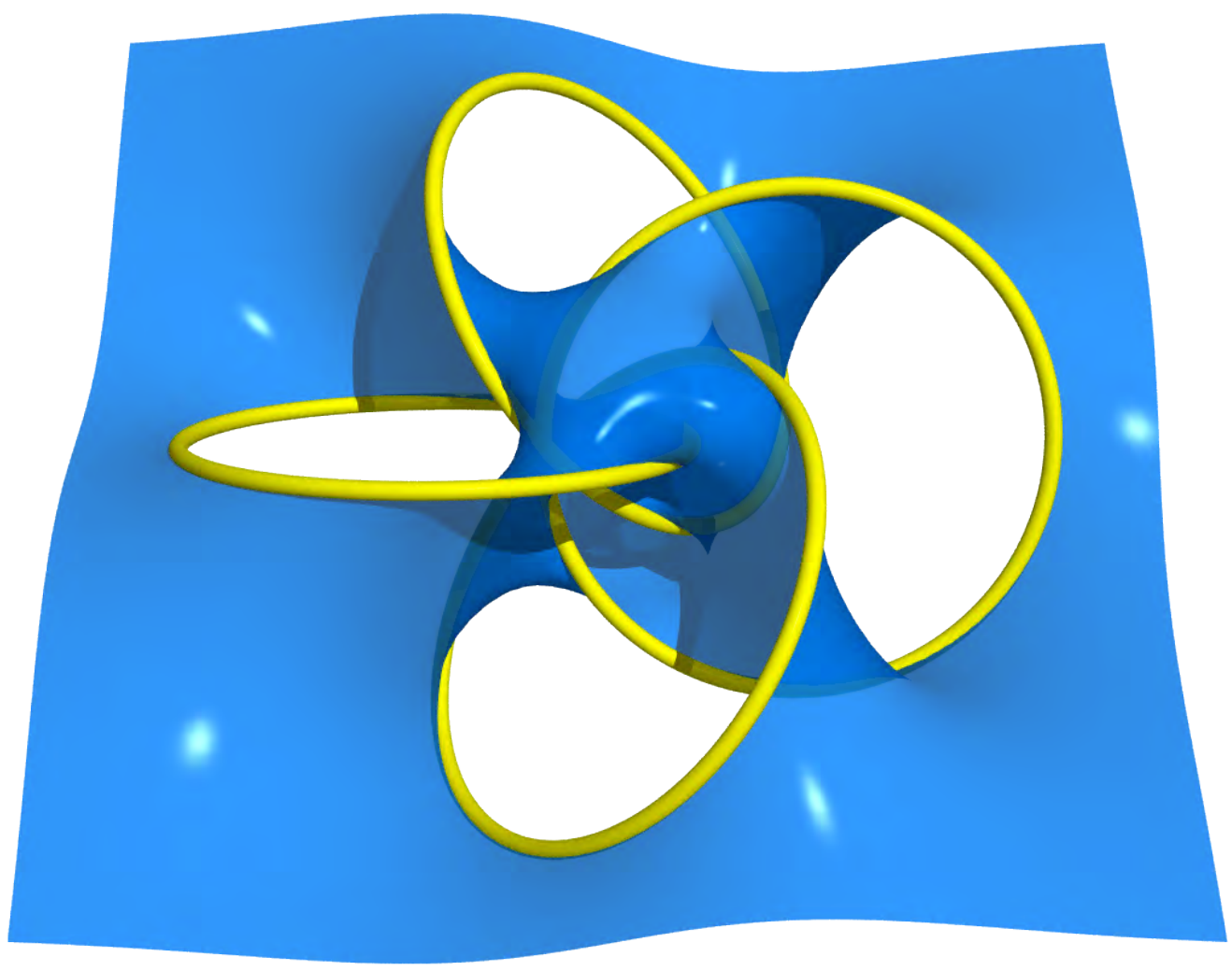

Four Hopf circles, bounding a surface. ©

## *Dante,* La Divina Commedia *and the* 3*-sphere*

This section is, of course, not necessary for the rest of the book.

It has been argued by Mark Peterson[97] that Dante's universe, as described in the *Divine Comedy*, is homeomorphic to a 3-sphere. Even though I am not fully convinced that "it is clear that Dante invents the notion of manifold[98]", I like this cosmological vision. At least it answers a question that all of us asked our parents when we were children: what happens when we reach the boundary of the universe ? Well, Dante's universe is a compact 3-manifold without boundary!

[97] M. A. Peterson. Dante and the 3-sphere. *American Journal of Physics*, 47:1031–1035, 1979.

[98] M. A. Peterson. The geometry of *paradise*. *Math. Intelligencer*, 30(4):14–19, 2008.

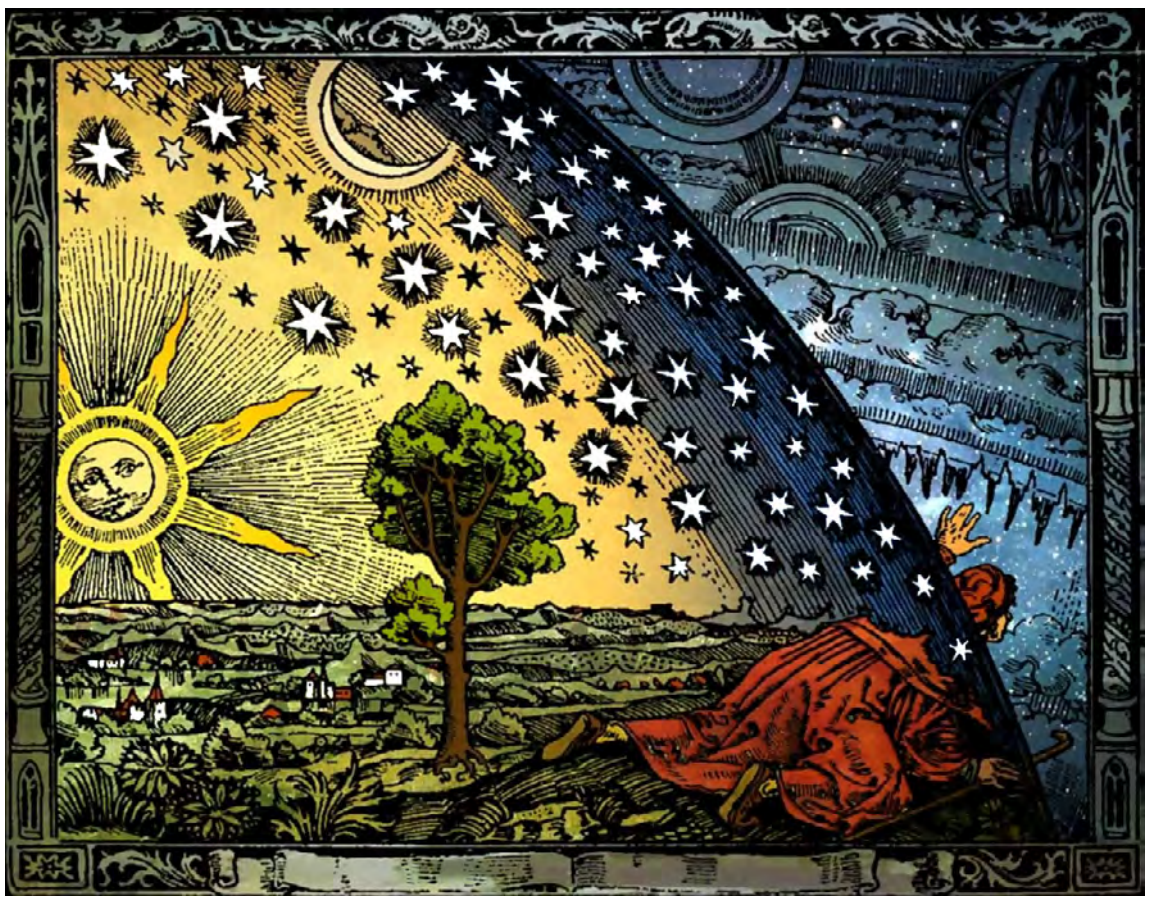

Flammarion engraving (1888). ©

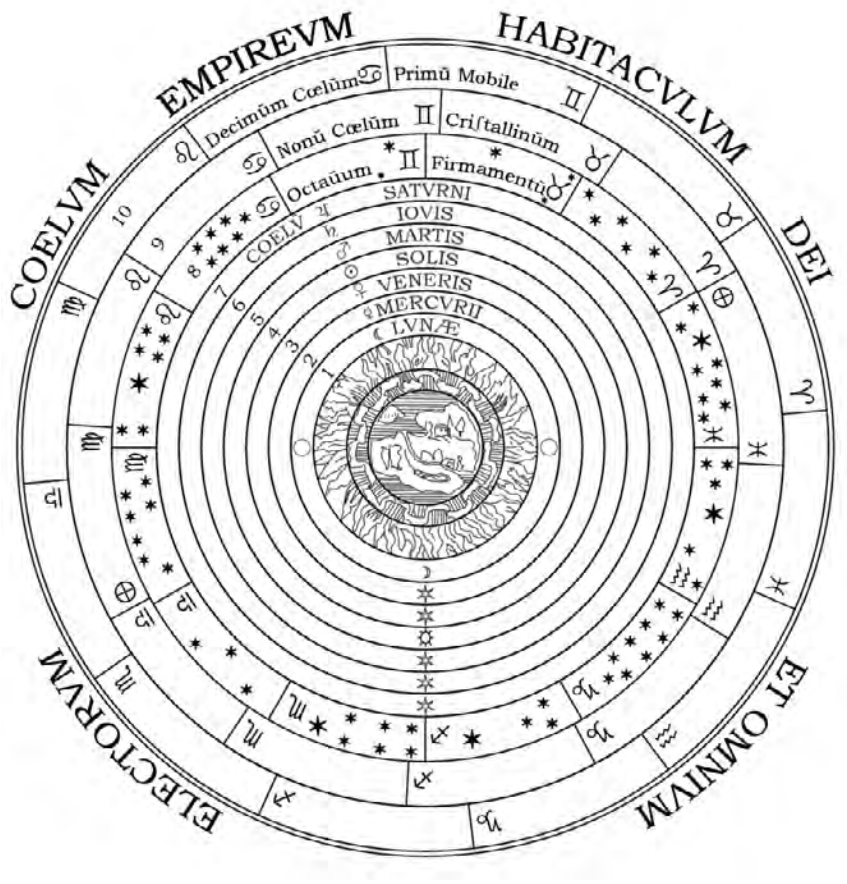

Aristotle Ptolemy geocentric system. ©

Let me recall that the world inherited from the ancient Greeks is finite[99]. Basically, in the Aristotle-Ptolemy system, the Earth is fixed at the center of the universe and is surrounded by seven celestial spheres, each one carrying a "planet": the Moon, Mercury, Venus, Sun, Mars, Jupiter and Saturn. An eighth sphere carries the *fixed stars*. Finally, a ninth sphere, called *Primum Mobile*, serves as a container for the full system and generates the motion of the other spheres. The sensible world is therefore a 3-dimensional ball, whose boundary is the Primum Mobile. Beyond this boundary begins the realm of the so-called *Empyrean* whose nature is not very clear. According to Aristotle[100] "it is clear that [it contains] neither space, nor void, nor time".

[99] A. Koyré. *From the Closed World to the Infinite Universe*. Johns Hopkins University Press, 1957.

Of course the Moon and the Sun are not planets but they turn around the Earth... in the geocentric system.

[100] Aristotle. *On the Heavens*, volume I, 9, 278b–279a. transl. by J. L. Stocks.

The *Divine Comedy* is a long poem written in 1320 which offers a fascinating description of the christian world in the 14th century[101]. Dante tells us about his journey through Hell, Purgatory and Paradise, as an allegory to the salvation of souls. In the final part, Paradiso, his muse Beatrice helps him to visit successively the nine celestial spheres. When he arrives on the Primum Mobile, he can contemplate the world from its boundary, with the tiny Earth at the center. Suddenly, he turns back and discovers that the Empyrean has exactly the same structure as the sensible world. It consists of the same number of spheres which are now centered at God. These angelic spheres, symmetric to the celestial spheres, have the following names (going from God to the Primum Mobile): Seraphim, Cherubim, Thrones, Dominations, Virtues, Powers, Principalities, Archangels, Angels.

*Dante's universe is therefore the union of two 3-balls, glued along the Primum Mobile. Therefore, it is a 3-sphere!* ⊡

See[102] for a much deeper discussion on the medieval vision of the universe and [103] for more about poetry and mathematics.

[101] Dante Alighieri. *Divina Commedia*. Digital edition, Columbia University. There are many editions and translations. I recommend this beautiful 2017 edition, with English translation and commentaries.

[102] J. Grzybowski. *Cosmological and Philosophical World of Dante Alighieri: "The Divine Comedy"*. Peter Lang GmbH, 1st new edition, 2015.

[103] R. Osserman. *Poetry of the universe, From the Divine Comedy to Riemann and Einstein*. Anchor, 1995.

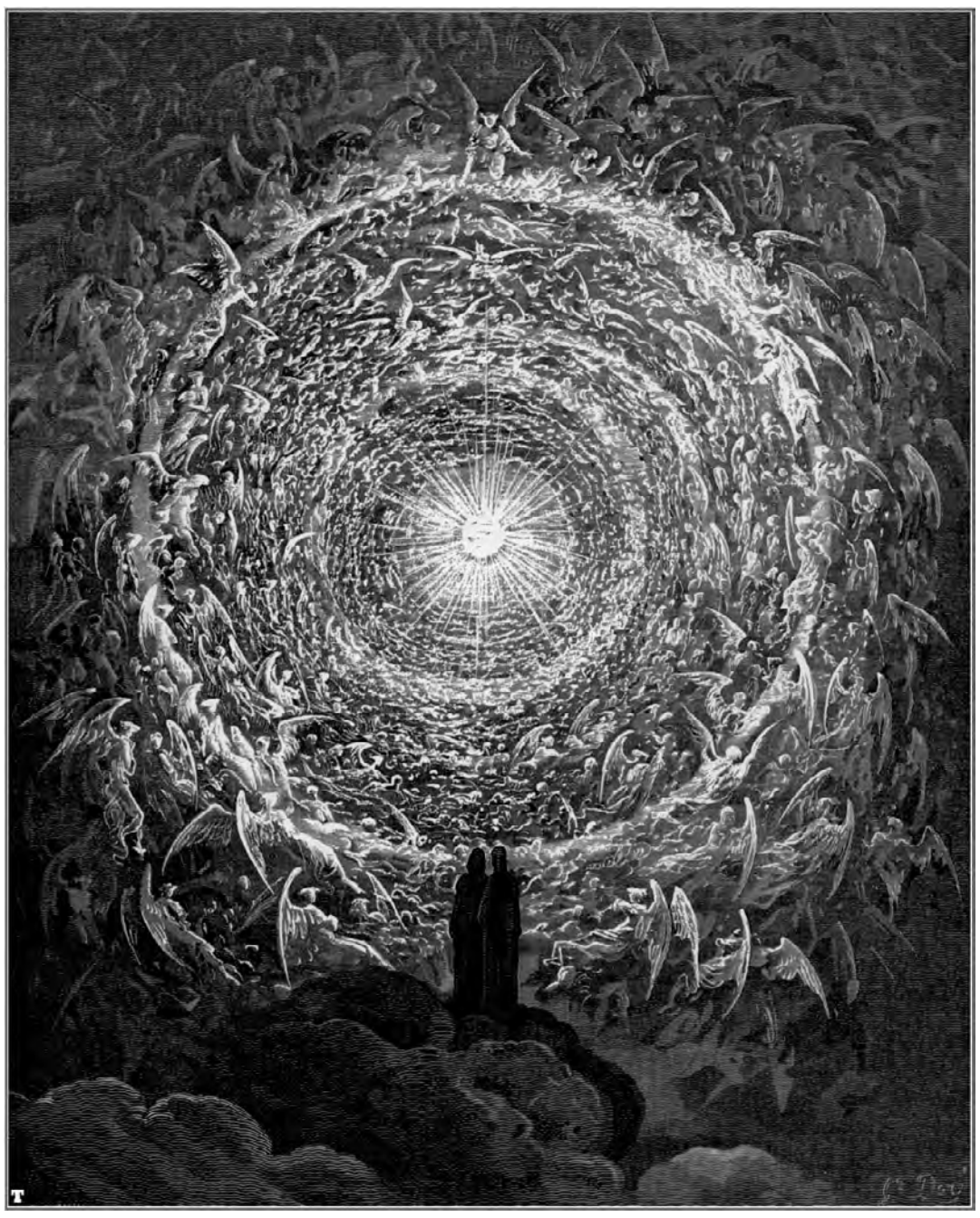

Rosa Celeste: Dante and Beatrice gaze upon the highest Heaven, The Empyrean (engraving by Gustave Doré). ©

The last verses of the Comedy:

*"All'alta fantasia qui mancò possa ; ma già volgeva il mio disiro e il velle, l'amor che move il sole e l'altre stelle."*

*"To the high fantasy here power failed; but already my desire and will were rolled — even as a wheel that moveth equally — by the Love that moves the sun and the other stars."*

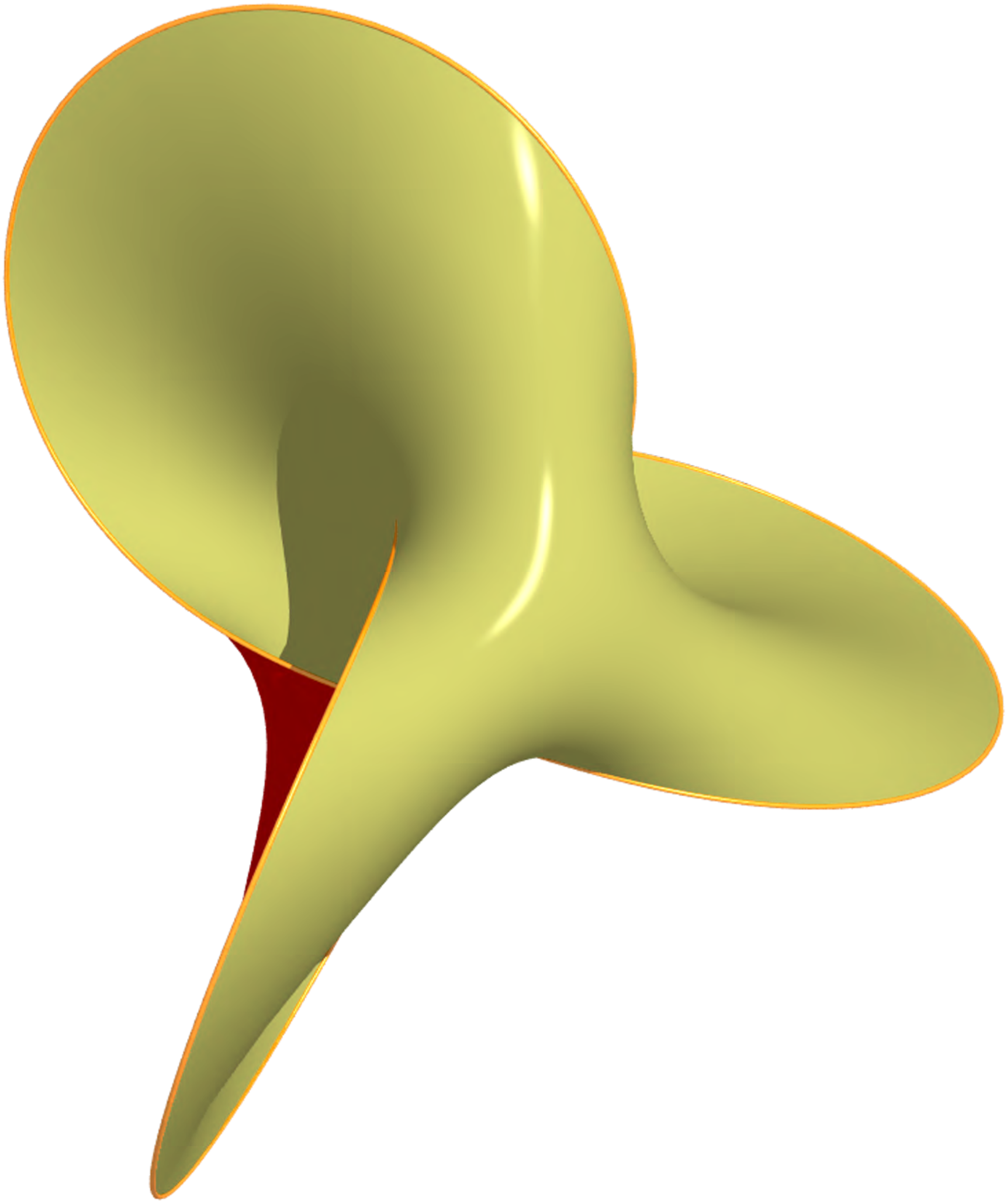

A page of the Milnor open book associated to the curve $y^2 - x^3 = 0$.

# *The cusp and the trefoil*

THE LOOSE PURPOSE OF OUR PROMENADE is to describe the topology of singularities of *real* analytic curves. As explained earlier, a shortcut through the *complex* domain might possibly shed some light on our "real" discussion. In any case, in this book we are more keen on detours than shortcuts.

For a full description of the topology of singularities of *complex* algebraic curves, I strongly recommend the excellent 721 page book[104] by Brieskorn and Knörrer. However:

> Un petit livre est rassurant.

as Jules Tannery wrote in the preface of a very concise and beautiful introduction[105] to Galois theory. Following this advice, I will limit myself to the basic features of the theory.

My only goal in this chapter is to convince the reader that the local topology of a singularity in the complex domain is incredibly rich.

[104] E. Brieskorn and H. Knörrer. *Plane algebraic curves*. Modern Birkhäuser Classics. Birkhäuser/Springer Basel AG, Basel, 1986.

[105] H. Vogt. *Leçons sur la résolution algébrique des équations*. Librairie Nony, 1895.

## *The link of a singularity*

The idea of intersecting a complex analytic curve $F(x,y) = 0$ by a small sphere $\mathbb{S}^3_\varepsilon = \{(x,y) \in \mathbb{C}^2 \mid |x|^2 + |y|^2 = \varepsilon^2\}$ is probably very old. The first paper explaining this construction is due to Brauner[106], published in 1928, following an idea of his PhD advisor Wirtinger in 1905. See[107] for an inspiring presentation of the historical development of these ideas.

[106] K. Brauner. Zur Geometrie der Funktionen zweier komplexer Veränderliche. *Abh. Math. Sem. Univ. Hamburg*, 6(1):1–55, 1928.

[107] M. Epple. Branch points of algebraic functions and the beginnings of modern knot theory. *Historia Math.*, 22(4):371–401, 1995.

It is important to recall that a curve over the complex numbers has dimension 1 over $\mathbb{C}$ and hence dimension 2 over $\mathbb{R}$, so that a complex curve is a real surface. This constant balance between curves and surfaces is one of the charms of the theory.

We already looked at the simplest case $F(x,y) = y - \lambda x$, leading to the Hopf fibration. We now look at the second significant example: the *cuspidal singularity*, defined by $F(x,y) = y^2 - x^3$.

How should we choose the small radius $\varepsilon$? Could we use some other hypersurface, like for example an ellipsoid? The answer is that under very mild assumptions, all these intersections define the *same* topological object, *up to homeomorphisms*. The case of $y^2 - x^3 = 0$ is particularly simple. Consider the following linear flow on $\mathbb{C}^2$:

$$\phi^t(x,y) = \left(e^{2t}x, e^{3t}y\right) \quad (t \in \mathbb{R}).$$

The space $\mathcal{O}$ of orbits of $\phi^t$ in $\mathbb{C}^2 \setminus \{(0,0)\}$ is homeomorphic to $\mathbb{S}^3$. Indeed, along such an orbit, the norm $|x|^2 + |y|^2$ is strictly increasing and each orbit intersects the sphere exactly once. The same argument could be used with an ellipsoid centered at the origin, or with our "square sphere" $\max(|x|,|y|) = \varepsilon$, or with many other hypersurfaces.

Now observe that the flow $\phi^t$ preserves our curve, whose equation is $y^2 - x^3 = 0$, so that the curve defines *canonically* a subset $K$ of $\mathcal{O}$. Identifying $\mathcal{O}$ with $\mathbb{S}^3_\varepsilon$, we realize that, up to homeomorphisms, the intersection of the curve with $\mathbb{S}^3_\epsilon$ is indeed independent of $\varepsilon$ and that we could as well use the square sphere.

So, let us intersect $y^2 - x^3 = 0$ with $\max(|x|,|y|) = \varepsilon$. If $\varepsilon < 1$, this is a parameterized by $\theta \in \mathbb{R}/2\pi\mathbb{Z}$

$$x = \varepsilon \exp(2i\theta) \quad ; \quad y = \varepsilon^{3/2}\exp(3i\theta)$$

in the solid torus (where $|x| = \varepsilon$ and $|y| \le \varepsilon$). This is the *trefoil knot*, seen as the $(3,2)$ *torus knot*: it is drawn on a standard torus of revolution in 3-space and goes three times around the meridian as it goes twice along the parallel.

There are many excellent books on the topology of knots. I recommend the "petit livre" by Sossinsky[108] and the very visual book by Kauffman[109].

In order to understand the topology of the cuspidal curve $y^2 - x^3 = 0$ in a small *ball* $|x|^2 + |y|^2 \le \varepsilon^2$, it suffices to note that all concentric spheres intersect the curve on such a trefoil. If follows that in a small ball, our curve is homeomorphic to the topological *cone* over the trefoil knot. The trefoil is a circle (embedded

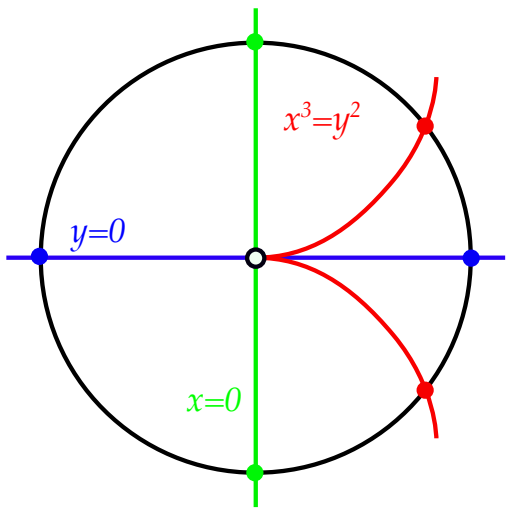

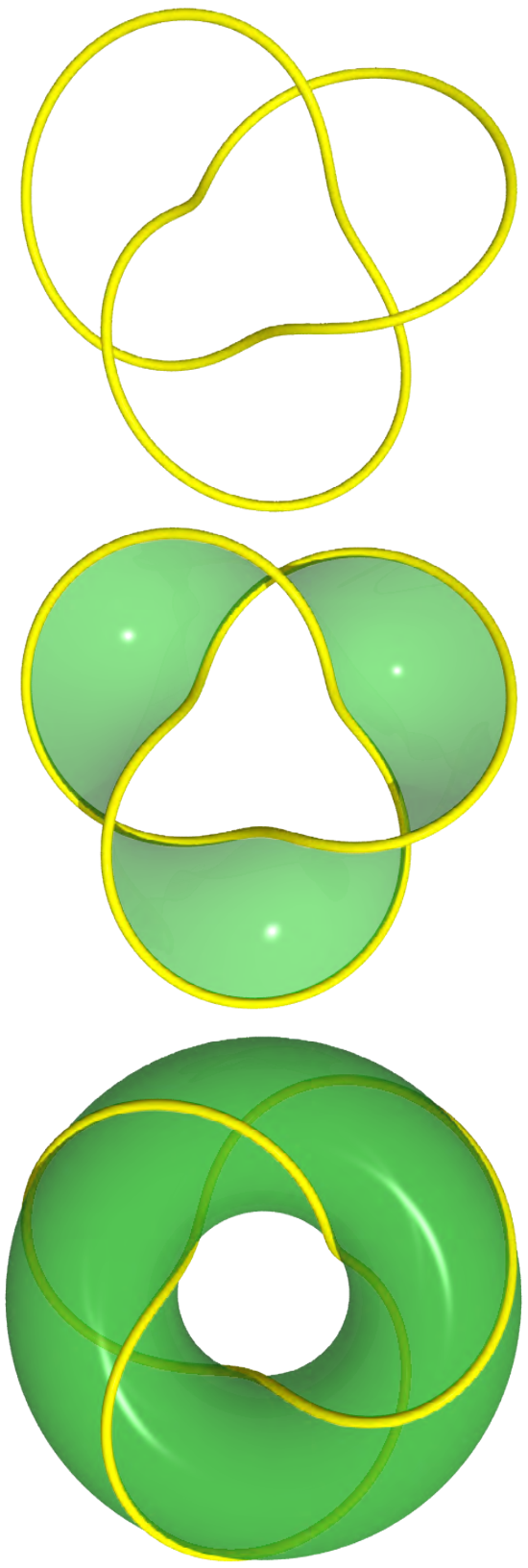

The trefoil knot. ©

108 A. Sossinsky. *Knots. Mathematics with a twist*. Harvard University Press, Cambridge, MA, 2002.

109 L. H. Kauffman. *On knots*, volume 115 of *Annals of Mathematics Studies*. Princeton University Press, Princeton, NJ, 1987.

in a knotted way). So the cone is a disc embedded in a tricky way into 4-space. Over the complex numbers the curve is topologically smooth, that is to say locally homeomorphic to a disc, but the embedding of this disc in $\mathbb{C}^2$ is knotted. This is a typical phenomenon that can be detected over the complex numbers and which is invisible over the reals, since the real curve $y^2 - x^3 = 0$, and indeed *every branch of a real analytic curve, is locally homeomorphic to a line in the plane*: this is what we called earlier Gauss's claim.

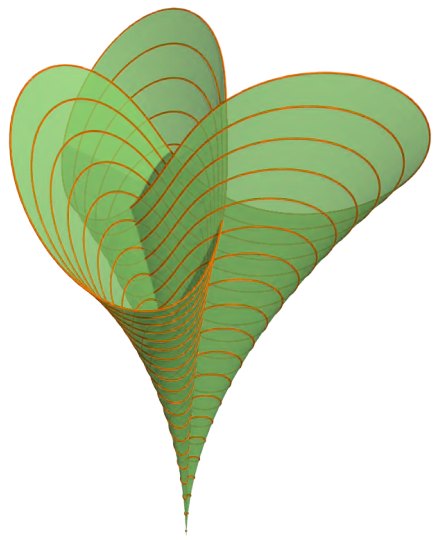

A local picture of a branch (projected in 3-space where the branch is not embedded).
©

## *Milnor's fibration*

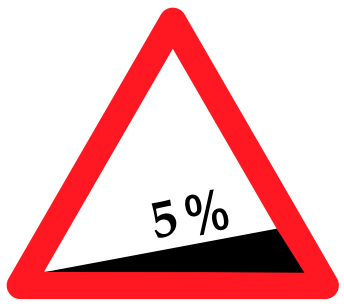

Let us describe the cuspidal curve in more detail and show a very special case of a general theorem of Milnor that will be presented later.

Consider the map

$$\mu : (x, y) \in \mathbb{S}^3 \mapsto y^2 - x^3 \in \mathbb{C}.$$

The inverse image of 0 is the trefoil knot. We want to look at the inverse image $\Sigma_\theta$ of a half line emanating from the origin, whose equation is $\arg(z) = \theta \in \mathbb{R}/2\pi\mathbb{Z}$. In other words, we look at the fibers of the map $\arg \circ \mu$ defined on the complement of the trefoil knot with values in the circle $\mathbb{R}/2\pi\mathbb{Z}$.

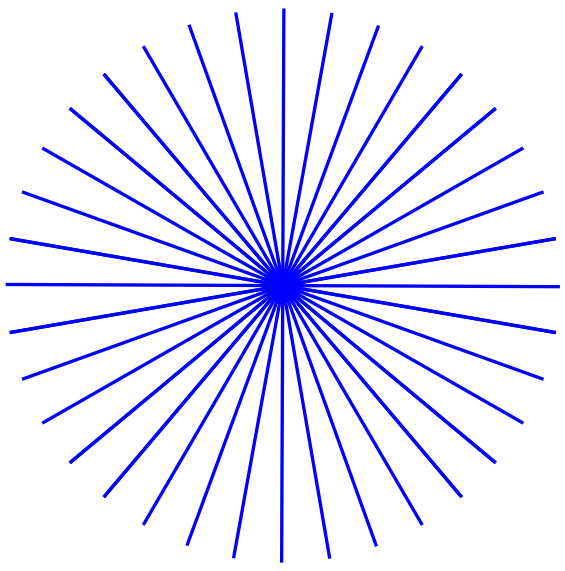

Radial half-lines with constant argument.

It is easy to see that $\arg \circ \mu$ is a submersion. Indeed, the flow

$$\psi^s(x, y) = \left(e^{2is}x, e^{3is}y\right)$$

preserves the spheres and satisfies $\arg(\mu \circ \psi^s) = \arg(\mu) + 6s$. It follows that the vector field associated to $\psi^s$ is not in the kernel of the differential of $\arg \circ \mu$. Observe that $\psi^s$ permutes the fibers $\Sigma_\theta$.

In the neighborhood of the trefoil knot, the situation is very easy to analyze. Our square sphere can still be used since $\arg \circ \mu$ is invariant under $\phi^t$ and we are in fact working in the orbit space $\mathcal{O}$. Parameterize a neighborhood of the trefoil by pairs $(\alpha, \zeta)$ with $\alpha \in \mathbb{R}/2\pi\mathbb{Z}$ and $\zeta$ a small complex number:

$$x = \varepsilon \exp(2i\alpha) \quad ; \quad y = \varepsilon^{3/2}\left(\exp(3i\alpha) + \exp(-3i\alpha)\zeta\right).$$

In these coordinates, $\arg \circ \mu$ is equal to $\arg(\zeta)$ to the first order. It follows that in the neighborhood of the trefoil, the $\Sigma_\theta$'s are surfaces which behave like the pages of a book around its binding.

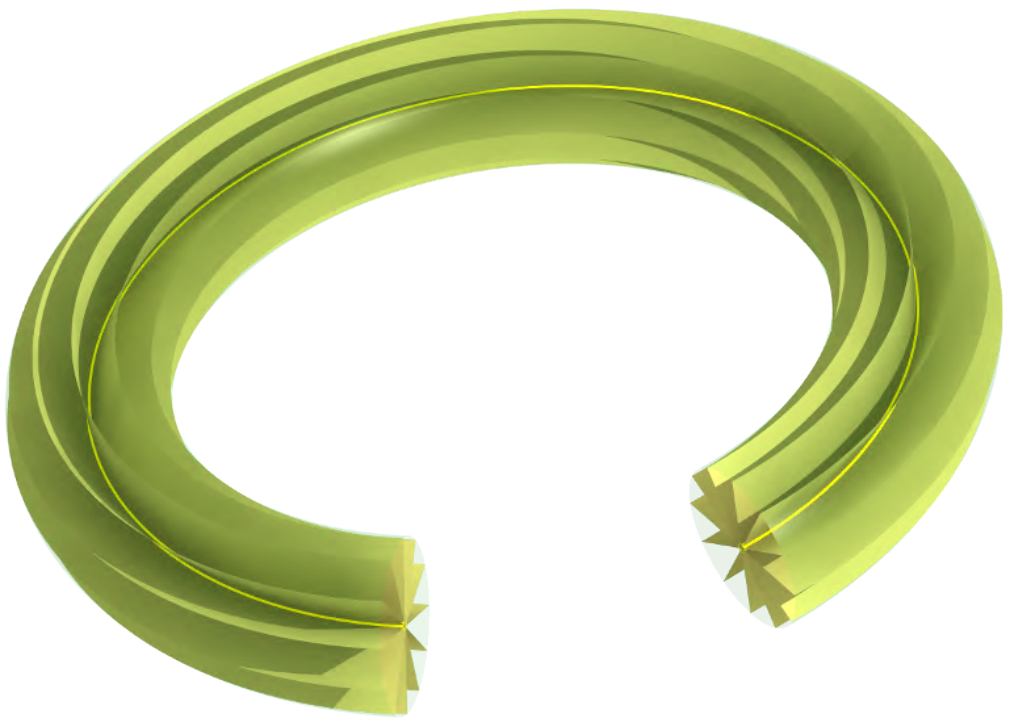

A strange book in which the pages are cyclically ordered and with no first page. A dream book that you read forever. At least the pages are orientable.

I will say that the trefoil knot is *fibered* or that its complement fibers over the circle. The fibers are disjoint *pages* whose closures in the 3-sphere all have the knot as their common boundary. Note that one page goes through the north pole in the 3-sphere, which is the center of the stereographic projection. This page, when projected in the Euclidean 3-space, is not compact.

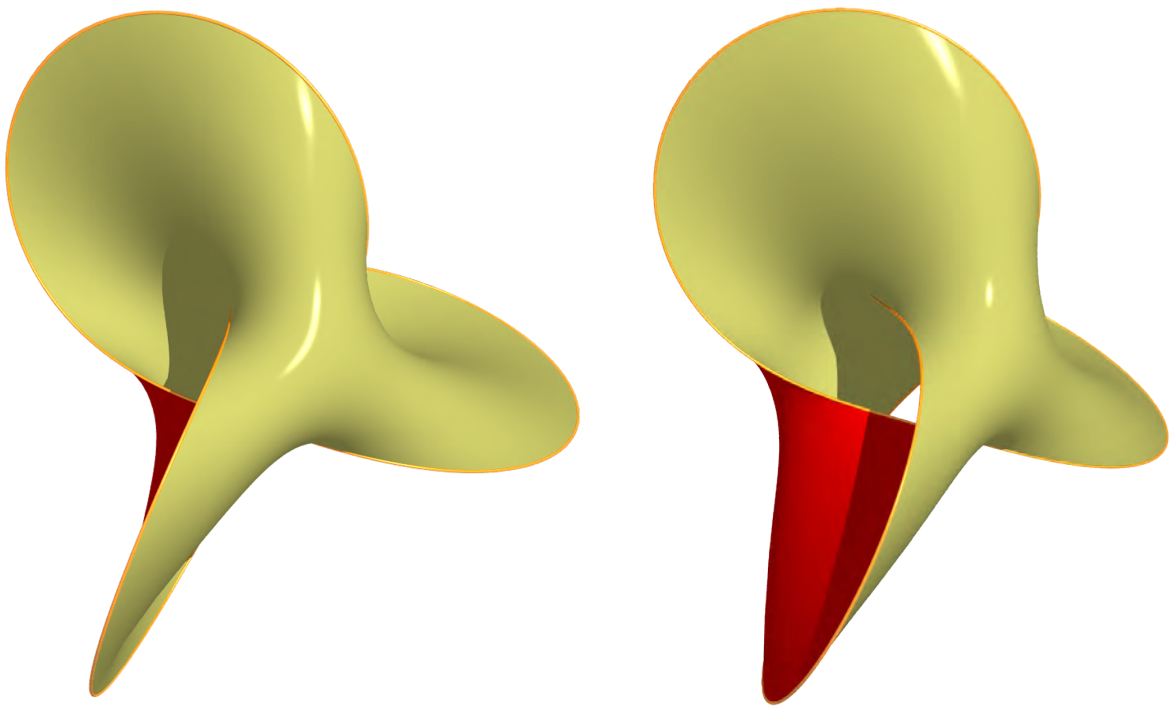

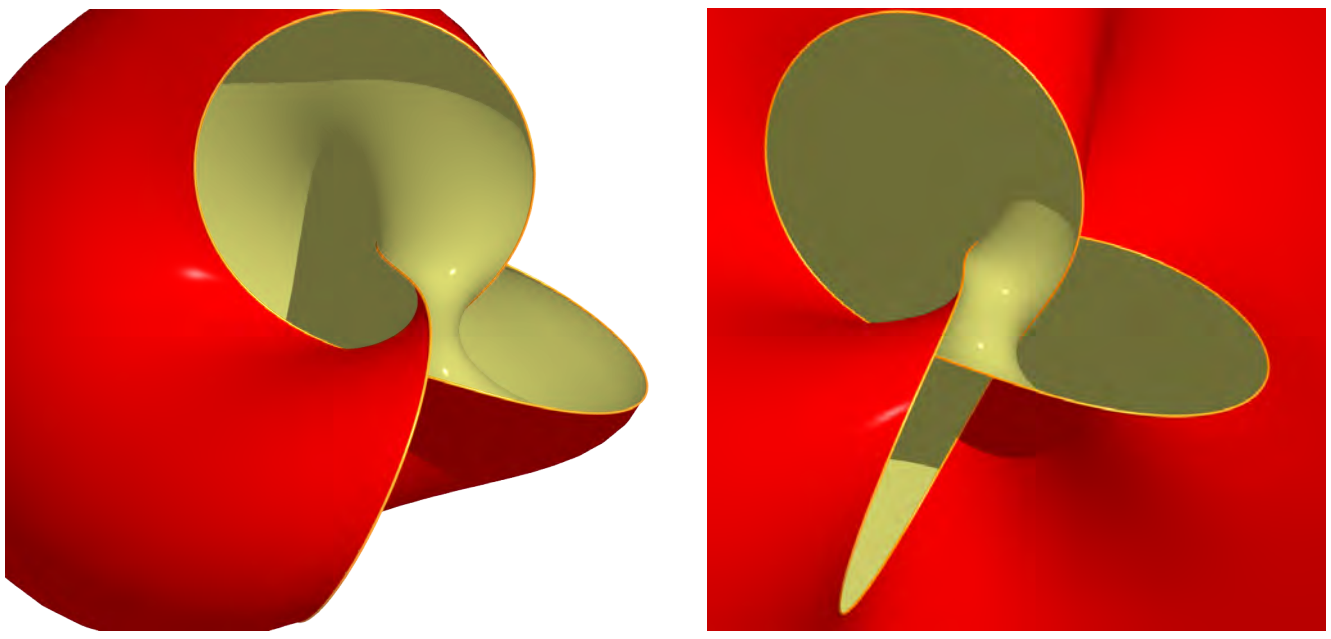

Some pages. ©

## *Monodromy*

The flow $\psi^s(x,y) = \left(e^{2is}x, e^{3is}y\right)$ permutes the pages of our book. More precisely, $\psi^s$ maps the page $\Sigma_\theta$ to the page $\Sigma_{\theta+6s}$. Note that $\psi^{2\pi}$ is the identity and that $\psi^{\pi/3}$ fixes globally each page, hence inducing on each page a homeomorphism of order 6, which is called the *monodromy* of the cusp.

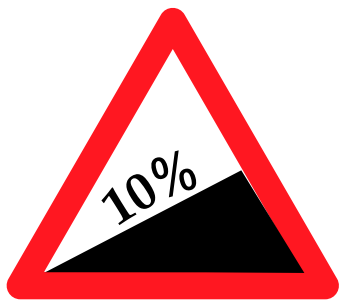

Our goal now is to describe the topology of the pages and the action of the monodromy.

By definition, a page $\Sigma$ is the set of $(x,y)$ in $\mathbb{S}^3$ such that the complex number $y^2 - x^3$ is in some half line, for instance the positive real axis $\mathbb{R}_+^* \subset \mathbb{C}$. Let $\mathcal{C}$ be the algebraic curve defined by $y^2 - x^3 = 1$ in $\mathbb{C}^2$. Think of the 3-sphere as the orbit space of the flow $\phi^t(x,y) = \left(e^{2t}x, e^{3t}y\right)$ acting on $\mathbb{C}^2 \setminus \{(0,0)\}$. The two real surfaces $\Sigma$ and $\mathcal{C}$ define the same object in this orbit space, so we work with $\mathcal{C}$. The action of the monodromy corresponds to

$$(x,y) \in \mathcal{C} \mapsto (\omega^2 x, \omega^3 y) \in \mathcal{C}$$

where $\omega = \exp(2i\pi/6)$ is a primitive 6-th root of unity.

The topology of $\mathcal{C}$ is easy to describe... if you know something about the genus of Riemann surfaces/algebraic curves. In $P^2(\mathbb{C})$, the homogenized cubic curve $y^2z - x^3 = z^3$ is a smooth elliptic curve intersecting triply the line at infinity in the point $[0:1:0]$. "Hence" the affine curve $\mathcal{C}$ is homeomorphic to a once punctured torus.

In a more down to earth way, one can proceed in the following manner. Set $Y = y^2 = 1 + x^3$, so that the map $(x,y) \in \mathcal{C} \mapsto Y$ is a six fold branched cover of $\mathbb{C}$, branched at 0 and 1, with order 2 and 3.

The points $(-1,0)$, $(\exp(2i\pi/6),0)$ and $(-\exp(2i\pi/6),0)$ are mapped to 0 and the points $(0,1)$ and $(0,-1)$ are mapped to 1.

Draw an arc in the complex plane connecting $Y = 0$ and $Y = 1$. It lifts to six arcs in $\mathcal{C}$.

Above $Y = 0$, there are 3 points, where the 6 arcs merge in 3 groups of 2. Above a small disc centered at 0, there are 3 *double plates*, as on the left part of the following picture. Above $Y = 1$, as in the second figure, there are 2 points, where the arcs merge in 2 groups of 3.

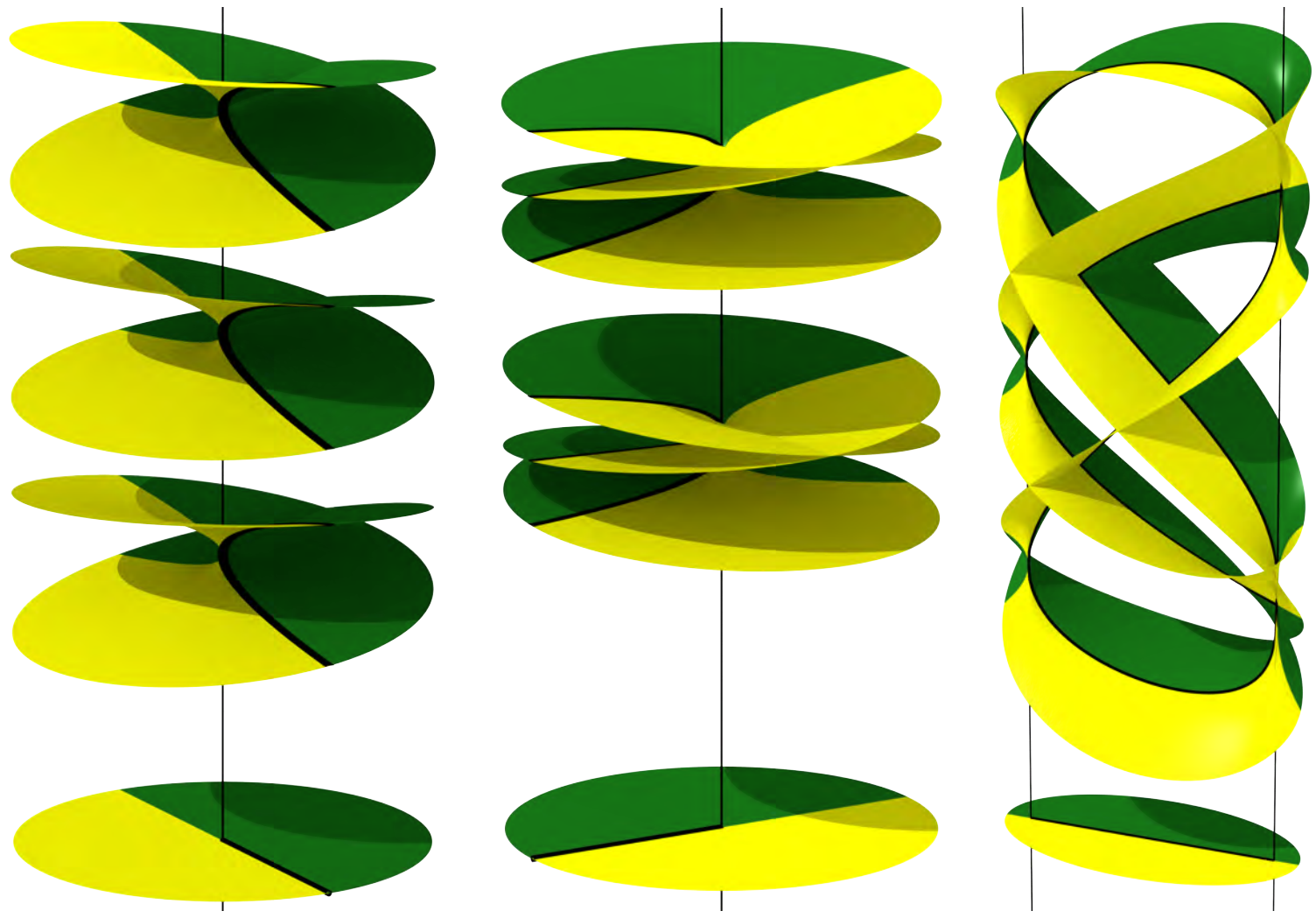

©

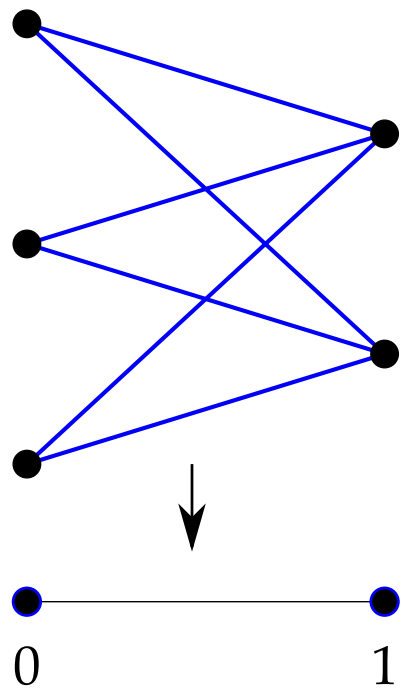

Since it is impossible to draw in the 4-dimensional $\mathbb{C}^2$, these pictures represent the graph of some appropriate combination of the real and imaginary parts of $\sqrt{Y} + \sqrt[3]{1-Y}$.

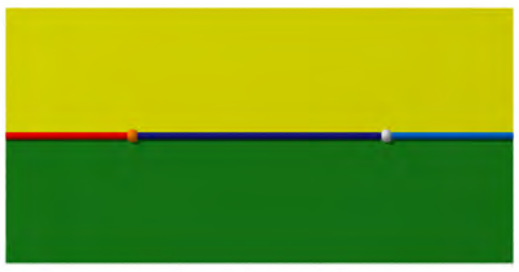

©

The combinatorics of the six arcs is represented in the right margin.

Cutting $\mathbb{C}$ along 3 arcs, from $-\infty$ to 0, from 0 to 1, and from 1 to $\infty$, we decompose $\mathbb{C}$ in two *triangles* where the imaginary part is positive or negative (yellow and green on the picture). These are indeed triangles with vertices at $0, 1$ and $\infty$. In $\mathcal{C}$, this produces in total 18 arcs, and $12 = 6 \times 2$ triangles.

Another way of seeing the same picture is the following. Identify opposite sides of a regular hexagon by translations. This defines a torus. Deleting the center, we get a punctured torus.

From the center of the hexagon, draw the 6 segments going to the vertices and the 6 heights to the sides. Our torus is now decomposed in 12 triangles, having in total 18 sides. The six roots of unity act by rotations on the (punctured) hexagon, permuting the triangles exactly as in the case of $\mathcal{C}$.

In summary, *each page of the book associated to the cusp is a punctured torus as above, and the monodromy is simply a rotation by 1/6th of a full turn.*

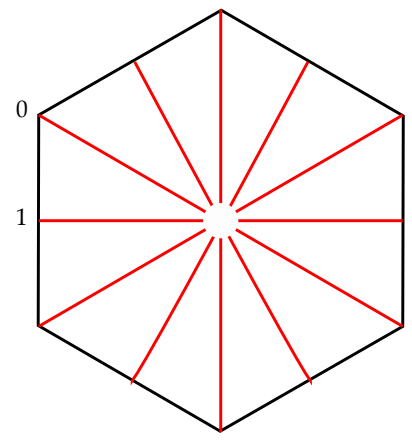

The more economical presentation of a torus from a square is more usual but this presentation with a hexagon is even more beautiful. Observe that the six sides define three arcs in the torus and the six vertices define two points in this torus.

## *Torus knots*

Most of what has been seen for the cusp $y^2 = x^3$ extends to a general curve $F(x, y) = 0$. This will require some work, but there is at least one family of examples where there is no extra work. Let $p, q$ be two relatively prime positive integers and let us look at the curve $y^p - x^q = 0$. We can assume $q > p$.

Just as before, let us look at the intersection with the square sphere.

$$x = \varepsilon \exp(ip\theta) \quad ; \quad y = \varepsilon^{q/p} \exp(iq\theta).$$

This is a $(p, q)$ torus knot $K_{p,q}$, drawn on a standard torus in 3-space, going around $p$ times the parallel and $q$ times the meridian.

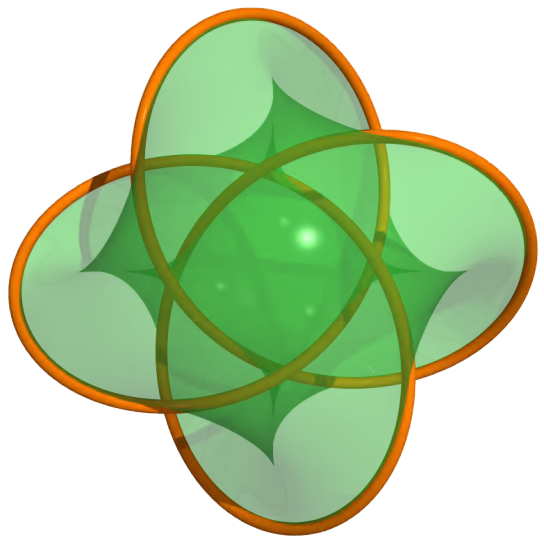

A $(3, 4)$ torus knot. ©

Exactly for the same reason, there is an open book decomposition and a fibration over the circle. Any page is homeomorphic to the affine algebraic curve $y^p - x^q = 1$, whose topology can be described in the same way. Set $Y = y^p$ and look at this curve as spread over $Y$, branched over 0 and 1. Over a point $Y$ different from $0, 1$, there are $pq$ points. Over 0 (resp. 1), there are only $q$ (resp. $p$) points, but each with multiplicity $p$ (resp. $q$). Replace the hexagon by a $pq$-gon and the situation is the same. There are $2pq$ triangles ($pq$ of each color) and $3pq$ edges. There are $p$ vertices above 0, $q$ above 1 and one above infinity.

This gives an Euler-Poincaré number equal to

$$p + q + 1 - 3pq + 2pq = 2 - 2\frac{(p-1)(q-1)}{2}.$$

Each page is now a punctured surface of genus $(p-1)(q-1)/2$.

I sketch now a proof of a fundamental fact: topology recovers a good part of the algebraic curve $y^p = x^q$.

I claimed in the preface that I tried to write a book that I could have understood myself as an undergraduate. I fear that this might not be the case for the end of this chapter. If this is too sketchy, just skip it!

**Theorem.** *If some homeomorphism of the 3-sphere sends the torus knot $(p_1, q_1)$ to $(p_2, q_2)$, then the sets $\{p_1, q_1\}$ and $\{p_2, q_2\}$ are equal.*

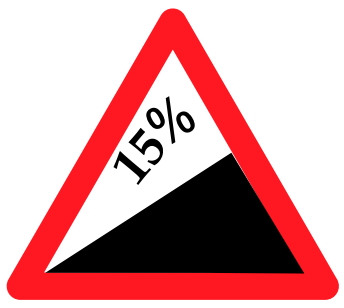

The proof will require some basic algebraic topology. From the topology of torus knots, we construct an algebraic gadget that will enable us to reconstruct $p, q$.

The complement of the torus knot $(p, q)$ in the 3-sphere is an open 3-manifold. Its most primitive invariant is its fundamental group, denoted by $\Gamma_{p,q}$. The key point is to extract algebraically $p$ and $q$ from this group. We will prove that $\Gamma_{p_1,q_1}$ is isomorphic to $\Gamma_{p_2,q_2}$ only if $\{p_1, q_1\} = \{p_2, q_2\}$.

Observe first that the map $\arg \mu = \arg(y^p - x^q)$ from the complement of the torus knot $(p, q)$ in the 3-sphere to $\mathbb{R}/2\pi\mathbb{Z}$ induces a surjective homomorphism

$$\lambda : \Gamma_{p,q} \to \pi_1(\mathbb{S}^1) \simeq \mathbb{Z}$$

between fundamental groups. Indeed, consider the following loop in the unit sphere defined for $t \in [0, 1]$ by

$$x(t) = \frac{\sqrt{2}(1 + \zeta(t))}{2|(1 + \zeta(t)|)} \exp(2ip\pi t) \quad ; \quad y = \frac{\sqrt{2}}{2} \exp(2iq\pi t).$$

Beginners are strongly encouraged to look at the remarkable website by Henri Paul de Saint Gervais dedicated to Analysis Situs.

For small values of $\zeta$, this is a small loop going around the $(p, q)$ knot. The argument of $\mu$ on this loop is close to

$$\pi + \arg \zeta(t) + 2\pi pqt$$

so that we can choose $\zeta(t) = \varepsilon \exp(i(1 - 2pqt)\pi)$ to make sure that the image of this loop by $\lambda$ is 1.

In the *first step*, we show that this homomorphism, up to sign, is the *only* surjection of $\Gamma_{p,q}$ onto $\mathbb{Z}$. It follows that the kernel of $\lambda$ only depends on the topology of the knot.

In the *second step*, we analyze the abelianization of the kernel of $\lambda$, show that it is a finitely generated free Abelian group, and find its rank. This will enable us to recover $\{p, q\}$ from the group $\Gamma_{p,q}$ as required.

Any map from $\mathbb{S}^1$ to $\mathbb{S}^3$ extends to the unit disk $D^2$. If this extension is an embedding, then $\gamma$ is a trivial knot. However, any knot is the boundary of an embedded oriented surface (of higher genus): this is called a Seifert surface.

The first step could be explained in a variety of ways, more or less sophisticated, most of them based on the so-called *Lefschetz duality*. Suppose some closed orientable manifold $X$, for example a circle, is embedded in some sphere, for example of dimension 3. Then the homology of the complement of $X$ does not depend on the way $X$ is embedded in the sphere. In particular, the homology of the complement of a knot in the 3-sphere is the same as in the case of a trivial knot, so this homology is simply isomorphic to $\mathbb{Z}$ in degree 1.

Note that the complement of an unknotted circle in the 3-sphere is homeomorphic to $\mathbb{R}^2 \times \mathbb{S}^1$.

I could present the same fact in the following way. Let $\gamma$ be a smooth loop in $\mathbb{S}^3 \setminus K_{p,q}$. Since the sphere is simply connected, $\gamma$ is the boundary of some smooth map $D \to \mathbb{S}^3$ which may not be an embedding. Put this disk in general position with $K_{p,q}$, so that the intersections between $K_{p,q}$ and $D$ are transversal. Count the number of intersections between the disk $D$ and $K_{p,q}$, the counting being algebraic, taking orientations into account. This number is the *linking number* $lk(\gamma)$. It turns out that it only depends on the homology class of $\gamma$ in $\mathbb{S}^3 \setminus K_{p,q}$. This follows from the fact that a surface with no boundary in $\mathbb{S}^3$ has a trivial algebraic intersection with any closed curve.

More precisely this is the linking number of $\gamma$ and $K_{p,q}$. This will be discussed later in this book.

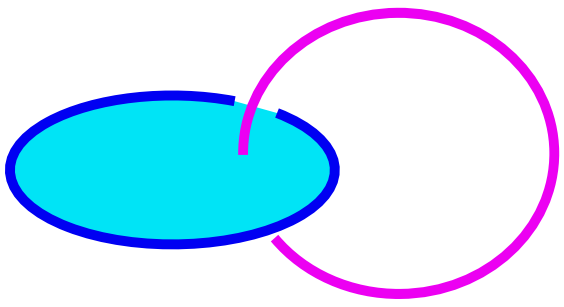

Linking number.

Therefore it defines a homomorphism

$$lk : H_1\left(\mathbb{S}^3 \setminus K(p,q), \mathbb{Z}\right) \to \mathbb{Z}$$

which is onto. Now, if $\gamma$ is in the kernel of $lk$, this means that $+$ and $-$ signs in the intersection can be coupled. Dig holes in $D$, around the intersection points, and connect their boundaries in pairs with tubes, in order to construct a surface whose boundary is still $\gamma$ and which does not intersect $K(p,q)$ anymore. Therefore the elements of the kernel of $lk$ are homologous to zero. In other words, $lk$ is an isomorphism. Finally recall that the first homology group is the abelianization of the fundamental group, so that any homomorphism $\pi_1\left(\mathbb{S}^3 \setminus K\right) \to \mathbb{Z}$ factors through $lk$. It follows in particular that $lk$ coincides (up to sign) with the previously defined $\lambda$. This is the first step.

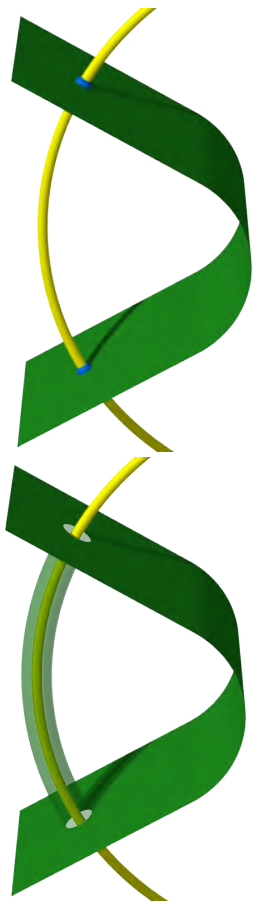

Removing two intersection points.

*We now proceed to the second step.* Denote the kernel of $lk$ by $G(p,q)$. It is therefore the fundamental group of some Galois covering of $\mathbb{S}^3 \smallsetminus K(p,q)$ whose group of automorphisms is infinite cyclic $\mathbb{Z} = \Gamma_{p,q}/\ker lk$. This covering is clearly the product $\Sigma \times \mathbb{R}$ of a page with $\mathbb{R}$ and the group of deck transformations is simply generated by

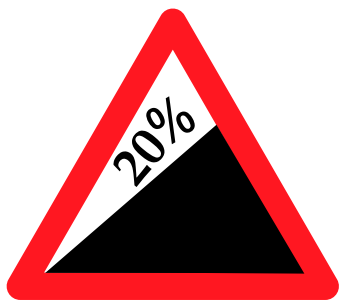

$$(p,t) \in \Sigma \times \mathbb{R} \to (M(p), t + 2\pi) \in \Sigma \times \mathbb{R}$$

where $M$ denotes the monodromy map. It follows that $G(p,q)$ is the fundamental group of a page $\Sigma$. We already described the topology of a page. Since $G(p,q)$ is not abelian, it might be easier to make it abelian. Denote by $H(p,q)$ this abelianization, which is nothing more than the first homology of a page.

I recall that $p$ and $q$ are relatively prime.

Let us describe this abelian group $H(p,q)$ and the action of $M$. In $\Sigma$ there is a graph containing $pq$ arcs, obtained by lifting the arc connecting 0 and 1. It contains $q$ vertices over 0 and $p$ vertices over 1. Recall that $\Sigma$ is obtained from a closed triangulated surface by deleting a vertex which is common to all triangles. Therefore, the punctured surface $\Sigma$ can be deformed to the union of all the edges opposite to this vertex which is our graph with $pq$ edges. This graph is usually called a *complete bipartite graph*. This produces a very simple 1-complex which computes $H(p,q)$.

The abelian group of 1-chains is freely generated by arcs $c_{i,j}$ where $i \in \mathbb{Z}/p\mathbb{Z}$ and $j \in \mathbb{Z}/q\mathbb{Z}$. The abelian group of 0 chains is generated by $p$ points $a_i$ and $q$ points $b_j$ with $i \in \mathbb{Z}/p\mathbb{Z}$ and $j \in \mathbb{Z}/q\mathbb{Z}$. The boundary operator $\partial$ sends $c_{i,j}$ to $b_j - a_i$. Finally the monodromy $\mathbb{Z}/pq\mathbb{Z} \simeq \mathbb{Z}/p\mathbb{Z} \times \mathbb{Z}/q\mathbb{Z}$ acts in an obvious way on the indices $i,j$.

The homology $H(p,q)$ fits into an exact sequence

$$0 \longrightarrow H(p,q) \longrightarrow \mathbb{Z}^p \otimes \mathbb{Z}^q \xrightarrow{\partial} \mathbb{Z}^p \oplus \mathbb{Z}^q \longrightarrow \mathbb{Z} \longrightarrow 0$$

which is equivariant with respect to actions of $\mathbb{Z}/p\mathbb{Z} \times \mathbb{Z}/q\mathbb{Z}$ at each level. The generator $M$ of monodromy is associated to the action of $(1,1)$. Tensoring by $\mathbb{R}$ to get vector spaces and linear maps, it follows that the dimension of $H(p,q) \otimes \mathbb{R}$ is $pq - (p+q) + 1$, i.e. $(p-1)(q-1)$. The characteristic polynomial of the action $M_\star$ of $M$ on $H(p,q)$ can even be computed using the exact sequence:

$$P(X) = \frac{(X^{pq}-1)(X-1)}{(X^p-1)(X^q-1)}.$$

Observe that the roots of $P$, eigenvalues of $M_\star$, are the $pq$-th roots of unity, minus the $p$-th and $q$-th root, plus 1. From this spectrum, the values of $p,q$ can be extracted.

The proof of the theorem is finished. From the fundamental group $\Gamma$ of the (complement of the) knot $K(p,q)$, construct its first derived group $\Gamma_1 = [\Gamma,\Gamma]$, which is also, as we have seen, the kernel of $lk$. Then, make $\Gamma_1$ abelian and define the group $\Gamma_1/[\Gamma_1,\Gamma_1]$. Now consider some element $g$ in $\Gamma$ with $lk(g) = \pm 1$ and the conjugation by $g$ on $(\Gamma_1/[\Gamma_1,\Gamma_1]) \otimes \mathbb{R}$. The values of $p,q$ can be obtained from the eigenvalues of this linear map. ⊡

This algebraic trick is actually a very general and powerful technique and is not restricted to knots. Given any group $\Gamma$, look at the action of the abelianization $\Gamma_{ab} = \Gamma/\Gamma_1 = \Gamma/[\Gamma,\Gamma]$ by conjugation on the abelianization $(\Gamma_1/[\Gamma_1,\Gamma_1]) \otimes \mathbb{R}$. This defines a family of commuting automorphisms whose conjugacy classes are invariants of the group $\Gamma$. One speaks of the *Alexander module* of $\Gamma$. This is one of the most primitive invariants of a group.

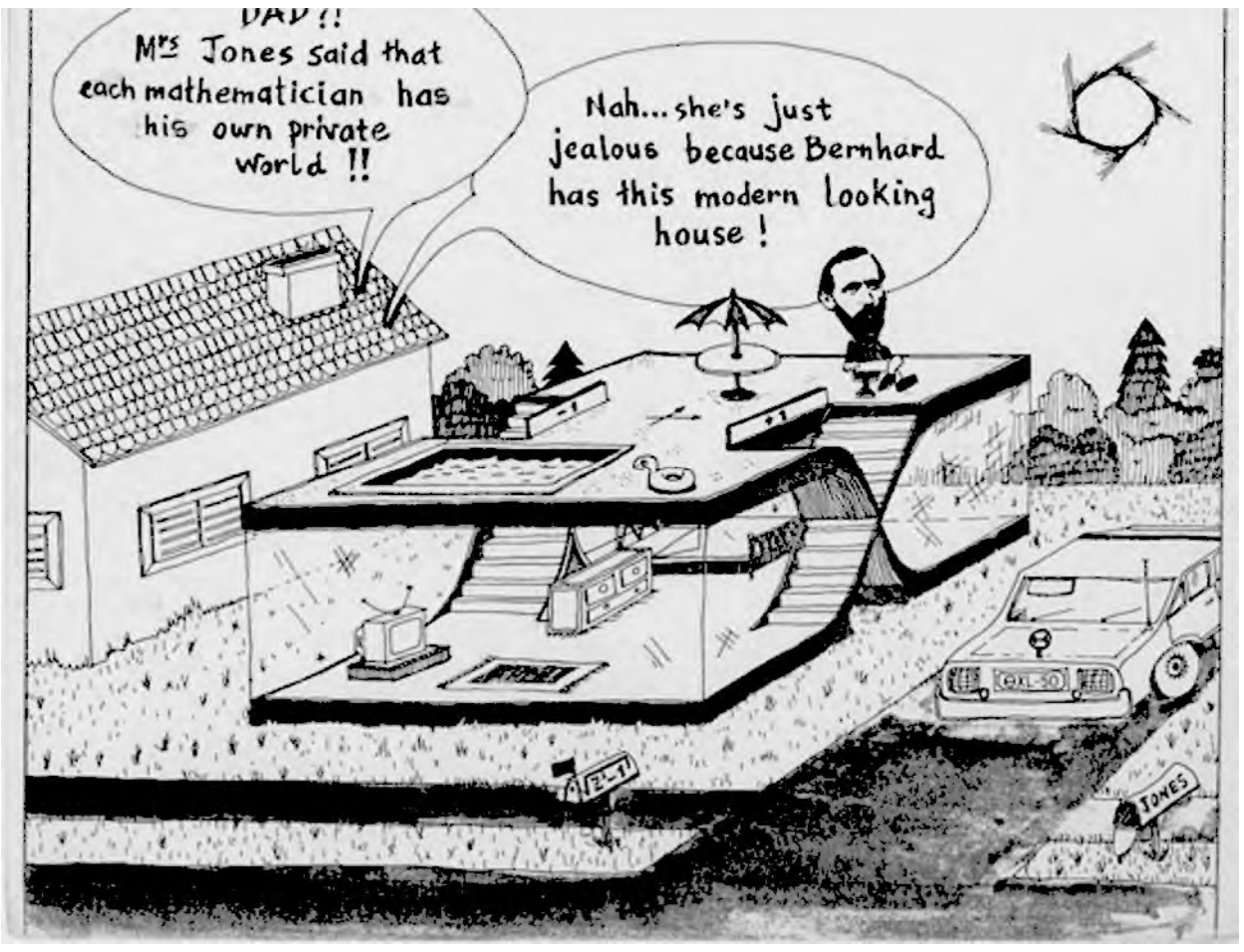

On André Nachbin's website.

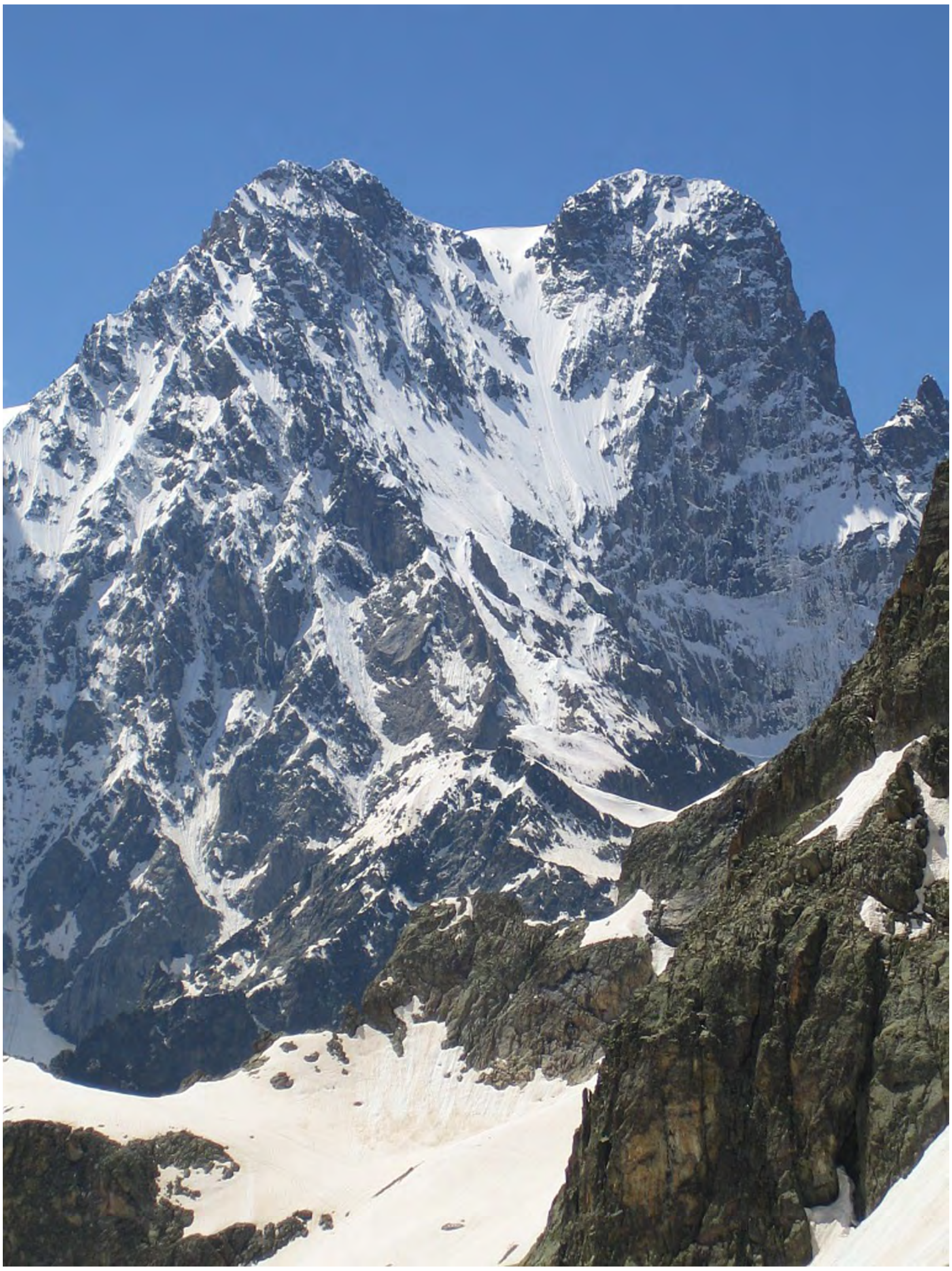

Pointe Puiseux. ©

# *Victor Puiseux, at last!*

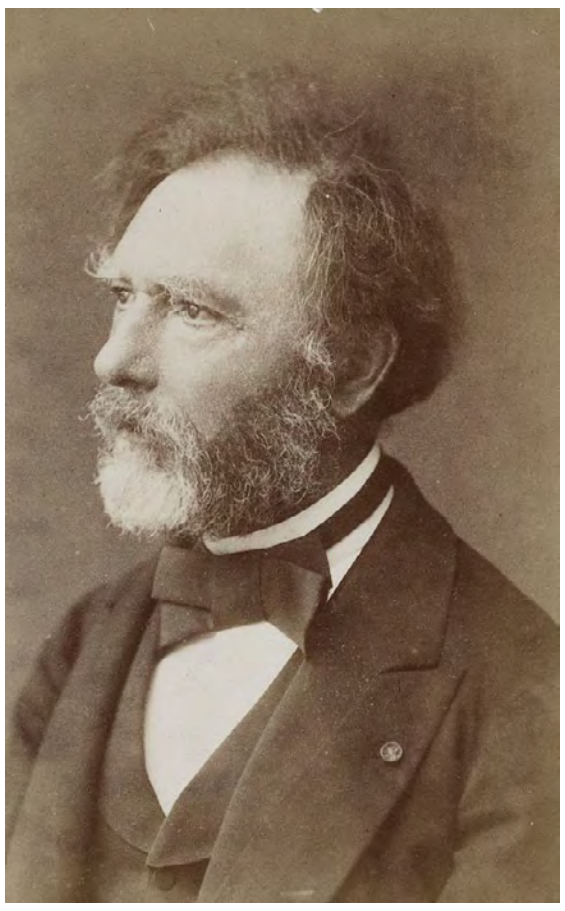

Victor Puiseux (1820–1883) ©

The name of Puiseux already appeared several times in this book. The reader may be anxious to know what he actually did. Unfortunately, the "well known Puiseux theorem" is not due to him but, as we have seen, to Newton, with some posterior help from Cramer. One could argue that neither Newton nor Cramer proved the convergence of the associated series but this convergence can be easily proved, for example using the *calcul des limites* of Cauchy.

However, Puiseux approached the problem of the local structure of singularities in a totally different way and his contribution is fundamental. In this chapter, I would like to explain his point of view. Unfortunately, it would be useless to stick to his original presentation.

Strange fate for a mathematician: he is "famous" for a theorem that was known long before him, and that we understand today much better than he did, using techniques that came long after him.

Fortunately, Puiseux is even more famous among alpinists since the highest peak of the Mount Pelvoux (3,946 m), in the Massif des Écrins, is called *pointe Puiseux*. He reached this peak on August 9, 1848. Unfortunately, it is not even sure that this was "a first" since Captain Durand claimed that he reached the summit 18 years earlier. Eternal second?

## *Puiseux's topological approach*

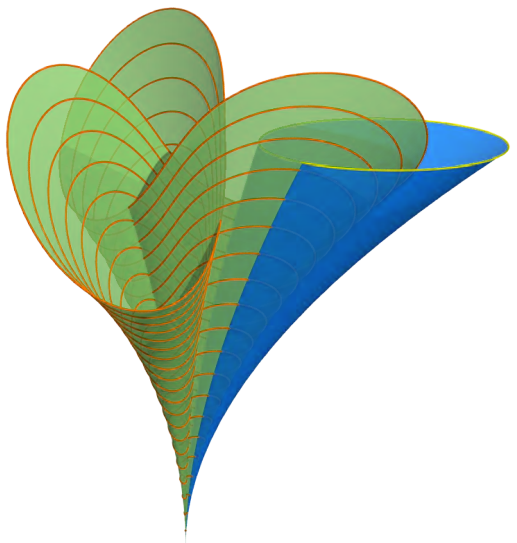

Two branches can be linked.

Let us recall what is usually called Puiseux's theorem.

**Theorem.** *Let $F(x,y)$ be a nonzero holomorphic function defined in the neighborhood of the origin in $\mathbb{C}^2$ and such that $F(0,0) = 0$. Then, there exists a finite number of holomorphic functions $g_1, \ldots, g_n$ defined in the neighborhood of $0 \in \mathbb{C}$ and positive integers $m_1, \ldots, m_n$ such that the curve $F(x,y) = 0$, again in the neighborhood of $(0,0)$, is the union of $n$ branches $t \mapsto (t^{m_i}, g_i(t))$ (for $i = 1, \ldots, n$) (plus, possibly, the $y$-axis). Moreover these branches are injective and they only intersect at the origin.*

We have already discussed a pre-Puiseux proof, very algebraic in spirit, where one finds first the formal series $g_i$ before proving that they converge. Puiseux proposed a topological approach[110] in 1850, just before the great papers by Riemann introducing topological ideas in algebraic geometry. We should therefore "forgive" him since, of course, he could not express himself in terms of Riemann surfaces.

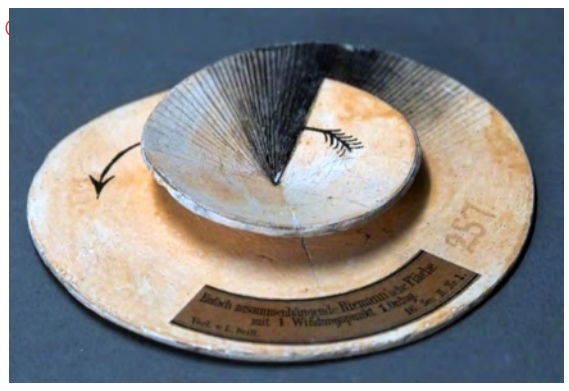

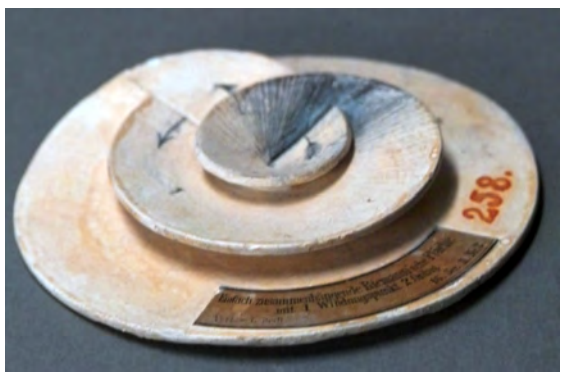

Two models from the Göttingen Collection of Mathematical Models and Instruments that the reader should definitely visit.

Let me sketch such a topological proof. Consider first $F(0,y)$. If this is identically zero, $F$ can be divided by some power of $x$ without changing the problem. We can therefore assume that the valuation of $F(0,y)$ (also called the multiplicity) is some positive integer $m > 0$. In particular $F(0,y)$ has an isolated zero at the origin (of multiplicity $m$). Choose some $\varepsilon > 0$ such that 0 is the only root of $F(0,y) = 0$ in $|y| \leq \varepsilon$. By a simple continuity argument, there is some $\eta > 0$ such that there is no root of $F(x,y) = 0$ on the solid torus $\{(x,y) \mid |x| \leq \eta\,; |y| = \varepsilon\}$. Dividing $x$ and $y$ by $\varepsilon$ and $\eta$, we assume that $\varepsilon = \eta = 1$.

[110] V. Puiseux. Recherches sur les fonctions algébriques. *Journal de Mathémagiques Pures et Appliquées*, 15:365–480, 1850.

Let us now make some assumption that will be analyzed in detail later on.

*Assume that the partial derivative $\partial F / \partial y$ does not vanish on the curve $F(x,y) = 0$, except at the origin.*

Do not forget that all this discussion is local so when I write "does not vanish", I mean "does not vanish in some neighborhood of the origin".

Denote by $\mathcal{C}^\star$ the punctured curve

$$\{(x,y) \in \mathbb{C}^2 \mid (x,y) \neq (0,0)\,; F(x,y) = 0\,; |x| \leq 1\,; |y| \leq 1\}.$$

The main assertion is that *the projection of $\mathcal{C}^\star$ onto the punctured disc $D^\star = \{x \mid |x| \leq 1\} \setminus \{0\}$ is a covering map.*

Again, do not miss Analysis Situs, by Henri Paul de Saint Gervais, available online.

Let me recall quickly the definition of covering maps and how they differ from a local homeomorphisms.

The word *étale* was introduced by French algebraic geometers. It means "stationary" and is frequently used to describe the surface of the sea, when at rest.

A continuous map $p : X \to Y$ is a *local homeomorphism* (sometimes called an étale map) if every point in $X$ has an open neighborhood $U$ such that $p(U)$ is open and the restriction of $p$ to $U$ is a homeomorphism onto $p(U)$.

A continuous map $p : X \to Y$ is a *covering map* if every point in $Y$ has an open neighborhood $V$ such that $p^{-1}(V)$ is a disjoint union of open sets $U_i$ such that the restriction of $p$ to each $U_i$ is a homeomorphism onto $V$.

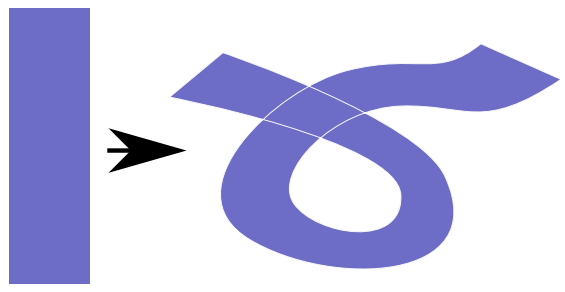

An étale map which is not a covering.

Clearly, a covering map is a local homeomorphism, but simple examples show that the converse is not true. One shows easily that a local homeomorphism is a covering space if it is *proper*.

A map is *proper* if the inverse image of a compact set is compact. Of course the converse is not true and a covering map might be non-proper (like for instance $t \in \mathbb{R} \mapsto \exp(it) \in \mathbb{S}^1$).

Let us now show our assertion that $\mathcal{C}^\star$ is a covering of the punctured disc. The fact that the projection is a local homeomorphism follows immediately from our assumption that $\partial F/\partial y$ does not vanish on $\mathcal{C}^\star$ and from the implicit function theorem. The properness of the projection is clear as well since a sequence of points on $\mathcal{C}^\star$ escapes from a compact set if and only if it converges to the origin.

The main theorem of covering space theory is that the connected covering spaces of a (locally simply connected) connected space are described, up to isomorphisms, by the subgroups of the fundamental group. For instance, connected covering spaces of $D^\star$ are isomorphic to some power map $x \in D^\star \mapsto x^m \in D^\star$, for some integer $m \geq 1$, or to the complex exponential map restricted to the half plane $\Re(x) \leq 0$.

Actually Puiseux used implicitly coverings when he described some loops followed by $x$ around the origin and the associated permutation of the values of $y$ satisfying $F(x,y) = 0$.

Choose some connected component $\mathcal{C}_0^\star$ of $\mathcal{C}^\star$. Since the covering $\mathcal{C}_0^\star \to D^\star$ has finite fibers, it is isomorphic to some covering $x \in D^\star \mapsto x^m \in D^\star$. Said differently, there is some *homeomorphism*

$$\phi : x \in D^\star \mapsto (x^m, g(x)) \in \mathcal{C}_0^\star.$$

This $\phi$ is clearly holomorphic on the punctured disc and we still have to show that it extends as a holomorphic function in the disc. This follows from the Riemann extension theorem: a bounded holomorphic function on a punctured disc is holomorphic in the full disc.

One more anachronism!

The theorem is proved, *under the assumption that the partial derivative $\partial F/\partial y$ does not vanish on the punctured curve $F(x,y) = 0$*, that will be discussed in the next paragraph.

## *Simple roots*

A holomorphic function of one complex variable $y$ and its derivative vanish simultaneously at some $y_0$ if and only if this zero is multiple. We therefore have to show that in Puiseux's theorem, one can always assume that $F$ has the property that, for $x_0$ small and nonzero, there are no small multiple roots of $F(x_0, y) = 0$.

It turns out that Puiseux did not consider general holomorphic functions $F(x,y)$ but polynomials in $x, y$. In this case, it is easy to deal with multiple roots and actually Puiseux dismisses the problem in one sentence (in a 135 page paper). Let us be just a little more careful than him.

Consider the *polynomial* $F$ as an element of $\mathbb{C}[x][y]$. This $F$ can be seen as a polynomial in one variable $y$ with coefficients in a factorial ring. We can write $F$ as a product of irreducible factors so that the curve $F(x,y) = 0$ is the union of the curves associated to these irreducible factors. We can therefore assume that $F$ is irreducible.

Suppose now that there is a sequence $(x_k, y_k)$ converging to $(0,0)$ with $x_k \neq 0$ and such that $y_k$ is a multiple root of $F(x_k, y) = 0$. Then the discriminant of the polynomial $F(x_k, y)$ is equal to 0. Therefore the discriminant of $F$, as an element of $\mathbb{C}[x]$, vanishes identically since it has an infinite number of roots. If the discriminant of some polynomial $P$ vanishes, the polynomial and its derivative have a common factor. This is impossible if $P$ is irreducible.

The discriminant of a polynomial is the resultant of this polynomial and its derivative.

Hence, if $F$ is irreducible in $\mathbb{C}[x][y]$, and if $x$ small and nonzero, then $F(x,y) = 0$ has no small multiple root as an equation in $y$. This is the ingredient that was missing for the proof of Puiseux theorem, *for a polynomial equation $F(x,y) = 0$*.

For a general holomorphic function $F(x,y) = 0$ (that, once again, Puiseux did not consider) there is still some work to be done.

## *Weierstrass's preparation theorem*

We already met Weierstrass's preparation theorem that we proved first in the context of formal series before establishing the convergence. The goal now is to prove the same theorem using complex analysis.

Let us recall the statement.

**Theorem.** *Let $F(x,y)$ be a nonzero holomorphic function defined in some neighborhood of the origin in $\mathbb{C}^2$. Then there exist $m$ holomorphic functions $a_0(x),\dots,a_{m-1}(x)$ defined in some neighborhood of $0 \in \mathbb{C}$, a holomorphic function $U(x,y)$ which is not vanishing at the origin, and an integer $r \geq 0$, such that*

$$F(x,y) = x^r U(x,y) \left(y^m + a_{m-1}(x)y^{m-1} + \dots + a_1(x)y + a_0(x)\right).$$

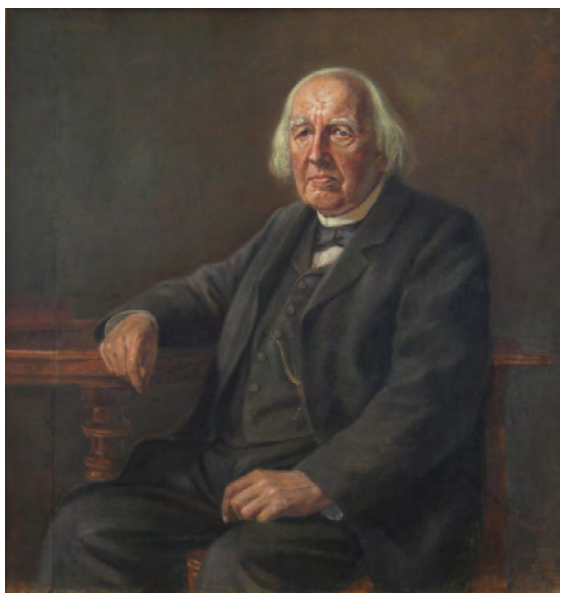

Karl Weierstrass (1815–1897). ©

This theorem is exactly what is needed. It states that up to non-vanishing functions, we can always assume that the function $F$ under study is a polynomial in the variable $y$, with coefficients in the ring $\mathbb{C}\{x\}$ of convergent series in $x$. The previous proof by Puiseux (that one can always assume that $\partial F/\partial y$ does not vanish identically on the curve, except at the origin), can therefore be reproduced word by word (replacing the ring of polynomials in $x$ by the ring of convergent series). Therefore, Puiseux's theorem is proved, using Weierstrass theorem.

I now present the standard analytical proof of Weierstrass.

Assume, after dividing $F$ by some $x^r$, that $F(x,y)$ does not vanish for $|x| \leq 1$ and $|y| = 1$. Fixing $x$ with $|x| \leq 1$, the function $y \mapsto F(x,y)$ has a finite number of zeros $y_1(x), y_2(x), \dots, y_m(x)$ in the unit disc, counted with multiplicity. The main difficulty is that it is impossible to choose these functions $y_i(x)$ as holomorphic functions of $x$, or even continuous, precisely because of the multivaluedness of the implicit $y(x)$ in $F(x,y) = 0$. However, we will show that all the symmetric functions of the $y_i(x)$'s are indeed holomorphic functions of $x$.

The simplest proof uses Cauchy formula. Let us evaluate

$$s_k(x) = \frac{1}{2i\pi} \int_{|y|=1} \frac{y^k F'_y(x,y)}{F(x,y)}\, dy.$$

The residue of $y^k F'_y(x,y)/F(x,y)$ as a function of $y$, at one of the roots $y_i(x)$, is the $k$-th power $y_i(x)^k$, so that $s_k(x)$ is the sum of the $k$-th powers of the roots. The integral shows clearly that $s_k(x)$ is a holomorphic function of $x$.

$s_0$ is the number of roots. Being an integer and a holomorphic function of $x$, it is constant. It was used implicitly a few lines above!

Since the $s_k$'s generate the symmetric functions, all symmetric functions of the $y_i(x)$ are holomorphic functions of $x$, in particular the elementary symmetric functions $a_i(x)$. By Viète's theorem, the polynomial

Another well known theorem of Newton.

$$y^m - a_{m-1}(x)y^{m-1} + \cdots + (-1)^{m-1}a_1(x)y + (-1)^m a_0(x)$$

vanishes exactly at the same points as $F$ with the same multiplicities, so that the quotient $U(x,y)$ does not vanish.

The Weierstrass preparation theorem, and Puiseux's theorem are proved. ⊡

## *Who proved Weierstrass's preparation theorem?*

My reader should have already guessed that the simple answer to this question is certainly not Weierstrass. Historians of mathematics know very well that questions like "who proved this first?" are far too naive, and frequently miss the point. It is nevertheless interesting to notice that two important mathematicians of the twentieth century, Henri Cartan[111] and Carl Siegel[112], wrote detailed papers trying to unfold the development of ideas around this theorem. Their papers are however not completely convergent. Let me only mention some steps.

111 H. Cartan. Sur le théorème de préparation de Weierstrass. In *Festschr. Gedächtnisfeier K. Weierstrass*, pages 155–168. Westdeutscher Verlag, Cologne, 1966.

– The fact that the symmetric functions of the roots of some holomorphic equation $F(x,y) = 0$, where $y$ is the unknown and $x$ a parameter, depend holomorphically on $x$ was known to Cauchy in 1831, with the proof that I presented.

112 C. L. Siegel. Zu den Beweisen des Vorbereitungssatzes von Weierstrass. In *Number Theory and Analysis (Papers in Honor of Edmund Landau)*, pages 297–306. Plenum, New York, 1969.

– Weierstrass published his proof in 1886 but mentions in a footnote that he had been lecturing on this theorem since 1860. Not surprisingly, he avoids as much as possible the use of Cauchy residues, but not completely, and works with series. His proof is only partially algebraic.

– The theorem is proved by Poincaré in his thesis, in 1879, with no mention to Cauchy. As usual, the word "proof" has

to be taken with great care in Poincaré's writing, and this is especially true in this early paper. Much later, for instance in his *Méthodes Nouvelles*, he referred to his thesis, without providing a better proof and without mentioning Weierstrass. Interestingly, Henri Cartan, one of the founding fathers of Bourbaki, does not mention Poincaré in his paper.

– In 1905, Lasker[113] provided a fully algebraic proof and deduced algebraic consequences for the rings of formal and convergent series.

– Siegel also emphasizes that according to him the shortest proof is due to Stickelberger[114] in 1887.

For a modern and elementary presentation of the theorem, see Ebeling's book[115]. For a careful description of the many variants of the theorem and additional historical comments, see Grauert and Remmert[116].

I cannot end this chapter without mentioning that there is a version of this theorem for $C^\infty$ functions, conjectured by Thom and proved by Malgrange[117]. But, that's another story[118]...

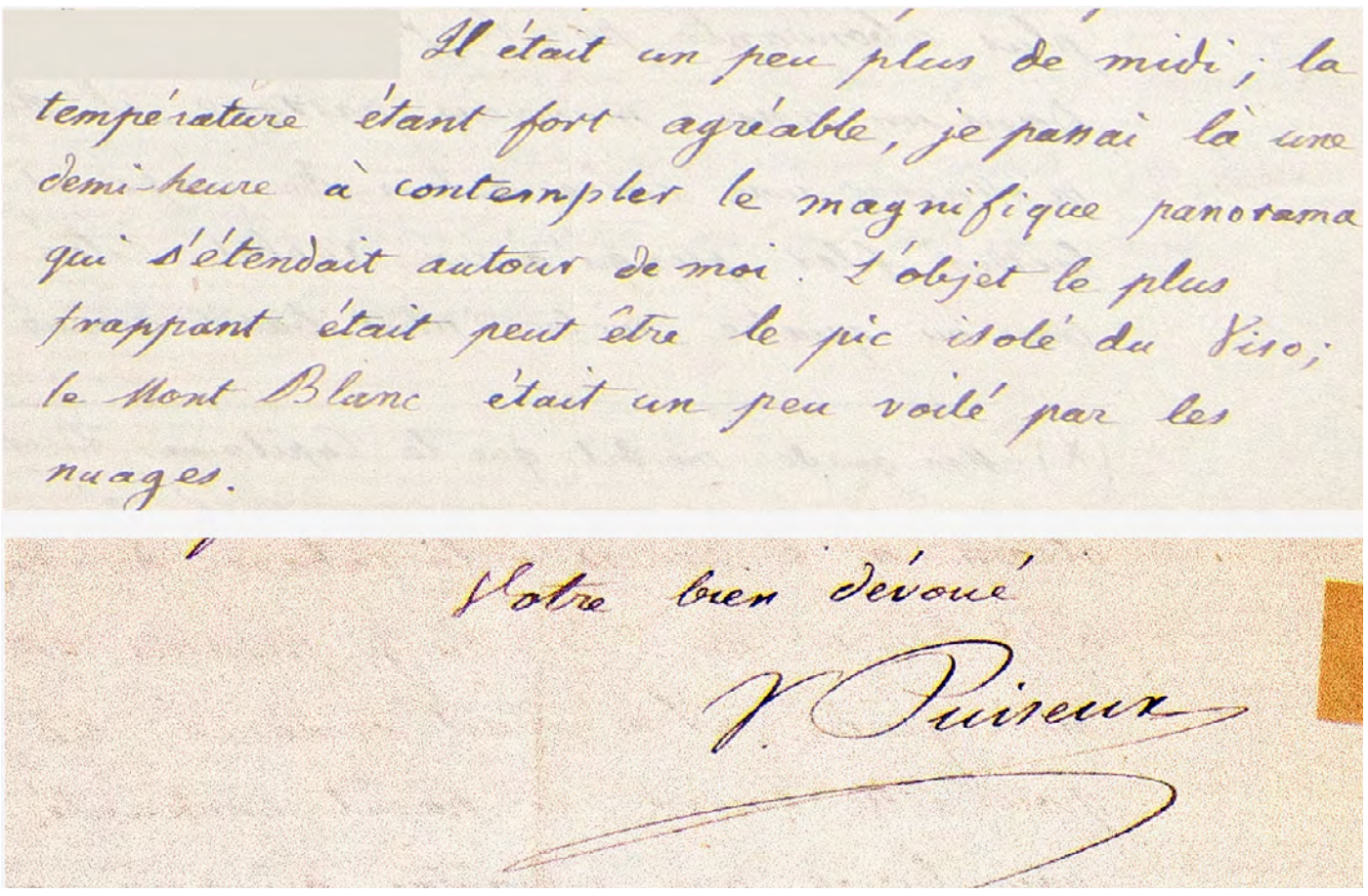
Il était un peu plus de midi; la
température étant fort agréable, je passai là une
demi-heure à contempler le magnifique panorama
qui s'étendait autour de moi. L'objet le plus
frappant était peut être le pic isolé du Viso;
le Mont Blanc était un peu voilé par les
nuages.

Votre bien dévoué

V. Puiseux

From a letter of Victor Puiseux describing his expedition to Mount Pelvoux: "Contemplating the magnificent panorama around me". Victor Puiseux could have been a model for Caspar David Friedrich when he painted the *Wanderer above the sea of fog*, but the Pelvoux expedition happened 30 years later! ©

[113] E. Lasker. Zur Theorie der Moduln und Ideale. *Math. Ann.*, 60(1):20–116, 1905.

[114] L. Stickelberger. Ueber einen Satz des Herrn Noether. *Math. Ann.*, 30(3):401–409, 1887.

[115] W. Ebeling. *Functions of several complex variables and their singularities*. GSM083. AMS, 2007.

[116] R. R. H. Grauert. *Analytische Stellenalgebren*. Die Grundlehren der mathematischen Wissenschaften 176. Springer-Verlag Berlin Heidelberg, 1971.

[117] B. Malgrange. *Ideals of differentiable functions*. Tata Institute of Fundamental Research Studies in Mathematics, No. 3. Tata Institute of Fundamental Research, Bombay; Oxford University Press, London, 1967.

[118] V. I. Arnold, S. M. Gusein-Zade, and A. N. Varchenko. *Singularities of differentiable maps. Classification of critical points, caustics and wave fronts*. Modern Birkhäuser Classics. Birkhäuser/Springer, New York, 2012. , Reprint of the 1985 edition.

"Il était un peu plus de midi, la température étant fort agéable, je passai là une demi-heure à contempler le magnifique panorama qui s'étendait autour de moi. L'objet le plus frappant était peut-être le pic isolé du Viso; le Mont Blanc était un peu voilé par les nuages."

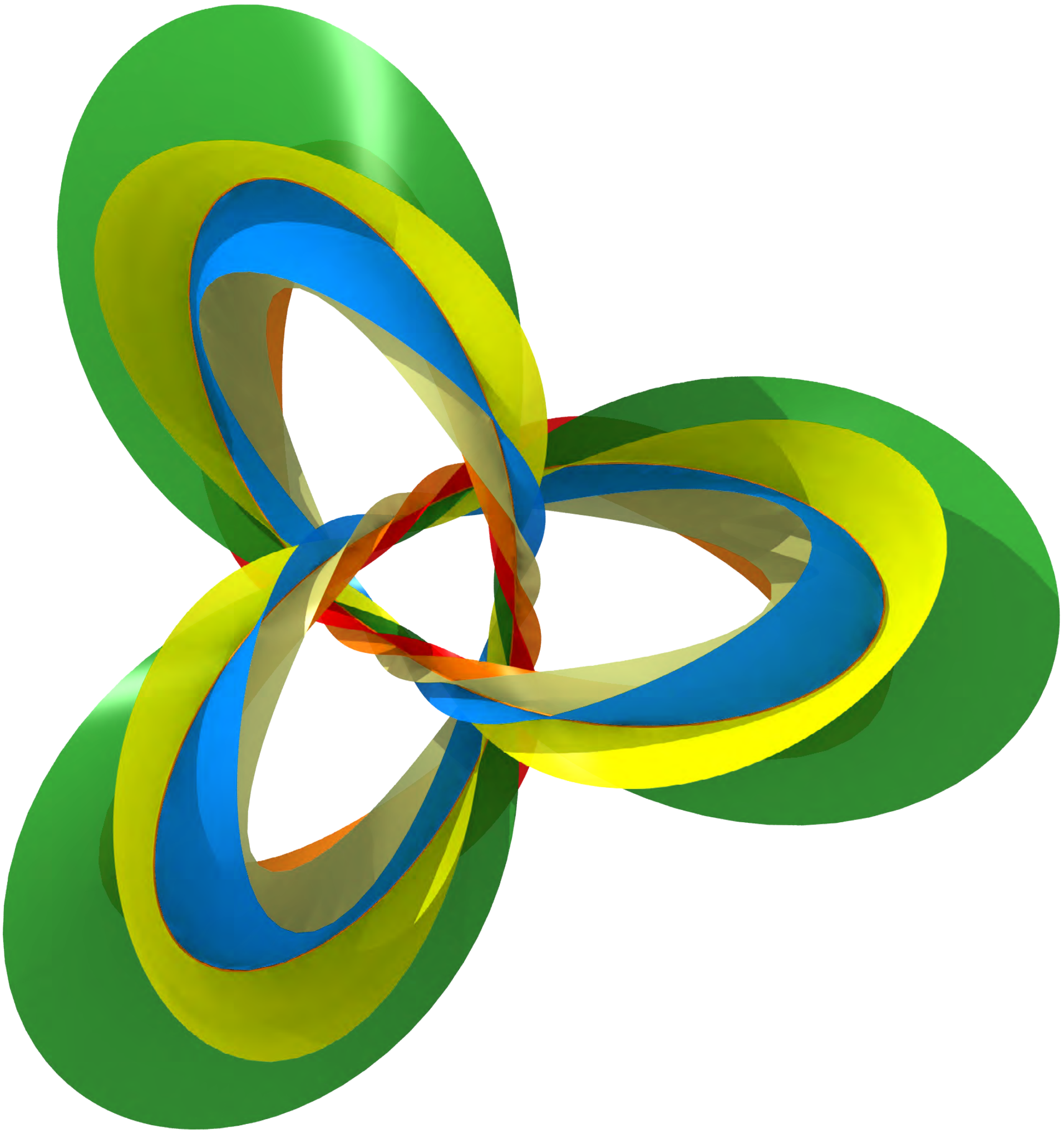

Milnor fibers of $x^3 - y^2$.

# *Jack Milnor and his fibration*

WHEN I ENTER A MATHEMATICAL LIBRARY, or when I navigate through the *Mathematical Reviews*, or simply when I *google*, I am frequently overwhelmed by the vastness of the mathematical world. Even topics that may look microscopic to the layman, like for instance the topology of algebraic curves, are actually huge territories whose exploration could easily require several lives. This feeling can be either depressing or intoxicating, depending on my mood ☹. In this petit livre, the best I can do is to describe one significant example, to mention some of the main results, and to refer to some of the (long) books proposing a complete discussion of the state of the art.

In any case, one single book should be emphasized as a gem and *has to be read by all students* interested by this topic: *Singular points of complex hypersurfaces*[119] by Milnor, a great master in the art of writing mathematics.

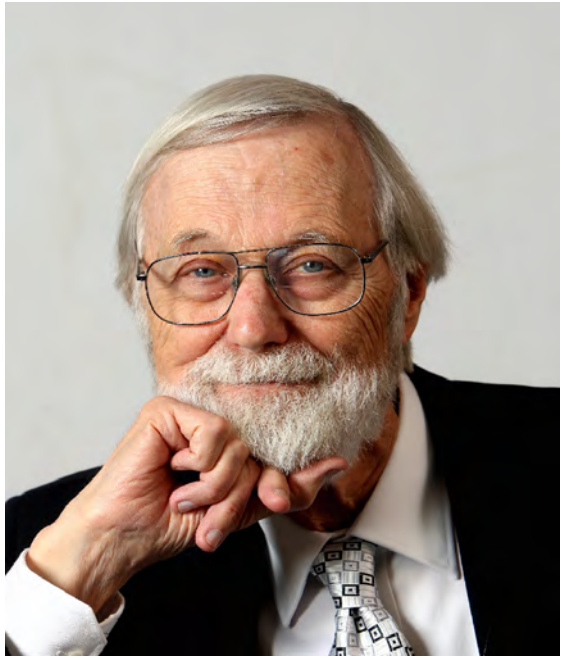

John Milnor. ©
French-speaking people are often surprised to learn that Jack is a nickname for John. The English translation of Jacques is James.

[119] J. Milnor. *Singular points of complex hypersurfaces*. Annals of Mathematics Studies, No. 61. Princeton University Press, Princeton, N.J.; University of Tokyo Press, Tokyo, 1968.

## *An example*

Look at the curve

$$F(x,y) = -x^{10} + x^9 + 6x^8y - 3x^6y^2 + 2x^5y^3 + 3x^3y^4 - y^6 = 0.$$

This $F$ has not been chosen at random. Any equation $F = 0$ can be solved using Puiseux series. In this example, I cheated and I started from the solution

$$y = x^{3/2} + x^{5/3}.$$

and I looked for the equation! We have

$$(y - x^{3/2})^3 = x^5.$$

Expanding and raising to a suitable power in order to eliminate rational exponents, we find indeed $F(x,y) = 0$. Actually, setting $x = x_1^2$ and $y = x_1^3(1 + y_1)$, as one should do using Newton's algorithm, the result factorizes, as it should:

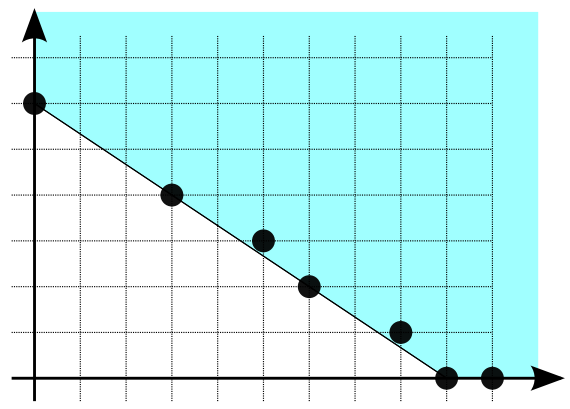

The Newton polygon of $F$.

$$F\left(x_1^2, x_1^3(1+y_1)\right) = -x_1^{18}\left(x_1 - y_1^3\right)\left(-8 + x_1 - 12y_1 - 6y_1^2 - y_1^3\right).$$

It follows that the zero locus of $F$ in the neighborhood of the origin contains exactly one branch $x_1 = y_1^3$ or

$$x = t^6 \quad ; \quad y = t^9 + t^{10}$$

so that $y = x^{3/2} + x^{5/3}$ as expected. For a general $F$, and even for a polynomial, we should expect an *infinite* Puiseux series but we will first look at this specific example.

Let us examine the *link of the singularity*, intersection of the curve $F = 0$ with a small sphere $\mathbb{S}^3_\epsilon$. The transversality of the intersection of the curve with small spheres is easy to see. Indeed the square of the norm

$$\phi : t \in \mathbb{C} \mapsto |t^6|^2 + |t^9 + t^{10}|^2 \in \mathbb{R}_+$$

is equivalent to $|t|^{12}$ for $t$ small, and the equation $\phi(t) = \epsilon^2$ defines a closed loop in $\mathbb{C}$, close to $|t| = \epsilon^{1/6}$, transversal to the radial lines. In other words, the intersection of $F = 0$ with a small sphere $\mathbb{S}^3_\epsilon$ is an embedded circle, i.e. a *knot*, which is the image by $\phi$ of this loop.

A *knot* is an embedding of the circle in the 3-sphere. A *link* is the disjoint union of finitely many knots. Two knots or links are considered equivalent if there is an orientation preserving homeomorphism of the sphere sending the first to the second.

Up to homeomorphisms of the sphere, this knot is independent of $\epsilon$. We could even use ellipsoids instead of spheres, or even our square sphere $\max(|x|, |y|) = \epsilon$. The detailed proof is technical and boring but the key idea is quite simple. Given two Euclidean norms $N_0, N_1$ in $\mathbb{R}^4$, we can intersect the curve $F(x,y) = 0$ with small spheres $N_0 = \epsilon$ and $N_1 = \epsilon$. This gives two knots in two manifolds homeomorphic to a sphere. We now construct a path of norms $\lambda N_1 + (1-\lambda)N_0$ (for $0 \leq \lambda \leq 1$) so that we actually have a continuous family of embedded circles in spheres

This non-proof is hiding a fundamental fact in differential topology. If we have a family of embeddings $i_\lambda : X \to Y$ (for $0 \leq \lambda \leq 1$) of some compact manifold $X$ in some other manifold $Y$, then there is an *isotopy*, that is to say a family of diffeomorphisms $\Phi_\lambda$, of $Y$ such that $i_\lambda = \Phi_\lambda \circ i_0$.

which therefore define "the same knot". A similar argument could be used for the "square sphere". See[120] for an illustration of possible mistakes that a naive beginner could make.

[120] F. Deloup. The fundamental group of the circle is trivial. *Amer. Math. Monthly*, 112(5):417–425, 2005.

Denote this knot by $K_F$. We will use the convenient square sphere $\max(|x|,|y|) = \epsilon$. The intersection with the curve is located in the solid torus $|x| = \epsilon$ and $|y| \leq \epsilon$, so that $|t| = \epsilon^{1/6}$. Let us rescale and set $X = x/\epsilon$ and $Y = y/\epsilon$. In particular, $X$ is in the unit circle and $Y$ in the unit disc. If $t = \epsilon^{1/6}\tau$ we get

$$X = \tau^6 \quad ; \quad Y = \epsilon^{1/2}\tau^9 + \epsilon^{2/3}\tau^{10}$$

where $\tau$ describes the unit circle.

*In all the following pictures, the solid torus $\mathbb{S}^1 \times D^2$ is drawn as the cylinder $[0, 2\pi[\times D^2$ and the two faces $\{0\} \times D^2$ and $\{2\pi\} \times D^2$ should be glued.*

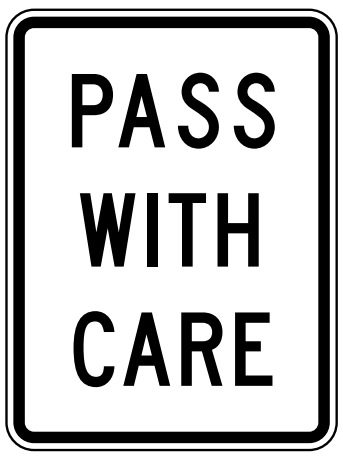

Take your time! Look at the pictures in this chapter with great care. This is not easy. ©

For every $X$ on the unit circle, there are exactly six values of $Y$, differing by multiplication of $\tau$ by some sixth root of unity. We say that the knot $K_F$ is in *braid form*: it intersects transversally all the discs $\{\star\} \times D^2$. Going around the circle, these six points are permuted in a way that will now be described.

©

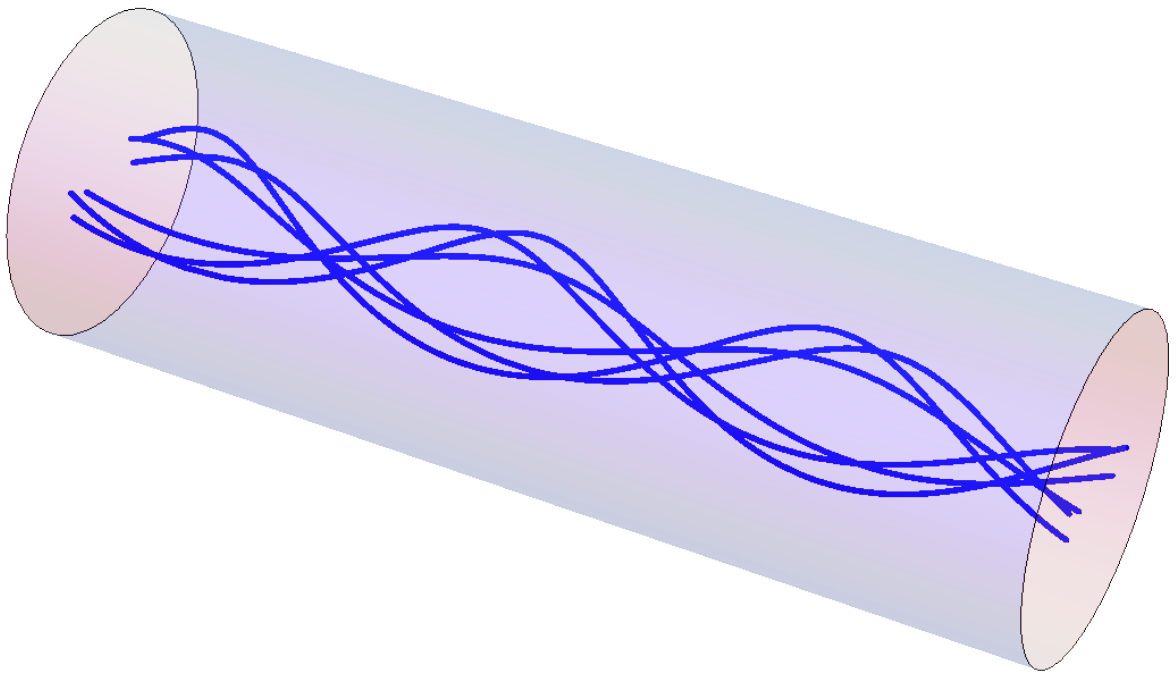

Observe that $\epsilon^{2/3}$ is small compared to $\epsilon^{1/2}$ for small $\epsilon$. Observe also that $\epsilon^{1/2}\tau^9$ takes only two values when one multiplies $\tau$ by a sixth root of unity. The knot associated to

$$X = \tau^6 \quad ; \quad Y_0 = \epsilon^{1/2}\tau^9$$

is simply our trefoil friend $x^3 = y^2$. When $X$ goes around the circle, the corresponding two points in $\{X\} \times D^2$ rotate by three half turns, producing the trefoil knot.

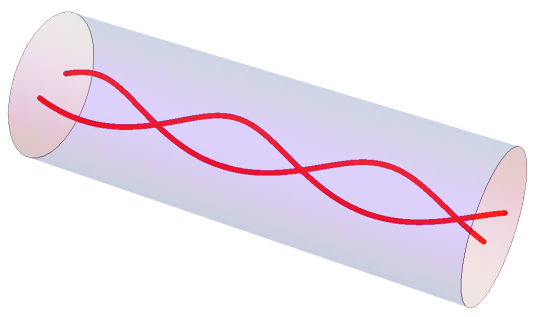

Consider $Y = \epsilon^{1/2}\tau^9 + \epsilon^{2/3}\tau^{10} = Y_0 + Y_1$ as a small perturbation of $Y_0$. Notice that

$$X = \tau^6 \quad ; \quad Y_1 = \epsilon^{2/3}\tau^{10}$$

gives three (very small) values of $Y_1$ for each value of $X$ on the unit circle. Hence the six points $Y$ on each disc $\{\star\} \times D^2$, come in two groups of three points. In other words, the knot $K_F$ lies in a thin tubular neighborhood of the trefoil knot and intersects small discs transversal to the trefoil in three points.

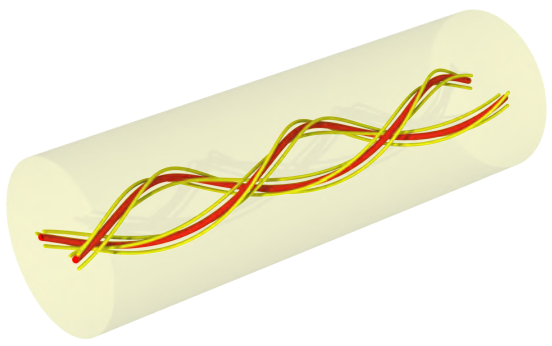

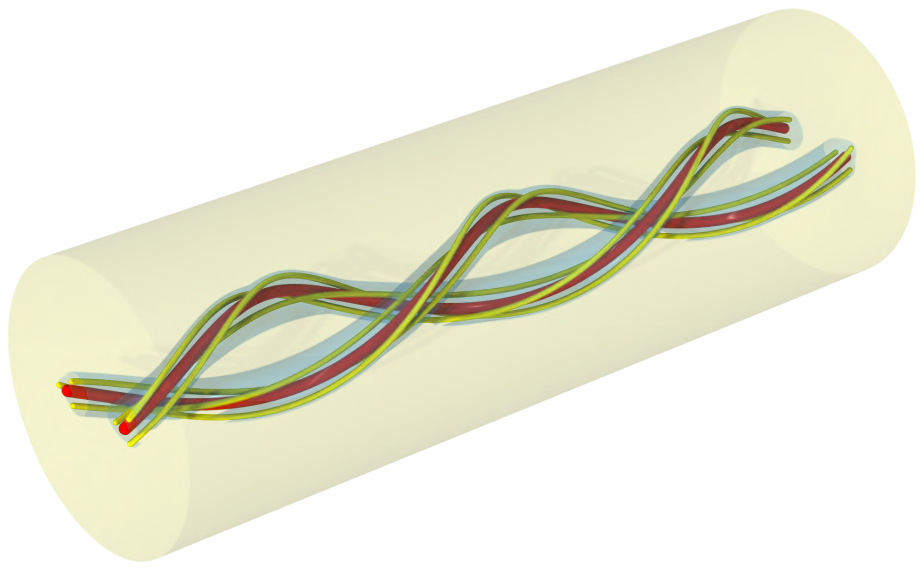

This tubular neighborhood of the trefoil can be parameterized in the following way.

$$(\mu, \zeta) \in \mathbb{S}^1 \times D^2 \mapsto \left(X = \mu^2, Y = \epsilon^{1/2}\mu^3 - 2\epsilon^{2/3}\mu^{-3}\zeta\right).$$

The coefficient 2 in front of $\epsilon^{2/3}\mu^{-3}$ is not important: its only role is to give enough thickness to the tube to contain our knot.

The circle $\mathbb{S}^1 \times \{0\}$, core of this solid torus, is mapped to the trefoil knot. One may ask why I chose $\epsilon^{1/2}\mu^3 - 2\epsilon^{2/3}\mu^{-3}\zeta$ and not simply $\epsilon^{1/2}\mu^3 + 2\epsilon^{2/3}\zeta$ which would also be a parametrization. The point is that, with this choice of coordinates, in this tubular neighborhood, $x^3 - y^2$ is equal to

$$\epsilon^3(\mu^2)^3 - \epsilon^2(\epsilon^{1/2}\mu^3 - 2\epsilon^{2/3}\mu^{-3}\zeta)^2$$

which is of the order of $4\epsilon^{19/6}\zeta$ for small $\zeta$, so that the argument of $x^3 - y^2$ is close to the argument of $\zeta$. Therefore *the Milnor fibers of $x^3 - y^2 = 0$ in the tubular neighborhood are close to the pages* $\arg\zeta =$ *constant.*

In these coordinates, we can relate the expressions $\epsilon^{1/2}\mu^3 - 2\epsilon^{2/3}\mu^{-3}$ and $\epsilon^{1/2}\tau^9 + \epsilon^{2/3}\tau^{10}$ for $Y$. Using $\mu = \tau^3$, our knot $K_F$ is the image of

$$\tau \in \mathbb{S}^1 \mapsto (\mu, \zeta) = (\tau^3, -\frac{1}{2}\tau^{19}) \in \mathbb{S}^1 \times D^2.$$

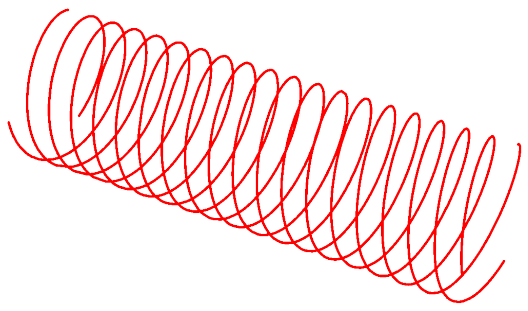

©

This is a $(19,3)$ torus knot.

This 19 might look strange. Note that there is a homeomorphism of the solid torus $\mathbb{S}^1 \times D^2$ which maps the $(19,3)$ torus knot to the $(19-3k,3)$-torus knot for any $k$ (for example to the $(1,3)$-knot, which is much simpler) but such a homeomorphism cannot be (homotopic to) the identity on the boundary and cannot be extended to the full sphere.

Indeed the map $(\tau, \zeta) \mapsto (\tau, \tau^k\zeta)$ is a homeomorphism which is twisting the solid torus.

It follows that $K_F$ is obtained by inserting a $(19,3)$ torus knot in a neighborhood of a $(3,2)$ torus knot. This is a typical example of an *iterated torus knot*. They are sometimes called *cable knots* since it resembles the construction of cables made of twisted strands that are braided together, or *satellites* turning around planets that rotate around the Sun.

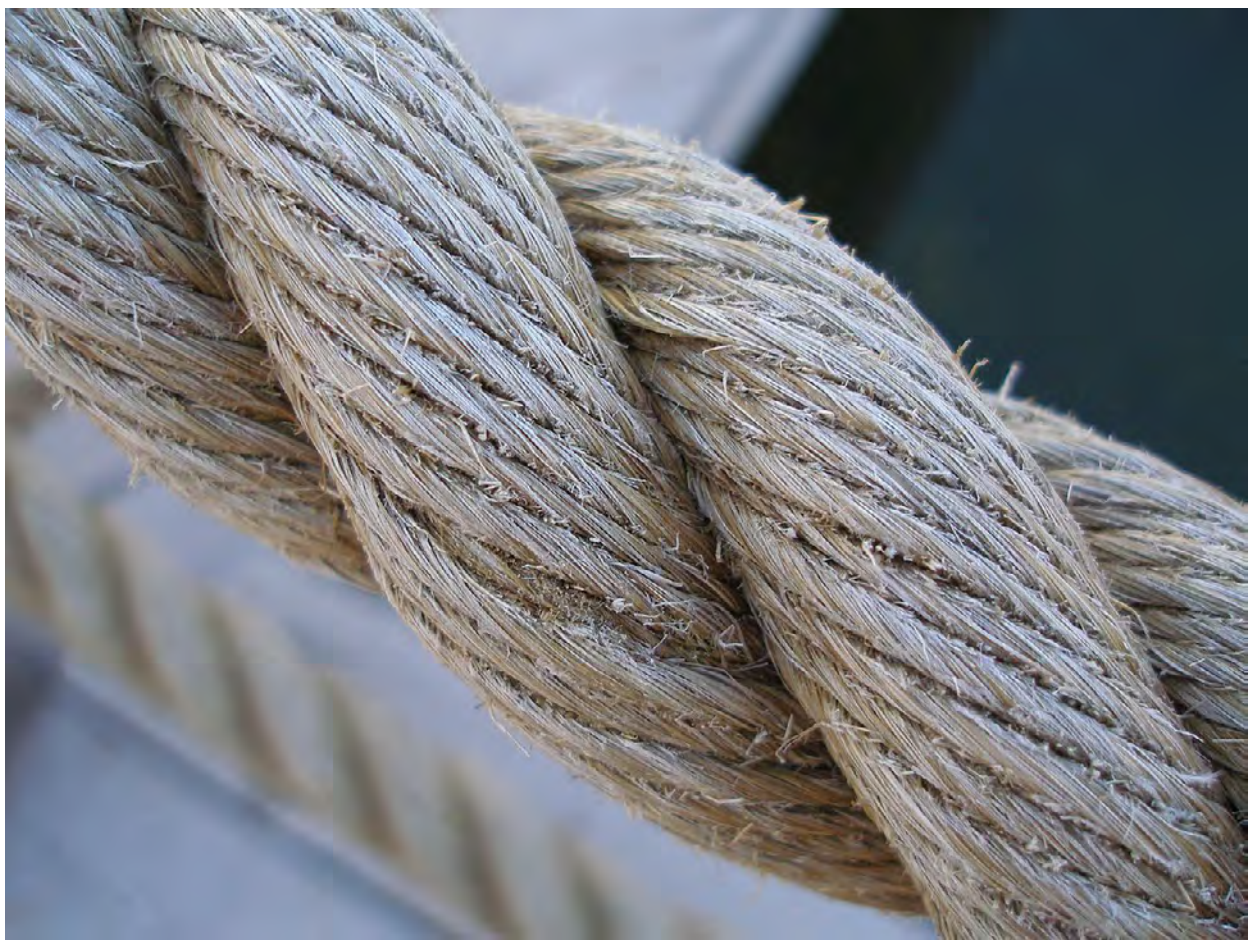

## *Milnor's fibration*

In the case of the trefoil, the level surfaces of $\arg(x^3 - y^2)$ fill the complement of the knot as the pages of a book whose bind-ing is the trefoil. We described the topology of those pages as punctured tori.

In his 1968 seminal book, Milnor showed that this is a general fact. He actually proved a theorem in all dimensions, but we limit ourselves to the (complex) dimension 2 case.

**Theorem.** *Let $F(x,y)$ be a nonzero reduced holomorphic function defined in the neighborhood of the origin of $\mathbb{C}^2$ such that $F(0,0)=0$. If $\varepsilon > 0$ is small enough, then*

*– the curve $F(x,y) = 0$ intersects transversally small spheres $\mathbb{S}^3_\epsilon$ along some link $L_\epsilon \subset \mathbb{S}^3_\epsilon$ whose topology is independent of $\epsilon$.*

*– the map*

$$(x,y) \in \mathbb{S}^3_\epsilon \setminus L_\epsilon \mapsto \arg(F(x,y)) = \frac{F(x,y)}{|F(x,y)|} \in \mathbb{S}^1$$

*is a locally trivial fibration. The closures of the fibers are compact surfaces whose boundaries all coincide with $L_\epsilon$. In a tubular neighborhood of the link, they look like an open book: the fibers are locally products of a radial segment in $D^2$ and a segment.*

A series $F(x,y)$ is *reduced* if it has no multiple factors in its decomposition in irreducible factors.

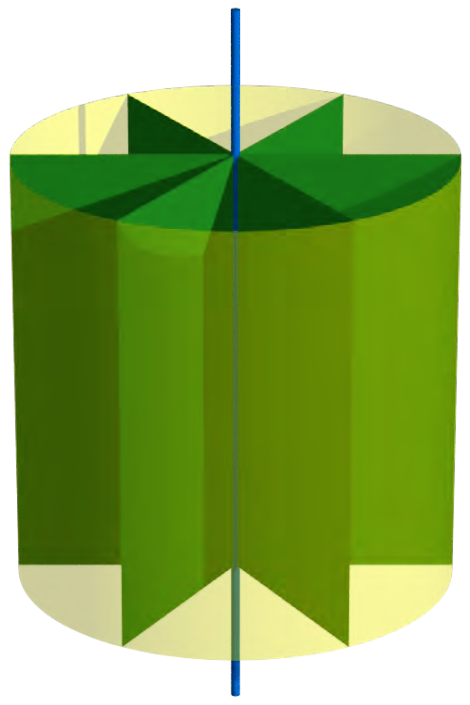

©

This theorem is the fundamental tool in the local study of singularities. However, I have to confess that for a long time, I had not looked at its proof and I was somehow convinced that it had to be elementary and straightforward. We have a very natural map to the circle; why shouldn't it be a fibration? I was wrong and the proof is indeed rather subtle. Amazingly, books dealing with this question are of two sorts. The first sort, arriving at the key point of the proof, in a very discreet way, just write "See Milnor, chapter 2". The second sort, arriving at the same key point, just copy almost word by word the content of "Milnor, chapter 2". Indeed, they are both right and it is difficult to do a better writing than "Milnor, chapter 2". My purpose here is not to innovate but to give some intuition on this theorem.

First observe that the theorem is true and elementary in dimension 1. Let $f(x)$ be some nonzero holomorphic function defined in the neighborhood of the origin of $\mathbb{C}$ and such that $f(0)=0$. Write $f(x) = ax^m + \cdots$ with $a \neq 0$ and look at the map

$$x \in \mathbb{S}^1_\epsilon \mapsto \arg(f(x)) = \frac{f(x)}{|f(x)|} \in \mathbb{S}^1.$$

This is obviously a covering map for $\epsilon$ small enough. Indeed, this is close in the $C^1$ topology to the covering map $x \mapsto \arg(a)\arg(x)^m$.

After this trivial case, let us consider a curve $F(x,y) = 0$. Assume that $F$ is reduced so that $\partial F/\partial y$ does not vanish on $F = 0$ (except at the origin). Instead of using the round sphere, we will use our square sphere $\max(|x|,|y|) = \epsilon$. For simplicity, assume that our curve only intersects the solid torus $T_1$ defined by $|x| = \epsilon$, $|y| < \epsilon$.

We want to show first that the argument of $F(x,y)$, restricted to this solid torus is a submersion (outside $F = 0$). In other words, given a point $(x,y)$ in $T_1$, we look for some tangent direction in $T_1$ along which the derivative of $\arg F$ is not zero. As a first attempt, we can try some vertical direction, fixing $x$. Then $F(x,y)$ changes according to the partial derivative $\partial F/\partial y$ so that $\arg F$ is indeed a submersion, *at least outside the zero locus* of $\partial F/\partial y$.

The *argument* of a nonzero complex number $z$ can be defined in several ways. It could be an element of $[0,2\pi[$, or of $\mathbb{R}/2\pi\mathbb{Z}$, or $z/|z|$, in the unit circle. In what follows, I always choose the most convenient way. I believe this will not create difficulties.

Now this zero locus $\partial F/\partial y(x,y) = 0$ is some other curve, only intersecting $F = 0$ at the origin. This new curve can be parameterized à la Puiseux, by $x = t^n$ and $y = f(t)$. By our trivial 1-dimensional case, the map

$$t \in \mathbb{S}^1_{\epsilon^{1/m}} \mapsto \arg F(t^m, f(t)) \in \mathbb{S}^1$$

is a covering map. Therefore, for points where $\partial F/\partial y(x,y) = 0$, we found some other direction in which the derivative of the argument does not vanish.

This argument is definitely not a complete proof of Milnor's theorem for several reasons.

The first is that we used a square sphere instead of a round one: this is not so serious and the argument could easily be adapted to the round sphere.

The second is that a submersion need not be a fibration, without some compactness assumption on the fibers. We have to study the local structure of our submersion close to the link $L_\epsilon \subset \mathbb{S}^3_\epsilon$. This is not difficult. The key point is that if $F(x,y) = 0$ and if we take some complex line in $\mathbb{C}^2$ passing through $(x,y)$ and transversal to the curve, then we can apply the trivial 1-dimensional case to analyze the argument $\arg F$ in the neighborhood of $(x,y)$ to get the local picture around the link.

Note that our simple presentation is limited to the dimension 2 case, and that Milnor's theorem holds in any dimension.

For an excellent presentation, I refer to Milnor, chapter 2 ☺.

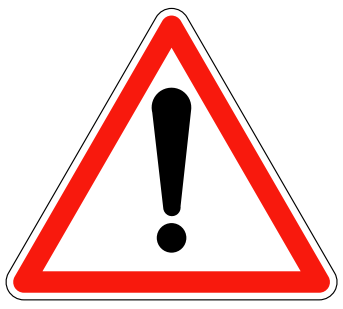

The end of this chapter requires a great attention, even though it is not necessary for the rest of the book.

## *Milnor's fibers in our example*

Let us come back to our example

$$F(x,y) = x^9 - x^{10} + 6x^8y - 3x^6y^2 + 2x^5y^3 + 3x^3y^4 - y^6 = 0.$$

The curve $F(x,y) = 0$ intersects a small sphere along a knot which is a satellite of the trefoil knot. We wish to describe the topology of the Milnor fibers $\arg F(x,y) = const$. If I ask my computer to draw one of these fibers, the resulting picture is the following.

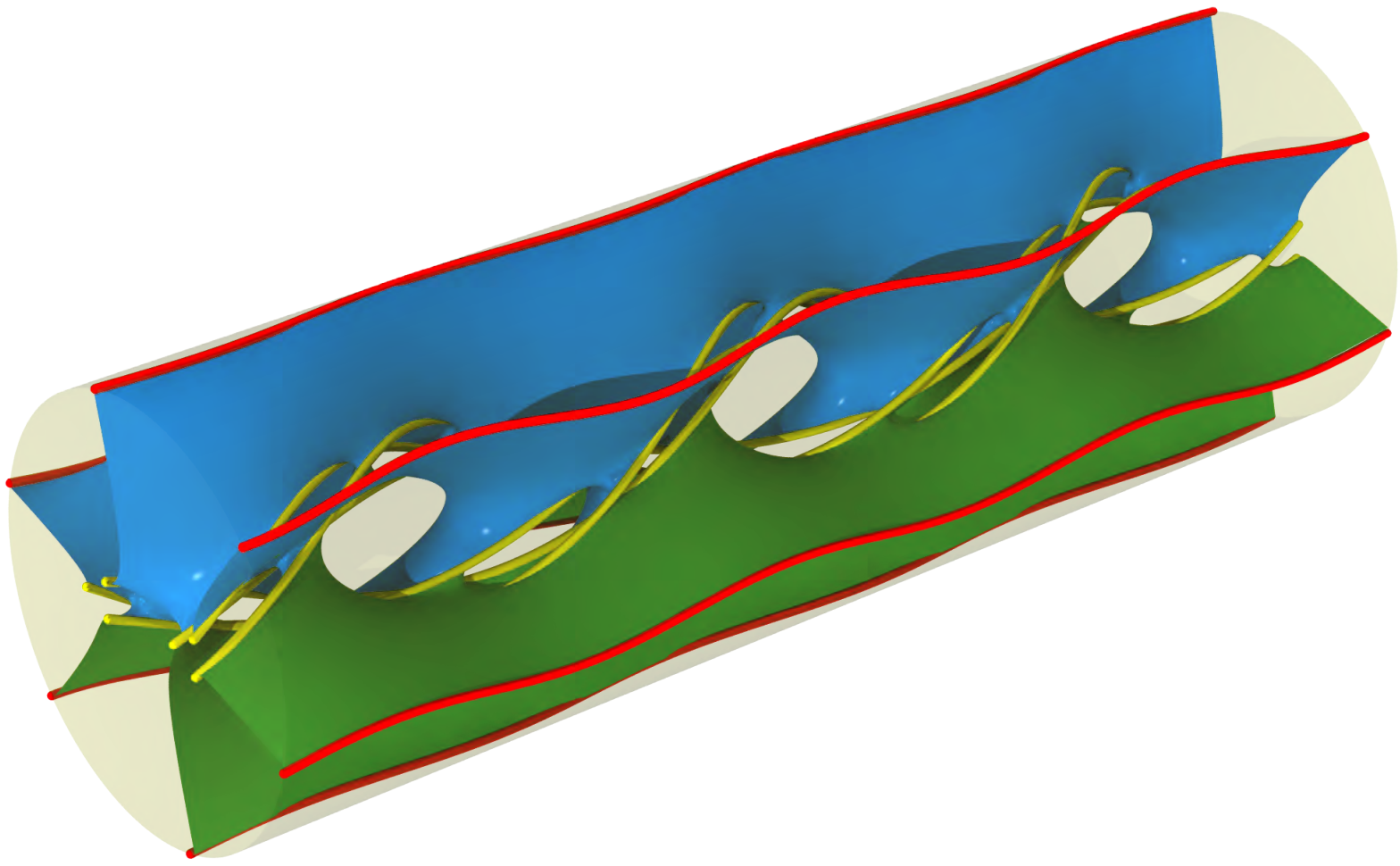

This is complicated and requires a careful analysis. Notice at least that this surface intersects the boundary of the solid torus along 6 curves (in red). The knot, represented in yellow, is also a boundary component. The colors blue and green of the two faces show that this surface is orientable. Keep in mind that there are 7 boundary components in the solid torus.

In order to understand this picture, let us look first at *the Milnor fibers of a* $p,q$ *curve* $x^p - y^q = 0$, where $p$ and $q$ are relatively prime (with $p > q$). We know that they are surfaces of genus $(p-1)(q-1)/2$ with a disc removed. Let us look at their position relative to our square sphere $\max(|x|,|y|) = \epsilon$. Since $p > q$, the intersection with $x^p - y^q = 0$ lies in the solid torus $T_1$ defined by $|x| = \epsilon$. On the boundary torus $|x| = \epsilon$ and $|y| = \epsilon$, the value of $x^p - y^q$ is very close to $y^q$ and the argument of $x^p - y^q$ is almost equal to $q$ times the argument of $y$. It follows that a Milnor fiber of $x^p - y^q$ intersects the boundary of $T_1$ along $q$ curves which are very close to $q$ parallels.

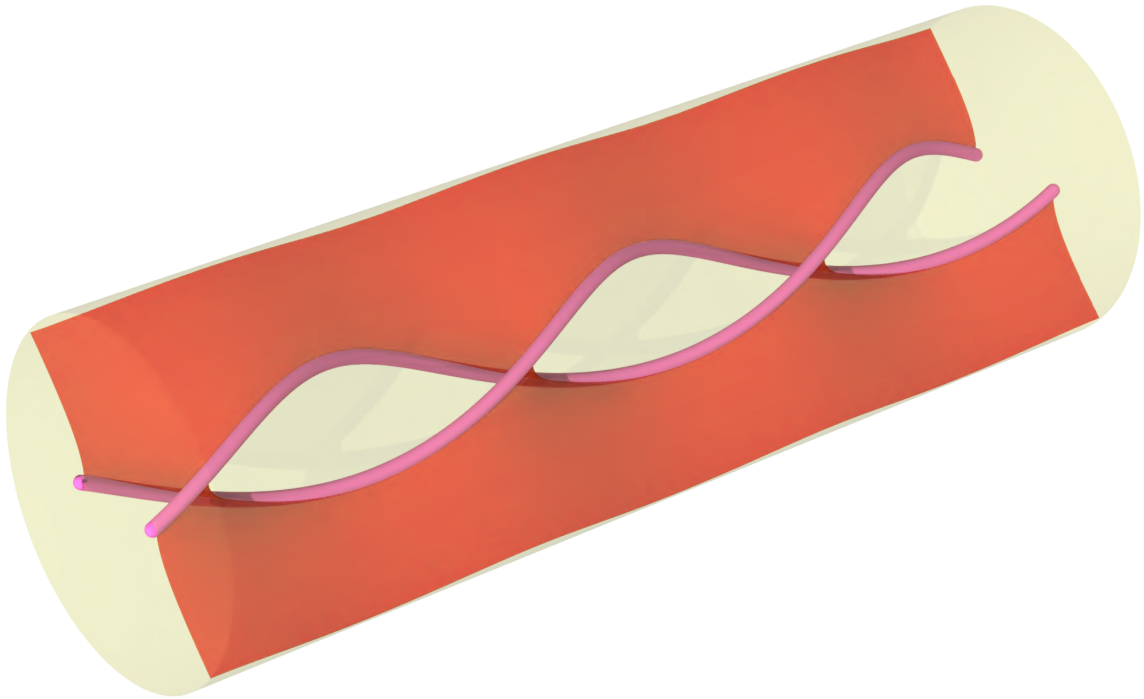

A Milnor fiber of trefoil knot, seen as a (3,2)-torus knot. The boundary of the surface is the knot and its intersection with the boundary of the cylinder consists of two parallels. ©

For the same reason, on the other torus $T_2$ defined by $|y| = \epsilon$, a Milnor fiber almost coincides with $q$ discs where the argument of $y$ takes $q$ values and $x$ describes the disc of radius $\epsilon$. In other words, the intersection of some Milnor fiber of $x^p - y^q$ with $T_1$ is a surface of genus $(p-1)(q-1)/2$ where $1+q$ discs have been removed. The boundary of the first removed disc is the torus knot, sitting inside $T_1$ and the $q$ other discs have boundaries $q$ circles on the boundary of $T_1$.

Let us now come back to our *more complicated example defined by* $F(x,y) = 0$.

Extracting the dominant terms in Newton's polygon, we get

$$F(x,y) = (x^3 - y^2)^3 - x^{10} + 6x^8y + 2x^5y^3.$$

Recall that we constructed a tubular neighborhood $\mathcal{T}_\epsilon$ of the trefoil knot, parameterized by $(\mu,\zeta) \in \mathbb{S}^1 \times D^2$ and in which $x^3 - y^2$ is of the order of $4\epsilon^{19/6}\zeta$. On the boundary of this solid

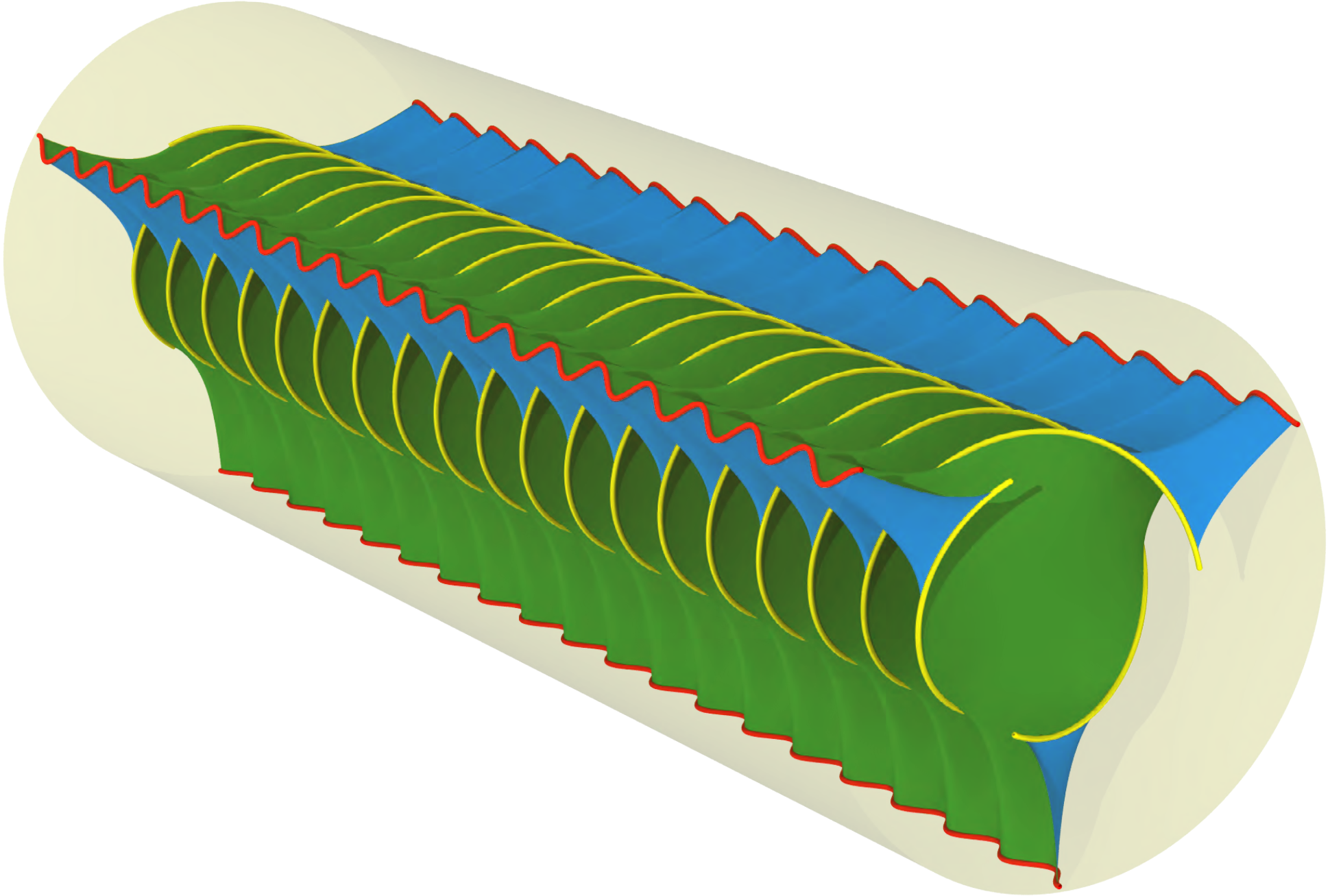

A Milnor fiber of a (19, 3)-torus knot. The boundary of the surface is the knot (in yellow) and its intersection with the boundary of the cylinder consists of 3 parallels (in red).

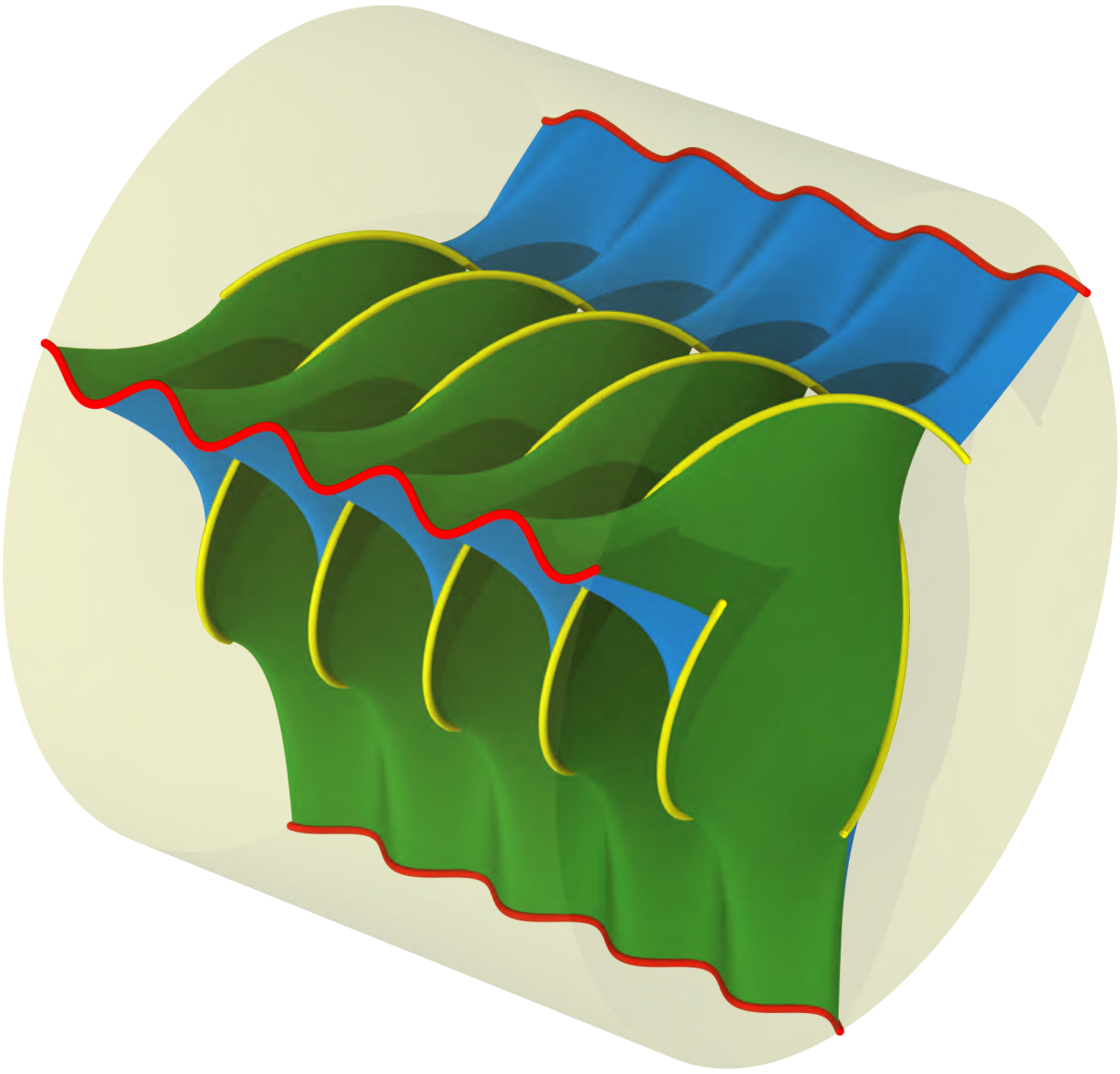

Here is a small slice to understand better the previous picture. The blue and green colors show that the surface is indeed orientable.

torus, where $|\zeta| = 1$, we have $|(x^3 - y^2)^3| \simeq 64\epsilon^{57/6}$ and $|x|, |y| \leq \epsilon$, so that $F(x,y)/(x^3 - y^2)^3$ is very close to 1, and the argument of $F(x,y)$ is close to $\arg(\zeta)^3$. In particular each Milnor fiber of $F$ *on the boundary* of $\mathcal{T}_\epsilon$ is very close to *three* parallels with $\arg\zeta = constant$.

This also holds outside $\mathcal{T}_\epsilon$: each Milnor fiber of $F$ *outside* of $\mathcal{T}_\epsilon$ is very close to *three* Milnor fibers of $x^3 - y^2 = 0$. Do not forget that the Milnor fibers of $x^3 - y^2 = 0$ are punctured tori. Their intersections with $T_2$ consist of two discs. Their intersection with $T_1$ have three boundary components, two being on the boundary of $T_1$ and the third being the knot itself. This is indeed what can be seen if we hide what is inside $\mathcal{T}_\epsilon$.

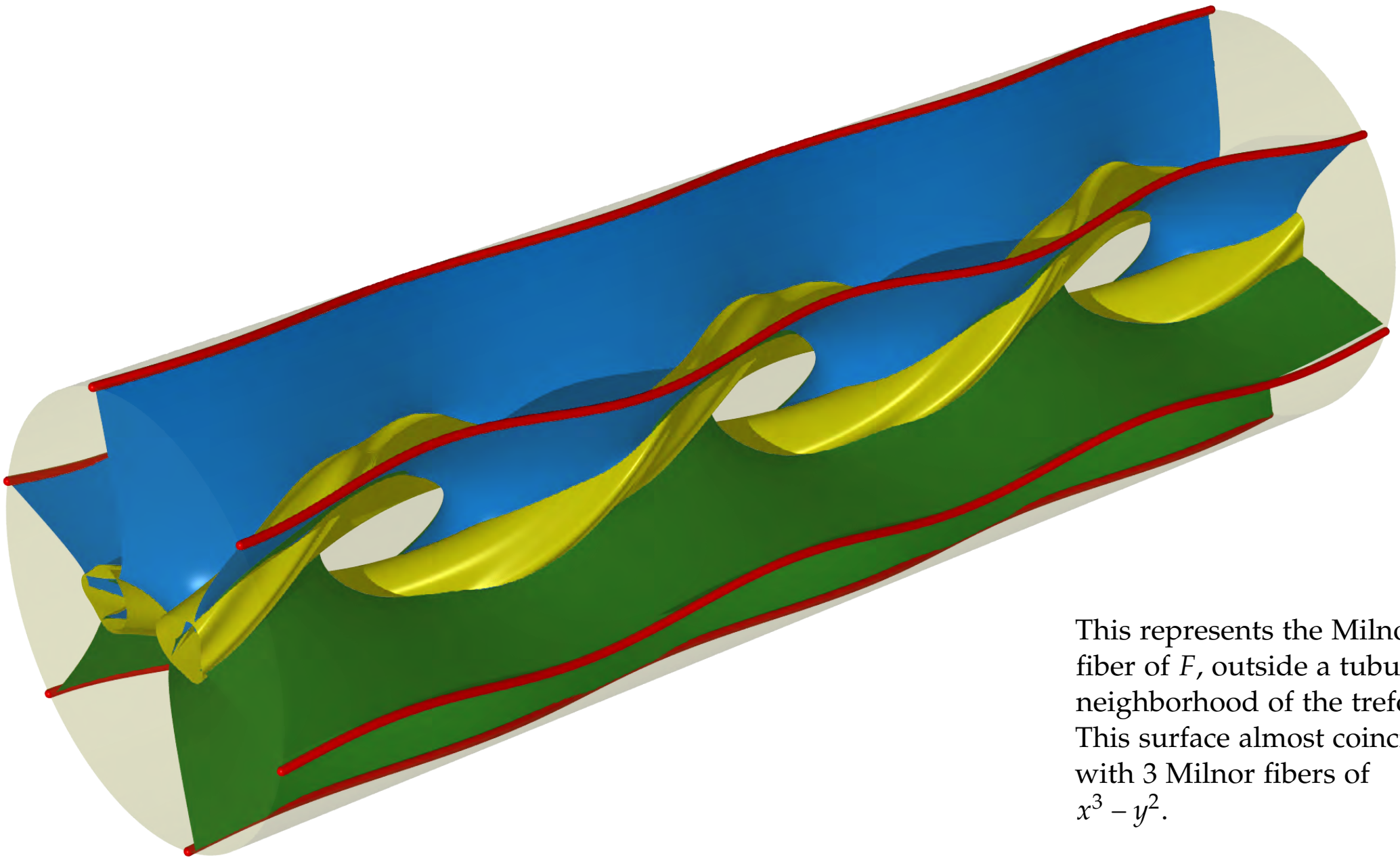

This represents the Milnor fiber of $F$, outside a tubular neighborhood of the trefoil. This surface almost coincides with 3 Milnor fibers of $x^3 - y^2$.

The interior of $\mathcal{T}_\epsilon$ is the realm of the $(19,3)$ knot that is inserted in the tube. Let us evaluate $F(x,y) = F(\epsilon X, \epsilon Y)$ inside the tube, in the coordinates $(\mu, \zeta)$:

$$F\left(\epsilon\mu^2, \epsilon^{3/2}\mu^3 - 2\epsilon^{5/3}\mu^{-3}\zeta\right).$$

This vanishes exactly on $\zeta = -\frac{1}{2}\mu^{19}$, to the first order. There-fore a Milnor fiber of $F$, inside $\mathcal{T}_\epsilon$, is close to a Milnor fiber of $y^3 = x^{19}$. It is a surface of genus $(3-1)(19-1)/2 = 18$ with 4 discs removed. One of the boundaries is the boundary of our knot, as it should be, and the three others are three parallels on the boundary of $\mathcal{T}_\epsilon$.

*In summary, a Milnor fiber of* $F(x,y) = 0$ *is homeomorphic to a closed orientable surface of genus* 18 *on which one performs three connected sums with a torus, and for which one finally deletes a disc.* It is a surface of genus 21. Quite complicated.

The rest of this chapter will be very vague: I cannot give more than a glimpse of the theory.

Splitting the complement of $K_F$ along one Milnor fiber $\Sigma$, we get a product $\Sigma \times [0,1]$. In order to reconstruct the complement of the knot, one should glue $\Sigma \times \{0\}$ to $\Sigma \times \{1\}$ using some diffeomorphism of $\Sigma$. This diffeomorphism, well-defined up to isotopy, is called the *monodromy* of the knot. The action of the monodromy on the first homology has a characteristic polynomial which is the *Alexander polynomial* of the knot. In our example, everything can be described in a rather concrete way.

Indeed, Milnor's fibration in the complement of *one* fiber is a fibration onto $[0,1]$, hence a trivial fibration since $[0,1]$ is contractible.

Let me just give the result. Our surface $\Sigma$ contains three closed curves $\gamma_i$ ($i = 1,2,3$) along which we performed the connected sum. Cutting along the $\gamma_i$'s, we get four components: $S, \Sigma_i$, where $S$ is a surface of genus 18 minus 4 discs and each $\Sigma_i$ is a punctured torus. Some monodromy map $\psi$ preserves the curves $\gamma_i$ and is a *Dehn twist* in some annulus around these curves. This means that $\psi$ in the neighborhood of these curves looks like the picture in the margin. This shows that the action of the monodromy on homology is periodic, but this would not be true in homotopy. The curves $\gamma_i$ are homologous to zero but not homotopic to zero.

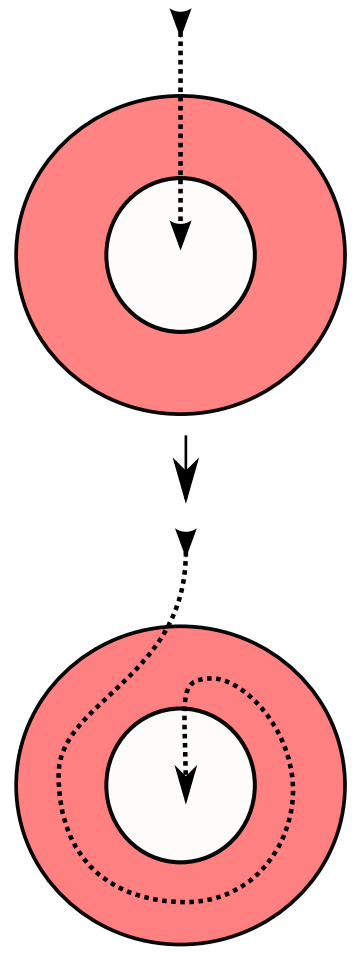

A Dehn twist. This homeomorphism is the identity on the boundary of the annulus, preserves the concentric circles, and twists them as shown.

If we cut open $\Sigma$ along the $\gamma_i$, we find the monodromies of $x^3 - y^2$, three times, and of $x^{19} - y^3$ once. It follows that the Alexander polynomial is the product of the cube of the polynomial for $x^3 - y^2$ and of the polynomial for $x^{19} - y^3$. Therefore, we get:

$$\frac{(X^6-1)^3(X-1)^3(X^{57}-1)(X-1)}{(X^2-1)^3(X^3-1)^3(X^{19}-1)(X^3-1)}$$

which is equal to

$$(1-X+X^2)^3(1-X+X^3-X^4+X^6-X^7+X^9-X^{10}+X^{12}$$
$$-X^{13}+X^{15}-X^{16}+X^{18}-X^{20}+X^{21}-X^{23}+X^{24}-X^{26}$$
$$+X^{27}-X^{29}+X^{30}-X^{32}+X^{33}-X^{35}+X^{36}).$$

## *The general case*

Let me only mention the most salient results.

The knots associated to a branch of a curve are always iterated torus knots.

The knots associated to two *irreducible* curves $F_1(x,y)=0$ and $F_2(x,y)=0$ are topologically equivalent through a homeomorphism of the 3-sphere if and only if the two associated branches have the same Puiseux characteristic invariant. Actually, one distinguishes these knots using the Alexander polynomials. This was proved a long time ago for the case of knots and even for the case of curves with two branches. The analogous fact, for non-irreducible curves, producing *links* consisting of *several* disjoint knots, was established much more recently.

The *monodromy* associated to general curves has been beautifully described by A'Campo[121].

I would say that the situation is now very well understood.

It is wise to stop our excursion here if we want to continue our promenade: there are other sites to visit. However, I would perfectly understand some frustration from the reader obliged to turn back on a path which seems to be (and which is) beautiful.

For much much much more on the topic, with a historical perspective, the reader should take a look at the wonderful survey by Weber[122] and at the already mentioned[123].

121 N. A'Campo. Sur la monodromie des singularités isolées d'hypersurfaces complexes. *Invent. Math.*, 20:147–169, 1973.

122 C. Weber. On the topology of singularities. In *Singularities II*, volume 475 of *Contemp. Math.*, pages 217–251. Amer. Math. Soc., Providence, RI, 2008.

123 E. Brieskorn and H. Knörrer. *Plane algebraic curves*. Modern Birkhäuser Classics. Birkhäuser/Springer Basel AG, Basel, 1986.

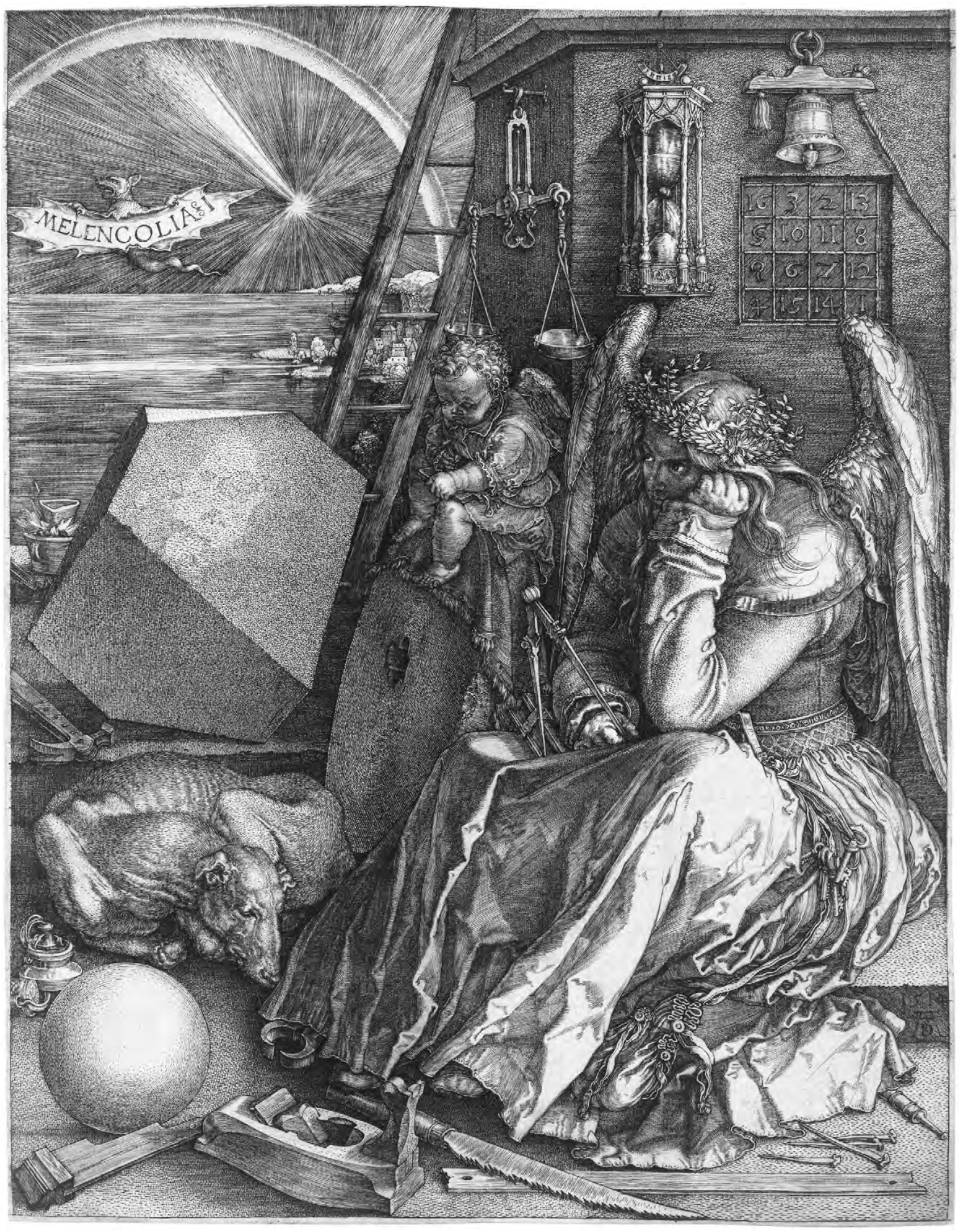

The famous engraving *Melencolia* by Dürer (1514). The polytope is not $K_5$ ! I recommend Günter Ziegler's article in *The Guardian* Dürer's polyhedron: five theories that explain Melencolia's crazy cube.

# *The Hipparchus-Schroeder-Tamari-Stasheff associahedron*

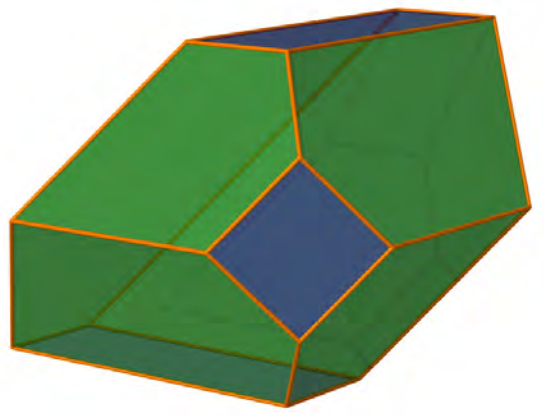

A model of $K_5$. ©

WE WILL FORGET ANALYTIC CURVES for a while and come back to trees, words, and combinatorics. There is a natural dictionary between the following three kinds of objects.

- *Binary rooted planar trees* with $n$ leaves.
- *Binary bracketing*s on a word of length $n$.
- *Partitions of a convex polygon* with $(n+1)$ edges (one of them being called *the root*) into $n$ triangles.

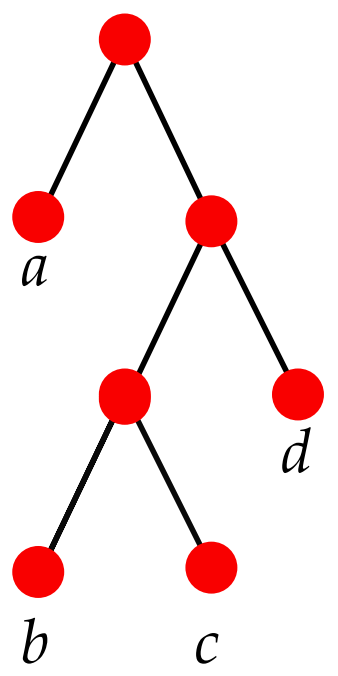

$a((bc)d)$

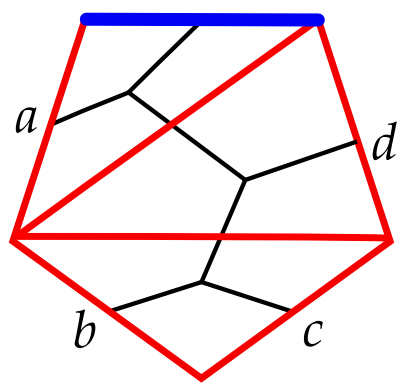

This is illustrated by the pictures in the margin.

We have also been concerned with planar rooted trees with $n$ leaves and such that every internal node has *at least* two children. These trees are associated with Schroeder bracketings on a word of length $n$ which are not necessarily binary. In terms of diagonals on a $(n+1)$-sided convex polygon, they correspond to collections of $k$ non-intersecting diagonals, with $0 \le k \le n$. The number of these objects is the $n$-th (small) Hipparchus-Schroeder number.

## *An abstract polytope*

We are going to construct a sequence of polytopes $K_n$, of dimension $n-2$, called the *associahedra*.

Draw an interval and label it with the unique rooted tree with one root having 3 children and no other node. Label its endpoints with the two rooted planar binary trees with 3 leaves. We get the following figure. This is $K_3$: just an interval.

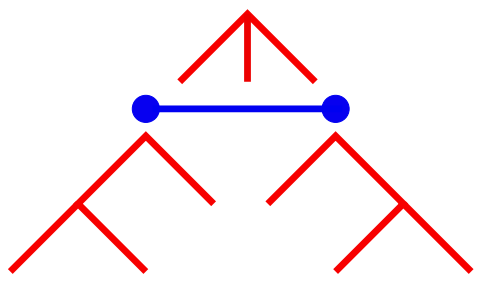

It is very tempting to connect two binary trees by an edge if one goes from one to the other by a local transition, as suggested in the previous picture. If you detect in some binary tree some sub-tree with 3 leaves, you delete it and you replace it by the other tree with 3 leaves: you defined an edge in the associahedron.

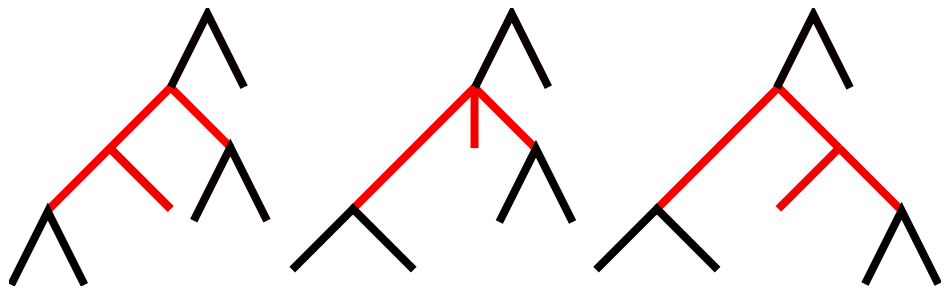

Let us draw a picture for $n = 4$. There are 5 binary trees with 4 leaves. We place them at the vertices of a pentagon. The 5 edges are labeled by the 5 planar trees with 4 leaves and exactly one 3-children node. There is one more planar tree with one root with 4 children. Place it in the center of the pentagon, as a label for the 2-dimensional face of the pentagon. This is $K_4$.

This suggests that one could define some *polytope* of dimension $n-2$ whose vertices are labeled by binary trees with $n$ leaves, whose edges are labeled by trees with a single 3-children node, etc. and whose unique top-dimensional face (of dimension $n-2$) is labeled by the single tree with one root and $n$ children (usually called a *corolla*).

Going to $n = 5$, we can still draw a picture.

It turns out that it is indeed possible to construct such a polytope for all values of $n$. The first problem is to give a precise definition of the word *polytope* in a *combinatorial context*. One would like a definition inspired by our geometrical intuition of a polytope *in Euclidean space*, but which should not take into account an embedding in some space.

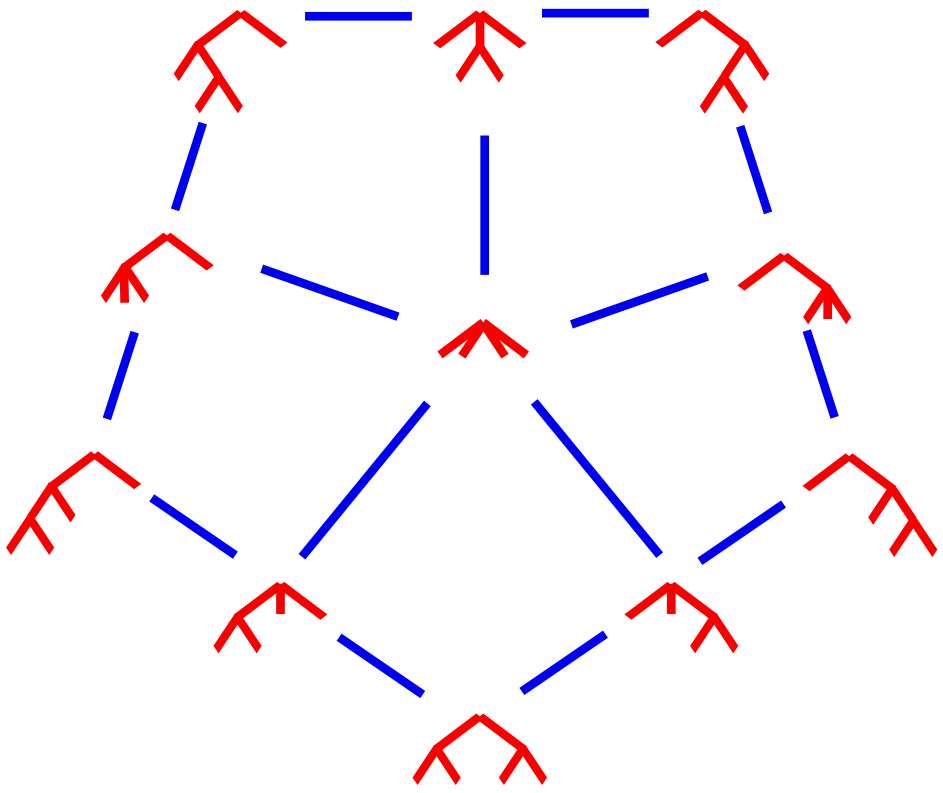

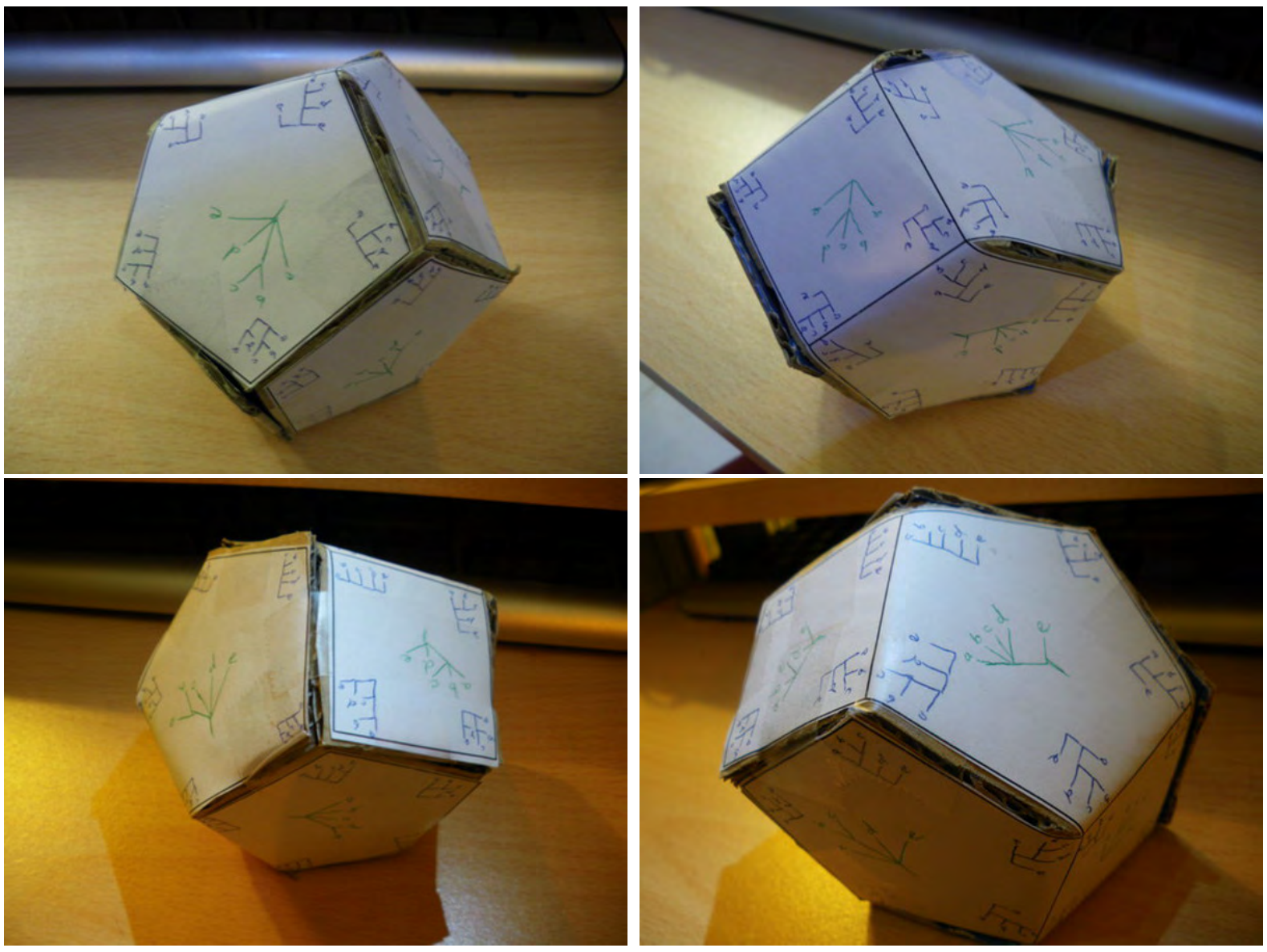

There is a well defined concept of *combinatorial polyhedron*, whose faces are segments, triangles, and simplices in general. Start with a set $V$ of points called *vertices* and select some subsets of $V$ which are called *faces* with only one condition: a subset of a face should be a face. If a face contains $k+1$ elements, one says that it is a $k$-dimensional simplex. This is a fairly easy definition but this is not suitable in our situation. For instance the 3-dimensional polytope $K_5$ above has 2-dimensional faces which are squares or pentagons, not triangles.

There are indeed several combinatorial (non-equivalent) definitions of *abstract polytopes* but we will not use them since our polytope will eventually be realized as a geometric object in Euclidean space. Nevertheless, an *abstract polytope* should at least be made of *faces* having some dimension, and there should be some partial ordering between faces, corresponding to the intuitive idea of adjacency. So, we will content ourselves with the definition of a partially ordered set $K_n$ of height $n-2$. This is very easy.

Let me define the height of a partially ordered set is the cardinality of a maximal totally ordered subset minus 1.

For simplicity, let us choose some convex polygon in the plane $\Pi_{n+1}$ with $(n+1)$ vertices, but the following construction is independent of the choice of this polygon. Let us also choose one side of the polygon, called the *root*.

A face of dimension $d$ of $K_n$ is by definition a set $F$ of $n-2-d$ non-intersecting diagonals in $\Pi_{n+1}$. The adjacency relation is defined using the *reverse* inclusion: say that a face associated to a subset $F_1$ is a *sub-face* of $F_2$ if $F_2 \subset F_1$. For instance, vertices of $K_n$, of dimension 0, correspond to partitions of $\Pi_{n+1}$ in $(n-1)$ triangles by $n-2$ diagonals. Using the root, these vertices are associated to planar binary trees, as desired.

A face of codimension $q$ of $K_n$ is associated with a rooted planar tree with $n$ leaves having exactly $q$ internal nodes (different from the root and the leaves), or, equivalently, having $q$ internal edges. Seen from the "tree point of view", one could say that the face associated to some tree $T_1$ is a subface of the one associated to $T_2$ if one obtains $T_2$ from $T_1$ by collapsing some edges.

For the time being, I only defined some partially ordered set. It would not be difficult to check that this does satisfy the axioms that define abstract polytopes... that I chose not to make explicit.

This is the *Hipparchus-Schroeder-Tamari-Stasheff associahedron*.

## *Some history*

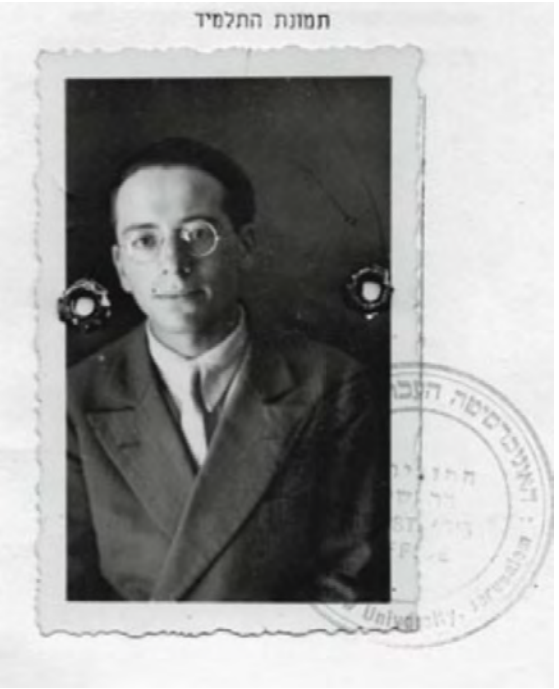

Dov Tamari (1911-2006).

As usual, giving a single name to a mathematical object is almost impossible

As we know, Catalan counted the number of vertices of $K_n$ and Hipparchus and Schroeder counted their faces.

Dov Tamari (formerly Bernhard Teitler) defined the combinatorial object in 1951 in his dissertation.

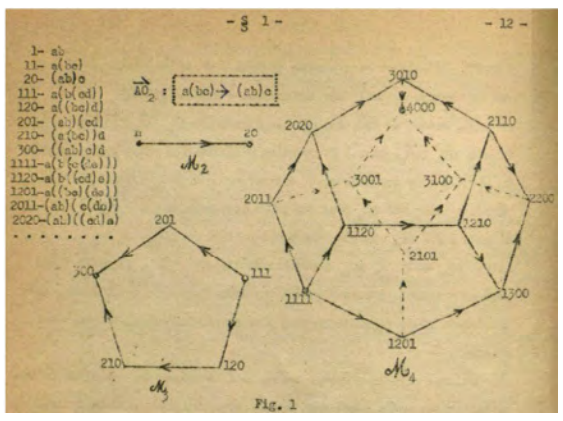

A figure from Tamari's dissertation.

I suggest reading the first chapter of the Tamari Memorial Festschrift[124] for a description of his motivation and biography (a "promenade" in Germany, Palestine, France, Israel, the USA, Brazil and the Netherlands, and across the twentieth century). One learns for instance that

> At least after 1948 Tamari opposed the injustices the Israelis did to the Palestinians, as well as discrimination directed against Jewish immigrants from Middle-Eastern countries, and these views were not at all widely accepted in those days.

[124] F. Müller-Hoissen, J. M. Pallo, and J. Stasheff, editors. *Associahedra, Tamari lattices and related structures. Tamari memorial Festschrift*, volume 299 of *Progress in Mathematical Physics*. Birkhäuser/Springer, Basel, 2012.

In 1963, J. Stasheff defined the same object, also in his dissertation, but in a very different topological context that will be discussed in some detail in the next chapter. He was not aware of the previous work of Tamari. The picture in the margin shows the *curved polytope* from his original paper.

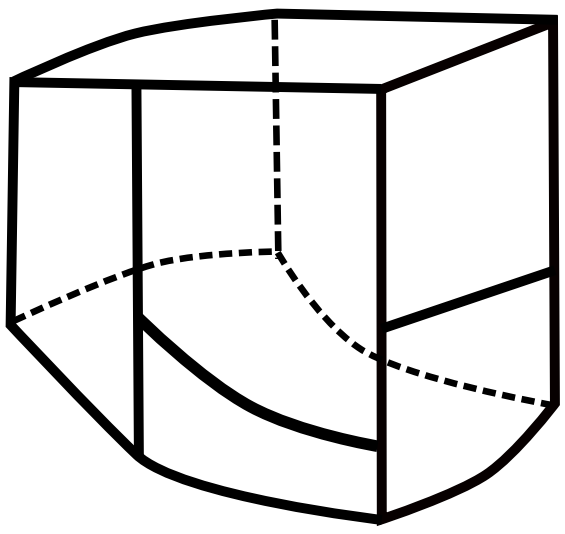

A figure from Stasheff.

The construction of a convex polytope *in some Euclidean space* was a very natural question. According to an anecdote, Milnor came to attend Stasheff's PhD defense with a cardboard model of $K_5$.

The name *associahedron* was coined by Kalai who asked Haiman if there is a geometric (non-abstract) convex polytope in $\mathbb{R}^n$ which realizes $K_{n-2}$. Haiman provided some construction in 1984 but did not publish it. A construction was published by Lee in 1989. Several authors provided other constructions. See the corresponding chapter by Ceballos and Ziegler in Tamari's Festschrift.

## *Loday's construction*

I now describe a beautiful construction due to Jean-Louis Loday[125] in 2004, of a convex polytope *in Euclidean space* whose faces (in the geometrical sense) realize precisely the combinatorics of the Hipparchus et al. associahedron.

[125] J.-L. Loday. Realization of the Stasheff polytope. *Arch. Math. (Basel)*, 83(3):267–278, 2004.

Consider a rooted planar binary tree $T$ with $n$ leaves, thought as a vertex of $K_n$. Label the leaves from 1 to $n$, from left to right. For every pair of leaves $i, j$, denote by $i \vee j$ the node of $T$ which is the least common ancestor of $i$ and $j$. For every integer $i$ such that $1 \leq i \leq n-1$ consider the node $i \vee (i+1)$ and denote by $\nu_l(i)$ (resp. $\nu_r(i)$) the number of its descendant leaves along its left (resp. right) branch. We associate to the tree $T$ the point

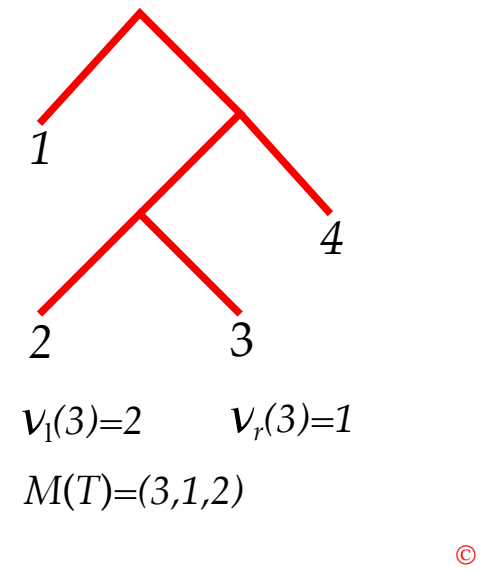

$$M(T) = (\nu_l(1)\nu_r(1), \nu_l(2)\nu_r(2), \ldots, \nu_l(n-1)\nu_r(n-1)) \in \mathbb{R}^{n-1}.$$

**Theorem.** *The convex hull of the set of points $M(T) \in \mathbb{R}^{n-1}$ where $T$ describes all planar rooted binary trees is a convex polytope whose combinatorics is precisely the one of the Hipparchus-Schroeder-Tamari-Stasheff associahedron.*

Let us show first that all the points $M(T)$ lie on the hyperplane of $\mathbb{R}^{n-1}$ whose equation is

$$x_1 + x_2 + \cdots + x_{n-1} = \frac{n(n-1)}{2}.$$

One way of proving this is to count the number of triples $(a, b, v)$ where $a < b$ are two leaves with $v = a \vee b$. Since $v$ is determined by $a, b$ this number is equal to the number of pairs $a < b$, equal to $n(n-1)/2$. Counting the same number according to the node $v$, we get the sum of the $\nu_l(i)\nu_r(i)$ from $i = 1$ to $n-1$. This proves the claim. ⊡

In order to prove Loday's theorem, we first identify the codimension 1 faces $F$ of $K_n$. They are labeled by (non-binary) trees having a single interior node. They are defined by two integers $1 \le p < p+q-1 \le n$ and are obtained by grafting a $q$-corolla at the $p$-th vertex of the $n-q+1$-corolla. The set $F_{p,q}$ of vertices of this face is the set of (rooted planar) *binary* trees which are obtained by grafting any rooted planar binary tree with $q$ leaves to the leaf numbered $p$ in any (planar rooted binary) tree with $(n-q+1)$ leaves.

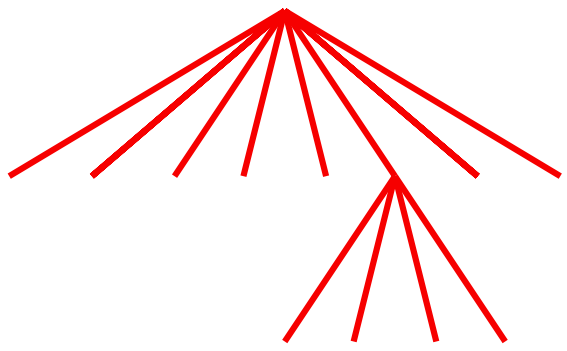

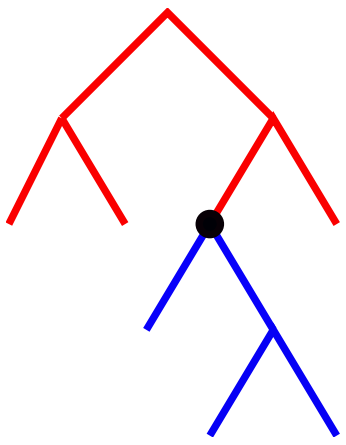

One could express the same thing in still a different way. A binary tree $T$ belongs to $F_{p,q}$ if and only the leaves $\{p,\dots,p+q-1\}$ are the descendants of a single node.

Define a linear function $l_{p,q}$ on $\mathbb{R}^{n-1}$ by:

$$l_{p,q}(x_1,x_2,\dots,x_{n-1}) = x_p + x_{p+1} + \cdots + x_{p+q-1}.$$

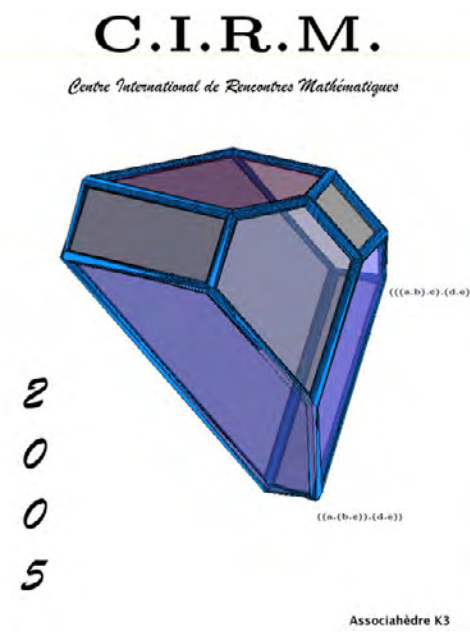

This is the CIRM season greeting card for 2005, representing K5. They even produced a T-shirt. ©

Let us evaluate $l_{p,q}(M(T))$ when $T$ is a vertex of the codimension 1 face $F_{p,q}$. We know that $T$ is a grafting of some rooted planar binary tree $T_1$ with $q$ leaves on the leaf numbered $p$ on some (planar rooted binary) tree $T_0$ with $n-q+1$ leaves. Clearly the $i$-th coordinate of $M(T)$ for $p \le i \le p+q-1$, is the $(i-p)$-th coordinate of $M(T_1) \in \mathbb{R}^{q-1}$ so that $l_{p,q}(T) = q(q-1)/2$ as seen earlier.

Suppose now that $T$ is not in $F_{p,q}$, so there is at least some leaf $i$ in the interval $\{p,\dots,p+q-1\}$ such that $i \vee (i+1)$ has some descendant outside it. While computing $l_{p,q}(M(T))$, one gets $q(q-1)/2$ if one counts only the descendants of $i \vee (i+1)$ which are inside the interval $\{p, p+q-1\}$. Any descendant falling outside the interval $p, p+q-1$ yields a greater sum.

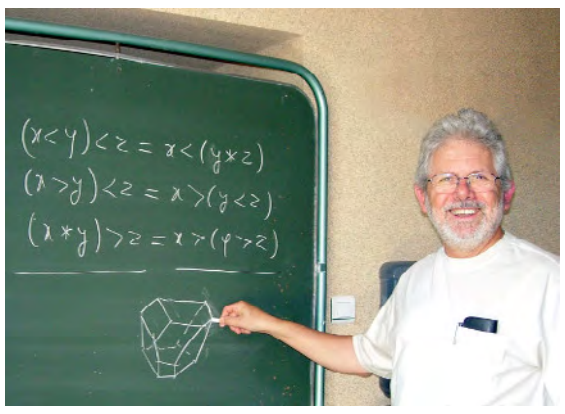

Jean-Louis Loday (1946-2012) lecturing on the associahedron. ©

The proof of the theorem is finished. Indeed for each codimension 1 face $F_{p,q}$ of $K_n$, the affine function $l_{p,q} - q(q-1)/2$ is zero on all $M(T)$ for all vertices $T$ of $F_{p,q}$ and positive on all $M(T)$ for vertices $T$ of $K_n$ which are not in $F_{p,q}$. In other words, we found supporting affine functions which show explicitly that the convex hull of the points $M(T)$ has indeed the combinatorics of our abstract polytope.

Jean-Louis Loday explained the discovery of this embedding in a nice online paper[126]. He frequently used Guillaume William Zinbiel as a pseudonym, due to his admiration for Leibniz.

[126] J.-L. Loday. Comment j'ai trouvé l'associaèdre.

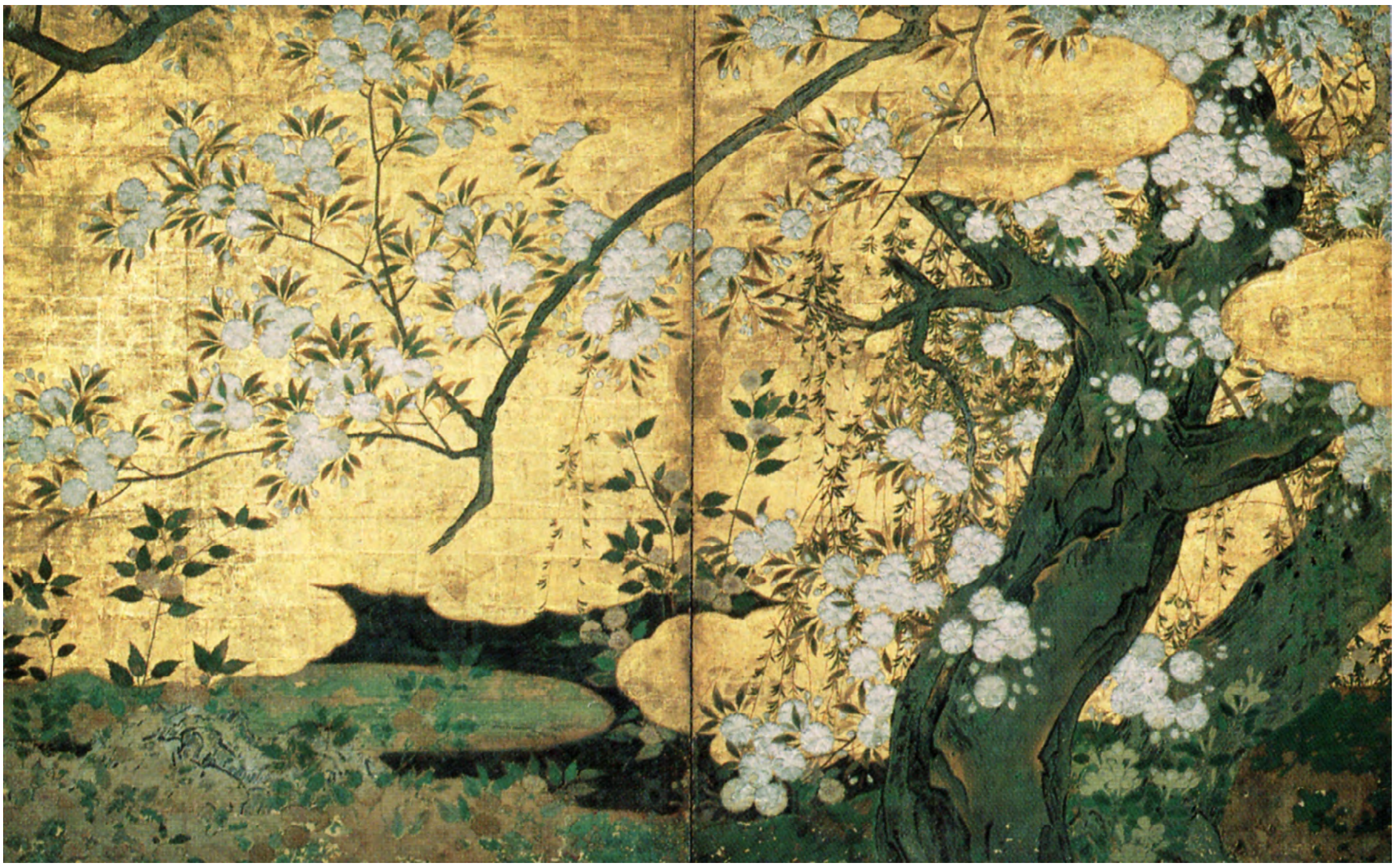

"Cherry Tree"
from Cherry and Maple,
Color Painting of Gold-Foil
Paper (1592).

# Jim Stasheff and loop spaces

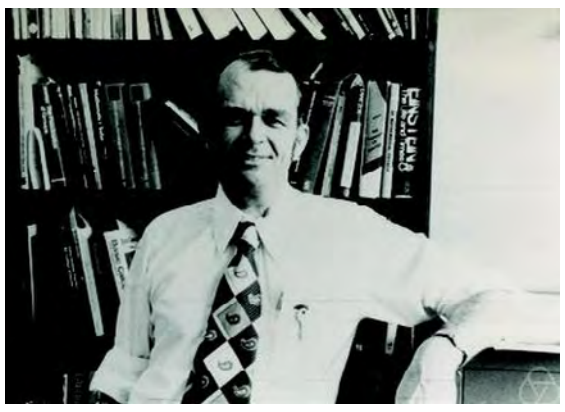

Jim Stasheff. ©

THERE IS NO NEED TO RECALL THE IMPORTANCE OF GROUPS in mathematics in general, and in topology in particular. One of the problems is that this concept is rather subtle in *homotopy theory* as explained in this chapter.

Recall that two continuous maps $f_0, f_1 : X \to Y$ are *homotopic* if there is a continuous map $F : X \times [0,1] \to Y$ (called homotopy) such that $F(x,0) = f_0(x)$ and $F(x,1) = f_1(x)$. Two topological spaces $X, Y$ have the same *homotopy type*, if there are *homotopy equivalences* $f : X \to Y$ and $g : Y \to X$, i.e. maps such that $f \circ g$ and $g \circ f$ are homotopic to the identity. Homotopy theory studies the *homotopy category* whose objects are topological spaces and whose arrows are homotopy classes of maps.

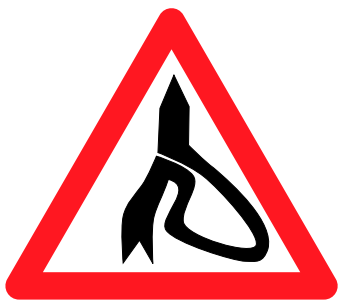

Detour! Strictly speaking, this chapter is not necessary for the rest of the book. It will serve as a motivation for the concept of operad, which is also not necessary, but sheds some light on the global picture. This is probably the most strenuous part of our promenade. ©

As a trivial example, a closed interval cannot be *homeomorphic* to a topological group since the set of its two end points is invariant under any homeomorphism while any group acts transitively onto itself by translations. Nevertheless, the interval is contractible: it has the same homotopy type as a point, which is a (trivial) group. In his 1961 dissertation (published in 1963[127]), Jim Stasheff addressed the question of determining which spaces have the homotopy type of a topological group.

## Topological groups, principal bundles

In this section, I give a very short overview of the role of topological groups in homotopy theory. My only purpose is to introduce enough terminology and basic facts to be able to explain

[127] J. D. Stasheff. Homotopy associativity of $H$-spaces. I, II. *Trans. Amer. Math. Soc. 108 (1963), 275-292; ibid.*, 108:293–312, 1963.

Stasheff's contribution. Once again, this is a huge territory and I have to refer to excellent books, like for instance, the one by... Milnor and Stasheff[128] dealing with the so-called characteristic classes.

[128] J. W. Milnor and J. D. Stasheff. *Characteristic classes*. Princeton University Press, Princeton, N. J.; University of Tokyo Press, Tokyo, 1974. Annals of Mathematics Studies, No. 76.

Let $G$ be *topological group*, i.e. a group equipped with a topology in such a way that the composition and inverse maps are continuous. A *principal $G$-bundle* is a *free action* of $G$ on some space $X$ with a "good quotient" $B$. Every point in $X$ should have a $G$-invariant neighborhood homeomorphic to a product $U \times G$ in which the $G$-action is just the action by translations in the second factor. Frequently, it is more convenient to think of the bundle as the projection map $p : X \to B$ on the space $B$ of $G$-orbits. A $G$-bundle $p' : X' \to B$ is *isomorphic* to $p$ if there is a $G$-equivariant homeomorphism between $X$ and $X'$ inducing the identity on $B$. We will say that the *total space* $X$ is *over* the *base* $B$ and that the inverse image of a point by $p$, which is a $G$-orbit, is *a fiber*.

Traditionally, in the context of principal bundles, groups act *on the right*.

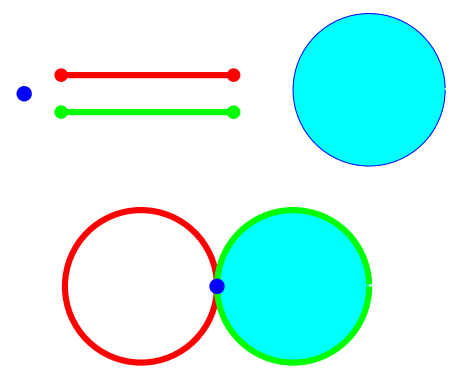

Take one point, two intervals and a disc. Glue the endpoints of the intervals to the point. This produces a figure eight. Glue the boundary of the disc to one component of the eight. You get a very simple example of a CW-complex. Of course, in general, the gluing maps might be much more complicated.

At this point, we should be very cautious about the kind of topological spaces that will be used. They should be Hausdorff and should not be too pathological. Usually, one restricts the study to *CW-complexes*. It is not my intention to give a precise description of these spaces. I will only mention that such a space $X$ is by definition an increasing union of subspaces $Sk_n(X)$, called their $n$-th *skeletons*. The $(n+1)$-st skeleton $Sk_{n+1}(X)$ is obtained from $Sk_n(X)$ by gluing some $(n+1)$-dimensional balls $B^{n+1}$ along some "attaching maps" $u : \partial B^{n+1} \to Sk_n(X)$. Hatcher's book[129] (freely available on the Internet) is an excellent reference.

[129] A. Hatcher. *Algebraic topology*. Cambridge University Press, Cambridge, 2002.

Principal bundles are fundamental objects in (differential) topology. For instance, given a smooth manifold $M$ of dimension $m$, look at the space $Fr(M)$ of pairs $(x, f)$ where $x$ is a point of $M$ and $f$ is a *frame* at $x$, in other words a basis of the tangent space $T_x(M)$. There is an obvious free action of the linear group $GL(m, \mathbb{R})$ on $Fr(M)$ and the map $p$ sending $(x, f) \in Fr(M)$ to $x \in M$ is a principal bundle.

Given a $G$-principal bundle $p : X \to B$ and a map $i : B_1 \to B$, one can *pull-back* $p$ to produce a principal bundle $p_1 : X_1 \to B_1$. Formally, $X_1$ is the subspace of $B_1 \times X$ consisting of couples $(b_1, x)$ such that $i(b_1) = p(x)$ and $p_1(b_1, x) = b_1$. For instance, if $i$ is an inclusion, $p_1$ is the restriction to "what is above $i(B_1)$ in $X$".

Here is an important example. Let $Gr_{k,n}$ be the space of linear subspaces of dimension $k$ in $\mathbb{R}^n$. This is a compact manifold, called a *Grassmannian manifold*. There is a *tautological* $GL(k,\mathbb{R})$ *bundle* over $Gr_{k,n}$ whose fiber over some subspace consists of the bases of that subspace. If a $k$ dimensional manifold $M$ is immersed in $\mathbb{R}^n$, then the differential of this immersion gives a map from $M$ to $Gr_{k,n}$. The pullback of the $GL(k,\mathbb{R})$ tautological bundle over $Gr_{k,n}$ is isomorphic to the frame bundle of $M$.

Let me mention two important facts concerning principal bundles.

– Any bundle with a contractible basis $B$ is trivial(izable), i.e. isomorphic to $B \times G$. See[130] for the history of this theorem.

– If $i, i' : B_1 \to B$ are homotopic, then the pull-back principal bundles $p_1, p_1'$ of $p$ by $i, i'$ are isomorphic.

These two properties show that the set of isomorphism classes of $G$-principal bundles over some space $B$ only depends on the homotopy type of $B$, and defines some contravariant functor on the homotopy category.

[130] M. Audin. Publier sous l'Occupation. I. Autour du cas de Jacques Feldbau et de l'Académie des sciences. *Rev. Histoire Math.*, 15(1):7–57, 2009.

## *Classifying spaces*

A $G$-principal bundle $p_G : E(G) \to B(G)$ is called *universal* if every principal $G$-bundle $p : X \to B$ is isomorphic to the pull-back of $p_G$ by some map $i : B \to B(G)$ which is unique, up to homotopy. Later, we will sketch a proof of the following.

**Theorem.** *For every topological group $G$, there exists a* universal *fiber bundle $p_G : E(G) \to B(G)$.*

In other words, there is a natural bijection between:

– (isomorphism classes of) $G$-principal bundles over some space $B$.

– Homotopy classes $[B, B(G)]$ of maps from $B$ to $B(G)$.

We will say that $B(G)$ is the *classifying space of* $G$.

Let me describe two important examples. Suppose first that $G$ is a *discrete group*. In such a situation, a $G$-principal bundle is nothing more than a Galois covering map with Galois group $G$.

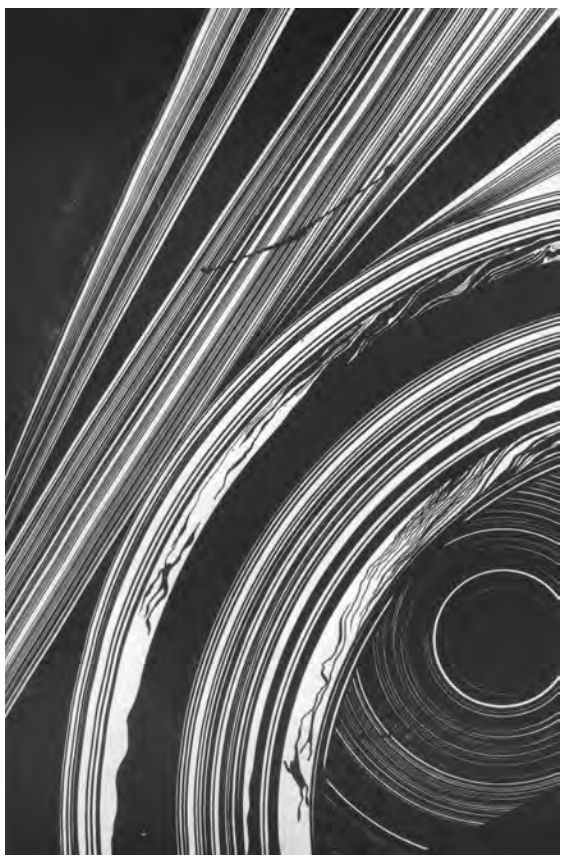

Spectral sequences and orbits of the action of groups, by A. Fomenko. ©

Covering spaces of a space $B$ are described by subgroups of its fundamental group. In this case $B(G)$ is an *Eilenberg-MacLane space* $K(G,1)$: its fundamental group is $G$ and its universal cover $E(G)$ is contractible (equivalently all higher homotopy groups of $K(G,1)$ are trivial). A $G$-Galois covering of some space $B$ is equivalent to a homotopy class of maps $B \to K(G,1)$.

For the theory of covering spaces, see Hatcher's book, or Analysis Situs.

*An important remark. In order to define the fundamental group of some space $X$, we need a base point $x \in X$. A notation like $f:(X,x) \to (Y,y)$ means that $f(x) = y$. When discussing homotopy of maps, I should mention explicitly if this homotopy preserves base points or not.* I should... but I will not! This would imply long and technical sentences and, as the reader has already noticed, this book is not a complete encyclopedia. I hope my reader will forgive this lack of precision.

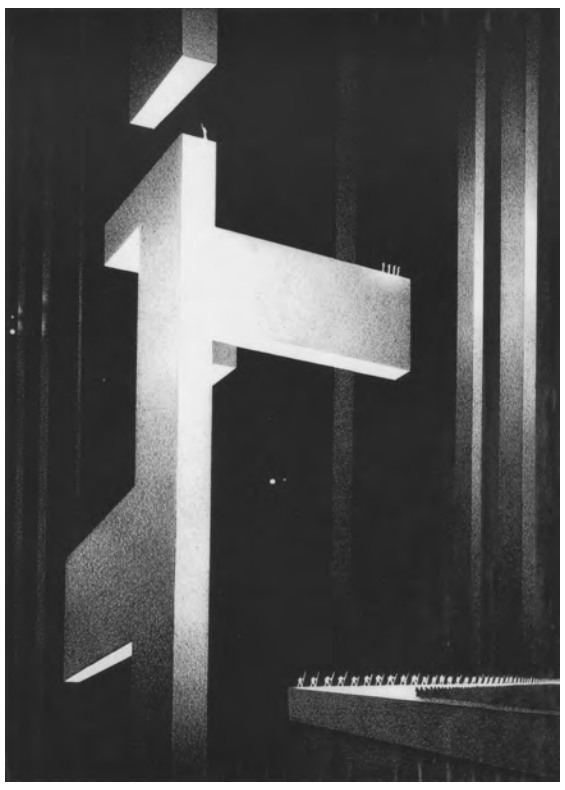

Polyhedra and simplicial chains 1973, A. Fomenko. ©

As a second example, consider the group $U(1)$ of complex numbers of modulus 1. For every $n$, there is an action of $U(1)$ on the unit sphere $\mathbb{S}^{2n-1}$ in $\mathbb{C}^n$. The element $\omega$ acts on $(z_1,\ldots,z_n)$ to produce $(\omega z_1,\ldots,\omega z_n)$. This defines a principal $U(1)$-bundle

$$p_n : \mathbb{S}^{2n-1} \to \mathbb{C}P^n$$

over the complex projective space. All these spheres and projective spaces are naturally nested, embedding $(z_1,\ldots,z_n)$ to $(z_1,\ldots,z_n,0,\ldots,0)$, so that we can define a principal $U(1)$-bundle whose total space is the infinite dimensional sphere and whose basis is the infinite dimensional projective space

$$p_\infty : \mathbb{S}^\infty \to \mathbb{C}P^\infty.$$

We will see that this is the universal bundle for $G = U(1)$. The key point is the following.

**Proposition.** *A $G$-bundle is universal if and only if its total space is contractible.*

The proof of this fundamental fact is a typical example of *obstruction theory*. Start with some $G$-bundle $p_G : E(G) \to B(G)$ such that $E(G)$ is contractible and let us show that it is universal. Consider some other $G$-bundle $p : E \to B$ and we want to show that it is the pullback of $p_G$ by some map $i : B \to B(G)$.

We construct $i$ on the skeletons of $B$ (which is, as always, a CW-complex), by induction on their dimensions. At each step, some continuous map has to be extended and the contractibility of $E(G)$ is precisely what is needed for this construction. See Milnor and Stasheff's book for the details and for the proof of the converse.

Check everything yourself, without opening too much the recommended books! Why is the sphere $S^\infty$ contractible?

In the example of $U(1)$ the classifying space $B(U(1))$ is $\mathbb{C}P^\infty$ since the infinite dimensional sphere is indeed contractible.

## *Milnor's join construction*

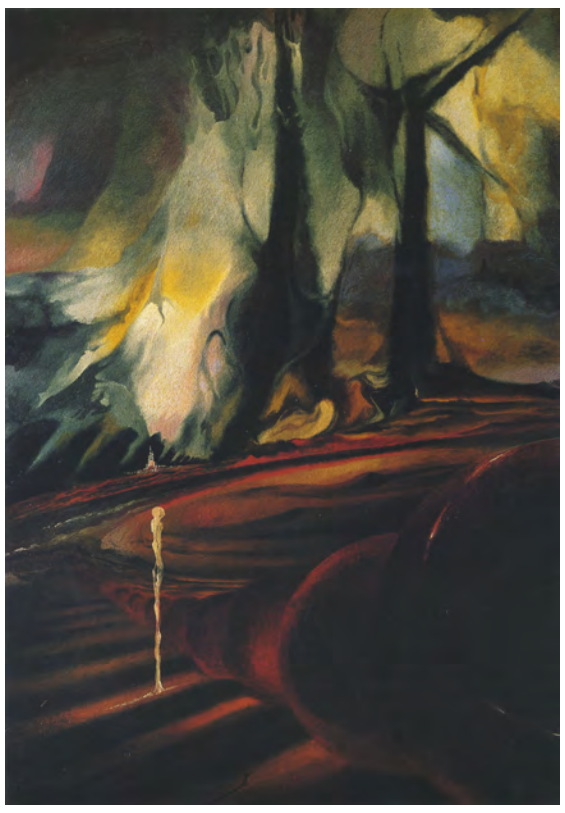

Cellular spaces, by A. Fomenko. ©

Milnor's construction of $B(G)$ is beautiful and easy[131]. Any topological group $G$ acts freely on itself but of course the group needs not be contractible. Therefore, we have to force the contractibility, preserving a free group action. Let $E(G)$ be the *simplex over* $G$. An element of $E(G)$ is by definition some finite formal barycentric combination of elements of $G$, that is to say a formal sum

$$\lambda_0 g_0 + \lambda_1 g_1 + \cdots + \lambda_k g_k$$

where $\lambda_i \geq 0$ and $\sum_i \lambda_i = 1$. This space is convex, hence contractible, and is equipped with a free action of $G$. The projection of $E(G)$ on its quotient $B(G)$ is therefore a classifying space. Et voilà !

[131] J. Milnor. Construction of universal bundles. II. *Ann. of Math. (2)*, 63:430–436, 1956.

We should be more careful in the definition of $E(G)$. Start with the disjoint union of products $G^{n+1} \times \Delta_n$ where

$$\Delta_n = \{(\lambda_0, \ldots, \lambda_n) \mid \lambda_i \geq 0 \text{ and } \sum_i \lambda_i = 1\}$$

is the standard simplex. Then, introduce an "obvious" equivalence relation generated by

$$\begin{aligned} &((g_0, \ldots, g_i, g_{i+1}, g_{i+2}, \ldots, g_n), (\lambda_0, \ldots, \lambda_i, 0, \lambda_{i+2}, \ldots, \lambda_n)) \in G^{n+1} \times \Delta_n \\ &\equiv ((g_0, \ldots, g_i, g_{i+2}, \ldots, g_n), (\lambda_0, \ldots, \lambda_i, \lambda_{i+2}, \ldots, \lambda_n)) \in G^n \times \Delta_{n-1} \end{aligned}$$

and we define $E(G)$ as the quotient space. This is *Milnor's join construction*. The name comes from the fact that "virtual connections" have been created, joining points in $G$.

## *Loops and their composition*

Given a space $B$ with a base point $\star \in B$, the *based loop space* $\Omega(B,\star)$ is the space of... based loops ☺, i.e. of continuous maps $\gamma : [0,1] \to B$ such that $\gamma(0) = \gamma(1) = \star$, equipped with the compact open topology. Given two based loops $\gamma_1, \gamma_2$, they can be *concatenated*. One possible definition is to set $\gamma_1 \bullet \gamma_2(t)$ as $\gamma_1(2t)$ for $0 \le t \le 1/2$ and $\gamma_2(2t-1)$ for $1/2 \le t \le 1$. This composition map

$$\Omega(B,\star) \times \Omega(B,\star) \to \Omega(B,\star)$$

is certainly not associative. In the composition $(\gamma_1 \bullet \gamma_2) \bullet \gamma_3$ we go along $\gamma_1$ when $t \in [0,1/4]$, then along $\gamma_2$ when $t \in [1/4,1/2]$ and finally along $\gamma_3$ when $t \in [1/2,1]$. This is not the same path as $\gamma_1 \bullet (\gamma_2 \bullet \gamma_3)$, even though these two loops are homotopic.

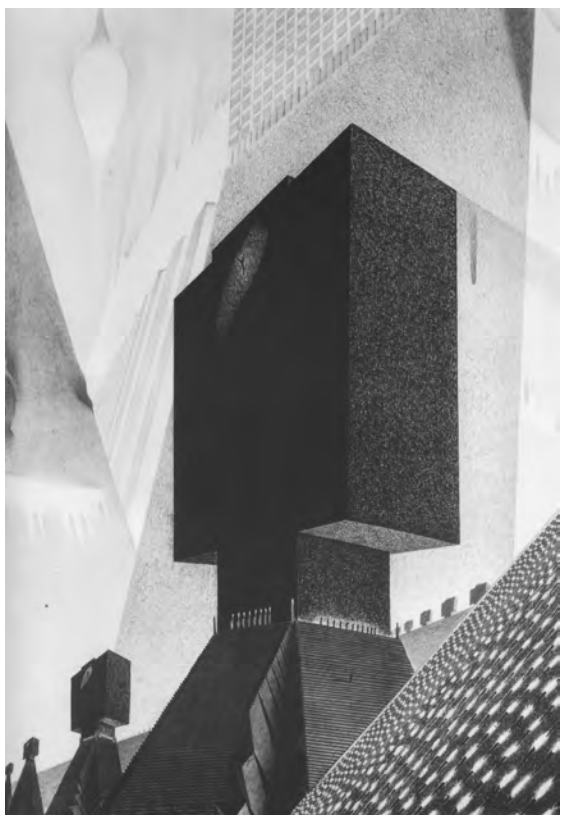

Simplicial complexes 1973, by A. Fomenko. ©

Changing a little bit the definitions, we get a loop space which is *strictly associative*, not only up to homotopy. Let us use the so-called *Moore loops*. Such a loop consists of some number $l \ge 0$ (thought as some length) and some continuous map $\gamma : \mathbb{R}_+ \to B$ such that $\gamma(0) = \star$ and $\gamma(t) = \star$ for $t \ge l$. There is a natural topology on the space of these fancy loops, denoted $\Omega_M(B,\star)$, which has the same homotopy type as $\Omega(B,\star)$. Given $(l_1,\gamma_1)$ and $(l_2,\gamma_2)$, their composition is defined as $(l_1 + l_2, \tilde{\gamma})$ where $\tilde{\gamma}(t) = \gamma_1(t)$ for $t \le l_1$ and $\gamma_2(t - l_1)$ for $t \ge l_1$. This is clearly associative.

There is another trick, less well known, to get another space, still with the same homotopy type as $\Omega(B,\star)$, which is now a *topological group*. This is due to Milnor (encore lui!) and described in a book[132] of Stasheff. We make a very mild assumption: $B$ is the geometric realization of a simplicial complex with a countable number of faces. Define a group $G(B)$ in the following manner. Start with the disjoint union of the $B^n$'s for $n \ge 0$. Think of an element of $B^n$ as a discrete path $b_1, \dots, b_n$ with $n$ steps, where a point, instead of following a continuous path, hops from point to point. Consider the equivalence relation where

$$(b_1, \dots, b_{i-1}, b_i, b_{i+1}, \dots, b_n) \equiv (b_1, \dots, b_{i-1}, b_{i+1}, \dots, b_n)$$

if $b_i = b_{i+1}$ or $b_{i-1} = b_i$.

[132] J. Stasheff. *H-spaces from a homotopy point of view*. Lecture Notes in Mathematics, Vol. 161. Springer-Verlag, Berlin-New York, 1970.

A *simplicial complex* is a combinatorial concept. It consists of a set $V$ whose elements are called *vertices*, and a family of finite subsets of $V$ whose elements are called *faces*. The only axiom is that a non-empty subset of a face is a face. Given a simplicial complex, there is an associated topological space called its *geometric realization*. It consists of functions $t : V \to [0,1]$ such that $\sum_{x \in V} t(x) = 1$ and such that $\{x | t(x) \neq 0\}$ is a face.

In the quotient space, define $G(B)$ as the subspace of classes which have a representative $(b_1, \ldots, b_n)$ such that $b_1 = b_n = \star$ and such that any two consecutive elements $b_i, b_{i+1}$ are in the same simplex (so that, mentally, we can connect them by a segment). The group structure is just concatenation. It is a simple exercise to check that this is indeed a topological group, with the same homotopy type as the loop space $\Omega(B, \star)$.

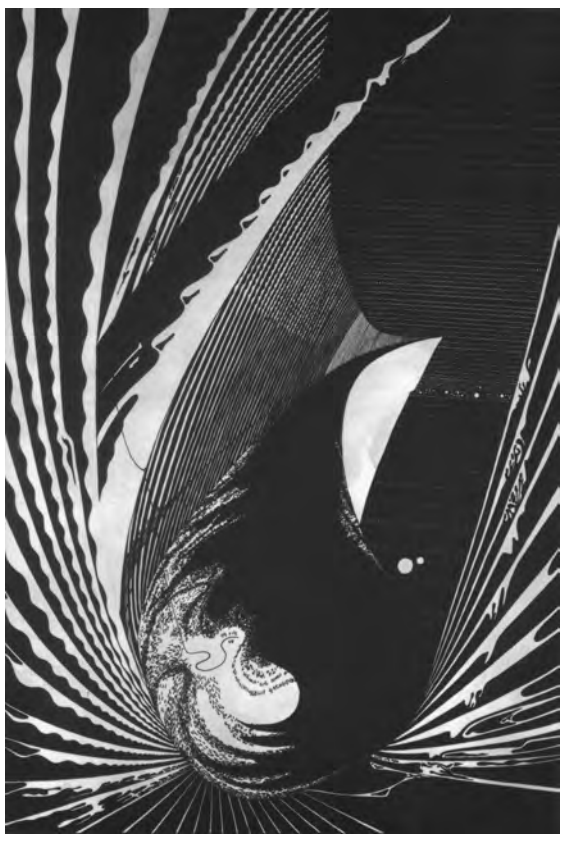

A fiber space, by A. Fomenko.

Well, this construction is not so complicated, but to keep in mind that the group that was produced is rather huge even if $B$ is very simple. This group is very seldom used "in practice".

Anyway, we should remember that a space $(B, \star)$ defines a useful $\Omega(B, \star)$ equipped with some concatenation map, which is not associative but that can be turned into an associative law or even into a topological group, at the cost of some topological contortions.

A final remark in this section:

**Proposition.** *Any topological group $G$ has the same homotopy type as the loop space of its classifying space.*

I only list the keywords in the proof in order to illustrate the kind of gymnastics required in this part of topology. If $(X, x)$ is a pointed space, its *suspension* $S(X, x)$ is obtained from $X \times [0, 1]$ by collapsing $X \times \{0\}$, $X \times \{1\}$ and $\{x\} \times [0, 1]$, to a single point. A map from $S(X, x)$ to some other space $(Y, y)$ is equivalent to a map from $(X, x)$ to the loop space $\Omega(Y, y)$. The following are equivalent:

– a homotopy class of maps from $(X, x)$ to $\Omega(B(G), e)$ ($e$ is the unit in $G$),

– a homotopy class of maps from $S(X, x)$ to $(B(G), e)$,

– an isomorphism class of a $G$-bundle over $X \times [0, 1]$ trivialized over $X \times \{0\}$, $X \times \{1\}$, and $\{x\} \times [0, 1]$,

– a path of $G$-bundles $p_t$ over $X$ and isomorphisms between $p_0$ and $p_1$ with the trivial bundle $X \times G$.

Observe that given a trivialized bundle $X \times G \to X$, the other trivializations are simply given by maps $X \to G$. Indeed an isomorphism from $X \times G \to X$ to itself sends $(b, g)$ to $(b, u(g)g)$ for some $u : B \to G$. Therefore the homotopy classes of maps from $(X,x)$ to $\Omega(B(G),e)$ are in canonical bijections with the homotopy classes of maps from $(X,x)$ to $(G,e)$. ⊡

Note that the fact that right and left translations commute is nothing more than associativity.

## *Stasheff's theorem on H-spaces*

A space $(X,\star)$ is called a *H-space* if it is equipped with a "multiplication"

$$m_2 : X \times X \to X$$

such that $m_2(x,\star) = m_2(\star,x) = x$.

*H* is in honor of Heinz **H**opf, and not of **H**omotopy.

The question studied by Stasheff is called the *recognition problem*. Is it possible to decide from $X$ and $m_2$ if there is some space $Y$ and some *homotopy equivalence* from $(X,x)$ to the loop space $\Omega(Y,y)$ that transforms $m_2$ in the concatenation of loops in $\Omega(Y,y)$? This is equivalent to the recognition of topological groups among $H$-spaces, up to homotopy.

This is indeed a fundamental question:

*What is the right concept of group in the homotopy category?*

As the reader has certainly guessed, the answer given by Stasheff will involve the associahedron introduced in the previous chapter.

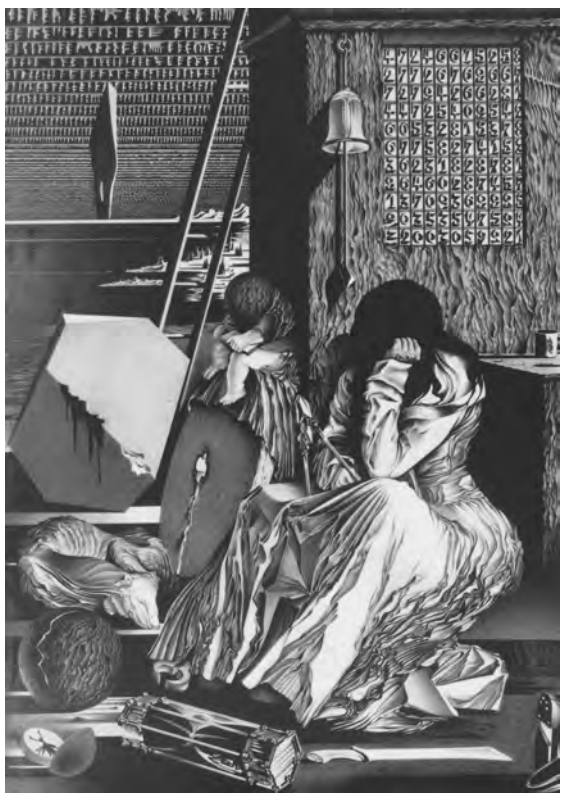

Anti-Durer - From the cycle - Dialogue with authors of the 16th century 1975, by A. Fomenko. ©

Suppose that $(X,x)$ has indeed the homotopy type of some loop space $\Omega(Y,y)$ and that $m_2$ is homotopic to concatenation. It is convenient to use the Moore loop space $\Omega_M(Y,y)$ with its *associative concatenation* $\mu$. The two maps

$$\begin{aligned}(x_1,x_2,x_3) \in X^3 \quad &\mapsto m_2(m_2(x_1,x_2),x_3) \in X \\ &\mapsto m_2(x_1,m_2(x_2,x_3)) \in X\end{aligned}$$

should be homotopic since $\mu$ is associative. This condition does not necessarily hold. If this happens, the $H$-space is associative *up to homotopy* and we say that the *H-space* $(X,m_2)$ *is an* $A_1$*-space*. The homotopy in this case is a map from $X^3 \times [0,1]$ to $X$ and the factor $[0,1]$ should be seen as the associahedron $K_3$.

Four terms define five maps $X^4 \to X$ associated to the five planar rooted binary trees with four leaves. These five trees can be seen as the vertices of a pentagon $K_4$. In the previous step, we have considered five maps $[0,1] \times X^4 \to X$. These maps agree on their boundaries and define a map $\partial K_4 \times X^4 \to X$. If $X$ has the homotopy type of a loop space $\Omega_M(Y)$, with strictly associative multiplication $\mu$, this map has to extend to the full pentagon $K_4 \times X^4 \to X$. If this happens, we say that $X$ *is an* $A_2$*-space.*

It should be clear that this picture continues in all dimensions.

We can at last state Stasheff's theorem:

**Theorem.** *A H-space* $m_2 : X \times X \to X$ *is homotopically equivalent to some loop space if and only if it is an* $A_\infty$*-space, i.e. if there exist coherent maps* $m_n : K_n \times X^n \to X$ *compatible with the faces of the associahedron* $K_n$ $(n \geq 1)$.

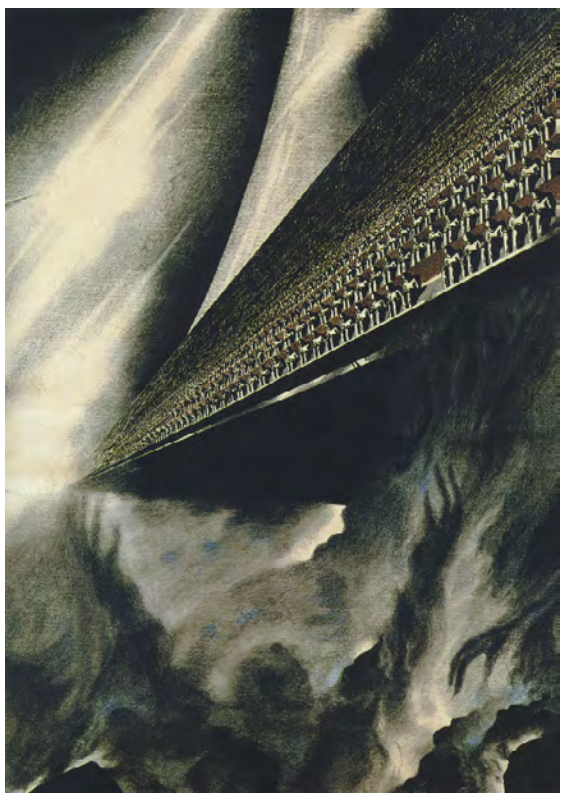

The method of killing spaces in homotopic topology, by A. Fomenko. ©

The necessary condition is clear by now. The most interesting part of the theorem is of course the sufficient condition, that I will *not* prove in the next section.

Start with some $A_\infty$-space $X$, with compatible maps $m_n : K_n \times X^n \to X$, and our purpose is to produce a space $Y$ whose loop space has the same homotopy type as $X$. We know that $\Omega(Y, y)$ has the same homotopy type as some topological group $G$, which in turn, has the homotopy type of $\Omega(B(G), \star)$. Therefore, it is tempting to choose $Y = B(G)$, but we don't know $G$.

We only know the collection of maps $m_n : K_n \times X^n \to X$ which are some kind of a substitute for a group structure. Therefore, our strategy is clear. We have to adapt Milnor's join construction of $B(G)$ to these more general $A_\infty$-structures. This project has been carried out by Stasheff.

Instead of starting with the disjoint union of $G^{n+1} \times \Delta_n$ and identifying points according to some "obvious" equivalence relation, start with the disjoint union of $K_n \times X^n$ and define some "obvious" equivalence relation in this disjoint union. This will produce a space $B(X)$ which is the *classifying space* of the $A_\infty$-space $X$.

The "only thing that we still has to show" is that, as expected, the loop space of $B(X)$ is a solution to our problem: that is to say a *delooping* of $X$. This is not easy and Stasheff proved it with some additional (minor) hypothesis on the topology of $X$.

## *Cherry trees*

In order to get some intuition behind this $B(X)$, let me describe briefly the *cherry trees* introduced by Boardman and Vogt[133]. Consider a rooted planar binary tree with $n$ leaves. If the $n-2$ internal edges are equipped with some length in $[0,1]$, we get a *metric* tree. This space of metric trees defines a cube $[0,1]^{n-2}$ for each binary tree.

[133] J. M. Boardman and R. M. Vogt. *Homotopy invariant algebraic structures on topological spaces*. Lecture Notes in Mathematics, Vol. 347. Springer-Verlag, Berlin-New York, 1973.

If one, or more, internal edges have length 0, these edges can be collapsed and the result is a rooted planar tree, which is not binary anymore, but whose internal edges still have a length. This produces some identifications along the boundary of those cubes. The set of these metric trees defines a cubical decomposition of $K_n$. Therefore $K_n$ can be viewed as a space of metric trees. For example the pentagon is decomposed in five squares. This presentation of $K_n$ as a space of metric trees enables us to define *grafting maps* $\iota_{k_1,\ldots,k_n}$ :

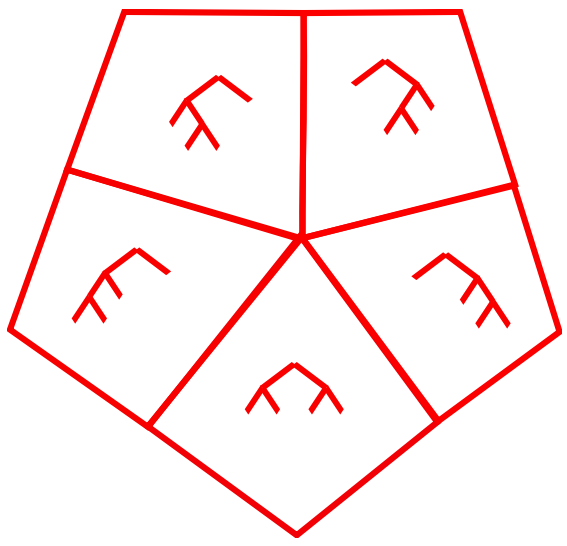

$$K_n \times \left(K_{k_1} \times K_{k_2} \times \cdots \times K_{k_n}\right) \to K_{k_1+\cdots+k_n}.$$

Simply attach metric trees at the leaves of a metric tree.

Now, we want to picture $K_n \times X^n$. Simply imagine that each leaf of our metric trees carries some element of $X$, thought as a cherry. This is a *cherry tree*: a metric tree with cherries on the leaves.

The classifying space $B(X)$ can be described using this terminology. Suppose that a cherry tree $T$ contains an edge $e$ which is *fully grown* of length 1. Cutting along $e$, the tree $T$ is decomposed in two metric trees. Let us denote by $T_1$ the part consisting of descendants of the endpoint of $e$: this is a cherry tree with $k \leq n$ leaves. The other tree, $T_0$, containing the root of $T$, is not quite a cherry tree since the newly created leaf, at the origin of $e$, is not equipped with a cherry. We can now evaluate $m_k$ on the cherry tree $T_1$ and deposit the result as a new cherry on the leaf of $T_0$ which was waiting for its cherry. This produces a new cherry tree.

By definition, the classifying space $B(X)$ is the quotient of the space of cherry trees by this operation of cutting fully grown edges and applying $m_k$ as explained.

I will *not* show that the loop space of $B(X)$ has the homotopy type of $X$. I have to admit that the references that I gave do contain proofs but are definitely not easy to read ☹.

The most readable reference that I know is[134].

[134] E. Hoefel, M. Livernet, and J. Stasheff. $A_\infty$-actions and recognition of relative loop spaces. 2013.

©

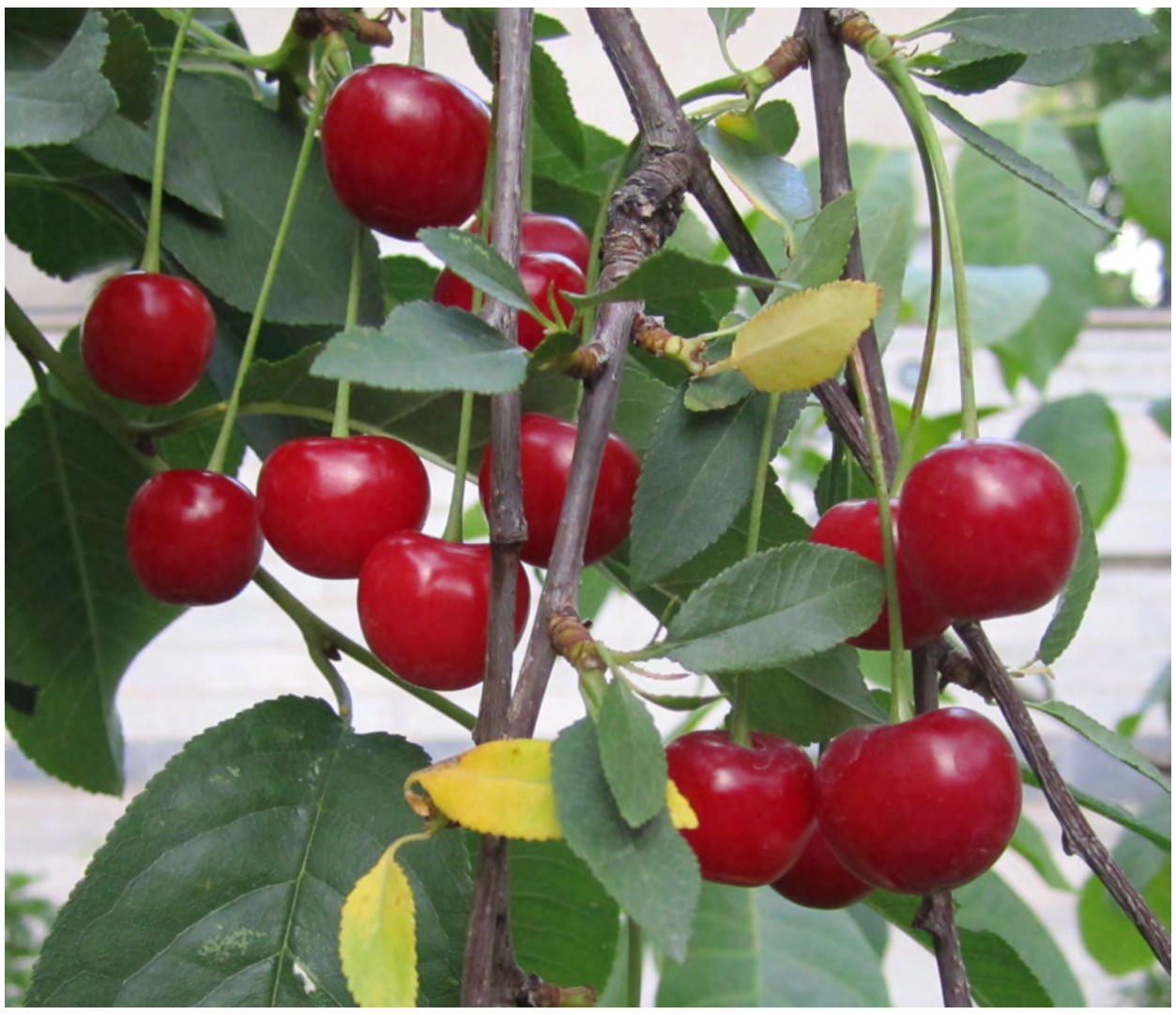

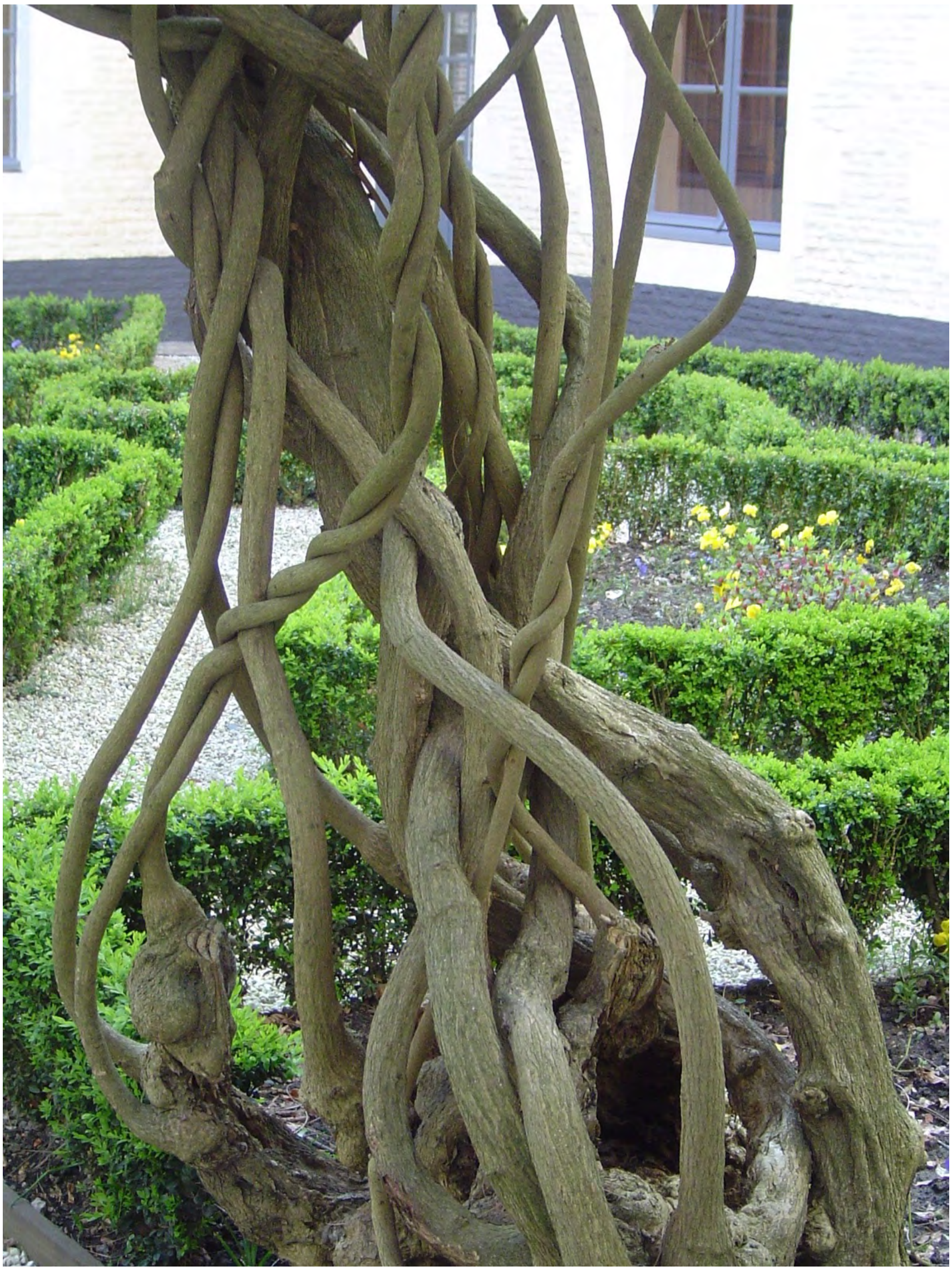

A braided tree,
illustrating an operad?

# *Operads*

Let me begin this chapter with a quote from Peter May[135].

> The name 'operad' is a word that I coined myself, spending a week thinking about nothing else. Besides having a nice ring to it, the name is meant to bring to mind both operations and monads. [...] What I did not foresee was just how flexible the notion would be, how many essentially different mathematical contexts there are in which it would play a natural role, how many philosophically different ways it could be exploited.

[135] J. P. May. Operads, algebras and modules. In *Operads: Proceedings of Renaissance Conferences (Hartford, CT/Luminy, 1995)*, volume 202 of *Contemp. Math.*, pages 15–31. Amer. Math. Soc., Providence, RI, 1997.

According to Wikipedia, another reason for this name is that May's mother was an opera singer. As almost all the concepts that we have met so far, operads "existed" much before their birth[136], or better to say, before they were baptized. May's definition is aimed at encapsulating many kinds of *operations*, most of them reminiscent of the grafting of trees that we already encountered.

[136] J. Stasheff. The pre-history of operads. In *Operads: Proceedings of Renaissance Conferences (Hartford, CT/Luminy, 1995)*, volume 202 of *Contemp. Math.*, pages 9–14. Amer. Math. Soc., Providence, RI, 1997.

A group is much more than a set equipped with some multiplication map satisfying some axioms. Groups only exist through their representations as automorphisms of "something". In the same way, operads only exist through their representations and we will not spend too much time on the abstract definitions.

In his talk *On teaching mathematics*, V. Arnold strongly criticizes the axiomatic approach to group theory, as it is usually taught in France. I fully agree with him.

An operad consists of

- sets $\mathcal{O}_n$ for $n \geq 0$ (thought as $n$-ary operations).
- some element $1$ in $\mathcal{O}_1$ called the unit,
- for every $n, k_1, \ldots, k_n$, an *operad operation*, i.e. a map

$$\mathcal{O}_n \times \left(\mathcal{O}_{k_1} \times \mathcal{O}_{k_2} \times \cdots \times \mathcal{O}_{k_n}\right) \to \mathcal{O}_{k_1+\cdots+k_n}$$

satisfying...some axioms. I don't want to write down the formulas expressing these axioms since I would be unable myself to read the formulas that I wrote. I prefer to give first an example (that the reader has probably already guessed) before describing the axioms in words.

Computer scientists taught us that some computable bijections $\mathbb{N} \to \mathbb{N}$ are easy to evaluate and have very complicated inverses. Writing a formula is usually easy and understanding it might be terribly complicated.

The example is given by planar binary rooted trees. Denote by $\mathcal{O}_n$ the set of planar binary rooted trees with $n$ leaves, and let 1 the tree with one leaf which is at the same time the root. The grafting operation that we used several times defines the easiest example of an operad.

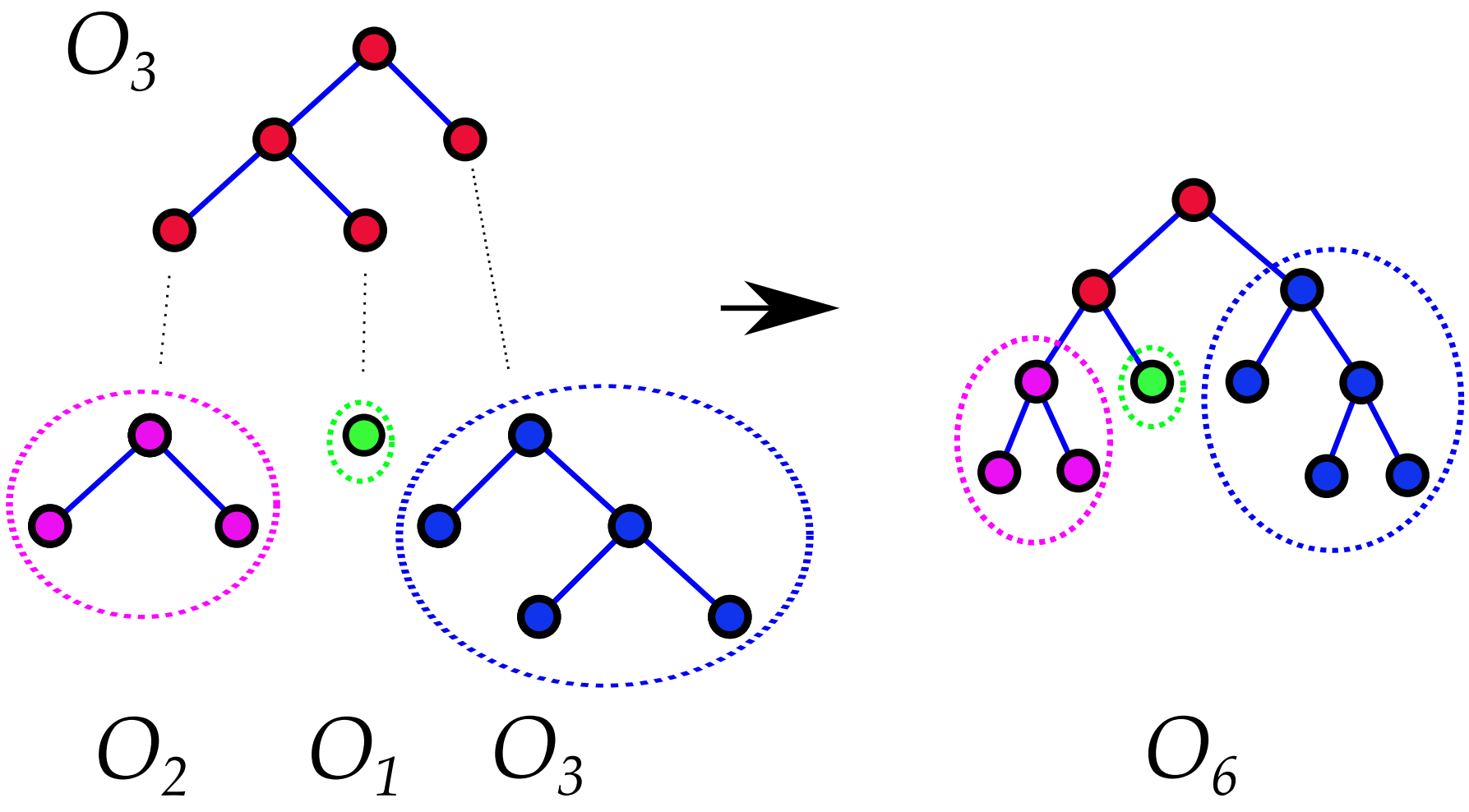

Now, what are the axioms for a general operad? Grafting 1 on some tree does not change the tree. If we graft some $B$ on some tree $A$ and then graft $C$ on the result, we could as well have grafted $C$ on $B$ and grafted the result on $A$. The operad axioms are nothing more than that, replacing the word *grafting* by the operad operation. Some unit condition and some sort of associativity.

Of course, we could also graft rooted planar trees which are not necessarily binary. For instance, we could use *pruned trees*, i.e. rooted planar trees such that every interior node has at least two children. This produces the *Hipparchus-Schroeder-Tamari-Stasheff operad*.

For a two page introduction to operads, see[137]. For a 27 page presentation, see[138]. For a more recent 634 page book on the same topic, see[139].

[137] J. Stasheff. What is ... an operad? *Notices Amer. Math. Soc.*, 51(6):630–631, 2004.

[138] F. Chapoton. Operads and algebraic combinatorics of trees, journal = Sém. Lothar. Combin. 58:Art. B58c, 27, 2007/08.

[139] J.-L. Loday and B. Vallette. *Algebraic operads*, volume 346 of *Grundlehren der Mathematischen Wissenschaften [Fundamental Principles of Mathematical Sciences]*. Springer, Heidelberg, 2012.

Here is another example of a naive operad. Choose some set $E$ and define $\mathcal{O}_n$ as the set of maps $E^n \to E$. In this example, 1 is the identity, and the operad operations are simply given by *substitution*. If you have a map $f : E^n \to E$ and $n$ maps $f_i : E^{k_i} \to E$, just replace $x_i$ by $f_i$ in $f(x_1, \dots, x_n)$ to get a map $E^{k_1+\dots+k_n} \to E$. This satisfies the axioms that we did not write down... This operad is denoted by $End(E)$.

An *algebra* over some operad $\mathcal{O}$ (also called a representation) is a set $E$ and some operad map from $\mathcal{O}$ to $End(E)$. In other words, each element of $\mathcal{O}_n$ is incarnated as an $n$-ary operation $E^n \to E$ in a *compatible way*.

We could work in many different categories. Instead of sets, we could use topological spaces, homotopy types, vector spaces etc.

I now describe more interesting examples of operads. More will come in the following chapters.

## *Permutations*

Recall from the beginning of this book that we studied some combinatorial questions concerning *pattern recognition*. This fits very well with the following operad.

$\mathcal{O}_n$ is the set of permutations of $\{1, \dots, n\}$. I wrote *"set"* and not *"group"* because we are not going to compose these permutations. Instead, we think of a permutation as a total order on $\{1, \dots, n\}$.

Suppose that we have total orderings $\sigma, \sigma_1, \dots, \sigma_n$ on $\{1, \dots, n\}$ and $\{1, \dots, k_1\}, \dots, \{1, \dots, k_n\}$. Write $\{1, \dots, k_1 + \dots + k_n\}$ as a disjoint union of $n$ consecutive intervals of sizes $k_1, k_2, \dots, k_n$. Order these intervals according to $\sigma$ and inside the $i$-th interval, order the elements according to the order $\sigma_i$. This produces a natural total ordering on $\{1, \dots, k_1 + \dots + k_n\}$. Clearly these *grafting operations* on orderings satisfy the axioms of an operad.

In the chapter dealing with separable permutations we noticed that any pruned planar rooted tree defines a permutation of the leaves.

In other words, this defines an operad map from the Hipparchus-Schroeder-Tamari-Stasheff operad to the Permutation operad. We also proved that this is an embedding, as the tree can be reconstructed from the permutation.

## *The free operad: Hipparchus-Schroeder again*

Consider a sequence of sets $(E_n)_{n \geq 1}$ and let us define the *free operad generated by* the $E_n$'s. We have to create sets $\mathcal{O}_n$'s whose elements are produced under the operad operations starting from elements of the $E_n$'s. Since we want a *free operad*, all these new elements should be assumed to be different, unless some use of the axioms implies that they are equal. It is not difficult to construct this free operad.

I am conscious of the fact that this is not a precise definition but I am reluctant to define this using initial objects in categories.

An element of $\mathcal{O}_n$ is a rooted planar tree with $n$ leaves, such that each node with $i$ children is equipped with some label belonging to $E_i$. The operad operations are again defined by grafting.

As an example, consider the case where all $E_n$'s are empty except $E_2$ containing one element. Then the free operad on one element "of degree 2" is the operad of rooted binary trees. An algebra over this operad is just a set with a binary operation.

As another example, let us consider the case where each $E_n$ contains a single element for $n \geq 2$. We get rooted planar trees, not necessarily binary, so that we are back to the Hipparchus-Schroeder bracketing. An algebra over this operad is just a set with an $n$-ary operation for each $n$.

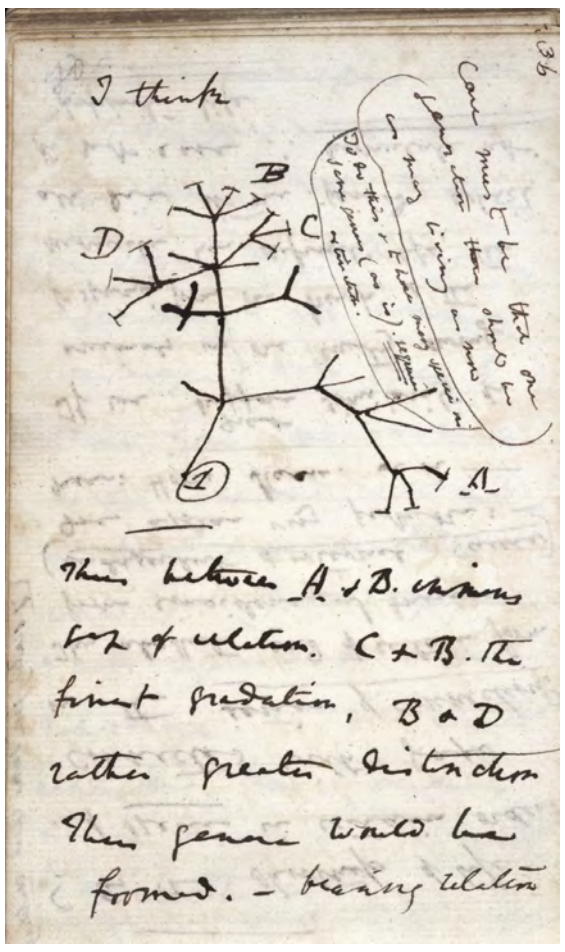

Charles Darwin's 1837 first diagram of an evolutionary tree, from one of his notebooks. ©

## *Symmetric and non symmetric*

Strictly speaking, we have dealt so far with the so-called *non-symmetric* operads. In many cases, there are actions of the symmetric groups $\Sigma_n$ on $\mathcal{O}_n$ which are compatible with the operad

operations

$$\mathcal{O}_n \times \left(\mathcal{O}_{k_1} \times \mathcal{O}_{k_2} \times \cdots \times \mathcal{O}_{k_n}\right) \to \mathcal{O}_{k_1+\cdots+k_n}.$$

In this case, the operad is called *symmetric*. Frequently, the adjective "symmetric" is omitted. I invite the reader to determine which of the previous examples are symmetric.

Find yourself the correct definition!

## *Small cubes and Stasheff again*

Recall the interpretation of the associahedra $K_n$ as spaces of metric trees, where each interior edge has some length in $[0,1]$. Grafting these trees produces maps

$$K_n \times \left(K_{k_1} \times K_{k_2} \times \cdots \times K_{k_n}\right) \to K_{k_1+\cdots+k_n}.$$

In other words, the real nature of the sequence of polytopes $K_n$ is that of an operad: *Stasheff's operad*. This is a topological operad since $\mathcal{O}_n$ is now seen as a topological space.

An algebra over $K_n$ is by definition a family of maps from $K_n \times X^n$ to $X$ which defines an operad homomorphism. Clearly, all definitions have been prepared in such a way that the *operad homomorphism* condition coincides with the definition of an $A_\infty$-space. Stasheff's theorem can therefore be restated.

**Theorem.** *A H-space $m_2 : X \times X \to X$ is homotopy equivalent to a loop space if and only if it extends as an algebra over the Stasheff operad.*

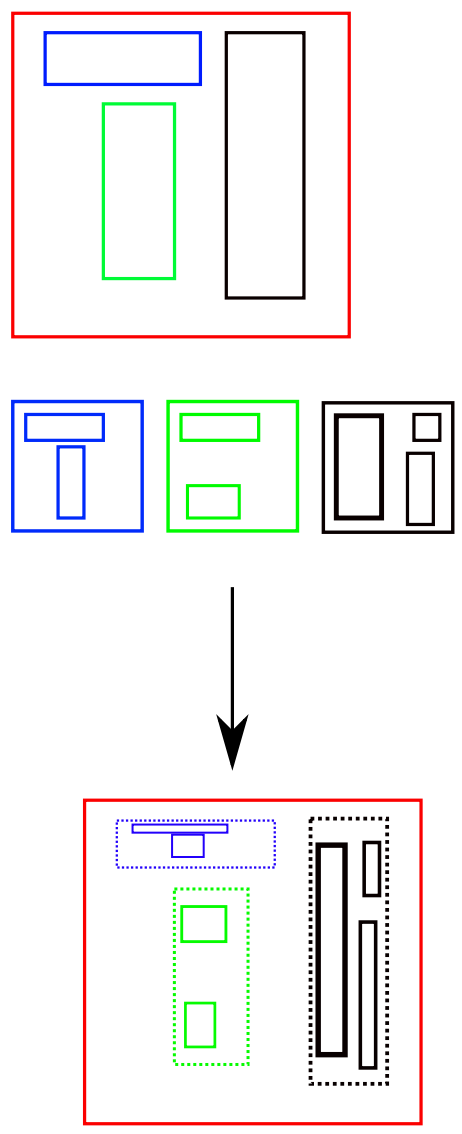

Boardman and Vogt transformed this statement still in another way and introduced the *little cubes operad*. Choose some dimension $d \geq 1$ and define a topological operad $Cub_d$ in the following manner.

- $Cub_d(n)$ is the space of $n$-tuples $(c_1, c_2, \ldots, c_n)$ of embeddings $[0,1]^d \to [0,1]^d$ such that the interiors of the images are disjoint. The $c_i$ should be affine and more precisely of the form $c_i(x_1, \ldots, x_d) = (\alpha_{1i}x_1 + \beta_{1i}, \ldots, \alpha_{di}x_d + \beta_{di})$ $(\alpha_{ij} > 0)$.
- The operad operations

  $$Cub_d(n) \times (Cub_d(k_1) \times \cdots \times Cub_d(k_n)) \to Cub_d(k_1 + \cdots + k_n)$$

  are "obvious". Simply insert the cubes as in the figure.

If $(Y, \star)$ is a topological space (in fact, as usual, a CW-complex) the $d$-loop space $\Omega^d(Y, \star)$ is the space of continuous (pointed) maps from the $d$ sphere to $Y$. We can also define $\Omega^d(Y, \star)$ as the space of maps $[0,1]^d \to Y$ which send the boundary of the cube to the base point.

The $d$-loop space of $Y$ appears naturally as an algebra over the little cubes operad. The operad operations

$$Cub_d(n) \times \Omega^d(Y, \star)^n \to \Omega^d(Y, \star)$$

are "obvious". Given $n$ little cubes $(c_1, c_2, \ldots, c_n)$ and $n$ elements $\gamma_i : [0,1]^d \to Y$ of $\Omega^d(Y, \star)$, we define a map from $[0,1]^d$ to $Y$. Inside the image of $c_i$, use the composition of $c_i$ and $\gamma_i$ and outside use the constant function sending everything to the base point.

In his review on the book by Markl, Shinder and Stasheff on operads[140], John Baez explains one of the motivations for operads.

> Most homotopy theorists would gladly sell their souls for the ability to compute the homotopy groups of an arbitrary space.

Indeed, Boardman, Vogt[141] and May[142] generalized Stasheff recognition theorem:

**Theorem.** *If a connected space $X$ is an algebra over the little cubes operad $Cub_d$, then it is homotopy equivalent to the $d$-loop space $\Omega^d(Y)$ of some space $Y$.*

## *More operads*

I believe my young reader has understood that operads occur almost everywhere in mathematics, at a foundational level. Maybe this great generality makes the theory a little bit too abstract?

Some mathematicians complain that the free group is too abstract to be a group and that it is just a bunch of words!

140 M. Markl, S. Shnider, and J. Stasheff. *Operads in algebra, topology and physics*, volume 96 of *Mathematical Surveys and Monographs*. American Mathematical Society, Providence, RI, 2002.

141 J. M. Boardman and R. M. Vogt. Homotopy-everything $H$-spaces. *Bull. Amer. Math. Soc.*, 74:1117–1122, 1968.

142 J. P. May. *The geometry of iterated loop spaces*. Springer-Verlag, Berlin-New York, 1972. Lectures Notes in Mathematics, Vol. 271.

To finish this conceptual chapter, let me give a few more examples. Enter a formula in a mathematical program, for instance the following one in *Mathematica*.

```
Sqrt[Sin[a + b + c]^2 + b^2 + c^2 + a/b]
```

$$\sqrt{\frac{a}{b} + b^2 + c^2 + \mathrm{Sin}[a+b+c]^2}$$

If you want to know how your computer "understands" this formula, just type the following.

```
TreeForm[Sqrt[Sin[a + b + c]^2 + b^2 + c^2 + a/b]]
```

You get. . . a tree.

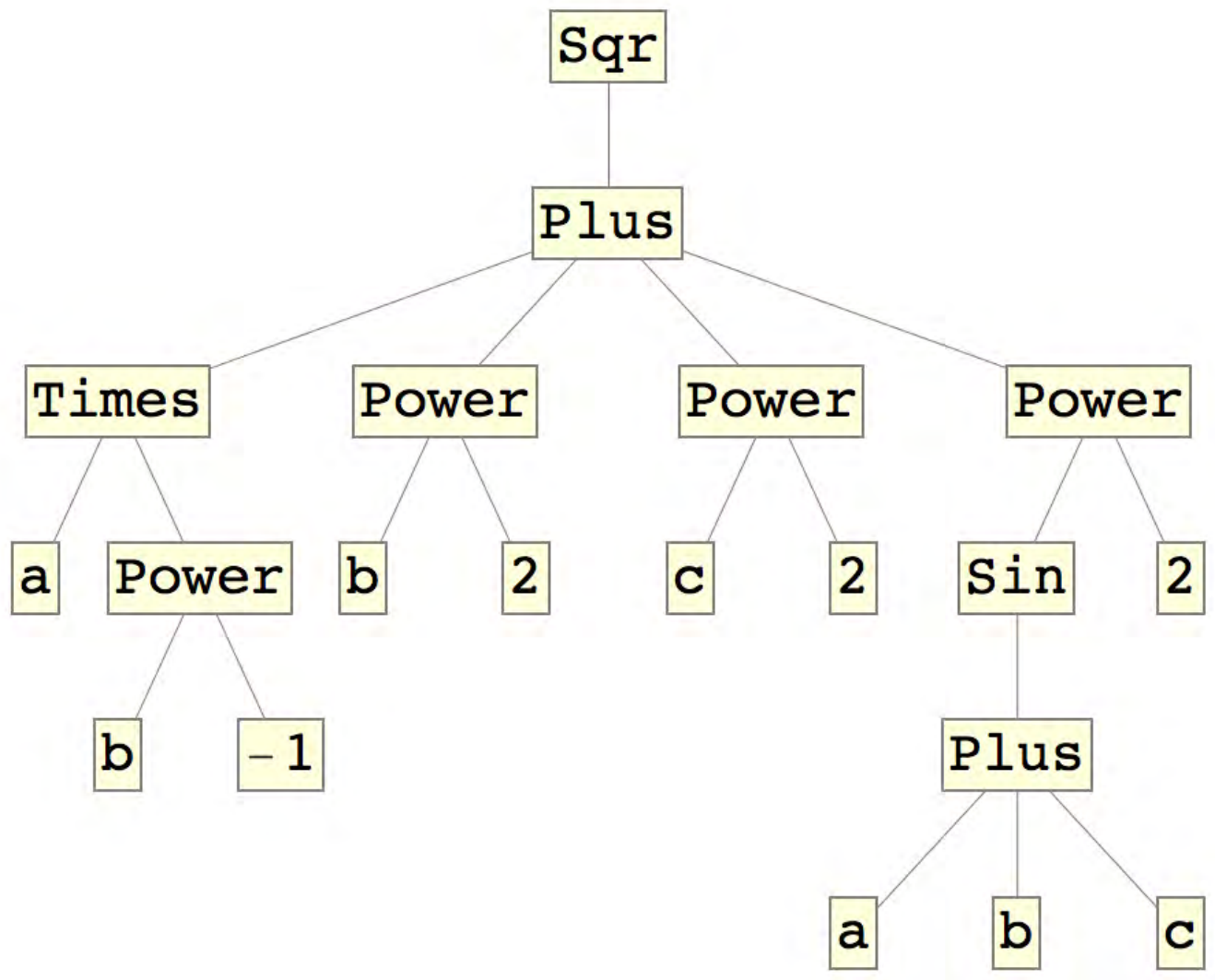

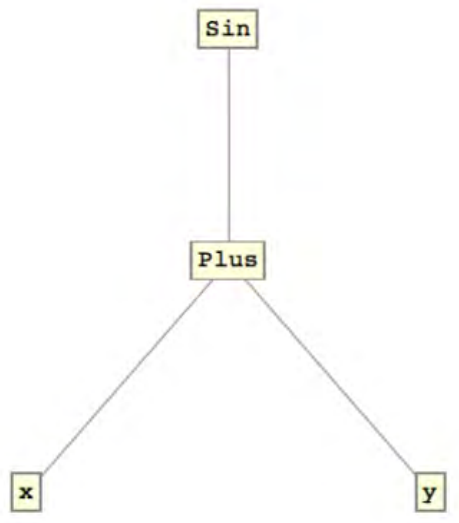

The nodes are labeled by operators, which could be $n$-ary for every $n$. The leaves are "atoms". Therefore the language of a software like *Mathematica* is actually some operad.

Note however that this operad is not free. For instance, the TreeForm of $\sin(x+y)$ is the tree in the margin.

If I substitute $a$ for $x$ and $\pi - a$ for $y$ and if I ask again for the corresponding tree:

```
TreeForm[Sin[x + y] /. x → a /. y → Pi - a]
```

I get the trivial tree with only one node labeled with 0. *Mathematica* "knows" that $\sin(\pi) = 0$ so that the tree has collapsed. In other words, the *Mathematica* operad is defined by *generators and relations*, which are built in. The user is allowed to add local rules and works in the corresponding quotient operad.

*Riemann surfaces* provide another good source of operads. Suppose that you have two Riemann surfaces $\Sigma_1, \Sigma_2$ *with boundary*. Note that a Riemann surface is canonically oriented and that this induces an orientation on the boundary. Suppose you have some orientation *reversing* diffeomorphism $f$ between some circle $S_1$ contained in the boundary $\partial\Sigma_1$ and some circle $S_2$ contained in $\partial\Sigma_2$. Glue the two surfaces along $f$ and you produced a new (oriented!) surface.

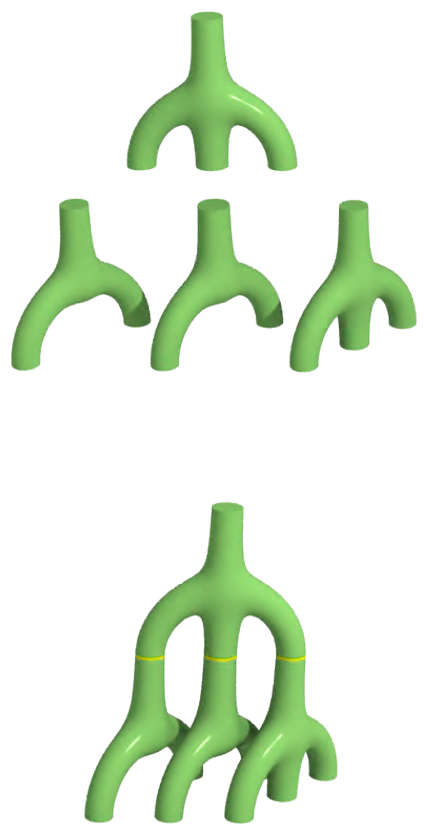

Gluing Riemann surfaces. ©

It turns out that this new surface is canonically a Riemann surface, i.e. is equipped with a structure of a 1-dimensional holomorphic manifold. This is easy to see if $f$ is real analytic since in this case it can be extended to a holomorphic diffeomorphism between small annuli, which can be used to define a holomorphic structure on the glued surface. You can also glue Riemann surfaces along non-analytic diffeomorphisms, but this is not important in our context.

Using this gluing operation defines an operad. An element of $\mathcal{O}_n$ is an isomorphism class of a compact Riemann surface with $(n+1)$ labeled boundary components, one being called entering and the $n$ others being exiting. Moreover, each boundary component is equipped with a diffeomorphism with the circle. Gluing surfaces along their boundaries, like in a *Lego* game, gives an example of an operad.

Several *functors* can be applied to operads in order to produce more operads. For instance, look at the operad $Cub_2$ of little squares. $Cub_2(n)$ has the homotopy type of the space of $n$ distinct points in a square. Its fundamental group $PB_n$ is called the *pure braid group*.

An element of $PB_n$ consists of $n$ little squares, numbered $1,2,\dots,n$ in a square (or a disc) which move along $n$ loops without intersecting each other. At the end of the loops, the squares came back to their initial positions and this last property is what is meant by *pure*.

The definition of the fundamental group requires a base point. Therefore, to define $PB_n$, one should choose some $n$-tuple of distinct points in the square, as "starting point" for our braids. I suggest that the reader finds himself these base points in order to define properly the operad maps.

Using fundamental groups, we get maps

$$PB_n \times \left(PB_{k_1} \times \cdots \times PB_{k_n}\right) \to PB_{k_1+\cdots+k_n}$$

and we get a group operad. This is not very complicated. An element of $PB_n$ gives rise to $n$ tubes in a cylinder $[0,1]^2 \times [0,1]$. The operad structure consists in inserting tubes into tubes.

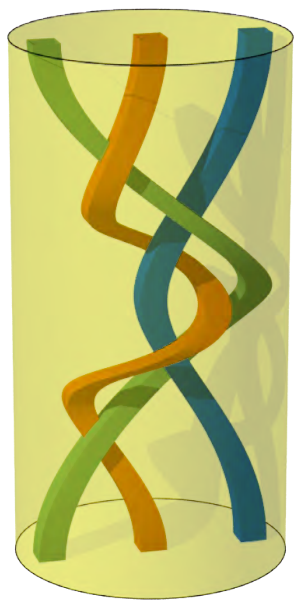

Three small squares follow a loop. ©

There is no reason to limit ourselves to the 2-dimensional case and to the fundamental group. Given a topological space $X$ and an integer $n$, the so-called *configuration space* $X^{[n]}$ is defined as the space of $n$-tuples of *distinct* points in $X$. If $X$ is a 2-dimensional disc, it turns out that the universal cover of this space is contractible so that, from the homotopy point of view, only its fundamental group is interesting: this is the pure braid group that we just introduced. However, if $X$ is a higher dimensional ball, $X^{[n]}$ is simply connected and we are tempted to describe its topology.

For instance if $X = \mathbb{R}^3$ and $n = 2$, the answer is easy: two distinct points $x, y$ in space are completely defined by their midpoint $(x+y)/2$ and the nonzero vector $x-y$. Therefore $(\mathbb{R}^3)^{[2]}$ has the same homotopy type as a 2-sphere. The situation is already more complicated for $n = 3$: three bodies in space...

A good approach is to study the homology or cohomology of these spaces not individually, for each value of $n$, but globally: the cohomology of the little cube operad.

I urge my reader to read the very accessible paper[143] which can serve as an entrance gate to operad theory. You will learn for instance that "The homology of the little $d$-cubes operad is the degree $d$ *Poisson* operad" (whatever that means).

[143] D. P. Sinha. The (non-equivariant) homology of the little disks operad. In *OPERADS 2009*, volume 26 of *Sémin. Congr.*, pages 253–279. Soc. Math. France, Paris, 2013.

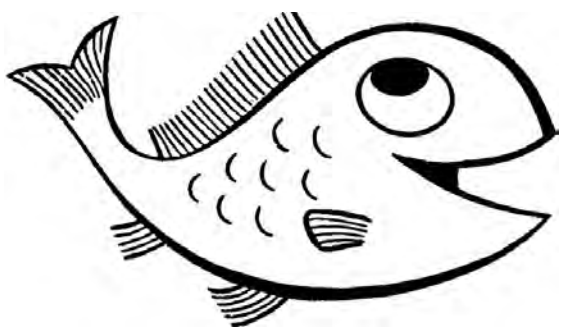

The Poisson operad. ©

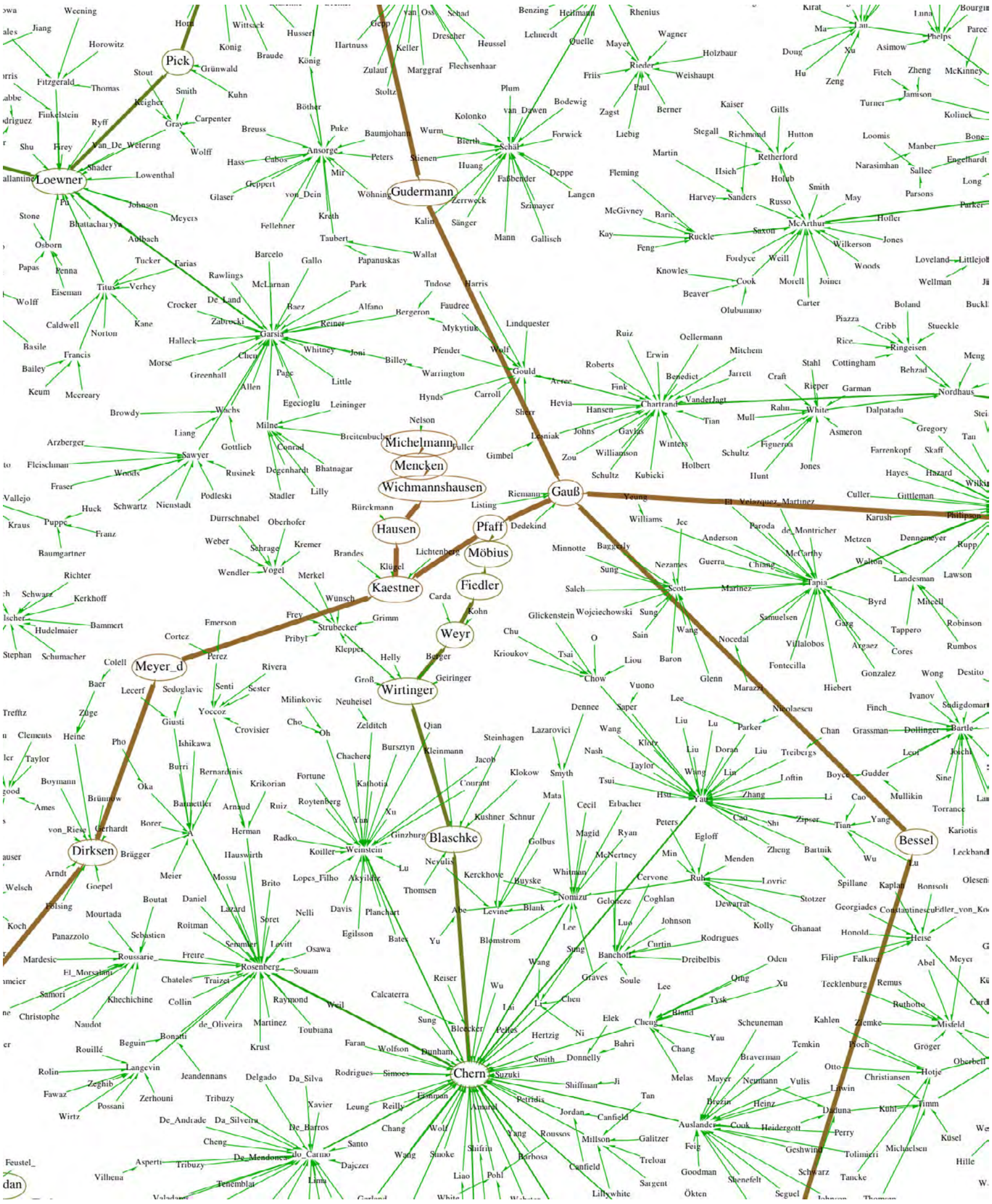

This is a local view, centered on Gauss, of the "tree of mathematicians" where one connects two mathematicians if one was the advisor of the other. ©

# *Singular operads*

We come back to our initial discussion on the relative position of the graphs of a family of polynomials.

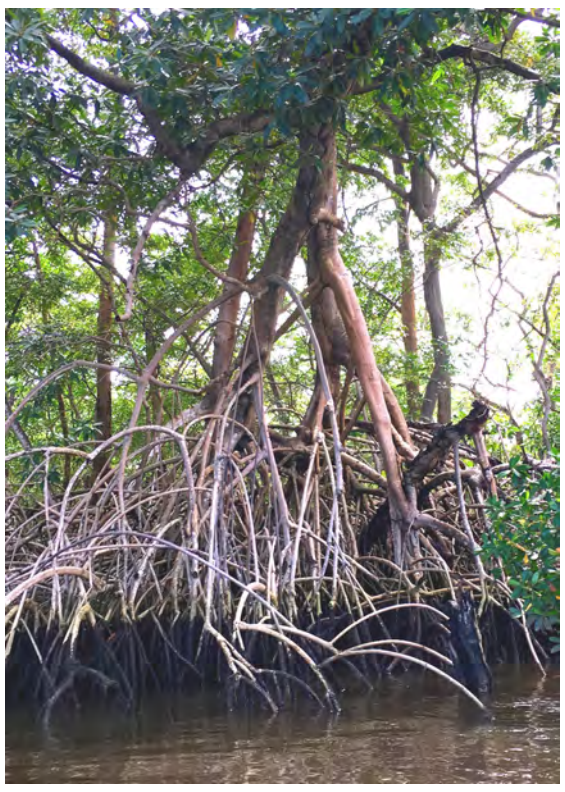

Roots or branches?
(Rio Preguiças, Brazil).

## *The real polynomials operad*

Let $(P_1, P_2, \ldots, P_n)$ and $(Q_1, Q_2, \ldots, Q_n)$ be two $n$-tuples of distinct polynomials in $\mathbb{R}[x]$ vanishing at the origin. I will say that they are *topologically equivalent* if for small nonzero values of $x$, the two $n$-tuples of real numbers $(P_1(x), \ldots, P_n(x))$ and $(Q_1(x), \ldots, Q_n(x))$ are ordered in the same way. Let $Pol_{\mathbb{R}}(n)$ be the (finite) set of equivalence classes of such $n$-tuples. Let us construct a very simple operad structure on $Pol_{\mathbb{R}}(n)$.

Suppose that we are given

- $(P_1, P_2, \ldots, P_n)$ (a representative of) an element of $Pol_{\mathbb{R}}(n)$,
- for each $i = 1, \ldots, n$, an element of $Pol_{\mathbb{R}}(k_i)$ given by (the class of) $(P_{i;1}, \ldots, P_{i,k_i})$.

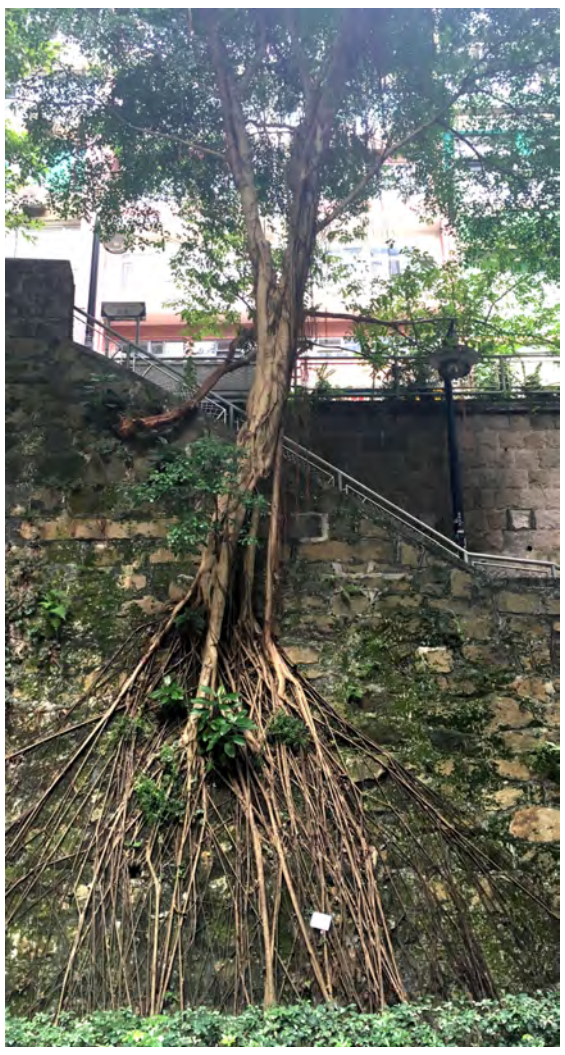

Roots and branches.
(Hong Kong).

We want to *graft* the $P_{i;j}$'s on the $P_i$'s. Just consider the $k_1 + \cdots + k_n$ polynomials

$$P_i(x) + x^{2N} P_{i;j}(x) \quad (1 \leq i \leq n \text{ and } 1 \leq j \leq k_i)$$

in the lexicographical order of $(i, j)$, where $N$ is some large integer.

Some explanations may be useful. The role of $x^{2N}$ is to make sure that the terms which are added to the $P_i$'s are much smaller than the differences $P_i - P_j$ $(i \neq j)$. To ensure this property, it suffices to choose $2N$ bigger than all $v(P_i - P_j)$ $(i \neq j)$. Fixing $i$, the graphs of the $k_i$ polynomials $P_i(x) + x^{2N}P_{i;j}(x)$ are very close to the graph of $P_i$. The *even* exponent $2N$ implies that, fixing $i$, the order between the $P_i(x) + x^{2N}P_{i;j}(x)$'s is the same as the order between the $P_{i;j}(x)$'s. Topologically, the graphs of the $P_i$'s have been transformed into some thin wedges in which the $P_{i;j}$'s have been inserted.

Recall that $v(P)$ denotes the valuation of $P$.

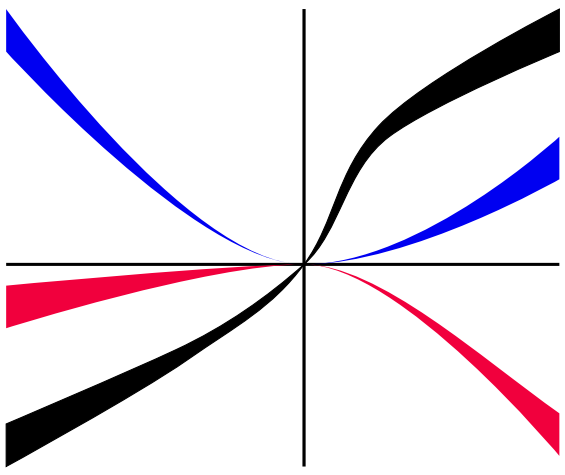

It should be clear that this is well defined and gives an operad structure on the $Pol_{\mathbb{R}}(n)$'s. This is a symmetric operad as the polynomials $(P_1, P_2, \ldots, P_n)$ can be permuted.

It should be equally clear, from the earliest chapters of this book, that this operad is very close to the (non-symmetric) operad of separable permutations.

Can you find an explicit relationship between the two operads?

## *The complex polynomials operad*

Let us play a similar game with *complex* polynomials. If $(P_1, P_2, \ldots, P_n)$ is some $n$-tuple of distinct polynomials in $\mathbb{C}[x]$, vanishing at the origin, we look at the following loops in $\mathbb{C}$ (for $1 \leq i \leq n$):

$$\gamma_i : \theta \in \mathbb{R}/2\pi\mathbb{Z} \mapsto P_i(\varepsilon \exp(\sqrt{-1}\theta)) \in \mathbb{C}.$$

Choose $\varepsilon$ very small. Then for every $\theta$ the $n$ points $\gamma_1(\theta), \gamma_2(\theta), \ldots$ are distinct. This defines a loop in the space of $n$-tuples of distinct points in the plane, that is to say an element of the *pure braid group* $PB_n$. To be precise, I should speak of a conjugacy class of a pure braid, since the initial points $\gamma_i(0)$ could be anywhere and the definition of the pure braid group requires some base point. This conjugacy class is independent of the choice of the small $\varepsilon$.

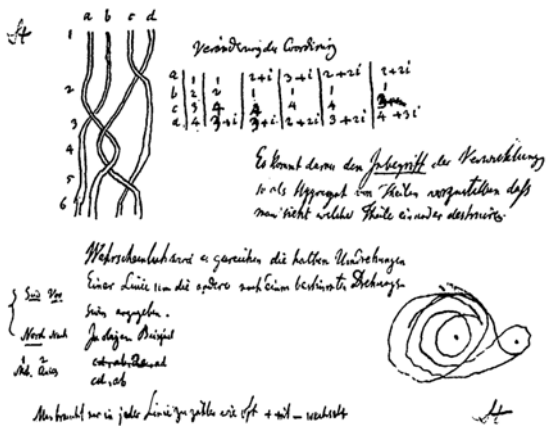

An unpublished manuscript of Gauss in which he starts the topological study of braids.

Say that $(P_1, P_2, \ldots, P_n)$ and $(Q_1, Q_2, \ldots, Q_n)$ are *topologically equivalent* if the corresponding braids are conjugate in $PB_n$. Denote by $Pol_{\mathbb{C}}(n)$ the set of equivalence classes.

Is this equivalent to the existence of a local homeomorphism of $\mathbb{C}^2$ in the neighborhood of the origin sending each complex curve $P_i = 0$ to the complex curve $Q_i = 0$?

Recall that the pure braid groups define one of our examples of operads, where the operations consist in inserting braids in tubular neighborhoods of the strands of a given braid. This suggests that $Pol_{\mathbb{C}}(n)$ could be a sub-operad of the $PB_n$. This is indeed the case as I show now.

We know that $\exp(-v(P_i - P_j))$ defines an ultrametric distance on $\{P_1, \dots, P_n\}$ which can be encoded by a rooted tree. The root corresponds to the full set $\{P_1, P_2, \dots, P_n\}$. Its children are the equivalence classes of the relation $v(P_i - P_j) \geq 2$ etc. Until we reach the singletons $\{P_i\}$ which are the leaves, labeled by $1, 2, \dots, n$.

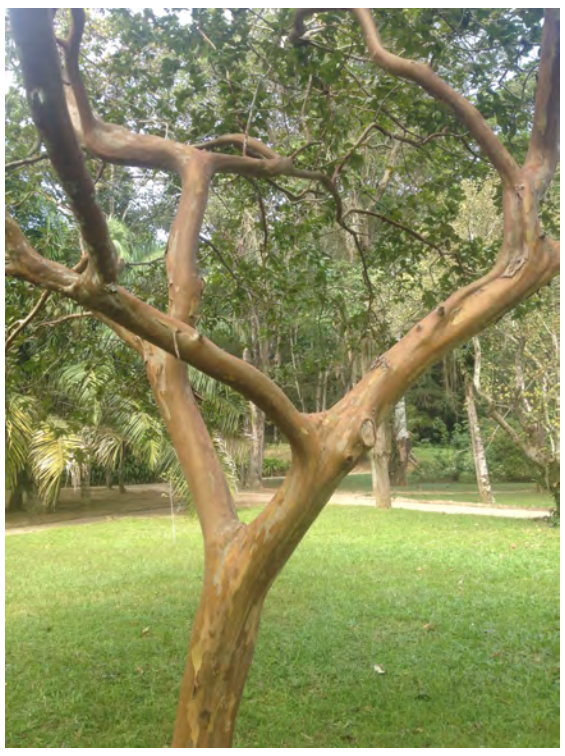

A non-planar tree in the Jardim Botânico, Rio de Janeiro.

There are two main differences with the case of real polynomials.

– There is no natural order structure on the nodes, so that our tree is not planar (after all, most trees in nature are not planar).

– In the real case, we did some pruning on the tree. We did that since for instance the pairs $(0, x)$ and $(0, x^3)$ are topologically equivalent *over the reals*: $x$ and $x^3$ have the same signs. But this is not true anymore in the complex domain: the braid associated to $(0, x^3)$ rotates 3 times unlike $(0, x)$ which rotates only once as $x$ describes the boundary of a small disc centered on 0.

The following (elementary) result gives a precise description of $Pol_{\mathbb{C}}(n)$.

**Theorem.** *Two $n$-tuples of complex polynomials $(P_1, P_2, \dots, P_n)$ and $(Q_1, Q_2, \dots, Q_n)$ are topologically equivalent if and only if they define the same rooted tree, or equivalently if $v(P_i - P_j) = v(Q_i - Q_j)$ for all $i, j$.*

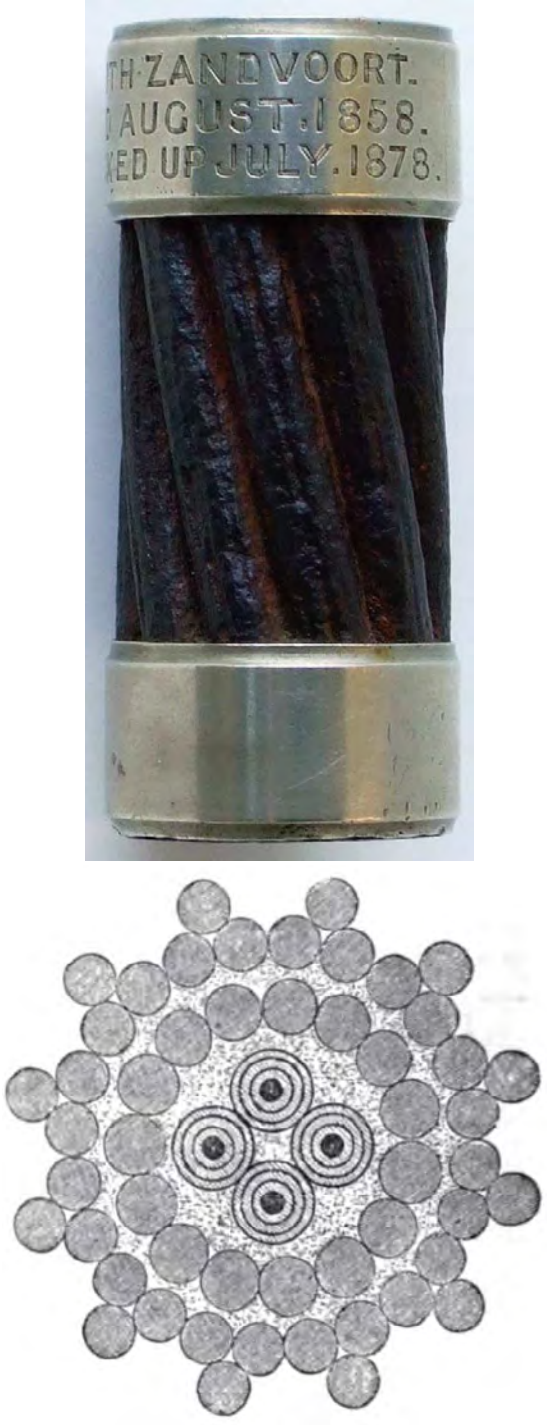

A 1858 England-Holland submarine telegraphic cable and its cross-section.

Let us first show that *the tree determines the braid*. Start from the root and descend $p$ edges until you reach the first node with $q \geq 2$ children. This means that the $p$-th Taylor polynomials of all $P_i$'s are all equal and that there are $q$ different $(p+1)$-st Taylor polynomials.

If $t_p$ denotes this common $p$-th Taylor polynomial, we subtract $t_p$ from all $P_i$'s without changing the corresponding braid. Therefore, $(P_1, P_2, \dots, P_n)$ comes in $q$ groups, where we place in the same group two $P_i$'s with the same Taylor polynomial of order $(p+1)$, of the form $a_\nu x^{p+1}$, for $\nu = 1, \dots, q$, where all the $a_\nu$'s are distinct complex numbers. When $x$ describes the circle of radius $\varepsilon$ these $q$ points $a_\nu x^{p+1}$ describe small circles under the same rigid rotation and they rotate by $p+1$ full turns. This (conjugacy class of this) braid in $PB_q$ only depends on $q$ and $p$ and not of the position of the $a_\nu$'s: we have $q$ points which simply rotate in block and describe $p+1$ full turns.

Around each $a_\nu \varepsilon^{p+1} \exp(\sqrt{-1}(p+1)\theta)$ draw a small disc which contains all the $P_i(\varepsilon \exp(\sqrt{-1}\theta))$ with $(p+1)$-st Taylor polynomial $a_\nu x^{p+1}$. Continue the process inside each of these discs, splitting again each group according to higher order Taylor polynomials.

The conclusion is that the braid is like a Solar system, moving along epicycles, à la Hipparchus. It consists of a group of $q$ small discs rigidly rotating by $p+1$ turns. Inside each disc, the picture is similar. And so on, until we arrive at the $\gamma_i(\theta)$'s. All these numbers $p,q$, for all these discs, are given by the combinatorics of the tree, so that the braid (always up to conjugacy) is indeed determined by the tree, in a very concrete way.

Encore lui!

I explain now how to construct the tree *from the braid*.

Choose two integers $1 \le i < j \le n$. A pure braid in $PB_n$ is a homotopy class of a loop of $n$-distinct points $(x_1, \ldots, x_n)$ in the plane. Forget about all points except $x_i$ and $x_j$: this defines some homomorphism from $PB_n$ to $PB_2$, for each $i,j$. The structure of $PB_2$ is very simple: it is isomorphic to $\mathbb{Z}$. When two distinct points move and come back to their original position, the vector $x_i - x_j$ is a loop in the punctured plane and has an index relative to 0. This is the isomorphism between $PB_2$ and $\mathbb{Z}$. This defines $n(n-1)/2$ homomorphisms $lk_{i,j} : PB_n \to \mathbb{Z}$, which simply express the number of turns of $x_i - x_j$. Two conjugate braids have the same images by $lk_{i,j}$.

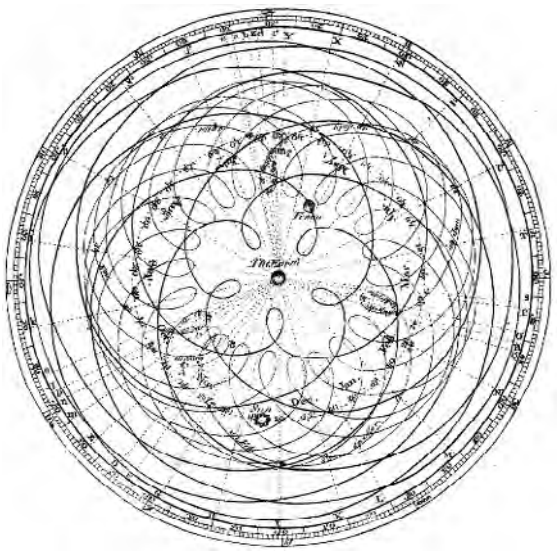

The apparent motion of the Sun, Mercury, and Venus from the Earth. ©

Let us come back to our braid defined by the $\gamma_i$'s. Evaluating $lk_{i,j}$ on this braid, we are counting the number of turns of $P_i(x) - P_j(x)$ when $x$ goes around the boundary of a small disc centered at the origin. This is obviously the valuation $v(P_i - P_j)$. Hence, the valuations $v(P_i - P_j)$ can be read off from the conjugacy class of the braid, as desired. ⊡

Conversely, given a rooted tree with $n$ leaves, it is easy to construct $n$ polynomials whose associated tree is the given one.

Our trees are non-planar. However, rooted non-planar trees can be grafted as soon as their their leaves are labeled from 1 to $n$. Given a rooted tree $T$ with $n$ leaves, labeled $1, \ldots, n$, and $n$ trees $T_1, \ldots, T_n$, we can graft $T_i$ on the leaf numbered $i$ of $T$. This

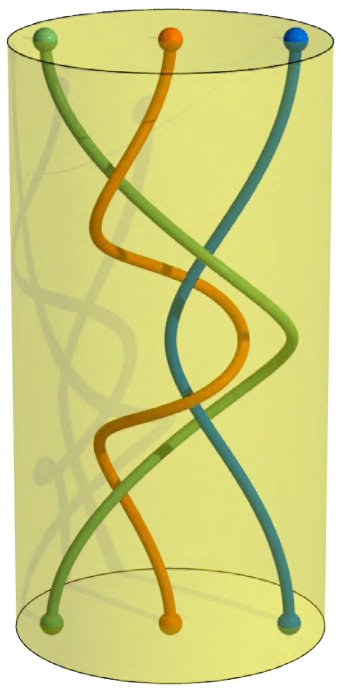

A pure braid with three strands. ©

defines the operad of *labeled rooted non-planar trees*.

Therefore, this grafting operation defines a natural operad structure on $Pol_{\mathbb{C}}(n)$. Moreover, the previous discussion shows that $Pol_{\mathbb{C}}(n)$ is a sub-operad of the pure braids operad.

This is a *symmetric* operad, since the polynomials can be permuted.

We can describe this structure in formulas, directly in terms of polynomials. Given some $n$-tuple $(P_1, \dots, P_n)$, define

$$\delta_i = \max_{j \neq i} v(P_i - P_j).$$

$\delta_i$ is the level at which the leaf $P_i$ is attached to the tree.

The operad action of $(P_1, \dots, P_n)$ on a family $(P_{i;j})$ $(1 \leq i \leq n$ and $1 \leq j \leq k_i)$ is defined as the $(k_1 + \dots + k_n)$-tuple of polynomials (in lexicographic order)

$$P_i(x) + x^{\delta_i} P_{i;j}(x) \ \ (1 \leq i \leq n \ \text{ and } \ 1 \leq j \leq k_i).$$

Check that this definition does realize the grafting of the associated labeled trees, as desired.

*In summary $Pol_{\mathbb{C}}(n)$ is isomorphic to the operad of labeled rooted non-planar trees and appears as a sub-operad of the pure-braid operad $PB_n$.*

## *An operad associated to complex singular curves*

For two complex polynomials $P_i(x)$ and $P_j(x)$ vanishing at 0, the valuation $v(P_i - P_j)$ is also called the *multiplicity* of intersection of the two smooth curves $y = P_i(x)$ and $y = P_j(x)$ at the origin. There is no surprise in this terminology since this is indeed the multiplicity, in the usual sense of the word, of the root 0 of $(P_j - P_i)(x) = 0$. The previous paragraph can also be re-interpreted in the following way.

The last two chapters of this book deal with the linking number in more detail.

The curves $y - P_i(x) = 0$ are smooth in $\mathbb{C}^2$. They intersect transversally a small sphere $\mathbb{S}^3_\varepsilon$ on a trivial knot. We looked at these knots in the square sphere $\max(|x|, |y|) = \varepsilon$ and we denoted them by $\gamma_i$. It turns out that the *linking number* of $\gamma_i$ and $\gamma_j$ is nothing more than their multiplicity of intersection. This is a topological interpretation of the multiplicity.

The linking number of two oriented knots $k_1, k_2$ in $\mathbb{S}^3$ is defined in the following manner. Choose an embedded oriented surface whose oriented boundary is $k_1$ and count the (algebraic) intersection number of this surface with $k_2$.

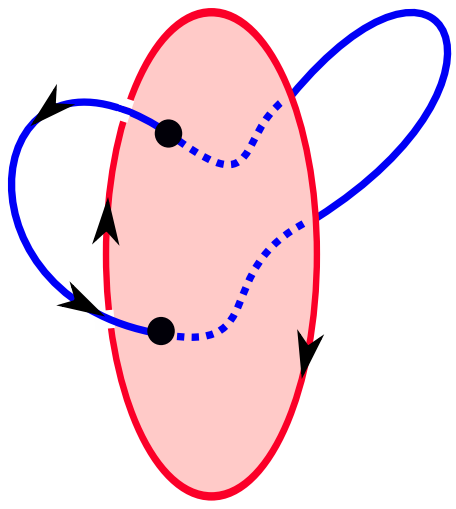

Linking number $= 1 - 1 = 0$.

In our simple case, consider the link of the smooth curve $y - P_i(x) = 0$ in the 3-sphere as the boundary of one of its Milnor fibers, where $y - P_i(x)$ is real and positive. In order to compute the linking number of the two curves $y - P_i(x) = 0$ and $y - P_j(x) = 0$ in the 3-sphere, we have to count the (algebraic) intersection of the knot $y - P_j(x) = 0$ with the Milnor fiber $y - P_i(x) \in \mathbb{R}_+$. Now $(P_j - P_i)(x) = ax^{v(P_i - P_j)} + \cdots$, so that the second knot intersects the Milnor fiber exactly $v(P_i - P_j)$ times (one should check that the intersections are positive). It follows that *the linking number between $\gamma_i$ and $\gamma_j$ is indeed the multiplicity $v(P_i - P_j)$*.

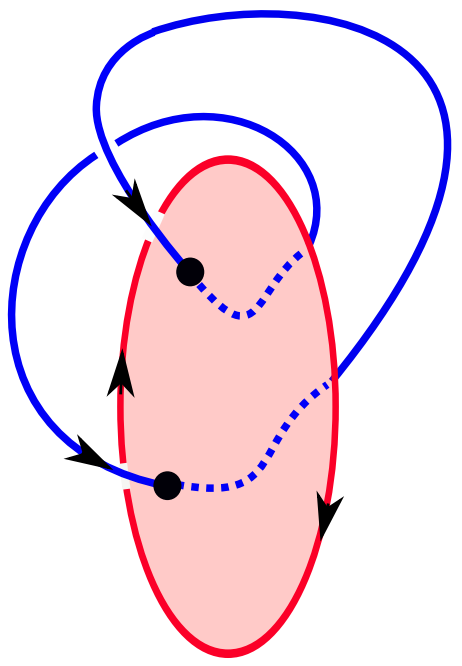

Linking number = 1 + 1 = 2.

It turns out that most of what has been said is true for branches of non-smooth curves. Let $F(x,y) = 0$ be a singular complex analytic curve, assumed reduced, admitting $n$ branches. Write $F = F_1 \cdots F_n$, a decomposition in irreducible factors. Each branch $F_i = 0$ admits a Puiseux parameterization:

$$t \in \mathbb{C} \mapsto (t^{m_i}, g_i(t)).$$

The integer $m_i$ is the *order* of the branch. The intersection of each branch with a small sphere is a knot $k_i$, which is not trivial if the branch is not smooth. The linking number $m_{ij}$ of $k_i$ and $k_j$ is the *multiplicity of intersection of the two branches*. The reader can view this as a definition, if his mind is topologically oriented. If she prefers algebra, she could proceed in the following way. Insert the parameterization of a branch in the equation of the another one and look at the multiplicity of the zero $t = 0$ of

$$t \in \mathbb{C} \mapsto F_j(t^{m_i}, g_i(t)).$$

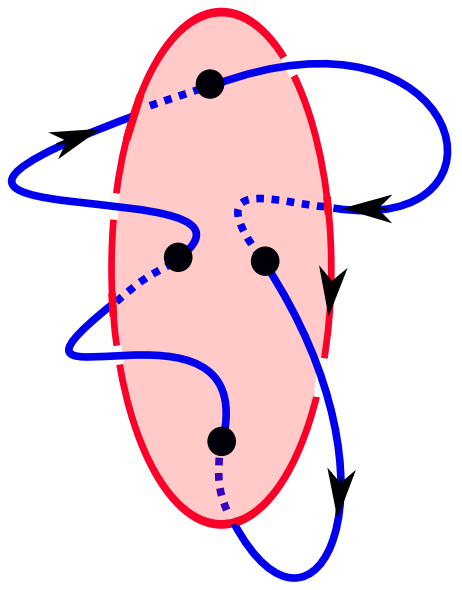

Linking number = 4.

Note the analogy of the two definitions: "algebraic multiplicity of intersection of two branches" and "linking number" of two knots. Both look asymmetric but are in fact symmetric.

It turns out that the multiplicities of intersection $m_{ij}$ between the branches always have the important properties that we noticed for smooth curves.

- $m_{ij}$ is a positive integer. This is easy.

Can you prove that a branch defines a trivial knot if and only if the branch is smooth?

- The $m_{ij}/m_i m_j$'s satisfy some sort of *ultrametric inequality*. In other words, for every $\varepsilon > 0$, the relation $m_{ij}/m_i m_j \geq \varepsilon$ is an equivalence relation in $\{1,\dots,n\}$ so that we can construct a tree, as before (except that the length of the edges are not integers but rational numbers).

The second item is due to Płoski[144] in 1985. See[145] for a modern presentation and application and [146] for an ample generalization.

We can now construct an *operad of complex singular curves*. Actually, as before, we need a labeling of the branches of our curve $y = P_1(x),\dots,y = P_n(x)$ (where the $P_i$'s are now the Puiseux series $g_i(x^{\frac{1}{m_i}})$). In other words, we consider the set of $n$-tuples of distinct Puiseux series, denoted $Curv_n$, as *complex curves with labeled branches*. We now define an operad structure on the $Curv_n$'s exactly as we did with polynomials. Define positive *rational numbers* $\delta_1,\dots,\delta_n$, associated to $(P_1,\dots,P_n) \in Curv_n$ by

$$\delta_i = \max_{j \neq i} v(P_i - P_j)$$

where $v$ is now the lowest rational exponent of a non-trivial term of a Puiseux series.

Given elements of $Curv_{k_1},\dots,Curv_{k_n}$, defined by the Puiseux series $P_{i;j}$ for $1 \le i \le n$ and $1 \le j \le k_i$, we can construct the element of $Curv_{k_1+\dots+k_n}$ as the $(k_1+\dots+k_n)$-tuple

$$P_i(x) + x^{\delta_i} P_{i;j}(x) \quad (1 \le i \le n \text{ and } 1 \le j \le k_i)$$

in lexicographic order. This defines a (symmetric) operad.

*Exercise*:

The intersection of $F(x,y) = 0$ with $|x| = \varepsilon$ defines a braid, which is not a pure braid anymore. Indeed each branch intersects $x = \varepsilon$ a number of times equal to its own order $m_i$. Unlike the previous case, the knot defined by each branch is not trivial.

Find a topological version of the previous operad. In other words, find a natural equivalence relation between curves $F(x,y) = 0$ (with labeled branches) in terms of the associated braids that they define. This equivalence relation should be such that the operad actually acts on equivalence classes.

[144] A. Płoski. Remarque sur la multiplicité d'intersection des branches planes. *Bull. Polish Acad. Sci. Math.*, 33(11-12):601–605 (1985), 1985.

[145] I. Abío, M. Alberich-Carramiñana, and V. González-Alonso. The ultrametric space of plane branches. *Comm. Algebra*, 39(11):4206–4220, 2011.

[146] E. García Barroso, P. González Pérez, and P. Popescu-Pampu. Ultra-metric spaces of branches on arborescent singularities. 2016.

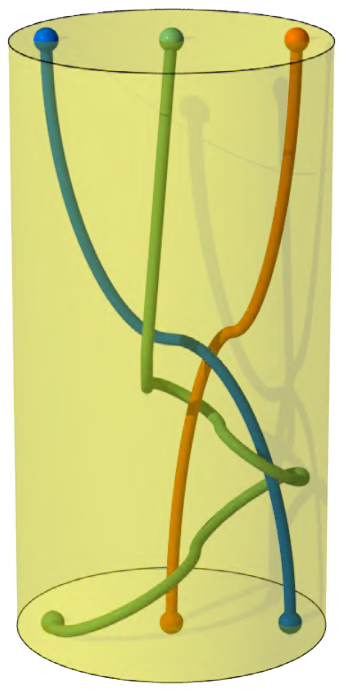

An "unpure" braid with 3 strands. ©

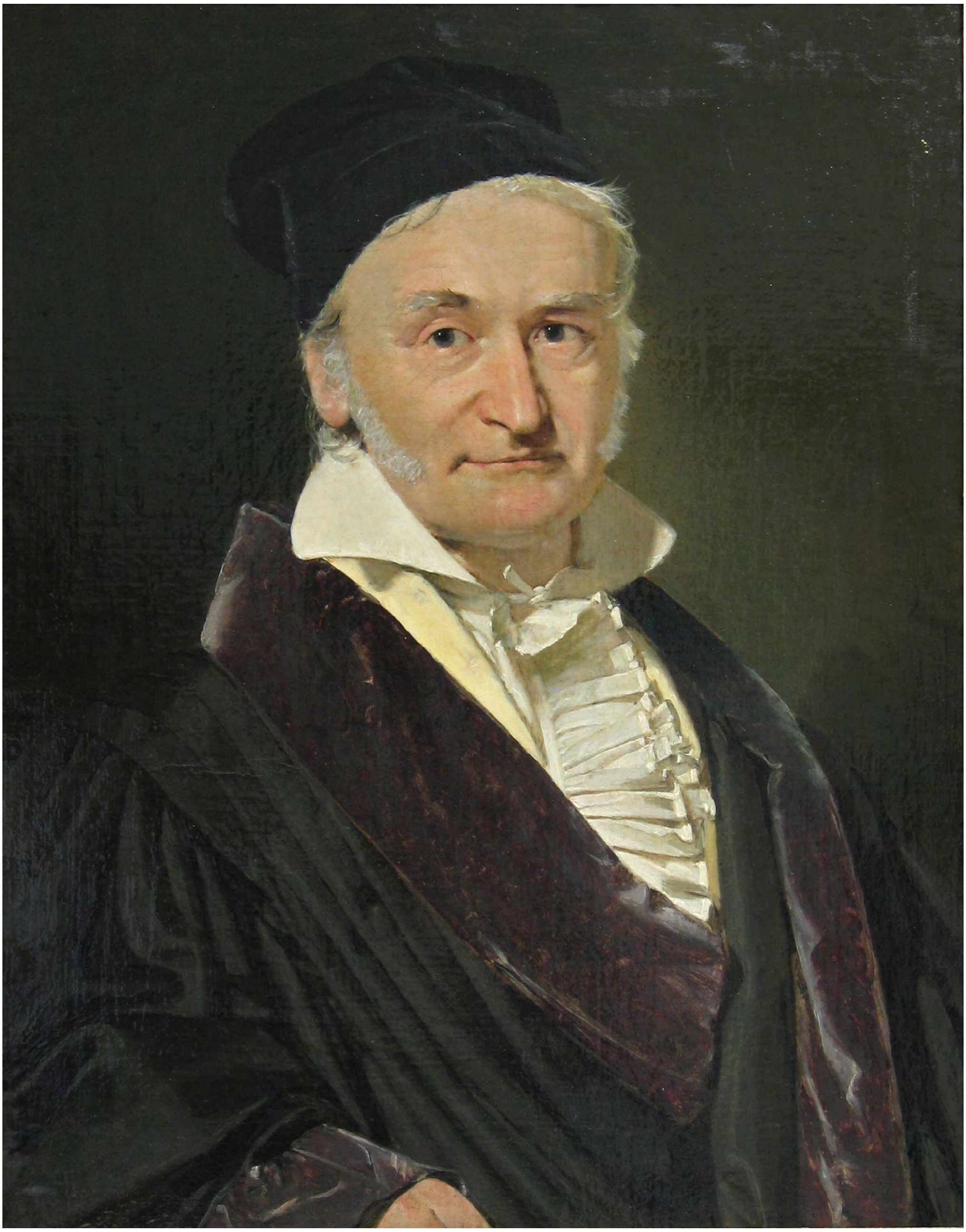

Gauss in 1840. ©

# *Gauss is back: curves in the plane*

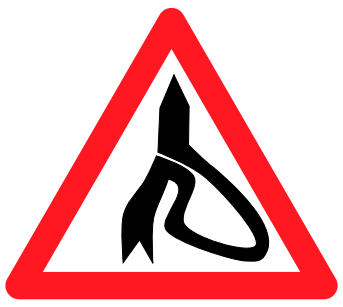

Detour! I like this chapter, but it is completely independent from the rest of the book.

Many great mathematicians, past or present, have enjoyed, or enjoy, drawing curves. As a quizz, I enclose some pictures and my reader should guess their authors.

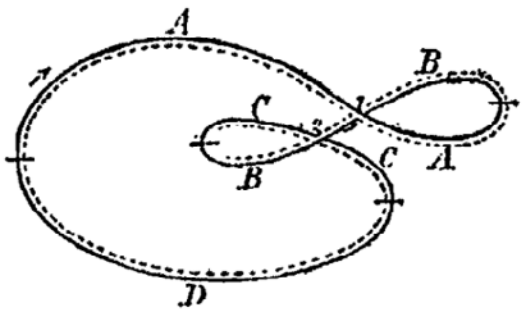

1

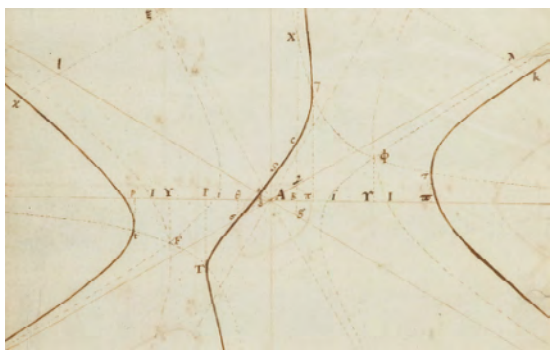

2

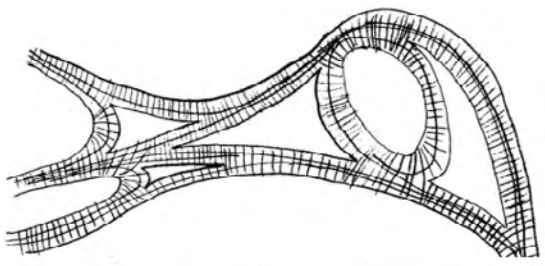

3

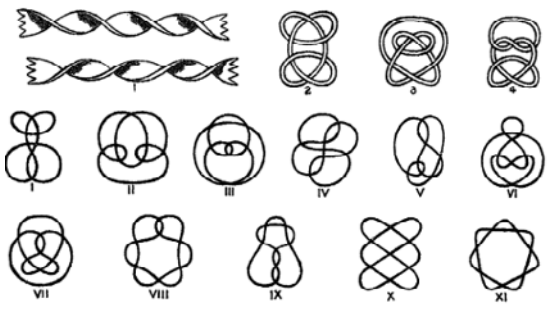

4

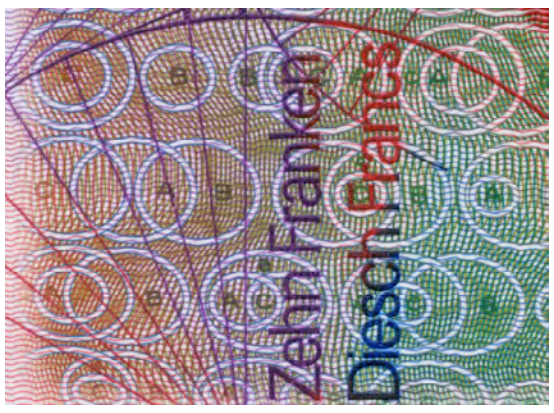

5

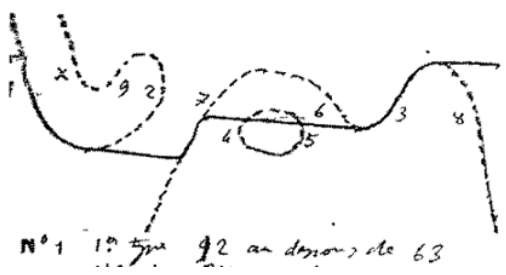

6

Solution: next page.

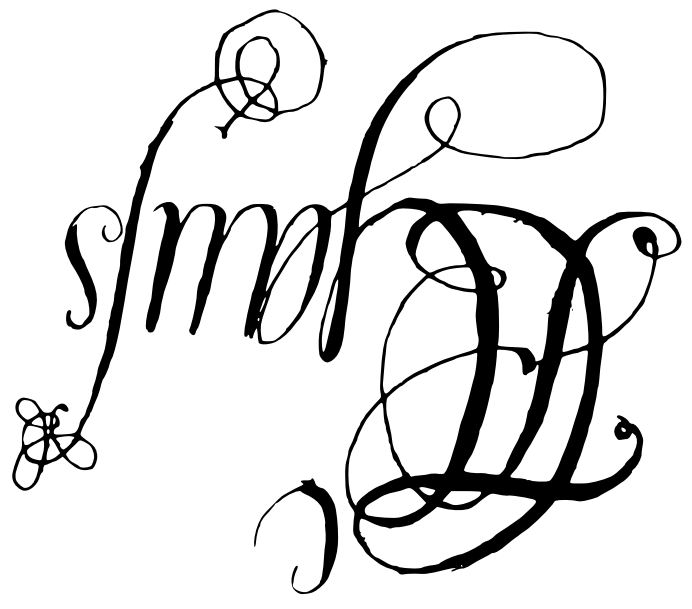

Gauss's signature in 1794 (he was 17).

**Picture 1** was drawn by *Gauss* in his Nachlass, around the problem which is discussed in this chapter.

**Picture 2** is a cubic curve drawn by *Newton* in his *Enumaratio Linearum tertii Ordinis* (1667). In modern terminology, Newton tries to classify these algebraic curves up to projective transformations. He speaks of the "generation of curves by shadows", which is a wonderful definition of projective geometry.

**Picture 3** is a drawing by *Thurston* in his *Notes on the topology and geometry of 3-manifolds* (1980) which has been very influential. This picture represents a *train-track*, an important concept which is useful in the description of the boundary of Teichmüller space.

**Picture 4** is a knot projection by *Tait* (1884) in a series of papers which represent the first serious attempt to classify knots.

**Picture 5** is from a Swiss banknote in honor of *Euler*. In his *Lettres à une princesse d'Allemagne*, he explained *Boolean* logic to a German young princess, using diagrams bounded by ellipses. He figured out all the possible configurations of intersections of four ellipses and the many possible cases are in the background of the banknote.

**Picture 6** is by *Poincaré* in an attempt to prove "his last theorem", on the number of fixed points of an area preserving diffeomorphism of the annulus. Actually, he did not prove it as he passed away a few months later. The theorem was proved soon after Poincaré's death by Birkhoff.

## *Gauss words*

In this chapter, we will review a beautiful question raised by Gauss about curves in the plane. Recall that his first proof of the fundamental theorem of algebra is based on the qualitative behavior of curves inside a disc.

Volume 8 of *Gauss's Werke* contains a few pages[147] (page 272 and 282-286) on immersed curves in the plane. We should be careful: these are the so-called *Nachlasse*, notes which remained unpublished during Gauss's life. Do not forget his motto *Pauca sed matura* (Few, but ripe). These pages should not be considered as an actual publication but more like a private draft.

[147] C. F. Gauß. *Werke. Band VIII*. Georg Olms Verlag, Hildesheim, 1973. Reprint of the 1900 original.

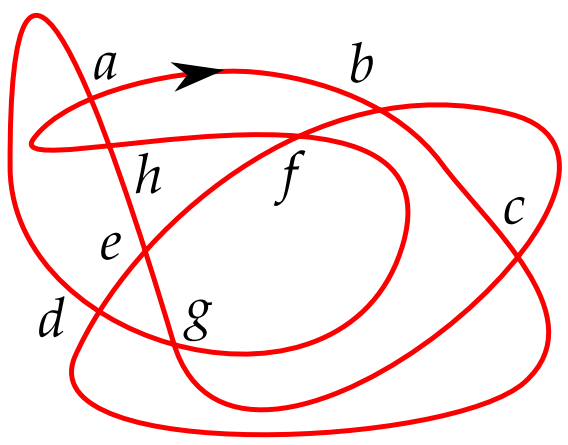

*abcdefbcgehadgfh*

Draw a closed generic oriented immersed curve $i : S^1 \to \mathbb{R}^2$ in the plane. Generic means that multiple points are only double points with two different tangents. Label the double points by $n$ symbols $a_1, \ldots, a_n$ (that Gauss called the *Knoten*). Going around the curve, we pass twice in each of these $a_i$'s and this defines a cyclic word of length $2n$ in which each of the $a_i$'s appears twice. *The closed curve therefore defines a chord diagram* with one chord for each $a_i$.

The question raised by Gauss is to recognize which chord diagrams arise from some planar curve. He first lists all possibilities for $n \leq 8$, by hand! Then he finds a *necessary* condition. Writing the word on a circle, between two occurrences of some $a_i$, there should be an *even* number of letters. In the example in the margin, the chord $d$ decomposes the circle in two intervals containing respectively 6 and 8 letters.

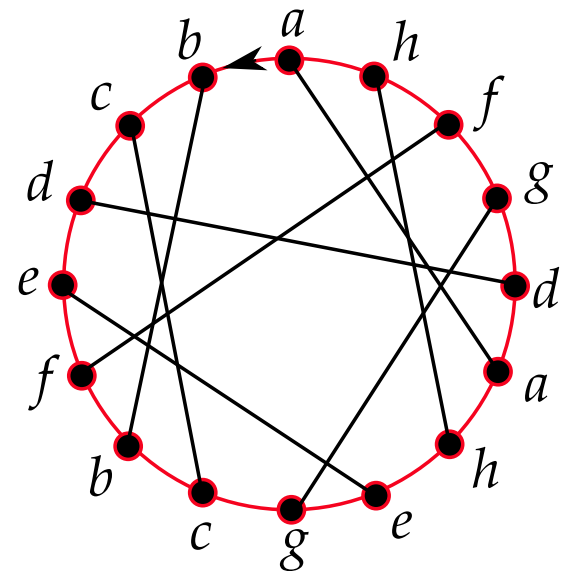

Some modern authors claim that this was a conjecture of Gauss and that he could not prove it. What a lack of respect! It seems clear to me that Gauss could prove it and did not take time to write it down in his private notebook.

One of the first theorems in topology, known to Gauss in his PhD, as we have already noticed, is that two closed curves in the plane intersecting transversally have an *even number of intersection points*. One of the possible proofs is to move the first curve by a generic path of translations so that at the end of the motion, there is no more intersection point. We then examine how intersection points appear or disappear in generic situations. It is not difficult to see that generically, points appear or disappear in pairs.

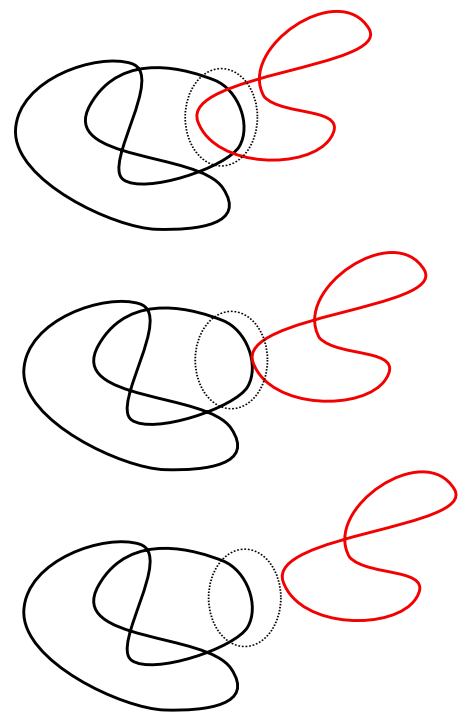

A closely related fact has been noticed by all pupils drawing doodles, during boring math classes. If you draw a generic closed curve in the plane, the connected components of the complement can be colored in black and white like in a checker board. Just use white for the component at infinity and for some other component connect it to infinity by some generic arc and color white if this arc intersects the closed curve in an even number of points, and black otherwise. This is coherent because of the previous observation.

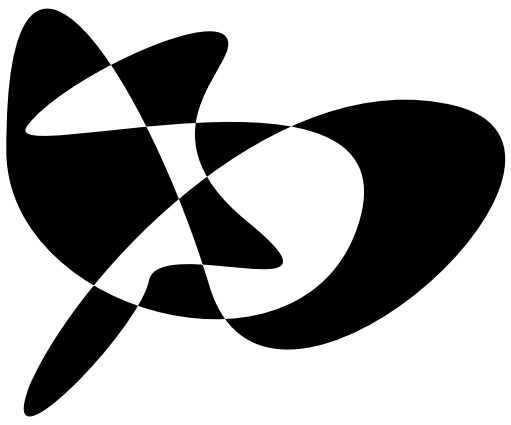

Now, Gauss's necessary condition is an easy corollary. Any chord decomposes the circle into two intervals $I_1, I_2$ which define two closed curves in the plane, say $\gamma_1$ and $\gamma_2$, starting from the same point. Slightly move these two curves to make them transverse. The number of letters in $I_1$ is equal to the number of intersection points between $\gamma_1$ and $\gamma_2$ plus twice the number of self-intersection of $\gamma_1$. Therefore it is even. On the picture in the margin, one of the two loops from $b$ to $b$ is slightly shifted and shown as a dotted blue loop.

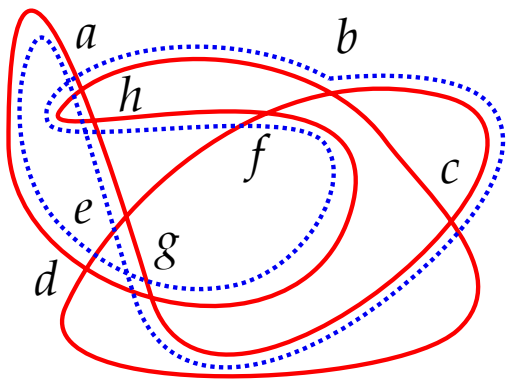

This necessary condition is not sufficient as Gauss knew very well.

## *Signed Gauss words*

Gauss's problem has been solved many times, in many different ways, in different mathematical communities, basically topological or combinatorial. This is in tune with Poincaré's quote:

> [...] il n'y a plus des problèmes résolus et d'autres qui ne le sont pas, il y a seulement des problèmes plus ou moins résolus.

I will only describe two solutions.

Let me present first a solution of a simpler problem, using a topological argument, mixing [148], [149] and [150]. If the plane is oriented, each double point of our generic curve defines two tangent vectors, so that one of them is "the first" and the other is "the second". Going around the curve, as we meet a double point, we mark a + sign if we pass through the first branch and a − sign otherwise. Hence, the cyclic word is now decorated with signs, or exponents, ±1. Each letter occurs twice, with different signs. Equivalently, we could think of an *oriented chord diagram* where each chord goes, say, from its + end to its − end. The *signed Gauss's problem* (that Gauss did *not* study) is the following. Given such a *signed word*, can one decide whether it is associated to some generic curve in the plane?

[148] N. Chaves and C. Weber. Plombages de rubans et problème des mots de Gauss. *Exposition. Math.*, 12(1):53–77, 1994.

[149] G. Cairns and D. M. Elton. The planarity problem for signed Gauss words. *J. Knot Theory Ramifications*, 2(4):359–367, 1993.

[150] G. Cairns and D. M. Elton. The planarity problem. II. *J. Knot Theory Ramifications*, 5(2):137–144, 1996.

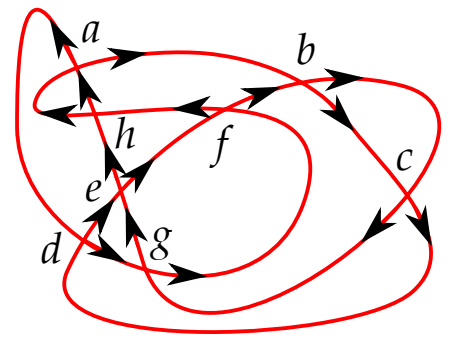

$a^+b^+c^-d^-e^+f^+b^-c^+g^-e^-h^+a^-d^+g^+f^-h^-$

For each symbol $a_1, a_2, \ldots, a_n$, draw a cross, as in the margin. Each cross has two entering sides and two exiting sides. Each cross has two arms, labeled + and −.

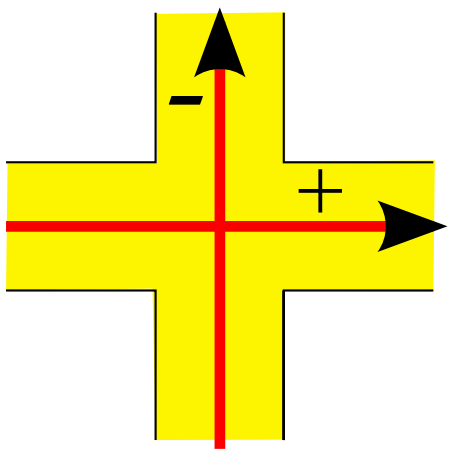

A signed cyclic word $w$ defines uniquely a way of glueing each exit side of each cross to some entrance side of some other cross. The result of this gluing operation is some oriented surface $S$ with boundary containing an immersed oriented curve $\gamma$. Going around this curve, we read precisely the signed word $w$. If $S$ is planar so that it can be embedded in the plane, our problem is solved since we constructed a curve in the plane. Conversely, if the word comes from some immersed curve in the plane, some neighborhood of its image is clearly a union of crosses, assembled as in $S$.

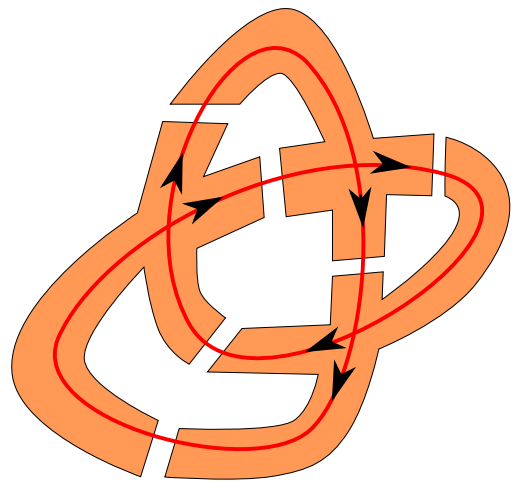

*Therefore $w$ is realizable by some immersed planar curve if and only if $S$ is planar.*

For the rest of this chapter, some familiarity with the basic theory of the homology of surfaces is necessary. This is a good opportunity to recommend the visual book by Fomenko[151]. Opening this book, my reader will immediately understand why I like it. A more standard book by Massey is very accessible[152].

Let $k$ be the number of boundary components of $S$. Let $\hat{S}$ be the surface obtained from $S$ by gluing a disc to each of its boundary components. The surface $S$ has the homotopy type of a graph with $n$ vertices and $2n$ edges. Hence, the Euler-Poincaré characteristic of $\hat{S}$ is $n - 2n + k$. Compact oriented surfaces without boundary are classified by their Euler-Poincaré characteristic. As a corollary, $S$ is planar if and only if $k - n = 2$.

This gives a very simple algorithm to decide if $w$ is realizable.

[151] A. Fomenko. *Visual geometry and topology*. Springer-Verlag, Berlin, 1994. Translated from the Russian by Marianna V. Tsaplina.

[152] W. S. Massey. *A basic course in algebraic topology*, volume 127 of *Graduate Texts in Mathematics*. Springer-Verlag, New York, 1991.

Start with $w$, glue the crosses, and count the number of components of the boundary: it should be $n+2$. This is essentially due to Carter[153].

[153] J. S. Carter. Classifying immersed curves. *Proc. Amer. Math. Soc.*, 111(1):281–287, 1991.

I now present another point of view. It follows from the classification of compact oriented surfaces with boundary that such a surface is planar if and only if any two closed transverse curves intersect in an even number of points. Indeed, as soon as the genus of a surface is $\geq 1$, it contains a punctured torus which contains two curves intersecting exactly in one point. Hence, to check whether the genus of $S$ is 0, it suffices to find a basis of its homology $H_1(S;\mathbb{Z}/2\mathbb{Z})$ modulo 2, and to compute the intersection.

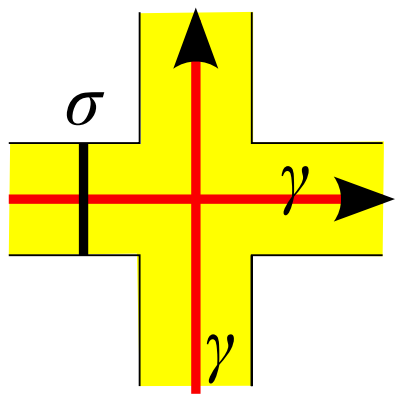

There is an easy way to find a basis of $H_1(S;\mathbb{Z}/2\mathbb{Z})$. As a preliminary observation, note that an ordered pair of points $(a^+, a^-)$ on the oriented circle defines an open interval, as we travel in the positive direction from $a^+$ to $a^-$. I will say that the elements of this interval are *between* $a^+$ *and* $a^-$. Be careful however: this interval is changed into its complement if the two points are switched.

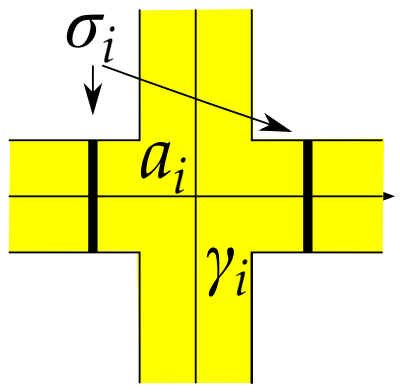

The original curve $\gamma$ is drawn on $S$ and therefore defines a homology class $[\gamma]$. Moreover, for each $i = 1, \ldots, n$, the interval from $a_i^+$ to $a_i^-$ on the circle defines a loop $\gamma_i$ on $S$ and a homology class $[\gamma_i]$ in $H_1(S;\mathbb{Z}/2\mathbb{Z})$. Note that when $\gamma_i$ enters a cross, it does not change direction with the exception of the cross labeled $a_i$, where it turns right. Said differently, the intersection of $\gamma_i$ with a cross different from the one labeled $a_i$ is either empty, or a straight segment, or two perpendicular segments.

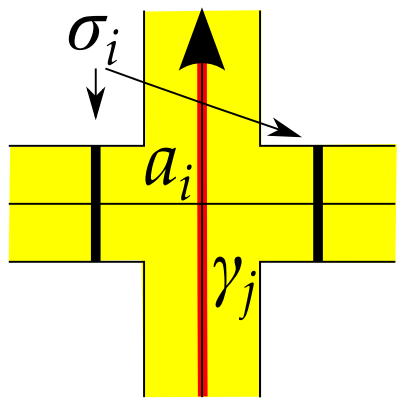

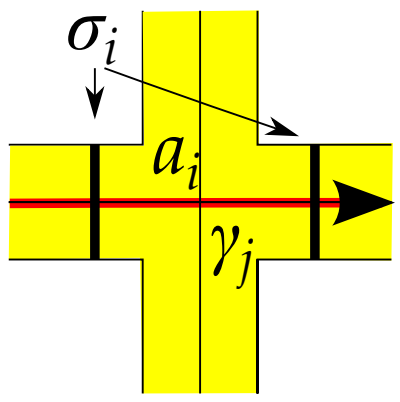

**Lemma.** *The classes* $[\gamma], [\gamma_i]_{1\leq i\leq n}$ *define a basis of* $H_1(S;\mathbb{Z}/2\mathbb{Z})$.

The surface $S$ has the homotopy type of a connected graph with $n$ vertices and $2n$ edges. The Euler-Poincaré characteristic is $-n$ and is equal to 1 minus the rank of $H_1(S;\mathbb{Z}/2\mathbb{Z})$. Therefore this rank is $(n+1)$. In order to prove the lemma we just have to show that the $[\gamma], [\gamma_i]_{1\leq i\leq n}$'s are linearly independent in $H_1(S;\mathbb{Z}/2\mathbb{Z})$.

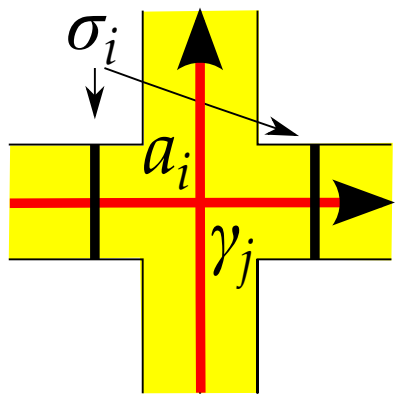

Any arc $\sigma$ in $S$ with endpoints in the boundary of $S$ defines a linear form in $H_1(S;\mathbb{Z}/2\mathbb{Z})$: just count intersection points with $\sigma$ (always modulo 2). For instance, choose $\sigma$ in some cross, as in the margin, in the previous page. Clearly $\gamma$ intersects $\sigma$ in only one point so that $[\gamma]$ is not trivial. In the cross associated to the letter $a_i$, let $\sigma_i$ be the sum of the two curves shown in the margin. The intersection of $\sigma_i$ with $\gamma_j$ is 0 if $i \neq j$ and 1 if $i = j$. The intersection of $\gamma$ and $\sigma_i$ is 0. Therefore one can use the $n+1$ linear forms $\sigma_i$ and $\sigma$ to show that $[\gamma], [\gamma_i]_{1\leq i \leq n}$ are linearly independent and the lemma is proved. ⊡

Now, the surface $S$ has genus zero if and only if the intersection numbers of the loops $\gamma, \gamma_i$ are all 0, modulo 2.

Observe that the self intersection of any closed curve on any orientable surface is 0 (over $\mathbb{Z}$ and hence over $\mathbb{Z}/2\mathbb{Z}$).

Why ?

Let us compute the intersection number of $\gamma$ and $\gamma_i$, denoted by $\gamma \cdot \gamma_i$. To make them transverse, move $\gamma$ slightly to its right in order to get a curve $\gamma'$ which is transversal to $\gamma_i$. Do not forget that the surface $S$ and the loops $\gamma_i$ are oriented. It follows that the intersection number of $\gamma$ and $\gamma_i$ is the number of letters between $a_i^+$ and $a_i^-$. We recover Gauss's necessary condition.

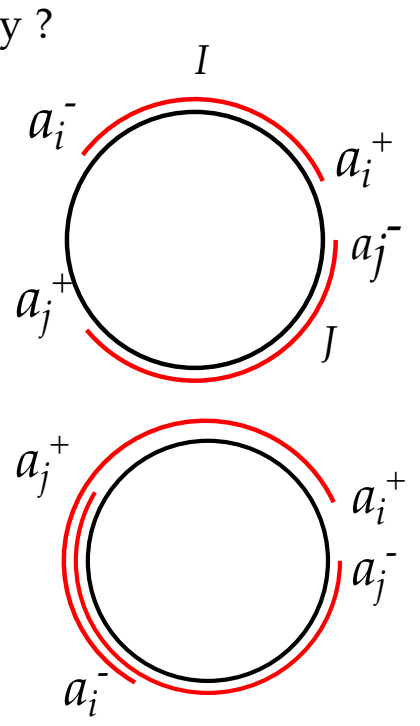

*Assume from now on that it is satisfied.*

Let us compute the intersection number modulo 2 of $\gamma_i$ and $\gamma_j$, denoted by $\gamma_i \cdot \gamma_j$.

*If the letters $a_i^\pm, a_j^\pm$ are* not *linked*, there are two *disjoint* intervals $I, J$ in the circle (or in the cyclic word) whose endpoints are $a_i^+, a_i^-$ and $a_j^+, a_j^-$ respectively. Since $\gamma_i$ can be replaced by $\gamma_i - \gamma$ in the computation of the intersection number, it follows that, in this unlinked case, $\gamma_i . \gamma_j$ is the number of letters in the word with one occurrence in $I$ and the second in $J$. This is therefore a second parity condition, necessary for the realizability of $w$.

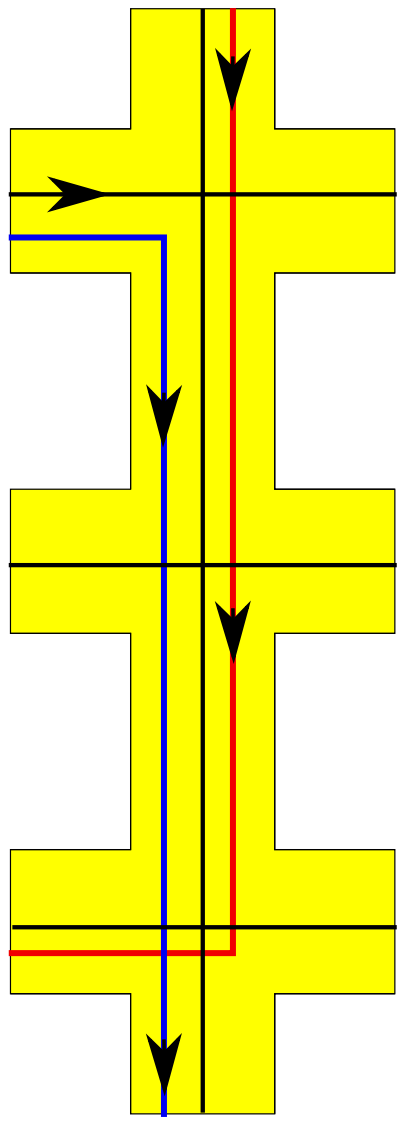

*If the letters $a_i^\pm, a_j^\pm$* are *linked*, the loops $\gamma_i$ and $\gamma_j$ are not transversal since they coincide on some non-trivial interval. Move $\gamma_i$ slightly to the right, to produce some $\gamma_i'$, and move $\gamma_j$ to the left, to get some $\gamma_j'$. The curves $\gamma_i'$ and$\gamma_j'$ are now parallel on this common part. We now count the intersection number of $\gamma_i'$ and $\gamma_j'$. Look at the picture.

There is one intersection in the cross $a_j$ and none in the cross $a_i$. The other intersections correspond to letters between $a_i^+$ and $a_i^-$ whose second occurrence is between $a_j^+$ and $a_j^-$.

Hence, *when $a_i^\pm, a_j^\pm$* are *linked, the intersection number $\gamma_i . \gamma_j$ is*

*equal to* 1 *plus the number letters between* $a_i^+$ *and* $a_i^-$ *whose other occurrence is between* $a_j^+$ *and* $a_j^-$.

We therefore get a very simple answer to the *signed* Gauss's problem.

**Theorem.** *A signed Gauss word is realizable by a planar immersed curve if and only if the following conditions are satisfied.*

1. *For every letter* $a_i$*, there is an even number of letters between* $a_i^+$ *and* $a_i^-$ *(Gauss's parity condition).*

2. *For every* $i, j$ *such that the letters* $a_i^\pm$ *and* $a_j^\pm$ *are not linked, let* $I, J$ *be the disjoint intervals whose endpoints are* $a_i^+, a_i^-$ *and* $a_j^+, a_j^-$ *(excluding these points). The number of letters in the word with one occurrence in* $I$ *and the other in* $J$ *is even.*

3. *For every* $i, j$ *such that the letters* $a_i^\pm$ *and* $a_j^\pm$ *are linked, the number of letters between* $a_i^+$ *and* $a_i^-$ *whose other occurrence is between* $a_j^+$ *and* $a_j^-$ *is odd.*

Equivalently, the first two necessary conditions can be expressed in terms of the chord diagram.

1/ Every chord is intersected by an even number of chords.

2/ Given two non-intersecting chords, the number of chords intersecting both of them is even.

## *Gauss's problem*

Let us come back to the original problem: *non-signed* words. Of course, we could cheat and try all the $2^n$ ways of choosing signs on the word. That might take a terribly long time. Even Gauss's computational force could have been beaten by $2^n$. Moreover this would not be very enlightening.

Note that Gauss's parity criterion is independent of the signs. The second condition, in the case where $a_i^\pm$ and $a_j^\pm$ are not linked, is also independent of the signs.

*We therefore assume that they are both satisfied for a non-signed word* $w$*.*

Let me introduce the so-called *interlace graph* $G(w)$. Its vertices are the integers $1, \ldots, n$, or the chords $a_i$, and there is an edge between two chords if they intersect.

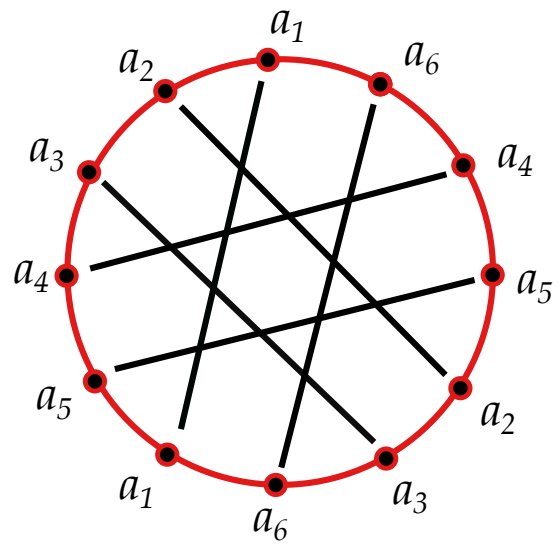

Here is an example from Cairns and Elton's paper:

$$w = a_1 a_2 a_3 a_4 a_5 a_1 a_6 a_3 a_2 a_5 a_4 a_6.$$

One checks easily that Gauss's parity condition is satisfied as well as the second condition in the non-linked case.

The interlace graph is represented in the margin.

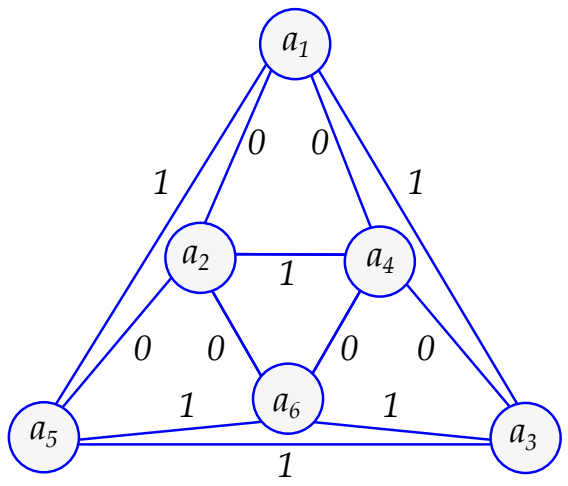

Choose some signed word $\overline{w}$ whose unsigned associated word is $w$, for instance

$$\overline{w} = a_1^+ a_2^+ a_3^+ a_4^+ a_5^+ a_1^- a_6^+ a_3^- a_2^- a_5^- a_4^- a_6^- .$$

Since we are back to signed words, we know how to compute the intersection numbers $\gamma_i.\gamma_j$ of the associated loops. This defines an element $f(e) \in \mathbb{Z}/2\mathbb{Z}$ for each edge $e$ of $G(w)$. For instance, the two chords 1 and 3 intersect. On the margin picture, there is no letter in the interval from $a_1^+$ to $a_1^-$ (in dotted green) whose other occurrence is between $a_3^+$ and $a_3^-$ (in dotted blue). Therefore $\gamma_1 \cdot \gamma_3$ is equal to $0 + 1$ modulo 2 and we write a 1 on the edge connecting $a_1$ and $a_3$ in the graph. Since we see 1 on some edges in our example, this *signed word* is not realizable by an immersed curve in the plane.

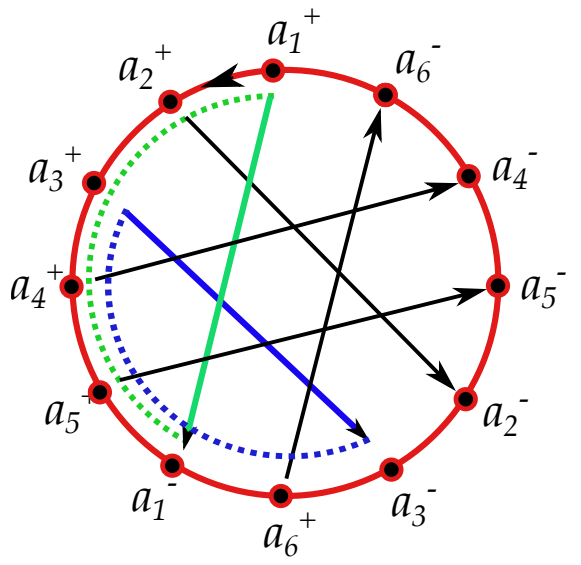

We now have to decide if there could exist a clever change in the signs in the letters so that all edges could be labeled by 0.

This labeling $f(e)$ can be thought as a 1-cochain $f$ (with values in $\mathbb{Z}/2\mathbb{Z}$) in the graph, hence a 1-cocycle (since there are no 2-faces in a graph). Let us examine how this cocycle changes when the signs on $\overline{w}$ are modified. A change of signs is defined by some function $u$ from $\{1, \ldots, n\} \to \mathbb{Z}/2\mathbb{Z}$, that we can think of as a 0-cochain in our graph. I claim that the new 1-cocycle $f'$ on $G(w)$, after the sign change associated to $u$, is simply equal to $f + du$, where $du$ is the *coboundary* of $u$. This $du$, evaluated on some edge $e$ is by definition the difference (or sum since we work modulo 2) of the values of $u$ at the two endpoints of $e$.

Let us begin by changing only the sign of *one letter*, say $a_k$. For every $i, j$ such that the letters $a_i^\pm$ and $a_j^\pm$ are linked, in other words for each edge of the graph $G(w)$, we have to compare two intersection numbers, for the two signed words $\overline{w}, \overline{w}'$ whose signs only differ on the letter $k$. Clearly, these intersection numbers are equal if $k$ is different from $i$ and $j$. It turns out that they differ by 1 (modulo 2), when $k = i$ or $k = j$. For instance, if $k = i$, we have to compare:

1/ the number of letters between $a_i^+$ and $a_i^-$ whose other occurrence is between $a_j^+$ and $a_j^-$, and

2/ the number of letters between $a_i^-$ and $a_i^+$ whose other occurrence is between $a_j^+$ and $a_j^-$.

Modulo 2, this difference is the number of letters between $a_j^+$ and $a_j^-$, *different from* $a_i$. This number is odd since we assumed that Gauss parity condition and $a_i^\pm$ and $a_j^\pm$ are linked. Hence, the effect of changing signs on a single letter $a_k$ is to change by 1 the labels on the edges of $G(w)$ with $a_k$ as an endpoint, and not to modify the other labels. So the formula $f' = f + du$ holds in this simple case. The general case follows since we can change signs one by one.

As a consequence, the object which is well defined, independently of the signs, is the *cohomology class* of $f$ in $H^1(G(w); \mathbb{Z}/2\mathbb{Z})$. This class is zero if and only if the cocycle is zero for some choice of the signs, i.e. if and only if the unsigned word $w$ is realizable by an immersed generic curve in the plane.

Finally, a cohomology class in a graph is trivial if and only if it is zero when evaluated on cycles. This gives a very efficient algorithm. Choose any signed word $\overline{w}$, compute the 1-cocycle, and sum its values on cycles in the interlace graph.

In our example, the cycle $a_1 \to a_3 \to a_5 \to a_1$ in $G(w)$ gives a total sum 1 so that the *unsigned* word $w$ is not realizable by an immersed curve in the plane.

See[154] for a history of the problem until 1972 and[155] for a more recent book.

[154] B. Grünbaum. *Arrangements and spreads*. American Mathematical Society Providence, R.I., 1972. Conference Board of the Mathematical Sciences Regional Conference Series in Mathematics, No. 10.

[155] C. Godsil and G. Royle. *Algebraic graph theory*, volume 207 of *Graduate Texts in Mathematics*. Springer-Verlag, New York, 2001.

## *The genus of a chord diagram*

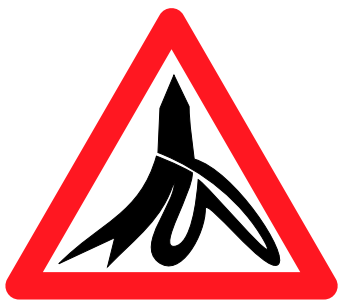

Detour in the detour.

Let $w$ be a diagram with $n$ chords. Take an annulus and glue $n$ bands on its boundary according to $w$, as in the figure next page. You get a surface $S(w)$ with boundary, which has some genus (which is by definition the genus of the closed surface obtained after gluing discs on each boundary component): this is *the genus* $g(w)$ *of the chord diagram*.

There is a nice way to compute this genus, due to Moran[156]. Consider the $n \times n$ matrix whose coefficients $a_{ij} \in \mathbb{Z}/2\mathbb{Z}$ are equal to 1 if the two occurrences of $a_i$ are linked with the two occurrences of $a_j$ and 0 otherwise. This is the *incidence matrix* of the interlace graph, modulo 2.

[156] G. Moran. Chords in a circle and linear algebra over GF(2). *J. Combin. Theory Ser. A*, 37(3):239–247, 1984.

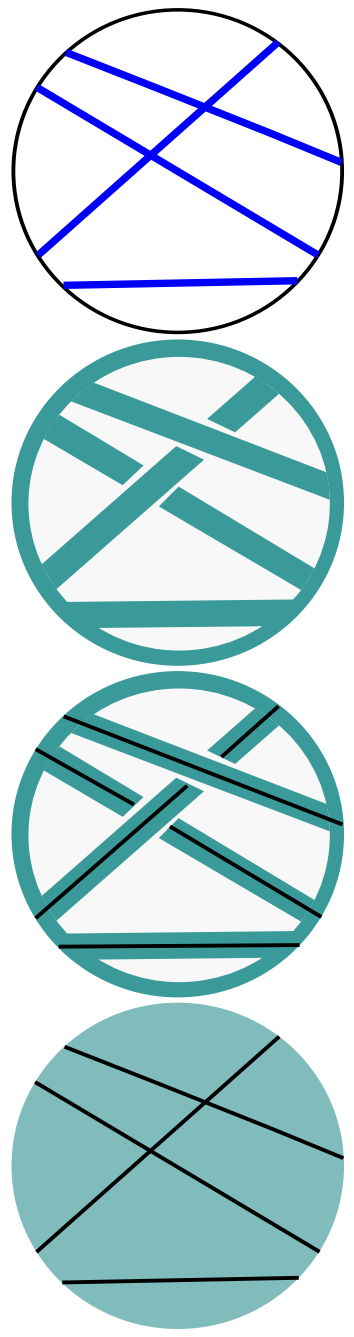

**Theorem.** *The genus of a chord diagram is half the rank of the incidence matrix modulo 2 of the interlace graph.*

The proof given by Moran is rather involved but can be presented in a simpler way.

The chords can be seen as arcs in $S(w)$.

We can also see these chords as arcs in the 2-dimensional disc.

Glue the disc and $S(w)$ along the outer circle to get a surface $S'(w)$. The two copies of each chord define loops $b_1, \dots, b_n$ generating the homology of $S'(w)$. In this basis the intersection of $b_i$ and $b_j$ is 0 if $i, j$ don't link and 1 if they link.

We now glue $k$ discs along the boundary components of $S'(w)$ in order to produce a closed oriented surface $\hat{S}(w)$. The embedding of $S'(w)$ in $\hat{S}(w)$ induces a surjection in the first homology modulo 2. Indeed any loop in $\hat{S}(w)$ can be homotoped away from the discs that have been added. However, this embedding does *not* induce an injection in homology. Indeed, when we glue discs along the boundary, boundary components die in homology since they are now... boundaries of discs. Nevertheless, any element in the kernel is in the kernel of the intersection form of $S(w)$ since it is homologous to a collection of curves parallel to the boundary. It follows that the intersection form on $S'(w)$, modulo its kernel, is isomorphic to the intersection form of $\hat{S}(w)$. For a closed oriented surface of genus $g$ the intersection form is non-degenerate of rank $2g$. ⊡

## *A theorem by Lovász and Marx*

There is a different solution to Gauss's problem. It fits very well with our description of separable permutations as those avoiding the two patterns 3142 and 2413. Interestingly, this theorem is published with no proof[157]. I hope the reader will enjoy finding the proof her(him)self.

[157] L. Lovász and M. L. Marx. A forbidden substructure characterization of Gauss codes. *Acta Sci. Math. (Szeged)*, 38(1–2):115–119, 1976.

Given a generic immersed planar curve there are two ways to delete a given double point, illustrated in the margin. In the first, the curve is split into two components, so that we can choose one of them.

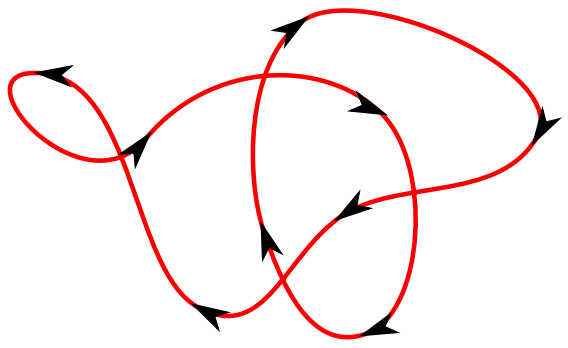

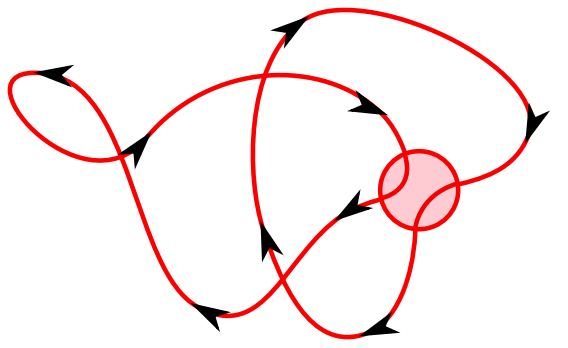

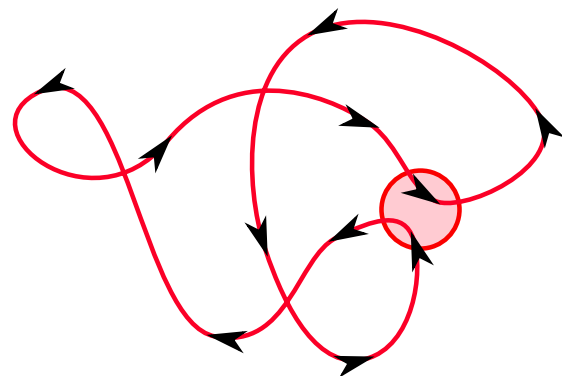

From the combinatorial point of view, these two operations can be expressed in the following way.

– Starting from a word $w = aUaV$, delete $a$ and consider the word $UV^{-1}$ (all words are written cyclically).

– Starting from a word $w = aUaV$, delete $a$ and all the letters that appear in $V$.

Therefore, each word $w$ produces some other words with less letters. The pictures in the margin show that the new words are realizable if the first was. Continue and produce new shorter words. These shorter words are said to be *contained* in the initial word $w$.

**Theorem.** *A Gauss word is realizable if and only if it does not contain the word $a_1a_2\cdots a_na_1a_2\cdots a_n$ for $n$ even.*

Note that the interlace graph of $a_1a_2\cdots a_na_1a_2\cdots a_n$ is the complete graph on $n$ vertices.

## *A Gaussian operad*

I define an operad structure on the set of Gauss words. More precisely, I will deal with generic *marked* oriented immersed curves in the plane $\mathbb{R}^2 \simeq \mathbb{C}$. *Generic* means as before that multiple points are only double points with two different tangents. *Marked* means that we have chosen one of the double points as the "starting" point and that the remaining double points are labeled from 1 to $n$.

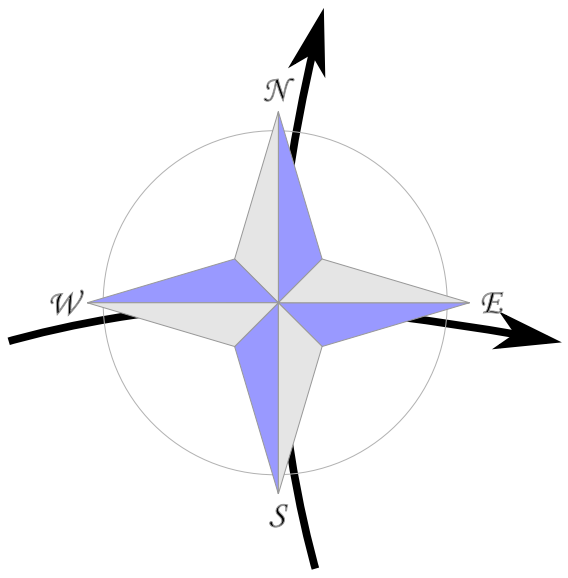

Since the plane and our curve are oriented, any double point defines a compass rose: the four intersection points with a small circle around the point can be labeled by the cardinal directions.

If the marked starting point is sent at infinity by inversion, we get a picture composed of two oriented *long curves*, $\gamma_b, \gamma_r : \mathbb{R} \to \mathbb{R}^2$ (blue and red) with the following properties.

- The blue (resp. red) curve $\gamma_b$ (resp. $\gamma_r$) goes from South to East (resp. from West to North). More precisely $\gamma_b(t)$ is equal to $it$ for large negative values of $t$ and to $t$ for large positive values of $t$ (resp. $t$ and $it$).

- $\gamma_b$ and $\gamma_r$ are immersed and the only multiple points of their union are transversal double points labeled $1, \ldots, n$.

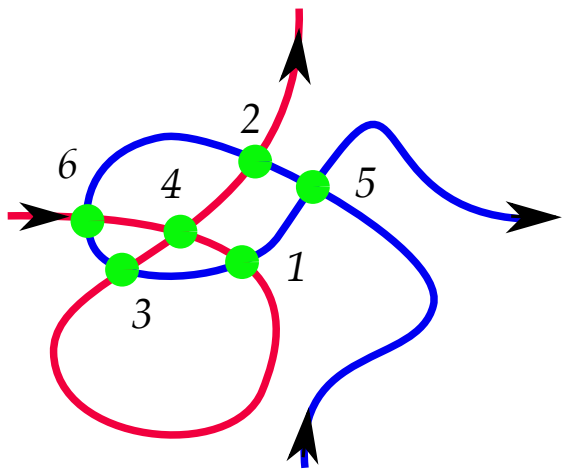

An element of $\Gamma_6$.

Denote by $\Gamma_n$ the set of pairs of curves satisfying these properties and having $n$ double points, up to orientation preserving diffeomorphisms of the plane.

There is a natural operad structure on the union of the $\Gamma_n$'s.

Take a pair $(\gamma_b, \gamma_r)$ of blue-red curves as above having $n$ labeled double points.

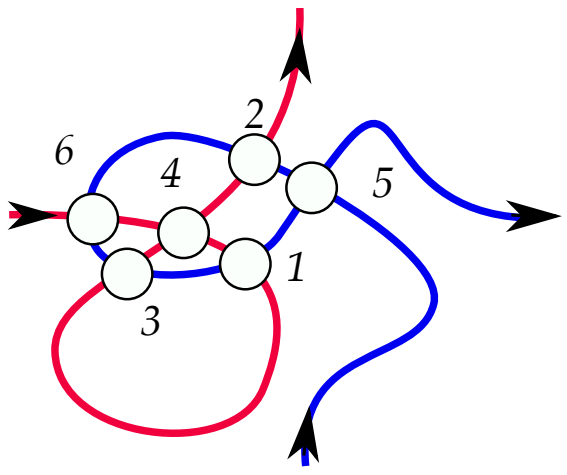

Choose $n$ pairs of blue-red curves $(\gamma_{b,i}, \gamma_{r,i})$ (for $i = 1, \ldots, n$), having $k_1, k_2, \ldots, k_n$ double points. Dig small discs around the double points of $(\gamma_b, \gamma_r)$. We would like to insert the $(\gamma_{b,i}, \gamma_{r,i})$'s into the disc with label $i$, respecting the cardinal directions. However, this is not possible. When we dig a hole, blue curves go from South to North and red curves go from West to East, so that this is not coherent with the South-East and West-North behavior of the blue and red curves $(\gamma_{b,i}, \gamma_{r,i})$ that we want to insert. It is easy to bypass this problem. Before inserting in the $(\gamma_{b,i}, \gamma_{r,i})$'s, it suffices to insert first a standard annulus containing oriented arcs switching North and East.

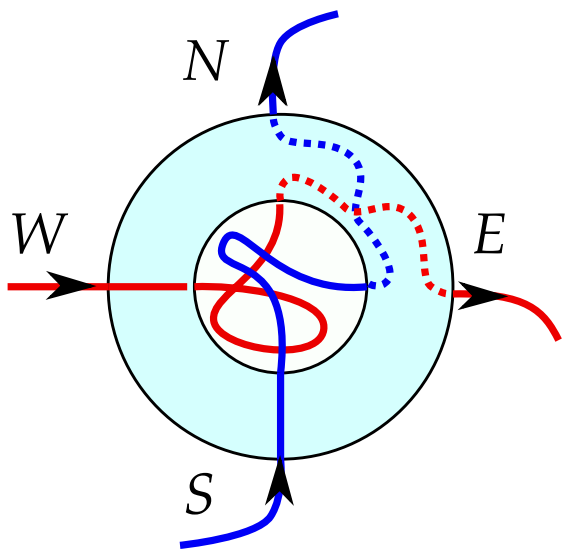

The result of this cut and paste operation is a pair of curves with $k_1 + \cdots + k_n$ double points.

This is the *Gaussian (symmetric) operad*. Can you find a generating system? relations among generators?

To conclude this chapter on some wide opening, I recommend the book[158] which is a remarkable and understandable introduction to Vassiliev knot invariants, where chord diagrams play a crucial role.

[158] S. Chmutov, S. Duzhin, and J. Mostovoy. *Introduction to Vassiliev knot invariants*. Cambridge University Press, Cambridge, 2012.

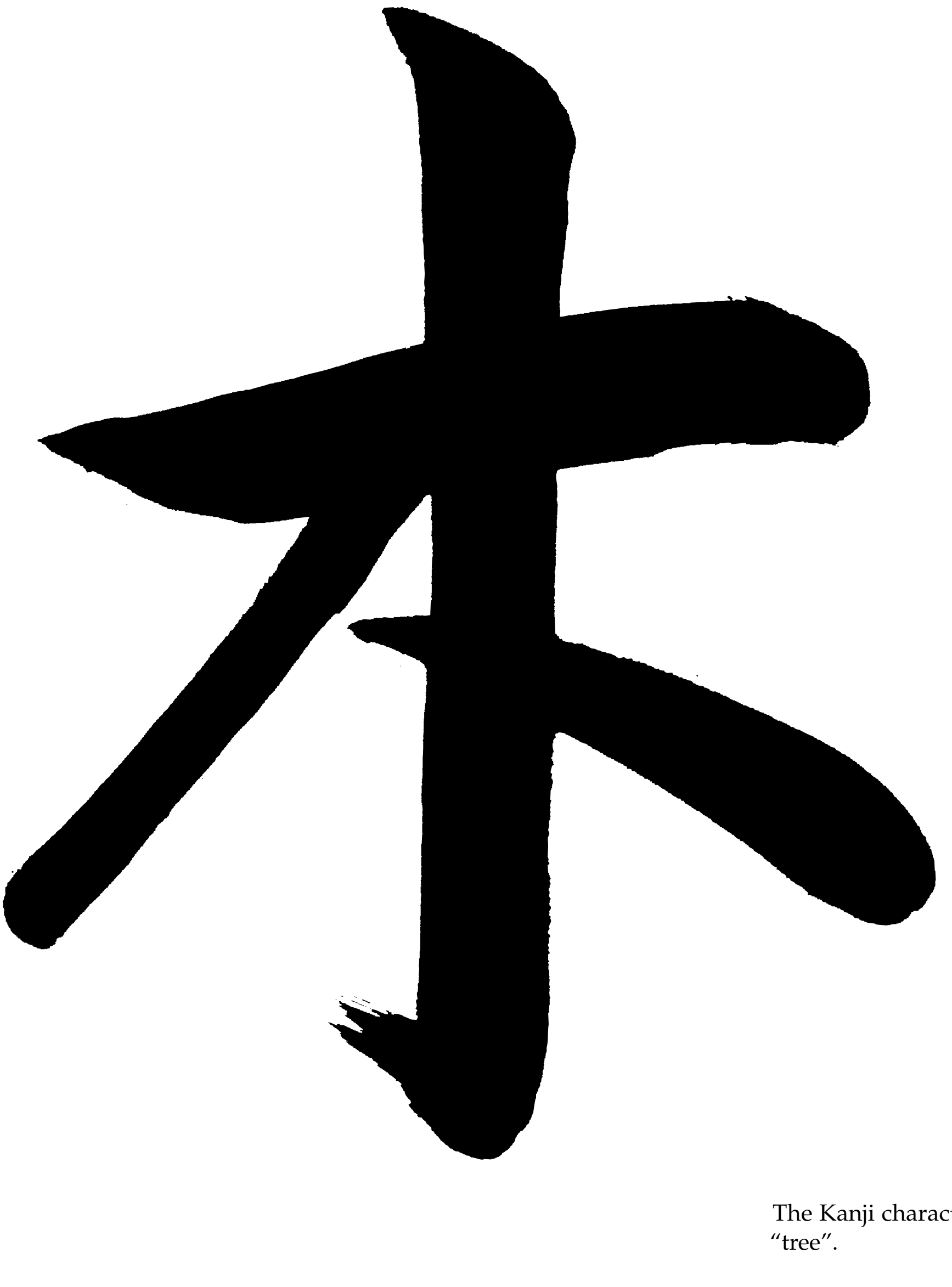

The Kanji character for “tree”.

# *Analytic chord diagrams: an algorithm*

In this chapter, we reach one of our goals: the algorithmic description of the chord diagrams that occur in the neighborhood of a singular point of a planar real analytic curve.

Recall that such a curve intersects a small circle around the singular point at an even number of points which come in pairs, each pair being associated to some real branch.

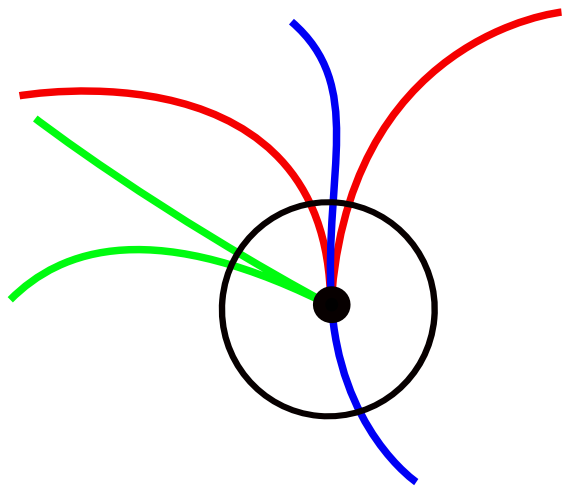

A curve with three branches.

One can think of this structure as a cyclic word of length $2n$ in which every letter occurs exactly twice (where the names of the letters are irrelevant). To be more pedantic (and precise), we are discussing fixed point free involutions on $\mathbb{Z}/2n\mathbb{Z}$ up to conjugacies by cyclic permutations. We can also draw $n$ chords in a circle.

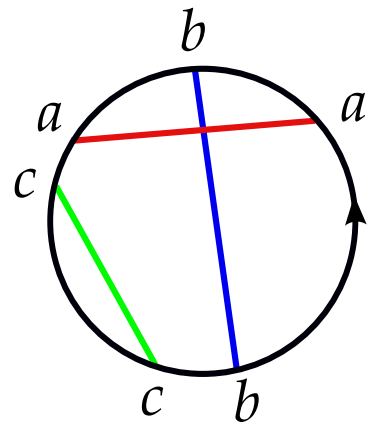

*abaccb*

The total number of these *chord diagrams* of length $2n$ has been studied in many papers. See for instance[159] with strong motivations from knot theory. The problem would be easy if, instead of a cyclic word, we look for standard (non-cyclic) words of length $2n$ in which every letter occurs exactly twice and in which the names of the letters are irrelevant. Indeed, write the first letter of the word and then choose any of the remaining $2n-1$ locations for the other letter which is identical to the first, then write the second letter in the first available free place and choose the other identical letter in any of the $2n-3$ remaining locations etc. Therefore the total number of these words is $(2n-1)\cdot(2n-3)\cdots 3\cdot 1$. These numbers are sometimes

Sometimes I skip the word "chord" in "chord diagram".

[159] A. Stoimenow. On the number of chord diagrams. *Discrete Math.*, 218(1-3):209–233, 2000.

called *double factorials* and denoted by $(2n-1)!!$. See[160] for a presentation of their combinatorial properties.

It would be tempting to divide $(2n-1)!!$ by $2n$ to take into account the cyclic permutations, but some words do admit symmetries and this is why the combinatorics is more subtle. In any case, it follows from these considerations that the number of chord diagrams of length $2n$ grows super-exponentially in $n$.

We will see that a very tiny proportion of chord diagrams are *analytic*, in the sense that they are associated to some singularity of a planar analytic curve.

[160] D. Callan. A combinatorial survey of identities for the double factorial. 2009.

Caution! The double factorial $(2n-1)!!$ is neither the factorial of the factorial nor an exclamation point!

Use Stirling's formula to show that $(2n-1)!!$ is equivalent to $\sqrt{2}\left(\frac{2}{e}\right)^n e^{n\log n}$ when $n$ tends to infinity. So $(2n-1)!!$ is indeed growing super-exponentially, but not much faster. For instance it is small when compared with $C\lambda^{n^{1+\varepsilon}}$ for any $\varepsilon > 0$ and $\lambda > 1$.

## *A necessary condition*

Recall that in the first chapter of this book we showed that for any separable permutation, there are two *consecutive* integers with consecutive images. This was the key point which enabled us to produce an algorithm deciding if a permutation is separable. I now prove a similar property for analytic chord diagrams.

I will say that a chord in a diagram is *solitary* if it connects two consecutive points of the diagram as in the picture on the left below. Two chords are *parallel* (resp. *antiparallel*) if they are as in the second (resp. third) picture, i.e. if the corresponding letters $a, b$ occur in the cyclic word as $\cdots ab \cdots ba \cdots$ (resp. $\cdots ab \cdots ab \cdots$). Finally, two chords as in the fourth picture constitute a *pitchfork* ($\cdots a \cdots bab \cdots$).

☺ Some theorems or lemmas are so famous that it seems to be forbidden to use words like "Theorem A, B" (Cartan), or "Fundamental lemma" (Ngo).

**Fundamental lemma ☺.** *Any analytic chord diagram contains a solitary chord, a pair of parallel or an antiparallel chords, or a pitchfork.*

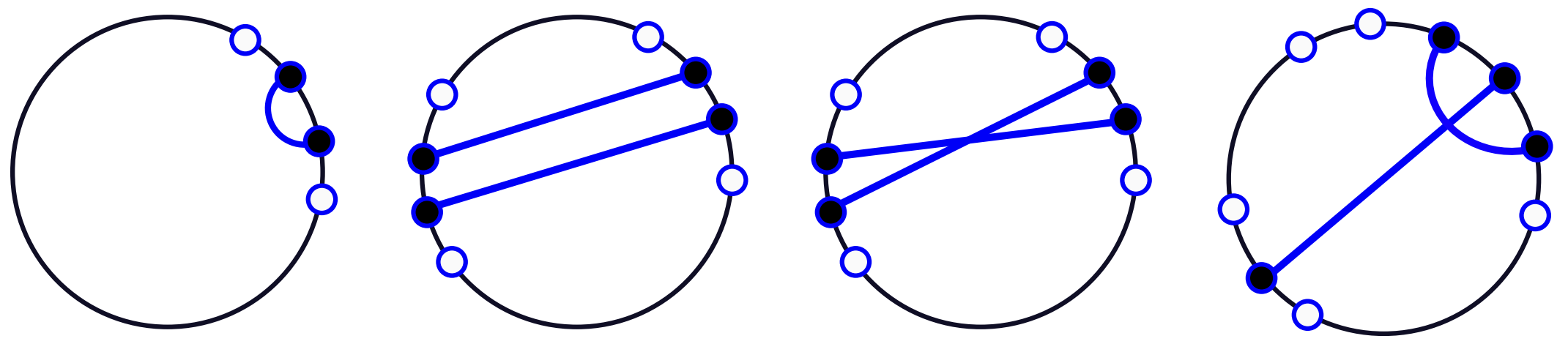

Let us observe that this implies immediately the theorem that was stated in the preface.

**Theorem.** *There is no singular analytic curve in the plane consisting of five branches intersecting a small circle as in the picture in the margin.* ⊡

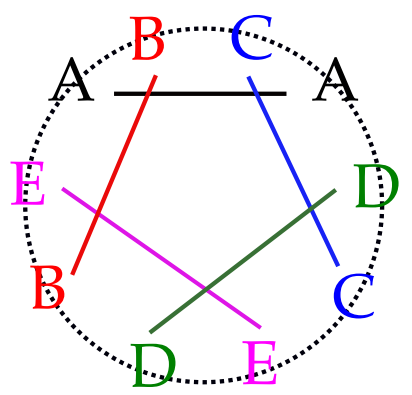

Impossible five branches.

## *Let us blow up*

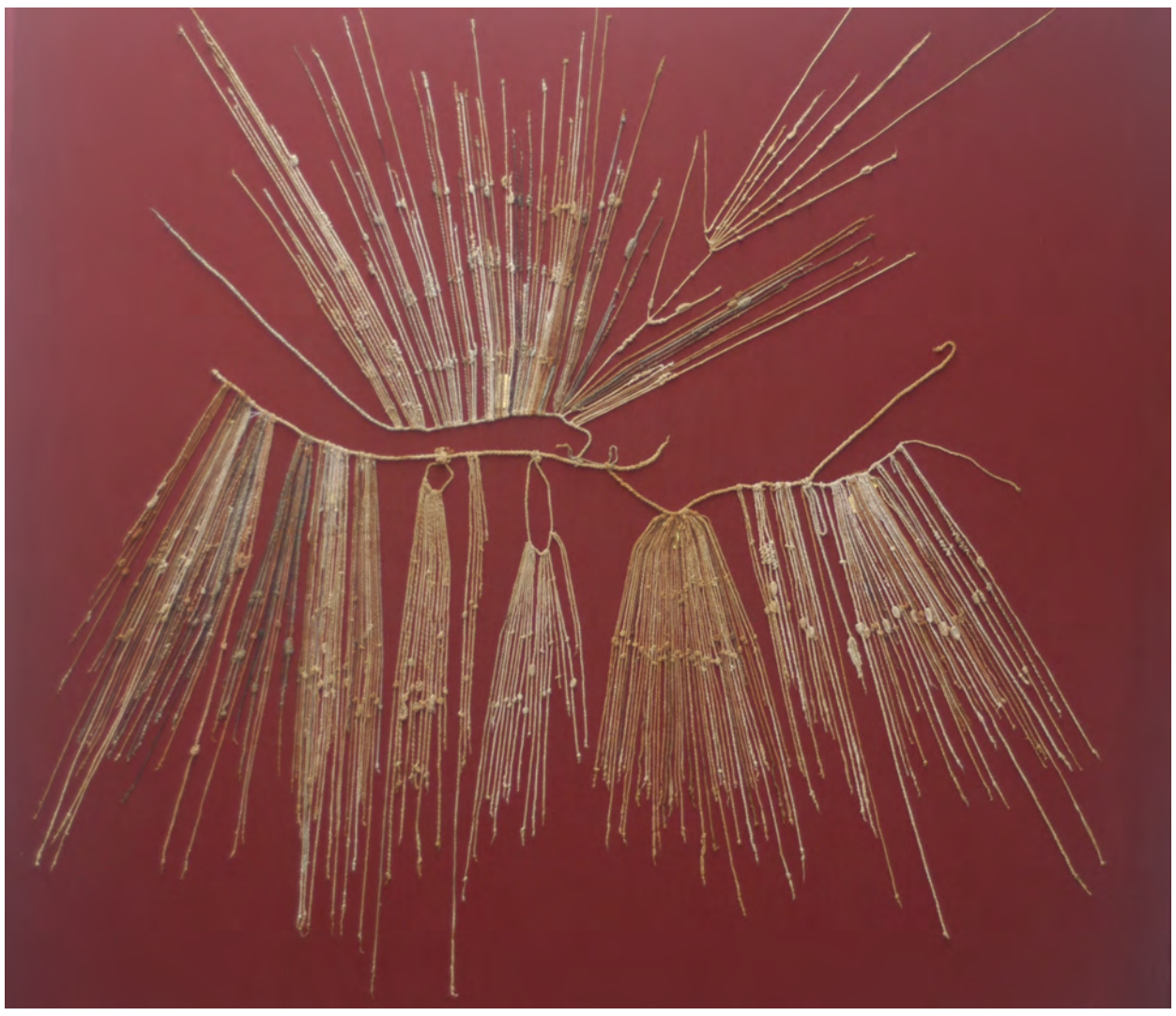

A beautiful *Quipu*: a knotted-string device that was used by the Incas for recording statistical information. Like a blown up projective line? (Centro Mallqui, Leymebamba, Peru)

Start with some singular point of some analytic curve in the (real) plane. Blow it up a first time. The result is a curve in some Moebius band, whose singular points are on the exceptional divisor, core of the band. If things go well, the singular point splits into several singular points, presumably simpler. Let us blow up all of them. It could happen that after one blow up, there is still a unique singular point on the divisor. Then, blow it up a second time. Let us continue the blowing up process as many times as necessary. We know that after some time, the singularity will be *resolved*. This means that the strict transform of the initial curve is now a collection of $n$ disjoint smooth analytic curves intersecting transversally the exceptional divisor.

This exceptional divisor is a union of real projective lines which are circles intersecting transversally. Consider the graph whose vertices are these projective lines and where an edge connects two vertices if the projective lines intersect. This graph is a tree as can be easily seen by induction. Indeed, in the inductive process of desingularization, at each step we blow up a point which can be either a smooth point of the exceptional divisor, or an intersection of two projective lines. In the first case, a new leaf is grafted to a tree and in the second case, an edge is split into two edges. The first projective line, coming from the first blow up, can be chosen as the root of this tree.

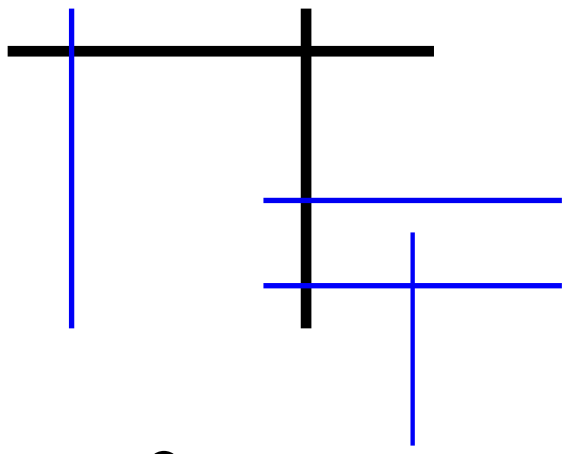

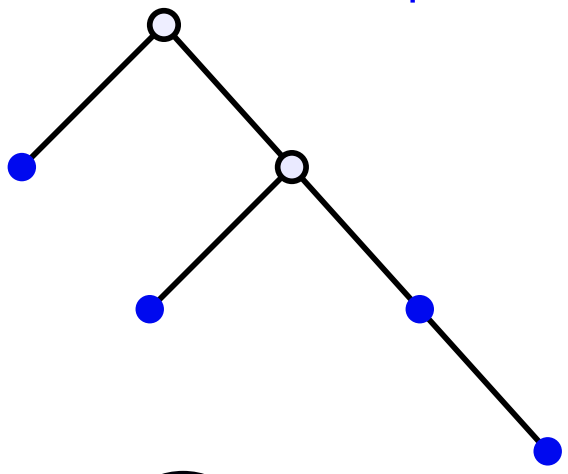

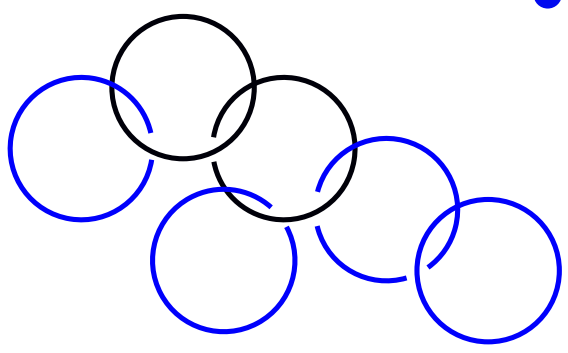

It will be convenient to blow up once more each of the $n$ points on the exceptional divisor, if necessary, introducing new projective lines, in order to make sure that at the end of the process each projective line contains *at most* one point of the strict transform. We can even suppose that those components of the divisor which meet the strict transform are leaves of desingularization tree.

Let us sum up. Given some analytic curve $\mathcal{C}$ defined in a neighborhood of $(0,0)$ in $\mathbb{R}^2$ by some equation $F(x,y) = 0$, we can construct the following objects.

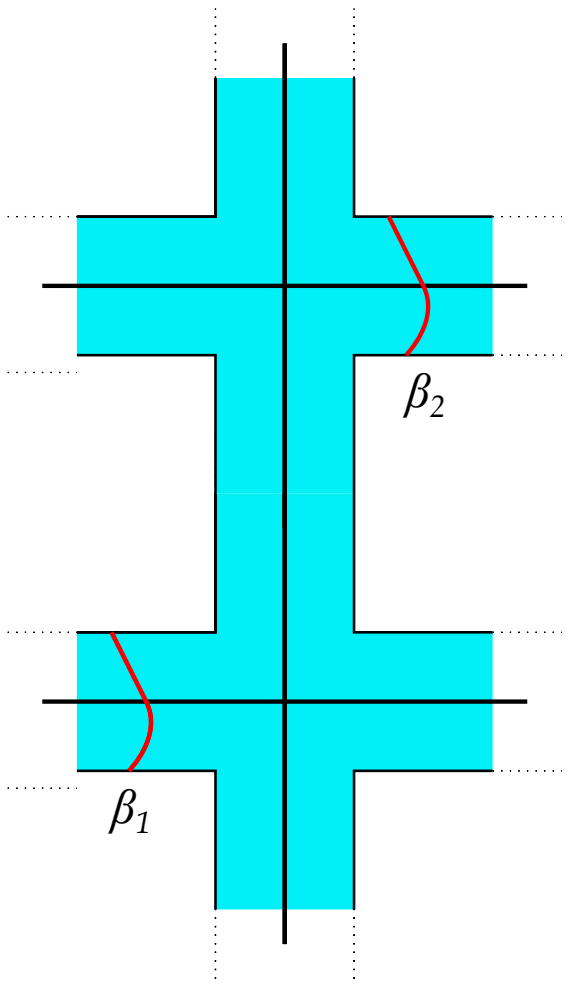

- A surface $S$ with connected oriented boundary.
- An exceptional divisor $E \subset S$, consisting of a certain number of circles intersecting transversally, each pair meeting at most once. The associated intersection graph is a rooted tree. The embedding $E \subset S$ is a homotopy equivalence.
- A finite disjoint union $\hat{\mathcal{C}}$ of *smooth* analytic arcs $\beta_1, \ldots, \beta_n$ in $S$ intersecting transversally $E$. The intersection of $\hat{\mathcal{C}}$ with of a component of $E$ is empty if this component is not a leaf of the tree, and contains at most one point if it is a leaf. We can assume moreover that $\hat{\mathcal{C}}$ is transversal to the boundary of $S$ and that each arc $\beta_i$ intersects the boundary in two points.
- A blowing down analytic map $\Psi : S \to \mathbb{R}^2$, collapsing $E$ to the origin, which is a diffeomorphism from $S \setminus E$ onto some small punctured disc, and which collapses $\hat{\mathcal{C}}$ to our singular curve $\mathcal{C}$.

Recall that each loop in $S$ can be orienting or disorienting. Let $\gamma$ be a closed immersed curve in a surface, passing once through a point $x$. When the surface is blown up at $x$, the self-intersection *modulo 2* of the strict transform of $\gamma$ (in the blown up surface) is equal to the self-intersection of $\gamma$ (in the original surface) plus 1. In the inductive construction, when a projective line appears for the first time in the exceptional divisor, it is the core of a Moebius band, of self-intersection 1. Later on, some of its points may be blown up. Each of these blowing ups permutes the orienting/disorienting status of a component of the divisor. The previous pictures in the margin (six lines, a tree with six vertices, and six circles) correspond to the same example that was analyzed in the chapter on necklaces: it is obtained after six blowing ups. I indicated in blue the components which correspond to Moebius bands.

Some of the components of $E$ intersect the desingularized curve $\hat{C}$: they define some of the leaves in the desingularization tree. Call those leaves *colored*. Observe that some leaves might be non-colored.

Note that if we choose some orientation of each component of $E$, the corresponding tree is planar so that the children of any node are linearly ordered. Changing the orientation reverses this order.

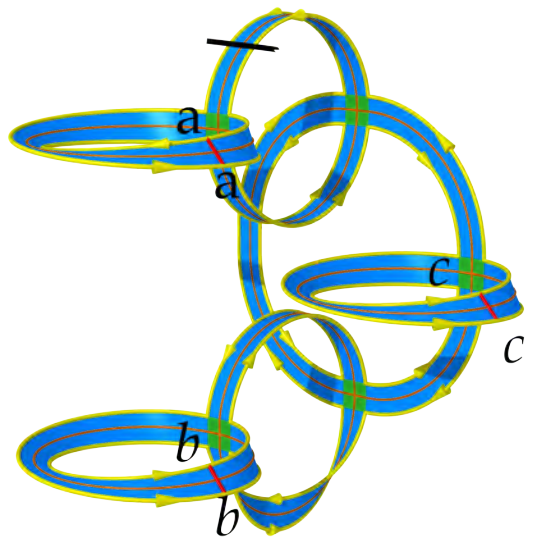

## *An example*

Look at the necklace in the margin. This is still the same object that was already discussed in the chapter on Moebius necklaces. Six blow ups produced six bands, two orientable and four non-orientable. The exceptional divisor consists of the six cores of the six bands. The desingularized curve consists of three red arcs, labeled $a, b, c$, each intersecting the boundary of $S$ in two points. On top, we see in black the strict transform of the $y$ axis. Going around the boundary of $S$ we can read the corresponding analytic chord diagram. Just follow the arrow and read *abacbc*. I must admit that it is not so easy to follow the arrows!

## *Proof of the fundamental lemma*

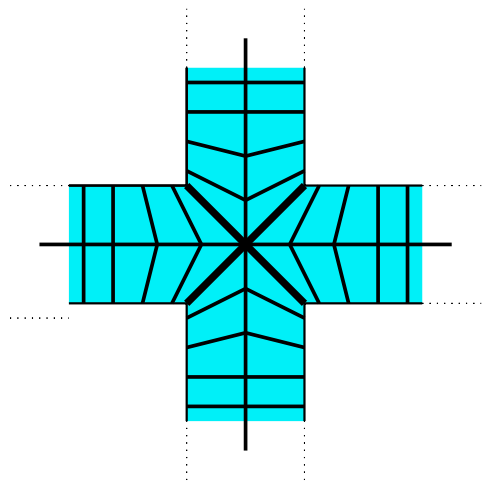

The fibers of $\pi$. Do not confuse $\pi$ (from $S$ to the divisor $E$) with $\Psi$ (from $S$ to a disc) which is blowing down $E$ to a point.

A first trivial observation is that the operation of deleting a chord in an analytic diagram transforms it into some other analytic diagram. It corresponds to deleting a branch.

Start with some chord diagram $w$ associated to some singularity of a planar analytic curve. Consider a desingularization tree as before. There is a projection $\pi$ of the surface $S$ onto the exceptional divisor $E$ which is a homotopy equivalence. For every point $x \in E$, the fiber $\pi^{-1}(x)$ is an arc connecting two points of the boundary if $x$ is a regular point, and a cross if $x$ is the intersection of two circles.

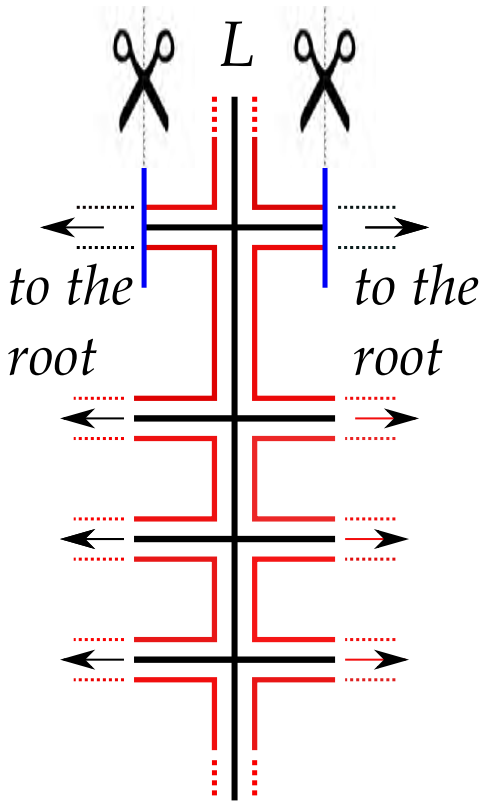

Let $L$ be some node of the tree, i.e. one of the projective lines that constitute the exceptional divisor $E$. There is a unique chain of nodes going from $L$ to the root. Cut two disjoint arcs in $S$ as in the figure, in order to disconnect $L$ from the root in $S$. The four endpoints of these arcs decompose the circle boundary of $S$ in four intervals. Two of them (colored red) correspond to "what is in $L$ or below $L$" in the tree. Going around the boundary of $S$ and reading the chord diagram, we therefore find two disjoint intervals of letters, below $L$, whose union is stable under the involution sending each occurrence of a letter to the other occurrence. Note that these intervals could be empty if there is no colored leaf below $L$. If there is a colored leaf below $L$, at least one of the two intervals is non-empty, but it could be the case that only one is not empty.

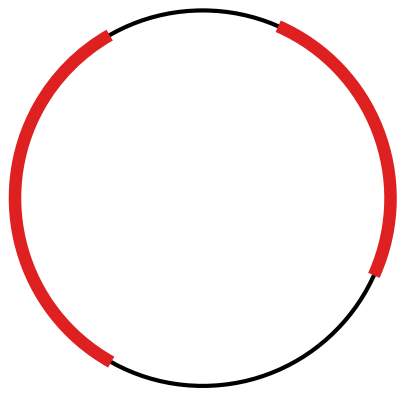

In summary, *every node $L$ in the tree defines a chord diagram $w(L)$ which is a sub-diagram of the original diagram $w$ and which is "connected" in the sense that its letters form one or two intervals in $w$.*

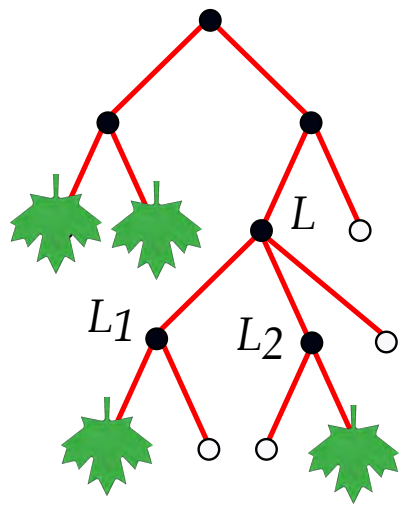

Colored leaves are green!

Think of a rooted tree as a genealogy tree, the root being the founding member of the family. Each node has a certain number of descendants, some of them being colored leaves. Let $L$ be one of the youngest members of the family having at least two colored leaves as descendants. Among the children of $L$, let $L_1, \ldots, L_k$ be the list of those having at least one colored descendant (ordered in this way along $L$). We have $k \geq 2$ since otherwise one of the children of $L$ would have at least two colored descendants. For the same reason, each $L_i$ has a unique colored descendant.

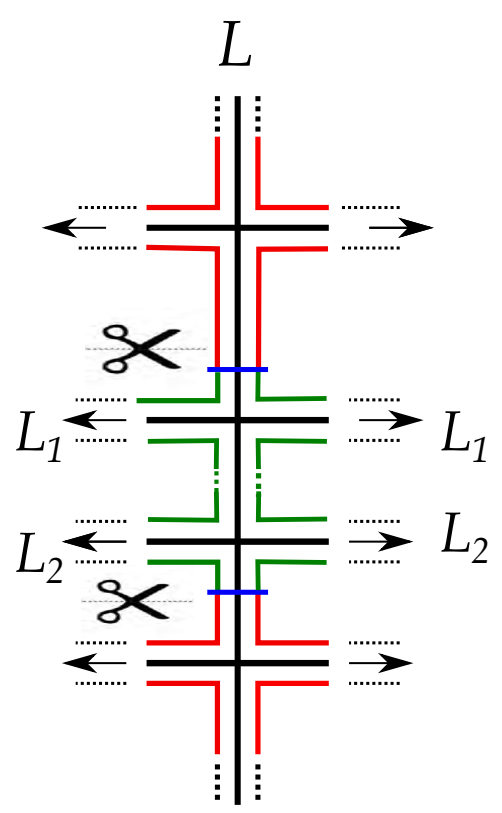

Now cut the surface $S$ as in the margin figure, disconnecting $L_1$ and $L_2$ from the root and from all other colored leaves. As before, this defines two (green) intervals on the boundary of $S$ whose union contains exactly four points of our initial diagram, associated to two chords. Two chords in two intervals can be organized in the following fifteen ways.

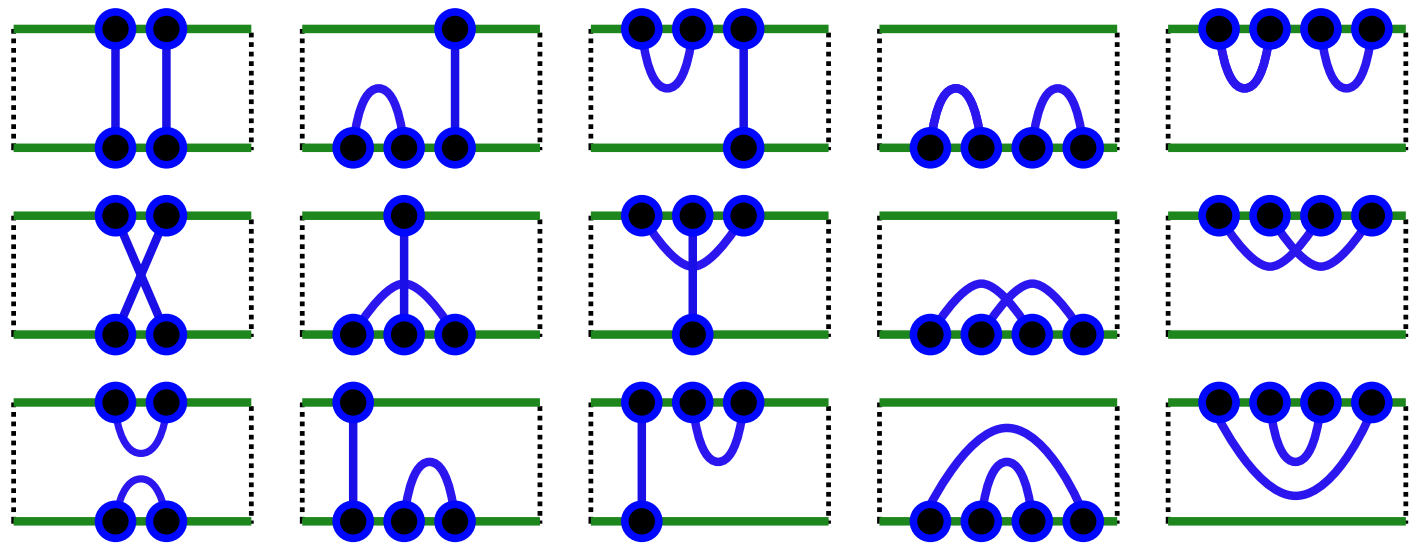

In each case, there is a solitary chord, a pitchfork, or a pair of parallel or antiparallel chords.

This is the end of the proof of the fundamental lemma. ⊡

## More non-analytic diagrams

We have observed that the deletion of some letters in an analytic chord diagram produces another analytic diagram. Let us say that a diagram is *basic non-analytic* if it is non-analytic but becomes analytic as soon as a single chord is deleted. Clearly a chord diagram is analytic if and only if it does not contain a *basic non-analytic chord diagram*. Recall that a permutation is separable if and only if it *does not* contain the Kontsevich permutation $(2,4,1,3)$ or its reverse permutation $(3,1,4,2)$, so that in this case there are only two basic non-separable permutations. The situation is more complicated in the case of chord diagrams.

**Theorem.** *There is an infinite number of basic non-analytic chord diagrams.*

Here is an example that will be denoted by $\mathcal{C}_n$ ($n \geq 5$). Consider the $2n$ points of $\mathbb{Z}/2n\mathbb{Z}$ ordered in a natural way on the circle. The chord diagram pairs $2k$ and $2k+3$ for $k = 1, \dots, n$. For $n = 5$, this is our previous example of non-analytic diagram with five chords. This diagram $\mathcal{C}_n$ ($n \geq 5$) is not analytic for the same reason as in the case $n = 5$. We still have to show that if one letter is deleted, the remaining diagram is analytic. For this, we need a *sufficient analyticity condition*. This will be provided by a very simple algorithm deciding if a chord diagram is analytic.

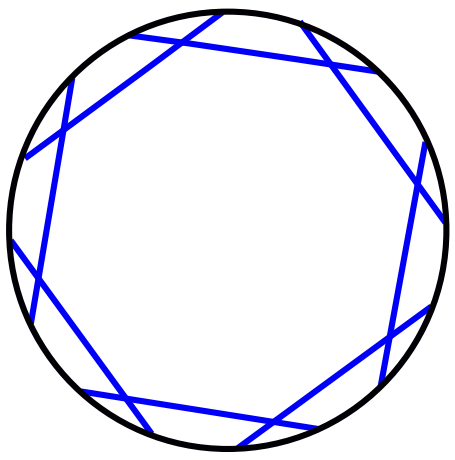

**Theorem.** *The following algorithm decides if a chord diagram is analytic:*

1. *If there is no solitary chord, no pitchfork, and no pair of parallel or antiparallel chords, the diagram is* not *analytic.*

2. *If there is a solitary chord, delete it and continue. If there is a pitchfork, delete the "small chord" and keep the handle. If there is a pair of parallel or antiparallel chords, delete one of them and continue.*

3. *If you end up with the empty diagram, the original one was analytic.*

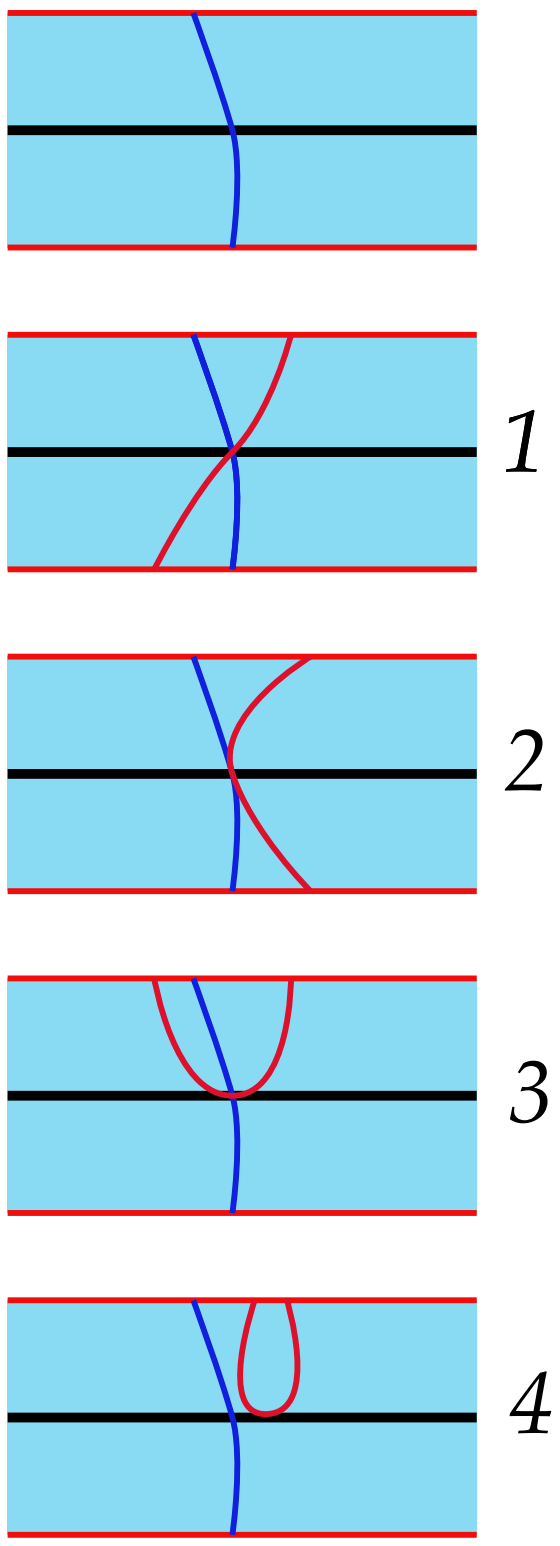

The proof is easy. We have to show that if $w$ is a diagram and if $\overline{w}$ is the new diagram with one less chord obtained after one step of the algorithm, then $w$ is analytic if $\overline{w}$ is analytic. The converse is easy: if $w$ is analytic, $\overline{w}$ is also analytic since it corresponds to deleting a branch. We now have to *add* an additional branch.

Choose some desingularization of some singular point on some surface $S$. A chord corresponds to some smooth arc $\gamma$ (in blue in the margin) connecting two points on the boundary, transverse to the divisor at some point $x$. Add a new analytic smooth (red) curve $\gamma'$ in $S$, also transverse to the divisor, very close to $\gamma$ and transverse to $\gamma$, as in 1/. The implosion of $S$ produces a new singular point with one more branch. Clearly, the new diagram has one more chord which is parallel or antiparallel to the initial chord, depending on the orientations on the boundary. Choose now $\gamma'$ as in 2/ with a quadratic tangency with $\gamma$ at $x$ and you get the other parallel or antiparallel situation.

In case you want to create a pitchfork with a given handle, just add a smooth curve $\gamma'$ close to $\gamma$ with a quadratic tangency with the divisor as in 3/.

Finally, if you want to add a solitary chord right after some given letter, proceed as in 4/. ⊡

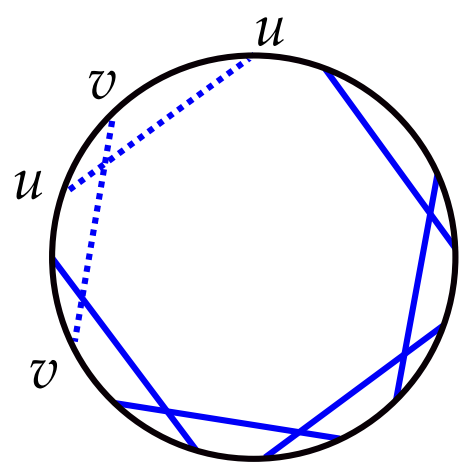

Let us test our algorithm on the previous diagram on $\mathcal{C}_n$ ($n \geq 5$). Deleting one chord, we get a non cyclic chain of $n-1$ chords. The first two chords $u, v$ define a pitchfork, so we may delete the first chord and continue until there is only one chord left. This diagram is therefore analytic. There is indeed an infinite number of basic non-analytic diagrams.

## *With a computer*

In order to count analytic chord diagrams, let us use a computer to test for small values of $n$. This is easy. We first list all possible words of length $2n$ in which each letter occurs twice. The only subtlety is to take into account the cyclic character of the word under consideration. Here is the result for $n \leq 7$.

In the following table:

| $n$ | 2 | 3 | 4 | 5 | 6 | 7 |
|---|---|---|---|---|---|---|
| Words | 3 | 15 | 105 | 945 | 10395 | 135135 |
| Chord diagrams | 2 | 5 | 18 | 105 | 902 | 9749 |
| Up to symmetry | 2 | 5 | 17 | 79 | 554 | 5283 |

– *Words* means "words of length $2n$ in which each letter occurs twice". The number of these words is the double factorial of $2n-1$.

– *Chord diagrams*, as we have defined them, are words up to cyclic permutations.

– The image of a chord diagram by a symmetry with respect to some line is another diagram, which may be the same diagram or not. The item "up to symmetry" counts the number of words up to these dihedral symmetries.

We then count the number of *analytic* diagrams. This is in principle not difficult, using the algorithm that was described earlier.

The On-Line Encyclopedia of Integer Sequences, A007769, A054499, A279207, A279208. Neil Sloane added the last two sequences based on a preliminary version of this book.

The result is:

| $n$ | 2 | 3 | 4 | 5 | 6 | 7 |
|---|---|---|---|---|---|---|
| Analytic diagrams | 2 | 5 | 18 | 102 | 817 | 7641 |
| Up to symmetry | 2 | 5 | 17 | 76 | 499 | 4132 |

It follows that *for $n \leq 4$, all diagrams are analytic.*

*Among the* 105 *diagrams with* 5 *chords, only the* 3 *examples in the margin are not analytic.*

The first diagram is already familiar, under the name $\mathcal{C}_5$. Let me denote the others by [symbol] and [symbol]. It was not difficult for me to guess the first but I must admit that I did not find the two others by hand but with a computer.

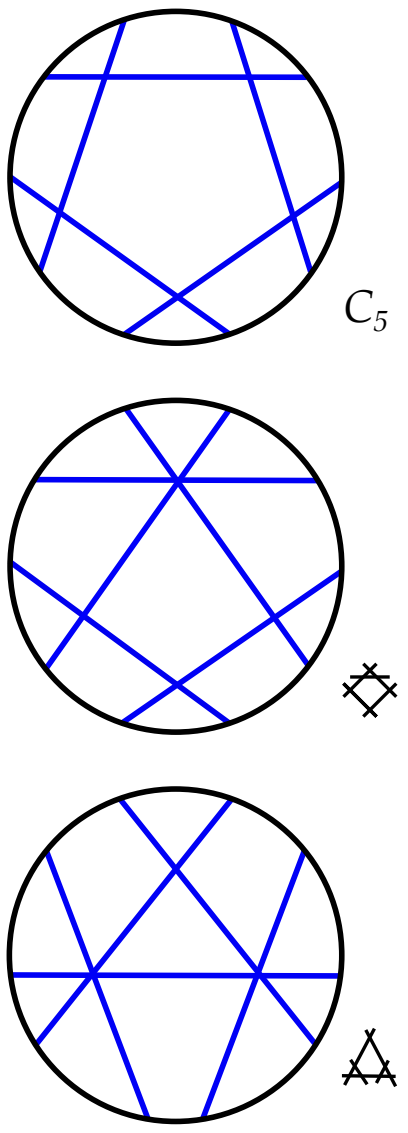

The 3 basic non-analytic chord diagrams with 5 chords.

*Among the* 902 *diagrams with* 6 *chords,* 85 *are not analytic.* However, the non-analyticity of most of them is due to the fact that one of their sub-diagrams is not analytic. There are only two diagrams with 6 chords which are *basic*: they are non-analytic and all their sub-diagrams are analytic.

Observe that the first one is the member $\mathcal{C}_6$ of the infinite family of basic diagrams $\mathcal{C}_n$. It corresponds to $\mathbb{Z}/12\mathbb{Z}$ where every even number $k \pmod{12}$ is connected to $k+3 \pmod{12}$. The second will be denoted by [symbol].

*Among the* 9749 *diagrams with* 7 *chords,* 2108 *are not analytic.* The only basic non-analytic example is $\mathcal{C}_7$.

In the next chapter, I will show that my computer did find *all* basic non-analytic diagrams.

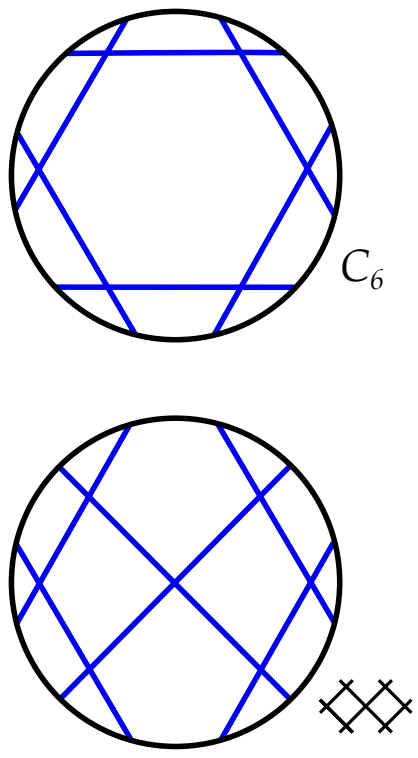

The 2 basic non-analytic chord diagrams with 6 chords.

## *Marked chord diagrams*

To conclude this chapter and to place it on a more general context, I would like to describe an operad structure which is behind the curtain. It will be convenient to introduce first a slight strengthening of the notion of analytic chord diagrams.

When we proved that our algorithm does decide if a diagram is analytic, the key point was the possibility of inserting a new branch. It turns out that more complicated singularities can also be inserted, as I explain now. Consider a desingularization of some curve $\mathcal{C}$ as before, so that we have a surface $S$, a divisor $E$, and a collection of smooth curves $\beta_1, \dots, \beta_n$ intersecting $E$ transversally in a finite number of points, $p_1, \dots, p_n$, where $n$ is the number of real branches. Choose one of these points, say $p_1$. Choose now some *other* singular curve $\mathcal{C}_1$, with $n_1$ real branches, and assume that it does not contain the $y$ axis. Delete $\beta_1$ and replace it by a copy of $\mathcal{C}_1$ in the surface $S$, in such a way that the $y$ axis for $\mathcal{C}_1$ is mapped into the divisor $E$ and the singular point of $\mathcal{C}_1$ is mapped to $p_1$. We can now blow down the union of this copy of $\mathcal{C}_1$ and $\beta_2, \dots, \beta_n$. The result is a new singular point, with $n + n_1 - 1$ branches: one of the branches of $\mathcal{C}$ has been replaced by a copy of $\mathcal{C}_1$.

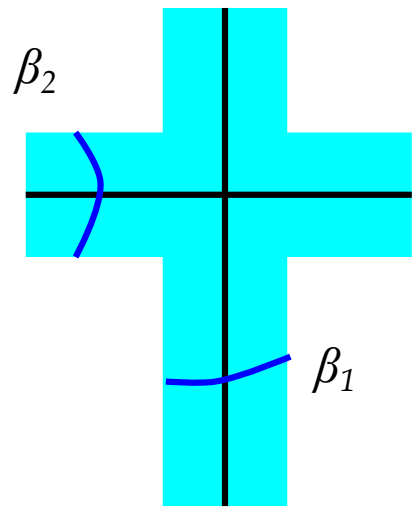

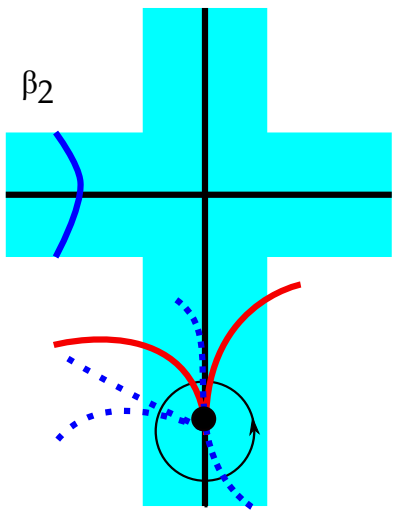

Let us examine the effect of this kind of operations on the associated chord diagrams. Looking at the diagram associated to $\mathcal{C}_1$, we see that the $y$ axis decomposes the word of length $2n_1$ in two components, *Left* and *Right*. In the new chord diagram with $2(n + n_1 - 1)$ letters, one pair of identical letters from the old diagram with $2n$ letters has been replaced by two intervals, which are *Left* and *Right*. We should be careful however, that in this process, the orders of the letters in *Left* and *Right* might have been reversed. Indeed, the two intersections of the oriented boundary of $S$ with $\beta_1$ might be of different signs. Moreover, the insertion of $\mathcal{C}_1$ in $S$ can be done in four ways since $S$ is not orientable and $E$ is not oriented.

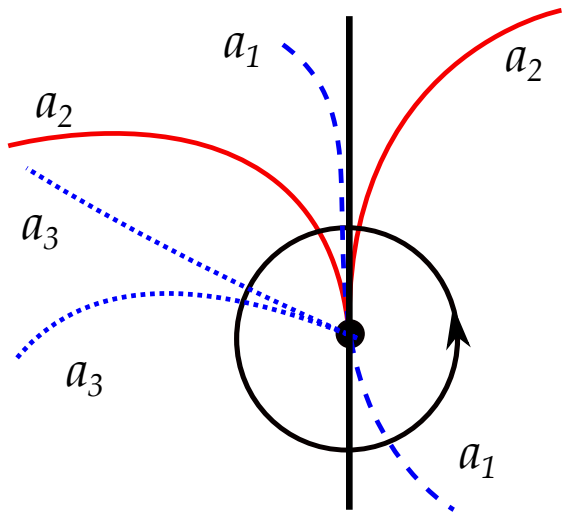

$Right = a_1 a_2$, $Left = a_1 a_2 a_3 a_3$

Note that *Left* or *Right* might be empty.

Of course, we can proceed in the same way with all other branches of $\mathcal{C}$, using other singular curves $\mathcal{C}_2, \dots, \mathcal{C}_n$.

All these remarks suggest the following definition.

In an equation $F(x, y) = 0$, we can replace $(x, y)$ by $(-x, y)$ or $(x, -y)$ or $(x, xy)$. The transformation

$$(x, y) \mapsto (x, xy)$$

preserves each vertical line, collapses the axis $x = 0$ to the origin, and reverses the orientation for $x < 0$. This is not a surprise: this is a blow down map. The square of this transformation preserves the orientation on each vertical line (for $x \neq 0$).

**Definition.** A marked chord diagram is a collection of $2n$ distinct points $a_1^{\pm 1}, \dots, a_n^{\pm 1}$ in the union of two opposite sides of a square $\{-1, 1\} \times [-1, 1]$ (up to orientation preserving homeomorphisms of each side).

Note the additional features if one compares with standard diagrams. Marked diagrams have a right and a left part. Moreover, each chord $a_i$ is now *labeled* with a number $i$ from 1 to $n$ and is *oriented* from $a_i^{-1}$ to $a_i^{+1}$.

Let me denote by $\mathcal{AMC}$ the set of those marked chord diagrams which are *analytic*, i.e. which arise from some analytic curve $F(x,y) = 0$ which does not contain the $y$-axis. Note that the analyticity of a marked diagram depends neither on the orientation of the chords nor on the labeling.

The role of the labelings and orientations is simply to give the relevant information about which marked diagram is inserted in each chord, and in which way.

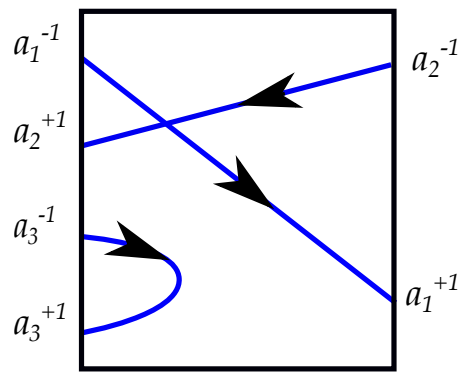

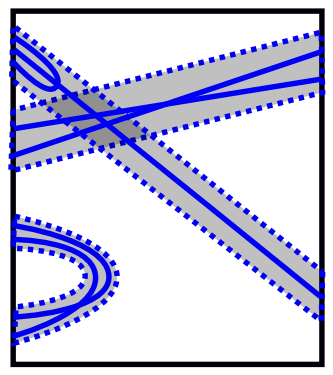

Now, I hope that the reader has guessed the operad structure on $\mathcal{AMC}$. Let $w$ be some analytic marked chord diagram with $n$ chords. Given $n$ analytic marked chord diagrams $w_1, \ldots, w_n$, with $k_1, \ldots, k_n$ chords, define the action of $w$ on $(w_1, \ldots, w_n)$ in the following way. Draw $w$ and thicken each chord $a_i$ of $w$, creating rectangles. Use the $a_i^{-1}$'s and the $a_i^{+1}$'s as the left and right sides of these new rectangles. Now, insert $w_1, .., w_n$ in these rectangles respecting the labels and the orientations. Rename the chords, from 1 to $k_1 + k_2 + \ldots + k_n$ using the lexicographic ordering. The result is another marked analytic chord diagram since this operation corresponds to the previously described insertion of analytic curves.

Hence, $\mathcal{AMC}$ has indeed a natural operad structure. For completeness, let us observe that the natural action of the symmetric groups, permuting the chord labels, shows that this is actually a symmetric operad.

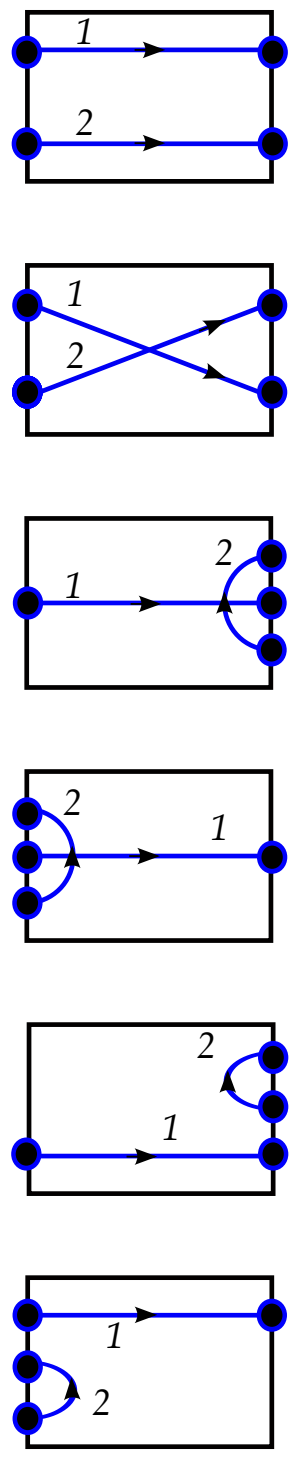

## *Let us bound the number of chord diagrams*

It would be great to have some precise information on the number $a_n$ of analytic diagrams with $n$ chords. For instance, an explicit formula for the generating series $\sum a_n t^n$ would give the exact exponential growth rate of $a_n$. Unfortunately, I was not able to compute this series ☹. In this section, I show that the fundamental lemma provides at least a reasonable bound.

Consider a finite planar rooted *binary* tree. Equip each of its interior nodes (including the root) with one of the six examples of marked diagrams with two chords represented in the margin. By recursive insertions of the diagrams of the siblings in the diagram of their parent, this produces a marked analytic diagram, and hence an analytic diagram, forgetting the labels, the orientations, and the two sides of the square.

I claim that all analytic diagrams with $n$ chords are produced with this recipe. This is true for $n = 2$ since both diagrams with 2 chords (linked of not linked) appear when one forgets the marking in the six examples. Now, let $w$ be some analytic diagram with $n+1$ chords, and apply the fundamental lemma. Therefore, we find $\cdots aa \cdots$, $\cdots ab \cdots ba \cdots$, $\cdots ab \cdots ab \cdots$ or $\cdots b \cdots aba \cdots$ in the diagram. In the case of $aa$, call $b$ the letter which comes *before* $a$ in the cyclic order. Our algorithm deletes $a$ and produces an analytic diagram $\bar{w}$ with $n$ chords, for which we can apply the induction. This means that $\bar{w}$ can be "clothed" with labels and orientations in such a way that it is produced by a binary tree, as above. Our diagram $w$ is obtained from $\bar{w}$ by replacing one chord by two chords. It is easy to check that our six examples are sufficient to realize this duplication using an insertion in the operad. Hence $w$ is constructed from a binary tree with $n+1$ leaves with the same recipe.

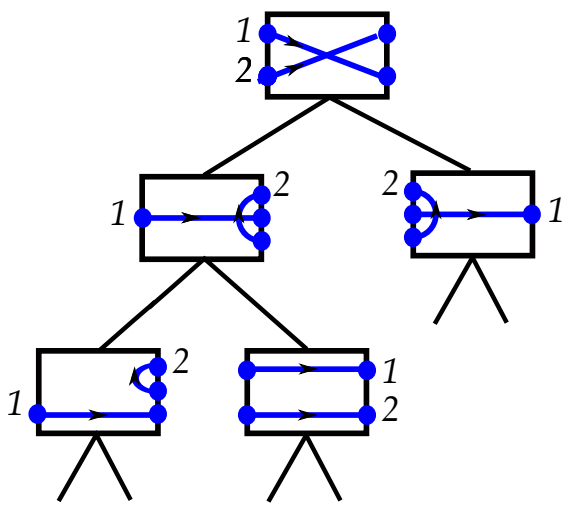

A rooted binary tree with $n$ leaves has $n-1$ interior nodes (including the root) so that there are $6^{n-1}$ possible labels on the interior nodes. The number of planar binary trees with $n$ leaves is given by the $(n-1)$-st Catalan number. Therefore, we get the following rough estimate.

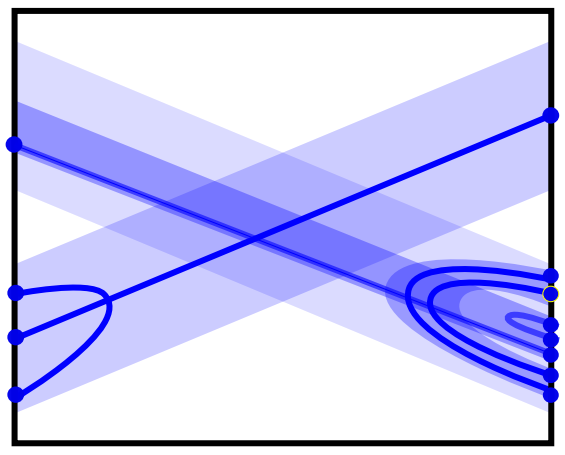

**Theorem.** *The number $a_n$ of analytic chord diagrams with $n$ chords is less than $6^{n-1}$ times the $(n-1)$-st Catalan number $C_{n-1}$.*

Recall that $\frac{1}{n}\log C_n$ converges to $\log 4$ when $n$ tends to infinity. Therefore

$$\limsup \frac{1}{n}\log a_n \leq \log(24).$$

Note that any permutation on $n$ letters can be seen as a diagram with $n$ chords, such that all its chords connect points on both sides of a square. In particular separable permutations produce analytic marked diagrams. This gives a lower bound for the growth of $a_n$ since we already counted separable permutations.

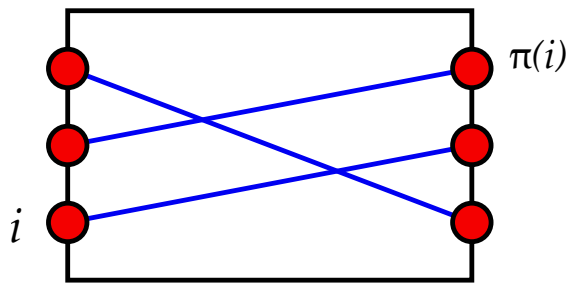

$$\liminf \frac{1}{n}\log a_n \geq \log(3+2\sqrt{2}).$$

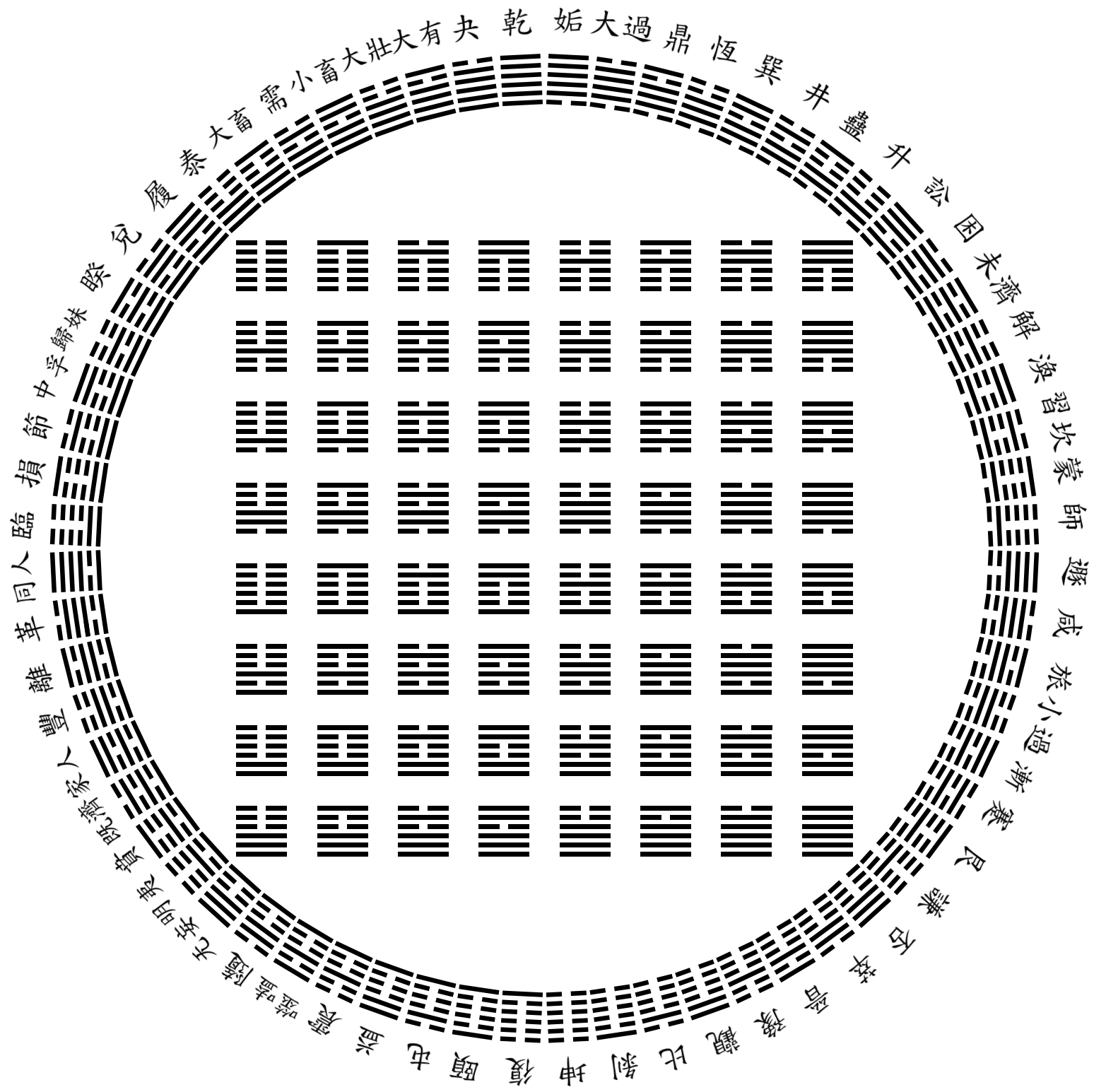

The $2^6 = 64$ Chinese hexagrams consist of six horizontal bars which can be either connected or disconnected. They appeared in the *I Ching* — the book of Changes — written more than 2500 years ago and are commonly used as a divination tool. Originally, they were ordered in a mysterious way, usually attributed to King Wen, that scholars are still trying to decipher. One thousand years ago, Shao Yong ordered them as shown in the picture, in a circle and in a square. In 1701, the jesuit Joachim Bouvet sent a copy of this configuration to Leibniz who explained it in terms of binary expansions and wrote one of the first systematic expositions of arithmetics in base 2. This is an interesting example of interaction between eastern philosophy and western science. I will discuss these I Chings a little bit more in the final section of this chapter. ©

# *Analytic chord diagrams: interlace graphs*

One of the pleasant aspects of random promenades is that they are full of surprises. Christopher-Lloyd Simon is an undergraduate student at *École Normale Supérieure de Lyon* and he kindly agreed to read the first draft of this book. While he was reading a preliminary version of the previous chapter he had the brilliant idea to transfer the discussion from chord diagrams to their associated *interlace graphs*. We already met this concept in our study of Gauss's words associated to generic immersed curves in the plane. Given a chord diagram, the set of vertices of its interlace graph is simply the set of chords, and edges connect linked (i.e. intersecting) chords. Not every graph comes from a diagram and a graph might come from several diagrams. Nevertheless, the interlace graphs coming from analytic diagrams turned out to be easy to analyze. The icing on the cake is that these graphs have been introduced forty years ago in a totally different context and are very well understood. Thanks to this new perspective, we will get the complete list of basic non-analytic chord diagrams.

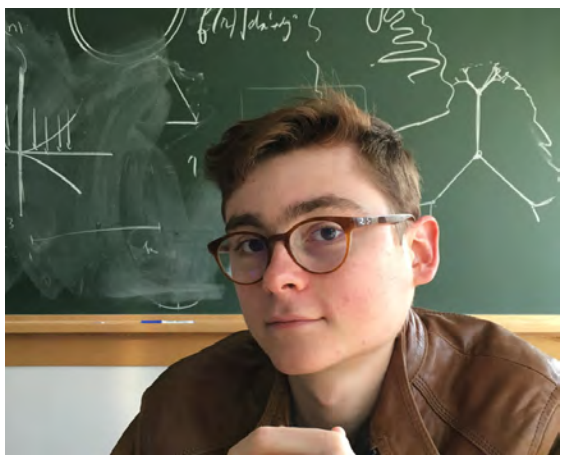

Christopher-Lloyd Simon.

## *Back to separable permutations*

In order to motivate what follows, let us revisit quickly the much simpler situation of polynomial interchanges (alias separable permutations) that we examined in the first chapters.

Let $\pi$ be a permutation of $\{1,\ldots,n\}$. The *permutation graph* $G(\pi)$ associated to $\pi$ has $\{1,\ldots,n\}$ as vertices, and an edge connects $i$ and $j$ if $\pi$ reverses the order of $(i,j)$. This is also the interlace graph of the associated marked chord diagram, with $n$ letters on each side.

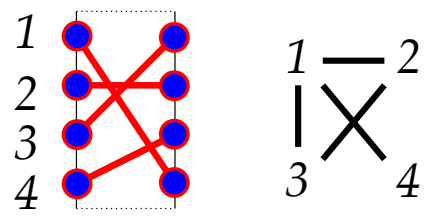

A permutation and its graph.

We know that if $\pi$ a polynomial interchange, at least two consecutive integers have consecutive images. The corresponding chords have therefore the property that any chord intersecting one of them intersects the other.

In terms of the graph $G(\pi)$, this suggests the following definition. Two vertices $x,y$ in a graph are called *twins* if they have the same neighbors (different from $x$ or $y$). They are called *true or false twins* depending on the existence of an edge connecting them. Two twins in a graph can be merged in a single vertex, producing a smaller graph with one less vertex.

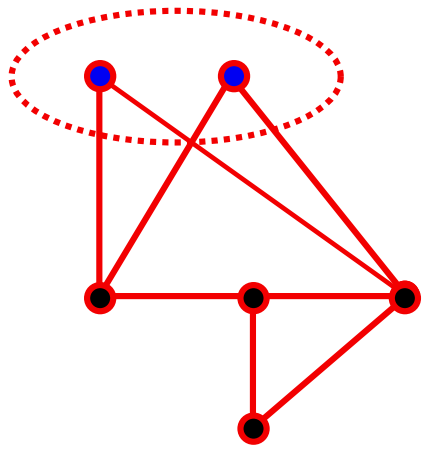

False twins. Observe that two isolated vertices are false twins. I am not responsible for the terminology which is classical and more convenient than dizygotic (or fraternal) twins. Il is also closer to the French "faux jumeaux".

Therefore the graph $G(\pi)$ coming from a polynomial interchange $\pi$ contains at least two twins, corresponding to two consecutive integers $i, i+1$ such that $\pi(i+1) = \pi(i) \pm 1$. Merging the twins in the graph amounts to merging the two elements $i, i+1$. Polynomial interchanges are characterized by the fact the iteration of this merging procedure eventually leads to the trivial permutation with $n = 1$.

**Definition.** A finite graph is called a *cograph* if it can be reduced to a trivial 1-vertex graph by merging twins successively.

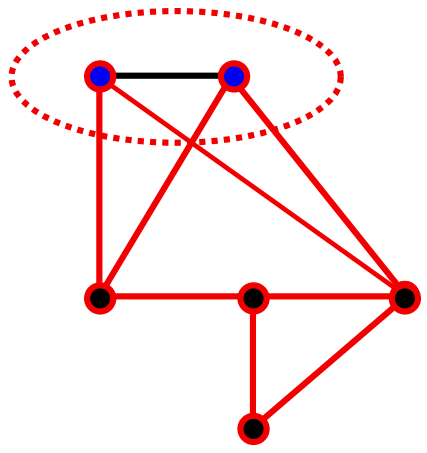

True twins, (identical or monozygotic), vrais jumeaux.

**Proposition.** *A permutation is a polynomial interchange if and only if its permutation graph is a cograph.*

I just explained why the permutation graph of a polynomial interchange is a cograph. To prove the converse, it suffices to show that if $G(\pi)$ is a cograph, there are two consecutive integers with consecutive images. The proof is by induction on $n$. If $i < j$ are false (resp. true) twins, the image by $\pi$ of the interval $\{i, i+1, \ldots, j\}$ is $\{\pi(i), \pi(i+1), \ldots, \pi(j)\}$ (resp. $\{\pi(j), \pi(j+1), \ldots, \pi(i)\}$). If $j \geq i+2$, the image $\pi(\{i, \ldots, j-1\})$ is also an interval and we apply the induction hypothesis to the restriction of $\pi$ to $\{i, \ldots, j-1\}$ so that one finds two consecutive integers with consecutive images. ⊡

The terminology "*co*graph" comes from the fact that the *com*plement of a cograph is also a cograph. A graph $G$ and its complement $\overline{G}$ have the same vertices and two vertices are adjacent in $\overline{G}$ if and only if they are not adjacent in $G$.

Cographs have been introduced in the 1970's under different names ($D^\star$-graphs, hereditary Dacey graphs, and 2-parity graphs: see[161] for references). They are not very difficult to describe. Let me list some of their properties, and leave the (elementary) proofs to the reader.

[161] A. Brandstädt, V. B. Le, and J. P. Spinrad. *Graph classes: a survey*. SIAM Monographs on Discrete Mathematics and Applications. Society for Industrial and Applied Mathematics (SIAM), Philadelphia, PA, 1999.

In what follows, all graphs are finite, with no loops and no multiple edges. A connected graph defines a metric space on its set of vertices. The distance between two vertices is, by definition, the length of the shortest path connecting them.

Graph theorists say that a subgraph $H$ of a graph $G$ is *induced* if any edge of $G$ connecting two vertices of $H$ is also an edge of $H$.

**Theorem.** *The following properties of a finite graph G are equivalent.*

1. *G is a cograph.*
2. *G is the permutation graph of some polynomial interchange.*
3. *Any two vertices in the same connected component of G are connected by a path of length at most* 2.
4. *There is no induced subgraph $P_4$ with four vertices as in the margin.*

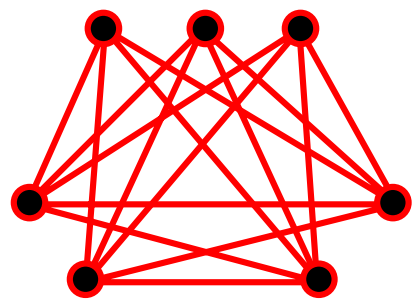

A cograph.

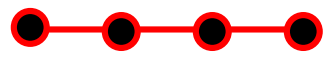

$P_4$

Note that the permutation graphs associated to the two forbidden Kontsevich's permutations $(2, 4, 1, 3)$ and $(3, 1, 4, 2)$ are both isomorphic to $P_4$.

## *Collapsible graphs*

Starting from a tree, you can strip off its leaves and do it again until the tree has been stripped completely naked. Let us say that a vertex in a graph is *pendant* if it is adjacent to a unique vertex. Any tree can be constructed by successive additions of pendant vertices, starting with the tree with only one vertex.

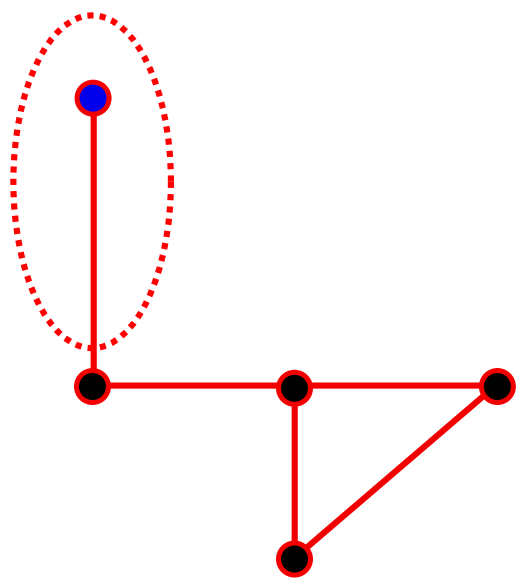

Pendant vertex.

**Definition.** A finite graph is *collapsible* if it can be reduced to a 1-vertex graph by applying two kinds of elementary operations: deleting a pendant vertex and merging twins.

If you are more constructive than destructive, you can express the same thing in another way. Start with the trivial graph with one vertex and apply two kinds of operations: adding a pendant vertex or creating a pair of twins. The second operation simply consists in duplicating a vertex and connecting the newly born twin to the rest of the graph as the original vertex was. Then, decide if you want true or false twins.

The key point is the following.

**Proposition.** *A chord diagram is analytic if and only if its interlace graph is collapsible.*

This will follow from the algorithmic description of analytic diagrams given in the previous chapter.

Before the proof, let me make an elementary remark, as an appetizer.

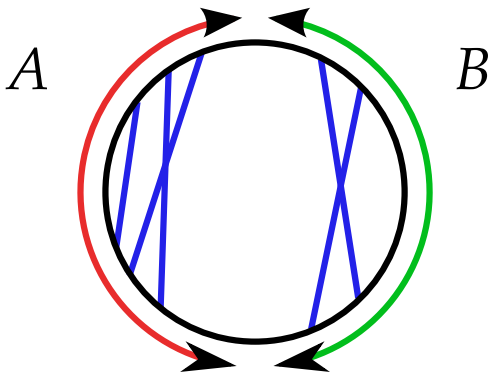

Let $w$ be a diagram and $A$ be a subset of its $2n$ letters on the circle. I will say that $A$ is *stable* under $w$ if any chord with one end in $A$ has its other end in $A$. Said differently $A$ is a sub-chord diagram $w_A$ of $w$.

Suppose that there is an *interval* $A$ which is stable under $w$, and let $B$ be its complement. Clearly the interlace graph $G(w)$ of $w$ is the disjoint union of the graphs $G(w_A)$ and $G(w_B)$ of $w_A$ and $w_B$. It follows that $G(w)$ is collapsible if and only if $G(w_A)$ and $G(w_B)$ are collapsible. Our algorithm shows that if $w_A, w_B$ are analytic so is $w$. Conversely, if $w$ is analytic, so are their sub-diagrams $w_A$ and $w_B$.

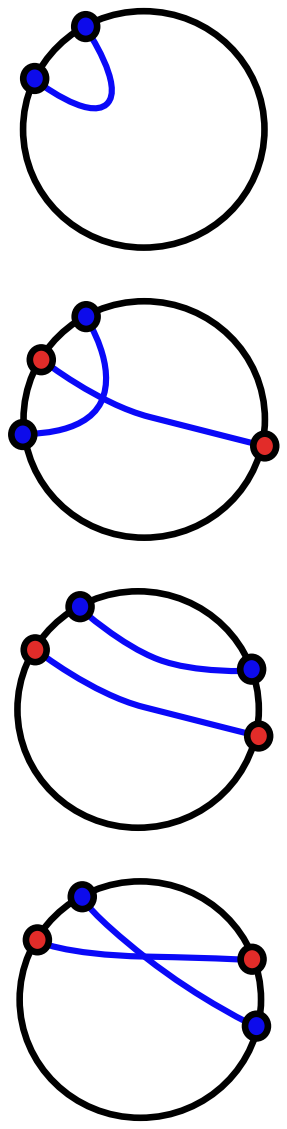

Let us prove now the proposition.

*Start with an analytic diagram $w$.* If two chords of $w$ are parallel or antiparallel, the associated vertices in the interlace graph are twins and our algorithm merges them. A pitchfork gives a pendant vertex in the graph and the algorithm deletes the short chord and keeps the handle. A solitary chord defines an isolated vertex, which is removed by the algorithm. It follows that the interlace graph associated to an analytic diagram is collapsible.

*For the converse, we show that every diagram $w$ whose interlace graph $G(w)$ is collapsible contains a solitary chord, or a pitchfork, or a pair of parallel or antiparallel chords.*

A solitary chord in a diagram corresponds to two consecutive identical letters $\cdots aa \cdots$ in the cyclic word.

A pitchfork corresponds to a subword of the form $\cdots aba \cdots$.

A pair of parallel (resp. antiparallel) chords corresponds to $\cdots ab \cdots ba \cdots$ (resp. $\cdots ab \cdots ab \cdots$).

Our proof will be by contradiction. Consider a possible counterexample $w$ to the previous assertion *with a minimal number of chords*. So $G(w)$ is collapsible and $w$ contains no solitary chord, no pair of parallel or antiparallel chords, and no pitchfork.

Since $G(w)$ is collapsible, there is a vertex $\alpha$ which is either isolated, or pendant, or is part of a pair of twins. Let $\overline{w}$ be the diagram obtained by deleting $\alpha$ from $w$.

Of course $G(\overline{w})$ is collapsible so that, by minimality, $\overline{w}$ contains a subword $\cdots aa \cdots$ or $\cdots aba \cdots$, or $\cdots ab \cdots ba \cdots$, or $\cdots ab \cdots ab \cdots$. The problem is that these are subwords of $\overline{w}$ and not of $w$, which also contains two copies of the letter $\alpha$, which could sneak into the above subwords.

Note that by minimality any interval which is stable under $w$ is either empty or everything.

A priori:

– $0, 1$ or $2$ letters $\alpha$ could sneak in the subword,

– the subword of $\overline{w}$ could correspond to a solitary chord, or a pitchfork, or to a pair of parallel or antiparallel chords,

– $\alpha$ could be isolated, pendant or twin, true or false, in $G(w)$.

That makes $3 \times 4 \times 4$ cases to examine! Fortunately, many cases can be studied simultaneously.

This case by case proof is not particularly pleasant. You can skip it if you wish, but if you do so, you should sympathize with me, who had to list all cases one by one.

1/ *If no letter $\alpha$* sneaks into the above subwords, there is no problem: our solitary chord or pitchfork, or pair of parallel or antiparallel chords in $\overline{w}$ have the same property for $w$ $\lightning$.

I use the symbol $\lightning$ to mean "contradiction".

2/ If $\alpha$ is isolated in $G(w)$, this means that no chord intersects $\alpha$. Therefore $\alpha$ decomposes the circle in two stable intervals, which have to be empty $\lightning$.

3/ *If two letters $\alpha$* sneak in, they cannot occur as consecutive letters since that would force the chord $\alpha$ to be solitary in $w$ $\lightning$.

So, still in this case, we have to look at

$$\cdots a\alpha b\alpha a\cdots, \text{or} \cdots a\alpha b\cdots b\alpha a\cdots, \text{or} \cdots a\alpha b\cdots a\alpha b\cdots$$

This produces respectively a pitchfork $(\alpha, b)$, a pair of parallel chords $(\alpha, b)$ or antiparallel chords $(\alpha, b)$ in $w$ ↯.

4/ *Inserting one* $\alpha$ *in a solitary chord* yields $\cdots a\alpha a\cdots$ which produces a pitchfork in $w$ with handle $\alpha$ ↯.

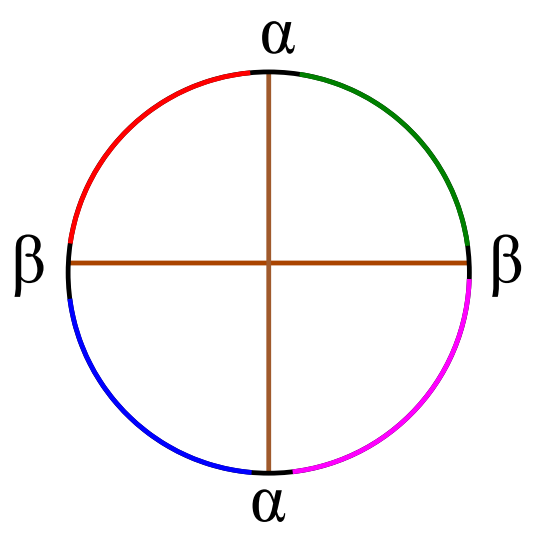

So far, we did not use the fact that $\alpha$ is pendant or twin. This will be used in the remaining cases, when a single $\alpha$ enters in a pitchfork or a pair of parallel or antiparallel chords of $\overline{w}$.

If $\alpha$ is pendant, let $\beta$ be the only chord in $w$ which intersects $\alpha$. If $\alpha$ has some twin siblings, let us denote one of them by $\beta$. The two chords $\alpha, \beta$ determine four intervals in the circle, *excluding* $\alpha, \beta$, that I will call *sectors*. If $\alpha$ is pendant, the union of two sectors which are on the same side of $\alpha$ is stable. If $\alpha, \beta$ are twins, the union of opposite sectors is stable.

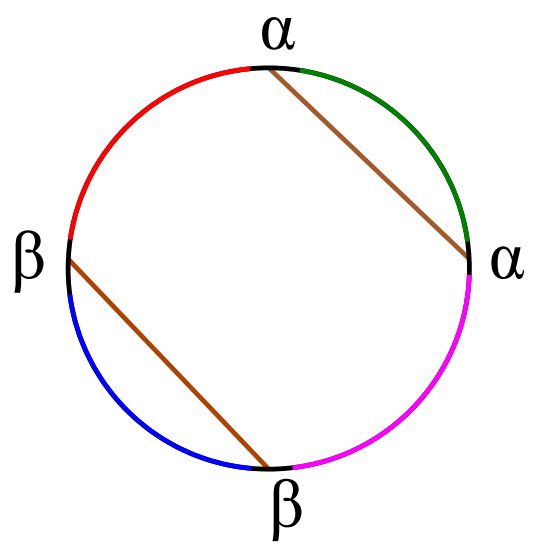

5/ *Suppose now that a single letter* $\alpha$ *enters* $\cdots aba\cdots$, *or* $\cdots ab\cdots ba\cdots$, *or* $\cdots ab\cdots ab\cdots$ *and that one of the two letters* $a, b$ *is equal to* $\beta$. Then, the two letters $\alpha$ and $\beta$ are consecutive in $w$. This implies that one of the sectors is empty.

5-1/ In the pendant case, this implies that the other sector, on the same side of $\alpha$, is a stable interval and therefore also empty. So, $\alpha, \beta$ is a pitchfork in $w$ with handle $\beta$ ↯.

5-2 In the twin case, this implies that the opposite sector is a stable interval, and therefore empty. So $\alpha, \beta$ is a pair of parallel or antiparallel chords in $w$ ↯.

6/ Finally, *suppose that a single letter* $\alpha$ *enters* $\cdots aba\cdots$, *or* $\cdots ab\cdots ba\cdots$, *or* $\cdots ab\cdots ab\cdots$ *and that none of the two letters* $a, b$ *is equal to* $\beta$.

6-1 Assume that $(a, b)$ are parallel or antiparallel chords in $\overline{w}$. Inserting one $\alpha$ in $ab$ yields $\cdots a\alpha b\cdots ba\cdots$ or $\cdots a\alpha b\cdots ab\cdots$. From the consecutive letters $ba$ (or $ab$), it follows that these two occurrences are on the same side of the chord $\alpha$. From the consecutive letters $a\alpha b$ in $w$, it follows that these other two occurrences of $a$ and $b$ are in different sides of $\alpha$. Hence the chord $\alpha$ intersects only one of the chords $a, b$.

6-1-1 If $\alpha$ is pendant, this forces $a$ or $b$ to be equal to $\beta$ ↯.

6-1-2 If $\alpha, \beta$ are twins, this is not possible $\lightning$.

6-2 Assume that $(a, b)$ is a pitchfork in $\overline{w}$. Inserting a letter $\alpha$ in $\cdots aba \cdots$ yields $\cdots a\alpha ba$ or $\cdots ab\alpha a$ so that the chord $\alpha$ should intersect the chord $a$.

6-2-1 If $\alpha$ is pendant, this forces $a = \beta$ $\lightning$.

6-2-2 If $(\alpha, \beta)$ are twins, the consecutive letters $a\alpha ba$ show that the two letters $a$ lie in opposite sides of $\alpha$, hence on opposite sectors. It follows that $b = \beta$ $\lightning$.

This finishes the proof. Ouf ! ⊡

Now, we have to understand the nature of collapsible graphs.

## *Collapsible, distance hereditary*

Collapsible graphs have been defined by several authors forty years ago, under different names, with very different motivations. We will see that these graphs are very close to being trees.

Howorka[162] defined *distance hereditary graphs* in 1977.

[162] E. Howorka. A characterization of distance-hereditary graphs. *Quart. J. Math. Oxford Ser. (2)*, 28(112):417–420, 1977.

**Definition.** A finite graph $G$ is *distance hereditary* if for every connected induced subgraph $H \subset G$, the distance between two vertices of $H$ *in* $H$ is equal to the distance between the same vertices *in* $G$.

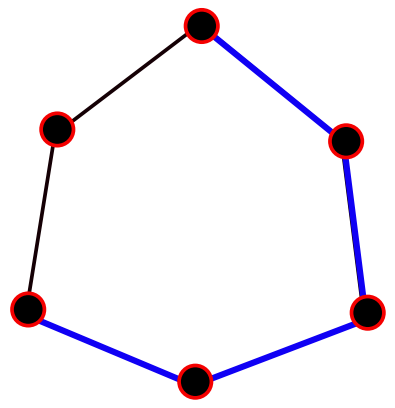

A cycle of length ≥ 5 is not distance hereditary.

For instance, a tree is *distance hereditary* and a cycle of length at least 5 is not. It suffices to choose $H$ as the induced subgraph defined by a path inside the cycle whose length is greater than one half of the length of the cycle (in blue on the picture).

Consider a finite graph and choose some length for each edge, which could be any positive real number. Define the length of a path as the sum of the lengths of its edges and the distance between two vertices as the smallest length of a path connecting them. One speaks of a *metric graph*.

We are looking for a characterization of metric spaces (usually called *metric trees*) arising in this way from *trees*. Here is the answer. Let $(E, d)$ be a finite metric space. Choose four points $x_1, x_2, x_3, x_4$ in $E$ and compute the sums of the lengths of the three pairs of diagonals:

$$d(x_1,x_2)+d(x_3,x_4)\,;\,d(x_1,x_3)+d(x_2,x_4)\,;\,d(x_1,x_4)+d(x_2,x_3).$$

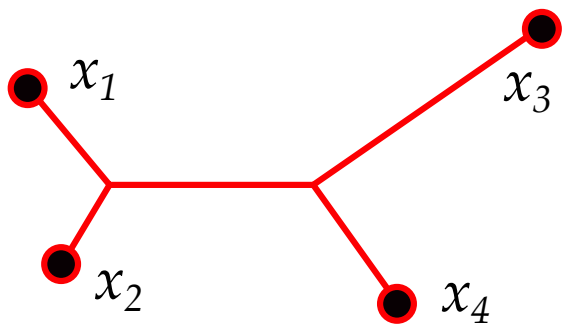

Let $s$ (resp. $m$, $l$) be the smallest (resp. medium, largest) of these three numbers: $s \le m \le l$. It turns out that a finite metric space is isometric to a subset of a metric tree if and only if $m = l$ for every quadruple of points. This is not difficult to prove and I leave it as an exercise $M25$. The lazy reader might see the proof in this short paper[163].

[163] P. Buneman. A note on the metric properties of trees. *J. Combinatorial Theory Ser. B*, 17:48–50, 1974.

We should be careful. A graph, where all edges have length 1, can be isometric to a subset of a metric tree without being itself a tree. Look at the example in the margin.

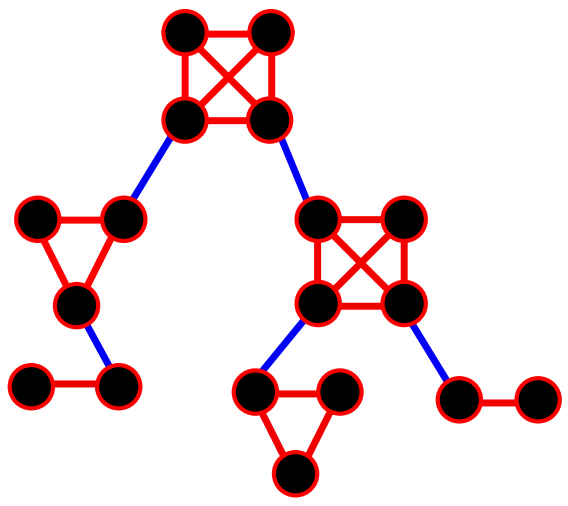

A block graph.

In graph theory, those graphs are called *block graphs*. In order to construct them, start with a tree, delete some of its vertices and replace them by *cliques*, i.e. finite graphs where all pairs of vertices are adjacent, as in the figure. I suggest that my reader proves that this is indeed a characterization of block graphs ($M15$ and, in case of emergency, see[164]).

[164] F. Harary. A characterization of block-graphs. *Canad. Math. Bull.*, 6:1–6, 1963.

In the 1980's, Gromov developed a geometric theory for *hyperbolic spaces* which had a very strong influence on combinatorial and geometric group theory (unfortunately not part of our promenade).

The definition is the following. A metric space $(E,d)$ is called *hyperbolic* if there exists some $\delta \ge 0$ such that for every quadruple of points as above, $m$ and $l$ are "almost equal", i.e. $l - m \le \delta$. Note that any finite metric space is trivially hyperbolic (for $\delta$ sufficiently big) so that this concept is only relevant for geometry *in the large*.

A metric space $(E,d)$ is *geodesic* if for every pair of points $(x,y)$ there exists an *isometric* embedding $i : [0, d(x,y)] \to E$ such that $i(0)) = x$ and $i(d(x,y)) = y$.

There are many equivalent formulations of this property, the most popular (for geodesic metric spaces) being that all geodesic triangles are *slim*. Consider three points $x, y, z$ and choose three geodesics $[x,y], [x,z], [y,z]$ connecting them. Every point in $[x,y]$ should be at some uniformly bounded distance from the union $[x,z] \cup [y,z]$, independently of the choice of $x, y, z$ (see the picture).

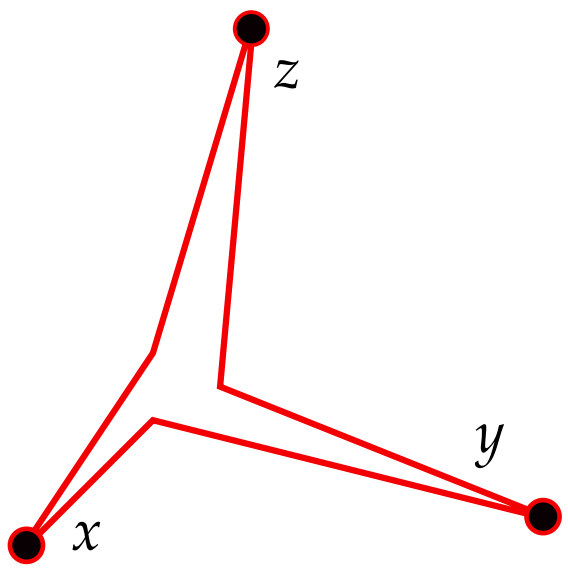

This concept is remarkably robust. For instance, the universal cover of a negatively curved compact Riemannian manifold is hyperbolic. These metric spaces are well approximated by trees, in a quantitative way. For more about this theory, the reader is encouraged to read[165].

In 1986, Bandelt and Mulder published a paper[166] proposing purely metrical characterizations of distance hereditary graphs, close to Gromov's hyperbolicity conditions.

**Definition.** A finite graph $G$ is *treelike* if for every 4-tuple of vertices $x_1, x_2, x_3, x_4$ two of the following three numbers are equal:

$$d(x_1,x_2)+d(x_3,x_4); \quad d(x_1,x_3)+d(x_2,x_4); \quad d(x_1,x_4)+d(x_2,x_3).$$

My reader has probably guessed that all these definitions turn out to be equivalent.

**Theorem.** *Let $G$ be a finite graph. The following properties are equivalent.*

1. *$G$ is the interlace graph of some analytic chord diagram,*
2. *$G$ is collapsible,*
3. *$G$ is distance hereditary,*
4. *$G$ is treelike,*
5. *$G$ does not contain a cycle of length at least five, or a house, a gem, or a domino, as an induced subgraph.*

The *house*, the *gem* and the *domino* are pictured below.

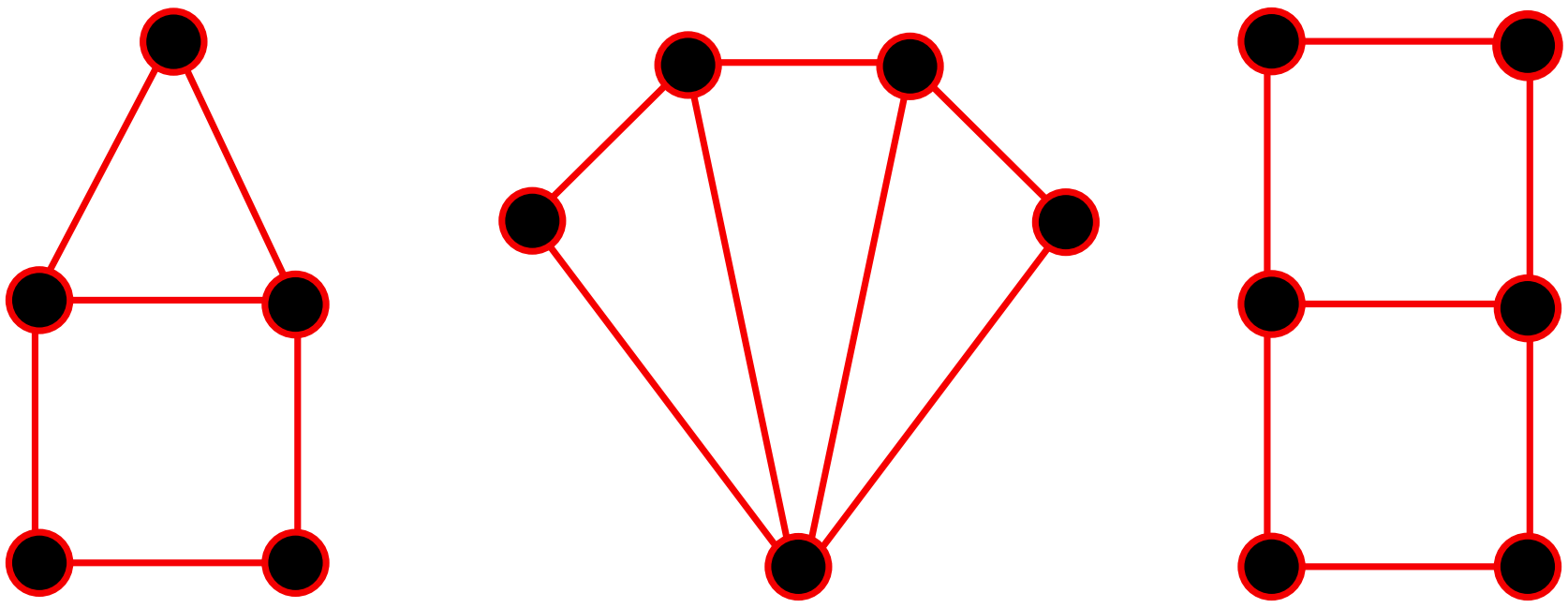

[165] É. Ghys and P. de la Harpe, editors. *Sur les groupes hyperboliques d'après Mikhael Gromov*, volume 83 of *Progress in Mathematics*. Birkhäuser Boston, Inc., Boston, MA, 1990.

[166] H.-J. Bandelt and H. M. Mulder. Distance-hereditary graphs. *J. Combin. Theory Ser. B*, 41(2):182–208, 1986.

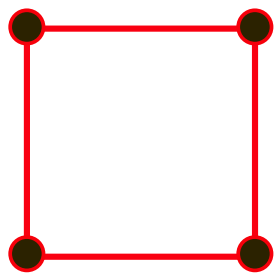

A cycle of length 4 is treelike but not a tree.

*Exercise*: Show that a treelike graph is hyperbolic in the sense of Gromov with $\delta$=2.

All the equivalences in the previous theorem (except of course the first item) are proved in the above mentioned papers. However, I will soon propose some elementary proofs.

It is now time to harvest the fruits of our labor and to get a very simple description of analytic chord diagrams.

You should not be surprised that the interlace graphs of , , are the house, the gem, and the domino.

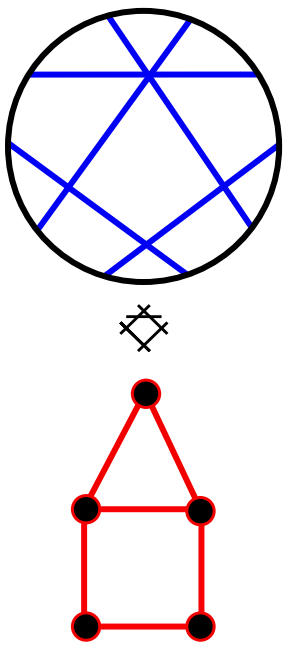

The house.

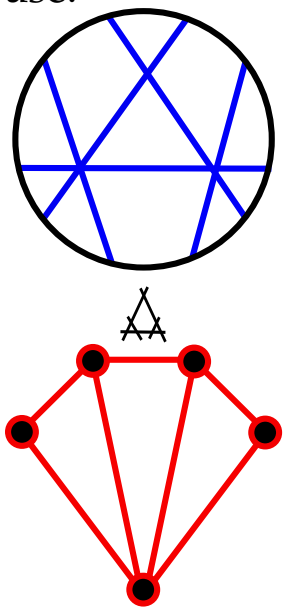

The gem.

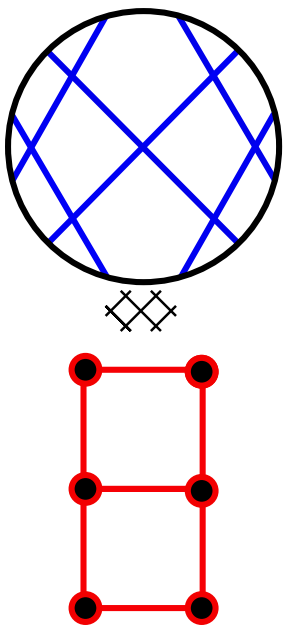

The domino.

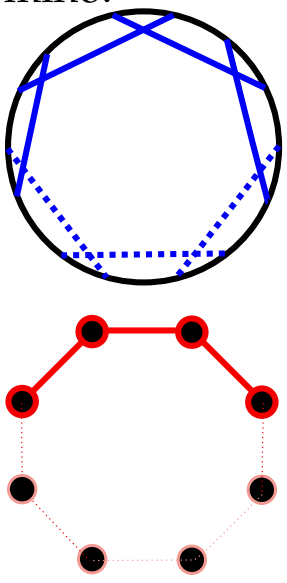

$\mathcal{C}_n$ and its interlace graph: the $n$-cycle.

*Exercise :* Show that , , are the only chord diagrams whose interlace graphs are the house, the domino and the gem.

In the same way, we have already described the non-analytic chord diagram $\mathcal{C}_n$ defined by $\mathbb{Z}/2n\mathbb{Z}$ ($n \geq 5$) where there is a chord connecting $2k$ and $2k+3$ (for $k = 1, \ldots, n$). Its interlace graph is a cycle of length $n$.

*Exercise :* Show that $\mathcal{C}_n$ is the only chord diagram whose interlace graph is a cycle of length $n$.

Finally, note that a sub-chord diagram defines an induced subgraph in the interlace graph. Therefore, we get a very satisfactory description of analytic chord diagrams. I print the following theorem in blue since it is a highlight in our promenade.

**Theorem.** *A chord diagram is analytic if and only if it does not contain , , or $C_n$ ($n \geq 5$) as a sub-chord diagram.*

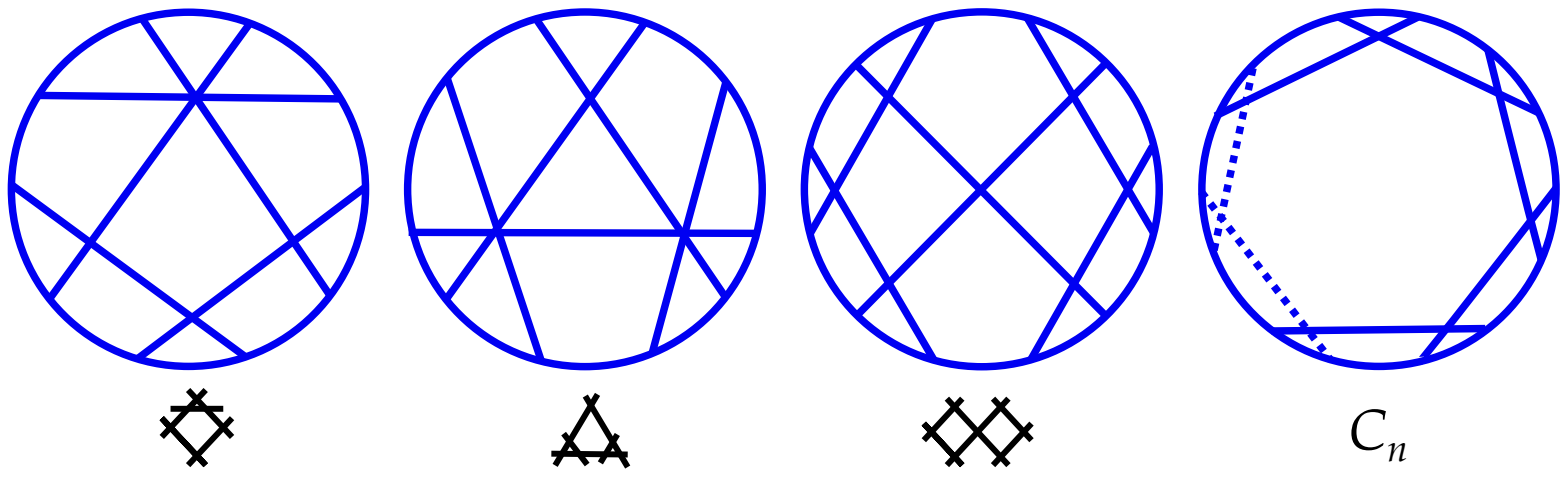

Note the complete analogy with our characterization of *polynomial interchanges* as the *separable permutations*, which are in turn precisely those permutations which do not contain Kontsevich's examples (2,4,1,3) and (3,1,4,2).

## *Some proofs*

I present now the proofs of the equivalences of the definitions in the previous section. They are mostly elementary and I suggest that the interested reader tries to prove them alone. It is important to draw pictures. In this specific case, it was probably more challenging to find the significant definitions than to prove their equivalence.

*No* 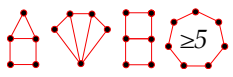 *is an induced subgraph* $\Longrightarrow$ *Distance hereditary.*

Let $H$ be an induced subgraph of a graph $G$. Connect two vertices $p, q$ of $H$ at distance $n$ *in* $H$ by a path $c = (x_0, x_1, \ldots, x_n)$ (with $p = x_0$ and $q = x_n$) in $H$. Two vertices $x_i, x_j$ are adjacent if and only if $i, j$ are consecutive since otherwise there would be a shortcut. In other words, the path $c$ is induced in $G$. It follows that in order to show that a graph is distance hereditary, we should prove that the distance between the endpoints of any induced path in $G$ is equal to the length of the path.

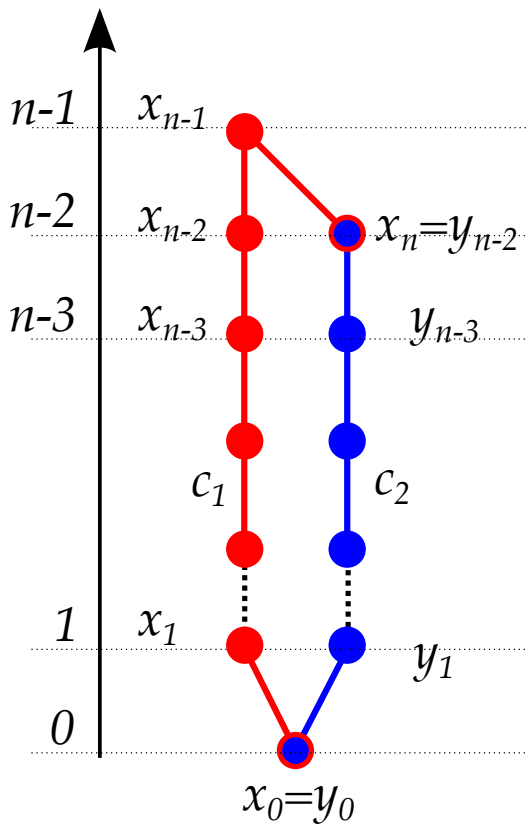

Suppose that no 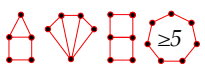 is induced in $G$. Choose an induced path $c_1 = (x_0, x_1, \ldots, x_n)$ and let us show, by induction on $n$, that the distance in $G$ between $x_0$ and $x_n$ is exactly $n$.

Connect $x_0$ to $x_n$ by a shortest path $c_2 = (y_0, y_1, \ldots, y_l)$ in $G$ (with $y_0 = x_0$ and $y_l = x_n$). Of course, $c_2$ is also induced and $d(y_0, y_i) = j$ for $0 \le i \le l$. By induction, $d(x_0, x_i) = i$ for $0 \le i \le n - 1$. It follows that $l$ is equal to $n-2$, $n-1$ or $n$. We show that the first two cases are not possible. Suppose that $l = n-2$ or $n-1$. By induction, we can assume that the two paths $c_1, c_2$ only intersect at their endpoints: any other intersection point could be used as the starting point of shorter paths $c_1'$ and $c_2'$.

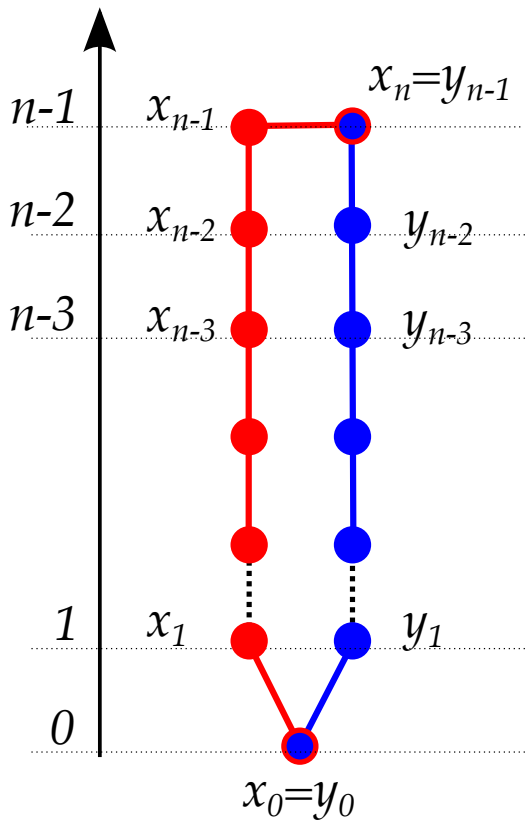

Draw a picture in the plane in such a way that the height of a vertex of $c_1$ or $c_2$ is the distance from $x_0$. The cases $l = n-2$ and $l = n-1$ are pictured in the margin. Vertices of $c_1$ are red and vertices of $c_2$ are blue. The union of $c_1$ and $c_2$ defines a cycle $c$ in $G$. The length of $c$ is at least 5. The cycle $c$ cannot be induced since there is no induced cycle of length $\ge 5$. Therefore there must exist diagonals connecting vertices of $c_1$ with vertices of $c_2$.

By the triangle inequality, the height difference of the two endpoints of a diagonal can only be $-1, 0, 1$. Moreover, diagonals connect points of different colors.

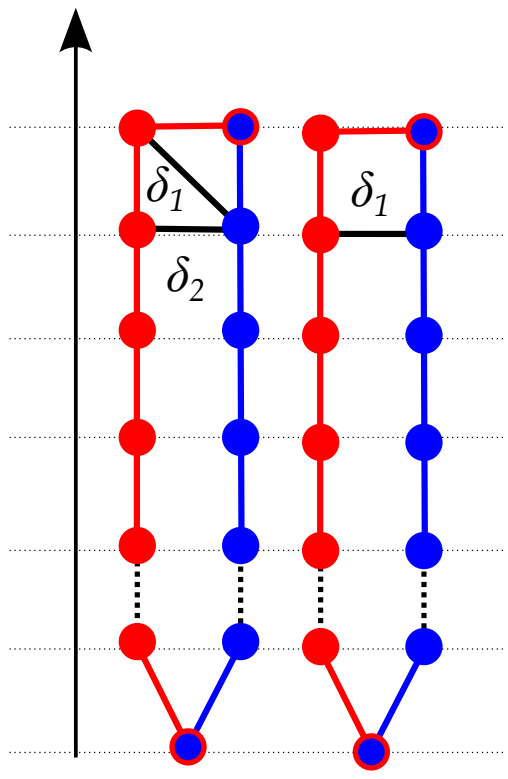

Let us order the diagonals $(x_i, y_j)$ from top to bottom, i.e. $(x_i, y_j)$ is before $(x_{i'}, y_{j'})$ if $j > j'$ or $j = j'$ and $i > i'$.

Now, try to construct the ladder, one diagonal at a time. The rule of the game is the following. You have to draw an ordered sequence of diagonals $\delta_1, \delta_2, \ldots$ respecting the conditions above and without creating any induced 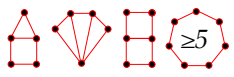. Note that the diagonal $\delta_k$ together with the part of $c$ which is above it defines a cycle. Any chord in this cycle has to be one of the previously chosen chords $\delta_1, \ldots, \delta_{k-1}$.

In the case $l = n-1$, there are only two possibilities for the first diagonal $\delta_1$. It could be $(x_{n-1}, y_{n-2})$ or $(x_{n-2}, y_{n-2})$. In the case $l = n-2$ there is only one possibility for $\delta_1$.

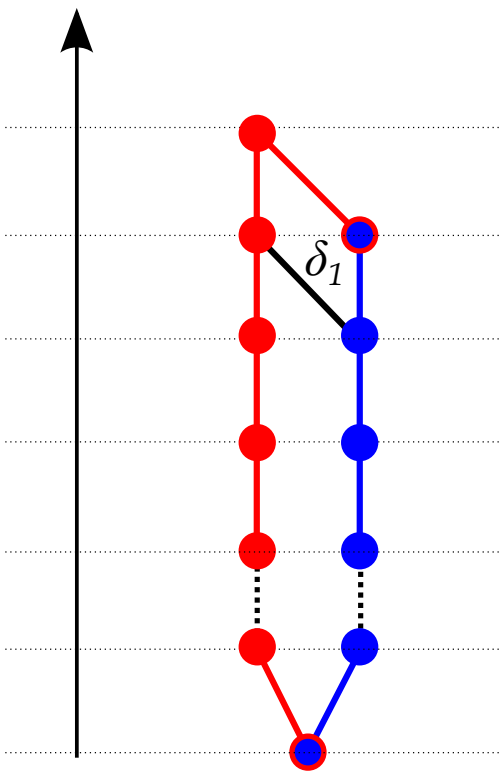

Then, try to select the second diagonal $\delta_2$, avoiding 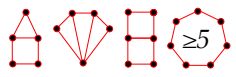. Only one of the three choices of $\delta_1$ allows you to do so.

Finally, try to draw the third diagonal, in the only case where you could draw $\delta_1, \delta_2$. It is not possible to continue without creating one of the forbidden graphs. ⊡

*Distance hereditary* $\Longrightarrow$ *Collapsible.*

Let me first make an elementary remark.

Choose some vertex $x$ in a connected distance hereditary graph $G$ and look at the largest $k$ such that the sphere $S_k$ in $G$ of radius $k$ and centered on $x$ is non-empty. Let $C$ be a connected component of $S_k$. If $C$ contains only one element, then it is a pendant vertex in $G$.

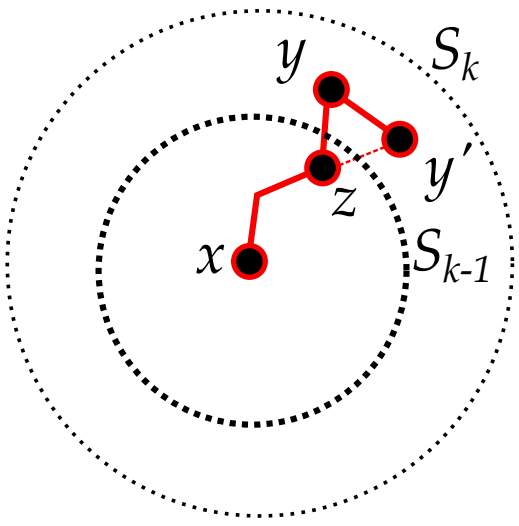

Choose two vertices $y, y'$ in $C$ which are adjacent in $C$. Choose some point $z$ which is adjacent to $y$, at distance $k-1$ from $x$. Choose a chain $c$ of length $k-1$ from $x$ to $z$ and call $c'$ the chain of length $k+1$ obtained by adding the edge between $z$ and $y$ and from $y$ to $y'$. Since the distance between $x$ and $y'$ is exactly $k$, this chain cannot be induced and $y'$ has to be adjacent to $z$. This implies that two points in $C$ are simultaneously adjacent or not to any point $z$ at distance $k-1$ from $x$. It follows that $C$ cannot contain an induced path $P_4$ of length 3, since together with $z$, it would produce a gem in $G$, which is not distance hereditary. Hence $C$ is a cograph and in particular contains a pair of twins. By the above observation, two twins in $C$ are twins in $G$. ⊡

*Collapsible* $\Longrightarrow$ *Treelike.*

Easy by induction. Take four points in a graph, delete a pendant vertex or merge two twins. One of the four points might be the vertex which has been removed. If this is the case replace it by the other end of the removed edge. Look at the corresponding points in the stripped graph (taking into account for instance the fact that two of our four points could be the two twins which have been merged). Apply the induction hypothesis. ⊡

*Treelike* $\Longrightarrow$ *No induced* 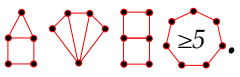.

Obvious since one checks easily that none of these examples of graphs are treelike. ⊡

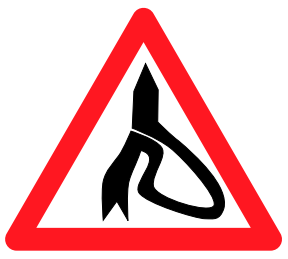

Another detour. The only purpose of this section is to describe the structure of collapsible graphs.

## *Completely decomposable graphs*

Graphs that can be stripped to a point by deleting only pendant vertices are trees. Graphs that can be stripped to a point by merging only pairs of twins are cographs.

My reader has probably guessed that collapsible graphs should not be far from being trees. This is indeed true as I explain now.

I prefer the terminology "decomposable" to "separable" which is also common in this area, since we already used the word "separable" for permutations.

Let $G$ be a finite connected graph. Suppose that its vertices have been partitioned in two parts $A_1$ and $A_2$. Let $B_1 \subset A_1$ (resp. $B_2 \subset A_2$) the set of vertices of $A_1$ (resp. $A_2$) which are adjacent to some vertex in $A_2$ (resp. $A_1$). *Suppose that* every *element of* $B_1$ *is adjacent to every element of* $B_2$. This condition is trivially satisfied if $A_1$ or $A_2$ contains only one (or zero!) element, so we assume that $A_1$ and $A_2$ contain at least two elements each. In this situation, the graph $G$ is called *decomposable* and the partition $A_1, A_2$ is a *split*. In order to keep track of this decomposition, let us create two graphs $G_1, G_2$ in the following way. The set of vertices of $G_1$ (resp. $G_2$) is $A_1$ plus one extra vertex $x_1$ (resp. $x_2$) called the *control vertex*. As for the edges of $G_1$ (resp. $G_2$), choose the edges of $G$ plus extra edges connecting $x_1$ (resp. $x_2$) to all elements of $B_1$ (resp. $B_2$).

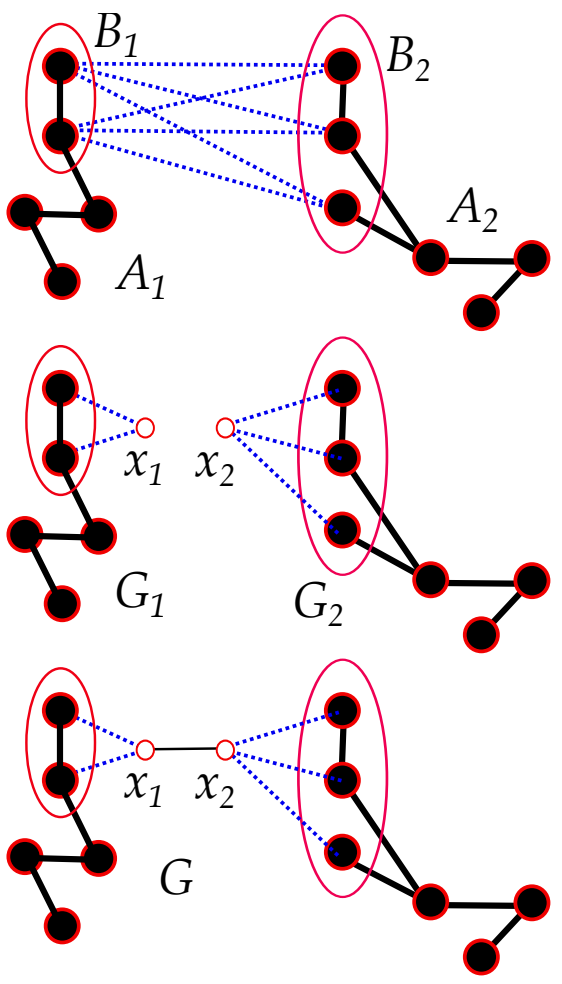

The graph $G$ can be reconstructed from $(G_1, x_1)$ and $(G_2, x_2)$ by an elementary *join* construction. Notice that the control points

$x_1, x_2$ are not vertices of $G$: they are only useful to define the edges connecting the two parts.

Note that when $A_2$ contains two elements, as in the margin, the graph $G$ has a pendant vertex or a pair of twins.

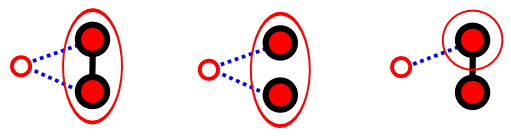

Hammer and Maffray[167] introduced another definition in 1987.

**Definition.** A finite graph is *completely decomposable* if every induced connected subgraph with at least four vertices is decomposable.

Note the analogy with trees. A finite connected graph is a tree if and only if every induced connected subgraph contains a *cut edge*, i.e. an edge that disconnects it.

[167] P. L. Hammer and F. Maffray. Completely separable graphs. *Discrete Appl. Math.*, 27(1-2):85–99, 1990. Computational algorithms, operations research and computer science (Burnaby, BC, 1987).

It is not hard to prove that completely decomposable (connected) graphs are precisely the collapsible (connected) graphs.

Indeed in the join construction, if $G_1$ and $G_2$ are collapsible, the same is true for $G$, so that completely decomposable graphs are collapsible, by induction.

Conversely, we have seen that a pendant vertex or a pair of twins gives rise to a decomposition. Therefore collapsible graphs are decomposable and even completely decomposable, since an induced subgraph of a collapsible graph is collapsible.

In order to give a precise description of completely decomposable graphs, let me state first the important *split decomposition* theorem for *general* connected finite graphs.

If a finite connected graph $G$ is decomposable, consider it as the join of $G_1$ and $G_2$ as before. Then, try to decompose $G_1$ and $G_2$ etc. until the resulting graphs become non-decomposable. The final result of this decomposition into "prime pieces" can be conveniently described by a *graph-labeled tree*, as explained below.

It consists of a tree $T$ where each internal node $x$ is equipped with a connected finite graph $G_x$. Moreover some bijection has been chosen between the vertices of $G_x$ and the edges getting out from the node $x$ in $T$. Assume that the valency of each node is at least 3. Given such a structure, we construct a graph $G(T)$ which is a "composition of the $G_x$'s controlled by $T$". The definition is the following.

The vertices of $G(T)$ are the leaves of $T$. In order to understand the edges of $G(T)$, let me just draw a picture, inspired by the paper of Gioan and Paul[168] who introduced this concept of *graph-labeled tree*. We see a tree with 16 leaves and 6 internal nodes, in pink. The associated graph, with 16 vertices, is drawn on the right.

[168] E. Gioan and C. Paul. Split decomposition and graph-labelled trees: characterizations and fully dynamic algorithms for totally decomposable graphs. *Discrete Appl. Math.*, 160(6):708–733, 2012.

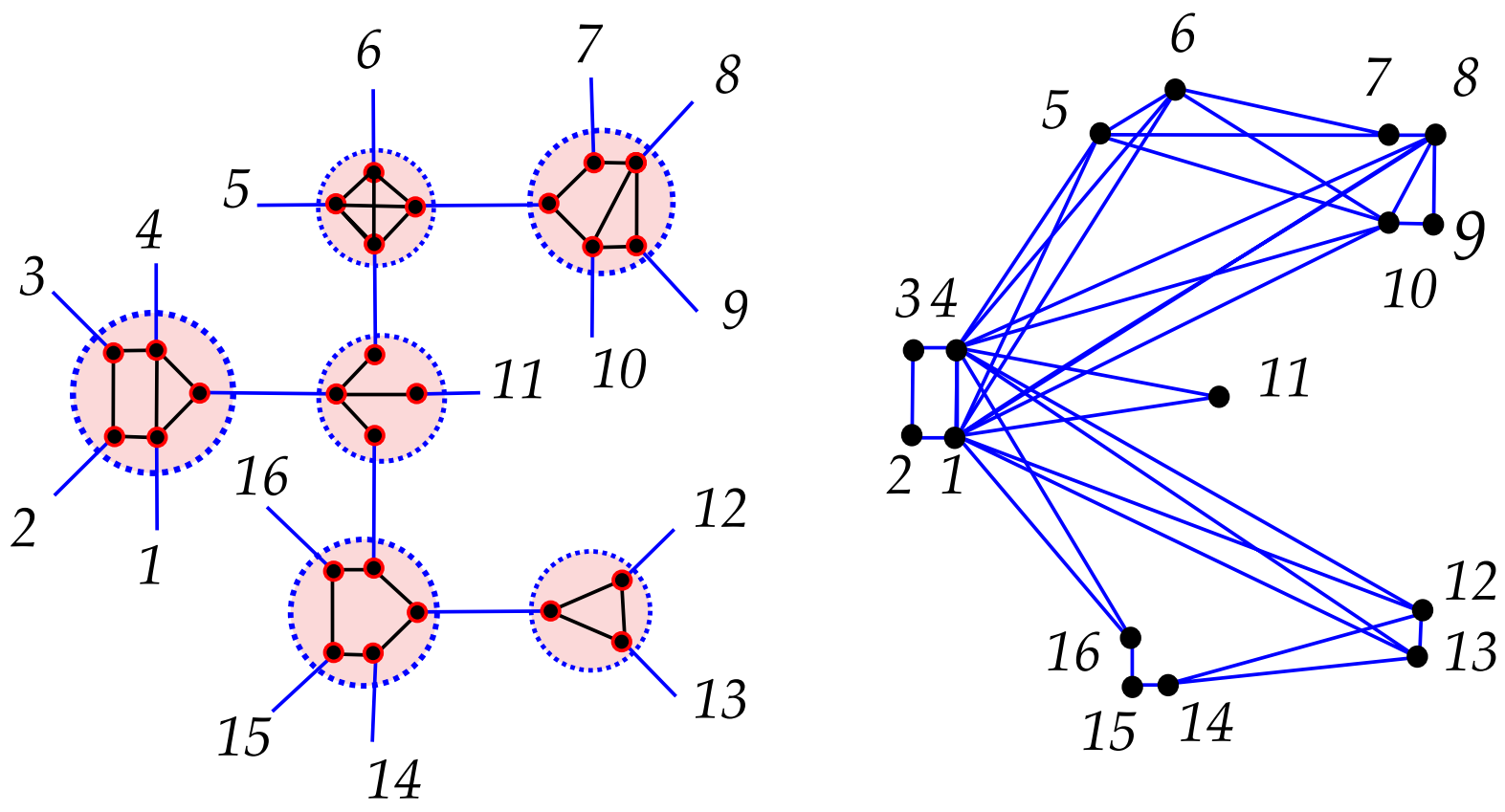

Choose two leaves of $T$ and connect them by the shortest path in the tree. For each node $x$ which is visited by this path, there is an entrance edge and an exit edge. In turn, these two edges define two vertices of $G_x$. Two vertices of $G(T)$, i.e. two leaves of $T$, are adjacent in $G(T)$ if for every node $x$ visited by this path, the two corresponding vertices of $G_x$ are adjacent in $G_x$. The vertices of the $G_x$'s generalize the two control vertices, as in the simple case where $T$ contains only one edge.

It seems that graphs can be labeled in the USA and labelled in other English speaking countries.

The main result, proved by Cunningham and Edmonds[169] in 1980 (and reformulated by Giona and Paul), is that *any finite connected graph is obtained by such a construction in which the $G_x$'s are indecomposable, in an essentially unique way*. The *existence* of this splitting is easy. The hard part is the "essential uniqueness" that I don't define since I will not need it.

[169] W. H. Cunningham and J. Edmonds. A combinatorial decomposition theory. *Canad. J. Math.*, 32(3):734–765, 1980.

Let us back to completely decomposable graphs. In this special case the $G_x$'s must *have at most* 3 *vertices*. Indeed, they are indecomposable and induced subgraphs of $G$, so that the claim follows from the definition of complete decomposability.

This gives a fairly precise geometric description of completely decomposable graphs. Take a tree such that every node has valency 3. For each node, choose a connected graph with 3 vertices (there are not too many choices!), and construct a graph-labeled tree as in the margin. All completely decomposable graphs are produced in this way.

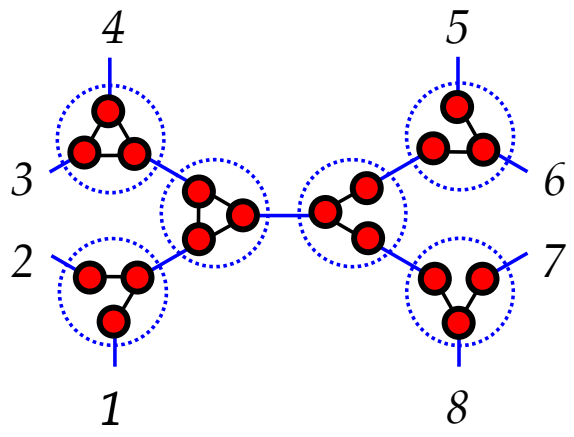

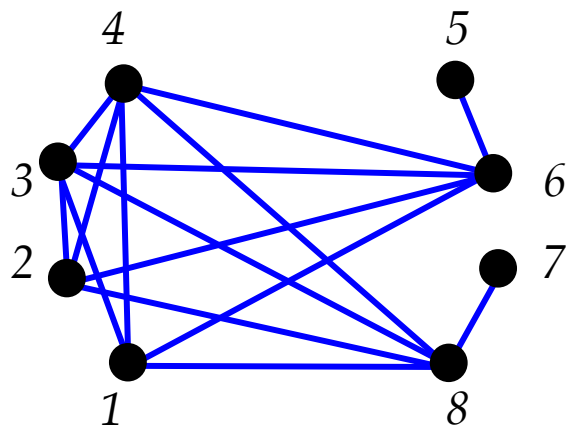

This is absolutely not a surprise. Indeed look at the third picture in the margin, showing three small graphs with three leaves. When you hook one of them to some (blue) leaf of a graph-labeled tree $T$, you get another graph-labeled tree $T'$ with one more leaf. If you examine the effect on the associated graph $G(T)$, you see that you have split a vertex in a pair of twins (true or false) or you have created a pendant vertex, depending on the cases. We are back to the original definition of collapsible graphs as graphs that can be constructed from a point by successive introductions of twins or of pendant vertices. We are also back to operads. This was indeed a detour!

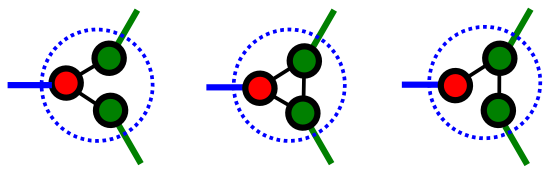

## *Computability*

There is an algorithm deciding in quadratic time (in $n$) if a graph of size $n$ is an interlace graph. This was proved in[170] after a long period of successive improvements (starting from a $n^9$ algorithm, in 1987).

[170] J. Spinrad. Recognition of circle graphs. *J. Algorithms*, 16(2):264–282, 1994.

Given a diagram with $n$ chords, constructing its interlace graph requires a time which is quadratic in $n$. Then, look for pendant vertices and twins and iterate the process $n$ times so that you decide in quadratic time if it is analytic.

## *An esoteric exercise*

The 64 hexagrams pictured on the first page of this chapter are traditionally grouped in 32 pairs of *complementary hexagrams*. Think of the *Yin and Yang*. To get the *dual of a hexagram*, just turn it upside down. In the case where the hexagram is symmetric, replace each connected line by a disconnected one and conversely.

Here are two examples of dual pairs.

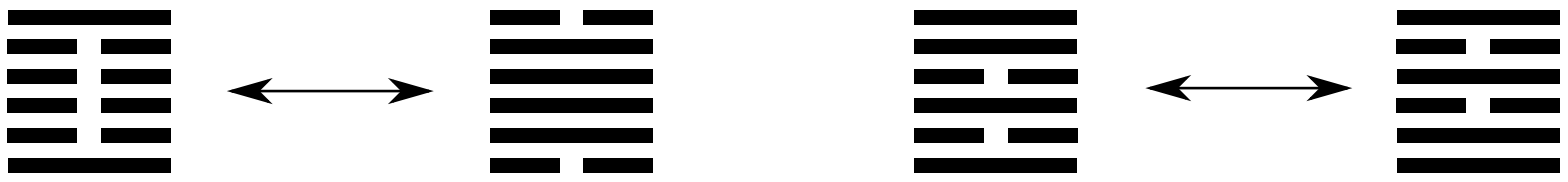

27 28 37 38

The numbers (27-28) and (37-38) are relative to the King Wen ordering. Many experts like to draw a segment between dual hexagrams. Working with Shao Yong's circular arrangement, this produces a diagram with 32 chords.

Will my reader have the patience to draw these chords and *decide whether or not this I Ching diagram is analytic?*

On November 14th 1701, Leibniz received a copy of Shao Yong's circular arrangement from the French Jesuit Joachim Bouvet who was living in China. Two years later, he published a remarkable paper on binary arithmetics[171] in the *Mémoires de l'Académie des Sciences*. According to him:

> These figures are perhaps the most ancient monument of science which exists in the world.

"Leibniz hoped that his astute analysis of the trigrams from the *I Ching* would awaken in China a deep appreciation for Western science and, ultimately, for Christianity [172]."

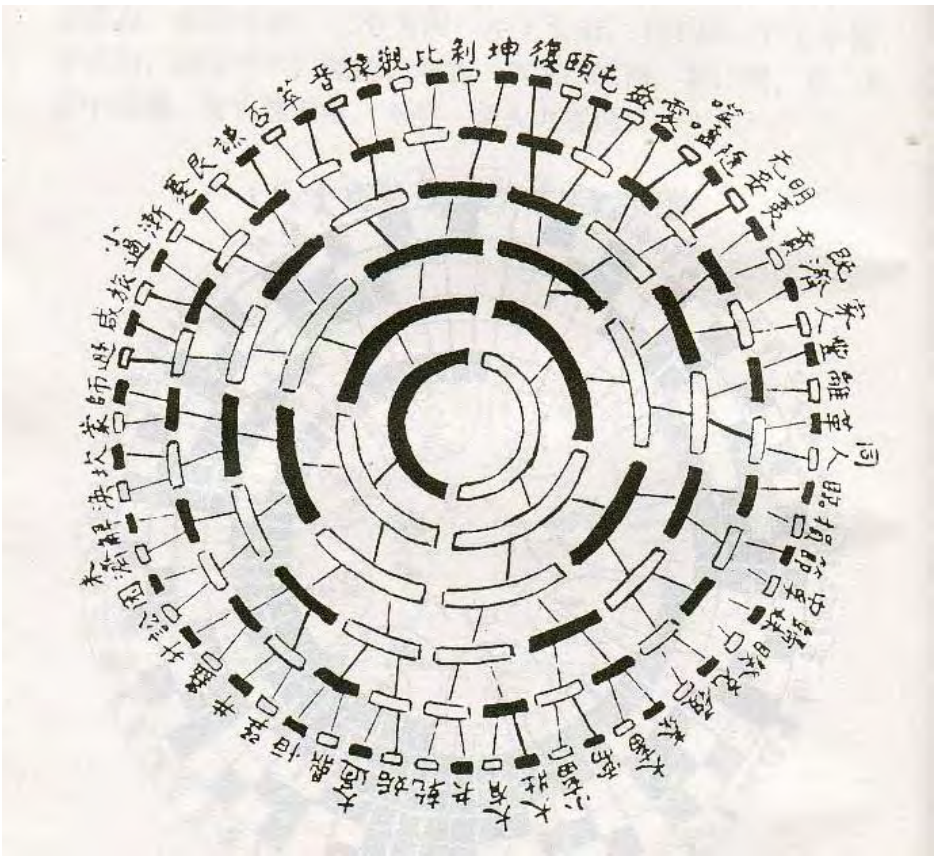

[171] G.-G. Leibniz. Explication de l'arithmétique binaire, qui se sert des seuls caractères 0 et 1 avec des remarques sur son utilité et sur ce qu'elle donne le sens des anciennes figures chinoises de Fohy.

"Ces figures sont peut-être le plus ancien monument de science qui soit au monde".

[172] D. Lach. Leibniz and China. *Journal of the History of Ideas*, 6(4):436–455, 1945.

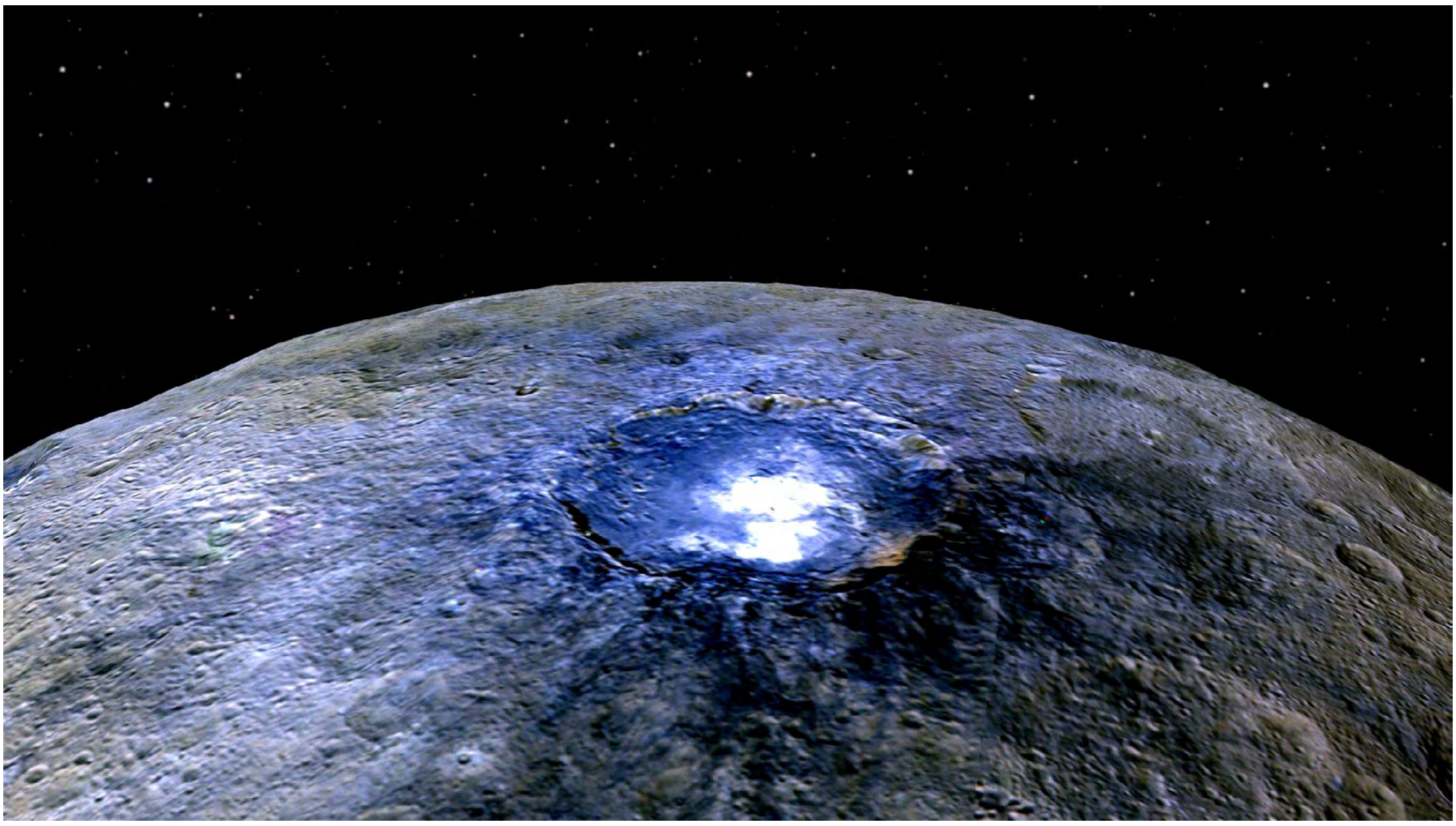

The dwarf planet Ceres, as seen by Dawn mission, in July 2016. The determination of its orbit was a spectacular achievement of Gauss. ©

# *Gauss, again: linking, magnetism and astronomy*

## *Gauss and linking numbers*

On January 22nd, 1833, Gauss wrote some enigmatic formula in his notebook[173].

ZUR ELECTRODYNAMIK. 605

[4.]

Von der *Geometria Situs*, die Leibnitz ahnte und in die nur einem Paar Geometern (Euler und Vandermonde) einen schwachen Blick zu thun vergönnt war, wissen und haben wir nach anderthalbhundert Jahren noch nicht viel mehr wie nichts.

Eine Hauptaufgabe aus dem *Grenzgebiet* der *Geometria Situs* und der *Geometria Magnitudinis* wird die sein, die Umschlingungen zweier geschlossener oder unendlicher Linien zu zählen.

Es seien die Coordinaten eines unbestimmten Punkts der ersten Linie $x, y, z$; der zweiten $x', y', z'$ und

$$\iint \frac{(x'-x)(dy\,dz'-dz\,dy')+(y'-y)(dz\,dx'-dx\,dz')+(z'-z)(dx\,dy'-dy\,dx')}{[(x'-x)^2+(y'-y)^2+(z'-z)^2]^{\frac{3}{2}}} = V$$

dann ist dies Integral durch beide Linien ausgedehnt

$$= 4m\pi$$

und $m$ die Anzahl der Umschlingungen.

Der Werth ist gegenseitig, d. i. er bleibt derselbe, wenn beide Linien gegen einander umgetauscht werden. 1833. Jan. 22.

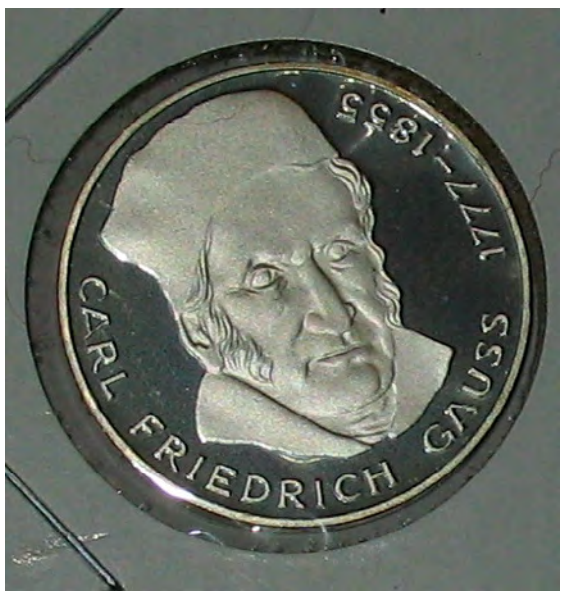

A coin from Germany, released in 1977, celebrating Gauss's 200th birthday. ©

[173] C. F. Gauß. *Werke. Ergänzungsreihe. Band V*. Georg Olms Verlag, Hildesheim, 1975. Briefwechsel: C. F. Gauss–H. C. Schumacher. Teil 3. [Correspondence: C. F. Gauss–H. C. Schumacher. Part 3], Edited by C. A. F. Peters, Reprint of the 1863 and 1865 originals.

*"Die Umschlingungen zweier geschlossener Lieben zu zählen."*

His purpose is to *"to count the linking number of two closed curves"*: an integer associated to two disjoint closed curves in 3-space which is invariant under deformation.

In 1833, Topology did not exist yet… Even the word would only appear in print twelve years later in a book by Listing. Leibniz had already coined the word *Analysis Situs* and was only dreaming about some science manipulating shapes just like algebra manipulates symbols. Gauss uses the terminology *Geometria Situs* and mentions Euler and Vandermonde as precursors.

Do not forget that these *Nachlasse* were not intended for publication. What would he have thought if he had been aware that his private drafts would become publicly available? This 1833 note was indeed published in 1867, after Gauss's death, and the editor included it in a volume dedicated to electromagnetism. That was a reasonable choice and a recent paper[174] does propose a good electromagnetic interpretation. Another paper[175] claims on the contrary that the formula has an astronomical origin and this paper seems just as credible. Who is right? Both, of course! Gauss was convinced by the deep unity of mathematics, and he would not build frontiers between mathematics, astronomy, physics etc. I will certainly not make a choice and I will present three parallel points of view: three definitions of the linking number of two disjoint closed curves in 3-space.

[174] R. L. Ricca and B. Nipoti. Gauss' linking number revisited. *J. Knot Theory Ramifications*, 20(10):1325–1343, 2011.

[175] M. Epple. Orbits of asteroids, a braid, and the first link invariant. *Math. Intelligencer*, 20(1):45–52, 1998.

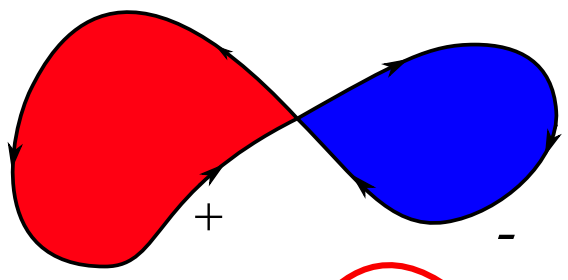

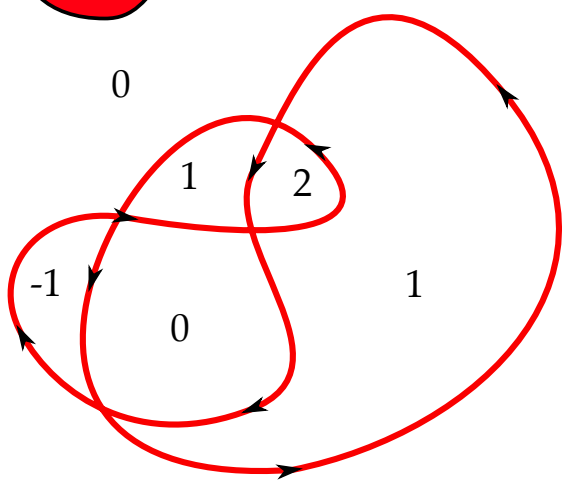

## *Geometry*

A closed oriented smooth *embedded* curve in the plane bounds a domain, which has some area. Let us say that this area is positive if the curve is oriented anti-clockwise and negative in the other case. That's easy. Now, if the curve is not embedded the situation is slightly more complicated, as for example in the figure eight curve in the margin. The left loop is oriented anti-clockwise and the right one clockwise, so that we are led to define the *signed area* as the algebraic sum of the two areas.

In the general case of an immersed curve with finitely many double points, we proceed in a similar way. The curve decomposes the plane into connected components. Let us equip the unbounded component with the coefficient 0. Now, equip each component with some integer with the convention that when we cross the curve positively this integer jumps by +1. In other words, a point moving on the curve in the positive direction sees a coefficient on its left equal to the coefficient on the right +1. It turns out that such a labeling exists and is unique. We then define the signed area of the curve as the linear combination

Try and prove the existence of such a labeling using the fact that the (algebraic) intersection number of two closed transverse oriented curves in the plane is 0. Let $c_1$ and $c_2$ be two transversal oriented curves in an oriented surface. Any intersection point of $c_1$ and $c_2$ has a sign $\pm 1$ depending on the orientation given by the pair of tangent vectors of $c_1$ and $c_2$ at this point. The sum of all these signs for all intersection points is the *algebraic intersection number* of $c_1$ and $c_2$. If $c_1, c_2$ are curves in the plane, this intersection number is 0. In modern terms, this follows from the fact that the homology of the plane is trivial. Gauss knew this fact. Can you produce a proof that he could have accepted?

of the geometric areas of the components with these integral coefficients. This definition is natural and is due to... Gauss.

Another definition comes from the fact that after all the *area below a curve* $y(x)$ is the integral of $y\,dx$. Consider the differential 1-form $\omega = -y\,dx$ in the plane and integrate it along the curve. I encourage the reader to check that these two definitions give the same number. En passant, note that the differential of $\omega$ is the 2-form $dx \wedge dy$, which is the area form. This is not a surprise for a 21st century mathematician but was far from obvious at the beginning of the 19th century.

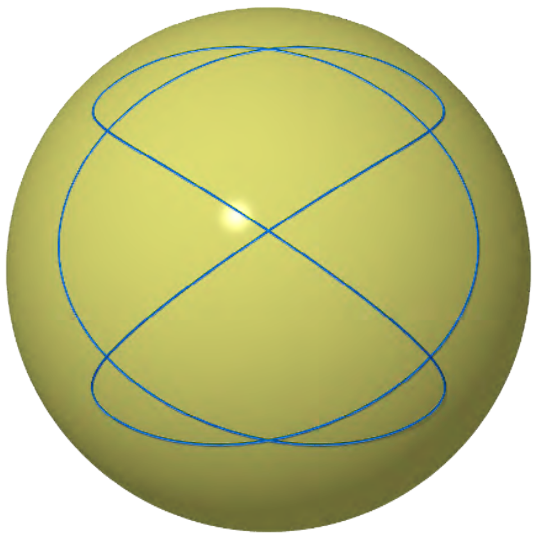

Consider now a closed oriented curve on the unit 2-sphere. Can we define the enclosed area? If the curve is embedded, there is no problem. The curve decomposes the sphere in two domains, one of them having the curve oriented as an anti-clockwise boundary, and we can define the area as the area of this domain. Now, if the curve is complicated, what can we do? We can still attribute numbers to the connected component of the complement with the same property as before, but they cannot be normalized by asking that some component at infinity has the weight 0, since there is no infinity. Therefore, all these integers are well defined up to the addition of the same integer to each component. The signed area enclosed by the curve is only defined up to the addition of an integral multiple of the area of the sphere, i.e. modulo $4\pi\mathbb{Z}$.

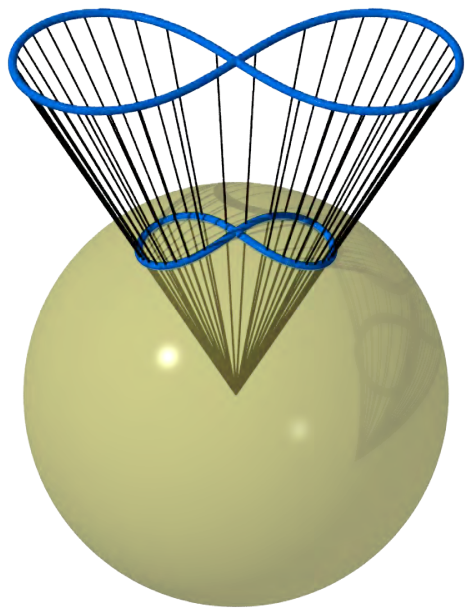

A closed curve in dimension 3 is the image of a map from $S^1$ to $\mathbb{R}^3$. In dimension 2, we can look at a map from $S^0$ to $\mathbb{R}^2$, i.e. two points $P, Q$ in the plane! The cone is now replaced by a triangle, and the angle function $A : \mathbb{R}^2 \setminus \{P, Q\} \to \mathbb{R}/2\pi\mathbb{Z}$ maps $x$ to the angle $\widehat{PxQ}$. Looking at the fibers of $A$ might remind you of secondary school.

A closed smooth curve in the sphere defines a cone in 3-space, with apex at the origin. The area of the curve is by definition the *solid angle* of the cone, hence defined modulo $4\pi\mathbb{Z}$. Just like an oriented angle in the plane is defined modulo $2\pi\mathbb{Z}$.

Suppose now that we have a closed oriented curve $\gamma$ in 3-space, not necessarily embedded. For every point $x$ *outside* $\gamma$, look at the solid angle $A_\gamma(x)$ of the cone with apex $x$ and based on $\gamma$. Its solid angle defines a function

$$A_\gamma : \mathbb{R}^3 \setminus \gamma \to \mathbb{R}/4\pi\mathbb{Z}.$$

Note that if $\gamma$ is a knot, i.e. if it is an embedded circle, the pre-image $A_\gamma^{-1}(\theta)$ is an orientable surface having $\gamma$ as its boundary for every regular value $\theta$ of $A_\gamma$. We have already met such surfaces under the name of Seifert.

Let us compute the differential $dA_\gamma$. Let $x, x'$ be two nearby points in space, both being away from $\gamma$. To compute $A_\gamma(x)$, we should translate $\gamma$ by $-x$, project radially the result onto the unit sphere and compute its signed area. The difference of areas $A_\gamma(x) - A_\gamma(x')$ is the signed area of the projection on the unit sphere of the annulus bounding the translations of $\gamma$ by $-x$ and $-x'$. Approximate $\gamma$ by some polygonal curve so that $A_\gamma(x) - A_\gamma(x')$ is approximated by the sum of signed areas of the projections of some parallelograms. Note that if $\delta x, \delta' x$ are two vectors in space, the volume of the pyramid with apex 0 and base

$$x, x + \delta x, x + \delta' x, x + \delta x + \delta' x$$

is

$$\frac{1}{3} det(x, \delta x, \delta' x).$$

If $\delta x$ and $\delta' x$ are very small, the corresponding solid angle is approximately obtained by dividing this value by the norm of $x$ cubed. Putting everything together, going to the limit, we get a formula for $dA_\gamma$ at the point $x$, on the vector $v$:

$$dA_\gamma(x, v) = \int_\gamma \frac{1}{||\gamma(t) - x||^3} det\left(\gamma(t), \frac{d\gamma}{dt}(t), v\right) dt.$$

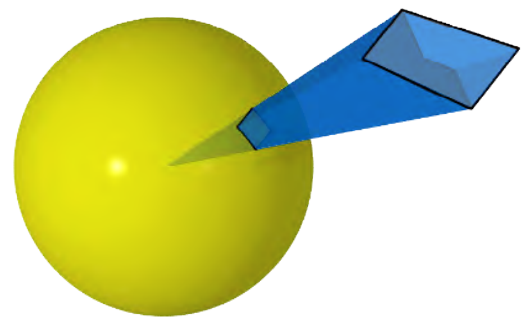

Solid angle. ©

We know the concept of *Cauchy index* of a closed curve $c$ in the punctured plane $\mathbb{C} \setminus \{0\}$: the number of "turns" around the origin when one goes along $c$. Follow the argument by continuity as you go around $c$ and then count the increase of the argument when you are back to the starting point. One could also use the differential form $\frac{1}{2i\pi} dz/z$ and integrate it on $c$. All this seems easy to today's students, but was not obvious for the founding fathers Gauss-Cauchy etc.

Now, do the exact same thing in 3-space, replacing the argument by the solid angle created by some closed curve $\gamma$. If a curve $\gamma'$ does not intersect $\gamma$, go around $\gamma'$ and look at the increase of the solid angle when you make the full turn (divided by $4\pi$). This index is called the *linking number* of $\gamma$ and $\gamma'$: this is an integer.

Using the formula for $dA_\gamma$, we get Gauss's formula for the linking number $link(\gamma,\gamma')$:

$$\frac{1}{4\pi}\iint \frac{1}{\|\gamma(t)-\gamma'(t')\|^3} det\left(\gamma(t)-\gamma'(t'), \frac{d\gamma}{dt}(t), \frac{d\gamma'}{dt'}(t')\right) dt\, dt'.$$

This is exactly what Gauss wrote in his notebook on January 22nd, 1833.

Note that the above formula shows that the linking is symmetric $link(\gamma,\gamma') = link(\gamma',\gamma)$, which was not obvious from the definition. This is what Gauss wrote:

> The value is symmetric: it remains the same when one interchanges the two curves.

*Der Werth ist gegenseitig, d. i. er bleibt derselbe, wen beide Linien gegen einander umgetauscht werden.*

Note also that if $\gamma$ and $\gamma'$ are deformed continuously in such a way that they don't intersect during the deformation, the linking number has to be constant: an integer cannot change continuously. This is the most important feature of the linking number: it is *invariant under deformation*.

## *Astronomy*

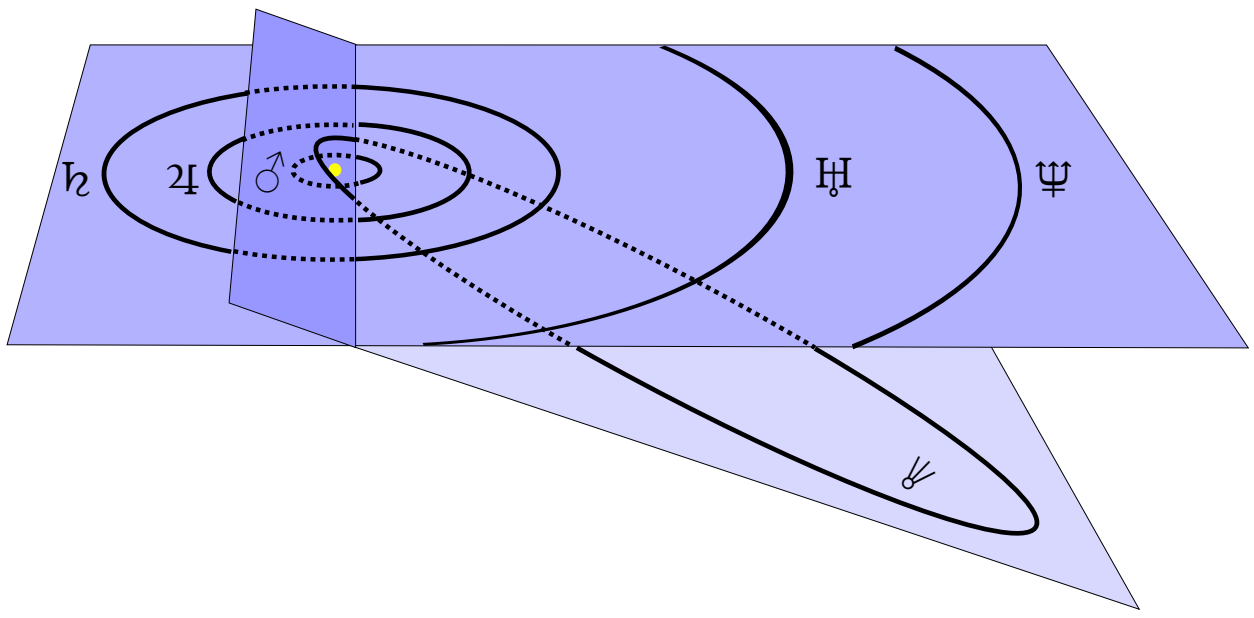

The orbit of Halley's comet, together with the orbits of Mars, Jupiter, Saturn, Uranus and Neptune.

The paper of Epple mentioned above proposes a possible approach to linking numbers. One of the first accomplishments of Gauss which made him famous was his determination in 1801 of the orbit of the dwarf planet Ceres, which had just been discovered. Suppose we observe a planet from a fixed position on our planet Earth. Where should we look in the sky? More precisely, let $\gamma$ be the trajectory of the Earth in fixed space (fixed with respect to the Sun) and let $\gamma'$ be the trajectory of the planet that we want to observe. For simplicity, I assume that $\gamma$ and $\gamma'$ are disjoint ☺. If the periods of rotation are rationally independent, the positions of the Earth and the planet *on their orbits* are independent random variables. Gauss calls *zodiacus* of the planet (relative to the Earth) the image of the map:

$$\varpi : (t, t') \in (\mathbb{R}/\mathbb{Z})^2 \mapsto \frac{\gamma(t) - \gamma'(t')}{||\gamma(t) - \gamma'(t')||} \in \mathbb{S}^2.$$

This is the zone in celestial sphere where the observer should look for the planet.

The integrand in Gauss's formula for the linking number is simply the Jacobian determinant of this map, so that the linking number is $1/4\pi$ times the signed area of the zodiacus. A modern mathematician knows that the integral of the Jacobian determinant of a map between two oriented manifolds of the same dimension is the topological degree of this map. Therefore the *linking number can also be defined as the degree of the zodiacus map* $\varpi$.

Of course, Gauss studied in detail the case of two ellipses in space.

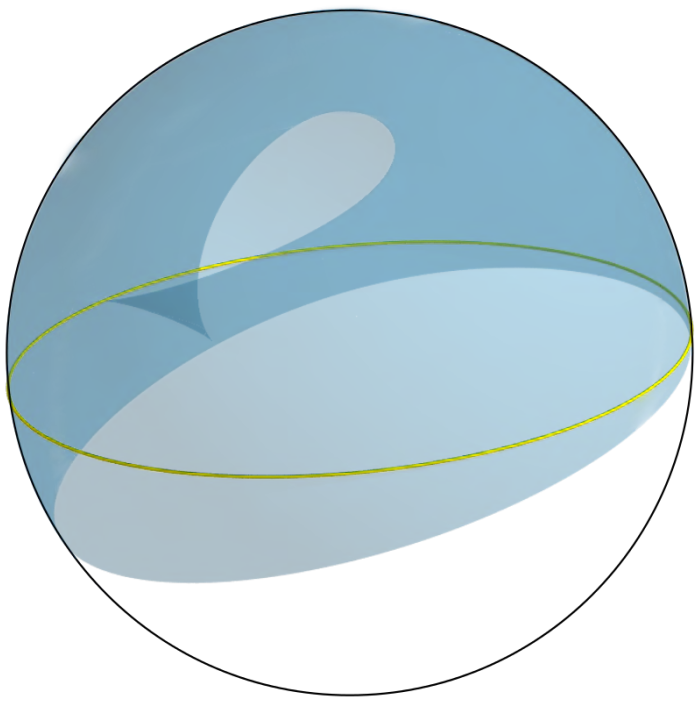

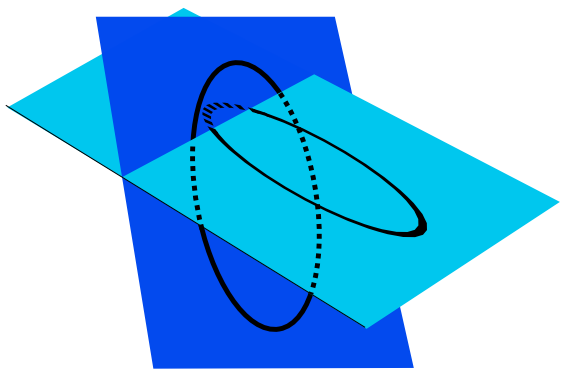

Two unlinked ellipses. The zodiacus is on the left.

*When the two ellipses are not linked*, like in the picture above, the zodiacus does not cover the full celestial sphere. The light blue zone in the zodiacus corresponds to points which are covered twice by $\varpi$. The darker zone, with two singular points, is covered four times. Compare with the usual picture (in the margin) of the perspective of a torus of revolution and the singularities which appear in its contour.

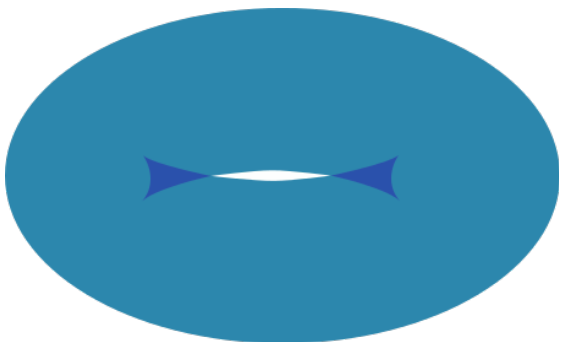

Projection of a torus on a plane.

*When the ellipses are linked*, like in the picture below, the zodiacus is the full sphere. The light blue zone in the zodiacus corresponds to points which are covered only once by $\varpi$. The darker zone, with four singular points, is covered three times.

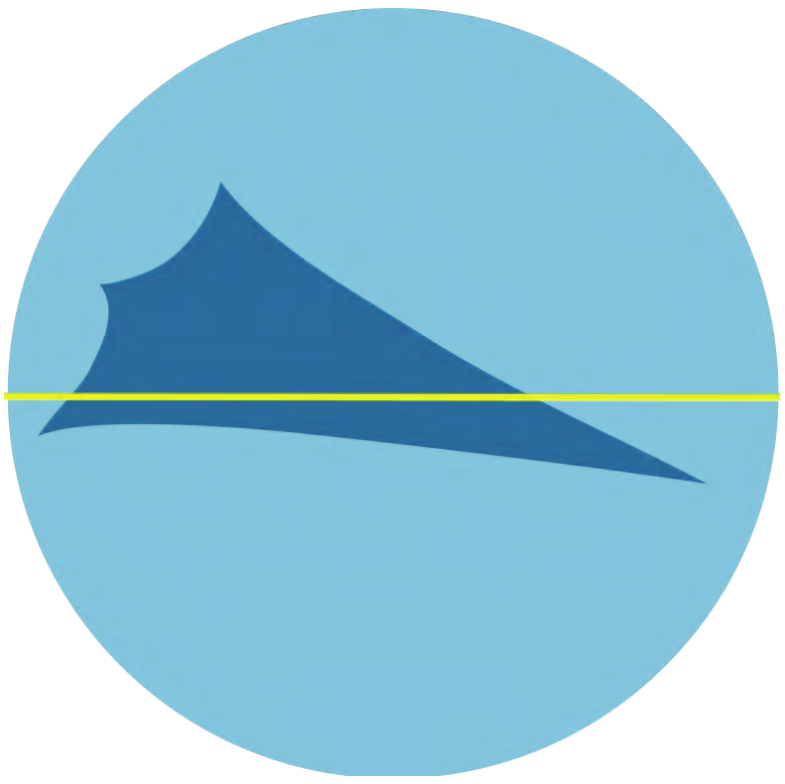

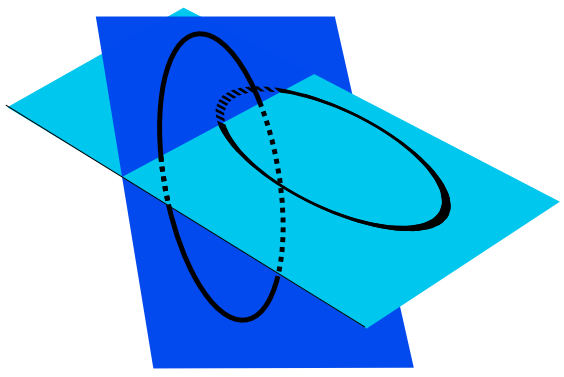

Two linked ellipses. The zodiacus is on the left.

Given a smooth map $f : M \to N$ between two compact oriented connected manifolds without boundary, there are several possible definitions for its *topological degree*. The first consists in choosing some volume form $vol$ on $N$, of total volume 1, and integrating its pull-back $f^\star vol$ on $M$. It is not hard to see that this is independent of $vol$. Indeed, if $vol'$ is another choice of volume form, $vol' - vol$ is an exact form, and therefore the integral of $f^\star vol' - f^\star vol$ is zero. From this definition, it is easy to see that this is invariant under deformation. Indeed, if two maps $f_0, f_1$ are homotopic, $f_0^\star vol - f_1^\star vol$ is an exact form. It is less easy to see that this degree is an integer.

A second definition consists in picking a regular value $y \in N$ of $f$ and looking at the finitely many pre-images $x_1, \dots, x_n$ in $M$. At each of these pre-images, the differential of $f$ preserves or reverses orientation, and we attribute them a + or - sign accordingly. The degree of $f$ is the sum of these signs. One has to show that this does not depend on the choice of the

Look again at the preceding pictures of the zodiacus of two ellipses and figure out the + or − signs.

regular value and that it is a homotopy invariant. This is proved brilliantly in Milnor's book[176]. One should also prove that the two definitions agree…One possibility is to use a sequence of volume forms on $M$ which converges to the Dirac mass at the regular value $y$.

[176] J. W. Milnor. *Topology from the differentiable viewpoint*. Princeton Landmarks in Mathematics. Princeton University Press, Princeton, NJ, 1997.

Let us use the regular value point of view to compute the degree of the zodiacus map $\varpi : \mathbb{R}^2/\mathbb{Z}^2 \to \mathbb{S}^2$. Choose the south pole of $\mathbb{S}^2$ as the point $y$. The pre-images of $y$ consist of the pairs of points $\gamma(t)$ and $\gamma'(t')$ such that $\gamma(t)$ is above $\gamma'(t')$: they have the same $x, y$ coordinates and the $z$-coordinate of $\gamma(t)$ is bigger than that of $\gamma'(t')$. The differential of $\varpi$ at such a point is easy to compute. It is non-degenerate if the projections of $\gamma$ and $\gamma'$ on the $(x, y)$ plane intersect transversally at the corresponding point. Its Jacobian determinant is *positive* (resp. *negative*) if the intersection of the projections is positive (resp. negative).

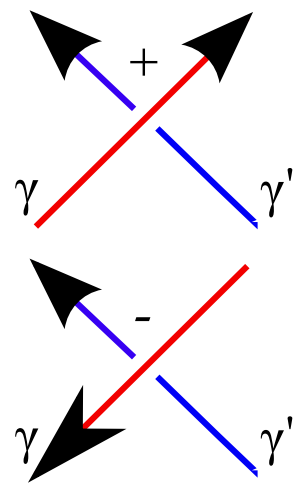

This leads to the combinatorial definition of the linking number, that Gauss obviously knew. Project the two curves $\gamma, \gamma'$ in a generic plane so that the two projections $\bar{\gamma}, \bar{\gamma}'$ intersect transversally. Mark an index +1 or −1 to each intersection point, according to whether the tangent vectors at $\gamma$ and $\gamma'$ define a positive or negative basis. Among those intersection points, select only those where $\gamma$ is *over* $\gamma'$. The sum of the corresponding signs is the linking number of $\gamma, \gamma'$.

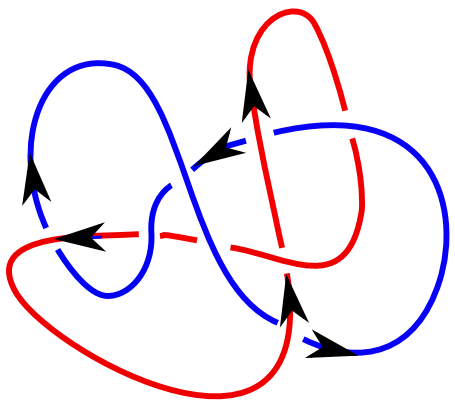

On the picture, the red curve passes 3 times over the blue one with signs $+1, +1, -1$. The linking number is 1.

The so-called *Whitehead link* in the margin has linking number 0 but this does not mean that the two components can be separated by some deformation[177]. Show that there is no complex algebraic curve with two branches such that the associated link is this link.

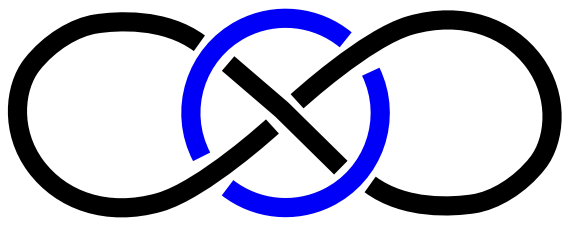

Whitehead link.

[177] D. Rolfsen. *Knots and links*, volume 7 of *Mathematics Lecture Series*. Publish or Perish, Inc., Houston, TX, 1990. Corrected reprint of the 1976 original.

## *Electromagnetism*

Gauss's formula is reminiscent of the Biot-Savart law in physics. An electric current generates a magnetic field. Suppose a closed wire $\gamma$ carries some steady current, with intensity $i$ and let $x$ be a point outside the wire. Then the magnetic field created at $x$ is

$$B(x) = \frac{\mu_0 i}{4\pi} \int_\gamma \frac{1}{||(\gamma(t) - x||^3} \left( (\gamma(t) - x) \wedge \frac{d\gamma(t)}{dt} \right) dt$$

where $\mu_0$ is the magnetic constant. This vector field is the dual of the closed 1-form $dA_\gamma$ (with respect to the Euclidean metric on the 3-dimensional physical space). It may also be interpreted as the gradient field of a local primitive of the 1-form $A_\gamma$. It follows that the circulation of the magnetic field on some loop $\gamma'$ is the same as the integral of $dA_\gamma$ on $\gamma'$, i.e. the linking number. Hence, the linking $link(\gamma, \gamma')$ is the circulation of the magnetic field created by a current.

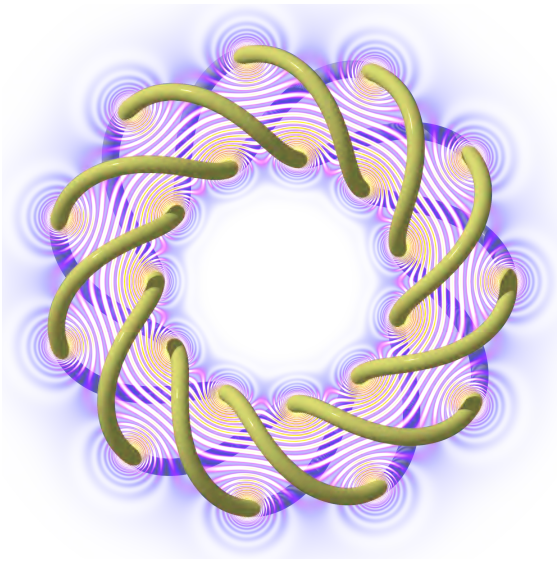

The magnetic field generated by a torus knot. ©

The above mentioned paper by Ricca and Nipoti gives an interesting reconstruction of what might have been the magnetic interpretation in Gauss's mind. Do not forget that, together with Weber, Gauss established the first telegraph transmitting messages across Goettingen.

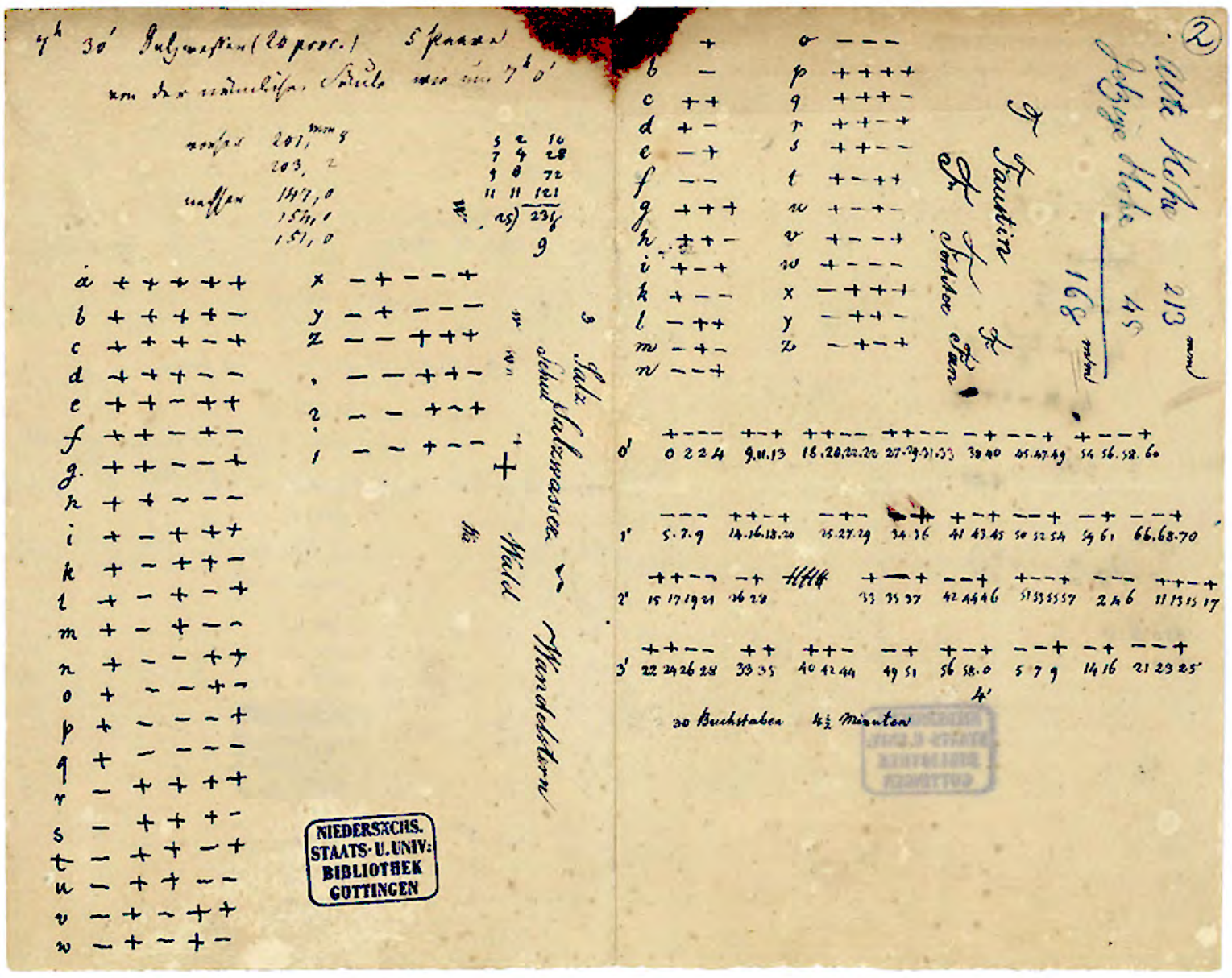

The Gauss-Weber code for their telegraph: combinatorics again! ©

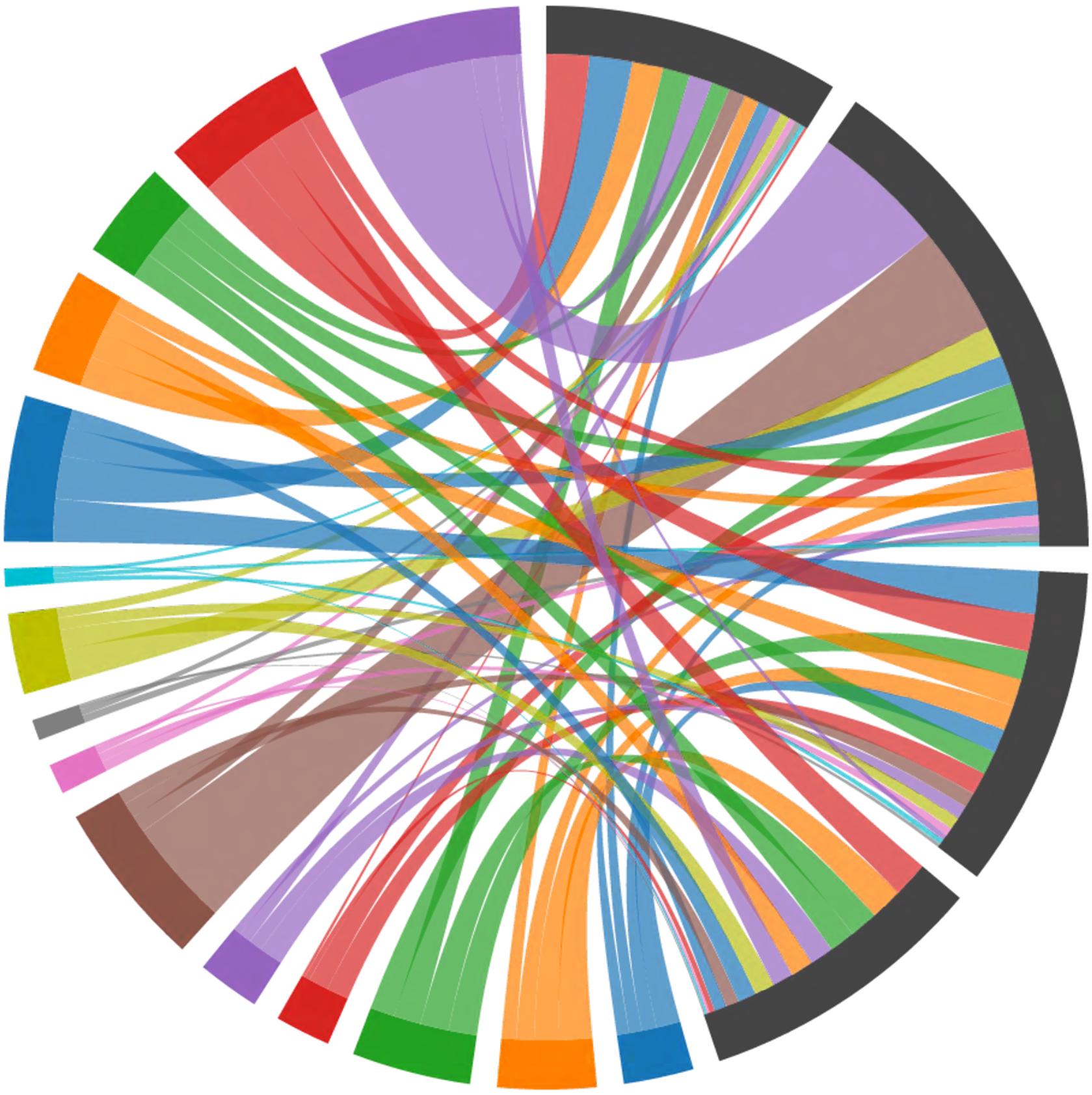

A chord diagram is also a graphical method of displaying the inter-relationships between data. Given a stochastic matrix $a_{ij}$ (i.e. $a_{ij} > 0$ and $\sum_j a_{ij} = 1$), one can think of $a_{ij}$ as the proportion of entity $i$ interacting with $j$. One draws $n$ intervals $I_1, ..., I_n$ around the circle, whose lengths $l_1, ..., l_n$ have to be determined, and bands connecting $I_i$ and $I_j$ with widths $a_{ij} l_i$. The compatibility condition can be expressed as $l_j = \sum_i a_{ij} l_i$. The existence of a solution is guaranteed by the Perron-Frobenius theorem. ©

# *Kontsevich is back: A universal invariant*

Is this promenade a loop homotopic to a point? We are back to our starting point: Maxim Kontsevich. This chapter is not a conclusion but an opening to a vast domain and shows that in mathematics it is possible to come back to very old ideas with a completely new perspective. I want to present a short introduction to a wonderful development in knot theory involving chord diagrams, in a 1993 paper of Kontsevich[178].

[178] M. Kontsevich. Vassiliev's knot invariants. In *I. M. Gelfand Seminar*, volume 16 of *Adv. Soviet Math.*, pages 137–150. Amer. Math. Soc., Providence, RI, 1993.

## *A new point of view on the linking number*

Let $\gamma_1, \gamma_2 : \mathbb{R}/\mathbb{Z} \to \mathbb{R}^3$ be two disjoint oriented closed curves. We know that the linking number of $\gamma_1$ and $\gamma_2$ is the topological degree of Gauss's *zodiacus map*, from the product of the two curves to the unit sphere. Kontsevich's formula will express the same number as a topological degree of a map *from an oriented* 1-*dimensional manifold (which is therefore a union of circles) to a circle*. The great advantage is that this new point of view enables us to define many more invariants.

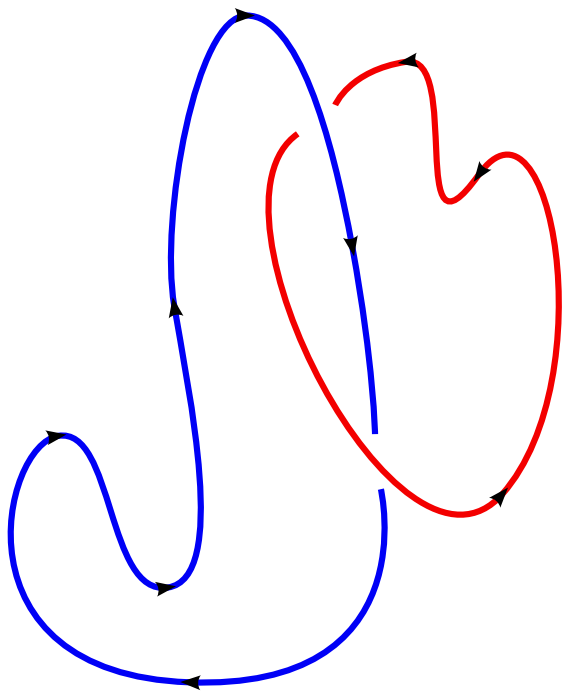

Let us think of the space $\mathbb{R}^3$ with coordinates $(x, y, t)$ as the product of the complex line $\mathbb{C}$ (with coordinate $z = x + iy$) and $\mathbb{R}$ (with coordinate $t$). Assume that our curves $\gamma_1, \gamma_2$ are *Morse*. This simply means that the projection onto the $t$-coordinate has a finite number of critical points and that the second derivative is not zero at these critical points. Assume that the critical values of the $t$-coordinates are all distinct.

Consider now the set of pairs of points on $\gamma_1, \gamma_2$ which have the same $t$-coordinate. Formally, this is the set

$$X = \{(s_1, s_2) \in (\mathbb{R}/\mathbb{Z})^2 \,|\, t(\gamma_1(s_1)) = t(\gamma_2(s_2))\}.$$

This is a smooth curve of the 2-torus. The only (easy) thing to check is that this is indeed the case in the neighborhood of critical points.

Look at this example. There are 8 critical values, decomposing the first curve in 18 strands and the second in 10 strands.

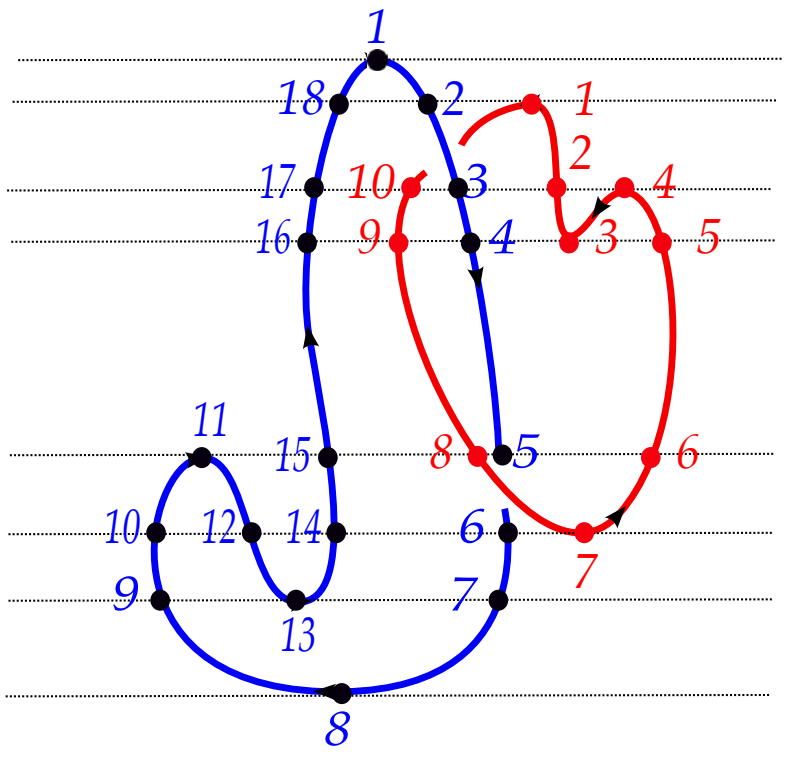

The submanifold $X$ is represented in the following picture.

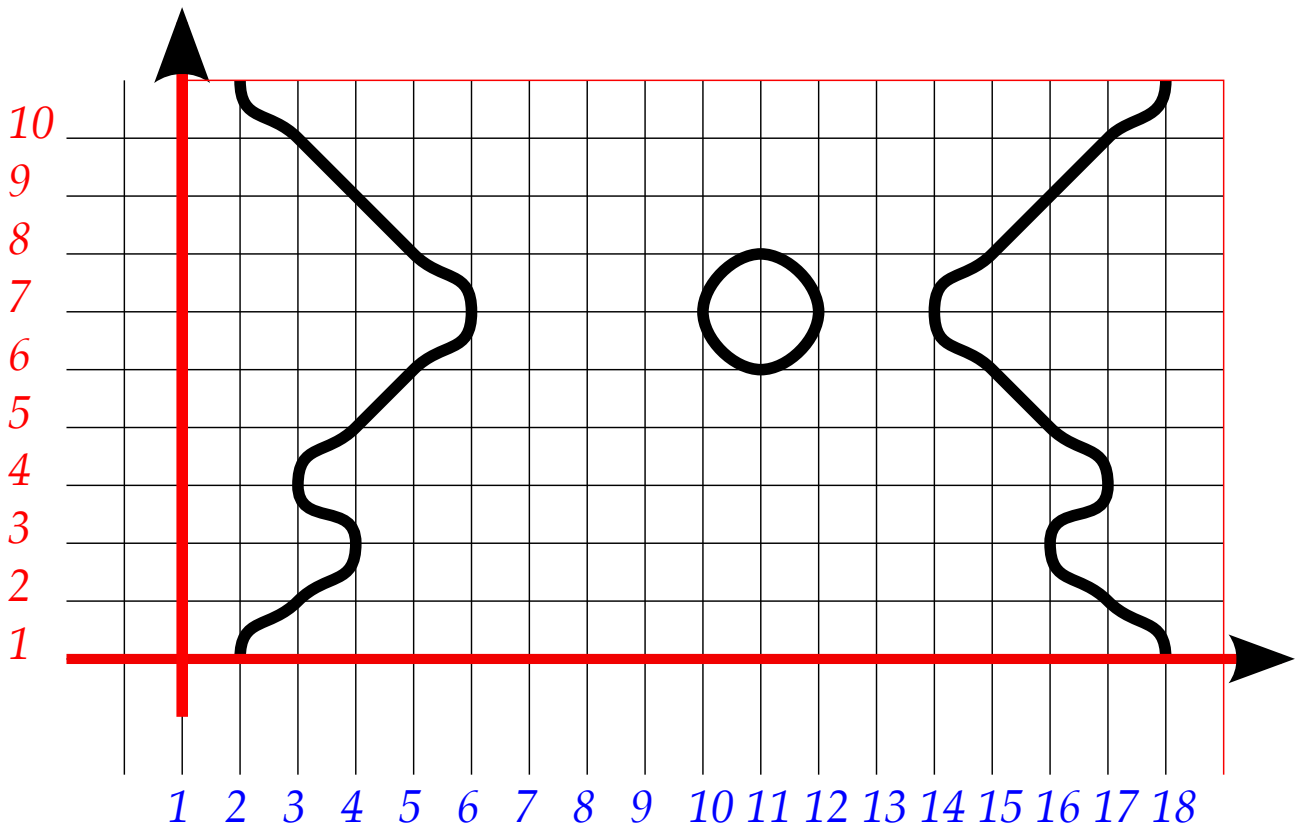

There is a canonical orientation on $X$. Choose a small interval $I$ in $X$, away from the critical values. This interval maps diffeomorphically onto some interval $I_1$ in $\gamma_1$ and onto some other interval $I_2$ in $\gamma_2$. A non-critical interval in $\gamma_1$ (or in $\gamma_2$) is equipped with two orientations, coming from the orientation of $\gamma_1$ (or $\gamma_2$) on the one hand, and from the $t$-coordinate, on the other hand. I will say that such an interval is *positive* if these two orientations agree and *negative* otherwise. Orient $I$ using increasing $t$ if $I_1$ and $I_2$ are both positive or both negative, and using decreasing $t$ otherwise.

Check that this does define an orientation on $X$.

Now, a point in $X$ defines two points $\gamma_1(s_1)$ and $\gamma_2(s_2)$ which project to *distinct* points in the complex plane $x + iy$. The argument of the difference defines a map $\pi : X \to \mathbb{S}^1$. This is a map between oriented 1-*dimensional* manifolds.

*I claim that the degree of* $\pi$ *is the linking number of* $\gamma_1$ *and* $\gamma_2$.

Let us prove this claim. The linking number is the topological degree of the map

$$\varpi : (s_1, s_2) \in (\mathbb{R}/\mathbb{Z})^2 \mapsto \frac{\gamma_1(s_1) - \gamma_2(s_2)}{\|\gamma_1(s_1) - \gamma_2(s_2)\|} \in \mathbb{S}^2$$

between oriented surfaces. The unit sphere $\mathbb{S}^2$ contains the *horizontal equator* $\mathbb{S}^1$ (where $t = 0$). The assumption that $\gamma_1$ and $\gamma_2$ are Morse with distinct critical values implies that $\varpi$ is transversal to $\mathbb{S}^1 \subset \mathbb{S}^2$. The inverse image $\pi^{-1}(\mathbb{S}^1)$ is $X$, by definition. The differential of $\varpi$ identifies the normal bundle of $X$, in the 2-torus, with the normal bundle of the equator, in the sphere. Our orientation convention on $X$ is such that this identification is positive.

We want to compare the two topological degrees of $\varpi$ and $\pi$. Take a regular value $v \in \mathbb{S}^1 \subset \mathbb{S}^2$ of $\pi$ and let $u$ be a point in its pre-image. Note that $v$ is also a regular value of $\varpi$. The sign of the Jacobian of the differential of $\pi$ at $u$ is the same as the sign of the Jacobian of $\varpi$ at $u$. It follows that the degrees of $\varpi$ and $\pi$ are equal. ⊡

We now get a new formula for the linking number, using Cauchy type indices. This is a special case of Kontsevich's theorem.

**Theorem.** *Slice $\gamma_1$ and $\gamma_2$ by horizontal planes passing through the critical points of the t-coordinate of $\gamma_1$ or $\gamma_2$. Between two consecutive planes, $\gamma_1, \gamma_2$ define a certain number of strands which are positive or negative. Choose one of the strands corresponding to $\gamma_1$, defined by some graph $(\zeta_1(t), t)$ (for $t_- \leq t \leq t_+$). Choose a strand $(\zeta_2(t), t)$ for the curve $\gamma_2$ (also for $t_- \leq t \leq t_+$). Compute the* amount of rotation

$$\epsilon \frac{1}{2i\pi} \int_{t_-}^{t_+} \frac{d\left(\zeta_1(t) - \zeta_2(t)\right)}{\zeta_1(t) - \zeta_2(t)}$$

*where $\epsilon$ is $+1$ if the two chosen strands have the same sign and $-1$ otherwise.*

*Sum all these numbers for all possible pairs of consecutive horizontal planes, and for all pairings of a strand for $\gamma_1$ and a strand for $\gamma_2$. The result is the linking number $lk(\gamma_1, \gamma_2)$.*

In the previous example, the 8 singular values define 7 intervals containing 2, 2, 2, 2, 4, 4, 2 blue strands, and 0, 2, 4, 2, 2, 0 red strands. Therefore there are

$$2 \times 0 + 2 \times 2 + 2 \times 4 + 2 \times 2 + 4 \times 2 + 4 \times 0 + 4 \times 0 = 24$$

pairings between strands, which corresponds to the number of intervals in $X$.

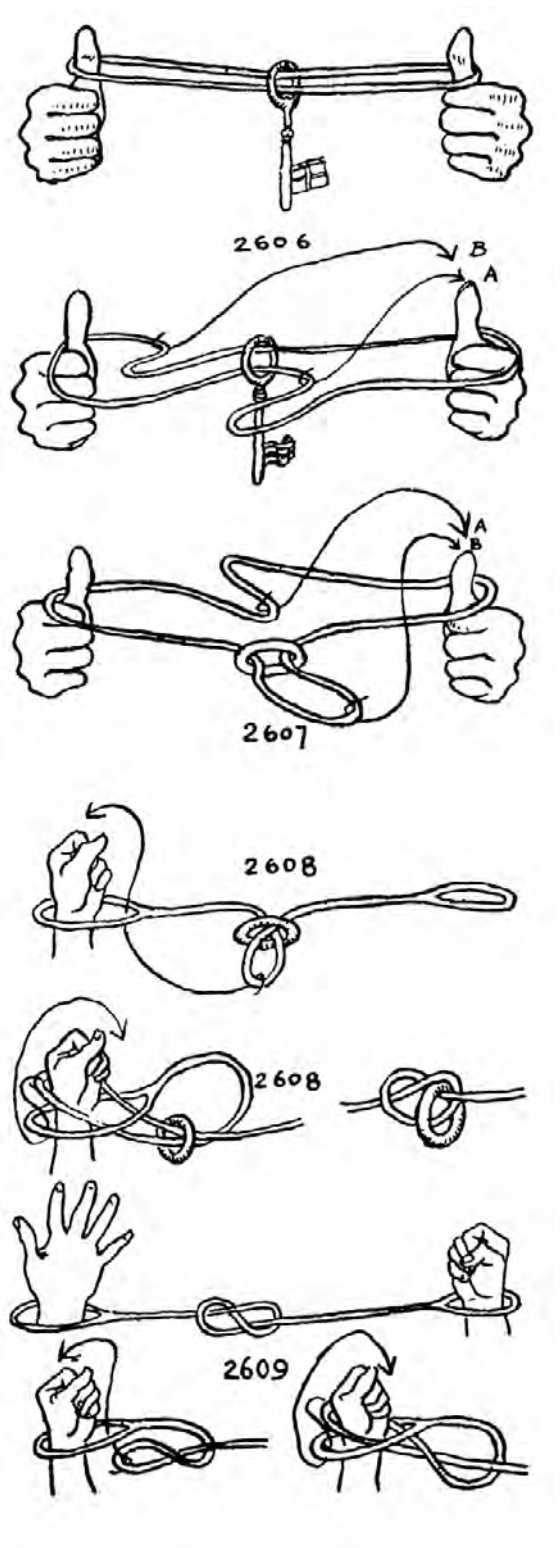

Some magic tricks trying to convince you that the linking number does not exist. From the marvelous "Ashley's book of knots" (1944), containing 7000 pictures representing 3800 knots. ©

## *The universal Kontsevich invariant of a knot with values in the chord algebra*

To conclude, I sketch the definition of an invariant associated to a knot with values in formal series with coefficients in chord diagrams. This is a brilliant idea of Kontsevich, from his famous 1993 paper.

Let $Chord(n)$ be the set of chord diagrams with $n$ chords.. As we have seen many times, they are sets of $2n$ points on an oriented circle, grouped in pairs, up to orientation preserving homeomorphisms of the circle. Denote by $\mathbb{C}[Chord]$ the vector space having the union $Chord$ of all $Chord(n)$'s as a basis. Its elements are therefore finite sums $\sum_{w \in Chord_n} \lambda_w . w$ where $\lambda_w = 0$ for all but a finite number of $w$. Consider $\mathbb{C}[Chord]$ as a *graded* vector space, the grading being given by $n$.

Denote by $\mathcal{A}$ the quotient of $\mathbb{C}[Chord]$ by the subspace generated by two relations, which might appear artificial at first sight:

— the *one term relation*. This means that any chord diagram obtained from the picture below by completing it in *any way* in the dotted part of the circle is declared to be 0 in $\mathcal{A}$. Said differently, every diagram containing a solitary chord is equal to 0 in $\mathcal{A}$.

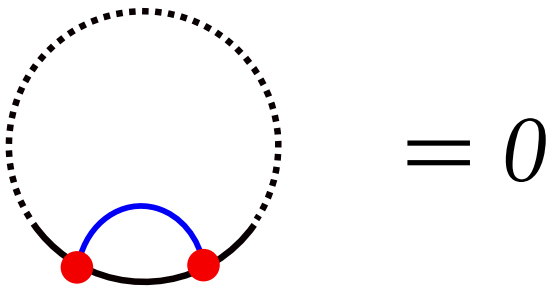

— the *four term relation*. Analogously, the dotted part of the circle can be completed in any way (but of course in the same way in the four constituents of the relation).

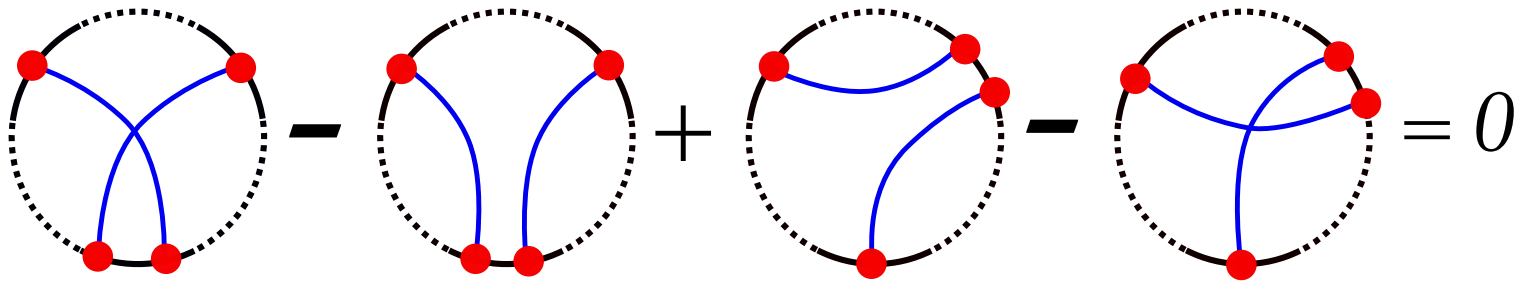

This vector space $\mathcal{A}$ is actually a graded *algebra* $\oplus_{n\geq 0}\mathcal{A}_n$. Two chord diagrams can be multiplied in the following way.

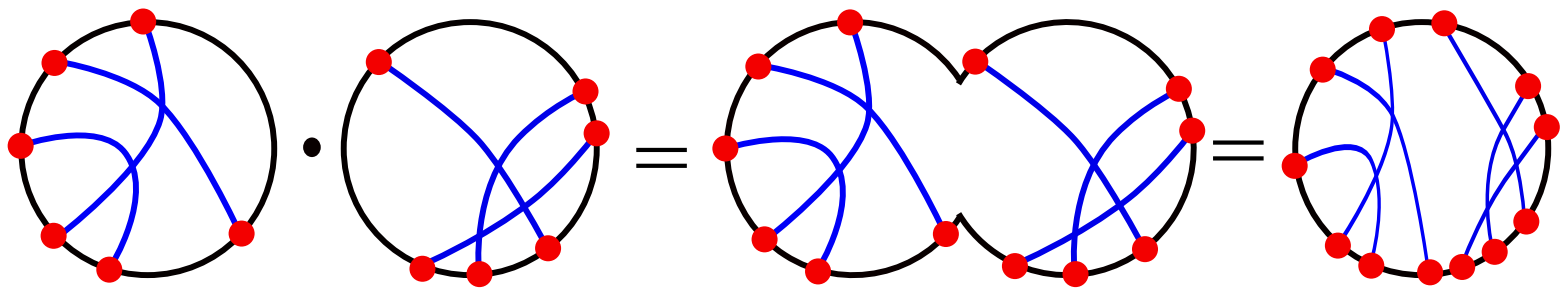

The four term relation is exactly what is needed to make sure that this operation is well defined, independently of the locus of the *connected sum*.

Consider the completion $\hat{\mathcal{A}}$, where we add infinite formal sums $\sum_{w\in Chord_n} \lambda_w.w$ with no condition on the numbers $\lambda_w$. Let us call $\hat{\mathcal{A}}$ the *chord algebra*.

I can now define the *Kontsevich universal invariant of a knot, with values in* $\hat{\mathcal{A}}$.

Let $\gamma$ be some knot in 3-space (assumed to be Morse).

Slice it by horizontal planes passing through the critical points of the $t$ coordinates. This decomposes the knot in a finite number of strands, which could be positive or negative, with respect to the orientation of the knot.

Choose some integer $n$. Consider the space of $2n$-tuples of distinct points $(p_1, q_1, \ldots, p_n, q_n)$ on the knot such that

$$t(p_1) = t(q_1) < t(p_2) = t(q_2) < \ldots < t(p_n) = t(q_n).$$

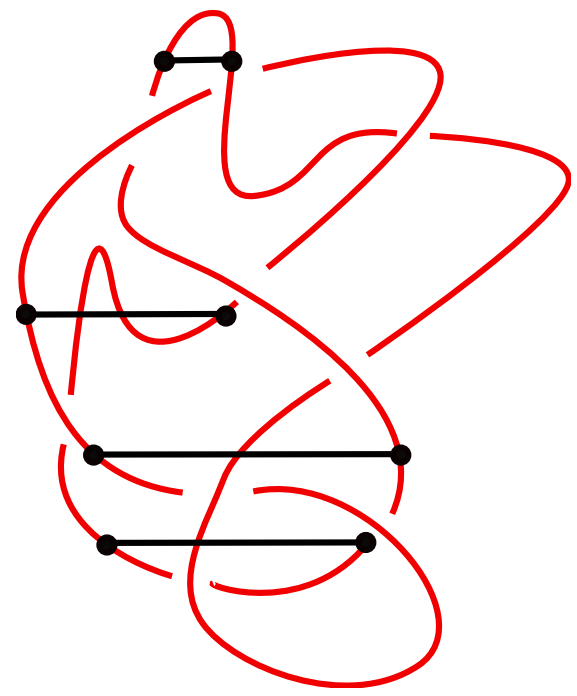

This is an $n$-dimensional submanifold $X_n$ with boundary of the $2n$-dimensional torus, canonically oriented by the orientation of the circle.

Note that any element of $X_n$ defines a chord diagram with $n$ chords.

There is a natural map $\varpi$ from $X_n$ to $(\mathbb{C}^\star)^n$. Indeed, if $p$ and $q$ are two distinct points on $\gamma$ with the same $t$-coordinate, their difference is a nonzero complex number. We therefore associate the $n$-tuple $(q_1 - p_1, \ldots, q_n - p_n) \in (\mathbb{C}^\star)^n$ to $(p_1, q_1, \ldots, p_n, q_n)$.

Consider now the (complex) differential $n$-form

$$\frac{1}{(2i\pi)^n}\varpi^\star\left(\frac{d\zeta_1}{\zeta_1} \wedge \ldots \wedge \frac{d\zeta_n}{\zeta_n}\right)$$

on $X_n$. Integrating it on each connected component of $X_n$, multiplying with the corresponding element of $Chord(n)$ and summing over all components of $X_n$, we get an element of $\mathcal{A}_n$. The formal sum of all these elements, for all values of $n$ defines an element of $\hat{\mathcal{A}}$: this is the *Kontsevich invariant of* $\gamma$, denoted by $Z(\gamma)$, which is an element of $\hat{\mathcal{A}}$.

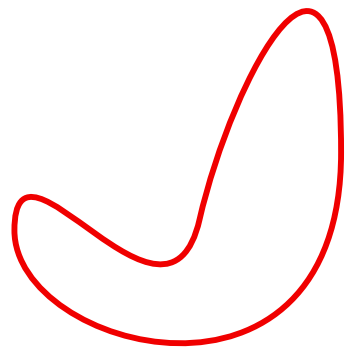

Strictly speaking, this is not yet an invariant! It turns out that this is only an invariant if the knot $\gamma$ is deformed among Morse knots, preserving the number of critical points. This is already a non-trivial fact.

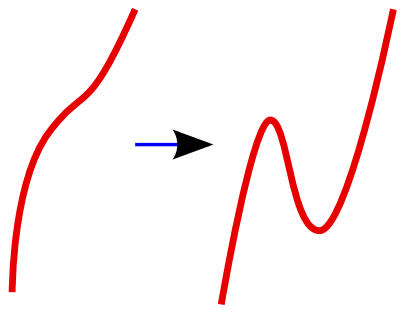

A general deformation of $\gamma$ could introduce a *hump*.

However, the change in this introduction of a hump can be completely described. Let $Z(H)$ be the invariant of the hump in the margin.

It can be shown that if the $t$-coordinate of a knot $\gamma$ has $2c$ critical points, the quotient

$$I(K) = Z(K)/Z(K)^{c/2} \in \hat{\mathcal{A}}$$

is an actual invariant of the knot $\gamma$, for any isotopy, that is any deformation of the knot, avoiding the creation of double points.

I still have to justify the division by $Z(K)^{c/2}$ in the algebra $\hat{\mathcal{A}}$. This is not difficult since it is easy to see that $Z(K)^{c/2}$ has the form $1 + a$ with $a$ of degree $> 1$, so that the inverse of $1 + a$ is $1 - a + a^2 - a^3 - \cdots$.

I proved absolutely nothing. I did not explain in which sense this invariant is universal. As a matter of fact, it is unknown whether two knots are equivalent if and only if they have the same invariant: that would be fantastic.

For a detailed presentation, I strongly recommend this article[179] and this book[180].

[179] S. Chmutov and S. Duzhin. The Kontsevich integral. *Acta Appl. Math.*, 66(2):155–190, 2001.

[180] S. Chmutov, S. Duzhin, and J. Mostovoy. *Introduction to Vassiliev knot invariants*. Cambridge University Press, Cambridge, 2012.

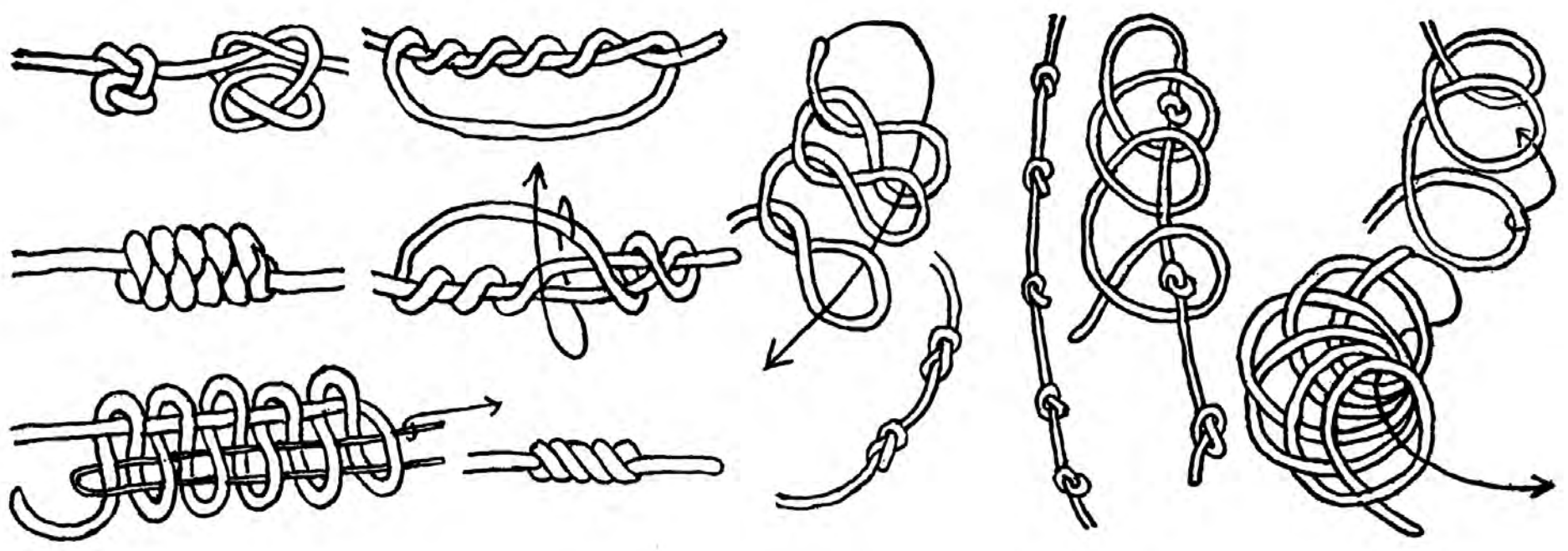

More from Ashley's book of knots.

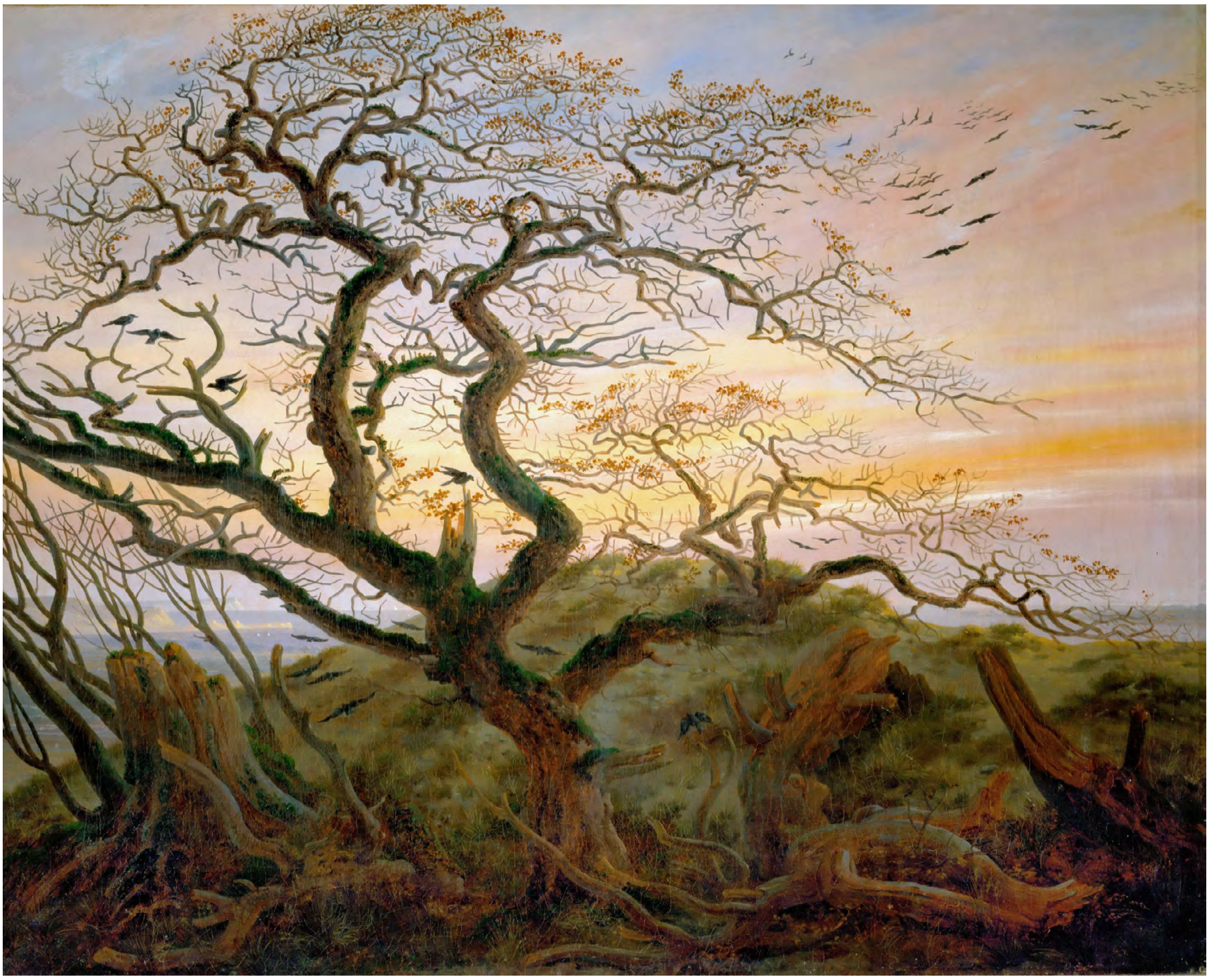

Caspar David Friedrich:
Tree of crows. ©

# *Postface*

Our promenade is over. We have wandered in quite a lot of mathematical forests. We have indeed seen many trees, and our travel was definitely not a geodesic path. My reader will hopefully want to travel more and to explore new territories in much more detail, maybe more seriously.

Since our stroll was some kind of closed loop which began with the romantic *Wanderer in the fog*, by Caspar David Friedrich, perhaps it is appropriate to now admire *The tree of crows* by the same artist, dated 1822. By this time, Gauss was dreaming about non-Euclidean geometry.

This painting has been chosen as a frontispiece for one of my favorite mathematical books[181] which also deals with trees, albeit very different from those seen during our promenade. The next destination for my reader?

[181] J.-P. Serre. *Trees*. Springer Monographs in Mathematics. Springer-Verlag, Berlin, 2003. Translated from the French original by John Stillwell, Corrected 2nd printing of the 1980 English translation.

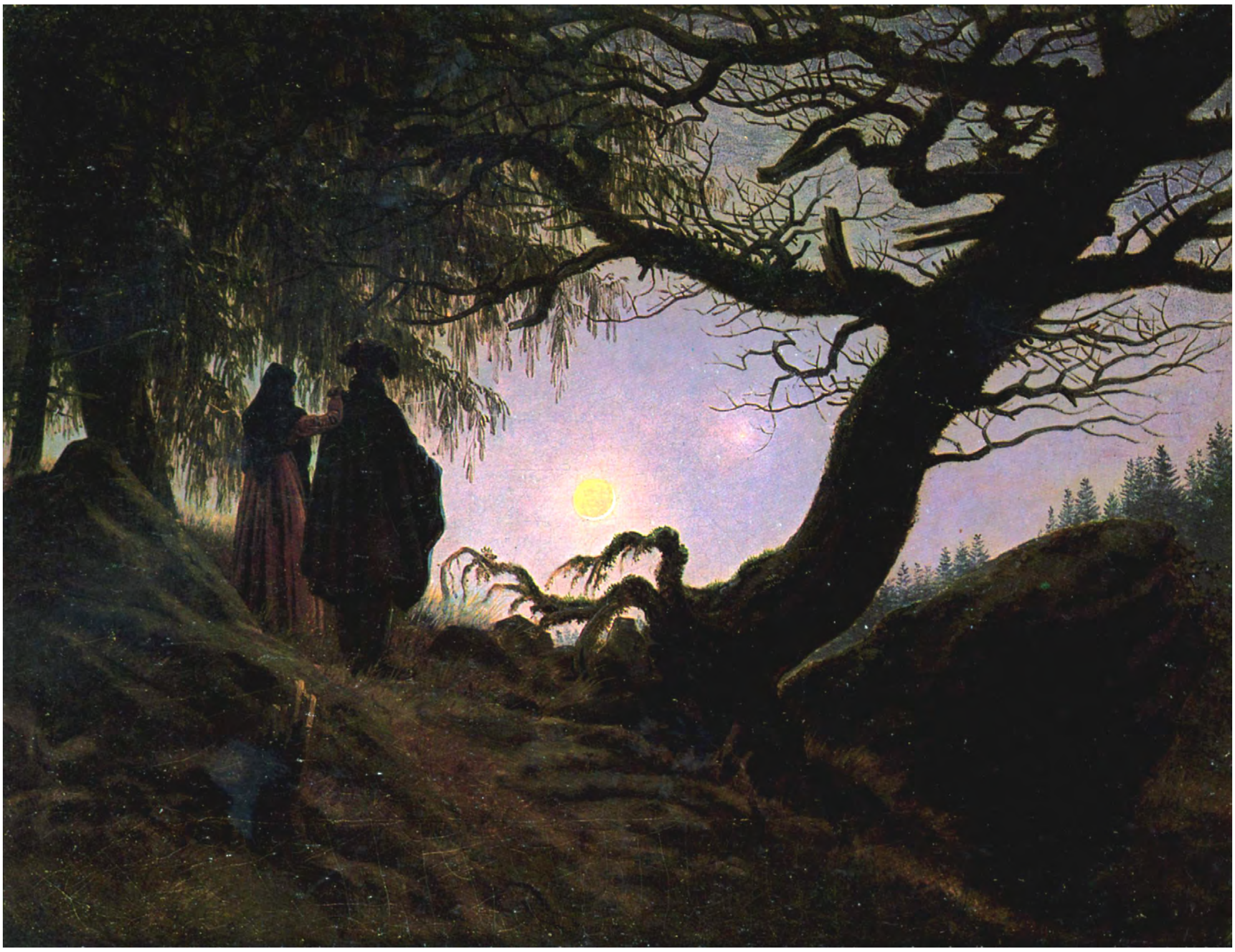

Caspar David Friedrich: Mann und Frau in Betrachtung des Mondes (Man and Woman contemplating the Moon) (1818-1824). ©

# *Acknowledgments*

Mathematics would be very sad if it could not be shared with others. Good friends and colleagues offered to read preliminary drafts of this book.

Thierry Barbot was one of my first PhD students and is now a renowned specialist in dynamical systems and lorentzian geometry. He spontaneously proposed to comment on an early version of the manuscript.

I met Grant Cairns long time ago when we were both students and we wrote two papers together. He is not specifically an expert on singularities but he has a wide interest in all aspects of mathematics. He likes promenades.

I hesitated to ask for the opinion of Pierre de la Harpe since I am well aware of his quality standards. His candid observations were fully appreciated.

Jos Leys is an artist and a geometer and we collaborated a lot, in particular in the production of mathematical movies. Many of the most sophisticated pictures are due to him.

Patrick Popescu-Pampu is one of the best experts of singularity theory, and loves the history of mathematics. We are both collaborators of Henri-Paul de Saint Gervais.

My colleague Bruno Sévennec knows everything in mathematics. He detected a large number of mistakes.

Christopher-Lloyd Simon agreed to play the role of the hypothetical "motivated undergraduate" trying to decipher these notes. The guinea-pig transformed into a collaborator. His suggestions improved the main results of the manuscript in a significant way.

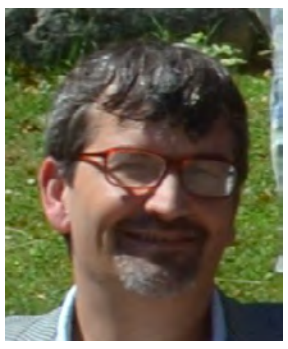

Thierry Barbot.

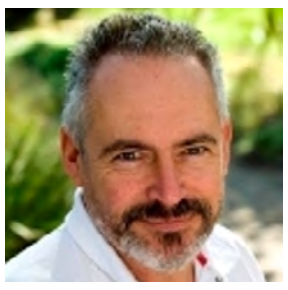

Grant Cairns.

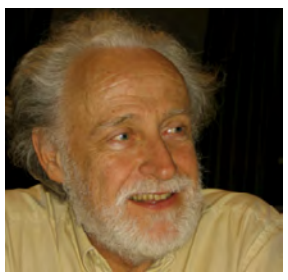

Pierre de la Harpe.

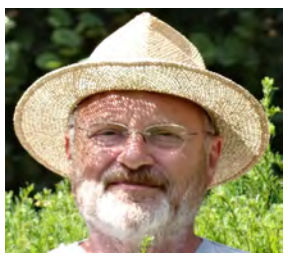

Jos Leys.

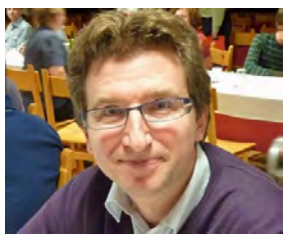

Patrick Popescu Pampu.

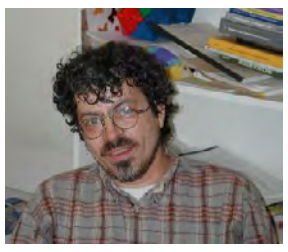

Bruno Sévennec.

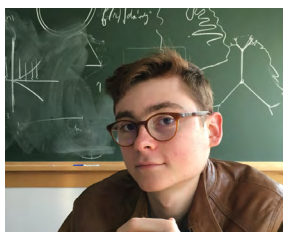

Christopher-Lloyd Simon.

A version of this text was available on *Arxiv* and on *AMS Open Math Notes* in November 2016. Since then two other former students examined the book in detail and sent me very interesting remarks.

Serge Cantat is usually dealing with algebraic geometry and dynamics in the complex domain.

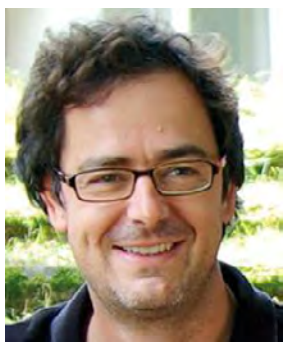

Serge Cantat.

Michele Triestino is more interested in probability, ergodic theory and theoretical physics, basically absent from this promenade!

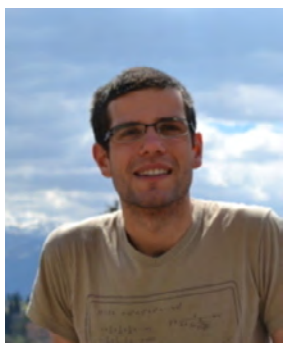

Michele Triestino.

I had the privilege to teach a graduate course on this topic in IMPA, Rio de Janeiro, in January and February 2017. Among the students who attended the classes, I would like to emphasize seven of them, who began a translation in Portuguese and Spanish: Aguiar Alina Sotolongo, Azul Fatalini, Carlos A. Gomes, Esteban Arreaga, Filipe Bellio da Nobrega, Gregory Cosac and Thiago Dourado. Dali Shen, a postdoctoral fellow at IMPA, also sent me a long list of remarks. I also received comments from Jean Barge, Brian Davey, Maxim Kontsevich, Manuel Ojanguren and Jean-Pierre Serre.

It is a great pleasure to thank them all.

I will not forget the anonymous reader who sent me chocolates because he liked the book. He identified himself as Paul Ynomial, rue du nœud de trèfle, Ker Lann.

I should probably not add the obvious comment that I am – of course – responsible for all remaining errors.

It is also my pleasure to thank three important institutions.

The *Centre National de la Recherche Scientifique* has been supporting me for many years.

The *Instituto de Matemática Pura e Aplicada*, Rio de Janeiro, provided an excellent scientific atmosphere and enough peace to write the first draft of this book.

The *Unité de Mathématiques Pures et Appliquées* de l'*École Normale Supérieure de Lyon* has been my home institution for a long time.

Étienne Ghys

Lyon, August 22, 2017

# *Image Credits*

The text and the illustrations without the symbol © have been produced by Étienne Ghys who has waived all copyright and related or neighboring rights. You can copy, modify, and distribute them, even for commercial purposes, all without asking permission.

The same *CC0* license applies to the images by Jos Leys on the following pages:

7, 7, 49, 112, 112, 116, 117, 124, 124, 142, 149, 152, 153, 153, 156, 158, 159, 176, 160, 161, 162, 163, 165, 170, 182, 179, 180, 180, 180, 181, 184, 185, 186, 186, 187, 200, 218, 219, 224, 227, 278, 283,

and the following documents from archive.org:

Gauss's dissertation: 6 — Colson's translation: 42 — Buffon's translation: 46 — Cramer's book: 63 — Cramer's book : 58 — Gauss figure: 72 — Serret's book: 84 — Serret's book: 85 — Divergentibus: 86 — Architecture book: 112 — Cable: 223 — Ashley: 288,

as well as the following from Wikipedia:

Wanderer: ii — Haeckel: 16 — Catalan: 21 — Catalan: 22 — Schroeder: 36 — Hipparchus: 36 — US signs: 37, 179 — Plutarch: 39 — Newton: 45 — Cramer: 61, — Dalembert: 77 — Euler: 87, — Hardy: 91 — Cauchy: 94 — Darwin finches: 109 — Moebius: 111 — Recycling: 117 — Olympic: 122, — Noether: 133 — Cremona: 137 — Paris: 147 — G. Doré: 155 — Puiseux: 169 — Weierstrass: 173 — Dürer: 190 — Cherry tree: 198 — Darwin: 214 — Epicycles: 224 — Gauss: 228 — Gauss's signature: 230 — Tree of crows: 292 — Mann und Frau: 294,

and

Four seasons: 9 — Tree of life: 10 — Mississipi: 34 — Circular Iching : 273 — Chord: 284.

The other illustrations are published under the Creative Commons licenses described below.

Signs: 3, 62, 87, 91, 134, 184, 199 — Knuth: 17, 19 — Mobile: 31 — : 32 — Cauldron: 118 — Cherry : 209 — Fish : 219 — IconPark : 242 — Nasa: 274.

Earth: 114 — Logo: 123 — Dimensions: 146 — Hopf: 149 — Flammarion: 154 — Geocentric: 154 — Pointe: 168 — Stasheff: 199.

I thank the following individuals and institutions for allowing me to use some of their documents

Maxim Kontsevich: 1 — Bib. Diderot ENS Lyon: 4 — Cambridge University library: 5, 45, 47, 48, 51, 52, 55, 57 —Museu de Arte da USP: 6, 120 — Frédérique Bassino: 33 — Pierre Gallais: 54 — U. Göttingen: 68, 283 — Scot Kolodziejski: 69, 275 — Institut de France: 78 — Gallica: 99 — Étienne Lecroart: 116 — Sylvie Pic: 118 — Laurent Bartholdi: 143, 119, 170 — Tadashi Tokieda: 125 — Ton Marar: 126 — Barry Sobel: 132 — Tony Philips: 148 — André Nachbin: 167 — Olivier Joseph: 175 — John Milnor: 177 — Rémy Oudompheng:193 — Robert Coquereaux: 197 — Eliane Loday: 197 —Anatoly Fomenko: 201, 202, 203, 204, 205, 206, 207 — Yifan Hu: 220, — Gary Urton: 245.

This book has been prepared using the "tufte-latex Document class" inspired by the work of Edward Tufte.

I limited the hyperlinks in the bibliography to freely accessible papers. I was told that almost all references are easily accessible on websites whose legality is questionable.

# A singular mathematical promenade

A stroll in the mathematical world. This is neither an elementary introduction to the theory of singularities, nor a specialized treatise containing many new theorems. The purpose of this little book is to invite the reader on a mathematical promenade. We pay a visit to Hipparchus, Newton and Gauss, but also to many contemporary mathematicians. We play with a bit of algebra, topology, geometry, complex analysis, combinatorics, and computer science. Hopefully motivated undergraduates and more advanced mathematicians will enjoy some of these panoramas.

Étienne Ghys is a CNRS senior researcher at École normale supérieure de Lyon.

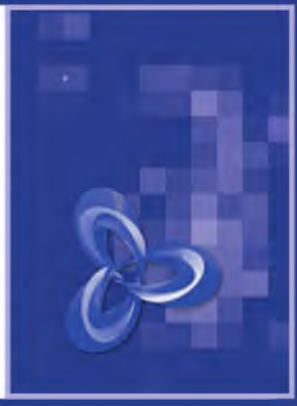